\documentclass[a4paper,10pt]{article}
\usepackage{amsmath}
\usepackage{caption}
\usepackage{graphicx} 
\usepackage{amsthm,amsfonts,amssymb,euscript,stmaryrd,graphics,color,amsmath,latexsym,marginnote, dsfont, math}

\usepackage[dvipsnames]{xcolor}
\usepackage[
    colorlinks=true,
    linkcolor=MidnightBlue,
    citecolor=BrickRed,
    urlcolor=RoyalBlue
]{hyperref}
\usepackage{enumitem}
\usepackage{todonotes}
\usepackage{thmtools}
\usepackage{hyperref}
\usepackage{cleveref}
\usepackage{tcolorbox}
\usepackage[left=2.2cm,right=2.2cm,top=2.2cm,bottom=2.2cm]{geometry}
\usepackage{cancel}
\usepackage[T1]{fontenc}
\usepackage{mathrsfs}
\newcommand{\Id}{{\rm Id}}
\usepackage{mathtools, upgreek}
\theoremstyle{plain}

\usepackage{amsmath}
\usepackage{accents, yfonts}
\usepackage{wrapfig}
\usepackage{soul}
\usepackage[intoc]{nomencl}  
\makenomenclature

\usepackage[T1]{fontenc}
\usepackage[utf8]{inputenc}
\usepackage[safe]{tipa}
\newcommand{\zetina}{{\text{\textctyogh}}}
\newcommand{\zeg}{{\text{\textctz}}}
\newcommand{\tipac}{{\text{\textctc}}}
\newcommand{\Lmap}{{\textswab{L}}}
\renewcommand{\eta}{\upeta}
\renewcommand{\psi}{\uppsi}

\newcommand{\zak}{{\mathfrak{Z}}}

\newcommand{\sign}{{\rm sign}}
\newcommand{\ov}{\overline}

\newcommand{\vr}{\varrho}
\newcommand{\vOmega}{{\bf \Omega}}
\newcommand{\opbw}{{{\rm Op}^{{\scriptscriptstyle{\mathrm BW}}}}}
\newcommand{\x }{\xi }
\renewcommand{\be}{\begin{equation}}
\newcommand{\ee}{\end{equation}}
\newcommand{\Opw}[1]{{{{\rm Op}^{{\scriptscriptstyle{\mathrm W}}}}\left(#1\right)}}
\newcommand{\Opbw}[1]{{{{\rm Op}^{{\scriptscriptstyle{\mathrm BW}}}}\left(#1\right)}}
\newcommand{\vOpbw}[1]{{{{\rm Op}_{\mathtt{vec}}^{{\scriptscriptstyle{\mathrm BW}}}}\left(#1\right)}}
\newcommand{\zOpbw}[1]{{{{\rm Op}_{\mathtt{out}}^{{\scriptscriptstyle{\mathrm BW}}}}\left(#1\right)}}
\newtheorem{theorem}{Theorem}[section]
\newtheorem{proposition}[theorem]{Proposition}
\newtheorem{lemma}[theorem]{Lemma}
\newtheorem{corollary}[theorem]{Corollary}
\newtheorem{definition}[theorem]{Definition}
\theoremstyle{definition}
\newtheorem{remark}[theorem]{Remark}

\renewcommand{\fA}{\mathsf{A}}

\renewcommand{\fC}{\mathfrak{C}}

\renewcommand{\fR}{\mathfrak{R}}
\renewcommand{\fP}{\mathfrak{P}}

\newcommand{\fs}{\mathfrak{s}}
\newcommand{\fv}{\mathfrak{v}}
\newcommand{\fc}{\mathfrak{c}}

\usepackage{bbm}

\newcommand{\intt}{{\rm int}}
\newcommand{\res}{{\rm res}}
\newcommand{\pperp}{\!\!\perp\!\!\!\!  \perp}

\newcommand{\scF}{\mathscr{F}}
\newcommand{\scG}{\mathscr{G}}
\newcommand{\scH}{\mathscr{H}}

\newcommand{\pare}[1]{( #1)}
\newcommand{\bra}[1]{[ #1]}
\newcommand{\sm}{\mathsf{m}}
\usepackage[normalem]{ulem}

\renewcommand{\bar}{\overline}

\newcommand{\bOmega}{{\boldsymbol{\Omega}}}

\newcommand{\bral}{[ \! [} 
\newcommand{\brar}{] \! ]}

\newcommand{\Set}{\Lambda}
\numberwithin{equation}{section}

\newcommand{\bPsi}{\mathbf{\Psi}}
\newcommand{\bUpsilon}{{\mathbf{\Upsilon}}}
\newcommand{\bXi}{\mathbf{\Xi}}
\newcommand{\wahlen}{\upzeta}
\newcommand{\Hilb}{\sfH}

\usepackage{microtype}

\usepackage{relsize}

\title{%
{\fontsize{14}{0}
\bfseries\scshape 
transfer of energy for  pure-gravity 
water waves \\ with constant vorticity}
}
\author{\fontsize{9.8}{12}\scshape beatrice langella\footnote{Università di Pisa, Department of Math, Largo Bruno Pontecorvo 5, 56127, Pisa, Italy\\ \textit{Email:} \texttt{beatrice.langella@unipi.it}}, \  alberto maspero\footnote{International School for Advanced Studies (SISSA), Via Bonomea 265, 34136, Trieste, Italy. \newline
	\textit{Emails:} \texttt{amaspero@sissa.it}, \texttt{shterra@sissa.it}}, \  federico murgante\footnote{	 Università degli studi di Milano, Department of Math. Federigo Enriques,
	Via Saldini 50, 20133 Milano, Italy \\ 
    and  \   Fields Institute, 222 College St, Toronto, Canada\newline 
	 \textit{Emails:} \texttt{federico.murgante@unimi.it},\texttt{
murgante@fieldsinstitute.ca}}, \ shulamit terracina$^\dagger$}
\date{}

\begin{document}

\maketitle

\begin{abstract}
We consider two-dimensional periodic gravity water waves with constant nonzero vorticity $\gamma$,  in  infinite depth and with periodic boundary conditions.
We prove that, if the
 characteristic wave number $\frac{\gamma^2}{g}$ is rational, 
 the system admits smooth small-amplitude solutions whose high Sobolev
 norms grow arbitrarily large while lower-order norms remain
 arbitrarily small,  thereby exhibiting  a genuine transfer of energy toward high
frequencies.
This yields
 the first rigorous construction of weakly turbulent solutions for a quasilinear  hydrodynamic wave system, in a regime where the flow remains smooth.
Moreover, the growth occurs simultaneously in the free surface and in the vertical component of the velocity at the interface, showing that the instability involves the full hydrodynamic evolution.

The proof relies on a new mechanism for generating energy cascades in quasilinear dispersive PDEs with sublinear dispersion and a nonlinear transport structure.
A central ingredient is to exploit quasi-resonances from 2-wave interactions to produce a transport operator that drives energy to high modes and causes Sobolev norm growth. A virial-type argument then shows that the resulting instability affects both the free surface elevation and the velocity field.


\end{abstract}

\begingroup
\hypersetup{linkcolor=black}
\tableofcontents
\endgroup
\section{Introduction}

The transition from smooth  motion to turbulence is one of the most ubiquitous and least understood phenomena in all of
fluid mechanics and mathematical physics. 
A central manifestation of this transition   is the energy cascade: 
the systematic transfer of energy from large spatial scales to small scales, or equivalently, from low to high Fourier modes, leading to the formation of increasingly fine structures. 
Analytically, a  well-established method to detect  such a transfer is to exhibit deterministic  solutions whose high Sobolev norms undergo arbitrarily large growth. 
In the context of {\em dispersive} partial differential equations,   such {\em weakly turbulent} solutions have been rigorously established only for a handful of semilinear models \cite{Kuk97.1,CKSTT, hani14, guardia_kaloshin, hani15, haus_procesi15, guardia_haus_procesi16, GHMMZ, Giu, GG, GS,
gerard_grellier0,gerard_grellier3, GL, GLPR, BE, BG}, leaving the problem completely open for the  quasilinear dispersive equations of classical hydrodynamics --which govern fluid motions where turbulent behavior is empirically observed.

Among hydrodynamical systems, the gravity water waves equations
occupy a distinguished position.
They describe the motion of the free interface between air and an inviscid,
incompressible fluid  under the action of gravity, 
and are the most accurate quasilinear PDE model  describing surface ocean waves. 
Pioneering formal computations by Hasselmann  \cite{Hasselmann1,Hasselmann2} and Zakharov
\cite{Zak65} in the 1960s identified resonant wave--wave interactions as the primary mechanism 
responsible for  energy cascades, providing  the basis for influential statistical theories of wave turbulence. Despite their impact on the physical understanding of ocean dynamics, 
and the recent progress on rigorous derivation of wave kinetic equations \cite{DH, DIP2}, 
the existence of  weakly turbulent solutions  in water waves  remains  open. 
This is the problem  we address in the present paper.

\smallskip
We consider the two-dimensional gravity water waves with constant vorticity $\gamma$ in infinite depth and periodic boundary conditions, under the action of gravity $g$. 
We prove that, whenever $\gamma^2/g \in \Q \setminus\{0\}$, the system admits smooth solutions that start with arbitrary small initial data and grow to arbitrary large size in high-Sobolev norm. 
The instability occurs entirely within the regime where the solution remains small in a low-Sobolev norm and the flow remains regular, thus it is not a byproduct of singularity formation or a loss of well-posedness, but a genuine 
   transfer of energy toward high frequencies within the full nonlinear dynamics. 
An important aspect of the result is that the instability manifests simultaneously in the free surface and in the vertical component of the velocity at the interface.

The difficulty in constructing weakly turbulent solutions for the water waves  is structural.
The water waves  are quasilinear PDEs with quadratic nonlinearities, 
and the cascade occurs on time scales substantially longer than
those provided by classical local well-posedness theory \cite{IFT}. This requires
a refined long-time analysis beyond standard 
methods, to guarantee the solution we construct {\em exists} over  enhanced time scales. 
Besides this, our approach departs from the by now classical methods in the literature on growth of Sobolev norms and is instead based on a new mechanism, which exploits the presence of 4-wave quasi-resonant interactions to single out a non-constant coefficients transport-type operator in the vector field. It is precisely this operator that drives the instability, effectively acting as a dispersive transport in frequency space. Once energy has reached higher frequencies, energy exchanges between the free surface and the velocity field provoke individual growth of all physically relevant quantities.
Let us now present precisely our result.

\smallskip
We study  the evolution of a two-dimensional
perfect and incompressible 
fluid with constant vorticity 
$\gamma \neq 0 $, under the action of gravity, which 
fills  the time dependent region 
$$ \cD_{\eta} := \big\{ (x,y)\in \T\times \R \ : \   y<\eta(t,x) \big\} \ , \quad   \T := \T_x :=\R/ (2\pi\Z)  \ , $$
with infinite depth  
and space periodic boundary conditions. 
The unknowns 
are the  free surface  
$$
\Gamma(t):= \{(x,y)\in \T \times \R \, \colon  y = \eta (t, x) \} $$
of $\cD_{\eta} $ and the incompressible 
 velocity field
\be\label{u.campo}
\vec{v} := \vec{v}(t,x,y):= \begin{pmatrix} v_1(t,x,y) \\ v_2(t,x,y) \end{pmatrix}\,,
\ee
evolving according to the  Euler equations
\be\label{euler}
\partial_t \vec v + \vec v \cdot \grad \vec v = - \grad P- g \vec{e}_y  \ , \quad 
{\rm div }\,  \vec v  = 0   \ , \quad 
{\rm rot }\,  \vec  v  = \gamma, \qquad \text{in } \cD_{\eta},
\ee
with $\vec{e}_y := \vect{0}{1}$. At the free surface, we impose   the kinematic and dynamic  boundary conditions
\be\label{kinematic}
\eta_t  = v_2  - v_1\, \eta_x \ , \qquad P = P_0  \quad \mbox{  at } y =\eta(t, x) \ , 
\ee 
that ensure that the fluid particles at the free surface remain on it along the evolution and that the pressure $P$ 
of the fluid at the free surface  equals 
the constant atmospheric pressure. At the bottom, we assume the boundary condition
\be\label{bottom}
v_2 \to 0 \ \ \  \mbox{ as } y \to - \infty \ , 
\ee
that imposes that the fluid  is horizontal at infinite depth.

For any $ s \in \R $, we  denote by 
$H^s(\T;\R)$ the Banach space of real-valued periodic Sobolev functions endowed with the norm  
$$
\norm{u}_{s} :=
\Big( \sum_{j \in \Z }  |u_j  |^2 \langle j \rangle^{2s} \Big)^{\frac12}  \,, 
$$
where  $ u_j  : = \frac{1}{2\pi} \int_{\T} u(x) e^{-\ii j x} \di x$ are the Fourier coefficients of $ u(x) $, and $\langle j \rangle := \max\{ 1, j\}$. 
Similarly we denote by $H^s(\T;\R^2)$ periodic functions with values in $\R^2$ and whose components have finite $\norm{\cdot}_s$ norm and with an abuse of notation, when $\vec{u}$ is in $H^s(\T; \R^2)$ we still denote by $\| \vec{u}\|_s$ the sum of the $\| \cdot\|_s$ norms of its components.

We begin by stating a simplified version of our main result,  focused on the velocity field at the free surface; 
we denote  by $\vec v(t,\cdot)\vert_{\Gamma(t)} := \vec v(t,x, \eta(t,x))$ the restriction of $\vec v(t,x,y)$ at the free surface $\Gamma(t)$, where $\vec{v}(t, \cdot)$ and $\eta(t, \cdot)$ are the solutions of the system \eqref{euler}--\eqref{bottom}.

\begin{theorem}\label{thm:main.intro}
Assume that the vorticity $\gamma$ in \eqref{euler} fulfills
\be \label{resonance.intro}
\frac{\gamma^2}{g} \in \Q  \setminus\{0\} \,. 
\ee
There exists $s_0 >\frac32$ such that for any $s \gg s_0$, any $0<\delta \ll 1$ and $K \gg 1$, 
there exist a  time  $T>0$ and 
 continuous in time classical solutions $\eta(t,x), \, \vec v(t,x, y)$  of \eqref{euler}--\eqref{bottom}, 
such that at the initial time 
    \be\label{explosion.intro}
    \| \eta(0,\cdot)\|_{{s+\frac12}} + 
      \| \vec v(0,\cdot)\vert_{\Gamma(0)}\|_{{s}} 
      \leq \delta  \,, 
      \ee
while, before the time $T$,  the free surface and the vertical component of the velocity field at the free surface undergo    Sobolev norm growth      
      \be\label{eta.u.intro}
    \sup_{t \in [0,T]}  \| \eta(t, \cdot)\|_{{s+\frac12}} \geq K   \,, 
     \qquad  
     \sup_{t \in [0,T]}  \norm{\vec{e}_y \cdot \vec v(t, \cdot)\vert_{\Gamma(t)} }_{{s}} \geq K  \,.
    \ee
    In addition, one has  the low-Sobolev norm  bound
    \be \label{small.in.small.intro}
 \| \eta(t,\cdot)\|_{{s_0+\frac12}} + 
      \| \vec v(t,\cdot)\vert_{\Gamma(t)}\|_{{s_0}} 
      \leq 2\delta   \,, \quad \forall \, t \in [0,T] \,. 
    \ee
\end{theorem}
This theorem is an immediate consequence of our main result   \Cref{thm:mainimprove}, stated  in the Zakharov variables,  and it will be proved immediately after it.

Condition \eqref{small.in.small.intro} ensures that the solution $(\eta(t), \vec{v}(t)\vert_{\Gamma(t)})$ 
exists and remains small over the time interval $[0, T]$ in the low $H^{s_0 + \frac 1 2}(\T) \times H^{s_0}(\T)$ norm, whereas, as guaranteed by the  inequalities in \cref{eta.u.intro}, by time $T$ the high Sobolev norms of the free surface and of the vertical component of the velocity field at the free surface become arbitrarily large. 
This implies that a genuine transfer of energy from low to high Fourier modes occurs over the time interval $[0, T]$.
Additionally, the smallness condition \eqref{small.in.small.intro} ensures the persistence of the Taylor sign condition
\be\label{taylor.cond}
- \nabla_{\vec n} P\big|_{\partial \mathcal D_\eta} \ge c_0 >0\,,
\ee
with $\vec n$  the outward unit normal to the fluid domain, 
which is a fundamental prerequisite for the local well-posedness of the free-boundary Euler equations \cite{CL00, Lind1, SZ08}.
While the sign condition \eqref{taylor.cond} 
is automatically verified  for  irrotational gravity waves with  flat or no  bottom  \cite{Wu0}, the rotational case  admits the possibility of a sign violation, producing   singular dynamics \cite{Ebin, Wu19, Su23}.
In our situation  (cf. \eqref{u.campo} and  \Cref{rem:taylor}) 
\be\label{taylor.explicit}
- \nabla_{\vec n} P\big|_{\partial \mathcal D_\eta} = \sqrt{1 + \eta_x^2} (1 + \sfa)\,, \qquad \sfa(t, x) := (\pa_t + v_1(t, x, \eta(t, x))\pa_x) v_2(t, x, \eta(t, x))\,,
\ee
so  condition \eqref{small.in.small.intro} ensures that the Taylor coefficient $\sfa$ remains  uniformly small  along the whole timescale considered in \Cref{thm:main.intro}.
Therefore,  the observed Sobolev norm inflation in \eqref{eta.u.intro} is a feature of the smooth nonlinear dynamics rather than a consequence of vorticity-induced singularity formation.
 Further comments will be made after \Cref{thm:mainimprove}.

\smallskip
The proof of Theorem \ref{thm:main.intro} introduces a novel paradigm for inducing energy transfers that is fundamentally distinct from the previous established mechanisms currently used in the study of semilinear Hamiltonian PDEs.
The first one, pioneered by Colliander-Keel-Staffilani-Takaoka-Tao \cite{CKSTT} for the cubic nonlinear Schr\"odinger equation on $\mathbb{T}^2$, relies on the construction of a finite-dimensional ``toy-model'' to drive energy to high-frequencies through a sequence of resonant clusters.
This approach has been  refined and extended  in  \cite{hani14, guardia_kaloshin, hani15, haus_procesi15, guardia_haus_procesi16, GHMMZ, Giu, GG, GS}. 
However, in all these works, the construction of weakly turbulent solutions has relied  on embedding the I-team toy-model \cite{CKSTT} into the nonlinear dynamics. While powerful, this mechanism is very rigid and so far could not be applied to any equation on tori which is not NLS-like (even for the quantum Euler equations \cite{GS},  the instability is established by mapping the system onto a semiclassical cubic NLS via the Madelung transform and subsequently exploiting the toy-model of \cite{CKSTT}).\\
The second mechanism,  pioneered by Gérard and Grellier \cite{gerard_grellier0} in the study of the Szegő equation, leverages certain  integrable  PDEs to characterize explicit families of unstable invariant objects—including quasi-periodic solutions \cite{gerard_grellier0, gerard_grellier3}, multi-solitons \cite{GL, GLPR}, and finite-dimensional resonant manifolds \cite{BE, BG}. These special solutions serve as the dynamical skeleton for constructing weakly turbulent trajectories. 
We also mention the unbounded growth of the energy density  proved for the integrable  binormal flow  by Banica-Vega \cite{BV}.
While mathematically profound, all these methods are intrinsically tethered to the integrable structure of the equations and, 
consequently, they provide little insight into the quasilinear dynamics of  water waves, which requires  a completely different analysis. 
\\
A first step toward generating Sobolev norm growth in a quasilinear setting was recently taken in \cite{MM} for a Schr\"odinger-Burgers equation. 
However, the analysis in \cite{MM} relied on a tailored, cubic nonlinearity specifically designed to enforce resonant interactions responsible for the cascade. 
While that model provided a conceptual proof-of-principle for Sobolev norm growth in a quasilinear system, it left open the fundamental question of whether such mechanisms are possible  in genuine, physical fluid equations with not engineered  nonlinearities.
In this paper we resolve this question by  retaining the full quasilinear nonlinearity of the free boundary Euler equations. 
Moreover, we establish a general and robust procedure that we believe can be applied to prove  growth of Sobolev norms for quasilinear PDEs  with sublinear dispersion and a nonlinear transport term.
We now present our main theorem.

\subsection{Main theorem}
We study the water waves problem in the 
\textit{Hamiltonian}  formulation 
\cite{Zak1,CrSu}, extended in
\cite{CIP} and 
\cite{Wh}
for constant vorticity fluids.  
In view of our incompressibility and constant vorticity assumptions \eqref{euler}, the velocity field $\vec v$  is the sum of the Couette flow $\begin{psmallmatrix} - \gamma y \\ 0 \end{psmallmatrix}$, which carries
all the  vorticity $ \gamma $, 
and an irrotational field, 
expressed as the gradient of a harmonic function $\Phi $, called the generalized velocity potential. 
Denoting by 
$\psi(t,x)$  the evaluation of the generalized velocity potential at the free interface
$ \psi (t,x) := \Phi (t,x, \eta(t,x)) $, one recovers
$ \Phi $ as 
the unique harmonic function 
$ \Delta \Phi = 0 $ in $ \cD_{\eta} $ with Dirichlet boundary condition
$ \Phi = \psi $ at $  y = \eta(t,x) $ and Neumann boundary condition
$\Phi_y ( t, x, y) \to  0  $ as $ y \to  - \infty $.  
Imposing  
the kinematic and dynamic boundary conditions \eqref{kinematic}, the 
evolution of the fluid \eqref{euler} is determined by the non-local quasi-linear 
equations 
\begin{equation} \label{eq:etapsi}
\left\{\begin{aligned}
 &   \partial_t \eta = G(\eta)\psi+ \gamma \, \eta \eta_x \\
&\partial_t\psi =  -g\,\eta  -\frac{1}{2} \psi_x^2 
+  \frac{1}{2}\frac{(\eta_x  \psi_x + G(\eta)\psi)^2}{1+\eta_x^2}
+ \gamma \, \eta \psi_x+ \gamma \partial_x^{-1} G(\eta) \psi 
\end{aligned}\right. \,, 
\end{equation}
where  $ g > 0 $ is the gravity constant 
and 
$ G(\eta)$ is the  Dirichlet-Neumann operator  
\begin{align}\label{DN}
G(\eta)\psi  := 
(- \Phi_x \eta_x + \Phi_y)\vert_{y = \eta(x)}  \,.
\end{align}
The equations \eqref{eq:etapsi} have a Hamiltonian structure, but not with respect to the standard symplectic form. This will be discussed in Section \ref{sec:basic},  where the Wahlén coordinates \cite{Wh}, which provide Darboux coordinates, will be introduced.
The quantity $\int_\T \eta(x) \, \di x$ is a prime integral of \eqref{eq:etapsi}
and then, with no loss of generality,  we restrict to interfaces  with zero average,
\be\label{zero.av} 
\int_\T \eta(x) \, \di x = 0  \,. 
\ee
The component $\eta $  of the solution of \eqref{eq:etapsi}
 will lie in a Sobolev space $H_0^{s+\frac{1}{2}}(\T)$ of periodic functions with zero mean. 
Moreover, the vector field on the right hand side of  \eqref{eq:etapsi} 
depends only on $ \eta $ and $  \psi - \frac{1}{2 \pi}\int_\T \psi(x) \, \di x  $
and therefore  $ \psi $ 
 will evolve in a 
homogeneous Sobolev space $ \dot H^{{s+\frac{1}{2}}}(\T)$ of periodic functions modulo constants. 
We also denote the horizontal and vertical components of the irrotational velocity field 
at the free interface by
\begin{align} 
\label{def:V}
& \sfV =  \sfV (\eta, \psi) :=  (\pa_x \Phi) (x, \eta(x)) = \psi_x - \eta_x \sfB \,,
\\
\label{form-of-B}
& \sfB =  \sfB(\eta, \psi) := (\pa_y \Phi) (x, \eta(x)) =  \frac{G(\eta) \psi + \eta_x \psi_x}{ 1 + \eta_x^2} \,. 
\end{align}
\smallskip
Note that, by \eqref{def:V} and \eqref{form-of-B}, the Dirichlet-Neumann operator in  \eqref{DN} is 
\begin{align}
    G(\eta)\psi= \sfB- \sfV \eta_x\,.
    \label{espr_DN}
\end{align}
For any $ s \in \R $,  we equip  the Sobolev spaces 
$$
 H^s_0 (\T; \C) := \Big\{ u(x) \in H^s (\T;\C) \, : \, \int_{\T} u(x) \, \di x = 0  \Big\} \, ,
\quad \dot H^s (\T; \C) = H^s (\T; \C) \slash \C \, ,
$$
with the norm 
$$
\| u \|_{H^s_0} := \| u \|_{\dot H^s} =
\Big( \sum_{j \in \Z\setminus\{0\} }  |u_j  |^2 | j |^{2s} \Big)^{\frac12}  \ . 
$$  
With an abuse of notation we will still use $\| \cdot \|_{s}$ instead of $\| \cdot \|_{H^s_0} $.
The quotient map induces 
an isometry between $ H^s_0 $ and $ \dot H^s $ and we shall often
 identify 
$ H^s_0 $ with $  \dot H^s $.

Following \cite{ABZ1}, we also  introduce the spaces
\be\label{Xs}
X^s:= H^{s+\frac12}_0(\T;\R) \times \dot H^{s+\frac12}(\T;\R) \times H^{s}(\T;\R) \times H^{s}(\T;\R)
\ee
endowed with the norms
\begin{equation}\label{norm:Xs}
    \| (\eta,\psi,\sfV,\sfB)\|_{X^s}:= \| \eta\|_{s+\frac12}+\| \psi\|_{s+\frac12}+\| \sfV\|_{s}+\| \sfB\|_{s}\,.
\end{equation}
 
The main result of the manuscript is the following one. 

\begin{theorem}
\label{thm:mainimprove}
    Let  $\gamma$ fulfill \eqref{resonance.intro}. There exists $\ts_0 > \frac32$  such that the following holds true.
    Let  $\ts> 5\ts_0$,  $0<\delta < 1$, $K >  1$   be given parameters. There exist  a time  $T>0$ and  a solution $(\eta, \psi, \sfV, \sfB) \in C^0([0,T]; X^\ts)$  of \eqref{eq:etapsi} 
    such that  at the initial time
    \be\label{initial.intro}
\norm{\eta(0)}_{\ts +\frac12} + \norm{\psi(0)}_{\ts +\frac12} + 
\norm{\sfV(0)}_\ts +\norm{\sfB(0)}_\ts  \leq \delta
    \ee
    whereas in the time interval $[0,\,T]$   one has 
\begin{equation}\label{explosion.improved}
        \sup_{t \in [0, T]} \norm{\eta(t)}_{\ts+\frac12} \geq K \,, \qquad 
        \sup_{t \in [0, T]}\norm{\sfV(t)}_\ts \geq K \,, \qquad 
        \sup_{t \in [0, T]}\norm{\sfB(t)}_\ts \geq K  \ . 
    \end{equation}
Finally, one has the low Sobolev norm bound
\be\label{low.norm.intro0}
        \sup_{t \in [0, T]} \| (\eta,\psi,\sfV,\sfB)(t)\|_{X^{\ts_0}} \leq 2 \delta\,.
\ee
\end{theorem}
Let us make some comments.

\medskip 

\noindent {\sc 1. \und{Simultaneous growth of the physical variables:}}
A crucial step towards the proof of  \eqref{explosion.improved} is to construct a 
 solution $(\eta, \psi, \sfV, \sfB)(t)$ of  \eqref{eq:etapsi} for which {\em at least one} variable among $\eta, \sfV, \sfB$ undergoes growth of Sobolev norms, i.e.  together with \eqref{low.norm.intro0},  it fulfills
\be\label{growthXs.intro}
   \| (\eta,\psi,\sfV,\sfB)(0)\|_{X^\ts} \leq \delta  \,, 
      \qquad 
\sup_{[0,\,T]} \| (\eta,\psi,\sfV,\sfB)(t)\|_{X^\ts} \geq K \,. 
       \ee       
We defer the discussion of how to prove  \eqref{growthXs.intro} to points 3--5 of this list, which rely on a different dynamical mechanism, and now provide an heuristic  explanation of how to exploit \eqref{growthXs.intro} and the wave-like structure of the water waves
to deduce the stronger estimate \eqref{explosion.improved}, that  proves that {\em each} variable actually grows.
The point is that in the complex variable 
$$
u \approx \frac{1}{\sqrt2}|D|^{-\frac14} \eta +  \frac{\im }{\sqrt2} |D|^{\frac14} \upomega ,  \quad (\mbox{cf. } \eqref{u.def})
$$
where $\upomega$ is Alinhac's good unknown \eqref{GU},
the water waves \eqref{eq:etapsi} take the form 
\be\label{linear.intro}
\pa_t u= - \im \Omega(D) u + \cN(u)  \ , 
\quad \Omega(D) = \sqrt{|D| + \frac{\gamma^2}{4}} + \im \frac{\gamma}{2} \Hilb \ , 
\ee
 $\Hilb$ being  the Hilbert transform in \eqref{def:Hilbert} and $\cN(u)$ a quadratic quasilinear nonlinearity (cf. \Cref{LemCompl}).
Roughly speaking, estimate \eqref{growthXs.intro} implies that $\|u(t)\|_{\ts+\frac34}$ becomes large, and hence at least one of $\|\Re\,u(t)\|_{\ts+\frac34}$ or $\|\Im\,u(t)\|_{\ts+\frac34}$ is large. By the norm equivalences $\|\Re\,u(t)\|_{\ts+\frac34}\approx \|\eta(t)\|_{\ts+\frac12}$ and $\|\Im\,u(t)\|_{\ts+\frac34}\approx \|\upomega(t)\|_{\ts+1}$, this corresponds to the growth of at least one of the two components. To deduce \eqref{explosion.improved}, it remains to show that {\em both} quantities become large. Heuristically, this follows from
 the oscillatory wave-like nature of \eqref{linear.intro}, which implies that $\Re\, u(t)$ and $\Im\,  u(t)$ are coupled and exchange energy over short amount of time.
 However, to make this heuristic rigorous is far from trivial,  since the 
full system is quasilinear and the nonlinear terms cannot be treated perturbatively, even on short time scales. Motivated by this picture, we instead adopt an energy-based approach, based on an “upside-down” virial argument, which captures the energy transfer indirectly by exploiting lower bounds on the second time derivatives of the norms of $\eta(t)$ and $\upomega(t)$, and that we describe better in \Cref{sec:scheme_proof}.

\medskip 

\noindent {\sc 2. \und{Enhanced time-scale:}}
A fundamental challenge in constructing  weakly turbulent solutions  for quasilinear systems is the disparity between the timescale required for energy transfer and the lifespan guaranteed by local theory.  Let us explain such disparity in our context.
We choose $0 < \epsilon \ll 1$ such that $K^{-1} \leq \epsilon^{\theta} \leq \delta$ for some small, positive $\theta$ (see \eqref{s0r} for the precise value), and we take a small-amplitude initial data of size
   \be\label{growthX.intro}
    \epsilon:=  \| (\eta,\psi,\sfV,\sfB)(0)\|_{X^{{\ts}_0}} \ll \delta, \qquad \| (\eta,\psi,\sfV,\sfB)(0)\|_{X^{{\ts}}} \leq \delta  \ . 
    \ee
 Standard local well-posedness results, combined with a modified energy method \cite{IFT,Wan}, ensure that the solution exists at least up  to the time  $T_{\rm wp} \sim \epsilon^{-2}$ over which the uniform bound $\sup_{t \in [0, T_{\rm wp}]} \| (\eta, \psi, \sfV, \sfB)(t) \|_{X^{\ts}} \lesssim \delta
   $ holds
   (actually \cite{IFT,Wan} adopts holomorphic coordinates, so we briefly sketch the cubic lifespan result in \Cref{sec:cubic} directly in the Zakharov variables). Consequently, no significant growth of Sobolev norms can manifest within time $T_{\rm wp}$.
Indeed, we prove that the instability in Theorem \ref{thm:mainimprove} develops on the logarithmically extended timescale
$$
T \sim \epsilon^{-2} \log(\epsilon^{-1}) \gtrsim \epsilon^{-2} \log(\tfrac{K}{\delta}) \ , $$ 
cf. \Cref{rem:enhanced}.
To bridge the gap between the time $T_{\rm wp}$ and the enhanced time scale $T$, we do not attempt to show that arbitrary initial data satisfying \eqref{growthX.intro} generate solutions existing up to time $T$--  which may be false. Instead, we restrict to a class of special initial data, which we call \emph{well-prepared}, whose $L^2$ mass is mostly concentrated on two low-frequency modes that generate only integrable $4$-wave resonances with all remaining modes. The existence of such modes relies on subtle arithmetic conditions, which we establish using number-theoretic arguments.

\medskip 

\noindent {\sc 3. \und{Resonance condition:}} 
The quantity $\gamma^2/g$, which has the dimension of an inverse length, determines the characteristic wave number arising from the balance between vorticity and gravity. 
Mathematically, the rationality assumption \eqref{resonance.intro} of $\gamma^2/g$ 
produces non-integrable $2$-wave resonances of the form
\be\label{2wave.intro}
\Omega_\tm(\gamma) = \Omega_\tn(\gamma) \ , \quad \tm < 0 < \tn \ , \quad \tm + \tn > 0   \ , \quad  \tm, \tn \in \Z\setminus\{0\}
\ee
among the linear frequencies of the system \eqref{eq:etapsi}, explicitly given by
\be\label{freq.intro}
\Omega_j(\gamma):= \sqrt{ g |j|+  \tfrac{\gamma^2}{4} } +  \tfrac{ \gamma}{2} \sign(j) \, ,  \quad j \in \Z\setminus\{0\}  \ . 
\ee
The  resonances \eqref{2wave.intro} have substantial dynamical consequences. We prove the existence of a complex normal form variable $z$ in which system \eqref{eq:etapsi} can be written as
	 \be \label{Z.eq.intro}
	 \begin{aligned} 
	 	&\pa_t z= 	 	- \im \Omega(D)z +
		 \Opbw{ \im  \langle \,\sfV\, \rangle (z;x) \xi}z +   \, \mbox{Remainder } 
	 \end{aligned}
	 \ee
where $\Omega(D)$ denotes the Fourier multiplier associated with \eqref{freq.intro}, the function $\langle \sfV\rangle(z;x)$ is real-valued and depends non trivially on the space variable $x$ (see \eqref{VresZ}), and the remainder consists of perturbative terms of higher order or smoother type.
The crucial feature of \eqref{Z.eq.intro} is the presence of the non-constant coefficient transport operator. This sharply contrasts with the normal forms obtained in \cite{BD, BFP, BMM2, MMS}, where the transport symbol has constant coefficients (at least at leading order). It also differs significantly from the normal form in \cite{CKSTT}. Indeed, the nonlinear vector field in \eqref{Z.eq.intro} is not Birkhoff-resonant: the  term 
$\Opbw{ \im  \langle \,\sfV\, \rangle (z;x) }z $ has phases of oscillations  given by
$$ 
\Omega_j(\gamma) - \Omega_{j +n-m}(\gamma) , \quad  j \in \Z\setminus\{0\}\ ,   \  \, 
\quad \forall \, m,n \quad \mbox{s.t.}\quad  \Omega_m(\gamma) = \Omega_n(\gamma) \ .
$$  
Formally, these terms could be removed by a Birkhoff normal form transformation. However, because of the quasi-linear nature of the equation, the required change of variables would be unbounded and ill-defined in Sobolev spaces. As a result, one encounters quasi-resonances rather than removable non-resonant interactions. \\
 In previous works on KAM theory \cite{BM, BBHM, BFM1, BFM2, FG, FMT} or long-time Cauchy stability \cite{BD, BMM, DIP, MMS, MRS}, the authors  have done ``whatever it takes'' 
 to put the transport term to {\em constant coefficients}
 (e.g. restricting to Cantor sets of parameters, 
 selecting the support of the initial datum 
 or   exploiting structural cancellations).
In our framework, by contrast, the non-constant coefficient transport is not an obstruction to be removed. Instead, we exploit its spectral properties as the driving mechanism of instability, as will be explained below.

It is worth noting that for non-rational values of $\gamma^2/g$
one expects a  different dynamical behavior characterized by enhanced stability. 
For instance, \cite{BFM2} proves that for a Borel set of vorticities $\gamma$ with asymptotically full Lebesgue measure, system \eqref{eq:etapsi} admits families of small-amplitude, time-quasi-periodic traveling wave solutions. These KAM solutions form islands of stability in phase space: they exist globally in time and their Sobolev norms remain uniformly bounded.

\medskip 

\noindent{\sc 4. \und{The initial datum:}} 
In the normal form variable $z$, the initial datum generating weakly turbulent solutions is supported on exactly four Fourier modes. Two of them belong to a set $\Lambda=\{\tm,\tn\}$ that supports a $2$-wave resonance of the form \eqref{2wave.intro}.
The remaining two modes lie at much larger frequencies, well separated from $\Lambda$.
The amplitudes and phases of these four modes are required to satisfy a number of structural conditions, detailed in \Cref{B}, that we shall explain better in the next section. A concrete example of such an initial datum is constructed in \eqref{z0}. 
Let us remark that the choice of the initial datum is not unique: for any $\gamma$ satisfying the resonance condition \eqref{resonance.intro}, we can exhibit infinitely many different choices of the set $\Lambda$. This, together with the fact that the structural conditions of Definition \ref{B} are open, amounts to infinitely many disjoint open neighborhoods of initial data leading to arbitrarily large growth of Sobolev norms.

\medskip 

\noindent {\sc 5. \und{The unstable transport operator:}}
\label{item:comment_transport}
The transport operator $\Opbw{\ii \langle \sfV \rangle(z; x)\xi)}$ appearing in \eqref{Z.eq.intro}
is at the core of our mechanism for  growth of Sobolev norms.
Indeed, we show that long-time controlled solutions are characterized, up to small remainders, by an oscillatory dynamics of the modes $z_\tm(t),\ z_\tn(t)$ supported on the set $\Set$, while on all other components of the solution, which remain smaller in a low-Sobolev norm, $\Opbw{\ii \langle \sfV \rangle(z; x)\xi)}$ is well approximated by the time independent, space dependent transport operator 
 \begin{equation}
  \Opbw{ \big(\tJ+\fv(x)\big) \xi}  \,, \quad   \tJ+\fv(x )=  \tV_\tm |z_\tm(0)|^2 + \tV_\tn |z_\tn(0)|^2 +
 2 \sfV^{(\rm res)}_{\tm, \tn} \, \Re \left(z_\tn(0) \, \bar{z_{\tm}(0)} \, e^{\im (\tn-\tm) x} \right)  \,
  \label{fH.intro}
\end{equation}
 where the coefficients $\tV_\tm$ and $\tV_\tn$ are defined in \eqref{VmVn.def} and $\sfV_{\tm, \tn}^{(\text{res})}$ in \eqref{Vres.coeff}, and depend only on the set $\Lambda$ and  the initial amplitudes of two Fourier modes of the initial data supported on  $\Lambda$.\\
The key mechanism generating  energy cascades lies in the spectral properties of such operator. 
If the smooth function $\tJ + \fv(x)$ changes sign -- a condition we are able to impose by carefully selecting the set $\Lambda$ and $z_\tm(0), z_\tn(0)$ -- 
 this operator possesses a non-trivial absolutely continuous spectrum, a feature that prevents the linearized waves from undergoing simple oscillations. Instead, it  creates  a direct pathway for energy to migrate toward higher and higher frequencies. 
 To rigorously prove quantitative Sobolev norm growth, we take inspiration from  the classical Mourre commutator theory \cite{Mourre}, and its application in growth of Sobolev norms for linear Schr\"odinger equations initiated in \cite{Mas21,Mas22}, and extend it to the nonlinear setting. 
 The no-sign condition on $\tJ + \fv(x)$ allows us to construct a conjugate operator $\fA$ -- supported on high frequencies -- that serves as an ``escape function'' for the dynamics, whose  commutator
$
\im [\fA,\Opw{(\tJ +\fv(x))\xi}]$
is strictly positive on large frequencies.
This ultimately proves an exponentially fast growth of high-Sobolev norms (that saturates the upper bounds given by the linear evolution \cite{MaRo}).

\medskip 

\noindent{\sc 6. \und{Comparison with  normal forms for water waves:}} 
Quasilinear normal form methods have recently emerged as a powerful tool in the study of the long-time dynamics of water waves. 
They have led to global/almost global existence results in the Euclidean setting \cite{GMS,GMS2, Wu, Wu11, IP15, AlDe1, IP18, IFRT, ifrTat, DIPP}, as well as long-time existence results on the torus \cite{BD, BFP, BMM2, IoP, DIP}.
In \Cref{sec:paraWW} we also perform a paradifferential normal form to study the long-time dynamics. 
Although we follow the approach  of \cite{BFP},  there are several significant differences that we highlight now.

First, in the presence of constant nonzero vorticity, $4$-wave resonances 
generate a cubic non-integrable normal form vector field, in contrast with the irrotational case of \cite{ZakD, BFP}, where the cubic normal form is integrable. This non-integrability obstructs the propagation of low-norm smallness in \eqref{small.in.small.intro}.
For this reason we give up to control {\em any} small amplitude solution, and we restrict to those whose initial data
have most of their $L^2$ mass 
concentrated on two low-frequency modes that generate only integrable $4$-wave resonances with all remaining modes.
 The existence of such modes relies on subtle arithmetic conditions, which we establish using number-theoretic arguments.
Exploiting this property and the Hamiltonian structure, we propagate \eqref{low.norm.intro0} on enhanced times. 

Second, unlike \cite{BFP},  in our case the  formal Birkhoff normal form of  \Cref{prop:cH^{(4)}_+}  does not coincide with the rigorous
paradifferential normal form of \Cref{thm:nf}. The difference is exactly  the transport-like term in \eqref{Z.eq.intro} which is absent from the formal normal form but appears in the rigorous one due to $4$-wave quasi-resonances. Hence the   formal evolution is qualitatively different from the  true one, in contrast with  \cite{BFP}.

Third, unlike \cite{BFP, IFT, BD, BMM2}  that rely on modified energy methods,  our normal form transformations are  genuine {\em changes of coordinates} in phase space, achieved via the theory of {\em admissible transformations} in \Cref{sez:admissible}.  This change is essential because the initial data leading to growth are naturally defined in the normal form variables.

\medskip 

\noindent{\sc 7. \und{Dynamics after the growth:}} 
The long-time asymptotic behavior of the solutions constructed in Theorem \ref{thm:mainimprove} beyond the timescale $T$ remains a profound open question. Given that global-in-time existence for the  gravity water waves  \eqref{eq:etapsi} in the periodic setting has not been established, one cannot exclude the possibility that the energy cascade eventually triggers a finite-time singularity. 
One possible way in which singularity can manifest is the phenomenon of wave breaking, characterized by the blow-up of the solution's gradient despite the boundedness of the interface elevation. 
In simplified models capturing the nonlocal dispersion of gravity waves, such as the Whitham and fractional Korteweg–de Vries equations,  sufficiently steep initial data lead to gradient blow-up and shock formation \cite{H1, Oh, Saut}. Similar phenomena occur in integrable shallow water models, such as the Camassa–Holm \cite{CE98} and Degasperis–Procesi equations \cite{LiuYin}.
A second breakdown is also physically plausible within the full free-boundary Euler equations, with the formation of splash and  splat singularities
 \cite{CCFGG, CoSh}, where the interface self-intersects at a point or along an arc.
Whether the quasilinear structure of the rotational water waves equations allows for the energy to be indefinitely redistributed across high frequencies, or whether the cascade  produces  a singularity in finite-time, 
constitutes a fundamental open problem at the intersection of wave turbulence and the theory of singularity formation.

\medskip 

\noindent {\sc 8. \underline{Small scales creation in Euler equation:}}
The study of small-scale creation and growth of high-order Sobolev norms for the incompressible Euler equations (without free boundary) has a long history.
In contrast with our setting, where the primary dynamical observables are the free surface $\eta$ and the irrotational velocity field $\nabla \Phi$, in the Euler equations without free boundary the fundamental quantity is the vorticity and its evolution.
Important progresses have been achieved in the last ten years. 
In the two dimensional case, Kiselev and Šverák \cite{KS14} proved the double–exponential growth of the vorticity gradient for the Euler flow on the disk, while Zlatoš  constructed examples with exponential growth on the torus \cite{Zl15} and  double-exponential growth  on the half-plane \cite{Zl26}. Analogous growth of the vorticity gradient has also been established in the free-boundary  Euler setting in  the pure capillary regime \cite{HLY2024}.
Very recently, Alazard and Said \cite{AS26} showed that for a dense $G_\delta$ subset of initial vorticities in $H^s_0(\mathbb{T}^2)$,  the  solutions undergo infinite time growth of Sobolev norms. 
For three-dimensional incompressible Euler flows, vorticity stretching obstructs global regularity, and whether it leads to finite-time singularities remains a central open problem. Finite-time blow-up has been proved by Elgindi in the  low-regularity  $C^{1,\alpha}$ space \cite{E17}.


\medskip 

\noindent {\sc 9. \underline{Scaling and Sign Symmetry:}}
Define, for $ \lambda >0$, 
$$
(\eta_\lambda(t,x), \psi_\lambda(t,x)):= ( \eta(\lambda t,x), \lambda \psi(\lambda t , x)) , \quad (\eta^\vee(t,x), \psi^\vee(t,x)):= ( \eta( t,-x),  \psi( t , -x))\,.
$$
If $(\eta,\psi)$ is a solution of \eqref{eq:etapsi} then $(\eta_\lambda, \psi_\lambda)$ is a solution of \eqref{eq:etapsi} with $(g,\gamma)\leadsto (\lambda^2 g , \lambda \gamma)$ and $(\eta^\vee, \psi^\vee)$ is a solution of \eqref{eq:etapsi} with $ (g,\gamma) \leadsto (g,- \gamma)$. 
For this reason, along the manuscript  \textbf{we assume $\mathbf{g=1}$ and $\mathbf{ \boldsymbol \gamma<0}$} without losing any generality.

\medskip 

Let us show   how \Cref{thm:mainimprove} directly implies  Theorem \ref{thm:main.intro}.
\begin{proof}[{\bf Proof of Theorem \ref{thm:main.intro}}]
Recall that solving equations \eqref{eq:etapsi} enables to reconstruct the velocity $\vec{v}$ and the pressure $P$ solving the original Euler system \eqref{euler}--\eqref{bottom} on the whole domain $\cD_\eta$  (see \cite{CIP, Wh}), with $\vec{v} = \nabla \Phi - \vect{\gamma y}{0}$. Then
by \eqref{def:V}, \eqref{form-of-B} one has $\vec{v}(t,\cdot)\vert_{\Gamma(t)} = (\sfV(t, \cdot)-\gamma \eta(t,\cdot),\ \sfB(t,\cdot))$.
So  the first of \eqref{initial.intro} yields 
$
\| \vec v(0,\cdot)\vert_{\Gamma(0)}\|_{H^\ts(\T)} \lesssim \delta$, 
whereas the third of \eqref{explosion.improved} gives
$
 \sup_{t \in [0,T]}  \norm{\vec{e}_y \cdot \vec v(t, \cdot)\vert_{\Gamma(t)} }_{H^\ts(\T)} \geq K$.
\end{proof}

\subsection{Scheme of the  proof and structure of the paper}\label{sec:scheme_proof}

The proof is organized around the extraction, from  the full water waves dynamics,  of  an effective system where an instability mechanism can be isolated and quantitatively analyzed. 
The whole proof combines a resonance analysis, a normal form construction adapted to a distinguished pair of modes, a bootstrap analysis for special solutions --that we call  long-time controlled-- and a Mourre's inspired positive commutator estimate. 
We now describe the main steps of the strategy.
\medskip

\noindent{\bf 1. Resonance analysis and selection of a resonant pair.}
The starting point of the analysis is the study of the linear frequencies $\Omega_j(\gamma)$ in \eqref{freq.intro}, carried out in Section \ref{sec:analisi.res}.
In Lemma \ref{lem:2wave} we   prove that, whenever $\gamma^2 \in \Q\setminus\{0\}$, there exist integers $m,n$ of opposite sign and different modulus for which non-integrable $2$-wave resonances of the  form $\Omega_m(\gamma)=\Omega_n(\gamma)$ occur. 
These resonances are the source of the transport term appearing in \eqref{Z.eq.intro}.
In contrast, Lemma \ref{3onde} establishes the absence of $3$-wave resonances, in fact for every value of $\gamma$. 
At order four, however, non-integrable resonances may appear. These constitute the rotational analogue of the Zakharov–Dyachenko resonances \cite{ZakD} arising in the irrotational case $\gamma=0$. Unlike the irrotational setting, the coefficients of the corresponding resonant monomials do not vanish, and this constitutes an obstruction to propagate smallness of low-Sobolev norms, as required in \eqref{low.norm.intro0}. 
It also plays a key role in governing the long-time nonlinear dynamics at low frequencies, which will be crucial for the instability mechanism.

To isolate the resonant effects, for a fixed vorticity $\gamma<0$ with $\gamma^2\in\Q$, we identify special two-mode sets 
$$
\Lambda=\{\tm,\tn\} \quad \mbox{ with } \tm<0<\tn \ \  \mbox{ and }\ \  \tn+\tm>0 \ , 
$$
which we call $\gamma$-good (see Definition \ref{g-good}), and  satisfy three key properties:
\begin{itemize}
\item[$(i)$] $\Lambda$ supports a $2$-wave resonance, namely $\Omega_\tm(\gamma)=\Omega_\tn(\gamma)$;
\item[$(ii)$]  any $4$-wave resonance involving at most two modes outside $\Lambda$ is integrable, in the sense that the corresponding Fourier indices are pairwise equal;
\item[$(iii)$]  the coefficients  $\tV_\tm$ and $\tV_\tn$ of the transport operator in \eqref{fH.intro} (explicitly defined in \eqref{VmVn.def}), which depend algebraically on $\tm$ and $\tn$, have opposite signs.
\end{itemize}
Property $(i)$  generates the resonant transport coefficient responsible for the instability.
Property $(ii)$, together with the Hamiltonian structure, shall  ensure that quartic interactions with two modes in $\Lambda$ and two modes outside $\Lambda$ 
preserve the $L^2$ norm and prevent energy leakage from $\Lambda$.
Property $(iii)$ guarantees that the effective transport coefficient $\tJ+\fv(x)$
in \eqref{fH.intro}
can change sign, which is essential for the Mourre-type positive commutator argument used later in the proof.

The construction of $\gamma$-good sets for arbitrary $\gamma^2\in\Q\setminus\{0\}$ is highly nontrivial. 
In \Cref{lem:goodset} we prove, via a delicate number-theoretic analysis of the dispersion relation, that for any $\gamma^2\in\Q\setminus\{0\}$ there exist infinitely many $\gamma$-good sets, and we provide an explicit construction. 
Subsequently, in Lemma \ref{lem:wres}, we prove that the arithmetic conditions \eqref{Omega*} and \eqref{G2} characterizing $\gamma$-good sets ensure that all $4$-wave resonances involving four, three, or two modes in $\Lambda$ are integrable.
This structural result allows us to effectively decouple the dynamics of the Fourier modes supported in $\Lambda$ from those supported outside $\Lambda$, and guarantees that the internal dynamics on $\Lambda$ is integrable up to very long times.  This last analysis is made rigorous by the strong $\Lambda$-normal form, that we now describe.
\medskip

\noindent{\bf 2. Strong $\Lambda$-normal form.}
In Section \ref{sec:formal} we exploit the Hamiltonian structure of the equation to compute the {\em formal} quartic Birkhoff normal form, with the aim of obtaining an explicit description of the nonlinear dynamics of the distinguished modes in $\Lambda$.
The procedure is as follows. 
After passing to the complex  variables \eqref{zak},  we eliminate all cubic monomials of the Hamiltonian by exploiting the absence of three-wave resonances.
At quartic order, however, the presence of non-integrable $4$-wave resonances prevents a complete reduction to an integrable normal form. Nevertheless, we prove that, fixing any  set $\Lambda$ which is  $\gamma$-good, 
 the quartic Hamiltonian can be reduced to a weaker notion of normal form, which we call the \emph{strong $\Lambda$-normal form} (see \Cref{def:wr}). According to this definition, $4$-wave resonances involving four, three, or two modes in $\Lambda$ are integrable and fulfill $L^2$-energy estimates (cf. \eqref{snf2}). 
As a consequence, the dynamics restricted to the modes in $\Lambda$ is approximately integrable over very long time scales, and there are no exchanges of energy between the modes in $\Lambda$ and those supported on $\Lambda^c$.
In addition, we shall
compute the coefficients of the quartic resonant monomials whose frequencies are entirely supported on $\Lambda$. This is the content
of \Cref{prop:cH^{(4)}_+}.  Some of the computations leading to the coefficients $\fa, \fb, \fc$ in \eqref{abc}  have been verified using Mathematica, and the 
 relevant  notebook \texttt{coeff\_abc.nb} can be downloaded at the page \url{https://git.sissa.it/amaspero/transfer-ww-vorticity}.

We stress that 
this Birkhoff normal form is only formal and does not converge,
unlike the 
rigorous 
paradifferential normal form that we describe now.

\medskip

\noindent{\bf 3. Paradifferential normal form.}
In \Cref{sec:paraWW}, we conjugate the water waves equations to  its rigorous paradifferential normal form. 
In contrast with the paradifferential normal forms developed in \cite{BD} (and later in \cite{BFP, BMM2}), our construction provides a genuine \emph{change of coordinates} in phase space: it defines a nonlinear transformation bounded from $H^s_\R(\T;\C^2)$ to itself provided a  low-norm $H^{s_0}_\R(\T;\C^2)$ is sufficiently small,  rather than a map acting only along curves of solutions. 
This is achieved through the introduction of \emph{admissible transformations} (defined in \Cref{sez:admissible}), which allow us to invert the normal form transformations generated as time–1 flows of paradifferential operators.
 In \Cref{thm:nf} we prove the existence of a normal form variable 
$z$, whose Sobolev norms are equivalent to those of the original variables, such that the system \eqref{eq:etapsi} can be written as
\be \label{Z.eq:intro}
	 \begin{aligned} 
	 	&\pa_t z= 	 	- \im \Omega(D)z +
		 \Opbw{ \im \langle \,\sfV\, \rangle (z;x) \xi }z
          +  
  X^{(\Set)}(z) + \text{ remainders}
	 \end{aligned}
	 \ee
where the resonant transport coefficient is given by 
	\be \label{VresZ.intro}
        \begin{gathered}
        \langle \, \sfV\, \rangle(z;x) := \sum_{j \in \Z_*} \sfV^{({\rm int})}_j |z_{j}|^2 +
		\sum_{\substack{m < 0 <n   \\ \Omega_m(\gamma) = \Omega_n(\gamma)}} 2\,  \sfV^{({\rm res})}_{m, n}\,  \Re \Big(  \bar{z_{m}(t)} \, {z_{n}}(t) \, e^{\im (n-m) x} \Big) \, ,
        \end{gathered}
		\ee 
 $X^{(\Set)}(z)$ is a smoothing, cubic vector field  in strong $\Lambda$-normal form (cf. \Cref{def:wr})
and  the remainders in \eqref{Z.eq:intro}
consist of higher-order perturbative terms which remain negligible on long time scales.
In particular, the paradifferential normal form reveals an effective transport operator governing the evolution of the high-frequency component of the solution, whose coefficients in 
\eqref{VresZ.intro}
are determined by the nonlinear dynamics of the modes supported on $2$-wave resonances.

Let us stress that, from the paradifferential normal form procedure alone, we can only deduce 
that $X^{(\Lambda)}$ is in {\em weak-$\Lambda$-normal form}. This means that all its monomials are resonant (see \eqref{wnf}), but they do not necessarily satisfy the 
$L^2$-energy estimate \eqref{snf2}. The underlying difficulty is that the paradifferential normal form procedure does not preserve the Hamiltonian structure of the water waves equations, so we cannot a priori guarantee that 
$X^{(\Lambda)}$ 
 is Hamiltonian.
In \Cref{prop.megliodemax2} we prove that $X^{(\Lambda)}$ actually satisfies the \emph{strong-$\Lambda$-normal form} property. This is achieved through an a posteriori identification argument, which allows us to identify $X^{(\Lambda)}$ with the  resonant cubic Hamiltonian vector field computed formally in \Cref{prop:cH^{(4)}_+}, and in particular to determine its dynamics on $\Lambda$. This argument builds on the one in \cite{BFP}, while accounting for the additional difficulty that in our setting the whole vector field of the paradifferential normal form \emph{does not} coincide with the formal one.

The normal form \eqref{Z.eq:intro}--\eqref{VresZ.intro} is not yet sufficient to produce growth of Sobolev norms. At this stage, it only guarantees that smooth initial data of size $\e$ generate solutions that exist up to times $T_\e \sim \e^{-2}$ and remain small on this time scale (cf. Proposition \ref{prop:e-2}).
To construct solutions exhibiting Sobolev norm growth, we must further analyze the resonant transport operator in \eqref{VresZ.intro}. 
Observe that such transport term depends  on the nonlinear dynamics of all pairs of $2$-wave resonant modes $z_m(t)$ and $z_n(t)$ with $\Omega_m(\gamma)=\Omega_n(\gamma)$. However, we show that it suffices to control the nonlinear  dynamics of the distinguished pair in $\Lambda=\{\tm,\tn\}$. This reduction follows from the strong-$\Lambda$ normal form combined with the notion of long-time controlled solutions, that we now describe.
\medskip

\noindent{\bf 4. Effective dynamics for long-time controlled solutions.}
In \Cref{sec:eff.eq} we exploit the strong-$\Lambda$ normal form to decouple the nonlinear dynamics of the resonant modes supported in $\Lambda$ from the remaining ones, and to analyze the integrable dynamics of the modes in $\Lambda$. 
This analysis applies not to arbitrary solutions, but to a special class of solutions, which we call \emph{long-time controlled} (see Definition \ref{A}).
These solutions are characterized by the following properties, given parameters $(s,\theta,T_*,\e)>0$:
 their initial data have most of the $L^2$ mass, of order $\e$, localized in  $\Lambda$, cf. \eqref{app.ass1}, and they are continuous curves from $[0,T_*]$ to $H^s$ satisfying the large a priori bound
$$
\sup_{0 \leq t \leq T_\star } 
\norm{z(t)}_s \leq  \e^{-\theta} \ . 
$$
The key structural property of long-time controlled solutions, proved in \Cref{prop:eff} and \Cref{cor:dyn.eff}, is that along  each long-time controlled solution, 
the normal form \eqref{Z.eq:intro} simplifies. 
In particular, 
up to perturbative remainders, the Fourier modes supported in $\Lambda$ evolve  according to the integrable vector field  \eqref{eff.sys}, so their nonlinear dynamics is, up to very long times, simply given by
\be\label{zmzn.intro}
\begin{aligned}
& z_{\tm}(t) \approx  e^{- \im t\big(\Omega_{\tm}(\gamma) +  2 \fa |z_\tm(0)|^2 + \fb|z_\tn(0)|^2\big)} z_{\tm}(0)\,,
\qquad 
z_{\tn}(t) \approx  e^{- \im t\big(\Omega_{\tn}(\gamma) +  2 \fc |z_\tn(0)|^2 + \fb|z_\tm(0)|^2\big)} z_{\tn}(0) \,,
\end{aligned}
\ee
where $\fa, \fb, \fc$ are the coefficients defined in \eqref{abc}.
Instead the modes outside $\Lambda$ --in a moving frame with constant speed--  evolve according to a perturbed transport equation:
\be\label{eff.intro}
\pa_t \zeg  = - \im \Omega(D)\zeg   + \im \Opbw{(\tJ + \fv(x)) \xi }\zeg  + \, \mbox{Remainder}\,,
\ee
where  the transport coefficient 
 $\tJ + \fv(x)$ in \eqref{fH.intro},
 obtained substituting \eqref{zmzn.intro} into \eqref{VresZ.intro} and discarding the monomials not supported in $\Lambda$, 
 is time-independent but space dependent. 
 Actually,  $\tJ + \fv(x)$  depends only on $\Lambda$ and on the initial values  $z_\tm(0), z_\tn(0)$.
The proof of this fundamental reduction relies crucially on the bootstrap argument in \Cref{lem:long.boot}, which makes essential use of the strong-$\Lambda$ normal form structure of the vector field, and properties $(i)$ and $(ii)$ of long-time controlled solutions.
This also allows to prove that the low  $H^{s_0}$-norm of long-time controlled solution remains of small size  on the time interval $[0, T_*]$.

\medskip
\noindent{\bf 5. Growth of Sobolev norms for long-time controlled solutions via a positive commutator argument.}
In \Cref{sec:mourre} we establish the second fundamental property of long-time controlled solutions. Any such solution $z(t)$ that  exists up to the enhanced time scale $T_*  \approx \e^{-2} \log(\e^{-1})$ and  has initial data satisfying additional structural conditions  -- that we  call {\em strongly well-prepared}, cf. \Cref{B} -- 
must undergo growth of Sobolev norms:
$$
\norm{z(0)}_{s} \lesssim \e^\theta \lesssim \delta\,, \qquad  \norm{z(T_*)}_s \gtrsim \e^{-\theta} \gtrsim K\ .
$$
This instability result is proved in \Cref{prop:instab}. The key mechanism underlying the growth is the non-constant coefficient transport operator in \eqref{eff.intro}. We choose the initial values $z_\tm(0)$ and $z_\tn(0)$ so that the function $\tJ+\fv(x)$ in \eqref{fH.intro} vanishes at some point. This is ensured by condition $(B2)$ in the definition of strongly well-prepared initial data (see \Cref{B}), which in turn relies on property \eqref{G3} of $\gamma$-good sets, guaranteeing that the coefficients $\tV_\tm$ and $\tV_\tn$ defined in \eqref{VmVn.def} have opposite signs.
 In this situation, the transport operator admits a conjugate operator $\fA$ (defined in \eqref{fA}) such that the commutator
$\im[\fA,\Opbw{(\tJ+\fv(x))\xi}]$ is positive on high-frequencies, as proved in \Cref{lem:pos.comm}.
This Mourre-type estimate allows us to derive the virial inequality
$$
\frac{\di}{\di t} \cA(t)\geq  \, c_1 \,  \e^2 \cA(t)
 -   C \e^{5-3\theta}\, \ , \quad \cA(t):= \la \fA \zeg, \zeg \ra
$$
for some $c_1, C >0$, 
which implies exponential growth of the virial functional $\cA(t)$ over time, provided that the initial quantity $\cA(0)$ is not too small. This last 
non-degeneracy condition depends only on the initial datum and is ensured by assumption $(B3)$ in the definition of strongly well-prepared data.

A fundamental strength of this approach is the robustness of positive commutator estimates: they remain stable under the remainders  of the effective equation \eqref{eff.intro} and yield growth of Sobolev norms of any long-time controlled solution with a strongly well-prepared initial data defined up to  $T_* \sim \e^{-2}\log(\e^{-1})$.

\medskip

\noindent{\bf 6. Growth of Sobolev norms for the normal form variable.}
 At this stage, we have shown that any long-time controlled solution with strongly well-prepared initial data existing up to the enhanced time scale
$T_*  \approx \e^{-2} \log(\e^{-1})$
must undergo Sobolev norm explosion. However, two issues remain open: we have not yet established the existence of strongly well-prepared initial data, nor the existence of long-time controlled solutions on enhanced time-scales.
 Section \ref{sec:growth.z} fills the gap.
First, in   Lemma \ref{lem:z0.wp},  we prove that strongly well-prepared initial data do indeed exist and construct an explicit example supported only of four modes: two modes in $\Lambda$, and two modes at a larger scale, cf. \eqref{z0}.

The central step is Proposition \ref{lem:cresce}, where we prove that {\em any} solution $z(t)$ of the normal form equation 
\eqref{Z.eq:intro} with a strongly well-prepared initial data
must undergo  growth of Sobolev norms. 
The argument  goes roughly as follows: 
consider any of such solutions and define $T$ as the maximal time such that on the interval $[0, T]$ such  solution remains below the threshold
$\norm{z(t)}_s \leq c\e^{-\theta}$ for some $c \in (0,1)$. 
We claim that $T$ must be finite. Indeed, if $T$ were larger than $\approx \e^{-2}\log(\e^{-1})$, then the solution would be long-time controlled up to that enhanced time scale. The instability result of the previous section would then force the $H^s$- norm to exceed the threshold $c\e^{-\theta}$ before time $T$, contradicting the very definition of $T$.
A key technical subtlety is the need to propagate the smallness of the lower $H^{s_0}$-norm  in order to preserve the Taylor sign condition and maintain the validity of the local Cauchy theory throughout the argument.

Finally, once growth for the normal form variable $z(t)$ has been established, we transfer the instability to the original water waves variables. Using the equivalence between the Sobolev norms of $z(t)$ and those of $(\eta(t),\upomega(t))$ --where $\upomega$ denotes Alinhac's good unknown (cf. \eqref{GU})-- and subsequently of the full state vector $(\eta,\psi,\sfV,\sfB)(t)$ (proved in \Cref{UZequiv} and \Cref{UXequiv}), we obtain the existence of a solution $(\eta,\psi,\sfV,\sfB)(t)$ satisfying \eqref{growthXs.intro}.

\medskip

\noindent{\bf 7. Separated growth for $\eta(t)$, $\sfV(t)$, $\sfB(t)$ and proof of the main Theorem.} 
Finally, in Section \ref{sec:growth.physical} we refine this conclusion and prove that all the relevant physical quantities undergo growth of Sobolev norms.
The argument exploits the wave-like structure of the paralinearized system \eqref{ParaWW} 
to show that  the free surface profile $\eta(t)$ and the good unknown $\upomega(t)$ exchange energy over short time intervals.
Technically, it is achieved by an ``upside-down virial argument'' that controls from below the second order derivatives in time of the high Sobolev norms of $\eta(t)$ and $\upomega(t)$. Let us explain it briefly. 
The outcome of  \Cref{sec:growth.z} 
is the existence of a  solution $(\eta(t),\upomega(t))$ and a time $T>0$ such that (cf. \eqref{eta_omega_gr})
$$
\| \eta(t)\|_{s-\frac14}+ \| \upomega(t)\|_{s+\frac14} \approx  \e^{-\theta}, \qquad \forall \, t \in [T-1,\, T].
$$
 By contradiction assume that only one variable has grown, say $\upomega(t)$, and that instead
\be\label{up.ass.intro}
\| \eta(t)\|_{s-\frac14}\ll\e^{-\theta} \ \ \  \forall \, t \in [T-1,\, T] \ . 
 \ee
Computing the  second time derivative of $\| \eta(t)\|_{s-\frac34}$ and using energy estimates, we show that it is  lower bounded by
$$
  \frac{\di^2}{\di t^2} \| \eta(t)\|_{s-\frac34}^2 
     \geq \|\upomega(t)\|_{s + \frac14}^2 \gtrsim \e^{-2\theta}  \ , \quad \forall \, t \in [T-1,\, T]
$$
This uniform positive  ``rapid acceleration'' of the norm $\| \eta(t)\|_{s-\frac34}$ eventually contradicts the boundedness assumption in \eqref{up.ass.intro}, implying that $\norm{\eta(t)}_{s - \frac 1 4}$ must itself grow.
 The same reasoning applies symmetrically to $\upomega(t)$, and is then transferred to the velocity components $\sfV$ and $\sfB$ using norm equivalences and paradifferential estimates.
 
A key feature of the approach is that it exploits the hyperbolic structure of the water waves system in Alinhac’s good unknown formulation. This allows us to rigorously capture the heuristic energy exchange described by  
\eqref{linear.intro}. The good unknown yields cancellations that prevent the loss of derivatives associated with the quasilinear nonlinearity (see \Cref{lem:A2s}), so that the nonlinear terms, while affecting the short-time dynamics, remain perturbative in the virial estimates.

\section{Basic properties of the water waves  and functional setting}\label{sec:fun}
This section contains basic properties of the water waves vector field and
the abstract functional setting used along the paper. 
In \Cref{sec:basic} we describe the Hamiltonian structure of the water waves, the analyticity properties of its vector field, and the linearized water waves  at the rest state.
In  Section \ref{sec:para} we present  definitions and results about paradifferential calculus
 following \cite{BD}
(with 
  the advantage to directly use the 
paralinearization of the 
Dirichlet-Neumann operator 
with multilinear expansions in \cite{BD}).
Note that, following \cite{BMM2}, we define the  notion of $ m$-operators (\Cref{def.m-op}) differently than \cite{BD}. 
In \Cref{sec:composition}, following \cite{BD, BMM2} we describe  composition theorems for paradifferential and $m$-operators.
In Section \ref{sez:admissible}, following \cite{MM}, we introduce a class of nonlinear transformations, \emph{admissible transformations}, which are characterized by being differentiable with respect to the internal variable and nonlinearly invertible. Differently from \cite{MM},  we need to allow such transformations to be unbounded. 
{This is due to the mismatch of  $\tfrac{1}{2}$ derivatives between the Hamiltonian complex Zakharov variable $\zak$, constructed from the canonical variables $(\eta,\wahlen)$, and the complex good unknown $U$, constructed from Alinhac’s good unknown $(\eta,\upomega)$.}
\\

\noindent {\bf Notation:}  Along the paper we deal with real parameters $s$, $s_0$, $\vr$ whose value may vary from line to line while still satisfying the relation:  
\begin{align*}
    s\geq s_0  \gg \vr.
\end{align*}
Moreover we use the following conventions for the sets of natural and integer numbers:
\begin{align}\label{integers}
\N:= \{ 1, 2, \ldots \}, \quad   \N_0:= \N\cup \{0\},  \quad \Z_*:= \Z \setminus \{0\}.
\end{align}

\noindent{\bf Function spaces.} We  denote $\dot L^2(\T;\C):= \dot H^0(\T;\C)$ and 
$ \dot L^2_r:=\dot L^2(\T;\R)= \dot H^0(\T;\R)$ 
the subspace of $\dot L^2(\T;\C)$ made by real-valued functions modulo constants.
Given $u, v \in \dot L^2(\T;\C)$ we define
\be\label{scpr12hom}
\la  u,  v \ra_{\dot L^2_r}:= \int_\T \Pi_0^\bot u(x)\, \Pi_0^\perp v(x) \, \di x  \, , \quad \text{respectively} \quad 
\la u, v \ra_{\dot L^2} := \int_\T \Pi_0^\bot u(x)\, \bar { \Pi_0^\bot v(x)} \, \di x \, ,
\ee
where $\Pi_0 u := \frac{1}{2 \pi} \int_\T u(x) \, \di x$  and $\Pi_0^\perp u := u - \Pi_0 u $ is the  projector onto the 
zero mean functions. 

For $s\in \R$ we shall denote with $H^s(\T;\C^2)$ the space of couples of complex valued Sobolev functions  in $H^s(\T; \C)$  and with 
\be\label{Hs}
\dot H^s_\R(\T;\C^2):=\Big\{ U= \vect{ u^+}{ u^-}\in \dot H^s(\T;\C^2)\colon  \ u^-= \bar{u^+}\Big\}, \qquad \dot L^2_\R(\T;\C^2):= \dot H^0_\R(\T;\C^2) \,.
\ee
Given $r>0$ we set $B_{s}(r)$ the ball of radius $r$ in $\dot{H}^s\left(\mathbb{T};\mathbb{C}^2\right)$ and $B_{s,\R}(r)$ the ball or radius $r$ in $\dot{H}^s_\R(\mathbb{T};\mathbb{C}^2)$.
Sometimes, abusing notation, we shall also denote by $B_s(r)$ the ball of radius $r$ in $H_0^s(\T;\R)$ and in $\dot H^s(\T;\R)$; it will be clear from the context to which ball we refer.
\\
Given an interval $ I\subset \R$ symmetric with respect to $ t = 0 $ and a Banach space $X$,  we use the standard notation $C(I,X)$ to denote the space of continuous functions with values in $X$. 
We denote by  $ B_{X}(I;r)$ the ball of radius $r$ in $ C(I,X )$.
 We simply write 
$B_s(I;r)$ for the ball of radius $r$ in $C(I,\dot{H}^s(\mathbb{T};\mathbb{C}^2))$ and $B_{s,\R}(I;r)$  the ball of radius  $r$ in $C(I,\dot{H}^s_{\R}(\mathbb{T};\mathbb{C}^2))$. \\
	Given $\ell \in \N_0$, we denote by $ W^{\ell, \infty}(\T) $ the space of 
	continuous functions $ u : \T \to \C $, $ 2 \pi$-periodic, 
	whose derivatives up to order  $\ell$ are in $L^\infty$, equipped with the norm 
	\be\label{Wlinf}
	\| u \|_{W^{\ell,\infty}} := 
	\sum_{j=0}^{\ell} \| \pa_x^j u \|_{L^\infty}.
	\ee
    
  \noindent{\bf Symmetries of operators and vector fields.} 
 Given a linear operator  $ A(U)$ acting on $L^2(\T;\C)$, 
we associate to it the linear  operator $\ov{A}(U)$  defined by the relation 
\begin{equation*}
\ov{A}(U)[v] := \ov{A(U)[\ov{v}]} \, ,   \quad \forall \, v: \T \rightarrow \C \, .
\end{equation*}
An operator $A$ is {\em real } if $A = \bar A$. 
We say that a matrix of operators acting on $L^2(\T;\C^2)$  is \emph{real-to-real}, if it has the form 
\begin{equation}\label{vinello}
R(U) =
\left(\begin{matrix} R_{1}(U) & R_{2}(U) \\
\ov{R_{2}}(U) & \ov{R_{1}}(U)
\end{matrix}
\right) \, , \quad \forall \, U \in L^2_\R (\T; \C^2)  \ . 
\end{equation}
 Note that a real-to-real matrix of operators $R(U)$  acts on the real space $L^2_\R (\T; \C^2)$. 
If $R(U)$ and $R'(U)$ are  real-to-real operators then also $R(U)\circ R'(U)$ is real-to-real.

Denoting  the translation operator by 
\be\label{tra}
[\tau_\varsigma u](x) :=u(x + \varsigma), \qquad \varsigma  \in \R,
\ee
we say that a  matrix  $R(U)$ as in \eqref{vinello}  is  translation 
invariant  if 
\begin{equation*}
\tau_\varsigma \circ R(U) = R(\tau_\varsigma U )\circ \tau_\varsigma \ , 
\ \forall \, \varsigma \in \R  \ .
\end{equation*}
Similarly, we will say that a vector field  
\be\label{rtr}
X(U): =  \vect{X(U)^+}{X(U)^-}  \quad \text{is real-to-real if} \quad
\bar{X(U)^+}=X(U)^- \,  \quad \forall \, U \in  L^2_\R (\T; \C^2) \,  
\ee
and translation   invariant if 
 \begin{equation}\label{X.gauge}
 \tau_\varsigma\circ X = X \circ  \tau_\varsigma  \ , 
\quad
\forall \, \varsigma \in \R 
 \,.
\end{equation}
The operator $ R(U)$ in \eqref{vinello} is translation invariant, then the  vector field
$X(U):= R(U)U$ is translation invariant as well.
\smallskip

\noindent{\bf Fourier expansion.}
Given a $2\pi$-periodic function $u(x)$ in $\dot L^2 (\T;\C)$, we expand it in Fourier series as 
\begin{equation}\label{fourierseries}
u(x)= \sum_{j \in \Z_*} u_j \,  e^{\im j x}, \quad u_j:= \frac{1}{2\pi}\int_{\mathbb{T}}u(x) e^{- \im j x }\,\di x \, .
\end{equation}
For functions in $H^s(\T;\C)$, $s\in\R$, we use the same Fourier expansion as in \eqref{fourierseries}, with the sum now taken over all $\Z$ instead of $\Z_*$. We adopt the standard notation $u_0 := \Pi_0 u$. 
We shall expand a couple of functions $ U\in \dot L^2(\T;\C^2)$ as 
\be \label{u+-}
U=\vect{u^+}{u^-}= \sum_{\sigma\in \pm} \sum_{j \in\Z_*} \tq^\sigma u_j^\sigma e^{\im\sigma j x},
 \quad 
   u^\sigma_j:=\frac{1}{2\pi}\int_{\mathbb{T}}u^\sigma(x) e^{- \im \sigma j x }\,\di x \,, \ \ \mbox{ where } \ \tq^+:= \vect{1}{0}\,, \ \tq^-:= \vect{0}{1} \,.
\ee
For $ \vec{\jmath} = (j_1,\dots, j_p) \in \Z_*^p$, $p \geq 1$, 
and $\vec{\sigma} = (\sigma_1,\dots,\sigma_p)\in \{\pm\}^p$ we  denote 
\begin{equation}\label{mod.e.lara}
    |\vec{\jmath}\;| := \max\{|j_1|.\dots, |j_p| \} \,,\qquad \langle \vec{\jmath} \; \rangle := \max\{ 1, |\vec{\jmath}\;|\}
\end{equation}
and
\begin{equation}\label{notationuvecjvecs}
u_{\vec{\jmath}}^{\vec{\sigma}}:= u_{j_1}^{\sigma_1}\cdots u_{j_p}^{\sigma_p} \, , 
\qquad \vec{\sigma} \cdot \vec{\jmath} := \sigma_1 j_1  + \dots+ \sigma_p j_p \,. 
\end{equation}
 We also denote by $\fP_p$ the set of indexes preserving the momentum
 \begin{equation}
 \label{mom1}
 \fP_p:= \left\{  ( \vec{\jmath} , \vec \sigma) \in \Z^p_* \times \{\pm \}^p \colon \quad  \vec \sigma \cdot \vec{\jmath} 
 = 0      \right\} \,. 
 \end{equation} 
  \noindent{\bf  Fourier representation of homogeneous operators and  vector fields.} 
In the sequel we shall encounter  matrices of linear operators,  translational invariant, of the form 
\be\label{Mupm}
M(U)= \begin{pmatrix}M^+_+(U)&M^-_+(U)\\M^+_-(U)&M^-_-(U)\end{pmatrix} , 
\ee
 depending on $U$ in a homogeneous way. We shall call them $p$-homogeneous if they are polynomials in $U$ of order $p$. 
We write them in  Fourier as
\be \label{smoocara0}
\begin{aligned}
M(U)V=\begin{pmatrix}
    (M(U)V)^+
    \\
    (M(U)V)^-
\end{pmatrix}
\,, \quad (M(U)V)^\sigma= 
  \sum_{ \sigma k = \vec \sigma_{p} \cdot \vec \jmath_{p}+ \sigma' j 
   } M_{\vec \jmath_{p}, j,k}^{\vec \sigma_{p}, \sigma',\sigma} u_{\vec \jmath_{p}}^{\vec \sigma_{p}} v^{\sigma'}_{j}  {e^{\im \sigma k x}}\, ,
\end{aligned}
 \ee 
where  the coefficients $M_{\vec \jmath_{p}, j,k}^{\vec \sigma_{p}, \sigma', \sigma} \in \C$  fulfill the  
the following symmetric property:
  for any permutation $ \pi $ of $ \{1, \ldots, p \} $, it results
\begin{equation}
\label{M.coeff.p}
M_{j_{\pi(1)}, \ldots,j_{\pi(p)}, j,k}^{ \sigma_{\pi(1)}, \ldots, \sigma_{\pi(p)},\sigma',\sigma} 
=  
M_{j_{1}, \ldots, j_{p}, j,k}^{ \sigma_{1}, \ldots, \sigma_{p},\sigma',\sigma} \, . 
\end{equation}
The operator $M(U)$ is real-to-real, according to the definition in \eqref{vinello}, if and only if its coefficients fulfill
\be\label{M.realtoreal}
\bar{M_{\vec \jmath_{p}, j,k}^{\vec \sigma_{p}, \sigma', \sigma}} = M_{\vec \jmath_{p}, j,k}^{-\vec \sigma_{p}, -\sigma', -\sigma} \ .
\ee
A $ (p+1)$-homogeneous vector field, which is   translation invariant (see \eqref{X.gauge}), can be expressed in Fourier as: for any $\sigma = \pm$,  
\be\label{polvect}
X(U)^\sigma
= \sum_{k\in \Z} X(U)^\sigma_k 	\, e^{\im\sigma k x} , 
\qquad X(U)^\sigma_k= 
\!\!\!\sum_{
	k \sigma = \vec{\sigma}_{p+1} \cdot \vec{\jmath}_{p+1}
	}\!\!\!\!\!\!\!\!\!
X_{\vec{\jmath}_{p+1}, k}^{ \vec{\sigma}_{p+1}, \sigma} \, u^{\vec{\sigma}_{p+1}}_{\vec{\jmath}_{p+1}} \,,
\ee
the last sum being in $(\vec \jmath_{p+1}, \vec \sigma_{p+1}) \in \Z_*^p \times \{\pm\}^p$, 
and with  coefficients  $ X_{\vec{\jmath}_{p+1}, k}^{ \vec{\sigma}_{p+1}, \sigma}\in \C$  satisfying 
the symmetry condition: for any permutation $ \pi $ of $ \{1, \ldots, p+1 \} $, 
\be\label{symmetric}
X_{\ j_{\pi(1)}, \ldots,j_{\pi(p+1)}, k}^{ \sigma_{\pi(1)}, \ldots, \sigma_{\pi(p+1)},\sigma} 
=  
X_{\ j_{1}, \ldots, j_{p+1}, k}^{ \sigma_{1}, \ldots, \sigma_{p+1},\sigma} \, . 
\ee
The constraint of the indexes in \eqref{polvect} can also be  written as $(\vec \jmath_{p+1}, k, \vec \sigma_{p+1}, - \sigma) \in \fP_{p+2}$ (recall \eqref{mom1}), and we shall often  use this notation.

If $X(U)$ is real-to-real, see \eqref{rtr}, then 
\be\label{X.real}
\bar{X(U)^+_k}=X(U)^-_k\quad \mbox{i.e.}  \quad\bar{X_{\vec{\jmath}_{p+1}, k}^{ \vec{ \sigma}_{p+1}, +}}=X_{\vec{\jmath}_{p+1}, k}^{- \vec{\sigma}_{p+1}, -} \, . 
\ee

\subsection{Basic properties of the water waves vector field}\label{sec:basic}
In this section, we gather some preliminary and basic properties of the water waves system \eqref{eq:etapsi}. 
By \cite{Zak1,CrSu,CIP,Wh} the equations \eqref{eq:etapsi} are the Hamiltonian system 
\begin{equation}\label{HamWW} 
\pa_t\vect{\eta}{\psi} = \cX_{\scH_\gamma}(\eta, \psi) = J_\gamma \begin{pmatrix}\nabla_{\eta} \scH_\gamma(\eta,\psi)\\ \nabla_{\psi} \scH_\gamma(\eta,\psi)\end{pmatrix} \qquad
\text{where} 
\qquad 
J_\gamma := \begin{pmatrix} 0 & \uno \\ -\uno & \gamma \pa_x^{-1} \end{pmatrix} 
\end{equation}
and $\scH_\gamma(\eta,\psi)$ is the real Hamiltonian
\be\label{H.gamma}
{\scH}_\gamma(\eta, \psi)  := \frac12 \int_\T \big(\psi \, G(\eta ) \psi+ g \eta^2\big)  \, \di x +   \frac{\gamma}{2} \int_{\T}  \big(-\psi_x \eta^2+ \frac{\gamma}{3} \eta^3\big)  \, \di x  \, .
\ee
The $L^2$-gradients $(\grad_\eta \scH_\gamma, \grad_\psi \scH_\gamma) $ in
\eqref{HamWW}  belong to (a dense subspace of)  
$ \dot L^2(\T)\times L^2_0(\T)$. 

It was proved in \cite{CM85,CSS97}(see also the monograph \cite[Theorem A.13]{LannesLivre} or the more recent paper \cite[Theorem 1.2]{BMV}), that 
there is $s_0>\frac32$ such that for any  $ \s \geq  s_0 $, 
there is $ r > 0 $ small enough such that  the Dirichlet-Neumann operator mapping    
\be\label{eq:DNan}
\eta \mapsto G(\eta) \, , \quad 
H^\s  (\T; \R) \cap B_{s_0} 
(r) \to {\cal L}(H^\s (\T; \R),H^{\s-1} (\T;\R)) \, , \quad \text{is analytic} \,  .
\ee
In view of  \eqref{eq:DNan} and the algebra properties of Sobolev spaces, for $ \s \geq{s_0}$, 
the water waves vector field   
in \eqref{HamWW} 
\be\label{wwana}
(\eta, \psi) \mapsto \cX_{\scH_\gamma}(\eta, \psi) 
= (X^{(\eta)}, X^{(\psi)}) \, , \quad 
B_{\s}(r)\times  \dot H^{\s}(\T;\R) \to 
H^{\s-1}_0(\T;\R)\times {\dot H}^{\s-1}(\T;\R) \, , 
\quad \text{is analytic} \, , 
\ee
possibly with a smaller $ r > 0 $, with quantitative estimates
\be\label{wanna}
\|X^{(\eta)}(\eta, \psi)\|_{\s - 1} \lesssim_\s \|\psi\|_{\s} + \|\eta\|^2_{\s}\,, \qquad \|X^{(\psi)}(\eta, \psi)\|_{\s -1} \lesssim_{\s} \|\eta\|_{\s} + \|\psi\|_{\s} \quad \forall \, (\eta, \psi) \in B_{\s}(r) \times B_{\s}(r)\,.
\ee

\smallskip 
\noindent{\bf Wahl\'en variables.} 
The variables $(\eta, \psi)$
are not Darboux coordinates, since the   Poisson tensor $ J_\gamma $ in 
\eqref{HamWW}  is not the canonical  one when $\gamma \neq 0$.
Wahl\'en noted in  \cite{Wh} that
under the linear change of variables 
\be \label{Whalen}
\begin{pmatrix} {\eta} \\ { \psi}\end{pmatrix} = \cW \begin{pmatrix} {\eta}\\ {\wahlen} \end{pmatrix} \, , \quad \cW := 
\begin{pmatrix} \uno & 0\\  \frac{\gamma}{2} \pa_x^{-1} & \uno \end{pmatrix} \, , \quad
\cW^{-1} = \begin{pmatrix} \uno & 0\\ - \frac{\gamma}{2} \pa_x^{-1} & \uno \end{pmatrix}  \, , 
\ee
the Poisson tensor  $ J_\gamma $ becomes 
$\cW^{-1} J_\gamma (\cW^{-1})^\top= J
$  the standard one, 
and the Hamiltonian system \eqref{HamWW} assumes the standard Darboux  form
\be \label{HamWW2}
\pa_{t}\vect{\eta}{ \wahlen } = J \vect{\nabla_{\eta} H( \eta, \wahlen)}{\nabla_{\wahlen}  H( \eta, \wahlen)}, \quad 
H( \eta, \wahlen):=\scH_\gamma \big( \eta, \wahlen+ \frac{\gamma}{2} \pa_x^{-1}\eta \big) \, .
\ee
Equivalently,   equipping  the real phase 
 space $\dot  L^2_{r} \times \dot L^2_r $ 
 with the scalar product 
 \be\label{real.bil.form.intro}
\Big\langle  \vect{v_1^+}{  v_1^-}, \vect{v_2^+}{ v_2^-} \Big\rangle_r:= \langle v_1^+, v_2^+ \rangle_{\dot L^2_r} +\langle v_1^-, v_2^- \rangle_{\dot L^2_r}
\ee
 and the symplectic form
\begin{equation}
\label{sympl.form}
\Omega_\R\left (\vect{\eta_1}{ \wahlen_1}, \vect{\eta_2}{ \wahlen_2} \right) := 
\Big\langle E_0
\vect{\eta_1}{\wahlen_1}, \vect{\eta_2}{\wahlen_2} \Big\rangle_r 
= - \langle \wahlen_1, \eta_2 \rangle_{\dot L_r} + \langle \eta_1, \wahlen_2 \rangle_{\dot L_r}  \,, 
\quad E_0 := \begin{pmatrix}
0 & -{\rm Id} \\
{\rm 
Id} & 0
\end{pmatrix} \,, 
\end{equation}
the vector field $X_H=J \vect{\nabla_{\eta} H( \eta, \wahlen)}{\nabla_{\wahlen}  H( \eta, \wahlen)}$ in \eqref{HamWW2} is the unique vector field satisfying 
\begin{equation}\label{defHS}
\Omega_\R\left( X_H(\eta, \wahlen), \vect{\breve \eta}{ \breve \wahlen}\right) = \di H(\eta, \wahlen)\left[\vect{\breve \eta}{ \breve \wahlen}\right] , \quad \forall \vect{\breve \eta}{\breve \wahlen} \in \dot  L^2_{r} \times \dot L^2_r \, . 
\end{equation}
The new Hamiltonian $ H$ is still translation invariant
so  is its Hamiltonian vector field. 
\\[1mm]
{\bf Linearized equation at the equilibrium.} 
The linearized equations  \eqref{HamWW2} in the Wahl\'en coordinates (see \eqref{Whalen}) at the equilibrium $ (\eta,\wahlen)=(0,0)$ are
obtained by conjugating the linearized equations \eqref{eq:etapsi} at $ (\eta,\psi)=(0,0)$, namely 
\be\label{eq:lin}
\pa_t \begin{pmatrix} {\eta}\\ {\wahlen} \end{pmatrix}  \!\! = 	\!
\cW^{-1} \! \begin{pmatrix} 0 & \!\!\!\! |D| \\ -1 & \!\!\!\! \gamma \Hilb  \end{pmatrix} \! \! \cW \!
\begin{pmatrix} {\eta}\\ {\wahlen} \end{pmatrix} 
= \!
\begin{pmatrix}\frac{ \gamma}{2} \Hilb  & \!\!\! |D| \\ -(1+ \frac{\gamma^2}{4} |D|^{-1} ) &\!\!\! \frac{ \gamma}{2}  \Hilb \end{pmatrix} \!\!\begin{pmatrix} {\eta}\\ {\wahlen} \end{pmatrix} 
\ee
where 
 $ D := \frac{1}{\im} \pa_x $
 \nomenclature{$D$}{the H\"ormander derivative $\frac{1}{\im} \pa_x$ } and the   Dirichlet-Neumann operator  
$|D|$ at the flat surface $\eta = 0$ is the Fourier multiplier 
with  symbol
 $ |\xi|$, whereas 
$\Hilb$
is 
 the Hilbert's transform 
\begin{equation}\label{def:Hilbert}
\Hilb u := -\im \ \sign(D)u \,, \qquad \Hilb[1]:= 0 \,.
\end{equation}
 We diagonalize 
system \eqref{eq:lin} introducing the complex variables 
\be\label{zak}
\begin{aligned}
 \zak:= \begin{pmatrix}{ \zetina}\\{ {\bar \zetina} }\end{pmatrix}:= 
   \Lmap^{-1}
   \begin{pmatrix} {\eta}\\ {\psi} \end{pmatrix}
   \quad 
 \, ,  \quad   \Lmap  := \cW \, \cM \,,    
\end{aligned}
\ee
where $\cW$ in \eqref{Whalen} and 
\be\label{defMD}
 \cM:= \frac{1}{\sqrt{2}} \begin{pmatrix} M^{-1}(D) & M^{-1}(D)\\ - \im M(D) & \im M(D) \end{pmatrix}, \qquad 
M(D):= \left( \frac{1+   \frac{\gamma^2}{4} |D|^{-1} }{|D| } \right)^{\frac{1}{4}} \,. 
\ee
The linear map  $\Lmap$ and its inverse $\Lmap^{-1}$ are explicitly given   by
{\small\begin{align}
\label{app_L}
&\Lmap^{-1}  = \frac{1}{\sqrt{2}} \begin{pmatrix}
    M(D)  - \im \frac{\gamma}{2} M^{-1}(D)\pa_x^{-1} & \im M^{-1}(D)
\\
M(D)  + \im \frac{\gamma}{2} M^{-1}(D)\pa_x^{-1} & -\im M^{-1}(D)
\end{pmatrix}
\colon H_0^\s(\T; \R) \times \dot H^\s(\T; \R) \to {H_\R^{\s - \frac14}(\T; \C^2)},
\\
\notag
&\Lmap =
\frac{1}{\sqrt{2}} \begin{pmatrix}
    M^{-1}(D)  & M^{-1}(D) \\
    - \im M(D) + \frac{\gamma}{2} M^{-1}(D)\pa_x^{-1}& 
     \im M(D) + \frac{\gamma}{2} M^{-1}(D)\pa_x^{-1}
\end{pmatrix}
\colon H_\R^\s (\T; \C^2) \to  H_0^{\s-\frac14}(\T; \R) \times \dot H^{\s+\frac14}(\T; \R) \,. 
\end{align}}
\normalsize
Using that  $\zak = \cM^{-1} \vect{\eta}{\wahlen}$, the variables $\zak$ solve the  diagonal linear system 
\be \label{diaglin}
\pa_t  \begin{pmatrix}{ \zetina}\\{ {\bar \zetina} }\end{pmatrix}= 
\cM^{-1} 
\begin{pmatrix}\frac{ \gamma}{2} \Hilb  & \!\!\! |D| \\ -(1+ \frac{\gamma^2}{4} |D|^{-1} ) &\!\!\! \frac{ \gamma}{2}  \Hilb \end{pmatrix} 
\cM  \begin{pmatrix}{ \zetina}\\{ {\bar \zetina} }\end{pmatrix} = 
-\im \vOmega(D)  \begin{pmatrix}{ \zetina}\\{ {\bar \zetina} }\end{pmatrix}, \qquad \vOmega(D):=\begin{pmatrix}\Omega(D) & 0 \\ 0 & -\overline{\Omega(D)} \end{pmatrix}
\ee
where, denoting by $\omega(D)$ the Fourier multiplier with symbol 
\be \label{omegaxi}
 \omega(\xi):=\sqrt{ |\xi| + \frac{\gamma^2}{4}  }  \,, 
 \ee
we put 
\be \label{relation} 
\Omega(D):= \omega(D) + \im \tfrac{ \gamma}{2} \Hilb , \qquad \overline{\Omega(D)}:= \omega(D) - \im \tfrac{ \gamma}{2} \Hilb  \, , 
\ee
and  $\Hilb$ is the Hilbert transform in \eqref{def:Hilbert}.
The identity in \eqref{diaglin}  can be checked by a direct computation 
using the identities
$M(D)^{-1}\circ \big( 1 +   \tfrac{\gamma^2}{4} |D| D^{-2} \big) \circ M^{-1}(D)= \omega(D)= M(D) \circ |D| \circ M(D)$.

System \eqref{diaglin} amounts to the  equation 
$\pa_t \zetina = - \im \Omega(D) \zetina$, 
which,   in Fourier coordinates 
$
\zetina(x) = 
 \sum\limits_{j \in \Z\setminus\{0\}} \zetina_j \,  e^{\im j x}$, where $ 
 \zetina_j := \frac{1}{2\pi} \int_\T \zetina(x) e^{- \im j x } \di x$ 
decouples in infinitely many harmonic oscillators
$$\pa_t \zetina_j = - \im \Omega_j(\gamma) \zetina_j \ , \quad  j \in \Z\setminus\{0\} \ , $$ 
where $  \Omega_j(\gamma) $ are the linear frequencies of oscillations 
\be \label{omegonejin}
\Omega_j(\gamma):= \omega_j(\gamma) +  \tfrac{ \gamma}{2} \sign(j) \, , 
\qquad 
\omega_j(\gamma):=\sqrt{  |j|+  \tfrac{\gamma^2}{4} } \, . 
\ee
\textbf{Hamiltonian in complex  variables $\zak$}.
In the variables $\zak$ in \eqref{zak}, the symplectic form $\Omega_\R$ in \eqref{sympl.form} is pulled-back to 
\be\label{sfc}
\Omega_\C\left(\vect{u_1}{ \bar u_1}, \vect{u_2}{ \bar u_2}\right) 
=  \Big\langle E_c \vect{u_1}{ \bar u_1}, \vect{u_2}{ \bar u_2} \Big\rangle_r\  \,, \quad
 E_c := \im E_0 \quad E_0 \mbox{ in } \eqref{sympl.form} \ , 
\ee
the real-valued Hamiltonian $ H(\eta, \wahlen)$ is   pulled-back  to the real valued Hamiltonian 
 $$\cH(\zetina, \bar \zetina):= H( {\cal M}(\zetina, \bar \zetina)) = \scH_\gamma(\Lmap^{-1} (\zetina, \bar \zetina)) \ , $$ and the Hamiltonian vector field in \eqref{HamWW2} is transformed in the 
Hamiltonian vector field   
\be\label{zetone:formale}
\pa_t \zak = 
\cX_{\cH}(\zak) = J_c \nabla \cH(\zak) = 
\vect{- \im \grad_{\bar \zetina} \cH}{ \im \grad_\zetina \cH} \ , \quad J_c := E_c^{-1} \,,
\ee
that satisfies the  characterization 
\be \label{complexvecham}
\Omega_\C( \cX_{\cH}, \cdot )= \di \cH(\zak)[\cdot ] \, .
\ee
In view of the analyticity in \eqref{wwana} (with $\s\leadsto \s-\tfrac14$) as well as the one of the map $\Lmap^{-1}$ in \eqref{app_L} (with $\s\leadsto\s-\tfrac54$), the push-forward vector field
\be\label{Xzak}
\zak \mapsto \cX_\cH(\zak) \equiv (\Lmap^{-1}_* X)(\zak) := \Lmap^{-1} X(\Lmap \zak), \quad B_{\s,\R}({r'}) \to H^{\s-\frac32}_\R(\T; \C^2)  \ , \quad \mbox{ is { analytic}} \,,
\ee
and Taylor expands as 
\be \label{iniziale_zak}
\pa_t \zak = 
\cX_{\cH}(\zak)
=
-\ii \vOmega(D) \zak + \hamvec{\cH_3}(\zak)+ \hamvec{\cH_4}(\zak)+ \hamvec{\cH_{\geq 5}}(\zak)\,,
\ee
where $\cH_3, \cH_4$ are the cubic and quartic terms of the Hamiltonian $\cH$, and $\cH_{\geq 5}$ collects the higher homogeneity terms.

Finally note that, 
 in the Fourier coordinates 
$ (\zetina_j)_{j \in \Z\setminus\{0\}} $,  
the symplectic form \eqref{sfc} reads, for any
$ U  = (u, \bar u)$, $ V = (v, \bar v)$, 
\be\label{Fex}
\Omega_\C\big( U, V\big) = 2\pi
\sum_{j \in \Z\setminus\{0\} } - \ii \ov{u_j} v_j + \ii u_j \ov{v_j} =-2\pi \im   \sum_{j \in \Z\setminus\{0\},\  \sigma \in \{\pm\} } \sigma u_j^{-\sigma} v_j^\sigma  \, , 
\ee
and consequently the Hamiltonian equation \eqref{zetone:formale}, in Fourier coordinates, are given by 
\be \label{hamvecco0}
\pa_t \zetina_k^\sigma = -\frac{\ii \sigma}{2\pi}  \pa_{\zetina_k^{-\sigma}} \cH(\zak) \, . 
\ee

\subsection{Paradifferential operators and $m$-operators}\label{sec:para}
In this section we present  definitions and results about paradifferential calculus following \cite{BD, BMM2}. We begin with the definition of symbols of paradifferential operators.

\smallskip
\noindent{\bf Symbols.} 
We define  the class of symbols which we will use along the paper. 
They correspond to the autonomous symbols of Definition 3.3 in \cite{BD}, where the dependence on time enters only through the function $U=U(t)$. In view of this, we do not need to keep track on the regularity indexes in time and we fix $K = K' =0$ with respect to Definition 3.3 of \cite{BD}.
\begin{definition}[Symbols]\label{def:sfr}
Let $ m  \in \R $, $N,\, \ell \in \N_0$, $p \in \N$, $s_0, r>0$.
\begin{enumerate}
\item{\bf H\"older symbols.}\label{Hsym} 
\nomenclature{$\Gamma^m_{W^{\ell,\infty}}$}{H\"older symbols, see \Cref{def:sfr}-\Cref{Hsym}}
We denote by $\Gamma^m_{W^{\ell,\infty}}$  the space of 
functions $ a : \T\times \R\to \C $,  $a(x, \xi)$, 
which are $C^\infty$ with respect to $\xi$ and such that, for any  $ \beta \in \N_0 $, 
there exists a constant $C_\beta >0$ such that
\begin{equation}\label{simb-pro}
\big\| \pa_\xi^\beta \, a(\cdot, \xi) \big\|_{W^{\ell,\infty}} \leq C_\beta \, \la \xi \ra^{m - |\beta|}  , \quad \forall \, \xi \in \R \,  . 
\end{equation}
We endow $\Gamma^m_{W^{\ell,\infty}}$ with the  family of norms defined, for any $n \in \N_0$, by
\begin{equation}
\label{seminorm}
\abs{a}_{m, {W^{\ell,\infty}}, n}:=  \max_{\beta\in \{0, \dots , n\} }\, 
\sup_{\xi \in \R} 
\, \big\| \la \xi \ra^{-m+|\beta|} \, \pa_\xi^{\beta} a(\cdot,  \xi) \big\|_{W^{\ell,\infty}} \, . 
\end{equation} 
		\item {\bf $p$-Homogeneous symbols.}\label{Hhom}
        \nomenclature{$\widetilde{\Gamma}^m_p$}{$p$-homogeneous symbols, see \Cref{def:sfr}-\Cref{Hhom}}
        We denote by $\widetilde{\Gamma}^m_p$ the space of  $p$-linear symmetric maps from $\left( \dot H^{\infty}\left(\mathbb{T};\mathbb{\C}^2\right)\right)^p$ to
        $ C^\infty(\mathbb{T}\times \R; \C) $ ,
		$ (x, \xi) \mapsto a_p(U_1, \dots, U_p ;x,\xi)$ defined by
		\be \label{espr.hom.sym}
	a_p(U_1, \dots, U_p; x, \xi):=  \sum_{\substack{\vec{\jmath}\in \Z_*^p \\ \vec{\sigma}\in \{ \pm\}^p}} 
a_{\vec \jmath}^{\vec\sigma}(\xi)
({U}_1)_{j_1}^{\sigma_1} \cdots ({U}_p)_{j_p}^{\sigma_p}\, e^{\ii (\vec{\sigma}\cdot \vec{\jmath}) x}\,,
		\ee
		where $a_{\vec{\jmath}}^{\vec{\sigma}}(\xi)$ are complex valued Fourier multipliers,  satisfying
        \begin{align}
        a_{\vec{\jmath}}^{\vec{\sigma}} := a_{j_1, \dots, j_p}^{\sigma_1, \dots, \sigma_p} = a_{j_{\pi(1)}, \dots, j_{\pi(p)}}^{\sigma_{\pi(1)}, \dots, \sigma_{\pi(p)}} \quad \text{ for any } \pi \text{ permutation of } \{1, \dots, p\}\,,
        \label{sym_sy}
        \end{align}
        and
        for some $	\mu\geq 0$, 
		\be \label{homosymbo}
		| \partial_\xi^\beta a_{\vec{\jmath}}^{\vec{\sigma}}(\xi) | \leq C_\beta \la \vec{\jmath}\; \ra^{\mu}  \langle \xi\rangle^{m-\beta},  \quad \forall \, \vec{\jmath}\in \Z_*^p,\, \vec{\sigma}\in\{\pm\}^p, \, \beta\in \N_0\,.
		\ee
		 We shall denote by 
        $$
        a_p(U;x, \xi):= a_p(U, \cdots, U; x, \xi)
        $$
        the polynomial symbol associated to the multilinear symmetric symbol.\\
		We  denote by $\widetilde{\Gamma}^m_0 $ the space of constant coefficients symbols $ \xi \mapsto a(\xi)$ which satisfy \eqref{homosymbo} with $\mu= 0 $, and by $\Sigma_p^{N} \wt \Gamma^{-m}$ the class of pluri-homogeneous symbols $\sum_{i=p}^{N} a_i(U;x, \xi)$ with $a_i\in \widetilde{\Gamma}^m_i $. For $p> N$ we mean that the sum is empty.
		\item  {\bf Non-homogeneous symbols.}\label{Hnon}
 \nomenclature{$\Gamma_{\geq p}^m[r]$}{non-homogeneous symbols, see \Cref{def:sfr}-\Cref{Hnon}}
We denote by $\Gamma_{\geq p}^m[r]$ the space of complex-valued functions  $ (U;x,\xi)\mapsto a(U;x,\xi) $,
		defined for $ U \in B_{s_0}(r)$ for some $s_0$ large enough 
        such that  for any $s\geq s_0$,  there is  $r':=r'(s)\in(0,r)$ such that  for  any $\beta \in \N_0$ the following holds. There is $C=C_{\beta, s}>0$ such that, for any $ U \in B_{s_0}\left(r' \right)\cap H^{s}\left(\mathbb{T};\mathbb{C}^2\right)$, and  $\ell \leq s - s_0$, one has the estimate
\begin{equation}\label{nonhomosymbo}
			\left\|\partial_\xi^\beta a\left({ U;\cdot,\xi }\right)\right\|_{W^{\ell,\infty}}  \leq C \langle \xi \rangle^{m-\beta} \| U\|_{{s_0}}^{p-1}\|U\|_{s} \, .
		\end{equation}
		  In addition, we require also  the translation invariance property
\begin{equation} \label{mome}
a\left( \tau_{\varsigma} U; x,\xi\right)= a\left( U; x+\varsigma, \xi\right),\quad \forall 
\varsigma\in \R \, , 
\end{equation}
where $\tau_\varsigma$ is the translation operator in \eqref{tra}.
\item {\bf Symbols.}\label{Hsymbols}
 \nomenclature{$\Sigma\Gamma^m_0[r,N]$}{the class of symbols, see \Cref{def:sfr}-\Cref{Hsymbols}}
We denote by $\Sigma\Gamma^m_0[r,N]$ the class of symbols of the form
\begin{equation}
\label{symbols}
a(U; x, \xi) = \sum_{i=0}^{N-1} a_i(U;x, \xi) + a_{\geq N}(U; x, \xi)\,, 
\end{equation}
where $ a_0(U; x, \xi) \equiv a_0 (\xi)\in \wt\Gamma^m_0$ is a Fourier multiplier, $a_i \in \wt \Gamma_i^m$ for $i=1,\dots, N-1$, and  $ a_{\geq N} \in \Gamma^m_{\geq N}[r]$.
For $q=1,\dots, N-1$, we denote by $\Sigma\Gamma^m_q[r,N]$ the class of symbols of the form \eqref{symbols} with $a_j = 0$ for $j\leq  q-1$.
We say that a symbol  $a(U;x,\xi) $ is \emph{real} if it is real-valued for any
		$ U \in B_{s_0,\R}(I;r)$. 
\item {\bf Functions.}\label{functions}
\nomenclature{$\wt \cF_p, \, \cF_{\geq p}[r], \, \Sigma \cF_p[r,N]$}{the class of $p$-homogeneous, non-homogeneous and general functions, see \Cref{def:sfr}-\Cref{functions}}
We denote by $\wt \cF_p$ (respectively $\cF_{\geq p}[r]$)  the subspace of $\wt \Gamma^0_p$ (respectively $\Gamma^0_{\geq p}[r]$)   made of those
symbols which are independent of $\xi$, by $\wt \cF^\R_p$ (respectively $\cF^\R_{\geq p}[r]$) the functions in $\wt \cF_p$ (respectively $\cF_{\geq p}[r]$)  which are real-valued, and by $\Sigma\cF_p[r,N]$ (respectively $\Sigma\cF^\R_{p}[r,N]$) the subspace of  $\Sigma\Gamma^0_p[r,N]$ (respectively $\Sigma\Gamma^0_{p}[r,N]$) made of those symbols which are independent of $\xi$. 
\end{enumerate}
\end{definition}
\begin{remark}\label{rmk.simbolini}
We now list some well known properties regarding symbols. 
\begin{itemize}
 \item  If $a$ is a symbol in $  \Gamma^m_{{W^{\ell,\infty}}} $ 
then $ \partial_x a  \in  \Gamma^{m}_{{W^{\ell-1,\infty}}}$ and  
$ \partial_\xi a \in   \Gamma^{m-1}_{{W^{\ell,\infty}}}$.
If  $b $ is a symbol in $  \Gamma^{m'}_{{W^{\ell,\infty}}}$  then 
$a b \in  \Gamma^{m+m'}_{{W^{\ell,\infty}}}$.
If $a \in \Sigma \Gamma^{m}_{ p}[r,N]$ and $b \in \Sigma \Gamma^{m'}_{q}[r,N]$, then  $ab  \in  \Gamma^{m+m'}_{\geq p+q}[r]$  and $\pa_x a \in\Sigma \Gamma^{m}_{ p}[r,N] $, $\pa_\xi a \in \Sigma \Gamma^{m-1}_{p}[r,N]$.

\item  $p$-homogeneous symbols in $ \widetilde \Gamma^m_p$ and non-homogeneous symbols in $\Gamma^m_{\geq p}[r]$  are actually functions with values in $\Gamma^m_{W^{\ell,\infty}}$ for some $\ell \in \N$, whose semi-norms \eqref{seminorm} are bounded, for some $\mu, s_0>0$, by
\begin{equation}\label{nonhomosymbo.homo}
|a_p|_{m, W^{\ell,\infty}, n} \leq C_n \, \| U\|_{{1}}^{p-1}\|U\|_{{\ell+\mu+1}}  \ , 
\quad
|a|_{m, W^{\ell,\infty}, n} \leq C_n \, \norm{U}_{{s_0}}^{p-1} \norm{U}_{s} \ , \ \ \ell\leq s - s_0  \ 
\end{equation}
and \eqref{espr.hom.sym} implies the translation invariance property \eqref{mome}.

\item If  $a$ is a $p$-homogeneous real-valued symbol in $\widetilde \Gamma^m_p$ then its Fourier coefficients (see \eqref{espr.hom.sym} satisfy
\begin{align}
   \ov{a_{\vec{\jmath}}^{\vec{\sigma}}}=a_{\vec{\jmath}}^{-\vec{\sigma}}, \qquad \text{for any } \vec{\jmath}\in\Z^p_*, \ \vec{\sigma}\in \{\pm\}^p.
   \label{reality_cond}
\end{align}
Moreover, $a(U;x,\xi)$ satisfies the translation invariant property \eqref{mome}. 

\item {If $a$ is a $p$-homogeneous symbol in $\wt \Gamma^m_p$, then it is also a non-homogeneous symbol in $\Gamma^m_{\geq p}[r]$ for any $r$.}
\end{itemize}
\end{remark}

\begin{remark}\label{rem:not.symm}
Sometimes we shall write a symbol $a_p(U; x, \xi)$ only in polynomial form
\be \label{not-symm}
		a_p(U;x,\xi):=\sum_{\substack{\vec{\jmath}\in \Z_*^p\\ \vec{\sigma}\in \{ \pm\}^p}} \wt a_{\vec{\jmath}}^{\vec{\sigma}}(\xi) \, u_{\vec{\jmath}}^{\vec{\sigma}} \, e^{\ii (\vec{\sigma}\cdot \vec{\jmath}) x}
		\ee
        with some Fourier multiplier coefficients $\wt a_{\vec{\jmath}}^{\vec{\sigma}}(\xi) $ not necessarily symmetric, but fulfilling the estimates
\eqref{homosymbo}. 
One obtains the symmetric coefficients $a_{j_1, \dots, j_p}^{\sigma_1, \dots, \sigma_p} $ in  expression \eqref{espr.hom.sym} by symmetrizing, i.e., denoting by $\cS_p$ the symmetric group of permutations of $\{1, \ldots, p\}$, writing
$$
a_{j_1, \dots, j_p}^{\sigma_1, \dots, \sigma_p} = 
\frac{1}{p!}\sum_{\pi \in \cS_p} \wt a_{j_{\pi(1)}, \dots, j_{\pi(p)}}^{\sigma_{\pi(1)}, \dots, \sigma_{\pi(p)}}  \ ;
$$
the new coefficients again fulfill  estimate \eqref{homosymbo}.
We shall use the  notation \eqref{not-symm} for example  in formulas \eqref{V2.exp0}-\eqref{V2.exp1} and  for the resonant  transport term in \eqref{VresZ}; the reason is that the transport term \eqref{VresZ} is perhaps the most important object of the paper, being the term responsible for the growth, and we prefer to express it in the simplest possible form.

\end{remark}

\smallskip

\noindent{\bf Paradifferential quantization.}
Given $p\in \N_0$ we consider   functions
  $\chi_{p}\in C^{\infty}(\R^{p}\times \R;\R)$ and $\chi\in C^{\infty}(\R\times\R;\R)$, 
  even with respect to each of their arguments, satisfying, for $0<\delta_0\leq \tfrac{1}{10}$\footnote{Actually the parameter $\delta_0>$ will be chosen in \Cref{sec:eff.eq}, see \eqref{s0r}.},
  \be\label{supp.chi}
\begin{aligned}
&{\rm{supp}}\, \chi_{p} \subset\{(\xi',\xi)\in\R^{p}\times\R; |\xi'|\leq\delta_0 \langle\xi\rangle\} \, ,\qquad \chi_p (\xi',\xi)\equiv 1\,\,\, \rm{ for } \,\,\, |\xi'|\leq \delta_0 \langle\xi\rangle / 2 \, ,
\\
&\rm{supp}\, \chi \subset\{(\xi',\xi)\in\R\times\R; |\xi'|\leq\delta_0 \langle\xi\rangle\} \, ,\qquad \quad
 \chi(\xi',\xi) \equiv 1\,\,\, \rm{ for } \,\,\, |\xi'|\leq \delta_0   \langle\xi\rangle / 2 \, . 
\end{aligned}
\ee
For $p=0$ we set $\chi_0\equiv1$. 
Moreover, we assume that 
$$ 
|\partial_{\xi}^{\ell}\partial_{\xi'}^{\beta}\chi_p(\xi',\xi)|\leq C_{\ell,\beta}\langle\xi\rangle^{-\ell-|\beta|} \, , \  \forall \, \ell\in \N_0, \,\beta\in\N_0^{p} \, , \quad
|\partial_{\xi}^{\ell}\partial_{\xi'}^{\beta}\chi(\xi',\xi)|\leq C_{\ell,\beta}\langle\xi\rangle^{-\ell-\beta}, \  \forall \, \ell, \,\beta\in\N_0 \, .
$$ 

\begin{definition}{\bf (Bony-Weyl quantization)}\label{quantizationtotale}
If $a_p(U;x,\xi)$ is a symbol in $\widetilde{\Gamma}^{m}_{p}$, 
respectively if $a \in \Gamma^m_{W^{M,\infty}}$ or $\Gamma^{m}_{\geq p}[r]$,
we set 
\be\label{regula12}
(a_p)_{\chi_{p}}(U;x,\xi) := \!\!\!\!\! \sum_{\substack{\vec{\jmath}\in \Z^p_*\\ \vec{\sigma}\in \{ \pm\}^p}}\chi_p (\vec{\jmath},\xi) a_{\vec{\jmath}}^{\vec{\sigma}}(\xi) u_{\vec{\jmath}}^{\vec{\sigma}} e^{\ii (\vec{\sigma}\cdot \vec{\jmath}) x},
\quad 
 a_{\chi}(x,\xi) :=\sum_{j\in \Z} 
\chi (j,\xi )\hat{a}(j,\xi)e^{\im j x}  \, ,
\ee
where in the last equality $  \hat a(j,\xi) $ stands for $j^{th}$ Fourier coefficient of $a(x, \xi)$ (or $a(U; x, \xi)$) with respect to the $ x $ variable. 
We define the \emph{Bony-Weyl} quantization of $a_p(U; x, \xi)$ or $ a(U; \cdot)  $ as 
\begin{align}\label{BW}
&\opbw(a_p({U};\cdot))v
= \!\!\!\!\!\!\!\!\!\sum_{\substack{(\vec{\jmath},j,k)\in {\Z^{p+2}_*}\\ \vec{\sigma}\in \{ \pm\}^p\\ \vec{\sigma}\cdot \vec{\jmath}+j=k}}\chi_p 
\left(\vec{\jmath},\frac{j+k}{2}\right)
 a_{\vec{\jmath}}^{\vec{\sigma}}\left(\frac{j+k}{2}\right) 
 u_{\vec{\jmath}}^{\vec{\sigma}} v_j
{e^{\im k x}} \, ,\\
&\opbw(a(U;\cdot))v
=
\!\!\!\!\!\! \sum_{(j,k)\in {\Z^2}} \!\!\!  \chi \left(k-j,\frac{j+k}{2} \right)
\hat{a}\left(U; k-j, \frac{k+j}{2}\right)v_j {e^{\im k x}} \, .\label{BWnon}
\end{align}
\end{definition}
\begin{remark}
Here are a few comments about paradifferential operators.
\begin{itemize}
    \item If $ \chi \Big( k-j, \tfrac{ k + j }{2}\Big) \neq 0 $
then $ |k-j| \leq \delta_0 \langle \frac{j + k}{2} \rangle  $
and therefore, for $ \delta_0 \in (0,1)$, 
\begin{equation}\label{rela:para}
\frac{1-\delta_0}{1+\delta_0} |k| \leq |j|
\leq \frac{1+\delta_0}{1-\delta_0}|k| \, , \quad \forall j, k \in \Z\, . 
\end{equation}
Analogously,
\begin{equation*}
    \chi_p\left(\vec{\jmath},\ \left \langle \frac{j + k}{2}\right \rangle \right) \neq 0 \ \Rightarrow |k-j|=|\vec{\sigma} \cdot \vec{\jmath}\;| \leq p|\vec{\jmath}\;| \leq p\,\delta_{0} \left \langle \frac{j+k}{2}\right \rangle\,,
\end{equation*}
and \eqref{rela:para} holds with $\delta_0 \leadsto p\, \delta_0$.
This relation shows that the action of a paradifferential  operator does not spread much the Fourier support of functions.
\item  If $ a$ is a homogeneous  symbol, the two definitions  of quantization in \eqref{BW} and \eqref{BWnon} differ by a  smoothing operator according to 
Definition \ref{omosmoothing} below. 
\item
Definition \ref{quantizationtotale} 
is  independent of the cut-off functions $\chi_{p}$, $\chi$,  
up to smoothing operators that we define below (see Definition \ref{omosmoothing}), see the remark at page 50  of  \cite{BD}. 
\item
Given  a paradifferential  operator
$ A = \Opbw{a(x,\xi)} $ it results
\be\label{A1b}
\bar A = \Opbw{\bar{a(x, - \xi)}} \, , \quad 
A^\top = \Opbw{a(x, - \xi)} \, , \quad
A^*= \Opbw{\bar{a(x,  \xi)}} \, , 
\ee
where $ A^\top $  and $ A^* $ denote respectively the transposed and  adjoint operator with respect to the complex, respectively real-valued,  scalar product 
of $  L^2(\T; \C) $ in \eqref{scpr12hom}. Moreover, one has $ A^* = \bar A^\top $. 
\item
 A paradifferential operator $A= \Opbw{a(x,\xi)}$ is {\bf real} (i.e. $A = \bar A$) if 
\be \label{realetoreale}
 \bar{a(x,\xi)}= a^\vee(x,\xi) \quad \text{ where} \quad a^{\vee}(x,\xi) := a(x,-\xi) \, .
 \ee
\item
A matrix of paradifferential operators $ \Opbw{A( x,\x)}$ is \textbf{real-to-real}, i.e. \eqref{vinello} holds, if and only if 
the  matrix of symbols $A(x,\x)$ has the form 
\begin{equation}\label{prodotto}
 A(x,\x) =\scalebox{0.8}{$
\left(\begin{matrix} {a}(x,\x) & {b}(x,\x)\\
{\ov{b^\vee(x,\x)}} & {\ov{a^\vee(x,\x)}}
\end{matrix} 
\right)=\left(\begin{matrix} {a}(x,\x) & 0\\
0 & {\ov{a^\vee(x,\x)}}
\end{matrix} 
\right)+\left(\begin{matrix} 0& {b}(x,\x)\\
{\ov{b^\vee( x,\x)}} & 0
\end{matrix} 
\right) \,   \, . 
$}
\end{equation}
Moreover, $\Opbw{A(x, \xi})$ is 
\textbf{self-adjoint} if and only if the matrix of symbols $A(x,\x)$ in \eqref{prodotto} satisfies 
\[
\ov{a(x, \xi)} = a(x, \xi)\,, \quad b(x, -\xi)=b(x, \xi)
\]
and \textbf{linearly Hamiltonian} if and only if it satisfies 
\begin{align}
\ov{a(x, \xi)} = - a(x, \xi)\,, \quad b(x, -\xi)= b(x, \xi)\,. 
\label{lin:ham_complex}
\end{align}
\end{itemize}
In view of \eqref{prodotto}, we will  denote matrices of real-to-real paradifferential operators according to following notations
\begin{equation}\label{vecop}
\begin{aligned} 
&\vOpbw{a(x,\xi)}:= \Opbw{\scalebox{0.8}{$
\begin{bmatrix} a(x,\xi)&0\\0& \bar{ a^{\vee}(x,\xi)}\end{bmatrix}
$}} \, , \  \ \ 
\zOpbw{b(x,\xi)}:= \Opbw{\scalebox{0.8}{$
\begin{bmatrix}0&b(x,\xi)\\ \bar{ b^{\vee}(x,\xi)}&0\end{bmatrix}
$}}\,.
\end{aligned}
\end{equation}
\end{remark}
Along the paper 
we shall use the following  result concerning the action of a paradifferential operator 
in Sobolev spaces. 
 We refer to 
\cite[Theorem A.7]{BMM} for the proof of $(i)$ and to  
 \cite[Proposition 3.8]{BD}   for the proof of $(ii)$, $(iii)$. 
\begin{theorem}{\bf (Continuity of Bony-Weyl operators)}
\label{thm:contS}
Let $m\in \R$, $p\in \N$ and $r>0$. Then:

$(i)$ Let  $ a \in \Gamma^m_{L^\infty} $.
Then $\Opbw{a}$ extends to a bounded operator 
 $H^s \to H^{s-m}$ for any $ s \in \R $  satisfying the estimate, for any $u \in H^s$,  
\begin{align} \label{cont00}
& \norm{\Opbw{a}u}_{{s-m}} \lesssim \, \abs{a}_{m, L^\infty, 4} \, \norm{u}_{{s}} \,.
\end{align}

$(ii)$ Let  $a\in \widetilde{\Gamma}_{p}^{m}$. 
There exists $ s_0 > 0 $ such that for any $s \in \R$,  
there is a constant $C>0$, depending only on $s$ and on \eqref{homosymbo} with $\ell=\beta=0$,
such that for any $U_1,\ldots,U_{p} \in \dot H^{s_0}(\T; \C^2)$ and $v \in \dot H^s(\T; \C)$, one has 
\begin{equation}\label{stimapar}
\|\Opbw{a(U_1,\dots, U_p;\cdot)}v\|_{{s-m}}\leq C\prod_{j=1}^{p}\|U_{j}\|_{{s_0}}\|v\|_{{s}} \, ,
\end{equation}
for $p\geq 1$, while for $p=0$ the estimate \eqref{stimapar} holds by replacing the right hand side with $C\|v\|_{{s}}$.

\smallskip
$(iii)$ Let $a\in \Gamma^{m}_{\geq p}[r]$.
There exists $ s_0 > 0 $ such that  for any $s \in \R$
there is a constant $C>0$  such that
for any $U \in B_{s_0}(r)$  one has 
\be\label{stimapar2}
\|\Opbw{a(U;\cdot)}\|_{\mathcal{L}(\dot {H}^{s},\dot {H}^{s-m})}\leq C\|U\|_{{s_0}}^{p} \, .
\ee
\end{theorem}
\noindent{\bf Classes of $m$-operators and smoothing operators.} 
We introduce $m$-operators and smoothing operators.  This is  a small adaptation of  \cite{BD, BMM2} where we consider only autonomous maps, where again the time dependence is only  through $U(t)$. 
In particular we put $K,K' = 0$ with respect to the notation in \cite{BD,BMM2}.
\smallskip

\noindent
 Given integers $(n_1,\ldots,n_{p+1})\in \N^{p+1}$, we denote by $\max_{2}\{n_1 ,\ldots, n_{p+1}\}$ 
the second largest among  $ n_1,\ldots, n_{p+1}$.

\begin{definition}[Classes of $m$-operators]
\label{def.m-op}
 Let $m \in \R$, $p \in \N_0$ and $r >0$.
\begin{enumerate}
\item  {\bf $p$-homogeneous $m$-operators}\label{Mhom}.
\nomenclature{$\wt \cM_p^m$}{$p$-homogeneous $m$-operators, see \Cref{def.m-op}-\Cref{Mhom}}
We denote by
$\wt \cM^{m}_p$ the class of 
 $(p+1)$-linear operators
$$(\dot H^{\infty}(\T;\C^{2}))^{p}\times \dot H^{\infty}(\T;\C)\ni  (U_{1},\ldots,U_{p}, v) \to M_p(U_1,\ldots, U_p)v \in \dot H^{\infty}(\T;\C)\, , $$    symmetric
 in $(U_{1},\ldots,U_{p})$,
with Fourier expansion
\begin{equation}
 \label{Mp}
M_p(U)v:=
M_p(U, \ldots, U)v= \!\!\!\!\!\!\!\!
 \sum_{\substack{  \vec\sigma_{p}\in \{\pm\}^{p} \\  k -j = \vec \sigma_{p} \cdot \vec \jmath_{p}} }
\!\!\!\!\!\!
  M_{\vec \jmath_{p}, j,k}^{\vec \sigma_{p}} \, u_{\vec \jmath_{p}}^{\vec \sigma_{p}}\,  v_{j}  \, {e^{\ii k x} }
 \end{equation}  
that satisfy the following. There exist $\mu\geq0$, $C>0$ such that 
 for any $(\vec \jmath_{p},j,k) \in \Z_*^{p+2} $,   $ \vec \sigma_p \in \{ \pm \}^{p} $,  one has 
\be\label{smoocara}
 |M_{\vec \jmath_p, j, k}^{\vec \sigma_p} |\leq C \, 
{\rm max}_2\{ \la j_1\ra,\dots, \la j_p\ra,\la j\ra \}^{\mu}\, \max\{ \la j_1\ra,\dots, \la j_p\ra,\la j\ra \}^{m}   \, . 
 \ee
 We denote by $\Sigma_p^{N} \wt \cM^{m}$ the class of pluri-homogeneous maps $\sum_{i=p}^{N} M_i(U)$ with $M_i(U)\in \wt \cM^m_i $. For $p> N$ we mean that the sum is empty.

\item  {\bf Non-homogeneous $m$-operators.}
\label{Mnom}
\nomenclature{$ \cM_{\geq p}^m[r]$}{non-homogeneous $m$-operators, see \Cref{def.m-op}-\Cref{Mnom}}
We denote by $\cM^{m}_{\geq p}[r]$ the class of   
 operators  $(U, v)\mapsto M(U) v $ defined on $B_{s_0}(r)\times \dot {H}^{s_0}(\T;\C) $ for some $ s_0 >0  $, 
  which are  linear in the variable $v $ and such that the following holds true. 
  For any $s\geq s_0$ there exist $C>0$ and 
  $r'=r'(s)\in]0,r[$ such that for any 
  $U\in B_{s_0}(r')\cap {H}^{s}(\T;\C^2)$, 
  any $ v \in  {H}^{s}(\T;\C)$, we have that
\begin{equation}
\label{piove}
\begin{aligned} 
& \|{ M(U)v }\|_{{s - m}} 
\leq C \left( \|{v}\|_{s}\|{U}\|_{s_0}^{p} 
 +\|{v}\|_{{s_0}}\|U\|_{{s_0}}^{p-1}\|{U}\|_{{s}} \right) \quad && \text{if} \quad p \geq 1\, ,\\
 & \|{ M(U)v }\|_{{s - m}} 
\leq C \left( \|{v}\|_{s}
 +\|{v}\|_{{s_0}} \|{U}\|_{{s}} \right) \quad &&\text{if} \quad p =0\,.
 \end{aligned}
\end{equation}
In addition, we require the translation invariance property: let $\tau_\varsigma$ be the translation in \eqref{tra}, then
\begin{equation}\label{def:R-trin}
M( \tau_\varsigma  U) [\tau_\varsigma v]  =  
\tau_\varsigma \big( M(U)v \big) \, , \quad \forall \varsigma \in \R \, . 
\end{equation}

\item {\bf $m$-Operators.} 
\label{Mm}
\nomenclature{$ \Sigma\cM_{q}^m[r,N]$}{ $m$-operators, see \Cref{def.m-op}-\Cref{Mm}}
Let $N \in \N$. We denote by $\Sigma\cM^m_0[r, N]$ the space of operators $(U,v) \to M(U)v$ of the form 
\begin{equation}\label{sum.maps}
M(U) = \sum_{q=0}^{N-1} M_q(U) + M_{\geq N}(U)\,,
\end{equation}
where  $M_q(U)$ is in $\wt \cM_q^m$ for $q\in \{0\,, \dots,\, N-1\}$, and  $M_{\geq N}(U)$ in $\cM^m_{\geq N}[r]$.\\
For $q=1, \dots, N$, we  denote  by $\Sigma\cM^m_q[r, N]$ the operators of the form \eqref{sum.maps} with $M_j = 0$ for all $j\leq q-1$. 
\end{enumerate}
\end{definition}
\begin{remark}\label{rmk:molto_bene}
We point out some comments about $m$-operators.
\begin{itemize}
    \item A $p$-homogeneous $m$-operator $M_p$ is a non-homogeneous $m$-operator. Indeed, \eqref{smoocara} implies the  quantitative estimate: 
for  $ s_0 \geq  \mu+1>0$, for any $s \geq s_0 $, 
any $U, U_1, \dots, U_p \in \dot H^s(\T;\C^2)$, and any $v \in \dot H^s(\T;\C)$
\be \label{smoothing}
\begin{gathered}
\| M_p(U)v \|_{{s-m}}\lesssim_s \| U \|_{{s_0}}^{p} 
	\| v \|_{ s}+ \| U \|_{{s_0}}^{p-1}\| U \|_{s} \|v \|_{{s_0}} \,,\\
  \| M_p(U_1, \dots, U_p)v \|_{{s-m}}\lesssim_s \prod_{j=1}^p\| U_j\|_{{s_0}}
	\| v \|_{ s}+ \sum_{j=1}^p \Big(\prod_{k \neq j} \|U_k\|_{s_0}\Big)\| U_j \|_{s} \|v \|_{{s_0}} \, .
\end{gathered}
\ee
See Lemma 2.8  and 2.9 in \cite{BMM2} for a proof and  note that the first of \eqref{smoothing} is actually \eqref{piove}.
\item {\bf (Paradifferential operators as $m$-operators)}
If $a(U;x,\xi)$ is a symbol in $ \Sigma \Gamma^{m}_{p}[r]$  
then the para-differential operator  $ \opbw(a(U;x,\xi))$ is an $m$-operator
$ \Sigma  \mM_{p}^m[r, N]$. This is a consequence of Theorem \ref{thm:contS}--$(ii)$\&$(iii)$. 
\item  We will meet vector fields of the form $X(U) = M(U)U$ where $M(U)$ is a matrix of $p$-homogeneous $m$-operators as in \eqref{Mupm}.
In this case, the relation between the Fourier coefficients  of the vector field in \eqref{polvect} and those of the $m$-operator in  \eqref{Mp} is given by 
\begin{equation}
	\label{simmetrizzata}
	X_{\ j_1, \ldots, j_{p}, j_{p+1}, k}^{\sigma_1, \ldots, \sigma_{p}, \sigma_{p+1},  \sigma} = 
	\frac{1}{p+1} \left(
	M_{\ j_1, \ldots, j_{p}, j_{p+1}, k}^{\sigma_1, \ldots, \sigma_{p}, \sigma_{p+1},  \sigma} + 
	M_{\ j_{p+1}, \ldots, j_{p}, j_1, k}^{\sigma_{p+1}, \ldots, \sigma_{p}, \sigma_1,  \sigma} +
	\cdots
	+
	M_{\ j_1, \ldots, j_{p+1}, j_{p}, k}^{\sigma_1, \ldots, \sigma_{p+1}, \sigma_{p},  \sigma} 
	\right) \, , 
\end{equation} 
namely, they are obtained symmetrizing with respect to the second last index $(j,\sigma')$  the coefficients  $M_{\ \vec \jmath_{p}, j, k}^{\vec \sigma_{p}, \sigma',\sigma}$.
\end{itemize}
\end{remark}
\smallskip

If $m  \leq 0 $ the  $m$-operators  are referred to as smoothing operators.
 \begin{definition}{\bf (Smoothing operators)} \label{omosmoothing}
Let $ \vr\geq0$,   $p\in \N_0$, $N \in \N$ and $q \in \{0, 1 \dots, N-1\}$.
We define the  $\varrho$-smoothing operators 
\be\label{smoothingoper}
\begin{aligned}
& \widetilde{\mathcal{R}}^{-\vr}_{p}:= \widetilde{\mathcal{M}}^{-\vr}_{p} \, ,
\quad
   \mathcal{R}^{-\vr}_{\geq p}[r]:= \mathcal{M}^{-\vr}_{ \geq p}[r]  \,, 
   \quad
 \Sigma\mathcal{R}^{-\vr}_q[r, N]:= \Sigma\mathcal{M}^{-\vr}_q[r, N]  \,.
 \end{aligned}
\ee
\end{definition}
\noindent Given an operator  $M(U)$ in $ \Sigma \mM_{p}^m[r,N]$   of the form \eqref{Mp} 
we denote, for $p\leq q \leq N$, by  
\be \label{pienne}
\cP_{\leq q}[ M(U)] :=  \sum_{j=p}^{q} M_j(U) \, ,
 \quad \text{resp.} \quad 
 \cP_{q}[ M(U)] := M_q(U) \, , 
 \ee 
 the projections on the pluri-homogeneous, resp. homogeneous, operators in 
 $\Sigma_p^q \widetilde \cM^m $ , resp. in $\wt\cM_q^m$. Given an integer $ p\leq p'\leq N$ we also denote 
\be \label{piennebot}
 \cP_{\geq p'}[ M(U)] :=  \sum_{q=p'}^{N-1} M_q(U) +  M_{\geq N}(U)\,.
\ee 
The same notation will be also used to denote 
 pluri-homogeneous/homogeneous components of symbols.

\smallskip 
 \noindent {\bf Spectrally localized maps.} Following \cite{BMM2}, we recall the notion of (non-homogeneous) spectrally localized maps. Such maps 
satisfy the same estimates of  paradifferential operators, see \eqref{bound:specloc}, and include 
paradifferential operators and also linear flows generated by paradifferential operators. The class of spectrally localized maps  is closed under compositions.
 \begin{definition}[{\bf Non-homogeneous spectrally localized maps}]\label{def:specloc}
 \nomenclature{$\cS_{\geq N}^0[r]$}{non-homogeneous spectrally localized maps, see \Cref{def:specloc}}
Let $N\in \N_0$. We denote by $ \cS_{\geq N}^0[r]$ the class of maps $(U,v) \mapsto B_{\geq N}(U)v$ defined on $B_{s_0}(r)\times L^2(\T;\C)$ for some $s_0>0$, which are linear in the variable $v$ and such that the following holds true. For any $ s\in \R $ there are $C>0 $ and $r'=r'(s) \in (0, r)$ such that for any $U\in B_{s_0}(r')$ and any $ v\in \dot H^s(\T;\C)$, we have that 
\begin{align}\label{bound:specloc}
\| B_{\geq N}(U)v\|_{s} \leq C \| U\|_{s_0}^N \| v\|_s.
\end{align}    
\end{definition}

\subsection{ Composition theorems}\label{sec:composition}
Recall  
$D_{x}:=\frac{1}{\ii}\pa_{x}$. 
The following is   Definition 3.11 in \cite{BD}. 
\begin{definition}{\bf (Asymptotic expansion of composition symbol)}
\label{def:as.ex}
Let   $ \varrho  \geq 0 $, $m,m'\in \R$, $r>0$ and $N\in \N$.
Consider symbols $a \in \Sigma\Gamma^{m}_{p}[r,N]$   and $b\in  \Sigma\Gamma^{m'}_{p'}[r,N]$, $p,p' \in \{0, 1, \dots, N-1\}$. For $U$ in $B_{s}(I;r)$ and $\varrho< s- s_0$(where $s_0>0$ is the parameter in \Cref{def:sfr} for $a$ and $b$),
we define the symbol
\begin{equation}\label{espansione2}
(a\#_{\varrho} b)(U;x,\x):=  \sum_{k=0}^{\varrho-1}  \frac{1}{2^k} \sum_{\ell+\beta=k}
\frac{(-1)^{\beta}}{\ell! \beta!} 
(\pa_\xi^\ell D_x^\beta a ) \cdot
(\pa_\xi^\beta D_x^\ell b ) \,.
\end{equation}
\end{definition}
\begin{remark}\label{rem_symbols}
One has the following properties:

\noindent
$ \bullet $
The symbol $ a\#_{\varrho} b $ belongs  to $\Sigma\Gamma^{m+m'}_{p+p'}[r,N]$ 
with estimate
\be \label{sharp.est}
| a \#_{\vr} b |_{m+m', W^{\ell-\vr,\infty}, n} \lesssim | a |_{m, W^{\ell,\infty}, n+\vr-1}| b |_{m', W^{\ell,\infty}, n+\vr-1}\,.
\ee
\noindent
$ \bullet $ Let $\vr \geq 2$, we have  that 
$ a\#_{\varrho}b =ab+\frac{1}{2 \ii }\{a,b\} $ 
up to a symbol in $\Sigma\Gamma^{m+m'-2}_{p+p'}[r,N]$,     
where 
\be\label{poisson}
\{a,b\}  :=  \pa_{\xi}a\pa_{x}b -\pa_{x}a\pa_{\xi}b  \in \Sigma \Gamma^{m+m'-1}_{p+p'}[r,N]
\ee
denotes the Poisson bracket. Moreover, if $a\in \Gamma^m_{W^{M,\infty}}$ and  $b\in \Gamma^{m'}_{W^{M,\infty}}$ then $\{a,b\}\in \Gamma^{m+m'-1}_{W^{M-1,\infty}} $ with estimate
\be \label{poi.est}
| \{a,b\} |_{m+m'-1, W^{M-1,\infty}, n} \lesssim | a |_{m, W^{M,\infty}, n+1}| b |_{m', W^{M,\infty}, n+1}.
\ee

\noindent
$\bullet$ Note that the terms of even (resp. odd) rank in the asymptotic expansion \eqref{espansione2} in the Weyl
quantization are symmetric (resp. antisymmetric) in $(a, b)$. Consequently, the terms of even rank
vanish in the symbol of the commutator and 
\begin{align}\label{commutator}
    a\#_\vr b- b\#_\vr a= \tfrac{1}{\ii} \{a,b\} + \Sigma \Gamma_{p+p'}^{m+m'-3}[r,N].
\end{align}
\end{remark}
 \smallskip

\noindent{\bf Composition of paradifferential operators.} The following proposition was proved in  \cite[Theorem A.8]{BMM} and  \cite[Proposition 3.12]{BD}.

\begin{proposition}{\bf (Composition of Bony-Weyl operators)} \label{teoremadicomposizione}
Let $m, m' \in \R$, $N\in \N$, $0\leq p,p'\leq N-1$, $\varrho \geq 0$ and  $r >0$.

\smallskip
\noindent
$(i)$ 
Let $a \in \Gamma^m_{W^{\varrho, \infty}}$, $b \in \Gamma^{m'}_{W^{\varrho, \infty}}$. Then  
\begin{align}
\label{comp01A}
\Opbw{a}\Opbw{b} 
& =  \Opbw{a\#_\varrho b} + \cQ(a,b)
\end{align}
where the linear operator $\cQ(a,b)$ is bounded ${H}^s \to {H}^{s-(m+m')+\varrho}$,
for any $s \in \R$, and  satisfies, for some $M=M(\varrho) >0$,
\begin{align}
\label{comp020}
\norm{\cQ(a,b)u}_{{s -(m+m') +\varrho}} \lesssim  \left(\abs{a}_{m, W^{\varrho, \infty}, M} \, \abs{b}_{m', L^\infty, M} + \abs{a}_{m, L^\infty, M} \, \abs{b}_{m', W^{\varrho, \infty}, M}  \right) \norm{u}_{s} \,.
\end{align}
When $\vr =2$, one can take $M = 7$.\\
\smallskip
\noindent
$(ii)$ If 
$a\in  \Sigma{\Gamma}^{m}_p[r,N] $ and  $b\in  \Sigma{\Gamma}^{m'}_{p'}[r,N]$. Then the operator $\cQ(a,b)$ in \eqref{comp01A} satisfies 
\begin{equation}
 R(U):= \cQ(a(U;\cdot),b(U;\cdot)) \in \Sigma {\mathcal{R}}^{-\varrho+m+m'}_{p+p'}[r, N]\,.
 \label{compsymbols}
\end{equation}
\end{proposition}
\begin{remark}
Some comments are in order:
\begin{itemize}
    \item  A direct computation shows that  
\begin{align}
    \ov{a\#_\vr b}=\ov{b}\#_\vr \ov{a}, \qquad \big(a\#_\vr b\big)^\vee= b^\vee \#_\vr a^\vee,
\label{prop:ov}
\end{align}
where $ a^\vee $ is defined in \eqref{realetoreale}. Then one has also $ \overline{a^\vee} \#_{\varrho} 
\overline{b^\vee}  = \overline{a \#_{\varrho} b}^\vee  $.
\item {\bfseries Composition of three Bony-Weyl operators.}
  It directly follows by \eqref{comp01A}--\eqref{compsymbols}  that, for symbols $a\in \Sigma \Gamma_{p_1}^{m_1}[r,N]$, $ b\in \Sigma \Gamma_{p_2}^{m_2}[r,N]$ and $c\in \Sigma \Gamma_{p_3}^{m_3}[r,N]$ one has
\begin{equation}
    \Opbw{a}\Opbw{b}\Opbw{c}= \Opbw{a\#_\vr b\#_\vr c}+  \cQ^{(\mathtt{t})}(a,b,c),\label{comp3fin}
    \end{equation}
    where 
    \begin{equation}
    \begin{aligned}
   &a\#_\vr b\#_\vr c:= \tfrac{1}{2}( a\#_\vr(b\#_\vr c)+ (a\#_\vr b)\#_\vr c), \\
   &\cQ^{(\mathtt{t})}(a, b,c):= \tfrac12\big(\Opbw{a}\cQ(b,c) + \cQ(a, b\#_\vr c)+ \cQ(a,b)\Opbw{c} + \cQ(a\#_\vr b,c)\big)
\end{aligned}
\label{sharp3}
\end{equation} 
and recalling $\cQ(\cdot,\cdot)$ of \eqref{comp020}. 
Combining \eqref{compsymbols} and \eqref{stimapar}-\eqref{stimapar2}, one easily gets that $\cQ^{(\mathtt{t})}(a, b,c)$ belongs to $\Sigma \cR^{-\vr+ m_1+m_2+m_3}_{p_1+p_2+p_3}[r,N]$.
Note that then
\begin{align} \label{comp3ordine}
  &   a \#_{\vr} b \#_{\vr} c - \left(a \#_{\vr} b\right) \#_{\vr} c \,, \quad  a \#_{\vr} b \#_{\vr} c -  a \#_{\vr} \left(b \#_{\vr} c \right) \in \Sigma \Gamma_{p_1 + p_2 + p_3}^{-\vr + m_1 + m_2 + m_3}[r, N]\,\\
  \label{magicordine}
    &   a \#_{\vr} b \#_{\vr} c +  c \#_{\vr} b \#_{\vr} a - 2abc \in \Sigma \Gamma_{p_1 + p_2 + p_3}^{ m_1 + m_2 + m_3-2}[r, N]
\end{align}
Moreover, the choice of the symbol in \eqref{sharp3} and the identities in  \eqref{prop:ov} lead to  the following algebraic properties:
\begin{align}
\overline{a\#_\vr b\#_\vr c}= \ov{c}\#_\vr \ov{b}\#_\vr \ov{a}, \qquad
\big(a\#_\vr b\#_\vr c\big)^\vee = c^\vee \#_\vr b^\vee \#_\vr a^\vee .
\label{prop:ov3}
\end{align}
\item \label{comm.op.vz} {\bfseries Commutators of matrices of real-to-real Bony-Weyl operators.} Let $a \in \Sigma\Gamma^m_{p}[r,N]$ and $b \in \Sigma \Gamma^{m'}_{p'}[r,N]$. Then, with the notation in \eqref{vecop}, one has 
  \be\label{commurule}
 \begin{aligned}
 &\left[\zOpbw{b},\vOpbw{a}\right]= \zOpbw{b\#_\vr \overline{a^\vee}- a\#_\vr b}+ R(U)\\
 &\left[\zOpbw{a},\zOpbw{b}\right]= \vOpbw{a\#_\vr \overline{b^\vee}- b\#_\vr \overline{a^\vee}}+ R(U)\\
 &\left[\vOpbw{a},\vOpbw{b}\right]= \vOpbw{a\#_\vr {b}- b\#_\vr {a}}+ R(U)
 \end{aligned}
 \ee
where  $R(U)$ are real-to-real matrices of  smoothing operators in $\Sigma \cR_{p+p'}^{-\varrho + m'+m}[r,N]$.
\item If the operators on the left hand side of \eqref{commurule} are linearly Hamiltonian (according to \eqref{lin:ham_complex}), so are the paradifferential operators on the right hand side, as one verifies using the definition \eqref{lin:ham_complex} of linearly Hamiltonian  and the properties \eqref{prop:ov} of $\#_\vr$.
  \end{itemize}
\end{remark}
\smallskip
\noindent We continue this section with 
the  paralinearization of the product of functions.
\begin{lemma}{\bf (Bony paraproduct decomposition)}
\label{bony}
Let  $f,g,h$ be functions in $H^\s(\T;\C)$  with  
$\s >\frac12$. Then 
\begin{gather}
f\, g = \Opbw{f}g + \Opbw{g}f +\cR(f,g)\label{bonyeq0}
\end{gather}
where $\cR(f,g)=R_1(f)g+R_2(g)f$ and, for $j=1, 2 $, $R_j$ is a homogeneous smoothing operator in $ \widetilde \cR^{-\vr}_{{1}}$ for any $ \vr \geq 0$. Moreover, the bilinear operator $\cR(f,g)$ satisfies the estimates
\begin{align}
    \| \cR(f,g)\|_{\s_1+\s_2-s_0}\lesssim \| g\|_{\s_1}\|f\|_{\s_2}, \quad \text{for any } \s_1+\s_2\geq 0 \ \text{and} \ g \in H^{\s_1}(\T;\C), \ f\in H^{\s_2}(\T;\C).
    \label{stima_BMM}
\end{align}
\end{lemma} 
\begin{proof}
    It is a classical result. One can see e.g. \cite[Lemma 7.2]{BD} and also \cite[Lemma 2.7]{BMM} for the proof of  \eqref{stima_BMM}.
\end{proof}
 
\noindent {\bf Composition of $m$-operators.} The following lemma, which is a consequence of \cite[Proposition 2.14]{BMM2}, shall be used below.
 \begin{lemma}\label{lem:general_composition}
 	Let $m,m', m_0\in \R$, $\vr\geq 0$, $r>0$, $ N\in \N$ and $p \in \{0,1,\ldots ,N-1\}$. 
 	Then:
 	\begin{enumerate}
    \item\label{comp-mappe}  If  $M(U) $ is  in 
$ \Sigma\mathcal{M}^{m}_{p}[r,N]$ and $M'(U)$ is  in
$ \Sigma\mathcal{M}^{m'}_{p'}[r,N] $ then the composition 
$ M(U)\circ M'(U)$  
is  in $\Sigma\mathcal{M}^{m+\max\{m',0\}}_{p+p'}[r,N]$.
\item \label{comp-mappe-int}
If  $M (U) $ is a homogeneous $ m$-operator in $  \widetilde{\mathcal{M}}_{p}^{m}$  
and $M^{(\ell)}(U)$, $\ell=1,\dots,p+1$, are matrices of $q_\ell$-homogeneous  $ m_\ell $-operators  in 
$ \wt \mM^{m_\ell}_{q_\ell} $ with  $m_\ell \in \R$, 
$q_\ell\in \N_0$,   
then  
$$ 
M(M^{(1)}(U)U, \ldots,  M^{(p)}(U)U)M^{(p+1)}(U) \in \wt\mM_{p+\bar q}^{m+ \bar m}
$$ 
with  $ \bar m:=\sum_{\ell=1}^{p+1} \max\{m_\ell,0\}$ and $\bar q:= \sum_{\ell=1}^{p+1}q_\ell$.
\item \label{comp-no-homo} {If  $M (U) $ is a homogeneous $ m$-operator in $  \wt{\mathcal{M}}_{p}^{m}$  
and $M^{(\ell)}(U)$, $\ell=1,\dots,p+1$, are matrices of non-homogeneous $m_\ell $-operators  in 
$\mM^{m_\ell}_{\geq q_\ell}[r] $ with  $m_\ell \in \R$, 
$q_\ell\in \N_0$,   
then  
$$ 
M(M^{(1)}(U)U, \ldots,  M^{(p)}(U)U)M^{(p+1)}(U) \in \mM_{\geq p+\bar q}^{m+ \und{m}}[r]
$$ 
with  $ \und{m}:=\max\{ m_1, \dots, m_{p+1}, 0\}$ and $\bar q:= \sum_{\ell=1}^{p+1}q_\ell$.
}
\item\label{comp-smo-pseudo}  If $ R({U}) \in \Sigma\cR^{-\vr}_{ p}\bra{r,N} $ and $ \mathsf{a}\pare{U;x, \xi}\in \Sigma\Gamma^m_{ p'}\bra{r,N} $, $ 0\leq m \leq \vr$, then
 		\begin{align*}
 			R({U})\circ \Opbw{\mathsf{a}({U;x, \xi}) } \in \Sigma \cR^{-\vr + {m}}_{p+p'}\bra{r,N},
 			&&
 			\Opbw{\mathsf{a}\pare{U;x, \xi} } \circ R\pare{U} \in \Sigma \cR^{-\vr +{m}}_{p+p'}\bra{r,N}\,.
 		\end{align*}
 \item \label{sym-comp-2} 
 Let $ \tp(\xi)$ be a Fourier multiplier in $ \wt \Gamma_0^{m_0}$ and $ \bM_{\geq 1}(U)$ be a real-to-real matrix of operators in $\Sigma\cM_{1}^{m_0}[r,N]$.
 If $a_p(U;x,\xi)$ is a $p$-homogeneous symbol in $ \wt \Gamma_p^m$ then 
 \begin{gather*}
 		\di_U a_p(U;x,\xi)[ -\ii  \mathtt{p}(D)U]= p \, a_p( -\ii  \mathtt{p}(D)U, U, \ldots U;x,\xi)\in \wt \Gamma_p^m, \\
         \di_U a_p(U;x,\xi)[\bM_{\geq 1}(U)U]= p\,  a_p(\bM_{\geq 1}(U)U,U, \dots , U;x,\xi)\in \Sigma \Gamma^{m}_{p+1}[r, N].
 		\end{gather*}
      Moreover 
      \begin{gather}
          \di_U \Opbw{a_p(U;x,\xi)}[  -\ii  \mathtt{p}(D)U]= \Opbw{\di_U a_p(U;x,\xi)[ -\ii  \mathtt{p}(D)U]}\\
          \di_U \Opbw{a_p(U;x,\xi)}[\bM_{\geq 1}(U)U]= \Opbw{\di_U a_p(U;x,\xi)[ \bM_{\geq 1}(U)U]}+ R_{\geq p+1}(U),
      \end{gather}
      where $ R_{\geq p+1}(U)$ is a smoothing remainder in $ \Sigma \cR^{-\vr}_{p+1}[r,N]$ for any $ \vr\geq 0$.
 		\item\label{smooth-comp-2}   If $Q_p(U)$ is a  $p$-homogeneous smoothing operator in $\widetilde{\cR}_p^{-\vr}$ 
   and $ \bM_{\geq 1}(U)$  a real-to-real matrix of operators in $\Sigma\cM_{1}^{m_0}[r,N]$
        then  
 		\begin{gather*}
 		\di_U Q_p(U)[ -\ii  \mathtt{p}(D)U]= p Q_p( -\ii  \mathtt{p}(D)U, U, \ldots U)\in \widetilde{\cR}_p^{-\vr+\max\{0,m_0\}}, \\
         \di_U Q_p(U)[\bM_{\geq 1}(U)U]= p Q_p(\bM_{\geq 1}(U)U,U, \dots , U)\in \Sigma \cR^{-\vr+\max\{0,m_0\}}_{p+1}[r, N].
 		\end{gather*}
 	\end{enumerate}	
 \end{lemma}
 \begin{proof}
     \Cref{comp-mappe} was proved in \cite[Proposition 2.14-Item $(i)$]{BMM2}. \Cref{comp-mappe-int} was proved in \cite[Proposition 2.14, $(ii)$]{BMM2}, see in particular the proof of (2.55) therein. 
     \Cref{comp-no-homo} follows simply combining estimate \eqref{smoothing} and  \eqref{piove}.
     \Cref{comp-smo-pseudo} was proved in \cite[Proposition 3.16]{BD}. \Cref{sym-comp-2} was proved in \cite[Proposition 2.14-Item $(iv)$]{BMM2}. Finally \Cref{smooth-comp-2} was proved in \cite[Proposition 2.14- Item $(ii)$]{BMM2}.
 \end{proof}
\begin{remark}\label{comp.plurihom}
     By item \ref{sym-comp-2} it follows that if 
     $a(U;x, \xi)$ is a pluri-homogeneous symbol in $\Sigma_{1}^{N}\wt\Gamma^m$ then one has
     $\di_U a(U;x,\xi)[ -\ii  \mathtt{p}(D)U]\in\Sigma_{1}^{N}\wt\Gamma^m$ and 
     $\di_U a(U;x,\xi)[\bM_{\geq 1}(U)U]\in \Sigma \Gamma^{m}_{2}[r, N]$.
 \end{remark}
 
We now establish the composition rules for symbols and smoothing remainders under a nonlinear map close to the identity. These rules will be applied at the end of \Cref{sec:paraWW} to replace the internal variables according to the transformation performed (see \Cref{subsec:proofthm}).
 \begin{lemma}\label{lem:sostitution}
 	Let $m, m'\in \R$, $\vr\geq 0$ and  $r>0$. Let   $\cF$ is a non-linear map of the form 
        \begin{align}
            \cF(Z)= Z + \bF_{\geq 1}(Z)Z\,, 
        \end{align}
        where $\bF_{\geq 1}(Z)$ is a real-to-real matrix of operators in $ \cM_{\geq 1}^{m'}[r]$.
 	Then:
 	\begin{enumerate}
 		\item \label{item:composizione_simbolo} If  $a_2(U)$ is  a real-valued $2$-homogeneous symbol in $\wt{\Gamma}_2^{m}$ and $a_{\geq 3}(U)$ is  a real-valued non-homogeneous symbol in $ \Gamma_{\geq 3}^m[r]$,  there is $r'>0$ such that 
        \begin{gather*}
             a_2(\cF(Z);x,\xi)=a_2(Z;x,\xi)+ b_{\geq 3}(Z;x,\xi)\quad \text{where} \quad b_{\geq 3}(Z;x,\xi) \in \Gamma_{\geq 3}^{m}[r] \quad \text{and is  real-valued}\,,\\
              a_{\geq 3}(\cF(Z);x,\xi) \in \Gamma_{\geq 3}^{m}[r']\quad \text{and is real-valued}\,.
        \end{gather*}
        Moreover, one has the following substitution formul\ae
            \begin{align}
\Opbw{a_2(U;x,\xi)}_{|U=\cF(Z)}&=\Opbw{a_2(Z;x,\xi)+ b_{\geq 3}(Z;x,\xi)}+ R_{\geq 3}(Z), \quad R_{\geq 3}(Z) \in \cR^{-\vr}_{\geq 3}[r]\,,\notag\\
       \Opbw{a_{\geq 3}(U;x,\xi)}_{|U=\cF(Z)}&= \Opbw{a_{\geq 3}(\cF(Z);x,\xi)}\,. \label{sost_non-omo}
 		\end{align}
 		\item\label{item:sostituzione_smoothing}  If $Q_2(U)$ is a  $2$-homogeneous smoothing operator in $\widetilde{\cR}_2^{-\vr}$ and $Q_{\geq 3}(U)$ is a non-homogeneous smoothing operator in $\cR^{-\vr}_{\geq 3}[r]$,  there is $r'>0$ such that 
 		\begin{gather}
 		     Q_2(\cF(Z))= Q_2(Z) + R_{\geq 3}(Z), \qquad R_3(Z) \in \cR_{\geq 3}^{-\vr + \max\{0,m'\}}[r]\,,\notag \\
             Q_{\geq 3}(\cF(Z)) \in \cR^{-\vr + \max\{0,m'\}}_{\geq 3}[r']\,.\label{sost_non-omo2}
 		\end{gather}
 	\end{enumerate}
 \end{lemma}
\begin{proof}
    The entire statement follows from \cite[Proposition 2.14, items (ii) and (iv)]{BMM2}; the only claims that remain to be proved are \eqref{sost_non-omo} and \eqref{sost_non-omo2}. To do so, we first use the estimate \eqref{piove} for $ \bF_{\geq 1}(Z)$ obtaining that for any $ Z \in B_{s_0+m', \R }(r')$, 
    \begin{align}
        \| \cF(Z)\|_{s_0}\leq C \| Z\|_{s_0}+ \| Z\|_{s_0}\| Z\|_{s_0+m'} \leq C r' (1+ r'). 
    \end{align}
    Choosing $ r'>0$ so that $ C r' (1+ r') <r $, the compositions in \eqref{sost_non-omo} and \eqref{sost_non-omo2} are well defined and the non-homogeneous estimates \eqref{nonhomosymbo} and \eqref{piove} for $a_{\geq 3}(\cF(Z);x,\xi)$ and $Q_{\geq 3}(\cF(Z))$ follow from the corresponding non-homogeneous estimates \eqref{nonhomosymbo} and \eqref{piove} for $a_{\geq 3}(U;x,\xi)$ and $Q_{\geq 3}(U)$.
 \end{proof}
 \subsection{Admissible transformations}\label{sez:admissible}
  In this section, following the approach in \cite{MM}, we introduce a class of $U$-dependent transformations $U \mapsto \bF(U)$, that we call {\em admissible},  that are differentiable with respect to the internal variable $U$. 
 
  We associate to each of these maps 
   three parameters: the {\em order} $(\nu, m)$ and the {\em gain} $\vr$.
 The parameter $ \nu$ is the loss of derivatives of the map  in the external variable, 
 namely for every sufficiently small and regular $U$, $\bF(U)$  is bounded as a map  $\dot H_\R^{s}(\T; \C^2) \to \dot H_\R^{s-\nu}(\T; \C^2)$. 
 The parameter $m \geq \nu$ measures  the loss of derivatives in the Taylor expansion of the map, namely the homogeneous components $\bF_j(U)$ are bounded as maps $\dot H_\R^{s+m}(\T; \C^2) \to \dot H_\R^{s}(\T; \C^2)$. 
 Such loss may be present even if the map $\bF(U)$ is bounded: this happens for example if the map is the time-1  flow of a paradifferential operator with strictly positive order symbol. 
 Finally, the parameter $\vr$ measures a  gain in regularity in the internal variable $U$ with respect to the external one, namely $\bF(U)\colon \dot H_\R^{s}(\T; \C^2) \to \dot H_\R^{s-\nu}(\T; \C^2)$ requiring only $U \in \dot H^{s-\vr}_\R(\T; \C^2)$. Such gain is necessary in order to prove that admissible maps are nonlinearly invertible, see \Cref{loc.inv}.
 The main property of these maps is that they are differentiable with respect the internal variable.
 Examples of such maps  are flows of paradifferential and smoothing operators, see Lemma \ref{lem:flow.ad} and Lemma \ref{flow.s.ad}.

 \begin{definition}[{\bf Admissible transformations}]\label{admtra}
Let
$$
0 \leq \nu \leq m \leq \vr  \,.
$$
An {\em admissible transformation} of order $(\nu,m)$ with \emph{gain} $\vr$
is a real-to-real matrix $\bF(U)$ of non-homogeneous $\nu$-operators in $\cM_{\ge 0}^{\nu}[r_0]$ for some $r_0>0$, such that there exists $s_0>0$ for which \Cref{def.m-op} is verified and the following holds:

\begin{enumerate}
\item[$(i)$] \emph{(}{\bf Gain \& Linear invertibility}\emph{)}  $\bF(U)$ is linearly invertible  and its inverse  $\bF(U)^{-1}$ is a  real-to-real matrix of  non-homogeneous   ${\nu}$-operators in $\cM_{\geq 0}^{{\nu}}{[r_0]}$, satisfying: 
 for any $s\geq s_0+\vr $ there is a constant $C:=C_{s}>0$  and $r=r(s) >0$ such that for any $ U\in B_{s_0,\R}(r)\cap \dot H^{s{-\vr}}_\R(\T; \C^2)$ and $V \in \dot{H}^{s+\nu}_\R(\T;\C^2)$ one has  
	\begin{align} 
		&\| \bF(U)V\|_{s}+ \| \bF^{-1}(U)V\|_{s}\leq C (\| V \|_{s+\nu} + \|U\|_{s{-\vr}} \|V\|_{s_0}) \,.  \label{lin.est.F}
		\end{align}
	\item[$(ii)$]  \emph{(}{\bf Expansion}\emph{)}  $\bF(U) - \uno$ is a real-to-real matrix of $m$-operators in $\Sigma\cM^m_1{[r_0, 3]}$ expanding as
	\be\label{esp:F}
	\bF(U) = \uno+ \bF_{1}(U) + \bF_{2}(U)+ \bF_{\geq 3}(U),
 \quad \bF_{q}(U)\in \wt \mM^m_q, \quad q=1,2, \quad  \bF_{\geq 3}(U) \in \mM_{\geq 3}^m[r_0].
	\ee
\item[$(iii)$] \emph{(}{\bf Differentiability}\emph{)}  for any $s \geq s_0 + \vr $, 
there is $r'=r'(s)>0$ such that
        the map 
		\be\label{ad:diff}
		B_{s-\vr, \R}({r'})  \ni U \mapsto \bF(U)  \in \cL\big(\dot H^{s+m}_\R(\T;\C^2), \, \dot H^{s}_\R(\T;\C^2)\big) =: {\mathbb{X}^{s,m}} \ \ \mbox{\rm is differentiable}, 
		\ee
         and one has  the quantitative bound: there is $C=C_s>0$ such that for any $U \in B_{s-\vr,\R}(r') $,  $\hat U \in \dot H^{s-\vr}_\R(\T;\C^2)$   and $V \in \dot H^{s+m}_\R(\T;\C^2)$ one has
        \begin{equation}\label{stima.dF}
          \norm{\di_U \bF (U)[\hat U] V}_{s} \leq C  \| \hat U\|_{s-\vr} \norm{V}_{s+m}   \,.
        \end{equation}
		In addition, for any $U \in B_{s-\vr,\R}(r') \cap \dot H^{s+m}_\R(\T;\C^2) $ and $\hat U, V \in \dot H^{s+m}_\R(\T;\C^2)$, we require the rough estimate
		\be
		\begin{aligned}\label{stima.d.adm}
   &    \|   \di_U  \bF_{\geq 3 }(U)[\hat U] V \|_{s}
		 \leq  C \| U\|_{s+m}^{2} \|\hat U\|_{s+m}\,  \| V\|_{s+m}.
		\end{aligned}
		\ee
\end{enumerate} 
\end{definition}
\begin{remark}\label{rem.b.ad}
\begin{enumerate}
    \item[(1)] 
  {If $\nu=0$}, in view of Lemma \ref{lem:general_composition}-Item \ref{comp-mappe}, $\bF(U)$  conjugates  any matrix of $0$-operators in $\mM^0_{\geq p}[r_0]$ to another matrix of $0$-operators in $\mM^0_{\geq p}[r_0]$, namely
  $
\bF(U) \bB_{\geq p}(U) \bF^{-1}(U)
  $
  is a matrix of $0$-operators  in $\mM^0_{\geq p}[r_0]$ for any matrix of $0$- operators $\bB_{\geq p}(U)$ in $\mM^0_{\geq p}[r_0]$.
  \item[(2)]  Property $(ii)$ implies that  for any $s \geq s_0 + \vr$  there exists $r''>0$ such that, for any $U \in B_{s_0, \R}(r'') \cap \dot{H}^{s+m}_\R(\T; \C^2)$ and $V \in \dot{H}^{s+m}_\R(\T; \C^2)$, one has the tame estimates
\begin{align}
		&\| \left[\bF(U)-\uno \right]V\|_{s}+ \| \left[\bF^{-1}(U)-\uno \right]V\|_{s}\leq C (\| U\|_{s_0} \| V \|_{s + m} { + \| U\|_{s + m}\|V\|_{s_0} )}\  \label{lin.est.F1} 
	\end{align}
	 and for $p=1,2$ 
\begin{equation}\label{tri.est.F2}
\|\di_U \bF_p(U)[\hat U]V \|_{{s}}
\lesssim_s \| U\|_{s+m}^{p-1}\| \hat U\|_{s+m}\|V\|_{s+m} 	 \,.
\end{equation}
Since $\bF_{\geq 3}(U)$ is differentiable by difference, the only nontrivial part of \eqref{stima.d.adm} is the quadratic bound in $U$.
\item[(3)] The expansion \eqref{esp:F} for $\bF(U)$ implies the corresponding expansion for $\bF(U)^{-1}$:
\be\label{esp:F_inv}
\bF(U)^{-1}=\uno - 	\bF_{1}(U)-\bF_{2}(U) + \bF_1(U)\bF_1(U)+ \breve \bF_{\geq 3}(U),
\ee
where 
$$\breve \bF_{\geq 3}(U):= -\bF(U)^{-1}\bF_{\geq 3}(U) + (\bF(U)^{-1} - \uno)(\bF_1(U) \bF_1(U)-\bF_2(U)) + {\bF(U)^{-1}[ \bF(U)- \uno -\bF_1(U)] \bF_1(U)}$$
   is a real-to-real matrix of $3m$-operators in $\mM_{\geq 3}^{3m}[r_0]$.
   \end{enumerate}
 \end{remark} 
We now prove that admissible transformations are closed by composition provided the gain of the first map is larger than the order $m$ of the second map.
\begin{lemma}[{\bf Composition of admissible transformations}]\label{lem:comp}
Let  $\bF(U)$ be admissible  of order $(\nu_1,m_1)$ with gain $ \vr_1$ and  $\bG(U)$ be admissible of order $({\nu}_2,m_2)$ with gain $\vr_2$.
If $\vr_2 > m_1$, then  the composition 
$\bF(U)\bG(U)$ is an admissible transformation of order $(\nu_1+\nu_2,m_1+m_2)$ with gain $ \vr:=\min\{ \vr_1, \vr_2-m_1\}$.
\end{lemma}
\begin{proof}
The proof follows the same lines as \cite[Lemma 2.13]{MM}. We include it here to keep track of the slightly different Taylor expansion and of the possibly nonzero parameter $\nu$ in \Cref{admtra}.

Let $r_{0}^{(1)}, r_{0}^{(2)}$ be such that $\bF(U) \in \cM_{\geq0}^{\nu_1}[r_{0}^{(1)}]$ and $\bG(U) \in \cM_{\geq0}^{\nu_2}[r_{0}^{(2)}]$
and  $s_0^{(1)}$, $s_0^{(2)}$   the regularity thresholds required  respectively for  $\bF(U)$ and $\bG(U)$. We immediately note that, by Lemma \ref{lem:general_composition}-Item \ref{comp-mappe}, both the map $\bF(U)\bG(U)$ and its linear inverse $\bG^{-1}(U) \bF^{-1}(U)$ are homogeneous $\nu$-operators in $\cM^{\nu}_{\geq 0}[r_0]$, with $\nu := \nu_1 + \nu_2$, $r_0 := \min\{r_{0}^{(1)},\ r_{0}^{(2)}\}$
and regularity threshold 
$\tilde{s}_0:= \max\{ s_0^{(1)}, s_0^{(2)}\}$.
\\

We now prove that  $\bF(U) \bG(U)$ satisfies Items $(i)$--$(iii)$ of Definition \ref{admtra} with 
\begin{align}
\und{s}_0:= \max\{ s_0^{(1)}, s_0^{(2)}\}+ \vr_2+m_2=\tilde{s}_0+ \vr_2+m_2.
\label{def:s_0new}
\end{align}

\noindent{\scshape Proof of Item $(i)$}. We verify  \eqref{lin.est.F}. 
For any $U\in B_{\tilde s_0,\R}(r)\cap \dot H^{s{-\vr}}_\R(\T; \C^2)$ and $V \in \dot{H}^{s+\nu}_\R(\T;\C^2)$ one has
\begin{equation*}
\begin{aligned}
    \|\bF(U)\bG(U)V\|_{s} 
    &\lesssim_s 
    \|\bG(U)V\|_{s+\nu_1} + \|U\|_{s-\vr_1} \|\bG(U)V\|_{\tilde{s}_0}
    \\
    &\lesssim_s \|V\|_{s+\nu_1+\nu_2} + \|U\|_{s + \nu_1-\vr_2} \|V\|_{\tilde{s}_0} + \|U\|_{s-\vr_1} ( \|V\|_{\tilde{s}_0+\vr_2+\nu_2} + \|U\|_{\tilde{s}_0} \|V\|_{\tilde{s}_0} )
    \\
    &\lesssim_s \|V\|_{s+\nu_1+\nu_2} + \|U\|_{s+ \nu_1- \vr_2} \|V\|_{\tilde{s}_0}   + \|U\|_{s-\vr_1}  \|V\|_{\tilde{s}_0+\vr_2+\nu_2}\,.
 \end{aligned}
\end{equation*}
Therefore, using that $ \vr\leq \min\{ \vr_1,\vr_2-\nu_1\}$  and $\nu_2\leq m_2$, the operator $\bF(U)\bG(U)$ satisfies \eqref{lin.est.F} with $\und{s}_0$ defined in \eqref{def:s_0new}. Reasoning as above one obtains that also the linear inverse $\bG^{-1}(U) \bF^{-1}(U)$ satisfies \eqref{lin.est.F} . \\
\noindent {\sc{Proof of Item $(ii)$}}.
 We set $m:=m_1+m_2$ 
and we verify the decomposition \eqref{esp:F}. One has
\be \label{esp:compo}
\bF(U)\bG(U)=\uno +\bF_1(U) + \bG_1(U)+ \bF_2(U) 
 + \bG_2(U) + \bF_1(U) \bG_1(U) +\bH_{\geq 3}(U)
\ee
where, by item {\ref{comp-mappe}} of Lemma \ref{lem:general_composition}, 
the remainder
\begin{equation}\label{def:H3}
\begin{aligned}
 \bH_{\geq 3}(U)&:=  \bG_{\geq 3}(U) + \bF_1(U)(\bG_2(U) + \bG_{\geq 3}(U))\\
 &\quad + \bF_2(U)(\bG_1(U) + \bG_2(U) + \bG_{\geq 3}(U))
+\bF_{\geq 3}(U)\bG(U)
\end{aligned}
\end{equation}
is a real-to-real matrix of operators in $\mM^{m}_{\geq 3}[{r_0}]$. Note that, as stated, item \ref{comp-mappe} of \Cref{lem:general_composition} guarantees that $\bH_{\geq 3}(U)$ satisfies \eqref{piove} for some \emph{a priori} implicit regularity threshold $s_0>0$. A direct inspection of the proof shows that such an $s_0$ can be chosen not larger than $\underline{s}_0$ in \eqref{def:s_0new}.

\noindent
{\sc{Proof of Item $(iii)$.}}
We first prove  that, for any $ {s} \geq \und{s}_0+\vr$, the map $U \mapsto \bF(U) \bG(U)$ is differentiable at $U \in B_{{s}-\vr, \R}(r) $ for some $r=r({s})>0$ sufficiently small,   and its differential is given by
 \be\label{dFG}
 \di_U \big( \bF(U) \bG(U)  \big)[\hat U]  = (\di_U \bF(U)[\hat U]) \, \bG(U)   
 +
\bF(U) \, (\di_U \bG(U)[\hat U])   \,. 
 \ee
 Indeed, fix $U\in B_{{s}-\vr, \R}(r)$, take $\hat U$ with  $\| \hat U \|_{{s}-\vr} \ll   r$ and put
  \begin{align*}
 \bQ(U,\hat U):=&\bF(U+ \hat U)\bG(U+ \hat U)-\bF(U)\bG(U)- \left((\di_U \bF(U)[\hat U]) \, \bG(U)+\bF(U) \, (\di_U \bG(U)[\hat U]) \right)\\   
 =&\left( \bF(U+ \hat U)- \bF(U)- \di_U \bF(U)[\hat U]\right)\bG(U+ \hat U)\\
 &+ \bF(U) \left( \bG(U+\hat U)-\bG(U)- \di_U \bG(U)[\hat U]\right)+\di_U \bF(U)[\hat U]\left( \bG(U+\hat U )-\bG(U)\right) \\
 =:& \bQ_1(U,\hat U) +\bQ_2(U,\hat U) + \bQ_3(U,\hat U).
\end{align*}
We show that for any $s \geq s_0 + \vr$ 
\be\label{est.Q}
\norm{\bQ_j(U,\hat U)}_{\mathbb{X}^{{s},m}}\lesssim \|{\hat U}\|_{{s}-\vr}^2 \ , \quad j=1,2,3 
\ee
proving formula \eqref{dFG}.
Consider first $\bQ_1(U,\hat U) V$ with $V\in \dot H^{{s} + m}_{\R}(\T; \C^2)$.
Using the differentiability of $ \bF(U)$,  estimate  \eqref{lin.est.F} for $ \bG(U+ \hat U)$ and that $\vr =\min\{ \vr_1,\ \vr_2 - m_1\}$,
we get that 
\begin{align*}
\norm{\bQ_1(U,\hat U) V}_{{s}} & \lesssim \norm{\bF(U+ \hat U) - \bF(U)- \di_U \bF(U)[\hat U]}_{\mathbb{X}^{{s},m_1}} \| \bG(U+ \hat U)V \|_{{s}+m_1}\\
  & \lesssim \| \hat U\|_{{s}-\vr_1}^2 (\|V\|_{{s}+m_1+{\nu_2}}+ \| U+ \hat U\|_{{s}-(\vr_2-m_1)}\| V\|_{s_0}) {\lesssim} \| \hat U\|_{{s}-\vr}^2 \|V\|_{{s}+m_1+{\nu_2}},
\end{align*}
proving \eqref{est.Q} for $j=1$ as $m = m_1+m_2\geq m_1+\nu_2$. \\
We now prove the estimate for $j=2$. Using \eqref{lin.est.F},  the differentiability of $ \bG(U)$, $\vr= \min\{\vr_2-\nu_1, \vr_1\}$, $\nu_1+m_2\leq m$ and $s\geq \und{s}_0+\vr\geq \tilde{s}_0+\vr_2+m_2$, we get 
\begin{align*}
\norm{\bQ_2(U,\hat U) V}_{s}  \lesssim &\norm{\left( \bG(U+\hat U)-\bG(U)- \di_U \bG(U)[\hat U]\right)V}_{{s}+\nu_1}\\
&+ \| U \|_{{s} -\vr_1} \| \left( \bG(U+\hat U)-\bG(U)- \di_U \bG(U)[\hat U]\right)V\|_{\tilde s_0} \\
  \lesssim & \| \hat U\|_{s+\nu_1-\vr_2}^2 \|V\|_{{s}+\nu_1+m_2}+ \| U\|_{{s}-\vr_1}\| \hat U\|_{{\tilde s_0}}^2\| V\|_{{\tilde s_0}+\vr_2+m_2} {\lesssim} \| \hat U\|_{{s}-\vr}^2 \|V\|_{{s}+m}\,,
\end{align*}
proving  also \eqref{est.Q} for $j=2$.\\
Consider now $j=3$. 
Applying first \eqref{stima.dF} for $\di_U \bF(U)[\hat U]$  with $ m \leadsto m_1 $, then 
writing $\bG(U + \hat U) - \bG(U) = \int_0^1 \di_U \bG(U + \tau \hat U)[\hat U] \di \tau $
and using  \eqref{stima.dF} for $\di_U \bG(U+\tau \hat U)[\hat U]$, $\tau \in [0,1]$ with $m\leadsto m_2$ and ${s}\leadsto {s}+m_1$ we get
\begin{align*}
    \norm{\bQ_3(U,\hat U)V}_{s} 
    {\lesssim}&   \|{\hat U}\|_{s-\vr_1}\norm{ \bG(U+\hat U)V-\bG(U)V}_{s+m_1}\\
    \lesssim&   \|{\hat U}\|_{s-\vr_1}\,    \| \hat U\|_{s+m_1-\vr_2}\,  \| V\|_{s+m} 
    {\lesssim}  \| \hat U \|_{s-\vr}^2 \| V\|_{s+m},
\end{align*}
proving  also \eqref{est.Q} for $j=3$. We conclude that  \eqref{dFG} holds.\\
Next we show that 
$\di_U \big( \bF(U) \bG(U)  \big) $ computed in \eqref{dFG}
 fulfills estimate \eqref{stima.dF}. Fix $\hat U \in \dot H^{s-\vr}_\R(\T; \C^2)$, $U \in B_{s-\vr,\R}(r')  $ and 
 $V \in \dot H^{s+m}_\R(\T; \C^2)$ and 
 consider the first term in the right hand side of \eqref{dFG}. We have, for any $s \geq \und{s}_0 + \vr$,
 \begin{align*}
    &  \norm{(\di_U \bF(U)[\hat U]) \, \bG(U)   V  }_{s} \stackrel{\eqref{stima.dF}}{\lesssim}  \|{\hat U}\|_{s-\vr_1} \norm{\bG(U)   V }_{s+m_1} \\
     \stackrel{\eqref{lin.est.F}}{\lesssim}  & 
      \|{\hat U}\|_{s-\vr_1} \left( \norm{V}_{s  +m_1+{\nu_2}} + \norm{U}_{s + m_1-\vr_2}\norm{V}_{{\tilde{s}_0}} \right)
     {\lesssim }
      \|{\hat U}\|_{s-\vr}  \norm{V}_{s +m},
 \end{align*}
 where in the last inequality we used also $ m_1+\nu_2\leq m$, $m_1-\vr_2\leq -\vr$ and $s+m\geq \tilde{s}_0$.
 The second term in \eqref{dFG} has an  analogous estimate, proving that $\di_U \left(\bF(U) \bG(U)\right)$ satisfies \eqref{stima.dF}.\\ 
Finally, we prove the rough estimate \eqref{stima.d.adm} for $\bH_{\geq 3}(U)$ in \eqref{def:H3}. First we compute its differential 
\[
\begin{aligned}
\di_U \bH_{\geq 3}(U)[\hat U]V&= \di_U \bG_{\geq 3}(U)[\hat U]V+ (\di_U \bF_1(U)[\hat U])(\bG_2(U) + \bG_{\geq 3}(U))V+ \bF_1(U)  (\di_U \bG_2(U)[\hat U] \\
&+ \di_U \bG_{\geq 3}(U)[\hat U])V + (\di_U \bF_2(U)[\hat U])(\bG_1(U) + \bG_2(U) + \bG_{\geq 3}(U))V+\bF_2(U)\di_U \bG_2(U)[\hat U]\\
&+ \bF_2(U)  ( \di_U \bG_1(U)[\hat U]  + \di_U \bG_{\geq 3}(U)[\hat U])V+ \di_U \bF_{\geq 3}(U)[\hat U] \bG(U)V + \bF_{\geq 3}(U)  \di_U \bG (U)[\hat U] V\,.
\end{aligned}
\]
Estimate \eqref{stima.d.adm} (with \(m=m_1+m_2\)) for $\di_U \bH_{\geq 3 }(U)[\hat U]V$ follows from the corresponding estimates for $\di_U \bG_{\geq 3}(U)[\hat U]V$, $\di_U\bF_{\geq 3}(U)[\hat U]V$, estimates \eqref{piove}, \eqref{smoothing} for \( \bG_{\geq 3}(U)\), \(\bF_{\geq 3}(U)\), $\bF_1(U)$, $\bF_2(U),$ $\bG_1(U)$, \( \bG_2(U)\) and \eqref{tri.est.F2} for \(\di_U\bF_1(U)[\hat U],\) \(  \di_U\bF_2(U)[\hat U]\), \(\di_U\bG_1(U)[\hat U]\), \(\di_U\bG_2(U)[\hat U]\).
\end{proof}
Next we prove a local invertibility property of the nonlinear map $U \mapsto \bF(U)U$ when $\bF(U)$ is an admissible transformation.  
\begin{lemma}[{\bf Invertibility of admissible transformations}]\label{loc.inv}
	Let  $\bF(U)$ be  a $(\nu,m)$-admissible transformation with gain $ \vr$, with 
    \begin{align}
    \nu + m \leq \vr  \,. 
    \label{restriction:vr_nu_m}
    \end{align}
    Define  the  nonlinear map
	$\mF(U):= \bF(U)U$, then the  following holds true:
	\begin{itemize}
	\item[$(i)$] There exists ${s'_0} \geq 0$ such that for any $s \geq {s'_0}$, the  map $\cF$ is locally invertible: namely there is $r'=r'(s)>0$ and $\mF^{-1}: B_{{s'_0},\R}(r')\cap \dot H^{s+\nu}_\R(\T;\C^2)\to \dot H^{s}_\R(\T;\C^2)$ such that 
	$$
    \mF\circ \mF^{-1}(V)=V \,, \qquad \mF^{-1}\circ \mF(U)=U\,,  \quad \forall \, U,\ V \in B_{s_0',\R}(r')\,.
	$$
    In particular,
    \begin{equation}
\label{stima.inv.adm}
\|\mF^{-1}(V)\|_{s}\leq C_s \|V\|_{s + \nu}, \quad \text{for any} \quad V\in B_{{ s'_0},\R}(r')\cap \dot{H}^{s+ \nu}_\R(\T;\C^2)  \,.
    \end{equation}
\item[$(ii)$]\label{item:inversa_adm} One has  $\cF^{-1}(V) = \bG(V)V$ with $\bG(V) $ a matrix of non-homogeneous ${\nu}$-operators in $\cM_{\geq 0}^{{\nu}}[r'']$ such that $\bG(V) - \uno \in \Sigma \cM_{{1}}^{{3m+\nu}}[r'', 3]$ for some $r'' >0$ and  expands  as
	\begin{equation}\label{exp.G}
	\bG(V) = \uno  - \bF_1(V) -\bF_2(V) + \bF_1(V) \bF_1(V) + {\bF_1(\bF_1(V) V)} +\bG_{\geq 3}(V) \ , 
	\quad \bG_{\geq 3}(V) \in \cM^{3m+\nu}_{\geq 3}[r''] \,,
	\end{equation}
where $\bF_1$ and $\bF_2$ are the homogeneous maps appearing in the expansion \eqref{esp:F} of $\bF$.
	\end{itemize}
\end{lemma}

\begin{proof}
The proof follows the same lines as \cite[Lemma 2.14]{MM}. We include it here to keep track of the possibly nonzero parameter $\nu$ in \Cref{admtra}.

Let $s_0, {r_0} >0$  the parameters given by Definition \ref{admtra} associated to  $\bF(U)$.\\
{\sc Proof of Item $(i)$}  
Let $\sigma_0:= s_0 + \vr$ and define $r:= \min\{r_0,\ r(\s_0),\ r(\s_0 + \nu + m),\ r'(\s_0),\ r'(\s_0 + \nu)\}$, where $r(s)$ and  $r'(s)$ are respectively the radii in \Cref{admtra}-$(i)$ and \Cref{admtra}-$(iii)$. 
We prove that there exists $r_1>0$ such that for any  $V \in B_{\sigma_0 +m+2\nu, \R}(r_1)$ there is a unique solution  $U=\mF^{-1}(V) \in B_{\s_0, \R}(r)$  of the equation 
	$$	V= \mF(U)= \bF(U)U \,. $$
  Then we show that if $V$ has  higher regularity in $ \dot H^{s+\nu}_\R(\T; \C^2)$, $s > \s_0 + m + 2\nu$, then $U \in\dot H^{s}_\R(\T; \C^2)$.\\
Exploiting the linear  invertibility  of $\bF(U)$, we recast 
$V=\bF(U)U$ 
 as the fixed point problem 
\begin{equation}\label{ift}
	\cG(U;V):=\bF(U)^{-1} V = U \,.
	\end{equation}
First we  show that 
for any $V \in  B_{\s_0+m+{2\nu}, \R}(r_1)$, the map $U \mapsto \cG(U;V)$ is a contraction on the ball  $B_{\s_0, \R}(r)$ provided $r_1>0$ is small enough. \\
{\underline{$\cG(U;V)$ maps the ball into itself.}} 
Let $V \in  B_{\s_0+m+{2\nu}, \R}(r_1)$ and $U \in  B_{\s_0, \R}(r)$. By 
\eqref{lin.est.F} and $ \nu\leq m$, we have
\[
\norm{\cG(U;V)}_{\s_0} \leq C \left( \norm{V}_{\s_0+\nu} +  \norm{U}_{s_0} \norm{V}_{s_0} \right) \leq C r_1 (1+r) \leq r\,,
\]
which is verified provided $r_1$ is sufficiently small.\\
{\underline{$\cG(U;V)$ is a contraction.}}
Again let $V \in  B_{\s_0+m+ 2\nu, \R}(r_1)$ and $U_1, U_2 \in  B_{\s_0, \R}(r)$.
By \eqref{stima.dF} one has  
\begin{align*}\bF(U_1) - \bF(U_2) = \int_0^1 \di_U \bF(\tau U_1 + (1-\tau)U_2)[U_1 - U_2]  \di \tau  \,, 
\end{align*}  
which  applying $\bF(U_1)^{-1}$ to the left and $\bF(U_2)^{-1}$ to the right  yields
\begin{equation*}
    \bF(U_1)^{-1} - \bF(U_2)^{-1} = -\int_0^1 \bF(U_1)^{-1} \, \di_U \bF(\tau U_1 + (1-\tau)U_2)[U_1 - U_2]\,  \bF(U_2)^{-1}  \di \tau  \,.
\end{equation*}
Exploiting such formula, recalling the relation $s_0=\sigma_0-\vr$,  applying \eqref{stima.dF} with $s\leadsto \s_0+\nu\geq s_0+\vr$, we get 
\begin{align*}
 \|\cG(U_1;V) & - \cG(U_2;V)\|_{\sigma_0}  \stackrel{\eqref{lin.est.F}}{\leq} C 
 \sup_{\tau \in [0,1]}\norm{\di_U \bF(\tau U_1 + (1-\tau)U_2)[U_1 - U_2]\,  \bF(U_2)^{-1} V}_{\s_0+\nu} \\
 & \ \ \ + C  \norm{U_1}_{s_0}
 \sup_{\tau \in [0,1]}\norm{\di_U \bF(\tau U_1 + (1-\tau)U_2)[U_1 - U_2]\,  \bF(U_2)^{-1} V}_{s_0} 
 \\
 & \leq  C \left(1 + \norm{U_1}_{s_0} \right) \, \norm{U_1- U_2}_{\s_0+\nu - \vr} \, \norm{\bF(U_2)^{-1}V}_{\s_0 + \nu + m} \\
 & \stackrel{\eqref{lin.est.F}}{\leq}
 C \left( 1 + \norm{U_1}_{s_0} \right) \, \norm{U_1- U_2}_{\s_0} \left( \norm{V}_{\s_0+2 \nu + m} +  \norm{U_2}_{\s_0 + \nu + m - \vr} \norm{V}_{s_0} \right) \\
  & \leq C \left(1 + \norm{U_1}_{\s_0} \right) \, \left(1 + \norm{U_2}_{\s_0 + \nu + m - \vr} \right) \norm{U_1- U_2}_{\s_0} \,  \norm{V}_{\s_0+2\nu + m } \leq \frac12 \norm{U_1- U_2}_{\s_0} \,,
 \end{align*}
where in the last step we chose $\norm{V}_{\s_0+2\nu + m } \leq r_1$ small enough, we used the hypothesis \eqref{restriction:vr_nu_m} to bound $\norm{U_2}_{\s_0 + \nu + m - \vr} \leq \norm{U_2}_{\sigma_0} \leq r$, and we denote by $C$ a constant which changes from line to line.  
By Banach fixed point theorem, for any 
$V \in  B_{\s_0+m+{2\nu}, \R}(r_1)$, 
there is a unique 
$U \in B_{\s_0,\R}(r)$ solving the fixed point problem \eqref{ift}, and so we put 
\be\label{sol.imp}
\cF^{-1}(V) := U \, \quad \mbox{ so that }\quad \cG(\cF^{-1}(V); V) =  \cF^{-1}(V) \quad \forall\,V \in B_{\wt s_0,\R}(r_1)\,, \quad \wt s_0 := \s_0 + m + 2\nu\,. 
\ee
So far we have shown that
$\cF\circ \cF^{-1}(V) = V$ for any $V \in B_{\wt s_0,\R}(r_1)$.
Now we show that
$\cF^{-1}\circ \cF(U) = U$ provided $U \in B_{\wt s_0 {+ \nu},\R}(\wt r)$ and $\wt r\leq r$ is small enough. 
First of all, note that $ \cF^{-1} \circ \cF(U) $ solves the fixed point equation \eqref{ift} with $ V \leadsto \cF(U) $ and $ U \leadsto \cF^{-1} \circ \cF(U) $.
Provided $ \mF(U) \in B_{\wt s_0, \R}(r_1) $, the map $ \cG(\,\cdot\,;\cF(U)) $ is a contraction. As a result, the associated fixed point problem admits a unique solution, which must therefore coincide with $ U $.
We prove now that $ \mF(U) \in B_{\wt{s}_0, \R}(r_1) $. Indeed, 
estimate
 \eqref{lin.est.F}, for some $C>1$, gives
$$
\| \mF(U)\|_{\wt s_0}=\| \bF\left(U\right)U\|_{\wt s_0}\leq C\|U\|_{\wt s_0 + \nu}\leq r_1
$$
for any $U\in B_{\wt{s}_0 {+ \nu} ,\R}(\wt{r})$, choosing $\wt r:= r_1/C$. Then
\be\label{1112:1651}
 \cF \circ \cF^{-1}(V) = V\,, \quad \cF^{-1} \circ \cF (U) = U \qquad \forall  \, V \in B_{\wt s_0, \R}(\wt{r})\,,   \quad \forall \, U \in B_{ \wt s_0 + \nu, \R}(\wt{r})\,.
\ee
\underline{Upgraded regularity.} We now show that for any $s\geq \wt{s}_0 + \vr - (m + \nu) $ (see \eqref{sol.imp}), there exists $r'=r'(s) \in (0, r_1)$ such that
if $V\in B_{\wt{s}_0,\R}(r') \cap  \dot H^{s + \nu}_\R(\T; \C^2)$,  then  $\mF^{-1}(V)$ belongs to
$ \dot H^s_\R(\T; \C^2)$   and
\begin{equation}\label{vs.us}
 \|\cF^{-1}(V)\|_s \lesssim_s \|V\|_{s + \nu}\,.
\end{equation}
Fix the largest $\underline{n}\in \N$ so that  $s\geq  \s_0+\underline{n}\vr$, where we note that, since $s \geq \wt s_0 + \vr - (m + \nu),$ one has $\und{n} \geq 1$.
Then, from equation \eqref{ift} and estimate \eqref{lin.est.F}, and using also that,  from the fixed point argument,  $U\in B_{\s_0,\R}(r)$, arguing inductively on $n$ 
from the previous estimate we get
$$
\norm{U}_{\s_0+n\vr} = \norm{\bF^{-1}(U)V}_{\s_0+n\vr} \lesssim_n \big( \norm{V}_{\s_0+n\vr + \nu} + \norm{U}_{\s_0+(n-1)\vr} \norm{V}_{s_0} \big)  \ , \quad n =1 , \ldots, \underline{n} \,. 
$$
This shows that  $U\in \dot H^{\s_0+\underline{n}\vr}_{\R}(\T;\C^2)$ and, since $V\in B_{\wt{s}_0,\R}(r')\subseteq B_{{s}_0,\R}(r') $, we get
\begin{equation}\label{junatarella}  
\norm{U}_{\s_0+n\vr} \lesssim_{n} \norm{V}_{\s_0+n\vr + \nu}, \quad n=1,\ldots, \underline{n}\,.  
\end{equation}
Finally, using $ s-\vr \leq \s_0+\underline{n}\vr$ and again that
 $U\in B_{\s_0,\R}(r)$, $V\in B_{\wt{s}_0,\R}(r') \cap \dot H^{s + \nu}_{\R}(\T; \C^2)$, we deduce that for $r' \ll 1$  (depending on $s$) one has
$$
\| U\|_s 
= 
\| \bF^{-1}(U) V \|_s 
\stackrel{\eqref{lin.est.F}}{\lesssim_s} \left(\| V \|_{s + \nu} + \| U\|_{s-\vr}\|V\|_{s_0}\right) \quad \stackrel{\eqref{junatarella}}{\Rightarrow} \quad \|U\|_s \lesssim_s \| V\|_{s + \nu}\,,
$$
proving \eqref{vs.us}. Then the statement in item $(i)$ follows with $s_0':= \wt{s}_0 + \vr -(m + \nu) \geq \wt{s}_0$.\\
{\sc Proof of Item $(ii)$}  It follows from \eqref{sol.imp}  and \eqref{ift} that 
\begin{equation}\label{cF-1}
\cF^{-1}(V) = \bG(V) V \ , \quad  \bG(V):= \bF^{-1}(\cF^{-1}(V))\,.  
\end{equation}
We now verify that $\bG(V) \in \cM^\nu_{\geq 0}[r_1]$. Using \eqref{lin.est.F} and \eqref{vs.us}, for any $s \geq s'_0$ (see \eqref{1112:1651}) we have
\begin{equation}\label{stima:GdiV}
\|\bG(V) Z\|_{s - \nu} = \|\bF^{-1}(\cF^{-1}(V))Z\|_{s-\nu} \lesssim_s \|Z\|_{s} + \|\cF^{-1}(V)\|_{s - \nu -\vr} \|Z\|_{s_0} \lesssim_s \|Z\|_{s} + \|V\|_{s-\vr} \|Z\|_{{s}_0}\,.
\end{equation}
Using that $s_0 \leq s'_0$, this gives that the second of \eqref{piove} holds with $m \leadsto \nu$ and $s_0 \leadsto s'_0$, therefore $\bG(V) \in \cM^{\nu}_{\geq 0}[r_1]$.
Next we  show that $\bG(V)$ expands as in \eqref{exp.G}. 
In view of \Cref{lem:general_composition}-\Cref{comp-mappe-int}, we set
\[
\begin{aligned}
&\breve \cF^{-1}(V) := \breve \bF_0(V)V +\breve \bF_1(V)V+ \breve \bF_2(V)V,\qquad \breve\bF_0(V):=\uno, \quad\breve\bF_1(V):= - \bF_1(V) \in \wt \cM_1^m, \\
&\breve \bF_2(V)V:=- \bF_2(V) + \bF_1(V) \bF_1(V)+ \bF_1(\bF_1(V)V)\in \wt \cM_2^{2m}\,.
\end{aligned}
\]
Then, using the expansion 
$\cF(U) = \bF(U)U =  U + \bF_1(U) U + \bF_2(U)U + \bF_{\geq 3}(U)U$
we get 
\begin{equation}
\begin{aligned}
  &(\breve \cF^{-1}\circ \cF)(U) = U  + \bF'_{\geq 3}(U)U \,, \quad \mbox{ with }\\
  &\bF'_{\geq 3}(U):= \sum_{p=0}^2\sum_{\substack{j_1,\ldots,j_{p+1}\in \{0,1,2,3\}\\ j_1+\ldots+j_{p+1}+p\geq 3}} \breve \bF_p(\bF_{j_1}(U)U, \ldots ,\bF_{j_{p}}(U)U)\bF_{j_{p+1}}(U) \,
  \end{aligned}
  \label{mille.m}
\end{equation}
where in the above sum we also used the notations $ \bF_0:= \uno$ and $\bF_3:= \bF_{\geq 3}$.
Since, by the first bullet of \Cref{rmk:molto_bene}, $p$-homogeneous $m$-operators in $\wt\cM_p^{m}$ are also non-homogeneous $m$-operators in $\cM_{\geq p}^{m}[r]$, we apply \Cref{lem:general_composition}-\Cref{comp-no-homo} to get $\bF'_{\geq 3}(U)\in \cM^{3m}_{\geq 3}[r]$.
Substituting $U = \cF^{-1}(V)$ in the previous formula and using \eqref{cF-1},  we obtain
$$
\begin{aligned}
\cF^{-1}(V) &\stackrel{\eqref{mille.m}}{=} (\breve{\cF}^{-1} \circ \cF)(\cF^{-1}(V)) - \bF'_{\geq 3}(\cF^{-1}(V))\cF^{-1}(V)\\
& = V - \bF_1(V) - \bF_2(V) V + \bF_1(V) \bF_1(V) V + {\bF_1(\bF_1(V)V)V} + \bG_{\geq 3}(V) V\,, 
\end{aligned}
$$
with $ \bG_{\geq 3}(V):= -
\bF'_{\geq 3}(\bG(V)V) \bG(V)$ (use also \eqref{cF-1}).
First, \eqref{vs.us} ensure that there is $r''>0$ such that $\bG_{\geq 3}(V)$ is well defined for any $V \in B_{s_0+\nu,\R}(r'')$.
Moreover, as $\bG(V) \in \cM_{\geq 0}^\nu[r]$, one checks that  the map $\bG_{\geq 3}(V) \in \cM^{{3m} + \nu}_{\geq 3}[r'']$  iterating estimate \eqref{piove}.
This proves the expansion in 
\eqref{exp.G}.
\end{proof}

\noindent{\bf Flows of paradifferential and smoothing operators.}
 We now consider linear flows  of the form 
 \be\label{flussoG}
\begin{cases}
\pa_\tau {\bf \Phi}^\tau(U)= \bG(\tau,U){\bf \Phi}^\tau(U)\\
{\bf \Phi}^0(U)=\uno
\end{cases} \  \ 
 \ee
where $\bG(\tau, U)$ is a matrix of paradifferential operators. 
In the next lemma we give conditions on $\bG(\tau, U)$ so that the   flow map ${\bf \Phi}^\tau(U)$  is an admissible transformation for any $\tau\in [0,1]$.
 
\begin{lemma}\label{lem:flow.ad}
Let  ${\bf \Phi}^\tau(U)$ be the flow defined in \eqref{flussoG}, for $\tau \in [0, 1]$. 
\begin{enumerate}
    \item[$(i)$] If $\bG(\tau, U) = \vOpbw{ \frac{\beta(U;x)}{1+\tau \beta_x(U;x)}\ii \, \xi}$ 
    with $\beta\in \wt \mF_q^\R$ {and $q \in \{1, 2\}$},  then 
  ${\bf \Phi}^\tau(U)$  is a $({0},3)$-admissible transformation with gain $ \vr$ for any $\vr \geq 3$.
\item[$(ii)$]  If  $\bG(U) = \vOpbw{h(U;\cdot)}+\zOpbw{ g(U;\cdot)}$ with $h,\, g\in  \Sigma\Gamma_1^0[r_0,3]$ {for some $r_0>0$,  and there exist $C, s_0>0$} such that, for any $U, \widehat U \in B_{s_0,\R}(r_0),$ the functions 
 $${(x, \xi) \mapsto \di_U h_{\geq 3}(U; x, \xi)[\widehat U] \in \Gamma^m_{L^\infty}}, \qquad {(x, \xi) \mapsto \di_U g_{\geq 3}(U; x, \xi)[\widehat U] \in \Gamma^m_{L^\infty}}
 \,,
 $$
 with
 \begin{equation}\label{du.g3}
 \big| \di_U h_{\geq 3}(U; x, \xi)[\widehat U]\big|_{0,L^\infty, 4}+\big| \di_U g_{\geq 3}(U; x, \xi)[\widehat U]\big|_{0,L^\infty, 4} \leq C \| U\|^2_{s_0} \| \widehat U\|_{s_0}\,,
 \end{equation} then 
 ${\bf \Phi}^\tau(U)$ is a $({0},0)$-admissible transformation with gain $ \vr$ for any $\vr \geq 0$.
\item[$(iii)$] If $\bG(U) =  \vOpbw{\ii  f(U;\cdot)}$
with $f \in \wt \Gamma_1^\frac12$ and  $f=\ov{f}$, then 
  ${\bf \Phi}^\tau(U)$  is a $({0},\frac32)$-admissible transformation with gain $ \vr$ for any $\vr \geq \frac 3 2$.
\end{enumerate}
\end{lemma} 
 \begin{proof} 
 {\sc Well posedness, gain and linear invertibility:} We set  
 $m_\star=1$ in case $(i)$, and $m_\star=0$ in case $(ii)$. Under the assumption $(i)$ or $(ii)$ or $(iii)$, 
 there exist $s_0>0$ and, for any $s \in \R$, a positive  $r = r(s)$, such that  for any 
$ U\in B_{s_0,\R}(r)$ and $V \in \dot H^s_\R(\T;\C^2)$, 
 the flow ${\bf \Phi}^\tau(U)$
 is well defined as well as its linear inverse for any $\tau \in [0,1]$  and satisfies: 
\be \label{est.flow}
\sup_{\tau \in [0,1]}\norm{{\bf \Phi}^\tau(U)^{\pm 1}V}_s \leq C_s \norm{V}_s  \ , \qquad 
\sup_{\tau \in [0,1]}\norm{({\bf \Phi}^\tau(U)^{\pm 1} - \uno)V}_s \leq C_s 
\norm{U}_{s_0}\norm{V}_{s+m_*}
\ee
 see e.g.  Lemma 3.16 of \cite{BMM2} (with $k = K' =K=0$, item $(i)$ for the first estimate and item $(iii)$ with $p=1$ and $N=0$ for the second one). 
The latter estimate implies both 
\eqref{lin.est.F} (with $\nu =0$) and the second of \eqref{piove}, showing that  ${\bf \Phi}^\tau(U)$ is a matrix of $0$-operators in $\cM_{\geq 0}^0[r_0]$ for some positive $r_0$.\\
It remains to prove that ${\bf \Phi}^\tau(U)$ satisfies \eqref{esp:F}, \eqref{stima.dF} and \eqref{stima.d.adm}, namely expansion and differentiability properties. To this aim, we observe that case $(i)$ with $q=2$ follows from  of \cite[Lemma 2.15]{MM}. 
Since the remaining cases follow by the same argument, we restrict ourselves to the proof of case $(i)$ with $q=1$ and of case $(ii)$.

\smallskip
 \noindent{\sc Expansion:} We prove  that  ${\bf \Phi}^\tau(U)$ expands as in \eqref{esp:F}.  First expand 
 \begin{equation}\label{G:exp}
 \bG(\tau,U)=\bG_1(U) + \bG_2(\tau, U)+  \bG_{\geq 3}(\tau,U) \in \Sigma \cM_1^{m_\star}[r_0, 3]\,, 
 \end{equation}
 with 
 \begin{align}
 \notag
& \bG_1(U) := \begin{cases}
 	\vOpbw{ {\beta(U;x)}\, \ii\xi},\\
 	\vOpbw{ h_1(U;x,\xi)} + \zOpbw{g_1(U;x,\xi)}\,,
 \end{cases}  \!\!\!\!\bG_2(\tau, U) := \begin{cases}
 	\vOpbw{ {-\tau (\beta \beta_x)(U; x)}\,\ii \xi},\\
 	\vOpbw{ h_2(U;x,\xi)} +\zOpbw{ g_2(U;x,\xi)}\,, 
 \end{cases}\\
 & \qquad \qquad \qquad \qquad \qquad  \bG_{\geq 3}(\tau, U):= \begin{cases}
 	\vOpbw{ \frac{\tau^2 (\beta \beta^2_x)}{1 + \tau \beta_x}(U; x)\,\ii \xi},\\
 	\zOpbw{ h_{\geq 3}(U;x,\xi)} + \zOpbw{ g_{\geq 3}(U;x,\xi)}\, .
 \end{cases}
 \label{Gj}
 \end{align}
By Picard iteration 
\begin{align}
    {\bf \Phi}^\tau(U) =& \uno+ \int_{0}^\tau \bG(\theta,U){\bf \Phi}^\theta(U)\di \theta\label{centrino}\\
    =&\uno + \int_{0}^\tau \bG(\theta, U) \di \theta + \int_{0}^\tau \!\!\!\!\int_{0}^{\theta} \bG(\theta, U) \bG(\varsigma, U) \di \varsigma \di \theta +
    \int_{0}^\tau \!\!\!\!\int_{0}^{\theta}  \bG(\theta, U) \bG(\varsigma, U)  ({\bf \Phi}^{\varsigma}(U) - \uno) \di \varsigma \di \theta\,,\nonumber
\end{align}
which using the expansion of $\bG(\tau, U)$ in \eqref{G:exp} gives 
\be\label{centrone}
\begin{aligned}
   {\bf \Phi}^\tau(U)&= \uno + 
    {\bf \Phi}_1(\tau, U) +  {\bf \Phi}_2(\tau, U) +  {\bf \Phi}_{\geq 3}(\tau, U) 
\end{aligned}
\ee
with
$$
{\bf \Phi}_1(\tau, U) := \tau \bG_1(U) \in \wt\cM_1^{m_\star} \ , \quad 
{\bf \Phi}_2(\tau, U) := \int_0^\tau \bG_2(\theta, U) \di \theta + \frac{\tau^2}{2} \bG_1(U) \bG_1(U)\in \wt\cM_2^{2m_\star}
$$
and 
   \begin{align}
 {\bf \Phi}_{\geq 3}(\tau, U)  &:= \int_0^\tau \bG_{\geq 3}(\theta, U) \di \theta + \int_0^\tau \!\!\!\!\int_0^\theta (\bG_2(\theta, U) + \bG_{\geq 3}(\theta, U)) \bG(\varsigma, U) \di \varsigma \di \theta 
 \label{Phigeq3}
 \\
   &+ \int_0^\tau \!\!\!\! \int_0^\theta \bG_1(U) (\bG_2(\varsigma, U) + \bG_{\geq 3}(\varsigma, U)) \di \varsigma \di \theta +
  \int_{0}^\tau \!\!\!\!\int_{0}^{\theta}  \bG(\theta, U) \bG(\varsigma, U)  ({\bf \Phi}^{\varsigma}(U) - \uno) \di \varsigma \di \theta \,.\nonumber
   \end{align} 
 By substituting ${\bf \Phi}^\tau(U)-\uno$ with \eqref{centrino}, using \Cref{lem:general_composition} and estimate \eqref{est.flow}, ${\bf \Phi}_{\geq 3}(\tau; U)$ in \eqref{Phigeq3}  belongs to 
$\mM^{3m_\star}_{\geq 3}[r_0]$. 
This proves the expansion \eqref{esp:F} with $m =  3m_*$. 

\smallskip
\noindent{\sc Differentiability:}
We prove now Item $(iii)$ of Definition \ref{admtra}.
The 
 differential ${\rm d}_U \Phi^\varsigma(U)[\widehat U]$ fulfills the variational equation 
\be 
\begin{cases}
\partial_\varsigma {\rm d}_U {\bf \Phi}^\varsigma (U)[\widehat U] =\bG(\varsigma,U) {\rm d}_U {\bf \Phi}^\varsigma(U)[\widehat U]
+\di_U \bG(\varsigma,U)[\widehat{U}] {\bf \Phi}^\varsigma(U) \\
{\rm d}_U {\bf \Phi}^0(U)[\widehat U]=0 
\end{cases}  \ , 
\ee
whose solution is given by the  Duhamel formula 
\begin{equation}\label{centrale}
{\rm d}_U {\bf \Phi}^\varsigma(U)[\widehat U] = {\bf \Phi}^\varsigma(U) \int_0^\varsigma  {\bf \Phi}^\tau (U)^{-1} \ \di_U \bG(\tau,U)[\widehat{U}]\   {\bf \Phi}^\tau (U)\, {\rm d}\tau  \,. 
\end{equation}
We claim that, for all choices of $\bG(\tau, U)$ in \eqref{flussoG}, $\forall s\in \R$, 
\begin{align}  
&\sup_{\tau \in [0,1]}\norm{\di_U \bG_1(\tau,U)[\widehat{U}]W }_{s}\!\!\!\!\lesssim \|\widehat{U}\|_{s_0}\norm{W}_{s+m_\star}\,, \quad
\sup_{\tau \in [0,1]}\norm{\di_U \bG_2(\tau,U)[\widehat{U}]W }_{s} \!\!\!\!\lesssim \norm{{U}}_{s_0}\| \widehat {U} \|_{s_0}\norm{W}_{s+m_{\star}}\,, \label{claimata} \\
&\sup_{\tau \in [0,1]}\norm{\di_U \bG_{\geq 3}(\tau,U)[\widehat{U}]W }_{s} \lesssim \norm{U}_{s_0}^2 \| \widehat{U} \|_{s_0}\norm{W}_{s+m_{\star}}\,.\label{claimata1}
\  
\end{align}
Hence, we have
$
\sup_{\tau \in [0,1]}\norm{\di_U \bG(\tau,U)[\widehat{U}]W }_{s}\lesssim \|\widehat{U}\|_{s_0}\norm{W}_{s+m_\star} $ which, 
together with \eqref{est.flow} gives that, for any $s\in \R $, any $U \in B_{s_0, \R}(r)$,  $\widehat U \in \dot H^{s_0}_\R(\T; \C^2)$  and $V \in \dot H^{s+m_\star}_{\R}(\T; \C^2)$ one has
\be\label{stima:derflow}
\sup_{\varsigma \in [0,1]} 
\norm{{\rm d}_U {\bf \Phi}^\varsigma(U)[\widehat U]}_s \lesssim_s 
\| \widehat {U} \|_{s_0}\norm{W}_{s+m_\star} \ ; 
\ee
hence, for any  $\vr \geq 0$ and $s \geq s_0 + \vr$, we obtain estimate  \eqref{stima.dF}.

Let us immediately also prove \eqref{stima.d.adm}.
Differentiating ${\bf \Phi}_{\geq3}(\tau; U)$ in \eqref{Phigeq3} and using estimates  
\eqref{claimata} and \eqref{claimata1}, together with 
\eqref{est.flow}, \eqref{stima:derflow},  we obtain
\eqref{stima.dF} for $\di_U {\bf \Phi}_{\geq 3}(\tau; U)$ with $m\leadsto 3 m_\star$.

\smallskip 
We conclude by   proving \eqref{claimata}--\eqref{claimata1}. First of all, consider case $(ii)$, namely $ \bG(\tau,U)= \vOpbw{ h} + \zOpbw{ g}$.
By assumptions $g(U; \cdot) = g_1(U; \cdot) + g_2(U;\cdot) + g_{\geq 3}(U; \cdot)$ with $g_q \in \wt \Gamma_q^0$,  $q=1,2$, and $g_{\geq 3} \in \Gamma_{\geq 3}^0[r]$ fulfilling \eqref{du.g3};  an analogous decomposition holds for $h$.
Then  \eqref{claimata1} immediately follows from  \eqref{du.g3} and the continuity \Cref{thm:contS} $(i)$.
Regarding \eqref{claimata}, one has 
\small
\begin{gather*}
 \di_U \bG_1(\tau,U)[\widehat{U}]=\vOpbw{ h_1(\widehat{U};x,\xi)} + \zOpbw{ g_1(\widehat{U};x,\xi)}\,, 
 \\ 
 \di_U \bG_2(\tau,U)[\widehat{U}]= \vOpbw{ 2h_2(\widehat{U}, U;x,\xi)} +\zOpbw{ 2g_2(\widehat{U}, U;x,\xi)}\,,   
\end{gather*}
\normalsize
and  again \eqref{claimata} follows from  Theorem \ref{thm:contS}.

Next we analyze  case $(i)$, i.e.  $ \bG(\tau,U)=  \vOpbw{ \frac{\beta(U;x)}{1+\tau \beta_x(U;x)}\ii \xi}$ with $\beta \in \wt\cF_1^\R$. 
Differentiating $\bG_1(U)$, $\bG_2(\tau; U)$ and $\bG_{\geq 3}(\tau; U)$ in \eqref{Gj} we get
\be\label{diff.Gj}
\begin{aligned}
& \di_U \bG_1(U)[\hat U]= \vOpbw{ \beta(\widehat{U};x)\ii \xi} \ , \\
& 
\di_U \bG_2(\tau; U)[\hat U]=  
-\tau\vOpbw{ 
\big(\beta(\widehat{U};x)\beta_x(U;x) + \beta(U;x)\beta_x(\widehat{U};x) \big)\ii \xi} \ , \\
& 
\di_U \bG_{\geq 3}(\tau, U)[\hat U] =
\tau^2  \vOpbw{ 
 \wt g_{\geq 3}(\tau, U, \widehat U; x) \ii \xi} \,,
\end{aligned}
\ee
\normalsize
where
$$
\wt g_{\geq 3}(\tau, U, \widehat U; x):= \frac{\beta(\widehat{U}; x) \beta_x^2(U;x) + 2\beta(U;x)\beta_x(U;x) \beta_x(\widehat{U};x)}{1+\tau \beta_x(U;x)}
- \tau \frac{\beta(U;x) \beta_x(U;x) \beta_x(\widehat{U};x)}{(1+\tau \beta_x(U;x))^2}
\in L^\infty(\T; \R)\, . 
$$
Now notice that $\beta(\widehat{U};x)  \in \wt \cF^\R_{1}$, $\beta(\widehat U; x) \beta_x(U;x) + \beta( U; x) \beta_x(\widehat U;x) \in \wt \cF^\R_{2}$ and 
$\sup_{\tau \in [0,1]} \norm{\wt g_{\geq 3}(\tau, U, \widehat U; x)}_{L^\infty}
\lesssim \|{\widehat{U}}\|_{s_0}\norm{U}_{s_0}^2.$ 
Then Theorem  \ref{thm:contS} applied to \eqref{diff.Gj} proves the estimates 
 \eqref{claimata} and \eqref{claimata1} for case $(i)$.
 \end{proof}
Next, we consider the flow map  generated by a matrix of  smoothing operators:
 \be\label{flussoR}
\begin{cases}
\pa_\tau {\bf \Phi}^\tau(U)= \bR(U){\bf \Phi}^\tau(U)\\
{\bf \Phi}^0(U)=\uno
\end{cases} \quad \text{where}\quad
 \bR(U) \in \wt \cR^{-\varrho}_q, \quad q\in \{1,2\}\,. 
 \ee
\begin{lemma}\label{flow.s.ad}
Let $\vr>0$ and  $\bR(U) \in \wt \cR^{-\varrho}_q$,  $q\in \{1,2\}$. The flow ${\bf \Phi}^\tau(U)$ in \eqref{flussoR} is a $(0,0)$-admissible transformation with gain $\vr$.  
\end{lemma} 
 \begin{proof}
The statement for $q=2$ was proved in \cite[Lemma~2.17]{MM}, so we only consider the case $q=1$. 
Arguing as in \cite[Lemma 2.17]{MM}, one proves that for any $s \ge s_0+\vr$,
  $ U \in { B_{s_0, \R}(r_s)} \cap H^{s-\vr}_{\R}(\T; \C^2)$ with a sufficiently small $r_s>0$  {and} $V \in H^{s}_\R(\mathbb{T};\C^2)$
\be\label{est:smotR}
\sup_{\tau \in [-1,1]} \norm{ {\bf \Phi}^\tau(U)V }_{s}\leq 2  C_s  \left(\norm{ V}_{s } + \norm{V}_{s_0} \norm{U}_{s-\vr} \right),
\ee
            proving \eqref{lin.est.F} with $\nu = 0$, and that ${\bf \Phi}^\tau(U) \in \cM_{\geq 0}^0[r]$.\\
\smallskip
{\sc Expansion:} We prove  that  ${\bf \Phi}^\tau(U)$ expands as in \eqref{esp:F}.
By \eqref{centrino},  $ {\bf \Phi}^\tau(U)$ expands as \eqref{centrone} with 
\begin{align*}
& {\bf \Phi}_1(\tau, U) := \tau \bR_1(U) \in \wt\cR_1^{-\vr} \ , \quad 
{\bf \Phi}_2(\tau, U) := \frac{\tau^2}{2} \bR_1(U) \bR_1(U)\in \wt\cR_2^{-\vr} \ ,\\
& 
 {\bf \Phi}_{\geq 3}(\tau, U)  := 
  \int_{0}^\tau \!\!\!\!\int_{0}^{\theta}  \bR_1(\theta, U) \bR_1(\varsigma, U)  ({\bf \Phi}^{\varsigma}(U) - \uno) \di \varsigma \di \theta  \in\cR_{\geq 3}^{-\vr}[r]\,
\end{align*}
proving \eqref{esp:F}. \\
\smallskip
{\sc Differentiability:} We prove now Item $(iii)$ of Definition \ref{admtra}.
The differential  ${\rm d}_U {\bf \Phi}^\varsigma(U)[\widehat U]$ fulfills \eqref{centrale} with  $ \bG(\tau,U)\leadsto \bR(U)$ and, since  
 $ \di_U \bR(U)[\widehat U] =  \bR( \widehat U)$,
 we replace  \eqref{claimata} with  the bound 
 $
\norm{ \di_U \bR(U)[\widehat U]W}_s \lesssim_s \|{\wh U}\|_{s-\vr} \norm{W}_{s-\vr}$ 
 obtained from \eqref{smoothing} with $m =- \vr$ and $s-m \leadsto s$.
 Then both \eqref{stima.dF}
and then 
 \eqref{stima.d.adm} follow easily.
 \end{proof}

\section{Analysis of Resonances}\label{sec:analisi.res}
In this section we study second, third and fourth order resonances of  the linear frequencies $\Omega_j(\gamma)$ defined in \eqref{omegonejin}. 
In \Cref{sec:2wave} we characterize the $2$-waves non integrable interactions, namely collisions of frequencies of the form $\Omega_m(\gamma) = \Omega_n(\gamma)$ with $m \neq n$.
We also show that, for any value of $\gamma < 0$, there are no $3$-waves resonances,
i.e. integers $j_1, j_2, j_3 \in \Z_*$  and signs $\sigma_1, \sigma_2, \sigma_3 \in \{\pm\}$ fulfilling
$$
 \sigma_1 \Omega_{j_1}(\gamma) + \sigma_2  \Omega_{j_2}(\gamma)+ \sigma_3 \Omega_{j_3}(\gamma)   = 0 \ , \qquad \sigma_1 j_1 + \sigma_2 j_2 + \sigma_3 j_3 = 0 \,.
$$
The latter property was already considered in \cite[Section 2A]{IFT}. We include its proof for the sake of completeness.

The analysis of $4$-wave resonances,
 namely    integers $j_1, j_2, j_3, j_4 \in \Z_*$ and signs $\sigma_1, \sigma_2, \sigma_3, \sigma_4 \in \{\pm\}$, such that 
$$
 \sigma_1 \Omega_{j_1}(\gamma) + \sigma_2  \Omega_{j_2}(\gamma) + \sigma_3 \Omega_{j_3}(\gamma)  +\sigma_4 \Omega_{j_4}(\gamma) = 0 \ , 
 \qquad 
 \sigma_1 j_1 + \sigma_2 j_2 + \sigma_3 j_3 + \sigma_4 j_4 = 0 \ , 
$$
is much more involved, due to the presence of families of non-integrable resonances.

\noindent Also  when the vorticity  $\gamma < 0$ there are families of non-integrable resonances.  We do not classify them, but, as an  example, choose   $\tm < 0 < \tn$ with  $\frac{3\tn + \tm}{8} \in \Z\setminus \{0, 2\tm \}$ and $\tm+\tn>0$ and put the vorticity   $\gamma = - \sqrt{ \tfrac{(\tm+\tn)^2}{2(\tn-\tm)}}$; then 
\begin{equation}\label{0301:0122}
    - \Omega_{j_1}(\gamma) - \Omega_{j_2}(\gamma) +  \Omega_{j_3}(\gamma)  +\Omega_{j_4}(\gamma) = 0  \quad \mbox{ provided }
    \ \  (j_1, j_2, j_3, j_4) = \Big(\tm, \tm, \frac{3\tn + \tm}{8}, 2\tm - \frac{3\tn + \tm}{8}\Big) \,. 
\end{equation} 
We bypass this difficulty by retaining a weaker notion of integrability. Given $\gamma<0$ with $\gamma^2\in\Q\setminus{0}$, we construct sets $\Lambda=\{\tn,\tm\}\subset\Z_*$, $\tn\neq\tm$, which we call $\gamma$-good, such that $\Omega_\tn(\gamma)=\Omega_\tm(\gamma)$ and no non-integrable $4$-wave resonance involves $\tn,\tm$ and modes outside $\Lambda$. We further impose additional algebraic conditions on $\Lambda$, which are crucial in Section \ref{sec:eff.eq} (see \Cref{sec:4wave}). In \Cref{sec:4wave} we prove that all $4$-wave interactions involving two modes in $\Lambda$ are integrable (see \Cref{lem:wres}).
Along the section we  shall often use the lower bounds
\begin{equation}\label{Omega>0}
    \Omega_j(\gamma) \geq \und{c}(\gamma):=\sqrt{1+\frac{\gamma^2}{4}}- \frac{|\gamma|}{2} > 0\,, \quad  \forall j \in \Z_* \ , 
    \quad \forall \ \gamma < 0 \,. 
\end{equation}

\subsection{ Analysis of 2 and 3-waves resonances} \label{sec:2wave}
We begin to analyze $2$-wave resonances of the form $\Omega_m(\gamma) = \Omega_n(\gamma)$. 
The first result characterizes for which values of the parameters such resonances occur. 
\begin{lemma}[{\bf $2$-wave resonances}] \label{lem:2wave}
Let $\gamma <0$ and  $\{\Omega_j(\gamma)\}_{j \in \Z_*}$ be the frequencies in \eqref{omegonejin}. Then the following holds true:
\begin{itemize}
   \item[$(i)$] {\bf Resonances:} Let $m, n \in \Z_*$ with $m \leq n $, then
   \begin{equation}\label{2res}
\Omega_{m}(\gamma)  = \Omega_{n}(\gamma) \    \ \quad \Leftrightarrow \quad
\begin{cases}
\mbox{ either }  \ \ n = m \ \ \mbox{ or }\\
m<0<n \ \ \mbox{ and } \ \ \gamma^2 = \dfrac{(m + n)^2}{2(n-m)}\ \  \mbox{ and } \ \ n+m >0
\end{cases} \,.
    \end{equation}
    In particular, for any $n \in \Z_*$ there exists at most one integer $m \in \Z_*$, $m \neq n$, such that $\Omega_m(\gamma) = \Omega_n(\gamma)$.
    \item[$(ii)$] {\bf Spectral gaps:} Let    $\gamma^2\in \Q\setminus\{0\}$. 
    There is a constant $c=c(\gamma)>0$ such that  for any $ j, k \in \Z_*$ 
    \begin{equation}\label{2Wave_low}
        \Omega_{j}(\gamma)  \not= \Omega_{k}(\gamma) \quad \implies \quad \big | \Omega_{j}(\gamma)  - \Omega_{k}(\gamma)\big| \geq \frac{c}{|j|^\frac32+ |k|^\frac32}  \,  \cdot 
    \end{equation}
\end{itemize}
\end{lemma}
\begin{proof}
$(i)$  The function $x \mapsto \Omega_x(\gamma)$ is injective on the sets $ \{ x>0\}$ and $ \{ x<0\}$ separately, as it is clear from the definition \eqref{omegonejin}.
If  $\sign(m)=\sign(n)$, such  injectivity implies that equation  $\Omega_{m}(\gamma)  = \Omega_{n}(\gamma)$ holds if and only if $n=m$. 

Now consider the remaining possibility, $ m<0<n$. 
In this case, recalling \eqref{omegonejin}, $\Omega_{m}(\gamma)  = \Omega_{n}(\gamma)$  means to solve
\begin{equation}\label{solve.me.mn}
\sqrt{ -m +  \frac{\gamma^2}{4} } 
= 
\sqrt{ n+  \frac{\gamma^2}{4} } + \gamma \,.
\end{equation}
Since $\gamma < 0$, in order to have a solution we clearly need $m+n >0$.
Solving for $\gamma$, we get the expression in \eqref{2res}.
Vice versa, any $m, n$ fulfilling the right hand side of \eqref{2res} yields to a 2-wave resonance.
Finally, note that for any $n>0$ and $\gamma <0$ equation \eqref{solve.me.mn} has at most one solution in $m$, which concludes the proof of Item $(i)$.

$(ii)$ Without loss of generality we assume $ j>0 $ and we distinguish two cases according to the sign of $ k$. 

First, if $k>0$, by the hypothesis  $  \Omega_{j}(\gamma)  \not= \Omega_{k}(\gamma)$  one has $ |j-k|\geq 1$ and 
\be\label{casoeasy}
\big | \Omega_{j}(\gamma)  - \Omega_{k}(\gamma)\big|= \big | \sqrt{  j  +\frac{\gamma^2}{4}} -  \sqrt{  k  +\frac{\gamma^2}{4}}\big|= \frac{|j-k|}{ \sqrt{  j  +\frac{\gamma^2}{4}} + \sqrt{  k  +\frac{\gamma^2}{4}}}\geq \frac{c(\gamma)}{ \sqrt{  j} + \sqrt{  k  }}\geq \frac{c(\gamma)}{|j|^\frac32+ |k|^\frac32}.
\ee
Instead, if $k<0$,  we define 
\be\label{k*}
k^*:= -k - 2\gamma \sqrt{-k+\frac{\gamma^2}{4}}+\gamma^2 >0 \,. 
\ee
Note that $k^*$ is not necessarily integer, but in any case
 $ \Omega_{k^*}(\gamma)= \Omega_{k}(\gamma)$ and $ |k^*|\leq c(\gamma) |k|$. By the hypothesis $ \Omega_j(\gamma)\not= \Omega_k(\gamma)$ we deduce that $ k^*\not=j$, and, arguing as in \eqref{casoeasy} and using also $|k^*|\leq c(\gamma) |k|$, we get 
$$
0\not=\big | \Omega_{j}(\gamma)  - \Omega_{k}(\gamma)\big|= \big | \Omega_{j}(\gamma)  - \Omega_{k^*}(\gamma)\big|=  \big |\sqrt{  j  +\frac{\gamma^2}{4}} -  \sqrt{   k^*  +\frac{\gamma^2}{4}}\big|\geq c(\gamma) \frac{| j-k^*|}{| j|^\frac12+ |k|^\frac12}\,.
$$
It remains now to lower bound $ |j-k^*|$. We claim that there exists $c(\gamma)>0$ such that 
\be\label{claim:jstar} 
| j-k^*| \geq \frac{c(\gamma)}{|j|+|k|} \,, \quad \forall j \in \Z_* \,. 
\ee
Indeed, if  $  2\gamma \sqrt{-k+\frac{\gamma^2}{4}}= j+k-\gamma^2$ then, by the definition of $k^*$ in \eqref{k*} and recalling $k<0$,
$$
0\not= |j-k^*| = 2 \big|  2\gamma \sqrt{-k+\frac{\gamma^2}{4}}\big|\geq 2 \gamma^2
$$
and so \eqref{claim:jstar} holds true.      

If, instead, $  2\gamma \sqrt{-k+\frac{\gamma^2}{4}}\not= j+k-\gamma^2$ then there exists $C(\gamma)>0$ such that
$$0\not=\big| j+k-\gamma^2-2\gamma \sqrt{-k+\frac{\gamma^2}{4}}\big|\leq C(\gamma) (| j| +|k|)
$$
and finally rationalizing the numerator
$$
0\not= |j-k^*|
= \frac{\big| (j+k)^2- 2\gamma^2(j-k)\big| }{\big| j+k-\gamma^2-2\gamma \sqrt{-k+\frac{\gamma^2}{4}}\big|}
\geq \frac{1}{C(\gamma)} \frac{\big| (j+k)^2- 2\gamma^2(j-k)\big| }{|j|+|k|}\geq \frac{1}{C(\gamma)(|j|+|k|)}
$$
where, to get the last lower bound, we used that $ j\not=k^*$ and we write $ \gamma^2=\frac{p}{q}$ to obtain  that $q(j+k)^2 - 2p(j-k)$ is a nontrivial integer, so 
$$
0\not= \big| (j+k)^2- 2\gamma^2(j-k)\big|\geq \frac{1}{q}\,.
$$
In any case, \eqref{claim:jstar} is proved, as well as \eqref{2Wave_low}.
\end{proof}

Next we consider $3$-waves resonances. 
Denote by  
\begin{equation}
 \label{res3}
 \fR_3(\gamma) := \left\{
 ( \vec{\jmath} , \vec \sigma) \in \fP_3 
  \colon 
  \quad 
 \sigma_1 \Omega_{j_1}(\gamma) + \sigma_2  \Omega_{j_2}(\gamma)+ \sigma_3 \Omega_{j_3}(\gamma)   = 0
   \right\}  \ , \qquad \fP_3 \mbox{ in }  \eqref{mom1} \  , 
 \end{equation} 
 the subset of 3-waves resonant indexes.  We prove
\begin{lemma}[{\bf Absence of $3$-waves resonances}]\label{3onde}
The set $\fR_3(\gamma) $ in \eqref{res3} is empty for any $\gamma \in \R$. 
In particular, there exists  $\tc_3(\gamma)>0$ such that 
\begin{equation}
   \abs{ \sigma_1 \Omega_{j_1}(\gamma) + \sigma_2  \Omega_{j_2}(\gamma) + \sigma_3 \Omega_{j_3}(\gamma) } \geq  \tc_3(\gamma) \ , 
   \qquad \forall\, (j_1, j_2, j_3, \sigma_1, \sigma_2, \sigma_3) \in \fP_3  \,. 
   \label{3waves_lower}
\end{equation}
\end{lemma}
\begin{proof}
  In case  $\sigma_1 = \sigma_2 = \sigma_3$, we have that   $ |\sigma_1 \Omega_{j_1}(\gamma) + \sigma_2  \Omega_{j_2}(\gamma) + \sigma_3 \Omega_{j_3}(\gamma) | \geq 3 \und{c}(\gamma) >0$ in view of \eqref{Omega>0}.
    So let us consider the case $\sigma_1 = \sigma_2 = +$, $\sigma_3 = -$. The momentum condition becomes
$j_1 +j_2 = j_3$.
To compute the lower bound of $\Omega_{j_1}(\gamma) + \Omega_{j_2}(\gamma) - \Omega_{j_3}(\gamma)$, we first notice that 
\begin{equation}\label{omega3.lb}
    \omega_{j_1}(\gamma) + \omega_{j_2}(\gamma) - \omega_{j_1+j_2}(\gamma)
    \geq 2\sqrt{1+\frac{\gamma^2}{4}} - \sqrt{2+\frac{\gamma^2}{4}}=: \tf(\gamma)  \,, 
    \quad \forall\, j_1, j_2 \in \Z_* \,.
\end{equation}
Indeed, by triangular inequality and the monotonicity of the function $j \mapsto \omega_{j}(\gamma)$ one has 
\[
\omega_{j_1}(\gamma) + \omega_{j_2}(\gamma) - \omega_{j_1+j_2}(\gamma) 
\geq
\sqrt{|j_1| +\frac{\gamma^2}{4}} + \sqrt{|j_2| +\frac{\gamma^2}{4}}- \sqrt{|j_1| + |j_2| +\frac{\gamma^2}{4}} =: \Phi(|j_1|, |j_2|)\,.
\]
Then, one can easily see that the $\Phi$ is monotone non-decreasing in both the variables and therefore the minimum is achieved at $(|j_1|, |j_2|)=(1,1)$ where it takes value $\tf(\gamma)$ (see \eqref{omega3.lb}).
\\
In the following we shall use the even, strictly positive function 
$$
\tc_3(\gamma):=\tf(\gamma) - \frac{|\gamma|}{2} >0 \qquad \forall \, \gamma \in \R\,.
$$
We consider  several cases:\\
 Case 1. $j_1, j_2 >0$. Then $j_3>0$ and
    \begin{equation*}
          \Omega_{j_1}(\gamma) + \Omega_{j_2}(\gamma) - \Omega_{j_3}(\gamma)  \stackrel{\eqref{omegonejin}}{=}  \omega_{j_1}(\gamma) + \omega_{j_2}(\gamma) - \omega_{j_1+j_2}(\gamma) + \frac{\gamma}{2}\stackrel{\eqref{omega3.lb}}{\geq}  \tc_3(\gamma)\,.
    \end{equation*}
Case 2. $j_1<0$, $j_2 >0$ and $j_1 + j_2 >0$. Then 
    \begin{equation*}
        \Omega_{j_1}(\gamma) + \Omega_{j_2}(\gamma) - \Omega_{j_3}(\gamma)  \stackrel{\eqref{omegonejin}}{=}  \omega_{j_1}(\gamma) + \omega_{j_2}(\gamma) - \omega_{j_1+j_2}(\gamma) - \frac{\gamma}{2}\stackrel{\eqref{omega3.lb}}{\geq} \tf(\gamma) - \frac{\gamma}{2}\geq \tc_3(\gamma) \,. 
    \end{equation*}
Case 3. $j_1<0$, $ j_2 >0$ and $j_1 + j_2 <0$. Then 
    \begin{equation*}
        \Omega_{j_1}(\gamma) + \Omega_{j_2}(\gamma) - \Omega_{j_3}(\gamma) \stackrel{\eqref{omegonejin}}{=} \omega_{j_1}(\gamma) + \omega_{j_2}(\gamma) - \omega_{j_1+j_2}(\gamma) + \frac{\gamma}{2} \geq \tc_3(\gamma) \,. 
    \end{equation*}
Case 4. $j_1, j_2<0$. Then $j_1 + j_2 <0$ and 
    \begin{equation*}
        \Omega_{j_1}(\gamma) + \Omega_{j_2}(\gamma) - \Omega_{j_3}(\gamma)  \stackrel{\eqref{omegonejin}}{=}  \omega_{j_1}(\gamma) + \omega_{j_2}(\gamma) - \omega_{j_1+j_2}(\gamma) - \frac{\gamma}{2} \stackrel{\eqref{omega3.lb}}{\geq} \tc_3(\gamma)\,.
    \end{equation*}
This concludes the proof, since all other cases can be obtained either by permuting ${(\sigma_1,j_1), (\sigma_2,j_2), (\sigma_3,j_3)}$ or by flipping the signs: $(\sigma_1,\sigma_2,\sigma_3) \mapsto (-\sigma_1,-\sigma_2,-\sigma_3)$.
\end{proof}

\subsection{$\gamma$-good sets}\label{sec:gamma-good}
We begin with the definition of  $\gamma$-good sets.
\begin{definition}[{\bf $\gamma$-good set}] \label{g-good} 
Let $\gamma < 0$, $\gamma^2 \in \Q$. A set $\Set\subset \Z_*$ of the form
    \begin{equation}
  \label{Lambda}
  \Set := \{ \tm, \tn\} \,,\quad  \tm < 0 < \tn \,, \quad \tm + \tn > 0
  \end{equation}  
  is said to be {\sc $\gamma$-good} if the following holds true:
  \begin{itemize}
      \item[$(G1)$] $(${\sc $2$-wave resonant interaction:}$)$ $\gamma$, $\tm$ and $\tn$ fulfill  the resonance relation 
      \be\label{Omega*}\tag{G1}
      \gamma = -\frac{\tm + \tn}{ \sqrt{2(\tn - \tm)}} \,,  \ \ \  \mbox{  so that } \ \ 
      \Omega_\tm(\gamma) = \Omega_\tn(\gamma)  = \sqrt{\frac{\tn-\tm}{2}}=: \Omega_* \,;
      \ee
      \item[$(G2)$] $(${\sc $4$-waves resonances exclusion:}$)$ The following rational numbers are not integers
      \be \label{G2}\tag{G2}
\frac{3\tn + \tm}{8} \not\in \Z_{{>0}}\,,  \ \  \quad  \ \ 
\frac{3\tm +\tn}{8} \not\in \Z_{<0} \,; 
      \ee
      \item[$(G3)$] $(${\sc Sign condition:}$) $ 
    Consider the polynomial $p(\lambda):=  2 \lambda^6-25 \lambda^5 +50 \lambda^4 +8 \lambda^3 -4 \lambda^2 +\lambda$; one has the sign condition
    \begin{equation}
        \label{G3}\tag{G3}
        p\left(\frac{\tm}{\tn} \right)> 0\,;
    \end{equation}
      \item[$(G4)$] $(${\sc Counting:}$)$ One has
      \be\label{G4}\tag{G4}
      \tn \neq -2 \tm\,.
      \ee
  \end{itemize}
\end{definition}
Let us motivate such algebraic conditions.
Condition (G1) guarantees the existence of a non-integrable $2$-wave resonance, see Lemma \ref{lem:2wave}. 
Condition (G2) removes certain $4$-waves resonances, for example those appearing in 
 \eqref{0301:0122}. 
Condition (G3) is needed to guarantee that  the two key coefficients in \eqref{VmVn.def} have opposite signs.
Although condition (G4) is  not strictly necessary, it simplifies the combinatorics of \Cref{prop:cH^{(4)}_+}, in particular it excludes the presence of the monomial $(v)$ in \eqref{combinatorial:poisson}. As it is harmless, we include it.\\

\noindent
Despite its look, condition $(G3)$ is easy to verify:
\\
\begin{minipage}{0.38\linewidth}
\centering
\includegraphics[width=\linewidth]{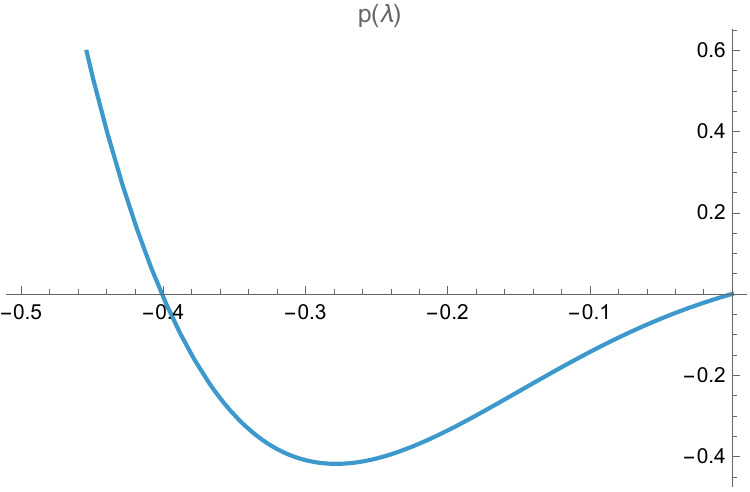}
\captionof{figure}{Graph of $p(\lambda)$.}
\label{plambda}
\end{minipage}
\hfill
\begin{minipage}{0.60\linewidth}
\begin{lemma}\label{lem:G3'}
There exists $\upsilon \in (2.493, 2.495)$ such that 
\be\label{G3'}
 \tn < -\upsilon\, \tm \quad \Leftrightarrow
 \quad  
 p\!\left(\frac{\tm}{\tn}\right) > 0 .
\ee
\end{lemma}

\begin{proof}
For $\lambda \leq 0$, the polynomial 
\[
p(\lambda)=2\lambda^{6}-25\lambda^{5}+50\lambda^{4}
+8\lambda^{3}-4\lambda^{2}+\lambda
\]
has exactly one negative root, that we denote by $-1/\upsilon$, lying in the interval $(-0.401,-0.4009)$. Moreover, $p(\lambda)>0$ in $(-\infty,0)$ if and only if $\lambda<-1/\upsilon$. See Figure~\ref{plambda}.
\end{proof}
\end{minipage}

\smallskip

Examples of $\gamma$-good sets $\Set$ for some $\gamma<0$ can be easily constructed by hand: for instance, the couple
$\Set = \{-2,\ 3\}$ is $\gamma-$good with $\gamma= -\frac{1}{\sqrt{10}}$.
Actually,  in the next result we show that for any value of the vorticity $\gamma<0$, provided $\gamma^2\in \Q$, there exist countable many 
$\gamma$-good sets. 
\begin{proposition}[{\bf Existence of $\gamma$-good sets}]\label{lem:goodset}
    For any $\gamma<0$ such that $\gamma^2\in \Q$, there exists a countable sequence of $\gamma-$good sets. More precisely, let $p,q \in \N$ be coprime numbers such that $\gamma^2 = \frac{p}{q}$.\\
    $1.$ If $p= 2^{\td} \tk$, with $2 \nmid \tk$, $\tk \in \N$ and $\td\geq 2$ even, then there exists $a_0>0$ such that the set $\Set_a := \{ (\tm_a, \tn_a) \}$ given by
    \begin{equation}
    \label{let.the}
        \tm_a := 2^{ \frac{\td}{2}} \tk a - \tk q a^2\,, \quad \tn_a := 2^{ \frac{\td}{2}} \tk a + \tk q a^2\,,
    \end{equation}
    is $\gamma-$good for any odd $a>a_0$.\\
    $2.$ If $p= 2^{\td} \tk$, with $2 \nmid \tk$, $\tk \in \N$ and $\td\geq 1$ odd, then there exists $a_0>0$ such that the set $\Set_a := \{ (\tm_a, \tn_a) \}$ given by
    \begin{equation}
    \label{number}
        \tm_a := 2^{ \frac{\td + 1}{2}} \tk a - 2\tk q a^2\,, \quad \tn_a := 2^{ \frac{\td + 1}{2}} \tk a + 2\tk q a^2\,,
    \end{equation}
    is $\gamma-$good for any odd $a>a_0$.\\
    $3.$ If $p$ is odd and $q \not \equiv_4 2$, then there exists $a_0>0$ such that the set $\Set_a := \{(\tm_a, \tn_a)\}$ given by
    \begin{equation}
    \label{theory}
    \tm_a := pa - pqa^2\,, \quad \tn_a := pa + pqa^2\,,
    \end{equation}
    is $\gamma-$good for any odd $a>a_0$.\\
    $4.$ If $p$ is odd and $q \equiv_4 2$, then there exists $a_0>0$ such that the set $\Set_a:= \{(\tm_a, \tn_a)\}$ given by
    \begin{equation}\label{run}
    \tm_a := \frac{p a}{2} -\frac{q pa^2}{4}\,, \quad \tn_a := \frac{pa}{2} + \frac{qpa^2}{4}\,,
    \end{equation}
    is $\gamma-$good for any odd $a> a_0$.
\end{proposition}

\begin{proof}
One easily sees that, in all cases,  $\tm_a<0 < \tn_a $ and $\tn_a + \tm_a >0$ provided $a$ is large enough. 
One also has $\tm_a, \tn_a \in \Z$. 
In the last case, just observe that when $q \equiv_4 2$, the number $\frac{q}{2}$ is odd. Hence $\tm_a$ and \(\tn_a\) in \eqref{run} are sums of two half–integers, which implies that they are integers.
Furthermore, in all cases \eqref{let.the}, \eqref{number},\eqref{theory}, \eqref{run} it is also immediate to check that $\tn_a \sim - \tm_a$ for $a$ large enough, which in particular implies that, for $a$ large enough, condition $(G4)$ is satisfied  and  $\tn < - \upsilon \tm$. By Lemma \ref{lem:G3'}, the latter condition implies that also condition $(G3)$ holds.
We now separately prove that conditions $(G1)$ and $(G_2)$ are satisfied in all four cases.\\
{\sc Case 1:} $p = 2^{\td} \tk$, with $\td\geq 2$ even. We have
$$
\Big(- \frac{(\tn_a + \tm_a)}{\sqrt{2(\tn_a-\tm_a)}}  \Big)^2 = 
\frac{(\tn_a + \tm_a)^2}{2(\tn_a-\tm_a)} = \frac{(2^{\frac{\td}{2} + 1} \tk a)^2}{4\tk q a^2} = \frac{2^{\td} \tk }{q} = \gamma^2\,,
$$
thus condition $(G1)$ is satisfied. We now check $(G2)$. Note that, since $p$ and $q$ are coprime, in this case $q$ is an odd number. Then for all odd $a$ we have
$$
3 \tn_a + \tm_a = 2^{\frac{\td}{2} + 2} \tk a + {2} \tk q a^2\,,
$$
where $2^{\frac{\td}{2} + 2} \tk a \equiv_8 0$ since $\td\geq 2$ and $\td$ even, and ${2} \tk a q^2 \not \equiv_8 0$, since $\tk, a$ and $q$ are all odd numbers. 
Thus $3\tn_a + \tm_a \not \equiv_8 0$, and $(3\tn_a + \tm_a)/8$ is not an integer number. Analogously, we have
\begin{equation}\label{3mn}
3 \tm_a + \tn_a = 2^{\frac{\td}{2} + 2} \tk a - {2}\tk q a^2 \equiv_8  - {2} \tk q a^2 \not \equiv_8 0\,,
\end{equation}
thus also $(3\tm_a + \tn_a)/8$ is not an integer. Thus condition $(G2)$ in this case is satisfied.\\
{\sc Case 2:} $p=2^{\td} k$, with $d\geq 1$ odd. Note that also in this case $q$ is an odd number, since $q$ and $p$ are coprime. Condition $(G1)$ holds since
$$
\Big(- \frac{(\tn_a + \tm_a)}{\sqrt{2(\tn_a-\tm_a)}}  \Big)^2 = 
\frac{(\tn_a + \tm_a)^2}{2(\tn_a-\tm_a)}= \frac{(2^{\frac{\td + 1}{2} + 1} \tk a)^2}{8\tk q a^2} = \frac{2^{\td} \tk}{q}= \gamma^2\,.
$$
To check $(G2)$, we observe that
$$
3 \tn_a + \tm_a = 2^{\frac{\td + 1}{2} + 2} \tk a + 4 \tk q a^2\,,
$$
where $2^{\frac{\td + 1}{2} + 2} \tk a \equiv_8 0$ since $\td + 1 \geq 2$ and $\td$ odd, and $4 \tk a q^2 \not \equiv_8 0$, since $\tk, a$ and $q$ are odd numbers. Thus $(3\tn_a + \tm_a)/8$ is not an integer number. For the same reason, one has
\begin{equation*}
3 \tm_a + \tn_a = 2^{\frac{\td + 1}{2} + 2} \tk a - 4 \tk q a^2 \equiv_{8}  - 4 \tk q a^2 \not \equiv_8 0\,,
\end{equation*}
thus $(3\tm_a + \tn_a)/8$ is not an integer number. Hence $(G2)$ in this case holds.\\
{\sc Case 3:} $p$ odd and $q\not \equiv_4 2$. In this case we have that
$$
\Big(- \frac{(\tn_a + \tm_a)}{\sqrt{2(\tn_a-\tm_a)}}  \Big)^2 = 
\frac{(\tn_a + \tm_a)^2}{2(\tn_a-\tm_a)} = \frac{(2p a)^2}{4 pq a^2}= \frac{p}{q} = \gamma^2\,,
$$
thus $(G1)$ is satisfied. We now check $(G2)$.\\
-- If $q \equiv_4 0$, namely $q=4\tilde q$ for some integer $\tilde q$, we have
\begin{equation}
\label{vediamo}
    3 \tn_a + \tm_a = 4pa + 2pqa^2 =  4 pa + 8 p\tilde q a^2 \equiv_8 4 pa \not \equiv_8 0\,,
\end{equation}
where the last passage holds since $p$ and $a$ are odd numbers.
For the same reason, one has
\begin{equation}
    \label{chi.sa}
    3 \tm_a + \tn_a = 4pa - 2pqa^2 =  4 pa - 8 p\tilde q a^2 \equiv_8 4 pa \not \equiv_8 0\,,
\end{equation}
 thus in this case $(G2)$ is satisfied.\\
-- If $q \equiv_4 1$, namely $q= 4 \tilde q + 1$ for some integer $\tilde q$, we have
$$
3 \tn_a + \tm_a = 4 pa + 2(4 \tilde q + 1) p a^2  \equiv_8 4 pa + 2 pa^2 \equiv_8 2pa (2 + a)\,.
$$
Now, since $p$ and $a$ are odd, $2pa(2 + a) \equiv_8 0$ would imply $2 + a \equiv_4 0$, which is not the case, since $a$ is odd. Thus $3 \tn_a + \tm_a \not \equiv_8 0$. Analogously, one proves
\begin{equation}
\label{magari}
3 \tm_a + \tn_a = 4 pa - 2(4 \tilde q + 1) p a^2  \not \equiv_8 0\,.
\end{equation}
This proves (G2).\\
-- If $q \equiv_4 3,$ namely $q =4 \tilde q + 3$ for some integer $\tilde q$, one has
$$
3 \tn_a + \tm_a = 4 pa + 2(4 \tilde q + 3) p a^2  \equiv_8 4 pa + 6 pa^2 \equiv_8 2pa (2 + 3a)\,,
$$
and since $p$, $a$ and $2 + 3 a$ are all odd numbers, one has $3\tn_a + \tm_a \not   \equiv_8  0$. Moreover, for the same reason
\begin{equation}
\label{puo.darsi2}
3 \tm_a + \tn_a = 4 pa - 2(4 \tilde q + 3) p a^2 \not \equiv_8 0\,,
\end{equation}
thus (G2) is proved. \\
{\sc Case 4:} $p$ odd and $q \equiv_4 2$. Then $q= 2\tilde q$, with $\tilde q$ an odd number. Then one has
$$
\Big(- \frac{(\tn_a + \tm_a)}{\sqrt{2(\tn_a-\tm_a)}}  \Big)^2 = 
\frac{(\tn_a + \tm_a)^2}{2(\tn_a - \tm_a)} = \frac{2(pa)^2}{2qpa^2} = \frac{p}{q} = \gamma^2\,,
$$
thus $(G1)$ is satisfied. We now check $(G2)$. One has
$$
3 \tn_a + \tm_a = 2 pa + \frac{q p a^2}{2} = 2 pa + \tilde q p a^2\,,
$$
which is an odd number, since $\tilde q, a$ and $p$ are odd, thus $3 \tn_a + \tm_a \not \equiv_8 0$. Analogously, one checks that
\begin{equation}\label{ancora.tu}
3 \tm_a + \tn_a = 2 pa - \frac{q p a^2}{2} = 2 pa - \tilde q p a^2 \not \equiv_8 0\,,
\end{equation}
being it also an odd number. Thus, $(G2)$ is satisfied also in this case.
\end{proof}

\subsection{ Analysis of 4-waves resonances}\label{sec:4wave}
In this section we study $4$-waves resonances. 
Define the set 
\begin{equation}\label{res1}
\fR_4:= \fR_4(\gamma) := \left\{
 ( \vec{\jmath} , \vec \sigma) \in \fP_4 : \ 
 \sigma_1 \Omega_{j_1}(\gamma) + \sigma_2  \Omega_{j_2}(\gamma) + \sigma_3 \Omega_{j_3}(\gamma)  +\sigma_4 \Omega_{j_4}(\gamma) = 0
   \right\}  \,, 
   \quad 
   \fP_4 \mbox{ in }  \eqref{mom1} \,, 
 \end{equation} 
 consisting in 4-waves resonant indexes.
Note that $\fR_4(\gamma)\neq \emptyset$, whatever the choice of $\gamma$: indeed, it contains always integrable resonances of the form    $(\ell, \ell, k, k)$ with corresponding signs $(+, -, +, -)$ and all their permutations. 
It also contains the Zakharov-Dyachenko non-integrable resonances when $\gamma =0$ and, for $\gamma < 0$, we exhibited in \eqref{0301:0122} an example of non-integrable $4$-waves resonance. 
 In what follows we want to study $4$-waves resonances when  certain indexes  among $(j_1, j_2, j_3, j_4)$ are constrained to belong to a set $\Lambda$. 
To formalize this, we give the following definition.
\begin{definition}Given a set $\Set \subseteq \Z$ and $n \in \{0, \ldots, 4\}$,  we  denote by 
$\fP_{\Set}^{(n)}$ the elements of $\fP_4$ (see \eqref{mom1}) having  exactly $n$ indexes {\em outside} the set $\Set$: 
\begin{equation}\label{momj}
\fP_\Set^{(n)} := \{ ( j_1, j_2, j_3, j_4 , \vec \sigma) \in \fP_4 \colon   \mbox{exactly } n\mbox{ indexes  among  } j_1, j_2, j_3, j_4 \mbox{ are  outside } \Set  \} \,. 
\end{equation}
We denote by $\fR_\Set^{(n)}(\gamma) $ the subset of $\fP_\Set^{(n)}$ made of $4$-waves resonances: with $\fR_4(\gamma)$ in \eqref{res1}, 
 \begin{equation}\label{resj}
\fR_\Set^{(n)}(\gamma) := \{ (  j_1, j_2, j_3, j_4 , \vec \sigma) \in \fR_4(\gamma)\colon  \mbox{exactly } n\mbox{ indexes  among }  j_1, j_2, j_3, j_4  \mbox{ are  outside } \Set  \} \,.
\end{equation} 
\end{definition}
In the next lemma we prove that, if the set $\Lambda$ is chosen to be $\gamma$-good according to \Cref{g-good}, the sets
$\fR_\Set^{(n)}$ with $n =0, 1, 2$ contain {\em only integrable resonances}. Recall that for any $\gamma <0$, $\gamma^2\in \Q$, by \Cref{lem:goodset} there exists countably many $\gamma$-good sets.
\begin{lemma}\label{lem:wres}
Let $\gamma <0$, $\gamma^2\in \Q$ and  $\Set = \{\tm, \tn\}$ be a $\gamma$-good  set of the form \eqref{Lambda}. 
Then the sets $\fP^{(n)}_\Set$, $\fR^{(n)}_\Set(\gamma)$,  defined in \eqref{momj} and \eqref{resj}, fulfill:
\begin{itemize}
\item[$(i)$] {The set $\fP^{(0)}_\Set \equiv \fR^{(0)}_\Set(\gamma)$ } and contains only integrable resonances: 
\be\label{R0}
\begin{aligned}
\fR_\Set^{(0)}(\gamma) & = \left\{ 
\big( \pi (\mathtt{k}, \mathtt{k}, \ell, \ell), \  \pi(+,-, +, -) \big), \colon \mathtt{k}, \ell\in \Set  \ , 
\ \  \pi \in \cS_4 \right\}  \ 
\end{aligned}
\ee
where $\pi$ is a permutation in the   symmetric group  $\cS_4$ of  permutations of  four symbols.
\item[$(ii)$] The set $\fR_\Set^{(1)}(\gamma) = \emptyset$. Moreover,
$\fP_\Set^{(1)}$ has finite cardinality and  there exists $c(\gamma) >0$ such that 
\begin{equation}
\label{lower.R1}
(\vec{\jmath}, \vec \sigma)  \in \fP_\Set^{(1)} \quad \Rightarrow
\quad
\abs{ \sigma_1 \Omega_{j_1}(\gamma) + \sigma_2  \Omega_{j_2}(\gamma) + \sigma_3 \Omega_{j_3}(\gamma)  +\sigma_4 \Omega_{j_4}(\gamma)  } \geq c(\gamma) \,.  
\end{equation}
\item[$(iii)$] The set $\fR_\Set^{(2)}(\gamma)$ contains only integrable resonances, i.e. 
\be
\begin{aligned}\label{fR2Lambda}
\fR_\Set^{(2)}(\gamma) & = \{ 
\big( \pi (\mathtt{k}, \mathtt{k}, \ell, \ell), \  \pi (+,-, +, -) \big)
 \colon  \ \mathtt k \in \Set, \  \ell  \in \Set^c , \  \pi \in \cS_4 \}  \,,
\end{aligned}
\ee
where again $\pi$ is a permutation in the   symmetric group  $\cS_4$ of  permutations of  four symbols.
\\
 Moreover, there exists $c(\gamma) >0$ such that 
{\small \begin{equation}
\label{lower.R2}
(\vec{\jmath}, \vec \sigma) \in \fP_\Set^{(2)}\setminus \fR_\Set^{(2)}(\gamma) \quad \Rightarrow
\quad
\abs{  \sigma_1 \Omega_{j_1}(\gamma) + \sigma_2  \Omega_{j_2}(\gamma) + \sigma_3 \Omega_{j_3}(\gamma)  +\sigma_4 \Omega_{j_4}(\gamma )} \geq \frac{c(\gamma)}{\max\limits_{a=1,\ldots, 4} ( \sqrt{|j_a|+\frac{\gamma^2}{4}}) } \,.  
\end{equation}}
\end{itemize}
\end{lemma}
\begin{proof} The proof involves certain computations that  have been verified using Mathematica; the relevant notebook \texttt{lemma3\_7.nb} can be found at the link  \url{https://git.sissa.it/amaspero/transfer-ww-vorticity}.

\smallskip
$(i)$ By the very definition of $\fP_\Set^{(0)}$  all indexes $j_1, j_2, j_3, j_4 \in \Set$. 
Without loss of generality, we reduce to study   one of the following cases. \\
\underline{Case I:} {$\sigma_1 = \sigma_2 = \sigma_3 = \sigma_4= +$. The momentum condition gives
$j_4 = - j_3 - j_2 - j_1$, meaning that $j_4 \in \{-3\tm,-2\tm-\tn,-\tm-2\tn,-3\tn\}$. 
Since $j_4 \in \Set$, this gives either $\tm = 0$, $\tn = 0$, $\tn +\tm = 0$, $\tm + 3\tn = 0$, $\tn = - 3 \tm$, all impossible  (the first three by  \eqref{Lambda} and the last one by \eqref{G3'}).}\\
\underline{Case II:}  $\sigma_1 = \sigma_2 = \sigma_3  = +$, $\sigma_4 = -$. The momentum condition gives $j_1 + j_2 + j_3 = j_4$, meaning that $j_4 \in \{ 3\tm, 
 2\tm + \tn, 2\tn + \tm, 3\tn\}$. Since $j_4 \in \Set$, this gives either $\tm = 0$, or $\tn = 0$, or $\tm + \tn =0$, or $3\tm = \tn$, or $3\tn = \tm$, all impossible in view of  \eqref{Lambda}.\\
\underline{Case III:} $\sigma_1 = \sigma_3 = +$, $\sigma_2 = \sigma_4 = -$. The momentum is $j_1 - j_2 + j_3 - j_4 = 0$. If $j_1 = j_2$ or $j_1 = j_4$, we are in the integrable case of \eqref{R0}. Otherwise, $j_4 \in \{2\tm - \tn, 2\tn - \tm  \} \cap \Set$, which yields $\tn = \tm$, a contradiction. \\
In conclusion, $\fP^{(0)}_\Set = \fR^{(0)}_\Set(\gamma)$ as claimed.\\

 $(ii)$ We can always assume that $j_1, j_2, j_3 \in \Set$ and $j_4 \in \Set^c$.
 Then, recalling \eqref{Omega*},
 \be\label{proof.lemma36.ii}
 \sigma_1 \Omega_{j_1}(\gamma) + \sigma_2  \Omega_{j_2}(\gamma) + \sigma_3 \Omega_{j_3}(\gamma)  +\sigma_4 \Omega_{j_4}(\gamma)  = (\sigma_1 + \sigma_2 + \sigma_3 )\Omega_*  + \sigma_4 \Omega_{j_4}(\gamma) \,. 
 \ee
 We analyze further cases.\\
 \underline{Case I:} $\sigma_1 + \sigma_2 + \sigma_3 = 3$.  The momentum is $j_1 + j_2 + j_3 + \sigma_4 j_4 =0$,  
 with  $j_1 + j_2 + j_3 \in \{ 3\tm, 2\tm + \tn, \tm + 2\tn, 3\tn \}$. We analyze every possibility:
 \begin{itemize}
     \item[I.a)]  Let $3\tm + \sigma_4 j_4 =0$. If $(j_4, \sigma_4) = (- 3\tm, +)$, then $3\Omega_* + \Omega_{-3\tm}(\gamma) > 4 \und{c}(\gamma) >0$ by \eqref{Omega>0}. If $(j_4, \sigma_4) = ( 3\tm, -)$, then, using the expressions in \eqref{omegonejin}, \eqref{Omega*}, 
     $$
     3\Omega_* - \Omega_{3\tm}(\gamma)  = -\frac{\sqrt{25 \tm^2-22 \tm \tn+\tn^2}+7 \tm-5 \tn}{2 \sqrt{2} \sqrt{\tn-\tm}}  \,, 
     $$
     which never vanishes for integers $\tm < 0<\tn$. 
     \item[I.b)] Let $2\tm + \tn + \sigma_4 j_4 = 0$. We can also assume that  $2\tm + \tn \neq 0$, otherwise $j_4 = 0$, not possible. The case $(j_4, \sigma_4) = (-2\tm - \tn, +)$ leads to 
     $3\Omega_* + \Omega_{-2\tm- \tn}(\gamma) > 4 \und{c}(\gamma) >0$ by \eqref{Omega>0}. In case  $(j_4, \sigma_4) = (2\tm + \tn, -)$ we have that 
     $$
     3\Omega_*- \Omega_{2\tm+ \tn}(\gamma) = 
     \begin{cases}
-\dfrac{\sqrt{-15 \tm^2+10 \tm \tn+9 \tn^2}+5 \tm-7 \tn}{2 \sqrt{2} \sqrt{\tn-\tm}}
& \mbox{ if } 2 \tm + \tn > 0\\
-\dfrac{\sqrt{17 \tm^2-6 \tm \tn-7 \tn^2}+7 \tm-5 \tn}{2 \sqrt{2} \sqrt{\tn-\tm}}
& \mbox{ if } 2 \tm + \tn < 0
     \end{cases}
     $$
  both of which  never vanish for integers $\tm < 0<\tn$. 
     \item[I.c)] Let $\tm + 2\tn + \sigma_4 j_4 = 0$. Again we can assume that  $\tm + 2\tn \neq 0$ to exclude $j_4 =0$.
      The case $(j_4, \sigma_4) = (-\tm - 2\tn, +)$ leads to 
     $3\Omega_* + \Omega_{-\tm- 2\tn}(\gamma) > 4 \und{c}(\gamma) >0$ by \eqref{Omega>0}. In case $(j_4, \sigma_4) = (\tm + 2\tn, -)$, we get
     $$
     3\Omega_*- \Omega_{\tm+ 2\tn}(\gamma) =
     -\frac{\sqrt{-7 \tm^2-6 \tm \tn+17 \tn^2}+5 \tm-7 \tn}{2 \sqrt{2} \sqrt{\tn-\tm}} 
     $$
          which never vanishes for integers $\tm < 0<\tn$. 
     \item[I.d)] Let $3\tn + \sigma_4 j_4 = 0$. 
     If $(j_4, \sigma_4) = (- 3\tn, +)$, then $3\Omega_* + \Omega_{-3\tn}(\gamma) > 4 \und{c}(\gamma) >0$ by \eqref{Omega>0}. 
 The case $(j_4, \sigma_4) = (3\tn, -)$ yields
     $$3\Omega_* - \Omega_{3\tn}(\gamma)  = 
     -\frac{\sqrt{\tm^2-22 \tm \tn+25 \tn^2}+5 \tm-7 \tn}{2 \sqrt{2} \sqrt{\tn-\tm}}
     $$
      which never vanishes for integers $\tm < 0<\tn$. 
 \end{itemize}
 \underline{Case II:} $\sigma_1 + \sigma_2 + \sigma_3 = 1$.  Without loss of generality  we assume $\sigma_1 = \sigma_3 = 1 = - \sigma_2$, so that the  momentum is $j_1 - j_2 + j_3 + \sigma_4 j_4 =0$ and 
 the resonant condition \eqref{proof.lemma36.ii} reduces to
 $$ \sigma_1 \Omega_{j_1}(\gamma) + \sigma_2  \Omega_{j_2}(\gamma) + \sigma_3 \Omega_{j_3}(\gamma)  +\sigma_4 \Omega_{j_4}(\gamma)  = \Omega_*  + \sigma_4 \Omega_{j_4}(\gamma) \,. $$
In addition, the number $j_1 - j_2 + j_3 \in \{ \tm, 2\tm - \tn, 2\tn - \tm, \tn \}$. We analyze every possibility:
 \begin{itemize}
     \item[II.a)] Let $\tm + \sigma_4 j_4 = 0$. Since $j_4 \in \Set^c$, the only possibility is $(j_4, \sigma_4) = (-\tm, +)$, but $\Omega_* + \Omega_{-\tm}(\gamma) > 2 \und{c}(\gamma) >0$ by \eqref{Omega>0}.
     \item[II.b)]  Let $2\tm - \tn + \sigma_4 j_4 = 0$. 
     If $(j_4, \sigma_4) = (\tn- 2\tm, +)$, we have
     $
     \Omega_* + \Omega_{\tn- 2\tm}(\gamma)> 2 \und{c}(\gamma) >0$ by \eqref{Omega>0}.
     In case $(j_4, \sigma_4) = (2\tm- \tn, -)$, we have
     $$
     \Omega_* - \Omega_{2\tm- \tn}(\gamma) = 
     -\frac{\sqrt{17 \tm^2-22 \tm \tn+9 \tn^2}+3 \tm-\tn}{2 \sqrt{2} \sqrt{\tn-\tm}}
     $$
      which never vanishes for integers $\tm < 0<\tn$. 
     \item[II.c)]
      Let $2\tn - \tm + \sigma_4 j_4 = 0$. 
       If $(j_4, \sigma_4) = (-2\tn+ \tm, +)$, we have
     $
     \Omega_* + \Omega_{\tm- 2\tn}(\gamma)> 2 \und{c}(\gamma) >0$ by \eqref{Omega>0}.
     In case $(j_4, \sigma_4) = (2\tn- \tm, -)$, we have
      $$
     \Omega_* - \Omega_{2\tn- \tm}(\gamma) = 
     -\frac{\sqrt{9 \tm^2-22 \tm \tn+17 \tn^2}+\tm-3 \tn}{2 \sqrt{2} \sqrt{\tn-\tm}}
     $$
 which never vanishes for integers $\tm < 0<\tn$. 
     \item[II.d)] Let $\tn + \sigma_4 j_4 = 0$. This case is analogous to II.a) and therefore  $\Omega_* + \sigma_4 \Omega_{j_4}(\gamma) \neq 0$.
 \end{itemize}
Observing that cases $\sigma_1 + \sigma_2 + \sigma_3 = -1,\ -3$ are the same as in Case I and II, this concludes the proof that  $\fR_\Set^{(1)}(\gamma) = \emptyset$. In addition, being  $\fP^{(1)}_\Set$ only finitely many cases, it shows \eqref{lower.R1}.\\

$(iii)$  Without loss of generality we   assume $j_1, j_2 \in \Set$, $j_3, j_4 \in \Set^c$. 
Then, recalling \eqref{Omega*},
\be\label{resIII}
\sigma_1 \Omega_{j_1}(\gamma) + \sigma_2  \Omega_{j_2}(\gamma) + \sigma_3 \Omega_{j_3}(\gamma)  +\sigma_4 \Omega_{j_4}(\gamma) = (\sigma_1  + \sigma_2 ) \Omega_* + \sigma_3 \Omega_{j_3}(\gamma)  +\sigma_4 \Omega_{j_4}(\gamma)\,.
\ee
We consider several cases.

\noindent \underline{Case I}: $\sigma_3 = \sigma_4 = +$.
Remark immediately that if 
$\sigma_1 = \sigma_2 = +$ then, by \eqref{Omega>0},
    $
 2\Omega_* +  \Omega_{j_3}(\gamma) +\Omega_{j_4}(\gamma) \geq 4\und{c}(\gamma)>0\,,
    $
    and similarly if $\sigma_1 = - \sigma_2 $, again by \eqref{Omega>0},
    $
 \Omega_{j_3}(\gamma) +\Omega_{j_4}(\gamma) \geq 2\und{c}(\gamma)>0\,.
    $\\
So from now on we consider 
\be\label{s1s2=-}
\sigma_1 = \sigma_2 = - \,, \qquad \sigma_3 = \sigma_4 = + \,. 
\ee
In this case  the momentum  and the resonant condition in \eqref{resIII} read 
\be\label{1911:1055}
j_1 + j_2 = j_3 + j_4 \ , \quad 
-2\Omega_* + \Omega_{j_3}(\gamma)  +\Omega_{j_4}(\gamma) \equiv  (\Omega_{j_4}(\gamma) - \Omega_*) + (\Omega_{j_3}(\gamma) - \Omega_*) \,. 
\ee
Note that, by \eqref{Lambda}, $|\tn| > |\tm|$. Therefore, if  both $|j_4|$ and $|j_3|$ are greater equal than $\tn+1$, by the monotonicity properties of the frequencies (see \eqref{omegonejin}) ande recalling $\gamma < 0$ we have $\Omega_{j_i}(\gamma) \geq \Omega_{|j_i|}(\gamma) \geq \Omega_{\tn + 1}(\gamma)$ for $i = 3,4$. Hence we conclude that
$$
 \Omega_{j_3}(\gamma)  - \Omega_* +\Omega_{j_4}(\gamma) - \Omega_* \geq 2(\Omega_{\tn+1}(\gamma)  - \Omega_*)> d(\gamma) := \frac{1}{2 \sqrt{\tn + 1 + \frac{\gamma^2}{4}}} > 0\,.
$$
From now on we let  $|j_3| \leq \tn$, the case $|j_4|\leq \tn$ being the symmetric one. Note that in this case, since $j_1, j_2 \in \Set$ and $|j_3| \leq \tn$, by the momentum condition also $j_4$ varies in a finite set; therefore, it is sufficient to prove that \eqref{resIII} never vanishes.\\
   Since $j_1, j_2 \in \Set$, using  \eqref{1911:1055} $j_4 +j_3 \in \{ 2\tm, \tm+\tn, 2\tn \}$. We consider each possibility.
    \begin{itemize}
        \item[I.a)] Let $j_4 = -j_3 + 2\tm$.  We consider other subcases according to the values of $j_3$, recalling $|j_3| \leq \tn$ and $\tm <0<\tn$, $\tn+\tm>0$.
        \begin{itemize}
            \item[I.a.1)]  If $ 1 \leq j_3 < \tn$,  ($j_3 = \tn$ is excluded because $j_3 \in \Lambda^c$), the resonant condition \eqref{resIII} with \eqref{s1s2=-} reduces to 
        \begin{align*}
    -2\Omega_* +&  \Omega_{j_3}(\gamma)  +\Omega_{2\tm - j_3}(\gamma)  
     = \frac{\sqrt{8(\tn - \tm) j_3 +17 \tm^2-14 \tm \tn+\tn^2}+\sqrt{8 j_3 (\tn-\tm)+(\tm+\tn)^2}+4 \tm-4 \tn}{2 \sqrt{2} \sqrt{\tn-\tm}}\,.
    \end{align*}
    The numerator is strictly increasing  as a function of $j_3$, being the 
    sum of  two strictly increasing  functions,  so it  can vanish only in one point.
    By a direct inspection, such point is 
$ j_3 = \dfrac{3\tn+\tm}{8} \in \Z_{{>0}}  $, 
   which is excluded by condition  (G2) of $\gamma$-good set, see \eqref{G2}. 
   Therefore the finitely many values $ \{-2\Omega_* +  \Omega_{j_3}(\gamma)  +\Omega_{2\tm - j_3}(\gamma) \}_{ 1 \leq j_3 \leq \tn-1} $ are all not zero. 
   \item[I.a.2)]  Next we consider the case  $2\tm < j_3 \leq -1$, with  $j_3 \neq \tm$ (since $j_3 \in \Set^c)$; in this case 
     \begin{align*}
    -2\Omega_* +&  \Omega_{j_3}(\gamma)  +\Omega_{2\tm - j_3}(\gamma)   = \frac{\sqrt{8(\tn - \tm) j_3 +17 \tm^2-14 \tm \tn+\tn^2}+\sqrt{-8 j_3 (\tn-\tm)+(\tm+\tn)^2}+6 \tm-2 \tn}{2 \sqrt{2} \sqrt{\tn-\tm}}  \,. 
    \end{align*}
    The derivative on the numerator as a function of $j_3$ is
    $$
 \frac{4 (\tn- \tm)}{ \sqrt{8 j_3(\tn - \tm) +17 \tm^2-14 \tm \tn+\tn^2}}-\frac{4 (\tn-\tm)}{\sqrt{-8 j_3 (\tn-\tm)+(\tm+\tn)^2}} \geq 0 \quad \Leftrightarrow
 j_3 \leq \tm  \ , 
 $$
 thus the numerator has a unique maximum at the point $j_3 = \tm$, at which  it vanishes. But the value $j_3  = \tm$  is  excluded since $j_3\in \Set^c$.    Therefore the finitely many values $ \{-2\Omega_* +  \Omega_{j_3}(\gamma)  +\Omega_{2\tm - j_3}(\gamma) \}_{ 2\tm < j_3 \leq -1} $ are all not zero.
    \item[I.a.3)]  Next consider the case  $-\tn  \leq   j_3 <  2\tm$ (note that $j_3 = 2\tm$ is not allowed, since it implies $j_4 =0$); then 
     \begin{align*}
    -2\Omega_* +&  \Omega_{j_3}(\gamma)  +\Omega_{2\tm - j_3}(\gamma)  = \frac{\sqrt{-8 j_3 (\tn - \tm) -15 \tm^2+18 \tm \tn+\tn^2}+\sqrt{-8 j_3 (\tn-\tm)+(\tm+\tn)^2}+4 \tm-4 \tn}{2 \sqrt{2} \sqrt{\tn-\tm}} \,.
    \end{align*}
     The numerator is strictly decreasing   as a function of $j_3$, being the 
    sum of the two strictly decreasing  functions,  so it  can vanish only at one point that, 
    by a direct inspection, results to be
   $j_3 = \dfrac{15\tm - 3\tn}{8} \in \Z$; however the constraint  $-\tn \leq j_3 \leq 2\tm $ forces  $3 \tm  + \tn \geq 0$, that contradicts \eqref{G3'} since $\upsilon < 3 $ and $\tm <0$ (and so it contradicts \eqref{G3}).
    Hence   the numbers $ -2\Omega_* +  \Omega_{j_3}(\gamma)  +\Omega_{2\tm - j_3}(\gamma) $ are different from zero for any finite allowed value of $|j_3| \leq \tn$. 
        \end{itemize}
          This concludes the analysis of Case I.a.
    \item[I.b)] Let $j_4 = -j_3 + \tm +\tn$. Remark that $j_3 \neq \tm + \tn$ otherwise $j_4 =0$. Again we consider subcases according to the values of $j_3$, recalling $|j_3| \leq \tn$ and $\tm <0<\tn$, $\tn+\tm>0$.
    \begin{itemize}
        \item[I.b.1)]   If $\tn +\tm < j_3 <  \tn$,  the resonant condition \eqref{resIII} with \eqref{s1s2=-} reduces to 
      \begin{align*}
    -2\Omega_* +&  \Omega_{j_3}(\gamma)  +\Omega_{\tm+\tn - j_3}(\gamma)  
        = \frac{\sqrt{8 j_3 (\tn-\tm)+(\tm+\tn) (9 \tm-7 \tn)}+\sqrt{8 j_3 (\tn-\tm)+(\tm+\tn)^2}+4 \tm-4 \tn}{2 \sqrt{2} \sqrt{\tn-\tm}} \,.
    \end{align*}
    The numerator is strictly increasing   as a function of $j_3$, being the 
    sum of the two strictly increasing  functions,  so it  can vanish only in one point that, 
    by a direct inspection, results to be 
    $j_3=\tn$, excluded since $j_3 \in \Set^c$.
    \item[I.b.2)]  If $ 0< j_3 < \tn +\tm$, we get
      \begin{align*}
    -2\Omega_* +&  \Omega_{j_3}(\gamma)  +\Omega_{\tm+\tn - j_3}(\gamma)   =\frac{\sqrt{-8 j_3 (\tn-\tm)-(7 \tm-9 \tn) (\tm+\tn)}+\sqrt{8 j_3 (\tn-\tm)+(\tm+\tn)^2}+2 \tm-6 \tn}{2 \sqrt{2} \sqrt{\tn-\tm}}
    \end{align*}
   The derivative of the numerator as a function of $j_3$ is 
   $$
   \frac{4 (\tn-\tm)}{\sqrt{8 j_3 (\tn-\tm)+(\tm+\tn)^2}} - \frac{4 (\tn-\tm)}{\sqrt{-8 j_3 (\tn-\tm)+(9 \tn-7\tm) (\tm+\tn)}} \geq 0 \ \ \Leftrightarrow \ \ j_3 \leq \frac{\tm + \tn}{2} \ , 
   $$
    thus the numerator has a unique maximum at the point $j_3 = \frac{\tm + \tn}{2}$, at which  it assumes the  value $$2\left(\sqrt{(5 \tn-3\tm) (\tm+\tn)}+\tm-3 \tn\right)$$ which is easily verified to be strictly negative.
 \item[I.b.3)]  If $ -\tn \leq  j_3 < 0$, $j_3\neq \tm$, we get
      \begin{align*}
    -2\Omega_* +  \Omega_{j_3}(\gamma)  +\Omega_{\tm+\tn - j_3}(\gamma)  
        & =\frac{\sqrt{-8 j_3 (\tn-\tm)-(7 \tm-9 \tn) (\tm+\tn)}+\sqrt{-8 j_3 (\tn-\tm)+(\tm+\tn)^2}+4 \tm-4 \tn}{2 \sqrt{2} \sqrt{\tn-\tm}}\,.
    \end{align*}
     The numerator is strictly decreasing   as a function of $j_3$, being the 
    sum of  two strictly decreasing  functions,  so it  can vanish only at one point that, 
    by a direct inspection, results to be $j_3=\tm$, that is excluded.
    \end{itemize}
  This concludes the analysis of Case I.b.

     \item[I.c)] Let $j_4 = -j_3 + 2\tn$. We consider other subcases according to the values of $j_3$, recalling $|j_3| \leq \tn$ and $\tm <0<\tn$, $\tn+\tm>0$.
     \begin{itemize}
         \item[I.c.1)] If $ 0 < j_3 < \tn$,  the resonant condition \eqref{resIII} with \eqref{s1s2=-} reduces to 
      \begin{align*}
    -2\Omega_* +&  \Omega_{j_3}(\gamma)  +\Omega_{2\tn - j_3}(\gamma)  
         =\frac{\sqrt{-8(\tn-\tm) j_3 +\tm^2-14 \tm \tn+17 \tn^2}+\sqrt{8  (\tn-\tm) j_3+(\tm+\tn)^2}+2 \tm-6 \tn}{2 \sqrt{2} \sqrt{\tn-\tm}} \,.
    \end{align*}
      The derivative on the numerator as a function of $j_3$ is
    $$
\frac{4 (\tn-\tm)}{\sqrt{8 j_3 (\tn-\tm)+(\tm+\tn)^2}}
-\frac{4 (\tn-\tm)}{\sqrt{-8 j_3 (\tn-\tm)+\tm^2-14 \tm \tn+17 \tn^2}}
\geq 0 \quad \Leftrightarrow
 j_3 \leq \tn \ , 
 $$
 thus the numerator has a unique maximum at the point $j_3 = \tn$, at which  it vanishes. But the value $j_3  = \tn$  is  excluded since $j_3\in \Set^c$.    Therefore the finitely many values $ -2\Omega_* +  \Omega_{j_3}(\gamma)  +\Omega_{2\tn - j_3}(\gamma)  $ for $j_3$ in the considered range are all not zero.
    \item[I.c.2)] If $ -\tn \leq  j_3 < 0$, $j_3 \neq \tm$, we get
      \begin{align*}
    -2\Omega_* + & \Omega_{j_3}(\gamma_*)  +\Omega_{2\tn - j_3}(\gamma_*)  \\
    &=
         \frac{\sqrt{-8(\tn - \tm) j_3 +\tm^2-14 \tm \tn+17 \tn^2}+\sqrt{-8 j_3 (\tn-\tm)+(\tm+\tn)^2}+4 \tm-4 \tn}{2 \sqrt{2} \sqrt{\tn-\tm}}
    \end{align*}
    The numerator is strictly decreasing   as a function of $j_3$, being the 
    sum of  two strictly decreasing  functions,  so it  can vanish only at one point that, 
    by a direct inspection, results to be
    $
    j_3= \dfrac{3\tm +\tn}{8} \in \Z_{< 0}$, 
    which is excluded by \eqref{G2}.
     \end{itemize}
     This concludes the analysis of Case I.c.
\end{itemize}
The analysis of Case I is concluded and we have proved that, however the signs $\sigma_1, \sigma_2$ are chosen,  there exists $c >0$ such that
\be\label{0301:0100}
\abs{(\sigma_1 + \sigma_2)\Omega_* + \Omega_{j_3}(\gamma)  +\Omega_{j_4}(\gamma) }\geq c \ , \qquad \forall j_3, j_4 \in \Set^c \,. 
\ee
\noindent \underline{Case II}:  $\sigma_3 = + $ and $\sigma_4 = -$.  This is more complex and we have several cases to check.\\
Case II.A)  $\sigma_1 = \sigma_2 = +$. The momentum  and the resonant condition in \eqref{resIII}  read
\be\label{CaseIImom} j_4-j_3 = j_1 + j_2  \in  \{2\tm, \tm +\tn, 2\tn \} , \ \ \ \ \ 
2\Omega_* + \Omega_{j_3}(\gamma) - \Omega_{j_4}(\gamma) \,. 
\ee
 We consider several subcases according to the signs of $j_3, j_4$.
\begin{itemize}
    \item[II.A.a)] Assume that  $j_3, j_4 > 0$.  Then, using \eqref{omegonejin}, 
    the resonance \eqref{CaseIImom} reads
    \begin{align}\label{res.caseIIa}
        2\Omega_* + \Omega_{j_3}(\gamma) - \Omega_{j_4}(\gamma)
        & 
        = 2\Omega_* +  \frac{j_3- j_4}{\sqrt{j_3 + \frac{\gamma^2}{4}} + \sqrt{j_4 +  \frac{\gamma^2}{4}}} \,.
    \end{align}
     If $j_3 \geq j_4 $ then  $2\Omega_* + \Omega_{j_3}(\gamma) - \Omega_{j_4}(\gamma) \geq 2 \und{c}(\gamma) >0$ by \eqref{Omega>0}. Therefore we restrict to  $0<j_3 < j_4$, and the momentum condition \eqref{CaseIImom} becomes   $j_4 - j_3\in  \{ \tm +\tn, 2\tn \}$ (namely we exclude $2\tm<0$).
     We treat the different subcases.
\begin{itemize}
    \item[II.A.a.1)] Let $j_4 - j_3 = \tm + \tn$. Then 
    \begin{align*}
        \eqref{res.caseIIa} & =  2\Omega_*  -  \frac{\tm + \tn}{\sqrt{j_3 + \frac{\gamma^2}{4}} + \sqrt{j_3 + \tm+\tn + \frac{\gamma^2}{4}}} 
        \geq 2\Omega_*  - \frac{\tm + \tn}{\sqrt{1 + \frac{\gamma^2}{4}} + \sqrt{1 + \tm+\tn + \frac{\gamma^2}{4}}}\\
        & \stackrel{\eqref{omegonejin}}{=}  2\Omega_* +\Omega_1(\gamma) - \Omega_{1+\tm+\tn}(\gamma) \\
        & =   \frac{\sqrt{2 (\tm+4) \tn+(\tm-8) \tm+\tn^2}-\sqrt{2 (\tm+4) \tn-\tm (7 \tm+8)+9 \tn^2}-4 \tm+4 \tn}{2 \sqrt{2} \sqrt{\tn-\tm}} > 0
    \end{align*}
    since 
    \begin{align*}
   \left( -4 \tm+4 \tn + \sqrt{2 (\tm+4) \tn+(\tm-8) \tm+\tn^2}\right)^2-(2 (\tm+4) \tn-\tm (7 \tm+8)+9 \tn^2) \\
   = 8(\tn-\tm)\left(\sqrt{2 (\tm+4) \tn+(\tm-8) \tm+\tn^2}-3 \tm+\tn\right) >0.
    \end{align*}
  \item[II.A.a.2)] Let    $j_4 - j_3 = 2\tn$. When  $j_3\geq R$ is sufficiently large, one has 
    \begin{align*}
    \eqref{res.caseIIa} & =  2\Omega_*   - \frac{2\tn }{\sqrt{j_3 + \frac{\gamma^2}{4}} + \sqrt{j_3 + 2\tn + \frac{\gamma^2}{4}}} > \Omega_* \stackrel{\eqref{Omega>0}}{>} \und{c}(\gamma) >0 \quad \mbox{ for } j_3 \geq  R \,.
    \end{align*}
    We now  consider the finitely many cases $0 < j_3 < R$, for which it is sufficient  to check that the expression in  \eqref{res.caseIIa} does not vanish. We have 
    \begin{align*}
    \eqref{res.caseIIa} & = \frac{-\sqrt{8 j_3 (\tn-\tm)+\tm^2-14 \tm \tn+17 \tn^2}+\sqrt{8 j_3 (\tn-\tm)+(\tm+\tn)^2}-4 \tm+4 \tn}{2 \sqrt{2} \sqrt{\tn-\tm}}  \,. 
    \end{align*}
    Assume by contradiction that  the numerator vanishes for some value of $j_3>0$,  then 
    \begin{align*}
        (-4 \tm+4 \tn -\sqrt{8 j_3 (\tn-\tm)+\tm^2-14 \tm \tn+17 \tn^2})^2 = 8 j_3 (\tn-\tm)+(\tm+\tn)^2\\
        \Leftrightarrow 
        8 (\tn-\tm) \left(\sqrt{-8 j_3 \tm+8 j_3 \tn+\tm^2-14 \tm \tn+17 \tn^2} +2(2\tn-  \tm) \right)  = 0
    \end{align*}
    which is never possible since $\tn \neq \tm$ and the two addends in the second factor are both strictly positive, absurd. 
\end{itemize}
\end{itemize}

\begin{itemize}
    \item[II.A.b)] Assume  $j_3> 0$ and $ j_4 <0$;  using \eqref{omegonejin}, 
    the resonance \eqref{CaseIImom} reads
       \be\label{CaseIIb}
       2\Omega_* + \Omega_{j_3}(\gamma) - \Omega_{j_4}(\gamma) 
        = 2\Omega_* +\gamma +  \frac{j_3+j_4}{\sqrt{j_3 + \frac{\gamma^2}{4}} + \sqrt{-j_4 +  \frac{\gamma^2}{4}}}  \,. 
        \ee   
In this case $j_4- j_3 <0$, so,  from the momentum condition \eqref{CaseIImom}, the only possibility is 
$j_4 - j_3 = 2\tm. $ 
This also confines  $0< j_3<-2\tm$, $j_3 \neq \tn$, so  we are reduced to finitely many cases, and we have only to check that \eqref{CaseIIb} does not vanish. 
Substituting $j_4 = j_3 + 2\tm$ we have 
 \begin{align*}
   \eqref{CaseIIb} 
    =
    \frac{-\sqrt{-8(\tn-\tm)j_3 - 16 \tm(\tn-\tm)+(\tm+\tn)^2}+\sqrt{8 j_3 (\tn-\tm)+(\tm+\tn)^2}+2 (\tn-3 \tm)}{2 \sqrt{2} \sqrt{\tn-\tm}} \,. 
    \end{align*}
 The numerator is strictly increasing   as a function of $j_3$, being the 
    sum of two strictly increasing  functions. So its minimum is attained at $j_3 = 0$, where it is strictly positive. Thus the functions $  2\Omega_* + \Omega_{j_3}(\gamma) - \Omega_{j_3+2\tm}(\gamma) $ never vanish when $0< j_3<-2\tm$.
\end{itemize}
 
\begin{itemize}
    \item[II.A.c)] Assume $j_3< 0$ and $ j_4 > 0$. In this case $j_4-j_3 >0$, so by the momentum condition \eqref{CaseIImom} $j_4 - j_3 \in \{ 2\tn, \tm + \tn\}$.
    \begin{itemize}
        \item[II.A.c.1)] Let $j_4 - j_3 = 2\tn$. Since $j_4>0$, $j_3$ is confined in the interval  $-2\tn < j_3 <0$, and again we are restricted to analyze finitely many cases. 
     The resonance in \eqref{CaseIImom} becomes 
        \begin{align*}
     2\Omega_* + \Omega_{j_3}(\gamma) - \Omega_{j_3+2\tn}(\gamma)
    =
   \frac{-\sqrt{8 j_3 (\tn-\tm)+\tm^2-14 \tm \tn+17 \tn^2}+\sqrt{-8 j_3 (\tn-\tm)+(\tm+\tn)^2}-2 \tm+6 \tn}{2 \sqrt{2} \sqrt{\tn-\tm}}    \,. 
    \end{align*}
 The numerator is strictly decreasing   as a function of $j_3$, being the 
    sum of two strictly decreasing  functions. So its minimum in the interval  $-2\tn < j_3 <0$ is attained at $j_3 = -1$, where it is strictly positive. 
     Thus the functions $  2\Omega_* + \Omega_{j_3}(\gamma) - \Omega_{j_3+2\tn}(\gamma) $ never vanish when $-2\tn < j_3 <0$.
     
 \item[II.A.c.2)] Let $j_4 = j_3 + \tn+\tm$. Since $j_4>0$, $j_3$ is confined in the interval   $-\tm-\tn < j_3 <0$, $j_3 \neq \tm$, and again we are reduced to study  finitely many cases. The resonance in \eqref{CaseIImom} becomes 
 \begin{align*}
     2\Omega_* + \Omega_{j_3}(\gamma) - \Omega_{j_3+\tn+\tm}(\gamma)
    =\frac{-\sqrt{8 j_3 (\tn-\tm)+(9 \tn- 7\tm) (\tm+\tn)}+\sqrt{-8 j_3 (\tn-\tm)+(\tm+\tn)^2}-2 \tm+6 \tn}{2 \sqrt{2} \sqrt{\tn-\tm}} \,. 
    \end{align*}
   The numerator is strictly decreasing   as a function of $j_3$, being the 
    sum of two strictly decreasing  functions. So its minimum in the interval  $-\tm-\tn < j_3 <0$ is attained at $j_3 = -1$, where it is strictly positive. 
     Thus the functions $  2\Omega_* + \Omega_{j_3}(\gamma) - \Omega_{j_3+\tn+\tm}(\gamma) $ never vanish when $-\tm-\tn < j_3 <0$.
    \end{itemize}
\end{itemize}

\begin{itemize}
    \item[II.A.d)] Assume $j_3, j_4< 0$. Then, using \eqref{omegonejin}, 
    the resonance \eqref{CaseIImom} reads
\be\label{CaseIId}
 2\Omega_* + \Omega_{j_3}(\gamma) - \Omega_{j_4}(\gamma)
        = 2\Omega_* + \frac{j_4 - j_3}{\sqrt{-j_3 +\frac{\gamma^2}{4}} + \sqrt{-j_4 + \frac{\gamma^2}{4}}} \,. 
\ee
If $j_4 - j_3>0$, clearly   $2\Omega_* + \Omega_{j_3}(\gamma) - \Omega_{j_4}(\gamma) \geq 2 \und{c}(\gamma) >0$ by \eqref{Omega>0}.
So we restrict to $j_4 - j_3<0$, which by the momentum \eqref{CaseIImom} reduces to $j_4 = j_3 + 2\tm$. If $j_3 \leq  -R$ sufficiently large, 
\begin{align*}
\eqref{CaseIId}  =  2\Omega_* + \frac{2\tm}{\sqrt{-j_3 + \frac{\gamma^2}{4}} + \sqrt{-j_3 -2\tm + \frac{\gamma^2}{4}}} > \Omega_* \stackrel{\eqref{Omega>0}}{>} \und{c}(\gamma) >0 \quad \mbox{ for } j_3 \leq - R  \,. 
\end{align*}
   We now  consider the finitely many cases $-R < j_3 < 0$, for which it is sufficient  to check that the expression in  \eqref{CaseIId} does not vanish.
   Recalling that $j_4 = j_3 + 2\tm$ we get
$$
       2\Omega_* + \Omega_{j_3}(\gamma) - \Omega_{j_3+2\tm}(\gamma)  = \frac{-\sqrt{-8 j_3 (\tn-\tm)+17 \tm^2-14 \tm \tn+\tn^2}+\sqrt{-8 j_3 (\tn-\tm)+(\tm+\tn)^2}-4 \tm+4 \tn}{2 \sqrt{2} \sqrt{\tn-\tm}}  \,. 
       $$
Assume by contradiction that the numerator vanishes for some value of $  j_3 <0$,  then 
\begin{align*}
(-4 \tm + 4 \tn + \sqrt{8 j_3 (\tm - \tn) + (\tm + \tn)^2})^2  = 17 \tm^2 + 8 j_3 (\tm - \tn) - 14 \tm \tn + \tn^2
\\    \Leftrightarrow 
(\tn -\tm) \sqrt{8 j_3 \tm-8 j_3 \tn+\tm^2+2 \tm \tn+\tn^2}-2 \tm \tn +2 \tn^2
 = 0
\end{align*}
which is never possible since all the addends are strictly positive. 
\end{itemize} 

This concludes the analysis of Case II.A); in particular  we have proved that  there exists $c >0$ such that 
    \be\label{CaseII.finest1}
    \abs{2\Omega_* + \Omega_{j_3}(\gamma) - \Omega_{j_4}(\gamma)} \geq c >0 \ ,  \ \ \forall j_3, j_4 \in \Set^c \,.
    \ee 
Case II.B)   $\sigma_1 = \sigma_2 = -$. The momentum  and the resonant conditions  
in \eqref{resIII}  read
\be\label{CaseIImom2} 
j_3-j_4 = j_1 + j_2  \in  \{2\tm, \tm +\tn, 2\tn \} , \ \ \ \ \ 
2\Omega_* + \Omega_{j_4}(\gamma) - \Omega_{j_3}(\gamma) \,, 
\ee
hence this is Case II.A with the role of $j_3$ and $j_4$ reversed.
The analysis is the same and we obtain
 \be\label{CaseII.finest2}
    \abs{2\Omega_* + \Omega_{j_4}(\gamma) - \Omega_{j_3}(\gamma)} \geq c >0 \,,  \ \ \forall j_3, j_4 \in \Set^c \,.
    \ee 
Case II.C) 
 $\sigma_1 = - \sigma_2$. The momentum and resonant conditions in \eqref{resIII} read 
 \be\label{CaseIIbis}
j_4 - j_3 = \sigma_1 (j_1 - j_2) \in \{ \tm - \tn, \tn - \tm, 0\}  \ , \qquad 
\Omega_{j_3}(\gamma) - \Omega_{j_4}(\gamma)  
 \ee
 We distinguish several subcases according to the signs of $j_3$ and $j_4$.
 \begin{itemize}
     \item[II.C.a)] Assume $j_3 = j_4$. Then the momentum \eqref{CaseIIbis} forces  $j_1 = j_2$, leading to an integrable quadruplet of indexes  of the form $(j_1, j_2, j_3, j_4) = (\tk, \tk, \ell, \ell)$, $\tk \in \Set$, $\ell \in \Set^c$,  $(\sigma_1, \sigma_2, \sigma_3, \sigma_4) = (\pm, \mp, +, -)$. 
     Taking all possible permutations one gets the set $\fR_\Set^{(2)}(\gamma)$ of \eqref{fR2Lambda}.
 \end{itemize}

\begin{itemize}
    \item[II.C.b)] Assume $j_3 j_4 >0 $, $j_3 \neq j_4$. Then $\abs{|j_3| - |j_4|}=|\tn - \tm|$ by the first of  \eqref{CaseIIbis}, and the resonance in \eqref{CaseIIbis} becomes 
\be\label{CaseII.bis.f}
\abs{\Omega_{j_3}(\gamma)  - \Omega_{j_4}(\gamma) }= 
 \frac{|\tn-\tm|}{\sqrt{|j_3| +\frac{\gamma^2}{4}}+\sqrt{|j_4| +\frac{\gamma^2}{4}}}\geq \frac{c}{\max\left\{|j_3|^{\frac12}, |j_4|^{\frac12} \right\}}\,. 
\ee
\end{itemize}

\begin{itemize}
    \item[II.C.c)] Assume  $j_3<0$ and $j_4>0$. Then $j_4-j_3 >0$ and the momentum \eqref{CaseIIbis} forces   $j_4 =j_3 + \tn - \tm$ and, being $j_4>0$,  $j_3$ to be confined in  $\tm - \tn < j_3 < 0$. Then 
\[
\Omega_{j_3}(\gamma)  - \Omega_{j_3 + \tn - \tm}(\gamma)  = \frac{-\sqrt{8 j_3 (\tn-\tm)+9 \tm^2-14 \tm \tn+9 \tn^2}+\sqrt{8 j_3 (\tm-\tn)+(\tm+\tn)^2}+2 \tm+2 \tn}{2 \sqrt{2} \sqrt{\tn-\tm}} \,. 
\]
The numerator, as a function of $j_3$, is strictly decreasing since it is the sum of  two decreasing functions.
 Moreover, one verifies that it vanishes only at the point $j_3 = \tm$, which is excluded since $j_3 \in \Lambda^c$. Hence the numbers  $\{\Omega_{j_3}(\gamma)  - \Omega_{j_3 + \tn - \tm}(\gamma)\}_{\tm-\tn < j_3 < 0}$ never vanishes.
\end{itemize}

\begin{itemize}
    \item[II.C.d)] Assume  $j_3>0$ and $j_4<0$. 
     Then $j_4-j_3 <0$ and the momentum \eqref{CaseIIbis} forces   $j_4 =j_3 + \tm - \tn$ and, since $j_4 <0$,  $j_3$ is  confined in  $0 < j_3 < \tn - \tm$. 
     Then 
\[
\Omega_{j_3}(\gamma)  - \Omega_{j_3+\tm - \tn}(\gamma) = \frac{-\sqrt{-8 (\tn-\tm) (j_3+\tm-\tn)+(\tm+\tn)^2}+\sqrt{8 j_3 (\tn-\tm)+(\tm+\tn)^2}-2 \tm-2 \tn}{2 \sqrt{2} \sqrt{\tn-\tm}} \,. 
\]
The numerator, as a function of $j_3$, is the sum of two strictly increasing functions, so it is strictly increasing and can vanish only in one point that, by a direct inspection, is $j_3 = \tn$. This contradicts $j_3 \in \Lambda^c$.  
Hence the numbers  $\{\Omega_{j_3}(\gamma)  - \Omega_{j_3 + \tm - \tn}(\gamma)\}_{0 < j_3 < \tn - \tm}$ never vanishes.
\end{itemize}
This concludes the analysis of Case II.C; in all cases we proved that  there exists $c >0$ such that
\be\label{0301:0058}
\abs{\Omega_{j_3}(\gamma)  - \Omega_{j_4}(\gamma) }
\geq \frac{c}{\max\left\{|j_3|^{\frac12}, |j_4|^{\frac12} \right\}}\ , \quad \forall j_3 \neq j_4  \  , \ \ j_3, j_4 \in \Set^c  \,. 
\ee
Then estimate \eqref{lower.R2} follows from \eqref{0301:0100}, \eqref{CaseII.finest1}, \eqref{CaseII.finest2}, \eqref{0301:0058}.
\end{proof}

\section{Strong $\Lambda$-normal Form and Identification}\label{sec:formal}
The goal of this section is to transform the water waves Hamiltonian \eqref{H.gamma} into its \emph{formal} quartic Birkhoff normal
form. In particular, we are going to proceed as follows: first we shall remove all the cubic monomials by exploiting the absence of three-wave interactions (see
Lemma \ref{3onde}). Due to the  presence of non integrable 4-wave resonances, we cannot reduce the quartic terms to integrable ones; however, we shall prove that, if the set $\Set$ is $\gamma$-good (recall Definition \ref{g-good}), the quartic terms of the Hamiltonian can be reduced to a suitable weaker notion of normal form, that we call strong $\Set$-normal form, see \Cref{def:wr}. In addition, we shall compute the coefficients of the quartic monomials whose frequencies are supported on $\Set$. This is the content of Section \ref{sec:bnf.formal}.

Note that at this level we do not guarantee the convergence of the Birkhoff normal form; however, in Subsection \ref{sec:identif}
we are going to exhibit an identification argument which will enable us to identify the resonant contributions of the
formal Birkhoff normal form that we exhibit in Subsection \ref{sec:bnf.formal} with the ones of the para-differential normal form that we are
going to compute in Section \ref{sec:paraWW}.


\smallskip 
\noindent{\bf Bracket of Hamiltonians and vector fields:} Given two functions $\scF$, $\scG \colon H^1(\T;\C^2) \to \R$, we define the {\em Poisson brackets}
\be\label{PP}
\{ \scF , \scG \} := \frac{1}{2\pi \im} \sum_{k \in \setminus \{0\} } \left( \pa_{\zetina_k} \scG \, \pa_{\bar{\zetina}_k} \scF - \pa_{\bar \zetina_k} \scG \pa_{\zetina_k} \scF \right) \,.
\ee
The {\em Hamiltonian vector field} of a function $\scF$ is
\be\label{hvf}
[\hamvec{\scF}]_k^\sigma := 
- \frac{\im \sigma }{2\pi} \pa_{\zetina_k^{-\sigma}} \scF \qquad \forall\, k \in \Z_*\,,\ \forall\, \sigma \in \{\pm\}\,.
\ee
With this notation and recalling the definition of the complex symplectic form \eqref{Fex}, one has $$\{\scF, \scG\} =\Omega_\C(\hamvec{\scG},\hamvec{\scF})= \di \scG [\hamvec{\scF}]  = \langle \grad \scG, \hamvec{\scF} \rangle.$$
The {\em Lie bracket} of two vector fields is defined as
\be\label{lb}
\bral X, Y \brar (\zak):=  \di Y(\zak)[X(\zak))] - \di X(\zak)[Y(\zak))] \,, 
\ee
and one has the identity
\be\label{lb.id}
\bral \hamvec{\scG}, \hamvec{\scF} \brar  = \hamvec{\{\scG, \scF\}}  \,.
\ee
\smallskip
\noindent{\bf Projections of quartic Hamiltonians and cubic vector fields.}
We introduce now projections of quartic Hamiltonians and cubic  vector fields on the sets $\fP_\Lambda^{(n)} $ and $\fR_\Lambda^{(n)}$ defined in \eqref{momj} and \eqref{resj}.
Recall that a real, translation invariant,  quartic Hamiltonian expands as 
\begin{align}
    H_4(U)=
    2\pi \sum_{(j_1,\ldots,j_4,\sigma_1,\ldots,\sigma_4) \in\fP_4}H_{j_1,j_2,j_3,j_4}^{\sigma_1,\sigma_2,\sigma_3,\sigma_4}u_{j_1}^{\sigma_1}
    u_{j_2}^{\sigma_2}u_{j_3}^{\sigma_3}u_{j_4}^{\sigma_4}, \qquad \ov{H_{j_1,j_2,j_3,j_4}^{\sigma_1,\sigma_2,\sigma_3,\sigma_4}}=H_{j_1,j_2,j_3,j_4}^{-\sigma_1,-\sigma_2,-\sigma_3,-\sigma_4}
\end{align}
and that any real-to-real cubic  vector field $X(U)$,  translation invariant,  expand in Fourier as (see \eqref{polvect})
\begin{equation}
\label{X3.fou}
X(U)^\sigma = 
\!\!\!\!\!\!
\sum_{(j_1,j_2,j_3, k,  \sigma_1,\sigma_2,\sigma_3, - \sigma) \in \fP_4} \!\!\!\!\!\!
X_{\ j_1, j_2, j_3, k }^{\sigma_1, \sigma_2, \sigma_3, \sigma} \ u_{j_1}^{\sigma_1} \, u_{j_2}^{\sigma_2} \, u_{j_3}^{\sigma_3}\, { e^{\im \sigma k x }}  \ , 
\quad 
X_{\ j_{\pi(1)}, \ldots,j_{\pi(3)}, k}^{ \sigma_{\pi(1)}, \ldots, \sigma_{\pi(3)},\sigma} 
=  
X_{\ j_{1}, \ldots, j_{3}, k}^{ \sigma_{1}, \ldots, \sigma_{3},\sigma} \,, 
\end{equation}
for any permutation $ \pi $ of $ \{1, 2, 3 \} $.
Given a 
subset $A \subseteq \fP_4$, we denote by $\Pi_A H$ and  $\Pi_A X$ the Hamiltonian and the vector field  obtained restricting the indexes to belong to $A$, namely 
\begin{equation}
\label{proiettore.sui.siti}
\begin{aligned}
(\Pi_A H)(U) &= 2\pi \sum_{(j_1,\ldots,j_4,\sigma_1,\ldots,\sigma_4) \in A}H_{j_1,j_2,j_3,j_4}^{\sigma_1,\sigma_2,\sigma_3,\sigma_4}u_{j_1}^{\sigma_1}
    u_{j_2}^{\sigma_2}u_{j_3}^{\sigma_3}u_{j_4}^{\sigma_4}\,,
\\
(\Pi_A X)(U)^\sigma &:= 
\!\!\!\!\!\!\!
\sum_{(j_1,j_2,j_3, k, \sigma_1,\sigma_2,\sigma_3, - \sigma)  \in A} 
\!\!\!\!\!\!\!
X_{\ j_1, j_2, j_3, k }^{\sigma_1, \sigma_2, \sigma_3, \sigma} \ u_{j_1}^{\sigma_1} \, u_{j_2}^{\sigma_2} \, u_{j_3}^{\sigma_3}\,  { e^{\im \sigma k x }}   \,. 
\end{aligned}
\end{equation}
\subsection{Formal Birkhoff normal form}\label{sec:bnf.formal}
We now introduce the definition of weak and strong-$\Set$ normal form.
\begin{definition}[{\bf $\Set$-normal form}] 
\label{def:wr}
Let $\gamma < 0$, $\gamma^2 \in \Q$ and 
$\Set$  be a $\gamma$-good set according to Definition \ref{g-good}.
Let $\fP^{(n)}_\Set$, $\fR^{(n)}_\Set(\gamma)$ be as in \eqref{momj} and \eqref{resj}.
A cubic, translation  invariant vector field $X(Z)$ is said to be in 
\begin{itemize}
\item {\em weak-$\Set$ normal form} if  all its  monomials with at most two indexes outside $\Set$ are resonant, i.e. 
\begin{equation}\label{wnf}
 \Pi_{\fP^{(n)}_\Set} X = \Pi_{\fR^{(n)}_\Set} X  \ , \quad n = 0,1,2 \,;
\end{equation}
\item {\em strong-$\Set$ normal form} if  in addition there are no resonant  monomials with exactly  one  index outside $\Set$, and those with exactly  $2$ indexes outside $\Set$ are integrable, i.e. 
\begin{align}
\label{snf1}
& \Pi_{\fP^{(0)}_\Set} X = \Pi_{\fR^{(0)}_\Set} X  \,,  \quad 
 \Pi_{\fP^{(1)}_\Set} X  = 0 \,, \\
 \label{snf2}
 & \Re \left(  [\Pi_{\fP_\Set^{(2)}}X(Z)]_j^+    \bar{z_j} \right)= 0 \quad \forall\, j \in \Z_* \,. 
\end{align}
\end{itemize}
\end{definition}

In the next proposition we put the Hamiltonian of water waves, written in the complex Zakharov coordinates $\zak$ (see \eqref{zak},  \eqref{zetone:formale}), into its formal  Birkhoff normal form, and show that the quartic term of the new Hamiltonian is
in strong-$\Lambda$ normal form, provided the set $\Lambda$ is chosen to be $\gamma$-good.
More importantly, we compute the coefficients of the quartic normal form 
$\Pi_{\fR_\Set^{(0)}} \cX_{\cH_4^+}$, see
\eqref{H4+.str0}.
\begin{proposition}[{\bf Formal normal form}]\label{prop:cH^{(4)}_+}
    Let $\cH$ be the Hamiltonian in \eqref{zetone:formale}. There exists a \emph{formal} symplectic transformation $\Phi$ such that
\begin{equation}\label{new.birkh}
        \cH \circ \Phi = \cH_2 + \cH_4^+ + {\cH}^+_{\geq 5}\,,
    \end{equation}
    where $\cH_2$ is the quadratic Hamiltonian
\be\label{H2}
       \cH_2(Z)\equiv \cH_2(z, \bar z):= 2\pi \sum_{k \in \Z_*} \Omega_k(\gamma) |z_k|^2 \,, \quad \Omega_k(\gamma) \mbox{ in } \eqref{omegonejin}  \,, 
\ee
$\cH_4^+$ is a homogeneous quartic Hamiltonian and ${\cH}^+_{\geq 5}$ vanishes of order $5$ at $\zak = 0$.

In addition, if $\Set =\{\tm,\, \tn\}$ is a $\gamma$-good set (cf. \Cref{g-good}), the cubic 
Hamiltonian vector field
$\cX_{\cH^+_4}$ is  in strong-$\Lambda$ normal form (cf. \Cref{def:wr}) and explicitly
\begin{equation}
\label{H4+.str0}
 (\Pi_{\fP^{(0)}_\Set} \cX_{\cH^+_4}) (Z) = 
\begin{pmatrix}
-\im 
\left(2\fa |z_\tm|^2 + \fb |z_\tn|^2 \right) z_\tm e^{\im \tm x} 
-\im 
\left( 2\fc |z_\tn|^2  + \fb |z_\tm|^2 \right) z_\tn e^{\im \tn x} 
\\
 \im 
\left(2\fa |z_\tm|^2 + \fb |z_\tn|^2 \right) \bar{z_\tm} e^{-\im \tm x} 
+\im 
\left( 2\fc |z_\tn|^2  + \fb |z_\tm|^2 \right) \bar{z_\tn} e^{-\im \tn x}  
\end{pmatrix}   \,, 
\end{equation}
where
\begin{equation}\label{abc}
    \begin{aligned}
        \fa &:= \frac{\tm^3 \left(17 \tm^3-9 \tm^2 \tn+3 \tm \tn^2-3 \tn^3\right)}{4 (\tm-\tn)^2 (\tn - 3 \tm)}\,, 
        \quad 
         \fc := \frac{\tn^3 \left(3 \tm^3-3 \tm^2 \tn+9 \tm \tn^2-17 \tn^3\right)}{4 (\tm-3 \tn) (\tm-\tn)^2} \,, \\
        \fb & :=  \frac{\tm \tn \left(-19 \tm^5+133 \tm^4 \tn-222 \tm^3 \tn^2+106 \tm^2 \tn^3-31 \tm \tn^4+\tn^5\right)}{4 (\tm-3 \tn) (\tm-\tn)^2 (3 \tm-\tn)} \,,
    \end{aligned}
\end{equation}
whereas
\begin{equation}
\Pi_{\fP_\Set^{(1)}} \cX_{\cH^+_4}=0, \qquad
 \Re\left( [\Pi_{\fP_\Set^{(2)}}\cX_{\cH_4^+}(Z)]^+_j    \, \bar z_j \right) = 0 \quad \forall \, j \in \Z_* \,,
\label{H4+.str}
\end{equation}
\end{proposition}
The rest of the section is devoted to the proof of \Cref{prop:cH^{(4)}_+}.\\
The first step is to compute the  water waves Hamiltonian $\scH_\gamma(\eta, \psi)$ in  \eqref{H.gamma}  up to quartic terms in the complex variables $
\zak
= \begin{pmatrix} \zetina \\ \bar \zetina \end{pmatrix}
$
defined in \eqref{zak}. This is done in the next lemma.
\begin{lemma}\label{lem:expham}
The water waves Hamiltonian 
in the complex linear Zakharov variables $\zak$ \emph{(}see \eqref{zak},  \eqref{zetone:formale}\emph{)} Taylor expands as 
\begin{equation}
\cH (\zak)= \cH_2(\zak) +  \cH_3(\zak)  +\cH_4(\zak) +  \cH_{\geq 5}(\zak) \,,
\end{equation}
where \\
$\bullet$ $\cH_2$ is the quadratic Hamiltonian in \eqref{H2};\\
$\bullet$  $\cH_3$ is  the cubic  real-valued Hamiltonian 
\begin{align}\label{H3.form}
\cH_3(\zak) \equiv \cH_3(\zetina, \bar \zetina) & = 2\pi \sum_{ \vec \sigma \cdot \vec k = 0 } 
H_{k_1, k_2, k_3}^{\sigma_1, \sigma_2, \sigma_3} \, 
\zetina_{k_1}^{\sigma_1} \zetina_{k_2}^{\sigma_2}\zetina_{k_3}^{\sigma_3} \,,
\\
 \label{H3.coeff}
 H_{k_1, k_2, k_3}^{\sigma_1, \sigma_2, \sigma_3} & := \frac12\left(k_2 k_3 \sigma_2 \sigma_3 + |k_2| \, |k_3| \right) \sigma_2 \sigma_3 \fm_{k_1} \fn_{k_2} \fn_{k_3} 
  + \sigma_2 \sigma_3  \frac{\gamma}{2} \sign(k_2) |k_3|  \fm_{k_1} \fm_{k_2} \fn_{k_3}
\end{align}
where the real, even numbers $\fm_j$, $\fn_j$ are defined as
\be\label{al.be.m}
\fm_j:= \fm_j(\gamma):= \frac{1}{\sqrt 2}
\left( \frac{4 j^2}{4 |j| + \gamma^2}\right)^{\frac14}
\,, \quad
\fn_j:= \fn_j(\gamma):= \frac{1}{\sqrt 2} 
\left( \frac{4 j^2}{4 |j| + \gamma^2}\right)^{-\frac14}
\,,
\quad 
\quad j \in \Z_* \,;
\ee
$\bullet$ 
 $\cH_{4}$ is the  quartic real-valued Hamiltonian
\begin{align}\label{H4.form}
\cH_4(\zak) \equiv \cH_{4}(\zetina, \bar \zetina)& := 2\pi \sum_{\vec \sigma \cdot \vec k = 0} 
H_{k_1, k_2, k_3, k_4}^{\sigma_1, \sigma_2, \sigma_3, \sigma_4} \, \zetina_{k_1}^{\sigma_1} \zetina_{k_2}^{\sigma_2}\zetina_{k_3}^{\sigma_3} \zetina_{k_4}^{\sigma_4}\,,\\
\notag
 H_{k_1, k_2, k_3, k_4}^{\sigma_1, \sigma_2, \sigma_3, \sigma_4} & := 
 \frac12 \sigma_3 \sigma_4 \left( |k_3| k_4^2 - |\sigma_2 k_2 + \sigma_3 k_3| \, |k_3| \, |k_4| \right) \fm_{k_1} \fm_{k_2} \fn_{k_3} \fn_{k_4} \\
 \label{H4.coeff}
 & + \frac{\gamma}{4} \sigma_3 \sigma_4 \left( \sign(k_3) k_4^2 + k_3 |k_4| - 2 |\sigma_2 k_2 + \sigma_3 k_3| \, \sign(k_3) |k_4|\right) \fm_{k_1} \fm_{k_2} \fm_{k_3} \fn_{k_4}\\
 \notag
 & + \frac{\gamma^2}{8} \sigma_1 \sigma_4 \left( k_1 \sign(k_4) - \sign(k_1) \sign(k_4) \, |\sigma_3 k_3 + \sigma_4 k_4| \right) \, \fm_{k_1} \fm_{k_2} \fm_{k_3} \fm_{k_4} \,.
\end{align}
\end{lemma}
\begin{remark}\label{rem.Hreal}
Hamiltonians of the form \eqref{H3.form}, \eqref{H4.form} are real-valued provided the coefficients satisfy
\be\label{real.cond}
\bar{H_{\vec k}^{\vec \sigma} } = H_{\vec k}^{-\vec \sigma}  \,. 
\ee
Actually, the coefficients \eqref{H3.coeff} and \eqref{H4.coeff} are real.
\end{remark}
 
\begin{proof}[{\bf Proof of Lemma \ref{lem:expham}.}] We divide the proof in two steps, writing first the Taylor expansion of the  Hamiltonian in the Wahl\'en coordinates $(\eta,\wahlen)$ and then passing to the complex Zakharov variables $(\zetina, \bar \zetina)$. 

\smallskip
\noindent{\sc Step 1: from $(\eta, \psi)$ to the  Wahl\'en variables $(\eta,\wahlen)$}.
The Hamiltonian $\scH_\gamma(\eta, \psi)$ in the  Wahl\'en variables \eqref{Whalen} is the 
 real Hamiltonian $H(\eta, \wahlen) := \scH_\gamma(\eta, \wahlen + \frac{\gamma}{2} \pa_x^{-1} \eta ) $ explicitly given in view of \eqref{H.gamma} by
\begin{equation}\label{Hwahlen}
\begin{aligned}
H(\eta, \wahlen) = \int_\T
\left( \frac12 \wahlen \cdot G(\eta)\wahlen + \frac12  \eta^2  
 +\frac{\gamma}{2}  \wahlen  \cdot G(\eta)[\pa_x^{-1}\eta] - \frac{\gamma^2}{8}  \eta \cdot  \pa_x^{-1}G(\eta)\pa_x^{-1}\eta  - \frac{\gamma}{2} \wahlen_x \eta^2  - \frac{\gamma^2}{12}  \eta^3 \right) \di x  \,. 
\end{aligned}
\end{equation}
We Taylor expand $H(\eta, \wahlen)$ in homogeneity up to quartic terms. 
We recall the expansion of the Dirichlet-Neumann  operator
\begin{equation}\label{DNexp}
G(\eta)\psi = |D| \psi + G_1(\eta) \psi + G_2(\eta)\psi + G_{\geq 3}(\eta)\psi \ , \quad D:= \frac{1}{\im}\pa_x \ ,
\end{equation}
where 
\begin{equation}\label{DNexp2}
G_1(\eta) := - \pa_x \eta \pa_x - |D| \eta |D|  \ , \quad 
G_2(\eta):= - \frac12 \left( D^2 \eta^2 |D| + |D| \eta^2 D^2 - 2 |D|\, \eta \,|D|\,\eta \,|D| \right)
\end{equation}
and where $G_{\geq 3}(\eta) \in \cM_{\geq 3}[r]$, see e.g. formula (2.5) of \cite{CS2} and Remark 3.2 of \cite{BFP}.

We  Taylor expand the Hamiltonian \eqref{Hwahlen} and, using the expansion \eqref{DNexp}, \eqref{DNexp2}, the  identities 
\begin{equation}
|D| \pa_x^{-1} = \Hilb  \,, 
\quad 
\pa_x^{-1} \Hilb = - |D|^{-1} \,, 
\qquad \Hilb \mbox{ in } \eqref{def:Hilbert} \,,
\label{identities.fou}
\end{equation}
 and recalling that $\eta(x)$ has zero average (cf.  \eqref{zero.av}), 
 we get 
\be
H = H_{2} + H_{3}  + H_{4}  + H_{\geq 5} \,,
\ee
where
\be\label{Hwahlen.taylor}
\begin{aligned}
H_{2}(\eta, \wahlen)&:=  \int_\T \left(\frac12 \wahlen (|D| \wahlen)  +\frac12 \eta^2 + \frac{\gamma}{2}  \wahlen (\Hilb \eta) +  \frac{\gamma^2}{8} \eta (|D|^{-1} \eta) \right) \, \di x \,, \\
H_{3}(\eta, \wahlen)&:= \int_\T \left( \frac12 \eta \wahlen_x^2 - \frac12 \eta (|D|\wahlen)^2 
 - \frac{\gamma}{2}  \eta (\Hilb \eta) (|D| \wahlen) + \frac{\gamma^2}{24} \eta^3  - \frac{\gamma^2}{8} \eta (\Hilb \eta)^2 
\right) \di x \,, \\
H_{4}(\eta, \wahlen)&:= \int_\T \left( - \frac12 \eta^2 (|D| \wahlen) (D^2\wahlen)  + \frac12 \eta (|D|\circ\eta \circ |D|\wahlen)(|D|\wahlen)  - \frac{\gamma}{4}\eta^2 (\Hilb\eta)(D^2\wahlen) + \frac{\gamma}{4} \eta^2 \eta_x (|D|\wahlen) \right) \di x\\
& + \int_\T \left(  \frac{\gamma}{2}\eta  (|D|\circ\eta\circ\Hilb \eta) (|D|\wahlen)+ \frac{\gamma^2}{8} \eta_x\eta^2  (\Hilb\eta) + \frac{\gamma^2}{8} (\Hilb\eta) \eta(|D|\circ \eta\circ \Hilb \eta) \right) \di x \,
\end{aligned}
\ee
and $H_{\geq 5}$ collects the at least quintic terms. 
Expanding $(\eta, \wahlen)$ in the corresponding 
Fourier variables (see \eqref{fourierseries}), we compute  $H_{3}(\eta, \wahlen)$ and $H_{4}(\eta, \wahlen)$ obtaining that 
\begin{equation}\label{H3.provvisoria}
    H_{3}(\eta, \wahlen)  = 2\pi \!\!\!\!\!\!\!\!\!\!\!\!
    \sum_{\substack{k_1, k_2, k_3 \in \Z_* \\ k_1 + k_2 + k_3 = 0}} \!\!\!\!\!\!\!\!\!\!
    \tA_{k_1, k_2, k_3} \eta_{k_1} \wahlen_{k_2} \wahlen_{k_3}
+
 \tB_{k_1, k_2, k_3} \eta_{k_1} \eta_{k_2} \wahlen_{k_3}
 +
  \tC_{k_1, k_2, k_3} \eta_{k_1} \eta_{k_2} \eta_{k_3}
\end{equation}
with
{\small \begin{equation}\label{Hwahlen3taylor.coeff}
        \tA_{k_1, k_2, k_3} := - \frac12 \left( k_2 k_3 + |k_2| \, |k_3| \right) \ , \quad 
   \tB_{k_1, k_2, k_3}:= \im \frac{\gamma}{2} \sign(k_2) |k_3|  \ , \quad
   \tC_{k_1, k_2, k_3} := \frac{\gamma^2}{8} \left( \frac13 + \sign(k_2)\, \sign(k_3) \right)
\end{equation}}
and
\begin{equation}\label{Hwahlen4taylor}
    H_{4}(\eta, \wahlen) = 2\pi \!\!\!\!\!\!\!\!\!\! \!\!\!
    \sum_{\substack{k_1, k_2, k_3, k_4 \in \Z_* \\ k_1 + k_2 + k_3+ k_4 = 0}} \!\!\!\!\!\!\!\!\!\!\!\!\!
    \tA_{k_1, k_2, k_3,k_4} \eta_{k_1} \eta_{k_2} \wahlen_{k_3} \wahlen_{k_4}
+
 \tB_{k_1, k_2, k_3,k_4} \eta_{k_1} \eta_{k_2} \eta_{k_3} \wahlen_{k_4}
 +
  \tC_{k_1, k_2, k_3, k_4} \eta_{k_1} \eta_{k_2} \eta_{k_3} \eta_{k_4}
\end{equation}
with
{\small
\begin{align}
\notag
   \tA_{k_1, k_2, k_3, k_4} &:=  \frac12 \left( |k_2 + k_3|\, |k_3| \, |k_4| -  |k_3| k_4^2   \right)  \,, \quad 
   \tB_{k_1, k_2, k_3, k_4} := \im \frac{\gamma}{4} \left( \sign(k_3) k_4^2 + k_3 |k_4| - 2 \sign(k_3)\,  |k_2 + k_3| \, |k_4|  \right) \,,\\
     \label{Hwahlen4taylor.coeff}
   \tC_{k_1, k_2, k_3, k_4} &:= \frac{\gamma^2}{8} \left(  \sign(k_4) k_1- \sign(k_1) \, \sign(k_4) \, |k_3 + k_4| \right)\,.
\end{align}
}
Note that, in view of the cancellation 
    $
    \tC_{k_1, k_2, k_3}+ \tC_{k_2, k_3, k_1} +\tC_{k_3, k_1, k_2}=0$ (see Lemma $3.1$ in \cite{Sulem}), 
we  have 
     $$
    \sum_{\substack{k_1, k_2, k_3 \in \Z_* \\ k_1 + k_2 + k_3 = 0} }
  \tC_{k_1, k_2, k_3} \eta_{k_1} \eta_{k_2} \eta_{k_3}=0\,,
    $$
    so that the Hamiltonian $H_3$ in \eqref{H3.provvisoria} reads
    \begin{equation}
    \label{Hwahlen3taylor}
     H_{3}(\eta, \wahlen)  = 2\pi \!\!\!\!\!\!\!\!\!
     \sum_{\substack{k_1, k_2, k_3 \in \Z_* \\ k_1 + k_2 + k_3 = 0}}
     \!\!\!\!\!\!\!
     \tA_{k_1, k_2, k_3} \eta_{k_1} \wahlen_{k_2} \wahlen_{k_3}
+
 \tB_{k_1, k_2, k_3} \eta_{k_1} \eta_{k_2} \wahlen_{k_3}\,.
    \end{equation}
    The expression in \eqref{Hwahlen3taylor} agrees with formula $(3.2)$ in \cite{Sulem}. 
    Moreover, to prove that \eqref{Hwahlen4taylor} agrees with formula $(4.11)$ in \cite{Sulem}, it is sufficient to symmetrize the coefficients and to use the momentum condition $ k_2 +k_3= -( k_1+k_4)$ for $\tB_{k_1,k_2,k_3,k_4}$.
\medskip

\noindent{\sc Step 2: from Wahl\'en variables $(\eta,\wahlen)$ to Zakharov  variables $(\zetina,\bar\zetina)$.} 
We now write $H_{n}(\eta, \wahlen)$, $n=2,3,4$ in \eqref{Hwahlen.taylor} in the  complex variables $({\zetina}, {\bar \zetina})$ defined in \eqref{zak}, namely we compute
\begin{equation}\label{HC}
\cH(\zetina, \bar \zetina):= H(\cM (\zetina, \bar \zetina)) = \cH_2(\zetina, \bar \zetina) +  \cH_3(\zetina, \bar \zetina)  +\cH_4(\zetina, \bar \zetina) +  \cH_{\geq 5}(\zetina, \bar \zetina) \,.
\end{equation}
Since, by  \eqref{zak},  
$ 
\eta = \frac{1}{\sqrt 2}  M^{-1}(D) (\zetina + \bar \zetina)$ and  $  \wahlen = - \frac{\im}{\sqrt 2} M(D) (\zetina - \bar \zetina)$ with $  M(D) $ in \eqref{defMD} (use that $\vect{\eta}{\zeta}=\cM \zak$), 
 the $j$-th Fourier coefficients $\eta_j, \wahlen_j$ are given, in terms of the Fourier coefficients $\zetina_j^+$, $ \zetina_j^-$ (recall \eqref{u+-}),  by 
\begin{equation}\label{etajzetaju}
\eta_j = \sum_{\sigma \in \{\pm\}}\fm_j  \zetina_{\sigma j}^\sigma \,, 
\quad
\wahlen_j = \sum_{\sigma \in \{\pm\}}  - \im \sigma \fn_j  \zetina_{\sigma j}^\sigma \,, 
\quad 
\fm_j, \fn_j \mbox{ in } \eqref{al.be.m} \,. 
\end{equation}
The quadratic Hamiltonian $ \cH_2(\zetina, \bar \zetina)$ is readily computed to be 
\eqref{H2}
recalling that 
the  complex Zakharov variables $\zak$ diagonalize the quadratic Hamiltonian $\cH_{\gamma}^{(2)}$, see \eqref{diaglin} and use also \eqref{hvf}.

Using \eqref{Hwahlen3taylor}, \eqref{Hwahlen4taylor}, \eqref{etajzetaju} and  the fact that $\fm_j = \fm_{-j}$ and $\fn_j = \fn_{-j}$ for any $j \in \Z_*$ we  also compute
\begin{equation*}
\cH_3(\zetina, \bar \zetina) :=2\pi  \sum_{ \vec \sigma \cdot \vec k = 0 } 
H_{k_1, k_2, k_3}^{\sigma_1, \sigma_2, \sigma_3} \, \zetina_{k_1}^{\sigma_1} \zetina_{k_2}^{\sigma_2}\zetina_{k_3}^{\sigma_3} \,,
\qquad 
\cH_4(\zetina, \bar \zetina)  := 2\pi \sum_{ \vec \sigma \cdot \vec k = 0 } 
H_{k_1, k_2, k_3, k_4}^{\sigma_1, \sigma_2, \sigma_3, \sigma_4} \, \zetina_{k_1}^{\sigma_1} \zetina_{k_2}^{\sigma_2}
\zetina_{k_3}^{\sigma_3} 
\zetina_{k_4}^{\sigma_4} \,,
\end{equation*}
where
\begin{align*}
 H_{k_1, k_2, k_3}^{\sigma_1, \sigma_2, \sigma_3} &:= 
- \sigma_2 \sigma_3 \fm_{k_1} \fn_{k_2} \fn_{k_3}  \tA_{\sigma_1 k_1, \sigma_2 k_2, \sigma_3 k_3}
- \im 
 \sigma_3 \fm_{k_1} \fm_{k_2} \fn_{k_3}  \tB_{\sigma_1 k_1, \sigma_2 k_2, \sigma_3 k_3}\ ,\\
  H_{k_1, k_2, k_3,k_4}^{\sigma_1, \sigma_2, \sigma_3, \sigma_4} &:= 
-  \sigma_3 \sigma_4 \fm_{k_1} \fm_{k_2} \fn_{k_3}\fn_{k_4}  \tA_{\sigma_1 k_1, \sigma_2 k_2, \sigma_3 k_3, \sigma_4 k_4} 
- \im 
 \sigma_4 \fm_{k_1} \fm_{k_2} \fm_{k_3} \fn_{k_4}  \tB_{\sigma_1 k_1, \sigma_2 k_2, \sigma_3 k_3, \sigma_4 k_4}+ \\
 &   \qquad \qquad 
+ \fm_{k_1} \fm_{k_2} \fm_{k_3} \fm_{k_4} \tC_{\sigma_1 k_1, \sigma_2 k_2, \sigma_3 k_3, \sigma_4 k_4}\,.
\end{align*}
The expressions \eqref{H3.coeff}, \eqref{H4.coeff} follow using formulae \eqref{Hwahlen3taylor.coeff}, \eqref{Hwahlen4taylor.coeff}.
\end{proof}
Next we transform the Hamiltonian $\cH$ in \eqref{HC} into its {\em formal} Birkhoff normal form.
\begin{proof}[{\bf Proof of \Cref{prop:cH^{(4)}_+}}]
The proof is divided in two steps: in the first one we remove the cubic part of the  Hamiltonian $\cH$ in \eqref{HC}, whereas in the second one we compute the restrictions of the new quartic Hamiltonian to the sets $\fR_\Set^{(n)}$, $n =0,1,2$.

\smallskip
\textsc{Step 1: Elimination of the cubic Hamiltonian:} 
We proceed to remove all cubic monomials, using the absence of $3$-waves resonances proved in \Cref{3onde}. 
In order to do that, we look for a symplectic transformation $\Phi = \Phi_{\cF_3}$, defined as the time-1 (formal) flow 
 $\Phi_{\cF_3} = \Phi^t_{\cF_3}\vert_{t = 1}$, generated by
 a cubic real  Hamiltonian $\cF_3$ of the form \eqref{H3.form}. 
A Lie expansion then yields
\be\label{lie3}
\cH\circ\Phi_{\cF_3} = \cH_2  + \{ \cF_3, \cH_2 \} + \cH_3 + 
\cH_4 + \frac12 \{\cF_3, \{\cF_3, \cH_2 \} \} + \{ \cF_3, \cH_3 \} +  \cH^+_{\geq 5}\,,
\ee
where $\cH^+_{\geq 5}$ are new quintic contributions. 
The Hamiltonian $\cF_3$ is then determined as the {\em unique} solution of the cohomological equation
\be\label{hom.eq3}
\{ \cF_3, \cH_2 \} + \cH_3  = 0 \,.
\ee
\begin{lemma}\label{lem:Fede}
   The unique solution of \eqref{hom.eq3} is given by the real Hamiltonian
   \begin{equation}\label{def:F3}
    \cF_3(z, \bar z):= 2\pi 
\sum_{\vec \sigma \cdot \vec k = 0} F_{k_1, k_2, k_3}^{\sigma_1, \sigma_2, \sigma_3 }
z_{k_1}^{\sigma_1} z_{k_2}^{\sigma_2} z_{k_3}^{\sigma_3}\,, 
\qquad 
F_{k_1, k_2, k_3}^{\sigma_1, \sigma_2, \sigma_3 }:=
 \frac{- H_{k_1, k_2, k_3}^{\sigma_1, \sigma_2, \sigma_3}}{\im ( \sigma_1 \Omega_{k_1}(\gamma) + \sigma_2 \Omega_{k_2}(\gamma) + \sigma_3 \Omega_{k_3}(\gamma) )} \,,
\end{equation}
where the coefficients $H_{k_1, k_2, k_3}^{\sigma_1, \sigma_2, \sigma_3}$ are defined in 
\eqref{H3.coeff}, and one has the real property
$\bar{F_{\vec k}^{\vec \sigma}} = F_{\vec k}^{-\vec \sigma}$. 
\end{lemma}
\begin{proof}
   Recalling the definition of $\cH_2$ in \eqref{H2}, we have 
   $\partial_{z_k^\sigma} \cH_2= 2\pi \Omega_k(\gamma)z_{k}^{-\sigma}$ whereas 
    \begin{align*}
        \partial_{z_k^{-\sigma}} \cF_3= 2\pi \sum_{\sigma_1k_1+\sigma_2k_2-\sigma k=0} \left(F_{k,k_1,k_2}^{-\sigma,\sigma_1,\sigma_2}+F_{k_1,k,k_2}^{\sigma_1,-\sigma,\sigma_2}+F_{k_1,k_2,k}^{\sigma_1,\sigma_2,-\sigma}\right)z_{k_1}^{\sigma_1}z_{k_2}^{\sigma_2}.
    \end{align*}
    Then, by \eqref{PP},
    $
        \{ \cF_3, \cH_2 \} 
        = \frac{2\pi}{\ii} \!\!\!\!\!\!\!\!\!
        \sum\limits_{\sigma_1 k_1+\sigma_2 k_2 +\sigma k=0}\!\!\!
\left(F_{k,k_1,k_2}^{\sigma,\sigma_1,\sigma_2}+F_{k_1,k,k_2}^{\sigma_1,\sigma,\sigma_2}+F_{k_1,k_2,k}^{\sigma_1,\sigma_2,\sigma}\right)(-\sigma \Omega_k(\gamma))z_{k_1}^{\sigma_1}z_{k_2}^{\sigma_2}z_{k}^{\sigma}$ 
    which  yields 
    \begin{align}\label{poiss.H2}
    \{ \cF_3, \cH_2 \}= 2\pi  \sum_{\sigma_1k_1+\sigma_2k_2+\sigma_3 k_3=0} F_{k_1,k_2,k_3}^{\sigma_1,\sigma_2, \sigma_3} \ii \left( \sigma_1\Omega_{k_1}(\gamma)+\sigma_2\Omega_{k_2}(\gamma)+\sigma_3\Omega_{k_3}(\gamma)\right)z_{k_1}^{\sigma_1}z_{k_2}^{\sigma_2}z_{k_3}^{\sigma_3}
    \end{align}
    relabeling the indexes. Defining $F_{k_1,k_2,k_3}^{\sigma_1,\sigma_2, \sigma_3}$
    as in \eqref{def:F3} solves \eqref{hom.eq3}.
   Such solution is unique since, by Lemma \ref{3onde}, the system
$\sigma_1 k_1 + \sigma_2 k_2 + \sigma_3 k_3 = 0$, $\sigma_1 \Omega_{k_1}(\gamma) + \sigma_2 \Omega_{k_2} (\gamma)+ \sigma_3 \Omega_{k_3} (\gamma) = 0
$ has no integer solutions whatever $\gamma$. 
\end{proof}
With this choice of $\cF_3$, the Hamiltonian \eqref{lie3} reduces to
\be\label{lie3.bis}
\cH\circ\Phi_{\cF_3} = \cH_2  + 
\cH_4^+ + \cH^+_{\geq 5}  \ , \qquad \cH_4^+:= \cH_4 +  \frac12 \{ \cF_3, \cH_3 \} \,,
\ee
with $\cH_3 $, $\cH_4 $ in   \eqref{H3.form}, \eqref{H4.form} and $\cF_3$ in \eqref{def:F3}.

\smallskip
\textsc{Step 2: strong-$\Set$ normal form:}  We claim that the Hamiltonian $\cH_4^+$ in \eqref{lie3.bis} fulfills 
\begin{align}
\label{HR0}
\Pi_{\fR_\Set^{(0)}}\cH_4^+(Z) = 2\pi \left(
\fa |z_{\tm}|^4 + \fb |z_\tm|^2 \, |z_\tn|^2 + \fc |z_\tn|^4  \right) \,, 
\quad \fa, \fb, \fc \mbox{ in } \eqref{abc} \,, \\
\label{HR12}
\Pi_{\fR_\Set^{(1)}} \cH_4^+ =  0 \,, \qquad 
\left\lbrace 
\Pi_{\fR_\Set^{(2)}} \cH_4^+ , |z_j|^2 \right\rbrace = 0 \quad  \forall \, j \in \Z\setminus \{0 \} \,. 
\end{align}
Then the claimed properties  \eqref{H4+.str0},  \eqref{H4+.str} 
of $X_{\cH_4^+}$ follow by computing the Hamiltonian vector field  using \eqref{hvf} and -- for the second of \eqref{H4+.str}-- exploiting  that, with $f_j(Z) = |z_j|^2$,
$$
0 = \left\lbrace \Pi_{\fR_\Set^{(2)}} \cH_4^+ , f_j \right\rbrace = \di f_j(Z)[X_{\Pi_{\fR_\Set^{(2)}} \cH_4^+}] = 
\bar z_j \, \Pi_{\fR_\Set^{(2)}}
X_{\cH_4^+}  + z_j \bar{\Pi_{\fR_\Set^{(2)}}
X_{\cH_4^+} } = 2   \Re\left( [\Pi_{\fP_\Set^{(2)}}X_{\cH_4^+}(Z)]^+_j    \, \bar z_j \right) \,.
$$
{\bf Projection on $\fR_\Set^{(0)}$.} 
In view of item $(i)$ of Lemma \ref{lem:wres}, the set $\fR_\Set^{(0)} = \fP_\Set^{(0)}$ contains only integrable resonances. Thus, the restricted Hamiltonian has the form \eqref{HR0}
for some coefficients $\fa, \fb, \fc$ that we shall now compute.\\

\noindent\underline{ Computation of $\fa$.}
First of all, consider the coefficients coming from the Hamiltonian $\cH_4$. 
 In view of the expression \eqref{H4.form} (which is not symmetrized),  the monomial $|z_\tm|^4$  has coefficient given by 
\be\label{tHmmmm}
2\pi\, \tH_{\tm}:= 2\pi \left(
H_{\tm , \tm , \tm , \tm }^{+, +, -, -} +  H_{\tm , \tm , \tm , \tm}^{+, -, +, -} + H_{\tm , \tm , \tm , \tm}^{+, -, -, +}
+ 	  H_{\tm , \tm , \tm , \tm}^{-, -,  +, +} +  H_{\tm , \tm , \tm , \tm}^{-, +, -, +} + H_{\tm , \tm , \tm , \tm}^{-, +, +, -} \right)\,.
\ee
Next we compute the coefficient of the monomial $|z_\tm|^4$ coming from  $ \frac12 \{ \cF_3, \cH_3 \}$. 
Such monomial can only come from the following Poisson brackets
$$
\{ z_\tm^2 \bar{z_{2\tm}} , \ \bar{z_\tm}^2 {z_{2\tm}} \} \,, 
\quad 
\{ z_\tm^2 {z_{-2\tm}} , \ \bar{z_\tm}^2 \bar{z_{-2\tm}} \}  \,, 
$$
where we kept track also of the momentum conservation.
Now note that  the coefficients of the monomials 
 $z_\tm^2 \bar{z_{2\tm}}$ and 
 $z_\tm^2 {z_{-2\tm}}$
 coming from a generic real cubic Hamiltonian $\cG_3(Z)= 2\pi \sum_{\vec{\s} \cdot \vec{k} = 0} G_{k_1  k_2 k_3}^{\s_1 \s_2 \s_3} z_{k_1}^{\s_1} z_{k_2}^{\s_2} z_{k_3}^{\s_3}$  are respectively 
\be\label{tHm2m}
2\pi\, \tG_{\tm}^{-}:=  2\pi \left(G_{\tm, \tm, 2\tm}^{+, +, -} + G_{\tm, 2\tm, \tm}^{+, -, +} + G_{2\tm, \tm, \tm}^{-, +, +} \right) \ , 
\qquad
2\pi\, \tG_{\tm}^{+}:= 2\pi\left( G_{\tm, \tm, -2\tm}^{+, +, +} + G_{\tm, -2\tm, \tm}^{+, +, +} + G_{-2\tm, \tm, \tm}^{+, +, +}\right)\,.
\ee 
Then, using also the  reality properties \eqref{real.cond} and \Cref{lem:Fede} and the fact that the coefficients $H_{\vec{k}}^{\vec \sigma}$ are real, we get 
\begin{align}\notag
 \frac12 \{ \cF_3, \cH_3 \} & = 
\frac{(2\pi)^2}{2}    \{ \tF_\tm^- \, z_\tm^2 \bar{z_{2\tm}}
+  \tF_{\tm}^{+} \, z_{\tm}^2 z_{-2\tm} ,  \  \tH_{\tm}^{-} \, \bar{z_\tm }^2 z_{2\tm} + 
\tH_{\tm}^{+}\, 
 \bar{z_\tm}^2 \bar{z_{-2\tm}}  \}+ c.c. + \mbox{other monomials} \\
 \label{2611:1709}
& =  2\pi  \textup{Re}\left( \frac{1}{\im}  \tF_{\tm}^{-}  \,  \tH_{\tm}^{-}  - \frac{1}{\im}  \tF_{\tm}^{+} \, \tH_{\tm}^{+} \right) |z_\tm|^4 
+ \mbox{ other monomials} \,,
\end{align}
where with the expression "other monomials" we mean a sum of terms which do not contribute to $|z_\tm|^4$.
In conclusion, in view of \eqref{tHmmmm} and \eqref{2611:1709}, the Hamiltonian $\cH_4^+$ has the monomial $2\pi \fa |z_\tm|^4$ with 
\be\label{a.coeff}
\fa = \tH_{\tm}
+  \textup{Im}\left(   \tF_{\tm}^{-}  \, \tH_{\tm}^{-}  -  \tF_{\tm}^{+} \,\tH_{\tm}^{+} \right)  \,.
\ee
\underline{ Computation of $\fc$.}  This is analogous to the computation of $\fa$. 
The Hamiltonian $\cH_4^+$ has the monomial $2\pi \fc |z_\tn|^4$ with  (using the notations \eqref{tHmmmm} and \eqref{tHm2m} with $\tm \leadsto \tn$)
\be\label{c.coeff}
\fc= \tH_{\tn}
+  \textup{Im}\left(   \tF_{\tn}^{-}  \, \tH_{\tn}^{-}  -   \tF_{\tn}^{+} \,\tH_{\tn}^{+} \right)  \,.
\ee
\underline{ Computation of $\fb$.} 
First of all, consider the coefficient of the  monomial $|z_\tm|^2 \, |z_\tn|^2$ coming from $\cH_4$; in view of the expression \eqref{H4.form}, it is given by 
\be\label{tHmmnn}
2\pi\, \tH_{\tm, \tn}:= 
2\pi \sum_{(\vec j, \vec \sigma) \in \cV_{\tm, \tn}}
H_{j_1 , j_2 , j_3 , j_4 }^{\sigma_1, \sigma_2, \sigma_3, \sigma_4}  ,
\qquad
\cV_{\tm \tn} := \left\{ \left(\pi(\tm, \tm, \tn, \tn), \ \pi(+,-, +, -) \right)\colon \pi \in \cS_4 \right\}\,,
\ee
with $\cS_4$ the symmetric group of permutations of four symbols.

Next we compute the coefficient of the monomial $|z_\tm|^2 |z_\tn|^2$  coming from the Poisson bracket  $ \frac12 \{ \cF_3, \cH_3 \}$. 
Such monomial can only come from the following Poisson brackets:
\begin{align}
    (i) &\ \{|z_\tn|^2 z_j^{\sigma}\,,\ z_j^{-\sigma}|z_\tm|^2 \} , \quad 
    (ii) \  \{z_\tn z_\tm z_j^{\s}\,,\ z_j^{-\s} \bar{z_\tn} \bar{z_\tm}\}\,, \quad 
    (iii) \ \{z_\tn \bar{z_\tm} z_j^{\s}\,,\ z_j^{-\s} \bar{z_\tn} z_\tm\}\,, \qquad j \neq \tm,\tn\,, \label{combinatorial:poisson}\\
    (iv)& \ \{z_\tn z_\tm z_\tn^{\s}\,,\ z_\tn^{-\s} \bar{z_\tn} \bar{z_\tm}\}\,, \quad
    (v) \ \{z_\tn z_\tm z_\tm^{\s}\,,\ z_\tm^{-\s} \bar{z_\tn} \bar{z_\tm}\}\,, \quad
    (vi)
    \ \{z_\tn \bar{z_\tm} z_\tn^{\s}\,,\ z_\tn^{-\s} \bar{z_\tn} z_\tm\}\,, \quad 
    (vii) \ \{z_\tn \bar{z_\tm} z_\tm^{\s}\,,\ z_\tm^{-\s} \bar{z_\tn} z_\tm\}\,.\notag
\end{align}
Now the monomial in case $(i)$ cannot appear due to momentum conservation and $j \neq 0$. Again in view of momentum conservation, the cases in $(ii)$ and $(iii)$ give rise to the Poisson brackets
\begin{equation}\label{poiss.poss}
\begin{gathered}
\{ z_{\tm} \bar{z_{\tn}}  {z_{\tn -\tm }}, \, \bar{z_{\tm}} {z_{\tn}} \bar{z_{\tn -\tm}} \} , \quad 
\{ z_{\tm} z_{\tn} {z_{-\tm - \tn}}, \, \bar{z_{\tm}} \bar{z_{\tn}} \bar{z_{-\tm - \tn}} \}  \ , 
\quad
\{ \bar{z_{\tm}} {z_{\tn}} z_{\tm - \tn}, \, z_{\tm} \bar{z_{\tn}} \bar{z_{\tm - \tn}} \} \ , \quad 
\{ \bar{z_{\tm}} \bar{z_{\tn}} z_{\tm + \tn}, \, z_{\tm} z_{\tn} \bar{z_{\tm + \tn}} \} \,,\\
\text{with} \quad \tn - \tm\,,\  -\tm -\tn\,, \ \tm - \tn\,, \ \tm + \tn \notin \Lambda\,.
\end{gathered}
\end{equation}
Furthermore, recalling that  the set $\Lambda = \{\tm,\ \tn\}$ is $\gamma$-good (see \Cref{g-good}) the cases $(iv), (vi), (vii)$ are not possible by momentum conservation, since $\tm <  0 <\tn$, $\tm + \tn >0$, and  by momentum conservation case $(v)$ requires either $\tn - 2 \tm = 0$, absurd due to $\tm < 0 <\tn$, or $\tn + 2 \tm = 0$, absurd since $\tn \neq -2\tm$ in view of \eqref{G4}.

Consider a generic,   real cubic Hamiltonian 
$\cG_3(Z)= 2\pi \sum_{\vec{\s} \cdot \vec{k} = 0} G_{k_1  k_2 k_3}^{\s_1 \s_2 \s_3} z_{k_1}^{\s_1} z_{k_2}^{\s_2} z_{k_3}^{\s_3}$; the  
coefficient of its monomial  $z_\tm^{\sigma_1} {z_\tn^{\sigma_2}} {z_{{j}}^{\sigma_3}} $, for $j \not\in \Set$, 
is
\begin{align}\label{tHmnj}
   2\pi  \tG_{\tm, \tn, {j}}^{\sigma_1, \sigma_2 , \sigma_3}:= 2\pi \left(
    G_{\tm, \tn, {j}}^{\sigma_1, \sigma_2 , \sigma_3}+
    G_{\tm, {j}, \tn}^{\sigma_1, \sigma_3 , \sigma_2}+
    G_{\tn, \tm, {j}}^{\sigma_2, \sigma_1 , \sigma_3}+
    G_{\tn, {j}, \tm}^{\sigma_2, \sigma_3 , \sigma_1}+
    G_{{j}, \tm, \tn}^{\sigma_3, \sigma_1 , \sigma_2}+
    G_{{j}, \tn, \tm}^{\sigma_3, \sigma_2 , \sigma_1}\right) \,. 
\end{align}
Then, using the  reality properties in \eqref{real.cond} and \Cref{lem:Fede},  the notation \eqref{tHm2m}, and the second of \eqref{poiss.poss},
we get 
\begin{align}\notag
 \frac12 \{ \cF_3, \cH_3 \} & =
\frac{(2\pi)^2}{2} 
\left\lbrace 
 \tF_{\tm, \tn, \tn-\tm}^{+,-,+} \, z_\tm \bar{z_\tn} {z_{\tn- \tm}} 
+\tF_{\tm, \tn, -\tm-\tn}^{+, +, +} \, z_{\tm} z_{\tn} z_{-\tm - \tn} 
+\tF_{\tm, \tn, \tm-\tn}^{-,+,+} \, \bar{z_\tm} {z_\tn} {z_{\tm- \tn}} 
+\tF_{\tm, \tn, \tm+\tn}^{-, -, +} \,\bar{z_{\tm}} \bar{ z_{\tn}} z_{\tm +\tn} \right.
\,, \\ \notag
& \qquad  \ \left.
 \tH_{\tm, \tn, \tn-\tm}^{-,+,-} \, \bar{z_\tm} {z_\tn} \bar{{z_{\tn- \tm}}} 
+\tH_{\tm, \tn, -\tm-\tn}^{-, -, -} \, \bar{z_{\tm} z_{\tn} z_{-\tm - \tn}} +\tH_{\tm, \tn, \tm-\tn}^{+,-,-} \, {z_\tm} \bar{z_\tn} \bar{{z_{\tm- \tn}}} 
+\tH_{\tm, \tn, \tm+\tn}^{+, +, -} \, z_{\tm} z_{\tn} \bar{z_{\tm + \tn}} 
\right\rbrace \\
\notag
& \quad + c.c.  + \mbox{ other monomials }
  \\ 
  \notag
& = - 2\pi \textup{Re}\Big( \frac{1}{\ii} \big( 
\tF_{\tm, \tn, \tn-\tm}^{+,-,+} \, \tH_{\tm, \tn, \tn-\tm}^{-,+,-}
+
\tF_{\tm, \tn, -\tm-\tn}^{+, +, +}\tH_{\tm, \tn, -\tm-\tn}^{-, -, -}+
\tF_{\tm, \tn, \tm-\tn}^{-,+,+}\tH_{\tm, \tn, \tm-\tn}^{+,-,-}
+\tF_{\tm, \tn, \tm+\tn}^{-, -, +}\tH_{\tm, \tn, \tm+\tn}^{+, +, -} \big)
\Big) |z_\tm|^2 \, |z_\tn|^2 
\\
\label{bbb2}
&\quad + \mbox{ other monomials}\,,
\end{align}
where again with the expression "other monomials" we mean a sum of terms which do not contribute to $ |z_\tn|^2 |z_\tm|^2$.

In conclusion, by \eqref{tHmmnn} and \eqref{bbb2},  the Hamiltonian $\cH_4^+$ has the monomial $2\pi \fb |z_\tm|^2 |z_\tn|^2$   with 
\begin{equation}\label{b.coeff}
    \fb=\tH_{\tm, \tn} - \textup{Im}\Big(  
\tF_{\tm, \tn, \tn-\tm}^{+,-,+} \, \tH_{\tm, \tn, \tn-\tm}^{-,+,-}
+
\tF_{\tm, \tn, -\tm-\tn}^{+, +, +}\tH_{\tm, \tn, -\tm-\tn}^{-, -, -}
+
\tF_{\tm, \tn, \tm-\tn}^{-,+,+}\tH_{\tm, \tn, \tm-\tn}^{+,-,-}
+\tF_{\tm, \tn, \tm+\tn}^{-, -, +}\tH_{\tm, \tn, \tm+\tn}^{+, +, -}
\Big)\,.
\end{equation}
The explicit values of $\fa, \fb, \fc$  in  \eqref{abc} are then computed using formulae
\eqref{H3.coeff}, \eqref{H4.coeff}, \eqref{al.be.m}, \eqref{omegonejin}, \eqref{tHmmmm}, \eqref{tHm2m}, \eqref{tHmmnn}, \eqref{tHmnj};
the computation is also done in the Mathematica notebook \texttt{coeff\_abc.nb} at the page \url{https://git.sissa.it/amaspero/transfer-ww-vorticity}.\\
{\bf Projection on $\fR_\Set^{(1)}$.} By Lemma  \ref{lem:wres}--$(ii)$ the set $\fR_\Set^{(1)}(\gamma) = \emptyset$;
 the first of \eqref{HR12} follows.\\
{\bf Projection on $\fR_\Set^{(2)}$.} By Lemma  \ref{lem:wres}--$(iii)$ the set $\fR_\Set^{(2)}(\gamma)$ contains only integrable monomials. Therefore, one has 
\be\label{HR2}
\Pi_{\fR_\Set^{(2)}} \left(\cH_4 +  \frac12 \{ \cF_3, \cH_3 \}  \right) = 2\pi 
\sum_{j \in \Z_*} \left( \td_j |z_\tm|^2 \, |z_j|^2 + \te_j |z_\tn|^2 |z_j|^2\right)\,,
\ee
for some real numbers $\td_j$, $\te_j$, $j \in \Z_*$;   the second of \eqref{HR12} follows. 
\end{proof}

\subsection{Identification argument} \label{sec:identif}
In this section we provide an \textit{a posteriori} identification argument for the normal form of the water waves \eqref{eq:etapsi} that will be constructed in the next section, in the spirit of \cite{BFP}.
 Even though normal form transformations are generally not unique in the presence of resonances, the absence of three-wave resonances (\Cref{3onde}) ensures that the resonant cubic normal form is uniquely determined, coinciding with those of the formal Hamiltonian $\cH_4^+$ of \Cref{prop:cH^{(4)}_+}. Although the transformation constructed in \Cref{thm:nf} does not coincide with the formal normal form transformation of \Cref{sec:bnf.formal}, the quadratic component of any near-identity normal form transformation is uniquely fixed by the corresponding homological equation (see \Cref{ide:B} below), again thanks to the absence of three-wave resonances.  More precisely, we have the following proposition.
\begin{proposition}[{\bf Identification argument}]\label{prop.ident}
Let $\zak$ be the complex Zakharov variable in \eqref{zak} solving the Hamiltonian system \eqref{iniziale_zak}. Let $\bUpsilon(\zak)$ be a $(\nu, m)$ admissible transformation in $ \cM^\nu_{\geq 0}[r]$ with gain $\vr$, for some $0\leq \nu \leq m$, $\vr\geq m+\nu$  and $ r >0$ (see Definition \ref{admtra}).
Assume that the variable
\begin{equation}
    \label{eq:b.adm}
    Z := \bUpsilon(\zak)\zak
\end{equation}
 solves an equation of the form 
\be\label{eq:zetone}
\pa_t Z = -\ii \vOmega(D) Z + X^+_3(Z)+ M_{\geq 3}(Z)Z\,,
\ee
for some cubic vector fields $ X^+_3$ and some $ M_{\geq 3}(Z)\in \cM^{\nu'}_{\geq 3}[r]$, $\nu'\geq 0$. Then the resonant  part of the cubic vector field $X^+_3$ is completely determined and given by 
   \be 
   \Pi_{\fR_4^{(n)}} X_3^+\equiv  \Pi_{\fR_4^{(n)}}\hamvec{\cH_4^{+}}, \qquad  \cH_4^{+} \mbox{ in } \eqref{new.birkh} \,, 
   \quad n = 0, \ldots, 4 \,. 
   \label{ide:res}
   \ee
\end{proposition}

\noindent{\bf Proof of \Cref{prop.ident}.}
Writing $\scF(\zak) := \bUpsilon(\zak) \zak$, by \eqref{eq:b.adm} we have $Z = \scF(\zak)$. Thanks to \Cref{loc.inv} the map $\scF$ is locally invertible and thus we shall use also the relation  $\zak = \scF^{-1}(Z)$.
Recalling that, by \eqref{zetone:formale}, $\zak$ solves the equation $\pa_t \zak = \hamvec{\cH}(\zak)$, and denoting by $Y$ the vector field at the right hand size of \eqref{eq:zetone}, we have
\begin{equation}\label{eq:pfwd}
Y(Z) = \scF^* \hamvec{\cH} (Z):= \left.\di \scF(\zak)[ \hamvec{\cH}(\zak)]  \right|_{\zak = \scF^{-1}(Z)}\,.
\end{equation}
The expressions that appear throughout the proof depend on the variables $\zak$ and $Z$, and are evaluated at $\zak = \scF^{-1}(Z)$. 
For brevity, we shall often write the variable $\zak$ without explicitly recalling its relation to $Z$, and we simply keep in mind the relation between their norms: using \eqref{stima.inv.adm} and \eqref{lin.est.F}, we have
\begin{align}
    \|\zak\|_s \lesssim_s \|Z\|_{s+\nu}, 
    \qquad 
    \|Z\|_s \lesssim_s \|\zak\|_{s+\nu}\,,
    \label{equiv:zak_Z}
\end{align}
for all $s \ge s_0$ sufficiently large and for all $Z, \zak \in B_{s_0, \R}(r)$, with $r>0$ small enough.
Considering the Taylor expansion
\be \label{B.exp}
    \bUpsilon(\zak) = \Id + \bUpsilon_1(\zak) + \bUpsilon_2(\zak) + \bUpsilon_{\ge 3}(\zak)\,, \qquad \bUpsilon_q(\zak) \in \wt{\cM}_{q}^m\,,\ \ q = 1,2\,, \quad \bUpsilon_{\geq 3}(\zak) \in \cM_{\geq 3}^m[r]\,,
\ee
given by \eqref{esp:F}, we now compute $\scF^* \hamvec{\cH}$. 
First of all, we observe that the difference $\zak - Z$ satisfies the estimate
\begin{equation}    \label{est:diff_zak_Z}
    \|\zak - Z\|_{\sigma}
    = \|\bUpsilon_1(\zak)\zak + \bUpsilon_2(\zak)\zak + \bUpsilon_{\ge 3}(\zak)\zak\|_{\sigma}
    \lesssim \|\zak\|_{\sigma+m} \|\zak\|_{s_0}
    \lesssim \|\zak\|_{\sigma+m}^2
    \lesssim_s \|Z\|_{\sigma+\nu+m}^2,
\end{equation}
obtained by applying \eqref{smoothing} to $\bUpsilon_1(\zak)$ and $\bUpsilon_2(\zak)$, \eqref{piove} to $\bUpsilon_{\ge 3}(\zak)$, and the equivalence \eqref{equiv:zak_Z}.
From the expansion of $\hamvec{\cH}$ in \eqref{iniziale_zak} and the expansion of $\bUpsilon(\zak)$ in \eqref{B.exp}, we have,
\begin{equation}
\label{eq:vec.f.2}
  \scF^* \hamvec{\cH} (Z)  = -\ii \vOmega(D) Z + Y_{\geq 2}(Z) + Y_3(Z) + Y_{\geq 4}(Z)\,,
\end{equation}
where, recalling the definition \eqref{lb} of Lie bracket of two vector fields,
\begin{align}
\label{vf.order2}
Y_{\geq 2}(Z) &=\left.\left(\hamvec{\cH_3}(\zak)+\comm{-\ii \vOmega(D) \zak}{ \bUpsilon_1(\zak)\zak} \right) \right|_{\zak = \scF^{-1}(Z)}\,,\\
Y_3(Z)&=\hamvec{\cH_4}(Z) +\comm{ -\ii \vOmega(D) Z}{\bUpsilon_2(Z)Z} +\di \left(\bUpsilon_1(Z)Z\right)[ \hamvec{\cH_3}(Z)]\,,\label{X3+}\\
Y_{\geq 4}(Z)&= \Big(\ii \vOmega(D) \bUpsilon_{\geq 3}(\zak)\zak+ \hamvec{\cH_{\geq 5}}(\zak)+ \di \left(\bUpsilon_1(\zak)\zak\right)\left[\hamvec{\cH_{\geq 4}}(\zak)\right] \notag\\
&+ \di\left( \bUpsilon_2(\zak)\zak\right)\left[\hamvec{\cH_{\geq 3}}(\zak)\right]
 + \di (\bUpsilon_{\geq 3}(\zak)\zak) [ \cX_{\cH}(\zak)] + Y_3(\zak)- Y_3(Z) \Big)\Big|_{\zak = \scF^{-1}(Z)}\,,\label{N4}
\end{align}
and where we are using the notations $\cH_{\geq 3}:= \cH_{3}+ \cH_{4}+ \cH_{\geq 5}$ and $\cH_{\geq 4}:= \cH_{4}+\cH_{\geq 5}$.
We claim that 
\be\label{Ygeq2=0}
Y_{\geq 2}(Z) \equiv 0 \,.
\ee
Indeed, the vector field in \eqref{eq:vec.f.2} coincides with the one in \eqref{eq:zetone}, which has no quadratic component. Therefore, the quadratic contribution of \eqref{eq:vec.f.2} must vanish identically. After substituting $\zak = Z - \bUpsilon_1(\zak)\zak - \bUpsilon_2(\zak)\zak - \bUpsilon_{\ge 3}(\zak)\zak$, the quadratic component becomes
\be\label{home.q0}
\cX_{\cH_3}(Z) + \comm{-\ii \vOmega(D) Z}{ \bUpsilon_1(Z)Z} \equiv 0 \,.
\ee
Therefore, such quadratic vector field  vanishes on whatever variable is evaluated, so also at $\zak$, proving the claim \eqref{Ygeq2=0}.
We now prove that the homological equation \eqref{home.q0} allows us to identify the quadratic component $\bUpsilon_1(Z)$ of the transformation $\bUpsilon(Z)$. 
\begin{lemma}\label{ide:B}
    One has 
    \be\label{B=-F3}
\bUpsilon_1(Z) Z = - \cX_{\cF_3}(Z)\,,
\ee
where $\cF_3$ is the Hamiltonian function of \Cref{lem:Fede} that solves the homological equation \eqref{hom.eq3}. 
\end{lemma}
\begin{proof}
    In view of \eqref{home.q0} and \eqref{hom.eq3}--\eqref{lb.id} and since $\cX_{\cH_2}(Z) = - \im \vOmega(D)Z$,  both $\bUpsilon_1(Z) Z$ and $-\hamvec{\cF_3}(Z)$ solve the linear equation
\be \label{home.q}
\comm{-\ii \vOmega(D) Z}{ \ \cdot\ } +\hamvec{\cH_3}(Z) = 0
\ee
which, by the absence of $3$-waves resonances  (cf. \Cref{3onde}), has a  unique solution; hence \eqref{B=-F3} follows.
\end{proof}

We now claim that $Y_{\geq 4}$ vanishes to at least fourth order at $Z=0$ (see Lemma \ref{lem:Y4} below). But then, comparing the homogeneous terms of degree $3$ in the two equal vector fields \eqref{eq:zetone}, \eqref{eq:vec.f.2}, we have $Y_{3}(Z) = X_3^{+}(Z)$ meaning that 
\be \label{x3+.y3}
 X^+_3(Z)=\hamvec{\cH_4}(Z) +\comm{ -\ii \vOmega(D) Z}{\bUpsilon_2(Z)Z} -\di \hamvec{\cF_3}(Z)[ \hamvec{\cH_3}(Z)]\,, 
\ee
where we also used the identification \eqref{B=-F3}. 
The following algebraic lemma allows us to identify and further simplify the vector field $\di \hamvec{\cF_3}(Z)[\hamvec{\cH_3}(Z)]$.
\begin{lemma}
Let $\cX$, $\cY$ be smooth vector fields fulfilling
$\comm{- \im \vOmega(D)\ \cdot }{\cX} = \cY  $
then
\be\label{bellina0}
\di \cX [\cY ]  = \frac12 \comm{\cY }{\cX } + \frac12 \comm{ -\ii \vOmega(D)\ \cdot }{\di \cX [ \cX ]}\,.
\ee
\end{lemma}
\begin{proof}
We start by writing  
    \be 
\di \cX [ \cY ]= \frac12 \comm{\cY}{\cX}+ \frac12\left( \di \cX[ \cY ]+\di\cY [\cX ]\right) 
\label{gat1}\,.
    \ee
  We now compute the second addend of the r.h.s. of \eqref{gat1}. Using the assumption $\comm{- \im \vOmega(D)\ \cdot }{\cX} = \cY  $, we get
    \be
 \di \cX[ \cY]+\di\cY [\cX ]= \di \cX\left[ \comm{-\ii \vOmega(D)\ \cdot}{\cX }\right] + \di\comm{-\ii \vOmega(D)\ \cdot}{\cX}[\cX ]\,.
\label{gat2}
    \ee
    We now apply the easily verified vector field identity 
    \be 
\di \comm{A}{B}[H]= \comm{\di A[H]}{B}+\comm{A}{\di B[H]}{+} \di A \left[\comm{B}{H} \right]-\di B \left[\comm{A}{H} \right]+ [ \di A, \di B] H
    \ee
    to the vector fields $A=-\ii \vOmega(D) \ \cdot \ $, $B= \cX $ and $H= \cX $, obtaining 
\begin{align}
    \di\comm{-\ii \vOmega(D)\, \cdot}{\cX}[\cX ]= &\bcancel{\comm{-\ii \vOmega(D)\cX}{\cX}}+\comm{-\ii \vOmega(D)\ \cdot }{\di \cX[\cX]}{-}\xcancel{\ii \vOmega(D)\comm{\cX}{\cX}}\notag\\
    &-\di \cX\left[ \comm{-\ii \vOmega(D)\, \cdot}{\cX }\right]+ \bcancel{\left[-\ii \vOmega(D), \di \cX\right] \cX}\notag\\
    =&\comm{-\ii \vOmega(D)\ \cdot}{\di \cX[\cX]}- \di \cX \left[ \comm{-\ii \vOmega(D) \ \cdot}{\cX }\right]\,,
    \label{gat3}
\end{align}
where in the last equality we used the identity 
\begin{align*}
    \comm{-\ii \vOmega(D)\cX}{\cX}+ \left[-\ii \vOmega(D)\,,\ \di \cX \right]\cX =0\,.
\end{align*}
Gathering \eqref{gat1}, \eqref{gat2} and \eqref{gat3} we eventually get \eqref{bellina0}. 
\end{proof}
We apply  this lemma with $\cX = \cX_{\cF_3}$ and $\cY = \cX_{\cH_3}$,  which satisfy the assumption by \eqref{home.q0}, \eqref{B=-F3}. Then we get  the identity
    \be\label{bellina}
\di \hamvec{\cF_3}(Z)[ \hamvec{\cH_3}(Z)]= \frac12 \comm{\hamvec{\cH_3}(Z)}{\hamvec{\cF_3}(Z)} + \frac12 \comm{ -\ii \vOmega(D)Z}{\di \hamvec{\cF_3}(Z)[ \hamvec{\cF_3}(Z)]}\,, 
    \ee
that, inserted in \eqref{x3+.y3} 
gives
$$
X_3^+(Z) = \hamvec{\cH_4} - \frac 1 2 \comm{\hamvec{\cH_3}(Z)}{\hamvec{\cF_3}(Z)} + \comm{-\ii \vOmega(D)Z}{\bUpsilon_2(Z) Z } - \frac12 \comm{ -\ii \vOmega(D)Z}{\di \hamvec{\cF_3}(Z)[ \hamvec{\cF_3}(Z)]}\,.
$$
Then, note that for any cubic map $\bF_2(Z)Z$ of the form $[\bF_2(Z)Z]^\sigma =  \sum\limits_{ \vec \sigma\cdot \vec{\jmath} = \sigma  k} 
 	\tF_{\vec{\jmath}, k}^{\vec \sigma, \sigma} \, z_{\vec{\jmath}}^{\vec \sigma}\, e^{\im \sigma k x }$ one has 
 		$$
 	\comm{ - \im \vOmega(D) Z }{\bF_2(Z)Z}^\sigma  =
\!\!\!\!	
 	 \sum_{\vec \sigma\cdot \vec{\jmath} = \sigma  k}   \!\!\!
 	 -\im\left( \sigma_1 \Omega_{j_1}(\gamma) +\sigma_2 \Omega_{j_2}(\gamma)+\sigma_3 \Omega_{j_3}(\gamma)  - \sigma \Omega_{k}(\gamma) \right) 
 	\tF_{\vec{\jmath}, k}^{\vec \sigma, \sigma} \, z_{\vec{\jmath}}^{\vec \sigma}\, e^{\im \sigma k x } \,;
 		$$
 		it then follows that
\begin{equation}
\Pi_{\mathfrak{R}_4} \comm{-\ii \vOmega(D)Z}{\ \cdot\ }\equiv 0\,,
\label{zero_on_res}
\end{equation}
where $\mathfrak{R}_4$ and $\Pi_{\mathfrak{R}_4}$ are respectively defined in \eqref{res1}, \eqref{proiettore.sui.siti}, and we obtain 
$$
\Pi_{\mathfrak{R}_4} X^+_3 = \Pi_{\mathfrak{R}_4}\left(  \hamvec{\cH_4} - \frac12 \comm{\hamvec{\cH_3}}{\hamvec{\cF_3}}\right) \stackrel{\eqref{lb.id}}{=} \Pi_{\mathfrak{R}_4} \left( \hamvec{\cH_4 - \frac 1 2\{ \cH_3,\ \cF_3\}}\right)   \stackrel{\eqref{lie3.bis}}{=} \Pi_{\mathfrak{R}_4}\ \hamvec{\cH_{4}^+}\,.
$$
Then \eqref{ide:res} follows since, by their very definition in \eqref{resj}, $\mathfrak{R}_\Set^{(n)} \subseteq \mathfrak{R}_4$ for any $n=0, \dots, 4$. To conclude the proof of \Cref{prop.ident}, it thus remains only to bound the remainder $ Y_{\geq 4}$ in \eqref{N4}.
\begin{lemma}\label{lem:Y4}
    There exist $s_0, \mu>0$ such that if $Z\in B_{s_0,\R}(r)$ then the remainder in \eqref{N4} is a quartic vector field satisfying  the bound
    \be 
\| Y_{\geq 4}(Z)\|_{\sigma}\lesssim_\sigma \| Z\|_{\sigma+\mu}^4\, , \quad \text{for any } \sigma \geq s_0 \,.
\label{est:quartic}
    \ee
\end{lemma}
\begin{proof}
By \eqref{Xzak} we have
\begin{equation}\label{est:ham.vec}
    \| \hamvec{\cH}(\zak)\|_\s\lesssim \| \zak\|_{\s+\frac32}, \quad \| \hamvec{\cH_{\geq 3}}(\zak)\|_\s\lesssim \| \zak\|_{\s+\frac32}^2, \quad \| \hamvec{\cH_{\geq 4}}(\zak)\|_\s\lesssim \| \zak\|_{\s+\frac32}^3, \quad \| \hamvec{\cH_{\geq 5}}(\zak)\|_\s\lesssim \| \zak\|_{\s+\frac32}^4\,.
\end{equation}
We start by bounding the term $\ii \vOmega(D) \bUpsilon_{\geq 3}(\zak)\zak$ in \eqref{N4}. Using that $ \vOmega(D) \in \wt \cM_0^{\frac12}$ and the bound \eqref{piove} for $\bUpsilon_{\geq 3}(\zak)\zak$, we get 
\begin{equation}
    \| \ii \vOmega(D) \bUpsilon_{\geq 3}(\zak)\zak \|_\s \lesssim \| \bUpsilon_{\geq 3}(\zak)\zak\|_{\s+\frac12} \lesssim \| \zak\|_{\s+m+\frac12}^4\stackrel{\eqref{equiv:zak_Z}}{\lesssim_s} \| Z\|_{\s+m+\nu+\frac12}^4\,.
\end{equation}
The estimate for $\hamvec{\cH_{\geq 5}}$ follows by the last of \eqref{est:ham.vec} and the equivalence \eqref{equiv:zak_Z}. The estimates of $\di (\bUpsilon_{1}(\zak)\zak) [ \hamvec{\cH_{\geq 4}}(\zak)]$, $\di (\bUpsilon_{2}(\zak)\zak) [ \hamvec{\cH_{\geq 3}}(\zak)]\,,$ and $\di (\bUpsilon_{\geq 3}(\zak)\zak) [ \hamvec{\cH}(\zak)]$ are similar; we only bound the last term.

Since $\bUpsilon(\zak)$ is a $(\nu, m)$-admissible transformation,  by  \eqref{stima.d.adm} one has
\[
\| \di (\bUpsilon_{\geq 3}(\zak)\zak) [ \hamvec{\cH}(\zak)]\|_\s= \|   \di \bUpsilon_{\geq 3}(\zak)[\hamvec{\cH}(\zak) ] \zak  + \bUpsilon_{\geq 3}(\zak) \hamvec{\cH}(\zak)\|_{\s}
		 \lesssim_\s
\|\zak\|_{\s+ m +\frac{3}{2}}^4 \stackrel{\eqref{equiv:zak_Z}}{\lesssim_\s}  \| Z\|_{\s+m+\nu+\frac{3}{2}}^4\,.
\]
     Finally, we bound the term $Y_3(\zak)- Y_3(Z)$. First note that $ Y_3(Z)= M_2(Z)Z$ where $ M_2(Z)$ is a homogeneous map in $ \wt \cM_2^{m+\frac32}$.  Using the tri-linearity of $Y_3(\,\cdot\,)= M_2(\,\cdot\,) \ \cdot\ $, its bound \eqref{smoothing}, and the bound \eqref{est:diff_zak_Z}, we get 
    $$
    \| Y_3(\zak)- Y_3(Z)\|_{\sigma}\lesssim_\sigma \| \zak \|_{\sigma+m+\frac32}^2 \| \zak-Z\|_{\sigma+m+\frac32}+\| \zak\|_{\sigma+m+\frac32} \| \zak-Z\|_{\sigma+m+\frac32}^2+ \| \zak-Z\|_{\sigma+m+\frac32}^3\lesssim_\sigma \| \zak\|_{\sigma+\nu+2m+\frac32}^4\,. 
    $$
    This proves \eqref{est:quartic} with $\mu:= \nu+2 m+\frac32$.
\end{proof}
The proof of \Cref{prop.ident} is then concluded. \qed

\section{ Paradifferential  Normal Form} 
\label{sec:paraWW}
From now on we consider \eqref{eq:etapsi} as a system
 on (a dense subspace of)    $  L^2_0(\T, \R)\times \dot L^2(\T, \R)$. Namely,  
 denoting  $ X_\gamma (\eta, \psi) $ the right hand side in \eqref{eq:etapsi}, 
 we consider  
\be\label{WWhomo}
 \partial_t (\eta, \psi) =  X_\gamma (\Pi_0^\bot \eta, \psi)  \,,
\ee
where $\Pi_0^\perp$ is the  $ L^2 $-projector onto the space of functions with zero average. 

The main result of this section is the following theorem, which states the existence of an {\em admissible transformation} (according to Definition \ref{admtra}) that, 
provided the vorticity $\gamma <0$ is resonant, i.e. $\gamma^2 \in \Q$,  puts the cubic water waves vector field  in {\em strong $\Set$-normal form} (recall Definition \ref{def:wr}), where
$\Lambda$ is any chosen $\gamma$-good set (cf. \Cref{g-good}).
\begin{theorem}[\bf Strong $\Set$-normal form of water waves]\label{thm:nf}
Let $\gamma < 0$, $\gamma^2 \in \Q$ and $\Set = \{\tm, \tn\}$ be a $\gamma$-good set according to Definition \ref{g-good}.
Let ${N \in \N}$, $N \geq 3$, and 
 $\vr \geq  3N + 28$.
There exist $s_0>\frac32, r >0$ and    
   a $(\frac12, {8})$ admissible transformation $\bUpsilon(\zak) \in \cM_{\geq 0}^{\frac12}[r]$ 
   with gain ${\und{\vr}:=\vr-3N - 2}$  
   such that
 if $(\eta,\psi)(t) \in B_{s_0, \R}(I;r)$ solves \eqref{eq:etapsi}, and $\zak=\zak(t)$ is defined in \eqref{zak}, then
the variable 
\be \label{Def:ZU}
 Z=\vect{z}{\bar{z}}:= \bUpsilon(\zak)\zak  
  \ee   
  with  $\zak = \Lmap^{-1}\vect{\eta}{\psi}$  the linear Zakharov variable defined in \eqref{zak}, solves 
	 \be \label{Z.eq}
	 \begin{aligned} 
	 	&\pa_t Z= 	 	- \im \vOmega(D)Z +
		 \vOpbw{ \im  \sm_{\geq 2}^{(\res)}(Z;x, \xi)}Z 
          +  
  X^{(\Set)}(Z) 
  +  {\bB_{\geq N}(Z) Z} + {\bR}_{\geq 3}(Z)Z
	 \end{aligned}
	 \ee
	  where $ \vOmega(D) $ is the matrix of Fourier multipliers in \eqref{diaglin},  $\sm_{\geq 2}^{(\res)}(Z; x, \xi)$  is a {\em real} symbol in $\Sigma\Gamma^{1}_{ 2}[r, 3]$ of the form 
       \begin{equation}\label{D.final}
           \sm_{\geq 2}^{(\res)}(Z; x, \xi):=  \big( \langle \,\sfV\, \rangle (Z;x) + {\sfV}_{\geq 3}(Z; x)\big) \xi + {\td}_{\geq 2}(Z; x) \omega(\xi) +  {\tf}_{\geq 2}(Z;x)\,\sign{(\xi)}+  \tg_{\geq 2}^{(-\frac12)}(Z; x, \xi)
       \end{equation}
      with $\omega(\xi)$ defined in \eqref{omegaxi} and\\
       $\bullet$  $ \langle \, \sfV\, \rangle(Z;x)$  is the  real-valued   function in $\wt \cF^\R_2$ defined by
		\be \label{VresZ}
        \begin{gathered}
        \langle \, \sfV\, \rangle(Z;x) := \sum_{j \in \Z_*} \sfV^{({\rm int})}_j |z_{j}|^2 +
		\sum_{\substack{m < 0 <n   \\ \Omega_m(\gamma) = \Omega_n(\gamma)}} 2\,  \sfV^{({\rm res})}_{m, n}\,  \Re \Big(  \bar{z_{m}(t)} \, {z_{n}}(t) \, e^{\im (n-m) x} \Big) \, ,
        \end{gathered}
		\ee 
        $\Omega_n(\gamma)$ in \eqref{omegonejin}, 
       and  for  $\tm < 0 < \tn$ and any $m < 0 < n$ 
        \be\label{Vres.coeff}
        \begin{aligned}
           &     { \sfV^{({\rm int})}_\tm:= \frac{2 \tm^2 \left(5 \tm^2-2 \tm \tn+\tn^2\right)}{{(3 \tm - \tn)(\tm - \tn) }} \ , 
          \qquad
          \sfV^{({\rm int})}_\tn:= -\frac{2 \tn^2 \left(\tm^2-2 \tm \tn+5 \tn^2\right)}{(3 \tn-\tm) (\tn-\tm)}} \\
& {\sfV^{({\rm res})}_{m,n} := 
{\frac{m n (m + n)^3}{2(n-m)\sqrt{m n(m-3n)(n-3m)}}}, \quad 
\text{ provided } \quad \gamma = - \frac{m+n}{\sqrt{2(n-m)}}
}
        \end{aligned}
        \ee
        and, more in general, for any $j \in \Z_*$
        \be
        \label{Vres.coeff2}
        {\sfV^{({\rm int})}_j := -\frac{\left(4 | j|  (\tn-\tm)+(\tm+\tn)^2\right) \left(| j|  (\tm+\tn)+j \sqrt{8 | j|  (\tn-\tm)+(\tm+\tn)^2}\right)}{2 (\tn-\tm) \sqrt{8 | j|  (\tn-\tm)+(\tm+\tn)^2}}
        }
         \ ; 
         \ee
$\bullet$ ${\sfV}_{\geq 3}(Z; x)$ is  a real-valued function in $\cF^{\R}_{\geq 3}[r]$;\\
$\bullet$ ${ \td}_{\geq 2}(Z;x), \tf_{\geq 2}(Z;x)$ are real-valued  functions in $\Sigma\cF_{2}^{\R}[r,3]$, and $\tg_{\geq 2}^{(-\frac12)} (Z; x, \xi)$  is a real-valued symbol in $\Sigma\Gamma_{2}^{-\frac12}[r,3]$;\\
$\bullet$ $X^{(\Set)}(Z)$ is a real-to-real cubic vector field   in strong-$\Set$ normal form (see Definition \ref{def:wr}), of the form 
\be\label{Xlambda}
X^{(\Set)}(Z):=  {\bR}_2^{(\Set)}(Z)Z\,, 
\ee
where ${\bR}_2^{(\Set)}(Z)$  a matrix of smoothing operators in $\wt \cR_2^{-\vr +3N+\frac32}$; 
 explicitly
\begin{equation}
\begin{gathered}
 (\Pi_{\fP^{(0)}_\Set} X^{(\Set)}) (Z) = 
\begin{pmatrix}
-\im 
\left(2\fa |z_\tm|^2 + \fb |z_\tn|^2 \right) z_\tm e^{\im \tm x} 
-\im 
\left( 2\fc |z_\tn|^2  + \fb |z_\tm|^2 \right) z_\tn e^{\im \tn x} 
\\
 \im 
\left(2\fa |z_\tm|^2 + \fb |z_\tn|^2 \right) \bar{z_\tm} e^{-\im \tm x} 
+\im 
\left( 2\fc |z_\tn|^2  + \fb |z_\tm|^2 \right) \bar{z_\tn} e^{-\im \tn x}  
\end{pmatrix}   , 
\quad 
\fa, \fb , \fc \mbox{ in } \eqref{abc} 
\  ,  \\
\Pi_{\fP_\Set^{(1)}} X^{(\Set)}=0, \qquad 
 \Re\left( [\Pi_{\fP_\Set^{(2)}}X^{(\Set)}(Z)]^+_j    \, \bar z_j \right) = 0 \ \ \forall j \in \Z_* \ ;
\end{gathered}
\label{Y.str}
\end{equation}
 $\bullet$   $\bB_{\geq N}(Z)$ is a real-to-real matrix of spectrally localized maps in $\cS^{0}_{\geq N}[r]$ and ${\bR}_{\geq 3}(Z)$ is a real-to-real matrix of smoothing operators in $\cR_{\geq 3}^{-\varrho+3N+28}[r]$\,.
\end{theorem}

\begin{remark} 
1. 
As in \Cref{rem:not.symm}, we can expand 
$\langle \, \sfV\, \rangle$ as in \eqref{espr.hom.sym}, with symmetric coefficients given by 
\begin{equation*}
    V_{j_1,j_2}^{\sigma_1,\sigma_2}= \begin{cases}
        \frac12\sfV^{({\rm int})}_{j_1}& \mathrm{if} \ j_1=j_2 \ \mathrm{and} \ \sigma_1\sigma_2=-1\\
        \frac12\sfV^{({\rm res})}_{m, n}& \mathrm{if} \ m:=\min(j_1,j_2)<0<\max(j_1,j_2)=:n, \ \Omega_{j_1}(\gamma) = \Omega_{j_2}(\gamma) \ \mathrm{and} \ \sigma_1 \sigma_2=-1\\
        0 &\mathrm{otherwise}\,.
    \end{cases}
\end{equation*}
2. The coefficients $\sfV^{({\rm int})}_j $ in \eqref{Vres.coeff2} depend on $\tn, \tm$ only through $\gamma = - \frac{\tm+\tn}{\sqrt{2(\tn-\tm)}} $.
\end{remark}

The rest of this long section is devoted to the proof of Theorem \ref{thm:nf}. Let us comment on the procedure:
\begin{enumerate}
    \item In Section \ref{subsec:pretras}, we perform some preliminary transformations in order to rewrite the original water waves equations \eqref{eq:etapsi} in paradifferential complex form, obtaining the system \eqref{complexo}. This is achieved in three steps: first, we rewrite the water waves equations in terms of Alinhac’s good unknown in \Cref{laprimapara}; next, we pass to Wahl\'en coordinates in \Cref{wahlen}; finally, we introduce complex variables in \Cref{LemCompl}.
    \item In Section \ref{subsec:quadraticNF} we put the system into its paradifferential normal form by removing the quadratic terms of the vector field and thus reducing it to a cubic one. This yields \Cref{thm:quadratic.nf}.
This reduction can be carried out for any value of the vorticity $\gamma$, thanks to the fact that, by \Cref{3onde}, there are no three-wave resonances. The procedure consists of several steps: in \Cref{diag} we block-diagonalize the operator at order $\frac12$, while in \Cref{diag.ord0} we remove the off-diagonal symbols up to homogeneity $N \in \N$.
This represents an important difference with respect to the normal form constructions in \cite{BD, BFP, BMM2}. Although in principle one could remove the off-diagonal symbols at degree of homogeneity, we deliberately stop at a finite homogeneity. The reason is that we want the transformation to be admissible. Since at this stage the transformation is given as a composition of flows of the form \eqref{flussoG}, with generators $G(U)= \zOpbw{g(U;\cdot)}$, we need, in view of \Cref{lem:flow.ad}, the generator to be homogeneous.
As a consequence, we remove only homogeneous terms. The order $N$ will be fixed equal to $3$ in the next section (see just above \eqref{s0r}), but we keep it abstract here for future applications.
    Next, in \Cref{ridlinord1}, we remove the linear terms from the transport equation, reducing it to a quadratic transport, in \Cref{lem:rid12} we remove the linear terms from the symbols of order $\frac12$ and $0$, and in \Cref{riduzione.negativi} the linear terms from the symbols of negative order.
    Then, in \Cref{poinc} we remove the quadratic terms from the smoothing remainder, completing the reduction of the quadratic vector field to a cubic one.
    \item In Section \ref{subsec:cubicNF} we perform a cubic resonant normal form. 
    We fix the vorticity $\gamma<0$ to a resonant value, i.e. $\gamma^2 \in \Q$.
    By \Cref{lem:goodset} there exist infinitely many $\gamma$-good sets, and we fix one $\Lambda := \{\tm, \tn\}$,  $\tm <0 <\tn $, fulfilling (G1)--(G4) of \Cref{g-good}. 
    In particular we have the 2-wave resonant interaction  $\Omega_\tm(\gamma) = \Omega_\tn(\gamma)$ (see \eqref{Omega*}). 
  In \Cref{thm:cubic.nf.1} we remove the non-resonant monomials from the quadratic transport, reducing to the resonant transport term $\langle \sfV \rangle(Z;x)$ in \eqref{VresZ}, and remove non-resonant monomials with at most two frequencies outside $\Lambda=\{\tm, \tn\}$ 
  from the cubic smoothing remainder.
This is possible thanks to the resonances analysis of \Cref{lem:wres}, valid for $\gamma$-good sets. 
    \item Finally, in Section \ref{subsec:proofthm}, we prove \Cref{thm:nf} via an identification argument based on \Cref{prop.ident} 
    and using the formal $\Lambda$-normal form computed in \Cref{prop:cH^{(4)}_+}.
\end{enumerate}

Before turning to the proof of \Cref{thm:nf}, we establish a local well posedness result for the Cauchy problem for $Z(t)$ solving equation \eqref{Z.eq}.  
\begin{lemma}[\bf Local existence of the Cauchy problem for $Z(t)$]\label{loc.ex}
    There exists $\sigma_0 >0$ and for any $\s \geq  \s_0$, there exists $\tr := \tr(\s)>0$ such that the following holds. Let $Z_0 \in B_{\s_0}(\tr) \cap H^\sigma_\R(\T, \C^2)$. Then there exist $T_{\rm loc}>0$, $C_\sigma>0$ and a unique solution $Z(t) \in C \left([- T_{\rm loc},\, T_{\rm loc}];\ B_{\s_0}(C_\s\tr) \cap H^\sigma_\R(\T, \C^2)\right)$ of \eqref{Z.eq} with initial datum $Z(0) = Z_0$. 
\end{lemma}

The proof of \Cref{loc.ex} is essentially based on the equivalence between the $H^s$-norm of $Z(t)$ in \eqref{Def:ZU} and  the $X^{s-\frac 3 4}$-norm of $(\eta(t), \psi(t), \sfV(t), \sfB(t))$, and is contained in Appendix \ref{app:cauchy.z}.

\smallskip

Finally, along the section we shall also use that  equation 
\eqref{zetone:formale} for the variable $\zak$ has the form
\be\label{claim:zak.M}
\pa_t \zak = \cX_\cH(\zak) =  - \im {\bf \Omega}(D)\zak + \bM_{\geq 1}(\zak) \zak \ , \quad 
 \bM_{\geq 1}(\zak)\in \Sigma \cM^{\frac32}_1[r,N],
\ee
for any $N$ in $\N$.
To prove this,   use that by \cite[Lemma 5.1]{BMM2} with the capillarity $\kappa = 0$, the water waves equations \eqref{eq:etapsi} have the form
\begin{equation}\label{det.map}
\pa_t \vect{\eta}{\psi} = \begin{pmatrix}
    0 & G(0) \\
    -1 & \gamma G(0) \pa_x^{-1}
\end{pmatrix}
\vect{\eta}{\psi} + \breve{\bM}_{\geq 1}(\eta, \psi) \vect{\eta}{\psi} \ , \qquad 
\breve\bM_{\geq 1}(\eta, \psi)  \in \Sigma \cM^{1}_1[r,N]
\end{equation}
for any $N$ in $\N$.
Then, in view of \eqref{zak} and \eqref{app_L}, 
$
\pa_t \zak = - \im {\bf \Omega}(D)\zak
 + \bM_{\geq 1}(\zak) \zak$ 
where 
$\bM_{\geq 1}(\zak):= \Lmap^{-1} 
\breve \bM_{\geq 1}
\left(\Lmap\zak\right) \Lmap
\in \Sigma \cM^{\frac32}_1[r,N]$.

\subsection{Preliminary transformations}\label{subsec:pretras}
We begin by introducing and studying some real-valued functions that will be used in the sequel.
\begin{lemma}\label{lem:B.V.a} 
Let  $\sfV, \sfB$ be the real-valued functions defined in \eqref{def:V}--\eqref{form-of-B}.
Define the real-valued function  
\begin{equation}
    \label{def:Vgamma}
    \sfV_\gamma := \sfV_\gamma(\eta,\psi;x) := \sfV(\eta,\psi;x) - \gamma \eta (x), 
\end{equation}
and the  real-valued {\em Taylor coefficient}
\begin{equation}\label{a.taylor}
   \sfa:=\sfa(x):= \sfa(\eta, \psi; x):= (\pa_t \sfB)(\eta, \psi;x)+ \sfV_\gamma(\eta, \psi; x) (\pa_x \sfB)(\eta, \psi; x)\,,
\end{equation}
where  
$(\pa_t \sfB)(\eta, \psi) \equiv \di \sfB(\eta, \psi)[ \cX_{\scH_\gamma}(\eta, \psi)]$ and 
$ \cX_{\scH_\gamma}$ is the water waves Hamiltonian vector field in \eqref{HamWW}. 
Then the following properties hold.
    \begin{itemize}
\item[$(i)$] There is $\s_0>3$ such that for any $\s \geq \s_0$, there is $r' = r'(\s)>0$ so that 
\be\label{BVa.ana}
(\eta, \psi) \mapsto (\sfB, \sfV, \sfV_\gamma, \sfa) \mbox{ is analytic }
B_{\s}(r')\times B_{\s}(r')\to \left(H^{\s-1}(\T; \R)\right)^3 \times H^{\s-2}(\T; \R)
\ee
with estimates
\begin{equation}
\begin{gathered}\label{a.v.b.estimate}
\norm{\sfB}_{\s-1} + \norm{\sfV}_{\s-1}\lesssim_\s \|\psi\|_{\s}\,, \quad \norm{\sfa}_{\s - 2} + \norm{\sfV_\gamma}_{\s-1} \lesssim_\s
 \|\eta\|_{\s} + \|\psi\|_{\s} \ ,  \\
   \forall (\eta, \psi) \in B_{\s}(r')\times  B_{\s}(r')\,.
\end{gathered}
\end{equation}
 \item[$(ii)$] Let $N \in \N$. There is $r>0$ such that   $\sfV, \sfB, \sfV_\gamma, \sfa$  belong to $\Sigma \cF^\R_1[r,N]$.
\item[$(iii)$] One has the Taylor expansions
\begin{gather}\label{esp.b.gamma.1.2}
   \sfB(\eta,\psi)= |D|\psi +  
   \left(-\eta\psi_{xx}-|D|(\eta|D|\psi) \right)
   + \sfB_{\geq 3}(\eta,\psi)\,, \\
   \label{esp.v.gamma.1.2}
     \begin{aligned}
    \sfV_\gamma(\eta, \psi) =& (\psi_x - \gamma \eta) -\eta_x |D| \psi + \sfV_{\gamma, \geq 3}(\eta, \psi)\,,  
     \end{aligned}
\end{gather}
with 
$ \sfB_{\geq 3}$ and 
   $\sfV_{\gamma,\geq 3}$  real-valued functions in $ \Sigma \cF^{\R}_{3}[r, N] $\,.
\item[$(iv)$] There is $\s_0> \frac72$ such that for any $\s \geq \s_0$, there is $r' = r'(\s)>0$ so that  if $(\eta, \psi, \sfV, \sfB) \in B_{X^{\sigma}}(I;r')$ is a solution of  \eqref{eq:etapsi}, then (recall \eqref{norm:Xs})
\be\label{stima.der.BVa}
 \norm{\pa_t \sfB}_{\s -\frac 3 2} + \norm{\pa_t \sfV}_{\s-\frac 3 2}+ \norm{\pa_t \sfV_\gamma}_{\s -\frac 3 2} + \norm{\pa_t \sfa}_{\s - \frac 5 2} \lesssim_\s 
 \|(\eta, \psi, \sfV, \sfB)\|_{X^{\s}}\,.
 \ee
    \end{itemize}
\end{lemma}
\begin{proof}
    We postpone the proof to \Cref{subsec:preliminary}, see page \pageref{proof.lemma5.3}.
\end{proof}
\begin{remark}\label{rem:taylor}
   The Taylor coefficient $\sf a(x)$ in \eqref{a.taylor} is exactly the same that appears in formula \eqref{taylor.explicit},  using $\vec u =  \vect{-\gamma y}{0}+ \grad \Phi$ and the definitions \eqref{def:V}, \eqref{form-of-B} and \eqref{def:Vgamma}. 
   Let us give a quick proof of \eqref{taylor.cond}.
On the free surface $y=\eta(t,x)$ the pressure satisfies $P=0$, that, differentiated in the outward unit normal 
$\vec n=\frac{(-\eta_x,1)}{\sqrt{1+\eta_x^2}}$ to the fluid domain, yields 
$  \partial_{\vec n} P = \sqrt{1+\eta_x^2}\, P_y$ on  $y=\eta(x)$. Then \eqref{taylor.cond} follows using the 
vertical component of the Euler equations \eqref{euler} (with $g=1$).
\end{remark}

It is convenient to write the  water-waves equations \eqref{eq:etapsi} in the variables $(\eta, \upomega)$ where $\upomega$ is the  “good unknown” of Alinhac, introduced by Alazard-Metivier \cite{AlM} and defined as 
\begin{equation}\label{GU}
\begin{pmatrix}{\eta}\\{\upomega} \end{pmatrix} = \cG(\eta,\psi) \begin{pmatrix}{\eta}\\{\psi} \end{pmatrix}:= \begin{pmatrix}{\eta}\\{\psi- \Opbw{\sfB(\eta,\psi;x)} \eta} \end{pmatrix} \,. 
\end{equation}
\begin{lemma}{\bf (Water-waves equations in $(\eta, \upomega)$  
variables)}\label{laprimapara} 
Let  $N \in \N $ and $ \varrho_1  \geq 1$. 
There exist $ \s_0 > \frac32$ and $r > 0 $  such that,
if $ (\eta, \psi) \in B_{{\s}_0, \R} (I;r)  $ solves \eqref{eq:etapsi}, then $(\eta, \upomega)$ defined in \eqref{GU} solve
\begin{align}
 \label{ParaWW}
\pa_t \begin{pmatrix} \eta \\ \upomega \end{pmatrix}= &\begin{pmatrix} 0 & |D| \\ -1 &   \gamma \Hilb \end{pmatrix} \begin{pmatrix} \eta \\ \upomega \end{pmatrix} \\
& + \Opbw {\begin{bmatrix} - \im  \xi\#_{\vr_1} \sfV_\gamma    & 0 \\ -\sfa  &  -\sfV_\gamma\#_{\vr_1}\im  \x  \end{bmatrix} + \begin{bmatrix} 0 & b^{(-1)}_{\geq 1}   \\ 0 
  &  \gamma \frac{1}{\ii \xi}\#_{\vr_1} b_{\geq 1}^{(-1)} \end{bmatrix}} \begin{pmatrix} \eta \\ \upomega \end{pmatrix}+ \bR(\eta,\psi) \begin{pmatrix} \eta \\ \upomega \end{pmatrix} \notag
\end{align}
where $\Hilb$ is the Hilbert transform in \eqref{def:Hilbert}, $  \sfV_\gamma $ and
$\sfa$ are the real-valued functions defined respectively in  \eqref{def:Vgamma} and in \eqref{a.taylor}, and
\begin{itemize} 
 \item the symbol $b^{(-1)}_{\geq 1}$ belongs to $\Sigma \Gamma^{-1}_1[r, N]$;  its homogeneous components $b_{p}^{(-1)} := \cP_{p} \big[b_{\geq 1}^{(-1)}\big]$ (cf. the bullet below \Cref{omosmoothing}) are real-valued and even in $\xi$; moreover for any $\sigma \geq \s_0$, there is $r' = r'(\s)>0$ so that  if $(\eta, \psi, \sfV, \sfB) \in B_{X^{\sigma}}(I;r')$ is a solution of  \eqref{eq:etapsi}, then  for any $t\in I$ one has estimate
 \be\label{b-1dert}
  |b^{(-1)}_{\geq 1}|_{-1, {W^{\s-\sigma_0,\infty}},M}  + | \pa_t b^{(-1)}_{\geq 1}|_{-1, {W^{\s-\sigma_0,\infty}},M}\lesssim_{M} \| (\eta,\psi,\sfV,\sfB)\|_{ X^{\s}} \,, \quad \forall M \in \N\,.
  \ee
 
\item  $\bR(\eta, \psi)$ is a matrix of real  smoothing  operators in 
$\Sigma\mathcal{R}^{-\varrho_1 + 1}_{1}[r, N]$; 
 moreover, for any $\sigma \geq \s_0$,  if 
 $(\eta, \psi, \sfV, \sfB) \in B_{X^{\sigma_0}}(I;r) \cap X^{\s}$ is a solution of  \eqref{eq:etapsi}, then
\begin{equation}\label{Rt_smooth}
\begin{aligned}
      \| \pa_t \left(  \bR(\eta, \psi) \vect{\eta}{\upomega} \right) \|_{\s + \vr_1 {- 3/2}} &\lesssim_{\s} \| (\eta, \psi, \sfV, \sfB)\|_{X^{{\s}_0}}
      \| (\eta, \psi, \sfV, \sfB)\|_{X^\s}\,.
\end{aligned}
\end{equation}
\end{itemize}
 Moreover, the function $\pa_t\upomega$ has the following properties: for any $\s \geq \s_0 + \frac 3 2$, if $(\eta, \psi, \sfV, \sfB) \in B_{X^{\s_0}}(I; r) \cap X^\s$ is a solution  of \eqref{eq:etapsi},
 then 
\begin{equation}\label{pa_t:omega}
        \|\pa_t \upomega\|_{\s-\frac12} \lesssim_\s \| (\eta, \psi, \sfV, \sfB)\|_{X^\s}\,.
    \end{equation}
\end{lemma}
\begin{proof}
We start with proving that $\pa_t \upomega$ satisfies \eqref{pa_t:omega}. By the very definition of $\upomega$ in \eqref{GU} and using \Cref{thm:contS}-Item$(i)$, one has
\begin{align*}
    \|\upomega\|_{\s + \frac 1 2} = \| \psi - \Opbw{\sfB} \eta\|_{\s + \frac 1 2} \lesssim_{\s} \|\psi\|_{\s+\frac 1 2} +  \|\sfB\|_{\s_0} \| \eta\|_{\s + \frac 1 2}\,,
\end{align*}
as well as 
\begin{align*}
    \|\pa_t \upomega\|_{\s -\frac 1 2} = \|\pa_t \psi - \Opbw{\pa_t \sfB} \eta - \Opbw{\sfB} \pa_t \eta\|_{\s - \frac 1 2} \lesssim_{\s} \|\pa_t \psi\|_{\s-\frac12} + \|\pa_t \sfB\|_{\s_0} \|\eta\|_{\s - \frac 1 2} + \|\sfB\|_{\s_0} \|\pa_t \eta\|_{\s - \frac 1 2}\,.
\end{align*}
Then combining the above estimate with estimates \eqref{a.v.b.estimate}, \eqref{stima.der.BVa}, and \eqref{pa_t:etapsi} on $\sfB$, $\pa_t \sfB$, and $(\pa_t \eta, \pa_t \psi)$ respectively, one gets for any  
 $(\eta, \psi, \sfV, \sfB) \in B_{X^{\sigma_0}}(I;r) \cap X^{\s}$  any $t \in I$,
 \begin{equation}\label{omega-small}
     \|\upomega\|_{\s + \frac 1 2}  \lesssim_\s \|(\eta, \psi, \sfV, \sfB)\|_{X^{\s}}\,, \qquad \|\pa_t \upomega\|_{\s - \frac 1 2}  \lesssim_\s \|(\eta, \psi, \sfV, \sfB)\|_{X^{\s}}\,,
 \end{equation}
proving in particular \eqref{pa_t:omega}.

\noindent We now prove that \eqref{ParaWW}--\eqref{Rt_smooth} hold.
Let us decompose the vector field in \eqref{eq:etapsi} as $\cX_{\scH_\gamma}(\eta, \psi) := \cX_1(\eta, \psi) + \gamma \cX_2(\eta, \psi) $, where $\cX_1$ contains all the terms independent of $\gamma$ and $\gamma \cX_2$ the linear terms in $\gamma$. 
Then, recalling the definition of $\cG$ in \eqref{GU}, $(\eta, \psi)$ solves \eqref{eq:etapsi}
 if and only if $(\eta, \upomega)$ solves
 \begin{equation}\label{om.vf}
     \partial_t \vect{\eta}{\upomega}
   = \cG \cX_1(\eta,\psi)+(\partial_t \cG)\, \cG^{-1} \vect{\eta}{\upomega} + \gamma\cG \cX_2(\eta,\psi) \,. 
 \end{equation}
  From Proposition 3.1 of \cite{BFP}, the paralinearization of the terms in \eqref{om.vf} which do not contain $\gamma$ is given by 
\begin{equation}
\label{non.sono.eta.om}
\cG \cX_1   
+(\partial_t \cG)\, \cG^{-1} \vect{\eta}{\upomega}
= \Opbw{ \begin{bmatrix}
    - \im \xi \#_{\vr_1} \sfV & |\xi| + b^{(-1)}_{\geq 1} \\
    - 1 -\pa_t \sfB  - \sfV \sfB_x  &  - \sfV\#_{\vr_1} \im \xi  
\end{bmatrix}}\vect{\eta}{\upomega} + \bR(\eta, \psi)\vect{\eta}{\upomega}\,,
\end{equation}
where $b_{\geq1}^{(-1)} \in \Sigma \Gamma_1^{-1}[r, N]$ and  $\bR(\eta, \psi)$ is a matrix of smoothing operators in $\Sigma \cR_1^{-\varrho_1}[r, N]$. 
The fact that the symbol $b^{(-1)}_{\geq 1}$ has  homogeneous components $b_{p}^{(-1)} := \cP_{p}[b_{\geq 1}^{(-1)}]$  real-valued and even in $\xi$ and satisfies \eqref{b-1dert} is proved in \Cref{lemma:para.bella}. 
The estimate \eqref{Rt_smooth} for the smoothing remainder  $\bR(\eta, \psi)$ in \eqref{non.sono.eta.om} follows arguing as in the proof of \eqref{pa_t:smooth}-\eqref{pa_t:smooth2}, and by additionally using the estimate \eqref{omega-small} for $\upomega$ and $\pa_t \upomega$. 
We now paralinearize $\cX_2$, which has the form
$$
\cX_2(\eta, \psi) :=  \begin{bmatrix}
    \eta \eta_x   \\
    \eta \psi_x + \partial_x^{-1} G(\eta)\psi
\end{bmatrix} \,.
$$
Using Lemma \ref{bony}, item $(i)$ of
Proposition \ref{teoremadicomposizione}, and the identity $ \eta_x= \Opbw{\ii \xi} \eta$, we get 
\be\label{ga1}
\begin{aligned}
 \eta \eta_x 
  &=
  \Opbw{ \im \xi \#_{\vr_1} \eta }\eta  + {{R}(\eta)\eta}
  \end{aligned}
\ee
where $ {R}(\eta) $ is a homogeneous smoothing operator in $ \widetilde{\mathcal{R}}^{-\vr_1}_1$, which using \eqref{piove} 
satisfies the following: there is $s_0>0$ and for any $\s \geq s_0$, there are $C>0$ and $r = r(\s)$ such that for any $\eta \in B_{s_0}(r) \cap H^{\s}(\T;\R)$, any   $v \in C^1(I; H^\s(\T;\C))$,
\begin{equation}\label{0.resti}
  \begin{aligned}
    \| \pa_t ( R(\eta)  v) \|_{\s + \vr_1-\frac 1 2} &\lesssim_{\s} \| R(\pa_t \eta) v \|_{\s + \vr_1-\frac 1 2} + \|  R(\eta) \pa_t v \|_{\s + \vr_1-\frac 1 2} \\
    &\lesssim_{\s} \|\pa_t \eta\|_{s_0}\| v\|_{\s-\frac12}  +  \| \pa_t \eta\|_{\s-\frac12}\| v\|_{s_0}+ \|\eta \|_{\s-\frac12}\| \pa_t v\|_{s_0} +  \|\pa_t v \|_{\s-\frac 1 2} \| \eta \|_{s_0} \,.
\end{aligned}  
\end{equation}
In particular, if $v = \eta$, using Lemma \ref{lem:tutto.X},  we get  provided $\s_0 \geq s_0 +1$,
\begin{align}
\label{un.resto}
      \|\pa_t (R(\eta)  \eta) \|_{\s + \vr_1-\frac 1 2} \lesssim_{s}  \|(\eta, \psi, \sfV, \sfB)\|_{X^{\s_0}}
      \|(\eta, \psi, \sfV, \sfB)\|_{X^\s}\,.
\end{align}
By Lemma  \ref{bony}, item $(i)$ of Proposition \ref{teoremadicomposizione}, \eqref{comp3fin}--\eqref{sharp3}, and by using \eqref{GU}, we have
\begin{equation}
\label{pezzo.1.om}
\begin{aligned}
    \eta \psi_x 
    &= \Opbw{\eta\#_{\vr_1} \im   \xi} \psi + \Opbw{\psi_x} \eta + R(\eta, \psi),\\
    &= \Opbw{ \eta\#_{\vr_1} \im   \xi} \upomega + \Opbw{(\eta\#_{\vr_1} \im   \xi)\#_{\vr_1} \sfB+\psi_x}\eta + R(\eta, \psi)+ R'(\eta,\psi)\eta 
\end{aligned}
\end{equation}
where $R(\eta, \psi)=R_1(\eta)\psi + R_2(\psi)\eta$ with $R_j$ homogeneous operators in $ \wt \cR_1^{-\vr_1}$  for $j=1,2$, and  $R'(\eta,\psi) = \cQ(\eta \#_{\vr_1} \ii \xi, \sfB)$ is in $\Sigma \cR_1^{-\vr_1 + 1}[r, N]$. We now check that $R(\eta, \psi),\ R'(\eta, \psi)$ satisfy estimates \eqref{Rt_smooth}. Indeed, using \eqref{piove} and Lemma \ref{lem:tutto.X}, arguing as to obtain \eqref{un.resto} one gets
\begin{equation}
\label{due.resti}
  \begin{aligned}
    \|\pa_t (R(\eta, \psi)  \eta) \|_{\s + \vr_1-\frac 1 2}
    \lesssim_{\s}  \|(\eta, \psi, \sfV, \sfB)\|_{X^{\s_0}}
    \|(\eta, \psi, \sfV, \sfB)\|_{X^\s}\,.
\end{aligned}  
\end{equation}
By Leibniz rule, one checks that $\pa_t \cQ(a, b) = \cQ(\pa_t a, b) + \cQ(a, \pa_t b)$; thus by estimate \eqref{comp020} for $\cQ$, 
estimate \eqref{a.v.b.estimate} for $\sfB$, estimate \eqref{stima.der.BVa} for $\pa_t\sfB$
and by Lemma \ref{lem:tutto.X}  we have 
\begin{align}
\| \pa_t (R'(\eta,\psi) \eta)\|_{\s+\vr_1-{\frac 3 2}} 
&\lesssim_\s
\| (\eta,\psi, \sfV,\sfB)\|_{X^{\s_0}}^2 \| (\eta,\psi, \sfV,\sfB)\|_{X^{\s}}\,,
\label{tre.resti}
\end{align}
which proves \eqref{Rt_smooth}.
Furthermore, by Lemma \ref{lemma:para.bella} (noting that $-\ii \sfV\xi -\frac12 V_x=-\ii \xi\#_{\vr_1} \sfV$) and recalling \eqref{GU}, we have 
\begin{align*}
G(\eta) \psi = &|D| \upomega + \Opbw{- \im \xi \#_{\vr_1} \sfV }\eta + \Opbw{b_{\geq 1}^{(-1)}} \upomega +R(\eta) \psi \,
\end{align*}
with   $b_{\geq 1}^{(-1)} \in \Sigma \Gamma^{-1}_1[r, N]$, $R(\eta)$ satisfying \eqref{pa_t:smooth}. Then, by    Lemma \ref{bony},
Proposition \ref{teoremadicomposizione} and recalling \eqref{def:Hilbert}, one has
\begin{equation}
\label{pezzo.2.om}
\begin{aligned}
  \partial_x^{-1} G(\eta) \psi
  =& \Hilb \upomega  + \Opbw{-\sfV} \eta + \Opbw{\frac{1}{\ii \xi} \#_{\vr_1}  b_{\geq 1}^{(-1)}} \upomega + \pa_x^{-1}(R(\eta)\psi) + R''(\eta, \psi) \eta + R'''(\eta) \upomega\,,
\end{aligned}
\end{equation}
where  $R''(\eta, \psi), R'''(\eta) $ are real smoothing operators in $\Sigma \cR^{-\vr_1}_{1}[r, N]$ having  the form (cfr. \eqref{comp01A})
$$
R''(\eta, \psi) = \pa_x^{-1} \cQ\left({\ii \xi},\,  \sfV\right), \qquad R'''(\eta)=\cQ\left(\frac{1}{\ii \xi},\, b_{\geq 1}^{(-1)}\right).
$$
Then  $R''$ and $R'''$ satisfy estimates \eqref{Rt_smooth} thanks to the properties of $\cQ(\cdot, \cdot)$ in   \eqref{comp020}, the bound for  
$\pa_t b_{\geq 1}^{(-1)}$ in   \eqref{b-1dert},  for $\pa_t \sfV$ in \eqref{stima.der.BVa}, for $\pa_t \eta$ in \eqref{pa_t:etapsi} and $\pa_t \upomega$ in \eqref{pa_t:omega}.

Therefore by \eqref{ga1}, \eqref{pezzo.1.om} and \eqref{pezzo.2.om}, and collecting all the estimates on the smoothing remainder terms, 
one has
$$
\gamma \cX_2 = \gamma \, \Opbw{\begin{bmatrix}
    \ii \xi \#_{\vr_1}  \eta   & 0 \\ 
   ( \eta  \#_{\vr_1} \im \xi )\#_{\vr_1} \sfB  +\psi_x  -\sfV  & -\im\,  \sign(\xi) + \eta \#_{\vr_1} \im \xi+ \frac{1}{\ii \xi}\#_{\vr_1} b_{\geq 1}^{(-1)}
\end{bmatrix}} \begin{pmatrix} \eta \\ \upomega \end{pmatrix}\, + \gamma \bR_1(\eta, \psi)\vect{\eta}{\upomega}\,,
$$
where $\bR_1$ is a matrix of real smoothing operators in $\Sigma\cR_1^{-\vr_1}[r, N]$ whose time derivative satisfies \eqref{Rt_smooth}. Using Lemma \ref{bony}, Proposition \ref{teoremadicomposizione}, and a straightforward computation, along with \eqref{def:V} and the identity 
$$
(\eta\#_{\vr_1} \ii \xi) \#_{\vr_1} \sfB- \sfB\#_{\vr_1} (\ii \xi \#_{\vr_1} \eta)= \eta \sfB_x - \sfB \eta_x \stackrel{\eqref{def:V}}{=}  \eta \sfB_x + \sfV - \psi_x 
$$
 we obtain that
\begin{align}
    \gamma \cG \cX_2   =& \gamma \, \Opbw{\begin{bmatrix}
    \ii \xi \#_{\vr_1} \eta   & 0 \\ \sfB_x \eta  & -\ii \sign(\xi) + \eta\#_{\vr_1}\ii  \xi + \frac{1}{\ii \xi}\#_{\vr_1} b_{\geq 1}^{(-1)}
\end{bmatrix}} \begin{pmatrix} \eta \\ \upomega \end{pmatrix}+  \bR'(\eta, \psi)\begin{pmatrix} \eta \\ \upomega \end{pmatrix}\,,
\label{metti.vort}
\end{align}
where
$$
\bR'(\eta, \psi) = \gamma \cG \bR_1(\eta, \psi) + \gamma \bR_2(\eta, \psi)\,, \quad \bR_2(\eta, \psi) := \begin{pmatrix}
    0 & 0 \\
    -\cQ(\sfB,\ \ii \xi \#_{\vr_1} \eta) & 0
\end{pmatrix}\,,
$$
is a matrix of real smoothing operators in $\Sigma\cR_1^{-\vr_1}[r, N]$.  Using Leibniz rule, estimate \eqref{a.v.b.estimate}, Lemma \ref{lem:tutto.X}, and the fact that $\bR_1(\eta, \psi) \vect{\eta}{\upomega}$ satisfies \eqref{Rt_smooth}, one has that $\cG \bR_1(\eta, \psi) \vect{\eta}{\upomega}$ satisfies \eqref{Rt_smooth}. Furthermore, the corresponding estimate for $\bR_2 (\eta, \psi) \vect{\eta}{\upomega}$ follows arguing as to obtain \eqref{tre.resti}. Thus also $\bR'(\eta, \psi) \vect{\eta}{\upomega}$ satisfies \eqref{Rt_smooth}.

Then the thesis follows by \eqref{om.vf} summing the vector fields in \eqref{metti.vort} and in \eqref{non.sono.eta.om}.
\end{proof}
\paragraph{Wahl\'en coordinates.} 
We now write  system \eqref{ParaWW}  in the Wahl\'en coordinates  defined in \eqref{Whalen}. 

\begin{lemma}\label{wahlen}
{\bf (Water-waves equations in Wahl\'en variables)}
Let  $N \in \N $ and $ \varrho_1  \geq 1$. 
There exist $ s_0 ,r > 0 $ such that,  
if $ (\eta, \psi) \in B_{s_0, \R} (I;r)  $ solves \eqref{eq:etapsi}
then 
\begin{equation}
 (\eta, \underline{\wahlen}) := \cW^{-1} (\eta, \upomega)      \ ,  \quad \  \cW^{-1}\, \mbox{ in } \eqref{Whalen}
\label{BISMA_wahlen}
 \end{equation}
  solves 
\begin{align}
\pa_t \begin{pmatrix}
\eta \\ \underline{\wahlen} 
\end{pmatrix}=&
\begin{pmatrix}
\frac{ \gamma}{2} \Hilb & |D| \\ 
-(1 + \frac{\gamma^2}{4}|D|^{-1} ) & \frac{ \gamma}{2} \Hilb \end{pmatrix} 
\begin{pmatrix}
\eta \\ \underline{\wahlen} 
\end{pmatrix} +  \Opbw {\begin{bmatrix}
- \im   \xi\#_{\vr_1} \sfV_\gamma   
& 0  
\\ 
-\sfa
& - \sfV_\gamma\#_{\vr_1} \im  \xi
\end{bmatrix}}\begin{pmatrix}
\eta \\ \underline{\wahlen} 
\end{pmatrix}\notag \\
&+ \Opbw{\begin{bmatrix}  \frac{\gamma}{2}b_{\geq 1}^{(-1)} \#_{\vr_1} \frac{1}{\ii \xi}&  b_{\geq 1}^{(-1)}  \\  \frac{\gamma^2}{4} \frac{1}{\ii \xi} \#_{\vr_1} b_{\geq 1}^{(-1)} \#_{\vr_1} \frac{1}{\ii \xi}&  
\frac{\gamma}{2} \frac{1}{\ii \xi} \#_{\vr_1}b_{\geq 1}^{(-1)}
\end{bmatrix}}\begin{pmatrix}
\eta \\ \underline{\wahlen} 
\end{pmatrix}+ \bR(\eta,\psi) \begin{pmatrix}
\eta \\ \underline{\wahlen} 
\end{pmatrix}  \label{WWW}
\end{align}
where $\Hilb$ is the Hilbert transform in \eqref{def:Hilbert},   $\sfV_\gamma \in \Sigma {\mathcal F}^\R_{1}[r, N] $ is defined in \eqref{def:Vgamma}, 
$\sfa \in \Sigma \cF^\R_{1}[r, N]$ is the Taylor coefficient defined in \eqref{a.taylor}, 
  $b_{\geq 1}^{(-1)} \in \Sigma \Gamma_{1}^{-1}[r, N]$ is the   symbol in \Cref{laprimapara} and 
  \begin{itemize}
\item the matrix of symbols in \eqref{WWW} satisfies,  for any $p=1,\ldots, N-1$, 
\begin{align}
 & \frac{1}{\ii \xi} \#_{\vr_1} b_{p}^{(-1)} \#_{\vr_1} \frac{1}{\ii \xi}\in \R \,,   &&\frac{1}{\ii \xi} \#_{\vr_1} b_{p}^{(-1)} \#_{\vr_1} \frac{1}{\ii \xi}=\big( \frac{1}{\ii \xi} \#_{\vr_1} b_{p}^{(-1)} \#_{\vr_1} \frac{1}{\ii \xi}\big)^\vee\,, \label{real:ham:sym0}\\
&\ov{b_{p}^{(-1)} \#_{\vr_1} \frac{1}{\ii \xi}}^\vee =b_{p}^{(-1)} \#_{\vr_1} \frac{1}{\ii \xi} &&\big(b_{p}^{(-1)} \#_{\vr_1} \frac{1}{\ii \xi}\big)^\vee =-\frac{1}{\ii \xi}\#_{\vr_1} b_{p}^{(-1)}\,.
 \label{real:ham:sym}
\end{align}
\item  $\bR(\eta,\psi) $ 
is a matrix of real 
smoothing operators in $ \Sigma \cR^{-\vr_1 + 1}_{1}[r, N]$. 
\end{itemize}

\end{lemma}

\begin{proof}
Recall that $\cW, \cW^{-1}$ are matrices of Fourier multipliers in $\wt\Gamma_0^0$.
By  \eqref{BISMA_wahlen} and \eqref{ParaWW}   one has 
\begin{align}
\pa_t \vect{\eta}{\underline{\wahlen}} & =  \cW^{-1} \begin{pmatrix} 0 & |D| \\ -1 & \gamma \Hilb \end{pmatrix} \cW \vect{\eta}{\underline{\wahlen}} \label{ilprimo} \\
& \ \ + 
\cW^{-1}  \Opbw{\begin{bmatrix} 
- \im  \xi\#_{\vr_1} \sfV_\gamma   & 0 
\\
-\sfa  &  -\sfV_\gamma\#_{\vr_1}\im  \x
\end{bmatrix}}\cW\vect{\eta}{\underline{\wahlen}} \label{ilprimo2} 
\\
& \ \ + \cW^{-1} 
\Opbw { \begin{bmatrix} 0 & b^{(-1)}_{\geq 1}   \\ 0
  &   \gamma \frac{1}{\ii \xi}\#_{\vr_1} b_{\geq 1}^{(-1)} \end{bmatrix}}
\cW\vect{\eta}{\underline{\wahlen}}+ \bR(\eta,\psi) \vect{\eta}{\underline{\wahlen}}  \label{ilprimo3}
\end{align}
where $b^{(-1)}_{\geq 1}   $ is given in Lemma \ref{laprimapara} and $\bR(\eta,\psi) $ is  a matrix of real smoothing operators in $ \Sigma \cR^{-\vr_1+1}_{1}[r, N]$.

We now compute the above conjugated operators 
applying the transformation rule
\begin{equation}\label{generic}
\cW^{-1}\begin{pmatrix} \cA& \cB\\ \cC& \cD\end{pmatrix}\cW= \begin{pmatrix} \cA+ \frac{\gamma}{2} \cB \pa_x^{-1} & \cB \\ \cC - \frac{\gamma}{2} \pa_x^{-1} \cA - \frac{\gamma^2}{4} \pa_x^{-1} \cB\pa_x^{-1} + \frac{\gamma}{2} \cD  \pa_x^{-1}& \cD - \frac{\gamma}{2} \pa_x^{-1}\cB\end{pmatrix}. 
\end{equation}
The operator at the right hand side of \eqref{ilprimo} is the Fourier multiplier given in \eqref{eq:lin}. 
Then, by  Proposition \ref{teoremadicomposizione} one has
\begin{align}  
\eqref{ilprimo2} 
 = &\Opbw{\begin{bmatrix}
    - \im  \xi\#_{\vr_1} \sfV_\gamma & 0 
     \\ -\sfa  & - \sfV_\gamma \#_{\vr_1}\im \xi
    \end{bmatrix}}+ \bR(\eta, \psi)
\, ,\label{conj2}
\end{align}
where  
$\bR (\eta,\psi) $  is a matrix of real smoothing operators in 
$ \Sigma \cR^{-\vr_1}_{1}[r, N]$. 
Finally, again by  Proposition \ref{teoremadicomposizione} and using \eqref{comp3fin}, \eqref{comp3ordine}, we deduce that 
\be  
\eqref{ilprimo3} = \Opbw {\begin{bmatrix}  \frac{\gamma}{2}b_{\geq 1}^{(-1)} \#_{\vr_1} \frac{1}{\ii \xi}&  b_{\geq 1}^{(-1)}  \\[0.5em]
\frac{\gamma^2}{4} \frac{1}{\ii \xi} \#_{\vr_1} b_{\geq 1}^{(-1)} \#_{\vr_1} \frac{1}{\ii \xi}&  
\frac{\gamma}{2} \frac{1}{\ii \xi} \#_{\vr_1}b_{\geq 1}^{(-1)}
\end{bmatrix}}+ \bR(\eta,\psi) \, ,  \label{conj3}
\ee
where 
$\bR(\eta,\psi) $ is a matrix of real smoothing operators in $  \Sigma \cR^{-\vr_1}_{1}[r, N]$. 

In conclusion, by \eqref{eq:lin}, \eqref{conj2}, \eqref{conj3}, 
system  \eqref{ilprimo}-\eqref{ilprimo3} has the form  \eqref{WWW}.
Finally, the matrix in \eqref{WWW} satisfies \eqref{real:ham:sym0}--\eqref{real:ham:sym} thanks to \eqref{prop:ov3}--\eqref{prop:ov}.
\end{proof}

\paragraph{Complex coordinates.}
We now introduce complex coordinates that diagonalize the linear part of  system \eqref{WWW}, by defining the variable
\be
U=\begin{pmatrix}{ u}\\{ {\bar u} }\end{pmatrix}:= \cM^{-1} \begin{pmatrix}{ \eta}\\{ {\underline{\wahlen}} } \end{pmatrix}
\stackrel{\eqref{BISMA_wahlen}, \eqref{zak}}{=} \Lmap^{-1}\begin{pmatrix}{ \eta}\\{ {\upomega} } \end{pmatrix}
\,, \qquad
\cM \mbox{ in } \eqref{defMD} \,, \quad \Lmap \mbox{ in } \eqref{app_L}  \,. 
\label{Ugrande}
\ee
Explicitly
\be\label{u.def}
u = \frac{1}{\sqrt{2}} \left(M(D) \eta + \im M^{-1}(D) \upomega - \frac{\gamma}{2}\im M^{-1}(D) \pa_x^{-1} \eta  \right)  \ , \quad M(D) \in \wt \Gamma^{-\frac14}_0 \mbox{ in } \eqref{defMD} \, ,
\ee
from which one  immediately deduces the  equivalence of norms: 
for some $C>1$ and  any   $\s \in \R$
 \be\label{eq:etaomegaUequiv}
C^{-1} \norm{U}_{\s+ \frac34} \leq \norm{\eta}_{\s+\frac12} + \norm{\upomega}_{\s+1} \leq C  \norm{U}_{\s +\frac34}\,.
 \ee
We prove the following result.
\begin{lemma}\label{LemCompl}
{\bf (Hamiltonian formulation of the water waves in complex coordinates)}
Let  $N \in \N $ and $ \varrho  \geq 0$. 
There exist $ s_0> \frac32$, $r > 0 $  such that, 
if $ (\eta, \psi) \in B_{s_0, \R} (I;r) $ is a solution of \eqref{eq:etapsi}, and $\zak=\zak(t)$ is defined in \eqref{zak}, then $ U$ defined in \eqref{Ugrande} solves 
\be \label{complexo}
\begin{aligned}
\pa_tU  =&  \begin{pmatrix}
    \frac{\gamma}{2}\Hilb  & 0 \\
    0 & \frac{\gamma}{2}\Hilb 
\end{pmatrix} U 
 +  \Opbw{\begin{bmatrix} -\ii \sfV_\gamma  \xi & -\frac{1}{4}(\sfV_\gamma)_x\\[0.3em]
-\frac{1}{4}(\sfV_\gamma)_x  & -\ii \sfV_\gamma  \xi 
 \end{bmatrix}+\begin{bmatrix} -\ii (1+ \frac{\sfa}{2}) \omega(\xi) & -\ii \frac{\sfa}{2}\omega(\xi)\\[0.3em]
\ii \frac{\sfa}{2}\omega(\xi) & \ii (1+ \frac{\sfa}{2}) \omega(\xi)
 \end{bmatrix}}U   \\
 & +\Opbw{\widetilde \bA_{\geq 1}^{(-\frac12)}(\zak)}U
 + \Opbw{\mathbf{B}_{\geq N}^{(-\frac 1 2)}(\zak)}U
  +\bR_{\geq 1}(\zak)U 
\end{aligned}
\ee 
where $\Hilb$ is the Hilbert transform in \eqref{def:Hilbert}, $\sfV_\gamma\in \Sigma {\mathcal F}^\R_{1}[r, N]$ is the real-valued function in \eqref{def:Vgamma}, $\omega(\xi)$ is the Fourier multiplier defined in \eqref{omegaxi}, $\sfa \in \Sigma \cF^\R_{1}[r, N] $ is the Taylor coefficient in \eqref{a.taylor},
$\zak= \Lmap^{-1} \vect{\eta}{\psi}$ is   the Zakharov variable defined in \eqref{zak} 
and  
\begin{itemize}

  \item $\tilde \bA_{\geq 1}^{(-\frac12)}(\zak)$ is a matrix of pluri-homogeneous symbols in $\Sigma_1^{N-1} \wt \Gamma^{-\frac 1 2}$
  (recall item \ref{Hhom} of Definition \ref{def:sfr}) which satisfies the linear Hamiltonian property \eqref{lin:ham_complex};
 \item $\mathbf{B}^{(-\frac 12)}_{\geq N}(\zak)$ is a real-to-real matrix of symbols in $\Gamma_{\geq N}^{-\frac 1 2}[r]$;
 \item  $\bR_{\geq 1}(\zak) $ is a  real-to-real matrix smoothing operators in $ \Sigma \cR^{-\vr}_{1}[r, N]$.
\end{itemize}
\end{lemma}
\begin{remark}\label{rem:sost_symbols}
  In formula \eqref{complexo}, it is convenient to keep the transport operator $\sfV_\gamma$ and the Taylor coefficient $\sfa$ expressed in the variables $(\eta,\psi)$, rather than in the complex variable $\zak$, unlike all other symbols and operators. Clearly, whenever needed, they can be rewritten as functions of $\zak$ by simply composing with $\Lmap$, namely $\sfV_\gamma(\Lmap\zak)$ and $\sfa(\Lmap\zak)$.
\end{remark}
\begin{proof}
Recall that $\cM,\ \cM^{-1},\ \Lmap,\ \Lmap^{-1}$ are matrices of Fourier multipliers with symbols in $\wt\Gamma_0^{\frac14}$ (cf. \eqref{defMD} and \eqref{app_L}).
We apply \Cref{wahlen} with a parameter $\vr_1 \ge 0$, to be specified at the end of the proof in terms of $\vr \ge 0$. 
Since $ (\eta, \underline{\wahlen}) $ solves \eqref{WWW}, the complex variable $ U$ in \eqref{Ugrande} solves 
\begin{align}
\pa_t U & = \cM^{-1} \begin{pmatrix}
\frac{ \gamma}{2} \Hilb & |D| \\ 
-(1 + \frac{\gamma^2}{4}|D|^{-1}  ) & \frac{ \gamma}{2} \Hilb \end{pmatrix} \cM U
\label{coniucomplessa}  \\ 
&\ \ +  \cM^{-1} 
\Opbw{\begin{bmatrix}
- \im   \xi\#_{\vr_1} \sfV_\gamma   
& 0  
\\ 
-\sfa
& - \sfV_\gamma\#_{\vr_1}\im  \xi
\end{bmatrix}} 
\cM U   \label{coniucomplessa1} \\
&\ \ + \cM^{-1} 
\Opbw{\begin{bmatrix}  \frac{\gamma}{2}b_{\geq 1}^{(-1)} \#_{\vr_1} \frac{1}{\ii \xi}&  b_{\geq 1}^{(-1)}  \\[0.5em]
\frac{\gamma^2}{4} \frac{1}{\ii \xi} \#_{\vr_1} b_{\geq 1}^{(-1)} \#_{\vr_1} \frac{1}{\ii \xi}&  
\frac{\gamma}{2} \frac{1}{\ii \xi} \#_{\vr_1}b_{\geq 1}^{(-1)}
\end{bmatrix}
} 
\cM U + \bR_1( \zak)U \label{coniucomplessa2}
\end{align}
where  $\bR_1(\zak)  := \cM^{-1}\bR(\Lmap\zak) \cM $ is a real-to-real matrix of smoothing operators in $ \Sigma \cR^{-\vr_1+\frac74}_{1}[r, N]$.
\\
The operator in the right hand side of \eqref{coniucomplessa} is computed in 
\eqref{diaglin}-\eqref{relation}, {as it gives the dispersion relation.}  

In order to compute the  conjugated operators  in \eqref{coniucomplessa1}-\eqref{coniucomplessa2}, we apply the following transformation rule,   
where we denote  by $M:= M(D)$ the Fourier multiplier in \eqref{defMD}:  
\begin{align}\label{generic2} 
&\scalebox{0.8}{$
\cM^{-1} \begin{pmatrix} \cA_1& \cA_2 \\ \cA_3 & \cA_4\end{pmatrix} \cM 
$}\\
&\stackrel{\eqref{zak}} = 
\frac12
\scalebox{0.8}{$ \begin{pmatrix} M\cA_1 M^{-1} + M^{-1} \cA_4 M + \ii M^{-1} \cA_3 M^{-1} - \ii M  \cA_2 M & M \cA_1 M^{-1}  - M^{-1} \cA_4 M + \ii M^{-1}  \cA_3 M^{-1} +\ii M \cA_2 M \\ 
M \cA_1 M^{-1}  - M^{-1} \cA_4 M - \ii M^{-1}  \cA_3 M^{-1} - \ii M  \cA_2 M &  M  \cA_1 M^{-1}  + M^{-1}  \cA_4 M  - \ii M^{-1}  \cA_3 M^{-1} +  \ii M  \cA_2 M \notag \end{pmatrix}
$}\,.
\end{align}
Therefore, 
using  \eqref{comp3fin}--\eqref{sharp3},  we have that
\begin{align*}
\eqref{coniucomplessa1}  &= 
 \frac{1}{2} \vOpbw{
      m_{\geq 1}^{(1)}- \ii a^{(\frac12)}_{\geq 1}}+
       \frac{1}{2} \zOpbw{
    m_{\geq1}^{(0)}- \ii a^{(\frac12)}_{\geq 1}} + \bR(\zak)
\end{align*}
where 
\begin{align}
m_{\geq 1}^{(1)}&:= M \#_{\vr_1} \Big( (-\ii \xi)\#_{\vr_1} \sfV_\gamma \Big)\#_{\vr_1} M^{-1}  + M^{-1}  \#_{\vr_1}  \Big(
\sfV_\gamma\#_{\vr_1} (-\ii  \xi)\Big)\#_{\vr_1} M \,, \label{uff.1}
\\
m_{\geq 1}^{(0)}&:= M   \#_{\vr_1} \Big((-\ii \xi)\#_{\vr_1} \sfV_\gamma \Big)\#_{\vr_1}  M^{-1}  - M^{-1}  \#_{\vr_1} \Big(\sfV_\gamma\#_{\vr_1} (-\ii  \xi) \Big) \#_{\vr_1} M\,, \label{uff.2}\\
a^{(\frac12)}_{\geq 1}&:= M^{-1}\#_{\vr_1} \sfa\#_{\vr_1} M^{-1} \label{def:a12}
\end{align}
and 
  $\bR(\zak) $  is a real-to-real matrix of smoothing operators in 
$ \Sigma \cR^{-\vr_1 + \frac54
}_{1}[r, N] $, obtained by expressing the smoothing remainder coming from \eqref{comp3fin} in terms of the complex variable $\zak$, through the change of variables \eqref{zak}.
By \eqref{prop:ov}--\eqref{prop:ov3}, the symbols in \eqref{uff.1}, \eqref{uff.2} and \eqref{def:a12} satisfy the linearly Hamiltonian properties (see \eqref{lin:ham_complex}):
$$
\ov{m_{\geq 1}^{(1)}} = - m_{\geq 1}^{(1)}\,,  \quad (m_{\geq 1}^{(0)})^{\vee} = m_{\geq 1}^{(0)}\, , \quad \ov{a^{(\frac12)}_{\geq 1}}=a^{(\frac12)}_{\geq 1}\, , \quad \big(a^{(\frac12)}_{\geq 1}\big)^\vee=a^{(\frac12)}_{\geq 1}\,.
$$
We further expand $m_{\geq 1}^{(1)}$,  $m_{\geq 1}^{(0)}$ and $a^{(\frac12)}_{\geq 1}$ using \eqref{espansione2}, \eqref{sharp3} and \eqref{magicordine} getting
\begin{equation}\label{nonsos2}
m_{\geq 1}^{(1)}= -2\ii\sfV_\gamma\, \xi + d_{\geq1}^{(-1)}\,, \quad m_{\geq 1}^{(0)} = - \frac{1}{2}(\sfV_\gamma)_x + r_{\geq1}^{(-1)}\, , \quad a^{(\frac12)}_{\geq 1}= \sfa M^{-2} + a^{(-\frac32)}_{\geq 1}=\sfa \omega(\xi) + a^{(-\frac12)}_{\geq 1},
\end{equation}
with
$d_{\geq1}^{(-1)}, r_{\geq1}^{(-1)} \in  \Sigma \Gamma^{-1}_{1}[r, N]$, $a^{(-\frac12)}_{\geq 1}\in \Sigma \Gamma^{-\frac12}_{1}[r, N]$ and $a^{(-\frac32)}_{\geq 1}\in \Sigma \Gamma^{-\frac32}_{1}[r, N]$, where to compute the second one we also used that 
$ (\sfV_\gamma)_x \xi (M \pa_\xi M^{-1} - M^{-1} \pa_\xi M)= \frac{1}{2}(\sfV_\gamma)_x $
 plus a symbol in $\Sigma \Gamma^{-2}_{1}[r, N]$ and that 
$$
M^{-2}(\xi) -\omega(\xi)=  -  \frac{\gamma^2}{4}  M^{-2}(\xi)|\xi|^{-1} \in \wt \Gamma_0^{-\frac12}\,.
$$
 By difference, $ d_{\geq 1}^{(-1)}$, $r_{\geq 1}^{(-1)} $ and $a^{(-\frac12)}_{\geq 1}$ satisfy the same algebraic properties 
 \eqref{nonsos2}
 as \( m_{\geq 1}^{(-1)} \),  \( m_{\geq 1}^{(0)} \) and $a^{(\frac12)}_{\geq 1}$, respectively.
Reasoning as above and using also \eqref{real:ham:sym}, we get
\begin{equation}
  \eqref{coniucomplessa2}  =  \frac{1}{2} \vOpbw{
      d^{(-\frac32)}_{\geq 1} + d_{\geq N}^{(-\frac 32)}}+
       \frac{1}{2} \zOpbw{
   r^{(-\frac32)}_{\geq 1} + r_{\geq N}^{(-\frac 32)}} + \bR(\zak)\,.\label{nonsos3}
\end{equation}
Here $ d^{(-\frac32)}_{\geq 1} = \sum_{p=1}^{N-1}d^{(-\frac32)}_{p} $, $r^{(-\frac32)}_{\geq 1}=\sum_{p=1}^{N-1}r^{(-\frac32)}_{p} $ are the
pluri-homogeneous symbols in $\Sigma_1^{N-1}\wt \Gamma^{-\frac32}$ defined as
\small
\be\label{dprp}
\begin{aligned}
    d^{(-\frac32)}_{p} &:= \frac{\gamma}{2} M\#_{\vr_1}\left( b_{p}^{(-1)}\#_{\vr_1} \frac{1}{\ii \xi}\right)\#_{\vr_1} M^{-1}+\frac{\gamma}{2}M^{-1} \#_{\vr_1} \left(\frac{1}{\ii \xi}\#_{\vr_1} b_{p}^{(-1)} \right)\#_{\vr_1} M\\
    & +\ii \frac{\gamma^2}{4}M^{-1}\#_{\vr_1} \left(\frac{1}{\ii \xi}\#_{\vr_1} b_{p}^{(-1)}\#_{\vr_1} \frac{1}{\ii \xi}\right)\#_{\vr_1} M^{-1}-\ii M\#_{\vr_1} b_{p}^{(-1)}\#_{\vr_1} M\\
    r_{p}^{(-\frac32)} &:= \frac{\gamma}{2} M\#_{\vr_1} \left(b_{p}^{(-1)}\#_{\vr_1} \frac{1}{\ii \xi}\right)\#_{\vr_1} M^{-1}-\frac{\gamma}{2}M^{-1} \#_{\vr_1} \left(\frac{1}{\ii \xi}\#_{\vr_1} b_{p}^{(-1)} \right)\#_{\vr_1} M\\
    & +\ii \frac{\gamma^2}{4}M^{-1}\#_{\vr_1} \left(\frac{1}{\ii \xi}\#_{\vr_1} b_{p}^{(-1)}\#_{\vr_1} \frac{1}{\ii \xi}\right)\#_{\vr_1} M^{-1}+\ii M\#_{\vr_1} b_{p}^{(-1)}\#_{\vr_1} M
\end{aligned}
\ee
\normalsize
 which satisfy, by \eqref{prop:ov3}, \eqref{real:ham:sym} and \eqref{real:ham:sym0},  the linearly Hamiltonian property (up to homogeneity $N-1$): for any $p=1, \dots, N-1$,
\begin{align*}
    \ov{d^{(-\frac32)}_{p}}=-d^{(-\frac32)}_{p}, \qquad \big(r^{(-\frac32)}_{p}\big)^\vee =r^{(-\frac32)}_{p}\,,
\end{align*}
and $d_{\geq N}^{(-\frac 32)}, r_{\geq N}^{(-\frac 32)}$ are symbols in $\Gamma_{\geq N}^{-\frac 3 2}[r]$,  defined as in \eqref{dprp} but with  $b_{p}^{(-1)} \leadsto  b_{\geq N}^{(-1)} $.
Moreover $\bR(\zak)$  is a real-to-real matrix of smoothing operators in
$ \Sigma \cR^{-{\vr_1}}_{1}[r, N] $.

In conclusion, by  \eqref{diaglin}-\eqref{relation}, \eqref{nonsos2}, \eqref{nonsos3}, 
computing the symbols at $ (\eta, \psi) = \Lmap \zak $ (see \eqref{zak}), 
and 
choosing $\vr_1:= \vr+ \frac74$,
we deduce that system 
\eqref{coniucomplessa}-\eqref{coniucomplessa2} has the form 
 \eqref{complexo} with 
 $\wt\bA_{\geq 1}^{(-\frac12)}(\zak) := \frac{1}{2} \vOpbw{
      d^{(-\frac32)}_{\geq 1} }+
       \frac{1}{2} \zOpbw{r^{(-\frac32)}_{\geq 1} } $
       and
       $
      \bB_{\geq N}^{(-\frac12)}(\zak):=  \frac{1}{2} \vOpbw{
       d_{\geq N}^{(-\frac 32)}}+
       \frac{1}{2} \zOpbw{r_{\geq N}^{(-\frac 32)}}$, and $\bR_{\geq 1}(\zak)$ is the sum of all the real-to-real matrices of smoothing operators.
\end{proof}

In the next lemma we prove that the map sending $\zak \to U$ has the form of an  admissible transformation.
\begin{lemma}\label{complex:good}
   Let $\wt \vr>1$. The complex variable $U$ in \eqref{Ugrande} reads
   \be\label{eq:complex_good}
U=  \pmb{\cG}_{\C}(\zak) \zak, \quad  \pmb{\cG}_{\C}(\zak):= \Lmap^{-1}\cG(\Lmap \zak)\Lmap,
   \ee
   where the real map $\cG$ is the good unknown of Alinhac defined in \eqref{GU} and $\Lmap^{-1}$ the linear operator of \eqref{zak}, and  $\pmb{\cG}_{\C}(\zak)$ is an admissible transformation of order $(\frac12, \frac12)$  and gain $\wt \vr $ according to \Cref{admtra}.
\end{lemma}
\begin{remark}
    The previous lemma highlights a loss of $\frac{1}{2}$ derivatives when passing from the Hamiltonian variable $\zak$ to Alinhac's good unknown. This phenomenon is well known in the theory of gravity waves. However, \Cref{complex:good,loc.inv} establish that the complex transformation $U = \pmb{\cG}_{\C}(\zak) \zak$—and consequently, Alinhac's good unknown itself—is \emph{nonlinearly} invertible.  
\end{remark}
\begin{proof}
Identity \eqref{eq:complex_good} follows by 
 \eqref{Ugrande}, \eqref{BISMA_wahlen}, \eqref{GU}, \eqref{zak}.
 We prove that $\pmb{\cG}_{\C}(\zak)$ fulfills the properties of admissible transformations. 

\noindent{\sc Gain \& linear invertibility:} computing $\pmb{\cG}_{\C}(\zak) $ explicitly one gets
\begin{gather*}
 \pmb{\cG}_\C(\zak)= \Id {-} \frac12 \begin{pmatrix}
    \ii & \ii \\ -\ii & -\ii 
\end{pmatrix}M^{-1}(D)\Opbw{\sfB(\Lmap\zak; x)}M^{-1}(D),\\  \pmb{\cG}_\C^{-1}(\zak)= \Id
{+} \frac12 \begin{pmatrix}
    \ii & \ii \\ -\ii & -\ii 
\end{pmatrix}M^{-1}(D)\Opbw{\sfB(\Lmap\zak; x)}M^{-1}(D)  
\end{gather*}
where $\sfB=\sfB(\eta,\psi)$ is defined in \eqref{form-of-B}  and $ M(D)$ is the  Fourier multiplier  in  \eqref{defMD}. 
Then estimate \eqref{lin.est.F}, with $\nu=\frac12$, follows from \eqref{stimapar} and the fact that $ M^{-1}(D)$ is of order $\frac14$. 

\noindent{\sc Expansion:} Using \eqref{esp.b.gamma.1.2} we expand
$
 \pmb{\cG}_{\C}(\zak)= \Id+ \pmb{\cG}_1(\zak)+\pmb{\cG}_2(\zak)+ \pmb{\cG}_{\geq 3}(\zak)$ with 
\begin{align*}
 & \pmb{\cG}_j(\zak):= \frac12 \begin{pmatrix}
    \ii & \ii \\ -\ii & -\ii 
\end{pmatrix}M^{-1}(D)\Opbw{\sfB_j(\Lmap\zak;x)}M^{-1}(D), \quad j = 1,2 \\  
& \pmb{\cG}_{\geq 3}(\zak):=\frac12 \begin{pmatrix}
    \ii & \ii \\ -\ii & -\ii 
\end{pmatrix}M^{-1}(D)\Opbw{\sfB_{\geq 3}(\Lmap\zak;x)}M^{-1}(D)
\end{align*}
where  $\sfB_1(\eta,\psi):= |D|\psi$, $ \sfB_2(\eta,\psi):=  -\eta\psi_{xx}-|D|\eta|D|\psi$ and, by \eqref{BVa.ana}, the map $B_{\s_0, \R}(r)\ni \zak \mapsto \sfB_{\geq 3}(\Lmap\zak)\in H^{\s_0 -\frac54}(\T;\C)$ is analytic.
Since $\cM^{-1}$ is a matrix of Fourier multipliers in $\wt \Gamma_0^{\frac14}$, the maps  $  \pmb{\cG}_1(\zak) \in \wt \mM_1^{\frac12}$, $ \pmb{\cG}_2(\zak) \in \wt \mM_2^{\frac12}$
and
$ \pmb{\cG}_{\geq 3}(\zak) \in  \mM_{\geq 3}^{\frac12}[r]$.

\noindent{\sc Differentiability:} By analyticity \eqref{BVa.ana}, 
 the map $B_{\s, \R}(r) \ni \zak \mapsto \sfB(\Lmap\zak) \in H^{\s -\frac54}(\T;\R)$ is differentiable  for any $\s \geq \s_0$.
 Moreover,  by \Cref{thm:contS}, for any $s \in \R$ the linear map $L^\infty(\T;\C) \ni f \mapsto \Opbw{f} \in \cL(H^s(\T;\C))$ is bounded hence analytic. 
Hence, provided  $s_0-\frac54 > \frac12$, the    composition 
 $B_{s_0, \R}(r) \ni  \zak \mapsto  \Opbw{\sfB(\Lmap\zak;x)} \in \cL(H^s(\T;\C))$ is differentiable, 
so is $B_{s_0, \R}(r)\ni
\zak \mapsto  \pmb{\cG}_\C(\zak)
\in  \mathbb{X}^{s,\frac12}$. 
 Thus, given an arbitrary  $\wt \rho \geq 1$ fixed, the differentiability \eqref{ad:diff} follows.

Finally we prove \eqref{stima.d.adm}. 
 By Cauchy estimates, for any $ \| \zak\|_{s_0}\leq \frac23 r$,   one has
\be\label{dBgeq3.est}
\| \di {\sfB}_{\geq 3} (\Lmap\zak)\|_{\cL(\dot H^{s_0};H^{s_0 -{5}\!/\!{4}})} \lesssim \sup_{\|\wt \zak- \zak \|_{s_0}< \| \zak\|_{s_0}/2}\frac{ \| {\sfB}_{\geq 3} (\Lmap \wt{\zak})\|_{s_0-\frac54}}{\| \zak\|_{s_0}/2}   \lesssim \|\zak\|^{2}_{s_0}\,,
\ee
yielding, by \Cref{thm:contS}, 
\begin{align*}
\norm{\di\pmb{\cG}_{\geq 3}(\zak)[\wh \zak]V}_s & \lesssim 
\norm{ \Opbw{\di \sfB_{\geq 3}(\Lmap\zak;x)[\Lmap\wh \zak] }}_{\mathbb{X}^{s+\frac14,0}}\norm{V}_{s+\frac12}
\lesssim 
|\di {\sfB}_{\geq 3} (\Lmap\zak)[\Lmap\hat \zak]|_{0, L^\infty, 4} \norm{V}_{s+\frac12}\\
& \lesssim 
\norm{\di {\sfB}_{\geq 3} (\Lmap\zak)[\Lmap\hat \zak]}_{s_0-\frac54} \norm{V}_{s+\frac12}
\stackrel{\eqref{dBgeq3.est}}{\lesssim}\norm{\zak}_{s_0}^2 \norm{V}_{s+\frac12}
\end{align*}
 after recalling that ${\sfB}_{\geq 3} (\Lmap\zak)$ does not depend on $\xi$.  
\end{proof}

\subsection{Quadratic paradifferential normal form}\label{subsec:quadraticNF}
The goal of this section is to put the vector field in \eqref{complexo} to its quadratic normal form, for any value of the vorticity $\gamma$. 
This is done by using paradifferential techniques to put in normal form the symbols, and one step of Poincar\'e-Birkhoff normal form to remove the quadratic components of the smoothing vector field. 
This is possible since, for any value of the vorticity $\gamma$,  there are no 3-waves interactions, see Lemma \ref{3onde}.
More precisely, the aim of this section is the following normal form result.
\begin{theorem}[{\bf Quadratic paradifferential normal form}]\label{thm:quadratic.nf}
Let $\gamma \in \R$ arbitrary and $N\in \N$, $N\geq 2$.  Let $\vr >{2N+\frac92}$. 
There exist $ s_0> \frac32$, $r > 0 $  and a $(0,{\frac{9}{2}})$- admissible transformation ${\bf T}(\zak)$
 with gain {$\vr-2N$} (recall Definition \ref{admtra}) 
    such that if $(\eta, \psi)\in B_{s_0,\R}(I;r)$ solves \eqref{eq:etapsi}, and $\zak=\zak(t)$ is defined in \eqref{zak}, then the variable $Y={\bf T}(\zak)U$, with $U$ in \eqref{Ugrande}, solves 
\begin{equation}\label{BNF1}
 \begin{aligned}
 \pa_{t}Y&= -\ii\bOmega(D)Y+ \vOpbw {\ii \sfD_{\geq 2}(\zak)}Y  +  \bB_{\geq N}(\zak) Y+  \bR_{\geq 2}(\zak)Y
 \end{aligned}
 \end{equation}
 where $\bOmega(D)$ is defined in \eqref{diaglin}-\eqref{omegaxi},  
 $\sfD_{\geq 2}$ is the real-valued symbol in $\Sigma\Gamma^{1}_{ 2}[r, N]$ given by
     \begin{equation}\label{def:D}
     \sfD_{\geq 2}(\zak; x,\x):= \big(\check{\sfV}_{2}(\zak;x) + \check{\sfV}_{\geq 3}(\zak; x)\big) \xi +  \tq_{\geq 2}(\zak; x) \omega(\xi) +  \tb_{\geq 2}(\zak;x)\,\sign{(\xi)}+ \tilde \tc_{\geq 2}^{(-\frac12)}(\zak; x, \xi)
     \end{equation}
     with $\omega(\xi)$ defined in \eqref{omegaxi} and\\
     $\bullet$ 
 $\check{\sfV}_{2}(\zak;x)$ is the real-valued, $2$-homogeneous function in $ \wt\cF^{\R}_{2}$, defined by
 \begin{equation}\label{V2.exp0}
    \begin{aligned}
        \check{\sfV}_{2} (\zak;x) & := \sum_{\substack{k_1, k_2 \in \Z_* \\ \sigma_1, \sigma_2 \in \{\pm 1\} }} V_{k_1, k_2}^{\sigma_1, \sigma_2}(\gamma) \,  \zetina_{k_1}^{\sigma_1}\zetina_{k_2}^{\sigma_2} e^{\im (k_1 \sigma_1 + k_2 \sigma_2) x}\,,
        \end{aligned}
        \end{equation}
        with coefficients given by
        \begin{equation}
\begin{aligned}\label{V2.exp1}
& V_{k_1, k_2}^{\sigma_1, \sigma_2}(\gamma):= 
    -\gamma \sigma_2 k_2 \left( \operatorname{sign}(\sigma_1 k_1 + \sigma_2 k_2) + \frac{\gamma^2 {(1-\delta(\sigma_1 k_1 + \sigma_2 k_2))}}{(\sigma_1 k_1 + \sigma_2 k_2)} + \sigma_1\operatorname{sign}( k_1) + \frac{\gamma^2 \sigma_1 }{ k_1}\right) \fm_{k_1} \fm_{k_2} \\
    &   + \sigma_1\left(  \sigma_1 k_1 (|\sigma_1 k_1 + \sigma_2 k_2| - |k_1| + \sigma_1 \sigma_2 k_1 \operatorname{sign}(k_2))+\sigma_2 \gamma^2 \frac{k_1^2+k_2^2}{k_2}
        + \gamma^2 \sigma_1 k_1 {\big( 1-\delta(\sigma_1 k_1+ \sigma_2 k_2) \big)} 
        \right)   \fe_{k_1} \fm_{k_2}\\
        & \ \ -\frac{\gamma}{2}\sigma_1 \sigma_2   \left( \sigma_1 \sigma_2  k_1 k_2 + | k_1 k_2| + 2 k_2^2\right)  \fe_{ k_1} \fe_{k_2}\, . 
\end{aligned}
\end{equation}
Here $\delta(j)$ denotes  the Kronecker delta, defined as
 \be\label{delta}
\delta(j):= 1  \mbox{\ \ if \ } j=0 \quad 
\mbox{ and } \quad \delta(j) := 0  \mbox{\ \ otherwise} \, ,
\ee
the numbers $\fm_j$ in \eqref{al.be.m}, while  
\be  \label{numerelli}
\fe_j:= \frac12\fm_j^{-1} + \gamma \frac{\fm_j}{2j} = 
\frac{1}{\sqrt{2}} \left(\frac{ 4|j| + \gamma^2}{4 j^2}\right)^{\frac14}  + \frac{1}{\sqrt{2}}\frac{\gamma}{2j}
\left( \frac{4 j^2}{4 |j| + \gamma^2}\right)^{\frac14} \ , \quad j \in \Z_*\,.
\ee 
$\bullet$ $\check{\sfV}_{\geq 3}(\zak;x)$ is a real-valued function in $\Sigma \cF^\R_{3}[r, N]$;
\\
$\bullet$ $\tq_{\geq 2}(\zak;x)$ and  $\tb_{\geq 2}(\zak;x)$ are real-valued functions in  $ \Sigma  \cF^\R_2[r, N]$ and  $\tilde \tc^{(-\frac12)}_{\geq 2}(\zak;x, \xi)$ is a pluri-homogeneous real-valued symbol in  $  \Sigma_2^{N-1} \wt \Gamma^{-\frac 1 2}$;
\\
$\bullet$ $\bB_{\geq N}(\zak)$ is a real-to-real matrix of spectrally localized,  bounded maps in $\cS_{\geq N}^0[r]$ and $ \bR_{\geq 2}(\zak)$ is a real-to-real matrix of smoothing  operators  
    in  $\Sigma\mathcal{R}^{-\vr+{2N+\frac32}}_{2}[r, N] $ (see \Cref{def:specloc}).
\end{theorem}

\begin{remark}
As anticipated in \Cref{rem:not.symm},  the coefficients $ V_{k_1, k_2}^{\sigma_1, \sigma_2}(\gamma)$ in 
 \eqref{V2.exp1} are 
  not symmetric; one can symmetrize them as described in the aforementioned remark.
\end{remark}

The rest of the  section is devoted to the proof of this theorem, which is finally achieved at pag. \pageref{proof.thm:quadratic.nf}.

\subsubsection{Block-Diagonalization at highest order}
In this section we  diagonalize the paradifferential operator of order $\frac12$ in \eqref{complexo}, that we denote with   
$\Opbw{\mathbf{A}_{\geq 0}^{(\frac 12)}}$, where
\begin{equation}
 \bA_{\geq 0}^{(\frac12)}:=  \ii \begin{bmatrix}
- (1 + \frac{\sfa(x)}{2}) &  -\frac{\sfa(x)}{2}  \\
\frac{\sfa(x)}{2}  & 1 + \frac{\sfa(x)}{2} 
\end{bmatrix} \omega(\xi)\,,  
 \quad \sfa(x) \mbox{ in } \eqref{a.taylor} \, \quad \omega(\xi) \mbox{ in } \eqref{omegaxi} \,.
\label{diag1}
\end{equation}
The eigenvalues of the matrix $\bA_{\geq 0}^{(\frac12)}$
are 
 $\,\pm \im \lambda( \zak; x)\omega(\xi)$, where 
\begin{equation}
\label{diag2}
\lambda( \zak; x) := \sqrt{\Big(1+\frac{\sfa(x)}{2}\Big)^2 -\Big(\frac{\sfa(x)}{2}\Big)^2}  = \sqrt{1 +\sfa(x)} \,,
\quad \lambda-1\in \Sigma\cF_1^{\R}[r, N] \, . 
\end{equation}
A matrix which diagonalizes $\mathbf{A}^{(\frac 12)}_{\geq 0}$
is
 \begin{equation} \label{diag4}
 \bH (\zak; x)  : = 
 \begin{pmatrix}
 h (\zak; x )& g(\zak; x) \\
 g(\zak;x ) & h(\zak; x ) 
 \end{pmatrix} \,, \quad
h :=  \frac{1+\frac{\sfa}{2}+ \lambda}{\sqrt{\big(1+\frac{\sfa}{2} + \lambda\big)^2 - \big( \frac{\sfa}{2}\big)^2}}  \, , \quad 
g:=   \frac{-\frac{\sfa}{2}}{\sqrt{\big(1+\frac{\sfa}{2}+ \lambda \big)^2 - \big( \frac{\sfa}{2}\big)^2}}  \,. 
 \end{equation}
 Note that both $h$ and $g$ satisfy the following identities:
\begin{align}
    (h^2+g^2) \, (1+ \frac{\sfa}{2})+hg\sfa=\lambda, \qquad 2 h g (1+ \frac{\sfa}{2}) 
  +(h^2 + g^2)   \frac{\sfa}{2}   = 0\ \,, \quad  h^2 - g^2 = \det \bH = 1 \, .
  \label{canc:diag}
\end{align}
Furthermore, $h-1$ and $g$ are symbols in $\Sigma \cF^{\R}_{1}[r, N]$ and the matrix $  \bH(\zak; x) $ fulfills 
 $\bH(\zak;x) - \mathds{1}\in \Sigma \cF^{\R}_{1}[r, N]$. In addition,
 $\bH(\zak;x)$ is  invertible with inverse
 \begin{equation}
 \label{diag5}
 \bH(\zak; x )^{-1} : = 
 \begin{pmatrix}
h (\zak; x )& - g(\zak; x) \\
-g(\zak; x ) & h(\zak; x )
 \end{pmatrix} \, 
 \end{equation}  
 and, using \eqref{canc:diag}, fulfills
 \begin{equation*}
\bH(\zak;x)^{-1} \,\im  \begin{bmatrix}
- (1 + \frac{\sfa(x)}{2}) &  -\frac{\sfa(x)}{2}  \\
\frac{\sfa(x)}{2}  & 1 + \frac{\sfa(x)}{2} 
\end{bmatrix}  \, \bH(\zak;x) = 
\begin{bmatrix}
- \im \lambda(\zak;x)  & 0  \\
0  & \im \lambda(\zak;x) 
\end{bmatrix}  \,. 
 \end{equation*}
We prove the following preliminary lemma.
\begin{lemma}
\label{anal:emmino}
  The real-valued function  
    \begin{align}
        m_{\geq 1}(\zak; x):=-\log{(h(\zak; x)+g(\zak; x))} \in\Sigma \cF^{\R}_1[r, N]\,, \quad h, g \ \text{in} \  \eqref{diag4}. 
        \label{emmino01}
    \end{align}
     Furthermore, there exist $s_0>\frac{11}{4}, r>0$ such that the map
     $B_{s_0, \R}(r)\ni \zak \mapsto m_{\geq 1}(\zak; \cdot)\in H^{s_0-\frac{9}{4}}(\T;\R)$ is analytic
     and its cubic part $m_{\geq 3} := \cP_{\geq 3} [ m_{\geq 1}]$ satisfies
     \be\label{d.m.geq.3}
        \| \di m_{\geq 3}(\zak)[\hat{\zak}]\|_{L^{\infty}} \lesssim \|\zak\|^{2}_{s_0} \|\hat \zak\|_{s_0} \quad \forall \hat \zak \in  H^{s_0}_\R(\T; \C^2)\,.
     \ee
\end{lemma}
     \begin{proof}
By the definitions of $h$ and $g$ in \eqref{diag4} and of $\lambda$ in \eqref{diag2}, $m_{\geq 1}$ 
in \eqref{emmino01}
depends on $(\eta, \psi)$ only through an analytic function of the Taylor coefficient $\sfa$ defined in \eqref{a.taylor}.
Therefore, by \eqref{BVa.ana}, for any $\s \geq \s_0$ there is $r = r(\s)$ so that 
 $B_\s(r)\times B_\s(r) \ni(\eta, \psi)\mapsto m_{\geq 1}(\eta, \psi)\in   H^{\s-2}(\T;\R)$ is analytic as well as  
 $(\eta, \psi) \mapsto \sfa(\eta, \psi)$.
The claimed analyticity of $\zak \mapsto m_{\geq 1}(\zak)$ then follows from \eqref{app_L} for the linear map  $ \zak \mapsto  (\eta, \psi)\equiv \Lmap \zak $. 
Moreover, the expansion claimed in \eqref{emmino01} follows from analyticity and the fact that $ m_{\geq 1}(0\,; x)\equiv 0$.

It remains to prove  estimate \eqref{d.m.geq.3}. By Cauchy estimates, for any $ \| \zak\|_{s_0}\leq \frac23 r$ and any $\hat \zak \in  H^{s_0}_\R(\T; \C^2)$,  using that $\|m_{\geq 3} (\zak)\|_{s_0 - \frac 9 4} \lesssim \|\zak\|_{s_0}^3,$ one has
$$
|\di {m}_{\geq 3} (\zak)[\hat \zak]|_{L^\infty}  \lesssim \| \di {m}_{\geq 3} (\zak)\|_{\cL(\dot H^{s_0};H^{s_0 -9\!/\!4})}  \| \hat \zak\|_{s_0} \lesssim \sup_{\|\wt \zak- \zak \|_{s_0}< \| \zak\|_{s_0}/2}\frac{ \| {m}_{\geq 3} (\wt \zak)\|_{s_0-\frac94}}{\| \zak\|_{s_0}/2} \| \hat \zak\|_{s_0}  \lesssim \|\zak\|^{2}_{s_0} \| \hat \zak\|_{s_0}\,.
$$
     \end{proof}

  We diagonalize the paradifferential operator of order $\frac12$.
 \begin{proposition}[{\bf Block-diagonalization at order $\frac12$}]
 \label{diag}
 Let $N\in\N$ and $\vr, \wt\vr>1$. 
There exist $s_0, r  >0 $ and a $(0,0)$-admissible transformation $\mathbf{\Psi}(\zak) \in \cM^0_{\geq 0}[r]$ with gain $\wt\varrho$ (recall Definition \ref{admtra}) such that if $(\eta, \psi)(t) \in B_{s_0, \R}(I;r)$ solves \eqref{eq:etapsi}, and $\zak=\zak(t)$ is defined in \eqref{zak}, then the variable
  $U_0:= \mathbf{\Psi}(\zak)^{-1}U$, 
  with $U$ solving \eqref{complexo},  solves the system
\begin{equation}\label{NuovoParaprod}
\begin{aligned}
\pa_t U_0& =  \vOpbw{ - \ii \frac{\gamma}{2} \sign(\xi)- \ii \sfV_\gamma(\zak;x)\xi  - \ii \lambda(\zak;x) \omega(\xi) -\ii \tilde d^{(-\frac12)}_{\geq 1}(\zak;x, \xi) } U_0+ \zOpbw{\tilde \tb_{\geq 1}^{(0)}(\zak;x,\xi)}U_0\\
&\quad +\Opbw{\bB_{\geq N}^{(0)}(\zak)}U_0 + \bR_{\geq 1}(\zak)  U_0
\end{aligned}
\end{equation}
where the real-valued function $\sfV_\gamma(\zak;x)\in \Sigma \mF^{\R}_{1}[r, N]$ is defined in \eqref{def:Vgamma}, the real-valued symbol $ \lambda(\zak;x) \in \Sigma \cF^{\R}_{0}[r, N]$ is defined in \eqref{diag2}, the Fourier multiplier $ \omega(\xi) \in \widetilde\Gamma^{\frac12}_0$ is defined in \eqref{omegaxi} and 
\begin{itemize}
\item $\tilde d^{(-\frac12)}_{\geq 1}(\zak;x, \xi)$  and $\tilde \tb_{\geq 1}^{(0)}(\zak;x, \xi)$ are pluri-homogeneous  symbols respectively in $\Sigma_1^{N-1} \wt\Gamma^{-\frac12}$ and $\Sigma_1^{N-1} \wt\Gamma^0$ which  satisfy the linear Hamiltonian property \eqref{lin:ham_complex}; 
\item $ {\bB}^{(0)}_{\geq N}(\zak; x,\xi)$  is a real-to-real matrix of   symbols  in $  \Gamma^{0}_{\geq N}[r]  $;
\item $\bR_{\geq 1}(\zak)$ is a real-to-real matrix of smoothing operators in $\Sigma \cR^{-\vr }_1[r,N]$.
\end{itemize}
\end{proposition}

\begin{proof}
By Lemma 4.2 in \cite{BFP} (see also the proof of Lemma 6.4 in \cite{BMM2} for an explicit expression of $m_{\geq 1}$),
the time 1 flow $\mathbf{\Psi}(\zak):=  \mathbf{\Psi}^{\tau}(\zak)|_{\tau=1}$  of 
\begin{equation}\label{def.m01}
\begin{cases}
\pa_\tau\mathbf{\Psi}^{\tau}(\zak)=
 \, \Opbw{\bM_{\geq 1}(\zak;x)}
 \mathbf{\Psi}^{\tau}(U)\\
\mathbf{\Psi}^0(\zak)={\rm Id} \, , 
\end{cases}\ ,  \quad 
\bM_{\geq 1}(\zak;x)
:=  \begin{bsmallmatrix}
 0 &  m_{\geq 1}(\zak;x) \\
  m_{\geq 1}(\zak;x) & 0
 \end{bsmallmatrix}\, ,  
\end{equation}
 and  $m_{\geq 1}$  in \eqref{emmino01} fulfills 
\begin{equation}\label{expo}
\mathbf{\Psi}(\zak) = \Opbw{\bH(\zak;x)}+ \bR(\zak) \, , 
\quad
\mathbf{\Psi}(\zak)^{-1} = \Opbw{\bH^{-1}(\zak;x)}+ \bR'(\zak) \ ,
\end{equation} 
where  $\bH(\zak;x)$ is defined in \eqref{diag4} 
 and  $\bR(\zak) , \bR'(\zak)  $ are 
 real-to-real
 matrices of smoothing operators in $ \Sigma\mathcal{R}^{-\vr}_{1}[r,N]$. 
 First we check that $\mathbf{\Psi}(\zak)$ preserves the real-to-real and linearly Hamiltonian structure. Recall that equation \eqref{complexo} has the form
 $$
 \partial_t U = \mathbf{Y}(\zak) U + \Opbw{\bB_{\geq N}^{(-\frac12)}(\zak)}U +\mathbf{R}_{\geq 1}(\zak) U\,,
 $$
where, recalling  
 \eqref{diag1} and \eqref{complexo},  $\bY(\zak)$ is the  real-to-real, linearly Hamiltonian  (cf. \eqref{lin:ham_complex})  matrix of paradifferential operators  
\begin{align}
\label{def:YOp} 
\bY(\zak) & = \Opbw{\begin{bmatrix}
     a & b \\ \ov{b}^{\vee} & \ov{a}^{\vee}
 \end{bmatrix}}\,, \quad  \ov{a} = - a\,, \quad b^\vee = b\,,
\end{align}
where the symbols $a$ and $b$ have the form
\begin{align*}
    a:= -\ii \frac{\gamma}{2}\sign{\xi}- \ii \sfV_\gamma \xi - \ii \big(1+ \frac{\sfa}{2}\big) \omega(\xi)+ d_{\geq 1}^{(-\frac12)}\,, \quad b:= -\frac14 (\sfV_\gamma)_x - \ii \frac{\sfa}{2}\omega(\xi)+ r_{\geq 1}^{(-\frac{1}{2})}\,,
\end{align*}
with  $d_{\geq 1}^{(-\frac12)}, \, r_{\geq 1}^{(-\frac12)} \in \Sigma_1^{N-1} \wt\Gamma^{-\frac12}$  the diagonal and out-diagonal symbols of the linearly Hamiltonian matrix $\wt\bA^{(-\frac12)}_{\geq 1}(\zak)$ in \eqref{complexo},  
 $\mathbf{B}_{\geq N}^{(-\frac 12)}(\zak)$ is a real-to-real matrix of symbols in $ \Gamma_{ \geq N}^{-\frac 1 2}[r]$, and $\mathbf{R}_{\geq 1}(\zak)$ a {real-to-real} matrix of smoothing operators in $  \Sigma\mathcal{R}^{-\vr}_{1}[r,N]$.
 
 Since $ U $ solves \eqref{complexo} then the variable $U_0 = \mathbf{\Psi}(\zak)^{-1} U $ solves
 \be
\begin{aligned}
\pa_t U_0   =&  \mathbf{\Psi}(\zak)^{-1}  \,  
\Big[ 
\bY(\zak) + \Opbw{\bB^{(-\frac 1 2)}_{\geq N} (\zak)}
 + \bR_{\geq 1}(\zak) 
  \Big] \mathbf{\Psi}(\zak)  \, U_0  
+   (\pa_t \mathbf{\Psi}(\zak)^{-1}) \mathbf{\Psi}(\zak)  \, U_0 \\
\stackrel{\eqref{expo}}{=}& \Opbw{\bH(\zak)^{-1} \#_\uprho \left(\begin{bmatrix}
     a & b \\ \ov{b}^{\vee} & \ov{a}^{\vee}
 \end{bmatrix} + \bB^{(-\frac 1 2)}_{\geq N}(\zak)\right) \#_\uprho \bH(\zak) } + (\pa_t \mathbf{\Psi}(\zak)^{-1}) \mathbf{\Psi}(\zak)  \, U_0 + \bR'_{\geq 1}(\zak) U_0\,,
\end{aligned}
\label{V1.eq}
\ee
where 
$\uprho = \vr +1$ and therefore, by \eqref{sharp3},
$\bR'_{\geq 1}(\zak)$ is a {real-to-real} matrix of smoothing operators in $\Sigma \cR^{-\uprho+1}_1[r,N] \equiv \Sigma \cR^{-\vr}_1[r,N]$.
We compute each term in \eqref{V1.eq}, using the expressions 
 for  $\bH(\zak)$ and $\bH(\zak)^{-1}$ in \eqref{diag4}, \eqref{diag5}. 
The first is 
\begin{equation}
\bH(\zak)^{-1} \#_\uprho \begin{bmatrix}
     a & b \\ \ov{b}^{\vee} & \ov{a}^{\vee}
 \end{bmatrix} \#_\uprho \bH(\zak)  = \begin{bmatrix}
        h & -g \\
        -g & h
    \end{bmatrix} \#_\uprho   \begin{bmatrix}
        a & b \\
\overline{b}^{\vee}& \overline{a}^{\vee}
    \end{bmatrix}  \#_\uprho \begin{bmatrix}
        h & g \\
        g & h
    \end{bmatrix} = \begin{bmatrix}
        d & r \\
\overline{r}^{\vee}& \overline{d}^{\vee}
    \end{bmatrix}\,,
    \label{der}
\end{equation}
with 
\begin{equation}
\begin{gathered}
     d = h \#_\uprho a \#_\uprho h - g \#_\uprho \ov{a}^{\vee} \#_\uprho g + h \#_\uprho b \#_\uprho g - g \#_\uprho \ov{b}^{\vee} \#_\uprho h\,, \\
     r = h \#_\uprho a \#_\uprho g - g \#_\uprho \ov{a}^{\vee} \#_\uprho h + h \#_\uprho b \#_\uprho h - g \#_\uprho  \ov{b}^{\vee} \#_\uprho g\,.
\end{gathered}
\label{cancelletti_rule}
\end{equation}
Using properties \eqref{prop:ov3} and $ \ov{a} = - a\,, b^\vee = b$, one also gets 
\begin{equation*}
    \ov{d} = - d  , \quad r^\vee = r \ , 
\end{equation*}
showing that the matrix of symbols in \eqref{der} is linearly Hamiltonian.
Then, by  \eqref{canc:diag}, the symbols $d$ and $r$ in \eqref{cancelletti_rule} expand as 
\begin{align}
d=-\ii \frac{\gamma}{2} \sign{\xi}-\ii \sfV_\gamma \xi-\ii \lambda \omega(\xi)+ d^{(-\frac12)}_{\geq 1}, \quad r=  {\sfV_\gamma(h_x g - g_x h)}- \frac14(\sfV_\gamma)_x+r^{(-\frac12)}_{\geq 1}\,,
\label{nuov:der}
\end{align}
where $\lambda$ is defined in \eqref{diag2} while $d^{(-\frac12)}_{\geq 1}$ and $r^{(-\frac12)}_{\geq 1}$ are symbols in $\Sigma \Gamma_{ 1}^{-\frac12}[r,N]$ which satisfy the linear Hamiltonian property \eqref{lin:ham_complex} by difference with $d$ and $r$.

We now compute $ \bH^{-1}(\zak) \#_\uprho \bB_{\geq N}^{(-\frac12)}(\zak)\#_\uprho \bH(\zak)$, where $\bB_{\geq N}^{(-\frac12)}(\zak)$ has the form 
$
\begin{bmatrix}
    d_{\geq N}^{(-\frac12)} & r^{(-\frac12)}_{\geq N}\\[0.3em] \big(\ov{r^{(-\frac12)}_{\geq N}}\big)^\vee &\big(\ov{d^{(-\frac12)}_{\geq N}}\big)^\vee
\end{bmatrix}$ 
for symbols  $d_{\geq N}^{(-\frac12)}$, $ r^{(-\frac12)}_{\geq N}$ in $ \Gamma_{\geq N}^{-\frac12}[r]$.
Since $ \bH^{-1}(\zak) \#_\uprho \bB_{\geq N}^{(-\frac12)}(\zak)\#_\uprho \bH(\zak)$ satisfies the same conjugation rule  \eqref{der} with $a\leadsto d_{\geq N}^{(-\frac12)}$ and $ b \leadsto r_{\geq N}^{(-\frac12)}$, we get that $\bH^{-1}(\zak) \#_\uprho \bB_{\geq N}^{(-\frac12)}(\zak) \#_\uprho \bH(\zak)$ is a matrix of symbols in $ \Gamma_{\geq N}^{-\frac12}[r]$.

Next we consider  $(\pa_t \mathbf{\Psi}(\zak)^{-1} )\mathbf{\Psi}(\zak)$. 
By using formula (A.4) of \cite{BFP} applied with $L=1$,
we have the following expansion for the time derivative
\begin{align}   \label{detLie}
    (\pa_t \mathbf{\Psi}(\zak)^{-1} )\mathbf{\Psi}(\zak)& =\Opbw{- \pa_t \bM_{\geq 1}(\zak)}
    \\\notag & \qquad
    + \int_0^1 (1-\theta)\mathbf{\Psi}^{-\theta}(\zak) \big[\Opbw{\bM_{\geq 1}(\zak)}, \Opbw{\pa_t \bM_{\geq 1}(\zak)}\big] \mathbf{\Psi}^{\theta}(\zak)\, \di \theta\,.
\end{align}
By \eqref{def.m01},  Lemma \ref{anal:emmino}  and \eqref{claim:zak.M},  we have
 \begin{equation} \label{Ham.terza:eq}
 \pa_t \bM_{\geq 1}(\zak) =
\begin{pmatrix}
0 & \!\! \!\!  \pa_t m_{\geq 1}(\zak)\\     \pa_t m_{\geq 1}(\zak)  & \!\! \!\!  0 
\end{pmatrix}\, \ , \quad \pa_t m_{\geq 1}(\zak)=\di m_{\geq 1}(\zak)[ \cX_{\cH}(\zak)] \in \Sigma\cF_1^\R[r,N] \ ,
\end{equation}
proving that $\pa_t \bM_{\geq 1}(\zak)$ is a linearly Hamiltonian matrix of functions in $ \Sigma \cF_1^\R[r,N]$. 
It remains to prove that the integral term in \eqref{detLie} is a {real-to-real matrix} of smoothing operators in $ \Sigma \cR^{-\vr}_1[r,N]$. Indeed, by Proposition \ref{teoremadicomposizione}, we have 
\begin{align}
    \left[ \Opbw{\pa_t \bM_{\geq 1}(\zak)}, \Opbw{{\bM_{\geq 1}}(\zak)} \right]= \overbrace{\Opbw{\pa_t \bM_{\geq 1}(\zak)\#_\uprho \bM_{\geq 1}(\zak)-\bM_{\geq 1}(\zak)\#_\uprho \pa_t \bM_{\geq 1}(\zak)}}^{=0} + \bR(\zak)
\end{align}
with $\bR(\zak)$ a {real-to-real} matrix of smoothing operators in $\Sigma \cR^{-\vr}_2[r, N]$,  where we used that $ \pa_t\bM_{\geq 1} (\zak)$ and $\bM_{\geq 1}(\zak)$ are commuting matrices of functions (thus depending only on $ x$). Since $ \mathbf{\Psi}^\theta(\zak)$ conjugates matrices of real-to-real smoothing operators into matrices of real-to-real  smoothing operators, it follows that 
\begin{align}
    (\pa_t \mathbf{\Psi}(\zak)^{-1} )\mathbf{\Psi}(\zak) =  \Opbw{\begin{bmatrix}
0 & \!\! \!\!  \pa_t m_{\geq 1}(\zak)\\     \pa_t m_{\geq 1}(\zak)  & \!\! \!\!  0 
\end{bmatrix}}+ \bR(\zak)\,,
    \label{detM:final}
\end{align}
 with $\bR(\zak)$ a {real-to-real} matrix of  smoothing operators in $\Sigma \cR_{1}^{-\vr}[r, N]$. By \eqref{nuov:der} and  \eqref{detM:final}, equation \eqref{NuovoParaprod} follows defining 
   \begin{align*}
       &\tb_{\geq 1}^{(0)}:= -\frac14 (\sfV_\gamma)_x + \sfV_\gamma(h_x g - g_x h)+r^{(-\frac12)}_{\geq 1}+\pa_t m_{\geq 1}, \quad \wt \tb_{\geq 1}^{(0)}:=\cP_{\leq N-1} \left(\tb_{\geq 1}^{(0)}\right) \\
       &\wt d^{(-\frac12)}_{\geq 1}:= -\ii \cP_{\leq N-1}\left[ \td^{(-\frac12)}_{\geq 1}\right]\\
       & \bB^{(0)}_{\geq N}(\zak):= \bH^{-1}(\zak) \#_\uprho \bB_{\geq N}^{(-\frac12)}(\zak)\#_\uprho \bH(\zak)+ \vOpbw{\cP_{\geq N}\left( \td^{(-\frac12)}_{\geq 1}\right)} + \zOpbw{\cP_{\geq N}\left(\tb_{\geq 1}^{(0)}\right)}\,.
   \end{align*}
It remains to check that the transformation $\mathbf{\Psi}(\zak)$ is $(0,0)$ admissible with gain $\wt\varrho$. This follows from Lemma \ref{lem:flow.ad}-$(ii)$, whose assumptions \eqref{du.g3} are satisfied thanks to  \eqref{d.m.geq.3}.
\end{proof}
\subsubsection{Block-Diagonalization at  lower order}
In the following proposition we cancel out the off-diagonal symbols up to order $-\vr$, for {an arbitrary but fixed} $\vr>1$ and up to homogeneity $N$. 
\begin{proposition}[{\bf Block-diagonalization up to  order $-\vr$ and homogeneity $N$}]
 \label{diag.ord0}
 Let $N\in \N$ and $\vr, \wt\vr>1$. 
There are $s_0, r  >0 $ and a $(0, 0)$ admissible transformation $\mathbf{\Psi}_0(\zak) \in {\cM^{0}_{\geq 0}[r]}$ with gain $\wt\vr$ such that if $(\eta, \psi)(t) \in B_{s_0, \R}(I;r)$ solves \eqref{eq:etapsi}, and $\zak=\zak(t)$ is defined in \eqref{zak}, then the variable
  $U_1:= \mathbf{\Psi}_0(\zak)U_0$, 
  with $U_0$ solving \eqref{NuovoParaprod},  solves the system
\begin{equation}\label{eq.diago.ord0}
\begin{aligned}
\pa_t U_1 &= \vOpbw{ - \ii \frac{\gamma}{2} \sign(\xi) - \ii \sfV_\gamma(\zak;x)\xi  - \ii \lambda(\zak;x) \omega(\xi) -\ii \tilde \td^{(-\frac12)}_{ \geq 1}(\zak; x, \xi) } U_1\\
& \qquad + \Opbw{{\bB}^{(0)}_{\geq N}(\zak; x,\xi) 
} U_1
+ \bR_{\geq 1}(\zak)  U_1
\end{aligned}
\end{equation}
where  the real-valued function $\sfV_\gamma(\zak;x)\in \Sigma \mF^{\R}_{1}[r,N]$ is defined in \eqref{def:Vgamma}, the real-valued symbol $ \lambda(\zak;x) \in \Sigma \cF^{\R}_{0}[r,N]$ is defined in \eqref{diag2}, the Fourier multiplier  $ \omega(\xi) $  in \eqref{omegaxi}, and 
\begin{itemize}
\item  $\tilde \td^{(-\frac12)}_{\geq 1}(\zak; x, \xi)$ is a pluri-homogeneous real-valued symbol in $  \Sigma_1^{N-1}\wt\Gamma^{-\frac12}$;

\item $ {\bB}^{(0)}_{\geq N}(\zak; x,\xi)$  is a real-to-real matrix of   symbols  in $  \Gamma^{0}_{\geq N}[r]  $;
\item $\bR_{\geq 1}(\zak)$ is a real-to-real matrix of smoothing operators in $ \Sigma \cR^{-\vr}_{1}[r,N] $.  
\end{itemize}
\end{proposition}
Proposition \ref{diag.ord0} is proved by the following iterative lemma. We anticipate that $\vr$ will be taken as an integer number later on in \eqref{paraN}.
\begin{lemma}\label{lem:indud}
For $ j = 0, \ldots,  2 \vr $, there are \\
\noindent $ \bullet $ 
a linearly Hamiltonian 
paradifferential operator (according to  \eqref{lin:ham_complex})
of the form  
 \begin{align}\label{sist2 j-th} 
\mathcal{Y}_{j}(\zak) :=\vOpbw{
-\ii \frac{\gamma}{2}  \sign(\xi) -\ii \sfV_\gamma \xi -\ii \lambda\omega(\xi) - \ii \wt\td^{(-\frac12)}_{\geq 1, j}}
+\zOpbw{\wt\tb_{\geq 1,j}^{(-\frac{j}{2})}}
\end{align}
where $\sfV_\gamma$, $\lambda$ and $\omega(\xi)$ are as in \Cref{diag.ord0}, the symbol
$\wt\td^{(-\frac12)}_{\geq 1, j}\in \Sigma_1^{N-1}\wt\Gamma^{-\frac12}$ is real-valued and the symbol $\wt\tb_{\geq 1,j}^{(-\frac{j}{2})}\in\Sigma_1^{N-1}\wt\Gamma^{-\frac{j}{2}}$ is even in $\xi$; 

\noindent $ \bullet  $  a real-to-real matrix 
${\bB}^{(0)}_{\geq N, j}(\zak; x,\xi)$
of symbols in $\Gamma_{\geq N}^{0}[r]$ and
a real-to-real  matrix $\bR_{\geq 1, j}(\zak) $  of smoothing operators in  
$\Sigma\mathcal{R}^{-\vr}_{1}[r,N]$;\\
$\bullet$ a  $(0,0)-$admissible transformation ${\bf \Phi}_{j}(\zak)
\equiv {\bf \Phi}_{j}^{\theta}(\zak)\vert_{\theta = 1}$
with arbitrary gain  $\wt\vr \geq 1$,  where   
${\bf \Phi}_{j}^{\theta}(\zak)$ is the time flow
\be\label{flow-jth}
\partial_{\theta} {\bf \Phi}^{\theta}_{j} (\zak) =  \zOpbw{ \wt m_{\leq (N-1)}^{(-\frac{j+1}{2})}(\zak; x, \xi)} {\bf \Phi}_{j}^{\theta}(\zak) \, ,
 \quad {\bf\Phi}_{j}^{0}(\zak) = {\rm Id} \, ,  
\ee
generated by  the even in $\xi$ symbol of negative order (recall \eqref{pienne})
\begin{equation}\label{generatore-jth}
\wt m_{\leq (N-1)}^{(-\frac{j+1}{2})} := \cP_{\leq (N-1)}\left[\frac{\wt\tb_{\geq 1,j}^{(-\frac{j}{2})}}{2\lambda}\right]\frac{1}{-\ii \omega(\xi)}\in
\Sigma_1^{N-1}\wt\Gamma^{-\frac{j+1}{2}} \,  , 
\end{equation}
such that, if $ W_j $, $ j = 0, \ldots,  2 \vr -1 $, solves 
\be\label{sist1 j-th}
\pa_{t}W_j=\left(\mathcal{Y}_{j}(\zak)+\Opbw{{\bB}^{(0)}_{\geq N, j}(\zak)}+\bR_{\geq 1, j}(\zak)\right)W_j\,, 
\ee
then the variable 
\begin{equation}\label{nuovavarjth}
W_{j+1}:={\bf \Phi}_{j}(\zak)W_{j}
\end{equation}
 solves a system of the form \eqref{sist1 j-th} with $ j + 1 $ instead of $ j $.
\end{lemma}
\begin{proof}
The proof proceeds by induction. 
\\[1mm]
{\bf Initialization.} System \eqref{NuovoParaprod} is \eqref{sist1 j-th} for $ j = 0 $ where  the 
paradifferential operator $ \mathcal{Y}_{0}(\zak)  $ has the form   \eqref{sist2 j-th} with 
$\wt\td_{\geq 1, 0}^{(-\frac12)}=\wt d_{\geq 1}^{(-\frac12)}$, 
$\wt\tb_{\geq 1, 0}^{(0)}=\wt \tb_{\geq 1}^{(0)}$, $\bB_{\geq N, 0}^{(0)}= \bB_{\geq N}^{(0)}$ and $\bR_{\geq 1, 0}= \bR_{\geq 1}$.
\\[1mm]
{\bf Iteration.} 
We now argue by induction. Suppose that $ W_j $ solves system \eqref{sist1 j-th} with operators
$ \mathcal{Y}_{j}(\zak)  $ of  the form   \eqref{sist2 j-th}, a real-to-real matrix of symbols $\bB_{\geq N, j}^{(0)}(\zak)\in \Gamma_{\geq N}^0[r]$ and a
real-to-real matrix of smoothing operators 
$\bR_{\geq 1, j} (\zak) $  in   $\Sigma\mathcal{R}^{-\vr}_{1}[r,N]$. 

Then the variable  $ W_{j+1} $ defined in  \eqref{nuovavarjth} solves the system  
\be\label{sis:j+1}
\begin{aligned}
\pa_{t}W_{j+1} = &  (\pa_t  {\bf \Phi}_{j}(\zak)) {\bf \Phi}_{j}(\zak)^{-1} W_{j+1}  +
 {\bf \Phi}_{j}(\zak)  \mathcal{Y}_j (\zak)   {\bf \Phi}_{j}(\zak)^{-1} W_{j+1}  \\
 & +
{\bf \Phi}_{j}(\zak) \Opbw{{\bB}^{(0)}_{\geq N, j}(\zak)} {\bf \Phi}_{j}(\zak)^{-1}  W_{j+1}  +{\bf \Phi}_{j}(\zak)\bR_{\geq 1, j}(\zak) {\bf \Phi}_{j}(\zak)^{-1}  
 W_{j+1}\,.
\end{aligned}
\ee
The operator
$(\pa_t {\bf \Phi}_{j}(\zak)){\bf \Phi}_{j}(\zak)^{-1}$, setting ${\rm Ad}_A[B]:=[A, B]$,
admits the Lie expansion  \cite[formula (A.4)]{BFP}
\be\label{pat.iter.lemma}
\begin{aligned}
    (\pa_t {\bf \Phi}_{j} (\zak) )({\bf \Phi}_{j}(\zak))^{-1}   
&=  \zOpbw{ \pa_t\wt m_{\leq (N-1)}^{(-\frac{j+1}{2})}} +
\sum_{q=2}^{L}\frac{1}{q!}{\rm Ad}_{\zOpbw{ \wt m_{\leq (N-1)}^{(-\frac{j+1}{2})}}}^{q-1}\Big[\zOpbw{\pa_t \wt m_{\leq (N-1)}^{(-\frac{j+1}{2})}} \Big] \\ 
& +\frac{1}{L!}\int_{0}^{1}
(1- \theta )^{L} 
{\bf \Phi}_{j}^{\theta}(\zak){\rm Ad}_{\zOpbw{ \wt m_{\leq (N-1)}^{(-\frac{j+1}{2})}}}^{L}\Big[ \zOpbw{ \pa_t\wt m_{\leq (N-1)}^{(-\frac{j+1}{2})}} \Big]({\bf \Phi}_{j}^{\theta}(\zak))^{-1} d \theta \, .
\end{aligned}
\ee
Since $\wt m_{\leq (N-1)}^{(-\frac{j+1}{2})}$ is a homogeneous symbol, by  \eqref{claim:zak.M} and  \Cref{comp.plurihom} (see also  \Cref{sym-comp-2} of \Cref{lem:general_composition}), its time derivative $\pa_t \wt m_{\leq (N-1)}^{(-\frac{j+1}{2})}$ fulfills, 
\begin{equation}\label{pat.tm}
    (\pa_t \wt m_{\leq (N-1)}^{(-\frac{j+1}{2})})(\zak)=\di \wt m_{\leq (N-1)}^{(-\frac{j+1}{2})}(\zak)[ \cX_{\cH}(\zak)] \in \Sigma\Gamma_1^{-\frac{j+1}{2}}[r, N]
    \ \quad  \mbox{ and } \quad (\pa_t \wt m_{\leq (N-1)}^{(-\frac{j+1}{2})})(\zak; x, \xi) \mbox{ even in } \xi \ ;
\end{equation}
In particular, $\zOpbw{\pa_t  \wt m_{\leq (N-1)}^{(-\frac{j+1}{2})}}$ is linearly Hamiltonian. 
\\
We recall (see \eqref{commutator}) that
$$
 \wt m_{\leq (N-1)}^{(-\frac{j+1}{2})}\#_{\vr}\pa_{t}  \wt m_{\leq (N-1)}^{(-\frac{j+1}{2})}-\pa_{t} \wt m_{\leq (N-1)}^{(-\frac{j+1}{2})}\#_{\vr}\wt m_{\leq (N-1)}^{(-\frac{j+1}{2})}
\in \Sigma\Gamma_{2}^{-(j+1)-1}[r,N]
$$ 
Therefore, by using the composition Proposition \ref{teoremadicomposizione}, in particular the matrix composition formula \eqref{commurule},  \eqref{pat.tm} and the last bullet at page \pageref{comm.op.vz},
 we have that 
${\rm Ad}_{\zOpbw{\wt m_{\leq (N-1)}^{(-\frac{j+1}{2})}}}[ \zOpbw{\pa_{t}\wt m_{\leq (N-1)}^{(-\frac{j+1}{2})}}]$ is a linearly Hamiltonian matrix of paradifferential operators 
with symbol in $\Sigma\Gamma^{-(j+1)-1}_{2}[r,N]$ plus a matrix of smoothing operators in $\Sigma\mathcal{R}^{-\vr}_{2}[r,N]$.
As a consequence, we deduce, for $k\geq1$,
\begin{align*}
{\rm Ad}^{k}_{ \zOpbw{\wt m_{\leq (N-1)}^{(-\frac{j+1}{2})}}}& [\zOpbw{\pa_t\wt m_{\leq (N-1)}^{(-\frac{j+1}{2})}}]=\opbw(\bB_{k}(\zak))+\bR_{k}(\zak)
\end{align*}
where 
$\bB_{k}(\zak)$ is a matrix of symbols in $ \Sigma\Gamma_{k+1}^{-\frac{j+1}{2}(k+1)-k}[r,N] $ and 
$\bR_{k}(\zak)$ a real-to-real matrix of smoothing operators in $\Sigma \mathcal{R}^{-\vr}_{k+1}[r,N]$.
As observed above, by the last bullet at page \pageref{comm.op.vz}, since the  paradifferenetial operators on the left are linearly Hamiltonian, so is the matrix of symbols $\bB_k(\zak)$.

By  taking   $ L $ in \eqref{pat.iter.lemma} large enough with respect to $ \vr $, by using also Lemma \ref{lem:general_composition} to estimate the integral in \eqref{pat.iter.lemma}, we obtain that 
$(\pa_t  {\bf \Phi}_1(\zak)) {\bf \Phi}_1(\zak)^{-1}$
is a matrix of
 paradifferential operators with symbols in $ \Sigma\Gamma^{-\frac{j+1}{2}}_{1}[r,N]$ plus a matrix of
 smoothing operators in $\Sigma\mathcal{R}^{-\vr}_{2}[r,N]$.
 \smallskip
 
 Next we consider $ {\bf \Phi}_{j}(\zak)\,  \mathcal{Y}_{j}(\zak) \,  {\bf \Phi}_{j}(\zak)^{-1}$. We  apply the Lie expansion \cite[formula (A.3)]{BFP}
 \begin{align}\label{Lie1Vec}
  {\bf \Phi}_{j}&(\zak)\,  \mathcal{Y}_{j}(\zak) \,  {\bf \Phi}_{j}(\zak)^{-1}  
  = \mathcal{Y}_j(\zak)
  + [ \zOpbw{\wt m_{\leq (N-1)}^{(-\frac{j+1}{2})}(\zak; \cdot)} , \mathcal{Y}_j(\zak) ] \\
  \notag
 &  +
  \sum_{q=2}^{L}\frac{1}{q!}{\rm Ad}_{\zOpbw{\wt m_{\leq (N-1)}^{(-\frac{j+1}{2})}}}^{q}[\mathcal{Y}_j(\zak) ]
 +  \frac{1}{L!} \int_{0}^{1}  (1- \theta)^{L} {\bf \Phi}_j^{\theta}(\zak) {\rm Ad}_{\zOpbw{\wt m_{\leq (N-1)}^{(-\frac{j+1}{2})}}}^{L+1}[\mathcal{Y}_j(\zak) ] ({\bf \Phi}_j^{\theta}(\zak))^{-1} 
\di  \theta \,.
\end{align}
 Let us first analyze the term 
$ \Big[ \zOpbw{\wt m_{\leq (N-1)}^{(-\frac{j+1}{2})}(\zak)} , \mathcal{Y}_j(\zak)   \Big]$ using the expression of $\mathcal{Y}_j(\zak)   $ in  \eqref{sist2 j-th}.
Applying \eqref{commurule} with $\vr\leadsto\vr+1$, we obtain
\begin{equation}\label{esp2-jth}
\begin{aligned}
\Big[ 
\zOpbw{\wt m_{\leq (N-1)}^{(-\frac{j+1}{2})} }, 
\vOpbw{ -\ii \frac{\gamma}{2} \sign{\,\xi} -\ii \sfV_\gamma\,\xi - \ii \lambda\omega(\xi)}
\Big] 
=\zOpbw{ 2\ii \wt m_{\leq (N-1)}^{(-\frac{j+1}{2})} \lambda\omega(\xi)}
\end{aligned}
\end{equation}
up to a linearly Hamiltonian matrix of paradifferential operators, with symbols in $\Sigma\Gamma^{-\frac{j+1}{2}}_{1}[r,N]$, and a real-to-real matrix of smoothing operators in $\Sigma\cR^{-\vr}_{{1}}[r,N]$.
Moreover, we have that
$
\Big[\zOpbw{\wt m_{\leq (N-1)}^{(-\frac{j+1}{2})}} , \vOpbw{ \ii \wt\td^{(-\frac12)}_{\geq 1, j}} \Big]
$
is a matrix of linearly Hamiltonian paradifferential operators with symbols in $\Sigma\Gamma^{-\frac{j+2}{2}}_{2}[r,N]$ up to a matrix of smoothing operators in $\Sigma\mathcal{R}^{-\vr}_{{2}}[r,N]$. In the same way one has that
$
\Big[ \zOpbw{\wt m_{\leq (N-1)}^{(-\frac{j+1}{2})}} ,\zOpbw{\wt\tb_{\geq 1,j}^{(-\frac{j}{2})}}\Big]
$
is a matrix of linearly Hamiltonian paradifferential operators with symbols
in $\Sigma\Gamma^{-\frac{2j+1}{2}}_{{2}}[r,N]$ up to a real-to-real matrix of smoothing operators in $\Sigma \cR_{2}^{-\vr}[r, N]$.
It follows that the off-diagonal symbols of order $ -\frac{j}{2} $ in the first line of \eqref{Lie1Vec} are 
of the form  $
\zOpbw{q_j(\zak;x,\x)}$
 with  
$$
q_j\stackrel{\eqref{sist2 j-th}, \eqref{esp2-jth}}{:=}\wt\tb_{\geq 1,j}^{(-\frac{j}{2})}  + 2\ii \wt m_{\leq (N-1)}^{(-\frac{j+1}{2})} \lambda \omega(\xi)\stackrel{\eqref{generatore-jth}}{=} 
\cP_{\geq N}\Big[\frac{\wt\tb_{\geq 1,j}^{(-\frac{j}{2})}}{\lambda}\Big]\lambda=: \tb_{\geq N,j+1}^{(-\frac{j}{2})} \in \Gamma^{-\frac{j}{2}}_{\geq N}[r]\, . 
$$
All together we find
    \begin{equation*}
    \begin{aligned}
     \cY_j(\zak) +   [\zOpbw{\wt m_{\leq (N-1)}^{(-\frac{j+1}{2})} } ,\mathcal{Y}_j(\zak) ]= & \vOpbw{-\ii \frac{\gamma}{2}  \sign(\xi)
-\ii \sfV_\gamma \xi -\ii \lambda\omega(\xi) - \ii \wt\td^{(-\frac12)}_{\geq 1, j}} +  \zOpbw{\tb_{\geq N, j+1}^{(-\frac{j}{2})}} \\
& + \Opbw{\ring{\bB}_{j+1}(\zak)} + \ring{\bR}_{j+1}(\zak)
\end{aligned}
    \end{equation*}
    with $\ring{\bB}_{j+1}(\zak)$ a linearly Hamiltonian matrix of symbols in $  \Sigma\Gamma^{-\frac{j+1}{2}}_1[r,N]$ and 
    $\ring{\bR}_{j+1}(\zak) $ is a 
    real-to-real matrix of 
smoothing operators in 
  $\Sigma\mathcal{R}^{-\vr}_{1}[r,N]$.
  \\
Now we consider the second line of \eqref{Lie1Vec}
and prove that it 
is  a matrix of linearly Hamiltonian paradifferential operators with symbols in $ \Sigma\Gamma^{-\frac{j+1}{2}}_{1}[r,N] $  (being the sum of commutators of linearly Hamiltonian operators, cf.  the last bullet at page \pageref{comm.op.vz}), 
and a matrix of
 smoothing operators in $\Sigma\mathcal{R}^{-\vr}_{1}[r,N]$.
Indeed, using Proposition \ref{teoremadicomposizione}, we deduce, for $k\geq2$, 
  \[
{\rm Ad}^k_{\zOpbw{\wt m_{\leq (N-1)}^{(-\frac{j+1}{2})}}}[\mathcal{Y}_j(\zak) ]
=\opbw(\widetilde{\bB}_{k}(\zak))+\widetilde{\bR}_{k}(\zak), \qquad 
\widetilde{\bB}_{k}(\zak)\in \Sigma\Gamma_{k+1}^{-\frac{j+1}{2}k}[r, N] \, \mbox{ linearly Hamiltonian}, 
\]
and $\widetilde{\bR}_{k}(\zak) $  a
real-to-real
matrix of smoothing operators in $ \Sigma\mathcal{R}^{-\vr}_{k+1}[r,N]$.\\
 Finally, consider the integral term in the second line of  \eqref{Lie1Vec}; provided  $L$ is  large enough, the iterated commutator is a real-to-real matrix of smoothing operators
 in $\Sigma \cR_{1}^{-\vr}[r,N]$, as well as its conjugation with the flow ${\bf \Phi}_j^\theta(\zak)$ 
 (use  Lemma A.2 in \cite{BFP}).
\\
Finally consider $ {\bf \Phi}_{j}(\zak) \Opbw{{\bB}^{(0)}_{\geq N, j}(\zak)}  {\bf \Phi}_{j}(\zak)^{-1} $; using again the 
Lie Expansion  
\cite[formula (A.3)]{BFP}
 \begin{equation*}
 \begin{aligned}
  {\bf \Phi}_{j}(\zak) \Opbw{{\bB}^{(0)}_{\geq N, j}(\zak)} &  {\bf \Phi}_{j}(\zak)^{-1}  
  = \Opbw{{\bB}^{(0)}_{\geq N, j}(\zak)}
  +\sum_{q=1}^{L}\frac{1}{q!}{\rm Ad}_{\zOpbw{\wt m_{\leq (N-1)}^{(-\frac{j+1}{2})}}}^{q}[\Opbw{{\bB}^{(0)}_{\geq N, j}(\zak)} ]\\
 &+  \frac{1}{L!} \int_{0}^{1}  (1- \theta)^{L} {\bf \Phi}_j^{\theta}(\zak) {\rm Ad}_{\zOpbw{\wt m_{\leq (N-1)}^{(-\frac{j+1}{2})}}}^{L+1}[\Opbw{{\bB}^{(0)}_{\geq N, j}(\zak)} ] ({\bf \Phi}_j^{\theta}(\zak))^{-1} 
\di  \theta \  , 
\end{aligned}
 \end{equation*}
applying Proposition \ref{teoremadicomposizione} and taking $L$ large enough, we deduce that 
\begin{equation*}
     {\bf \Phi}_{j}(\zak)  \Opbw{{\bB}^{(0)}_{\geq N, j}(\zak)}  {\bf \Phi}_{j}(\zak)^{-1}   = 
     \Opbw{\breve{{\bB}}^{(0)}_{\geq N, j}(\zak)} + \breve \bR_{\geq N}(\zak) \ ,
\end{equation*}  
with $\breve{{\bB}}^{(0)}_{\geq N, j}(\zak)$  matrix of symbols in $\Gamma_{\geq N}^0[r]$ and $\breve \bR_{\geq N}(\zak)$ a real-to-real matrix of smoothing operators in $\cR_{\geq N}^{-\vr}[r]$.

 We conclude that the variable $W_{j+1}$ solves a system of the form \eqref{sist1 j-th} with $j\leadsto j+1$ where \\
$\bullet$ $\cY_{j+1}(\zak)$  is the  matrix of paradifferential operators  
$$
\cY_{j+1}(\zak):= \vOpbw{-\ii \frac{\gamma}{2} \sign{\,\xi}-\ii \sfV_\gamma\, \xi -\ii \lambda\omega(\xi) - \ii \wt\td^{(-\frac12)}_{\geq 1, j+1} }
+ \zOpbw{\wt\tb_{\geq 1,j+1}^{(-\frac{j+1}{2})}}
$$
where the symbols $\wt\td^{(-\frac12)}_{\geq 1, j+1}\in \Sigma_1^{N-1}\wt\Gamma^{-\frac12}$  and  $\wt\tb_{\geq 1,j+1}^{(-\frac{j+1}{2})}\in\Sigma_1^{N-1}\wt\Gamma^{-\frac{j+1}{2}}$  are defined  by  
\begin{equation}\label{lem.iter.new.sym}
    \vOpbw{- \ii \wt\td^{(-\frac12)}_{\geq 1, j+1} }
+ \zOpbw{\wt\tb_{\geq 1,j+1}^{(-\frac{j+1}{2})}} :=  
 \vOpbw{- \ii \wt\td^{(-\frac12)}_{\geq 1, j} } + 
\cP_{\leq (N-1)} \big[ 
 \bP_j^{new}(\zak)
\big] \ 
\end{equation}
and we have denoted by 
$ \bP_j^{new}(\zak)$  the linearly Hamiltonian matrix of paradifferential operators in $\Sigma \Gamma_1^{-\frac{j+1}{2}}[r, N]$
\begin{equation*}
\begin{aligned}
    \bP_j^{new}(\zak):= \zOpbw{ \pa_t\wt m_{\leq (N-1)}^{(-\frac{j+1}{2})}} + \sum_{k=2}^L \frac{1}{k!} \Opbw{\bB_{k-1}(\zak)} + \Opbw{\ring{\bB}_{j+1}(\zak)}
    + \sum_{k=2}^L \Opbw{\widetilde{\bB}_k(\zak)}  \ ; 
\end{aligned}
\end{equation*}
$\bullet $ $\bB^{(0)}_{\geq N, j+1}(\zak)$ is the real-to-real matrix of symbols in $\Gamma_{\geq N}^0[r]$ defined by 
$$
   \Opbw{ \bB^{(0)}_{\geq N, j+1}(\zak)}:= 
    \Opbw{ \bB^{(0)}_{\geq N, j}(\zak)} +
    \zOpbw{\tb_{\geq N, j+1}^{(-\frac{j}{2})}}+
\Opbw{\breve{{\bB}}^{(0)}_{\geq N, j}(\zak)}   + 
    \cP_{\geq N} \big[ 
 \bP_j^{new}(\zak)
\big]
$$
$\bullet$ $\bR_{1, j+1}(\zak) $ is a
real-to-real matrix of smoothing operators in $ \Sigma\mathcal{R}^{-\vr}_{1}[r,N]$ collecting all the smoothing remainders. 

Note that  the right hand side of \eqref{lem.iter.new.sym} is linearly Hamiltonian, hence the symbol
$\wt\td^{(-\frac12)}_{\geq 1, j+1} $ is real-valued whereas 
$\wt\tb_{\geq 1,j+1}^{(-\frac{j+1}{2})}$ is even in $\xi$, 
proving the structure \eqref{sist1 j-th} at step $j+1$.

By Lemma \ref{lem:flow.ad}$-(ii)$,  the map ${\bf \Phi}_{j}^{\theta}(\zak)$ is a $(0,0)-$admissible transformation with arbitrary  gain $\wt\vr \geq 0$.
\end{proof}

\begin{proof}[{\bf Proof of Proposition \ref{diag.ord0}}]
We define 
\begin{equation}\label{finalPsi}
\mathbf{\Psi}_0:={\bf \Phi}_{2\vr -1}(\zak) \circ\cdots 
\circ{\bf \Phi}_{0}(\zak)
\end{equation}
where 
the maps ${\bf \Phi}_{j}(\zak)$, $j=0,1,\ldots,  2\vr-1 $
are defined in \Cref{lem:indud}.  
 The map $\mathbf{\Psi}_0(\zak)$ is
 a  $(0,0)-$ admissible transformation with arbitrary gain $\wt\vr>1$ as well as 
each ${\bf \Phi}_{j}(\zak)$ thanks to  \Cref{lem:comp}.
Lemma \ref{lem:indud} implies that 
if $U_0$ solves \eqref{NuovoParaprod} then 
 the function  $ U_1 := W_{ 2 \vr} = \mathbf{\Psi}_0(\zak)U_0 $
solves  system \eqref{sist1 j-th} with $j= 2\vr  $
which is  \eqref{eq.diago.ord0} with $ \wt\td_{\geq 1}^{(-\frac12)} := \wt\td_{1,  2\vr}^{(-\frac12)}$, ${\bB}^{(0)}_{\geq N}:= {\bB}^{(0)}_{\geq N, 2\vr}$
and
\[ 
\bR_{\geq 1}(\zak) := \zOpbw{\wt\tb_{\geq 1,  2\vr}^{(-\vr)}(\zak; x, \xi)}
+ \bR_{\geq 1,  2 \vr }(\zak) \, ,\qquad  \wt\tb_{\geq 1,  2\vr }^{(-\vr)}\in \Sigma_1^{N-1}\wt\Gamma^{-\vr} \, , 
\]
which is a {real-to-real matrix} of smoothing operator in $ \Sigma {\mathcal R}^{-\vr}_{1}[r,N]$.
\end{proof}
\subsubsection{Removal of diagonal symbols of homogeneity $1$}
The goal of this section is to remove all the components of homogeneity $1$ from the diagonal symbol 
$ - \ii \frac{\gamma}{2} \sign(\xi) - \ii \sfV_\gamma(\zak;x)\xi  - \ii \lambda(\zak;x) \omega(\xi) -\ii \tilde \td^{(-\frac12)}_{ \geq 1}(\zak; x, \xi) $
of the paradifferential operator in the first line of \eqref{eq.diago.ord0}. 
This is done in three steps:
\begin{enumerate}
    \item In Lemma \ref{ridlinord1}, via a paracomposition transformation,  we remove the terms of homogeneity $1$ from the transport term
    $-\im \sfV_\gamma(\zak; x) \xi$ and identify the new quadratic transport term  $\im \check{\sfV}_{2}(\zak;x)\xi$ that appears in equation \eqref{BNF1}--\eqref{def:D}. 
The algebraic procedure that allows to achieve this -- and compute explicitly the real-valued symbol $\check{\sfV}_{2}(\zak;x)$-- is achieved in \Cref{lem:homo.eq.trans} below.
    \item In Lemma \ref{lem:rid12} we remove the terms of homogeneity $1$ from the symbols of order $\frac12$ and $0$ appearing in the diagonal component; we prove that, at these orders, the symbols still have a product structure as
    $$
    \tq_{\geq 2}(\zak; x) \omega(\xi) \quad \mbox{ and } \quad \tb_{\geq 2}(\zak; x) \sign(\xi)
    $$
    giving rise to the terms of homogeneity greater than $2$ with the same notation in equation \eqref{def:D};
    \item Finally, in Lemma \ref{riduzione.negativi}, we remove the terms of homogeneity $1$ from the symbols of strictly negative order  appearing in the diagonal component; the remaining real-valued valued quadratic symbols give rise to the contribution $\tilde \tc_{\geq 2}^{(-\frac12)}(\zak; x, \xi)$ in \eqref{def:D}. 
\end{enumerate}

\paragraph{Normal form of transport term of homogeneity $1$.}
The goal of this section is to eliminate the homogeneity 1  component of the transport term $\vOpbw{-\im \sfV_\gamma(\zak; x) \xi}$ and identify the new quadratic component.
The  algebraic  equation allowing this procedure is analyzed in the following lemma.

\begin{lemma}\label{lem:homo.eq.trans}
Let $N\in \N$, $N \geq 2$. There are  $s_0, r>0$ such that the following holds true.
Let $(\eta,\psi)(t) \in B_{s_0, \R}(I;r)$ be a solution of \eqref{eq:etapsi}, 
and define $\zak=\zak(t)$ by \eqref{zak}.
Let $\sfV_\gamma$ be the real-valued function defined in \eqref{def:Vgamma}. 
 There is a real-valued homogeneous function $\beta \in \wt\cF^\R_1$ such that 
    \begin{equation}\label{v.gamma.1.2}
\beta_t - \sfV_\gamma  - \beta (\sfV_\gamma)_x + \sfV_\gamma \beta_x - \beta_x \beta_t = \check{\sfV}_2 + {\sfV}_{\geq 3}'
\end{equation}
  where $\check{\sfV}_2(\zak; x)$ is the real-valued function in $\wt \cF^\R_2$ defined in \eqref{V2.exp0}--\eqref{V2.exp1} and ${\sfV}_{\geq 3}'(\zak; x)  $ is a real-valued function in $ \Sigma \cF^{\R}_{3}[r, N] $.
\end{lemma}

\begin{proof}
 Recalling \eqref{esp.v.gamma.1.2}, we write 
 \be\label{Vg.exp}
 \sfV_{\gamma} = \sfV_{\gamma,1} + \sfV_{\gamma,2} + \sfV_{\gamma, \geq 3} \,, \quad  \mbox{ where } \ \
 \sfV_{\gamma,1}:= \psi_x - \gamma \eta \,, \ \ \ \sfV_{\gamma,2}:= - \eta_x |D| \psi \,. 
 \ee
   The function    $\beta \in \wt\cF^\R_1$ is then  chosen to remove the term of homogeneity exactly one from \eqref{v.gamma.1.2}, namely $\sfV_{\gamma,1}$. We now prove that in $(\eta, \psi)$ coordinates such $\beta$ is\footnote{Recall that $\Pi_0^\perp$ is the  $ L^2 $-projector onto the space of functions with zero average.}
    \begin{align} \label{def:beta}
   \beta(\eta, \psi; x):= 
    ( A \eta + \gamma \Pi_0^\bot \psi)   \,, 
    \quad A := -\left(\Hilb + \gamma^2 \partial_x^{-1 }  {\Pi_0^\bot}\right) 
   \,.
    \end{align}
     To verify that such $\beta$  solves \eqref{v.gamma.1.2}, we first compute $\pa_t \beta$.
    By using \eqref{DNexp}-\eqref{DNexp2}  and  the identities in \eqref{identities.fou}, we have that the 
     equations of motion \eqref{eq:etapsi}  Taylor expand  as\footnote{Recall that by \eqref{zero.av} $\Pi_0^\bot \eta = \eta$.} 
   \begin{align*}
       \eta_t &= |D| \psi - \partial_x \circ \eta \circ \partial_x \psi - |D| \circ \eta \circ |D| \psi + \gamma \eta \eta_x + F^{(\eta)}_{\geq 3}(\eta, \psi)\,,\\
       {\Pi_0^\bot} \psi_t& = -  \eta - \frac 1 2 {\Pi_0^\bot} \psi_x^2 + \frac 1 2 {\Pi_0^\bot}(|D|\psi)^2 + \gamma \Hilb \psi - \gamma \Hilb \circ \eta \circ |D| \psi + F^{(\psi)}_{\geq 3}(\eta, \psi)\,,
   \end{align*}
   with $F^{(\eta)}_{\geq 3}$ and $F^{(\psi)}_{\geq 3}$ real-valued functions in $\Sigma \cF^\R_{3}[r, N]$.
   Then a direct computation gives
    \begin{align}
    \beta_t (\eta, \psi; \cdot) &= A \eta_t + \gamma {\Pi_0^\bot} \psi_t 
    = \und{\beta}_{\,1} + \und{\beta}_{\, 2}+ F^{(\beta)}_{\geq 3}\,,
    \label{betat}
   \end{align}
   where
   \begin{align}
    \und{\beta}_{\, 1} &:= (A |D| + \gamma^2\Hilb) \psi - \gamma \Pi_0^\bot \eta \stackrel{\eqref{def:beta},\eqref{identities.fou}}{=} \psi_x - \gamma \eta \stackrel{\eqref{Vg.exp}}{=}
    \sfV_{\gamma, 1} \ , \label{def.und1}\\
   \label{def.und2}
       \und{\beta}_{\, 2} &:= - A (\partial_x \circ \eta \circ \partial_x \psi + |D| \circ \eta \circ |D| \psi ) + \gamma A (\eta \eta_x) -\frac{\gamma}{2 }{\Pi_0^\bot} \psi_x^2 +  \frac{\gamma}{2 }{\Pi_0^\bot} (|D|\psi)^2  - \gamma^{2}  \Hilb \circ \eta \circ |D|\psi\,
   \end{align}
   and $F^{(\beta)}_{\geq 3}$ is in $\Sigma \cF_3^\R[r, N]$. 
   Combining \eqref{betat} and \eqref{Vg.exp}, then equation \eqref{v.gamma.1.2} is fulfilled with   
\begin{gather}\label{undV2.def}
    \check{\sfV}_{2} := \und{\beta}_{\, 2} - \sfV_{\gamma, 2} - \beta(\sfV_{\gamma, 1})_x \in \wt \cF^\R_2 \,,\\
    {\sfV}_{\geq 3}' := -\sfV_{\gamma, \geq 3} + F^{(\beta)}_{\geq 3} - \beta\left( \sfV_{\gamma, 2} + \sfV_{\gamma, \geq 3}\right)_x + \beta_x \left( \sfV_{\gamma, 2} + \sfV_{\gamma, \geq 3} - \und{\beta}_{\, 2} - F^{(\beta)}_{\geq 3}\right) \in \Sigma \cF^\R_{3}[r,N]\,.\nonumber
\end{gather}
   \smallskip
\noindent{\sc Expression of $\check{\sfV}_2$ in terms of the $\zak$ variables.} To conclude we compute   $\check{\sfV}_{2}$  in \eqref{undV2.def} in terms of the $\zak$ variables, proving that it Fourier expands as in \eqref{V2.exp0}, \eqref{V2.exp1}.
First, substituting the expressions of $\und{\beta}_2$ and $\sfV_{\gamma,i}$, $i=1,2$, in \eqref{def.und2} respectively \eqref{Vg.exp} in the one for $\check{\sfV}_2$ in \eqref{undV2.def}, we obtain
\begin{align}
\label{tV2}
{\check{\sfV}}_2(\eta, \psi; \cdot) &:= - A (\partial_x \circ \eta \circ \partial_x \psi) - A ( |D| \circ \eta \circ |D| \psi) + \gamma A (\eta \eta_x) - \frac{\gamma}{2 } {\Pi_0^\bot } \psi_x^2 + \frac{\gamma}{2 }{\Pi_0^\bot}(|D| \psi)^2 \\
\notag
& \quad 
- \gamma^2 \Hilb \circ \eta \circ |D| \circ \psi
+  \eta_x |D| \psi
- (A \eta + \gamma {\Pi_0^\bot} \psi)(\psi_{xx} - \gamma \eta_x)\,. 
\end{align} 
We now write  each term in  Fourier series expanding
 $\eta(x) = \sum_{j\neq 0} \eta_j e^{\im j x}$, $\psi(x) = \sum_{j\neq 0} \psi_j e^{\im j x}$. 
First we have (recalling $\sign(0) = 0)$
\begin{align*} 
& \eta_x |D|\psi 
  -\gamma^2 \Hilb \circ \eta \circ |D| \circ \psi 
  = \sum_{j_1, j_2 \in \Z_*}\im  |j_1| \,  \left(  j_2 
 + \gamma^2     \sign(j_1 + j_2) \right) \, 
 \psi_{j_1} \eta_{j_2} e^{\im (j_1 + j_2) x}  \ , \\
& -\frac{\gamma}{2} {\Pi_0^\bot} \psi_x^2
 + \frac{\gamma}{2} {\Pi_0^\bot}(|D| \psi)^2 = 
  \sum_{\substack{j_1, j_2 \in \Z_* \\ j_1+j_2\not=0}} 
      \frac{\gamma}{2}\left(   j_1 j_2
      +  |j_1| |j_2| \right) 
      \psi_{j_1} \psi_{j_2} 
        e^{\im (j_1 + j_2) x}\,. 
 \end{align*}
Moreover, using  \eqref{def:beta} and \eqref{identities.fou}, 
 \begin{align*}
 \notag
  -A(\partial_x \circ & \ \eta \circ \partial_x \psi) 
  - A (|D| \circ \eta \circ |D| \psi) = 
    \left(|D| + \gamma^2\,{\Pi_0^\perp}\right)(\eta \psi_x)
  + \left(-\partial_x + \gamma^2\Hilb \right) (\eta |D| \psi)
  \\
        &= \sum_{\substack{j_1, j_2 \in \Z_* \\  {j_1 + j_2 \neq 0}} }
        \im \Big( 
        j_1 \left(|j_1 + j_2| + \gamma^2 \right)
         - |j_1| \left(  (j_1 + j_2)  + \gamma^2 \operatorname{sign}(j_1 + j_2)\right)
        \Big)
        \psi_{j_1} \eta_{j_2} e^{\im (j_1 + j_2) x} \ , 
\end{align*}
whereas, denoting by $\delta(\cdot)$ the Kronecker delta in \eqref{delta},
        \begin{align}
        \notag
        \gamma A (\eta \eta_x) 
        &= \sum_{j_1, j_2 \in \Z_*} -\gamma j_2 \left(\sign(j_1 + j_2) + \frac{\gamma^2 {(1-\delta(j_1+j_2))}}{(j_1 + j_2)}  \right)\eta_{j_1} \eta_{j_2} e^{\im (j_1 + j_2)x}\,.
\end{align}
Finally, we have
        \begin{equation*}
            \begin{aligned}
            &-\Big( A \eta + \gamma {\Pi_0^\bot} \psi\Big)(\psi_{xx} -\gamma \eta_x)
            = \sum_{j_1, j_2 \in \Z_*} \im \left(  j_1^2\left(\sign(j_2) +\frac{\gamma^2}{ j_2}\right) + \gamma^2 j_2\right) \psi_{j_1} \eta_{j_2} e^{\im(j_1 + j_2) x}\\
            & \quad +
            \sum_{j_1, j_2 \in \Z_*} \gamma j_2^2 \psi_{j_1} \psi_{j_2} e^{\im (j_1 + j_2)x} +
            \sum_{j_1, j_2 \in \Z_*} - \gamma j_2 \left( \sign(j_1) + \frac{\gamma^2}{ j_1}\right) \eta_{j_1} \eta_{j_2} e^{\im (j_1 + j_2) x}\,.
        \end{aligned}
        \end{equation*}
 Collecting all terms we get 
       \begin{equation}
           \label{puo.darsi}
          {\check{\sfV}}_2(\eta, \psi; x) = \sum_{j_1, j_2 \in \Z_*} \left(\tA_{j_1, j_2} \eta_{j_1} \eta_{j_2}  + \im \tB_{j_1, j_2} \psi_{j_1} \eta_{j_2} + \tC_{j_1, j_2} \psi_{j_1} \psi_{j_2} \right) e^{\im(j_1 + j_2)x}\,,
       \end{equation}
       where 
       \be\label{puo.darsi3}
       \begin{aligned}
      &  \tA_{j_1, j_2} := -\gamma j_2 \left( \operatorname{sign}(j_1 + j_2) + \frac{\gamma^2 {(1-\delta(j_1 + j_2))}}{(j_1 + j_2)} +\operatorname{sign}(j_1) + \frac{\gamma^2}{ j_1}\right)\,,\\
       & \tB_{j_1, j_2} := \left(  j_1 (|j_1 + j_2| - |j_1| + j_1 \operatorname{sign}(j_2))+ \gamma^2 \frac{j_1^2}{j_2}
        + \gamma^2(j_1 + j_2) 
        { -  \gamma^2 j_1 \delta(j_1+j_2) }
        \right)\,,\\
       & \tC_{j_1, j_2} := \frac{\gamma}{2} \left( j_1 j_2  + |j_1 j_2|+ 2 j_2^2\right)\,,
        \end{aligned}
        \ee
        and $\delta(j)$ the Kronecker delta in \eqref{delta}.
   As a final step, recall that $\vect{\eta}{\psi} = \Lmap \vect{\zetina}{\bar \zetina}$ (cf. \eqref{zak}), so using the formula for $\Lmap$ in the second of \eqref{app_L} and the one of $M(D)$ in \eqref{defMD} we obtain 
    \begin{equation}\label{epsi.fou}
    \begin{gathered}
        \eta_{j} = \sum_{\sigma \in \{\pm 1\}} \fm_{j} \zetina_{\sigma j}^\sigma\,, \quad 
        \psi_j = \sum_{\sigma \in \{\pm 1\}} \im \left(-\sigma \fn_j - \gamma\frac{\fm_j}{2j} \right) \zetina_{\sigma j}^\sigma\,, \qquad j \in \Z_*\,,
    \end{gathered}
    \end{equation}
       with $\fm_j, \fn_j$  in \eqref{al.be.m}. Then formulas \eqref{V2.exp0}--\eqref{V2.exp1} follow combining \eqref{puo.darsi}--\eqref{epsi.fou},   renaming $\sigma_1 j_1 =:k_1\,, \ \sigma_2 j_2 =:k_2$, and putting $\fe_j:= \fn_j + \gamma \frac{\fm_j}{2j} = \frac12\fm_j^{-1} + \gamma \frac{\fm_j}{2j} $.
\end{proof}
We now perform a paracomposition in order to remove the term $\sfV_{\gamma,1}(\zak;x)$ of homogeneity one. More precisely, we 
 prove the following result.
\begin{lemma}[{\bf Normal form  of  degree one homogeneous components in the transport term}]\label{ridlinord1}
Let $N \in \N$, $\vr> N\geq 2$ and $\wt\vr>3$. 
There exist $s_0, r  >0 $ and a $(0,3)$- admissible transformation $ {\bf \Psi}_1(\zak)$ with gain $\wt \vr$ (recall Definition \ref{admtra}) such that if $(\eta, \psi)(t) \in B_{s_0, \R}(I;r)$ solves \eqref{eq:etapsi}, and $\zak=\zak(t)$ is defined in \eqref{zak}, then the variable $U_2:= {\bf \Psi}_1(\zak)U_1$, with $U_1$ solving \eqref{eq.diago.ord0}, solves
\begin{equation}\label{brand.new.w}
\begin{aligned}
    \partial_t U_2  =  & \vOpbw{- \im 
    \frac{\gamma}{2}\sign(\xi) +
    \im \check{\sfV}_{2}(\zak;x) \xi + \im \check\sfV_{\geq 3}(\zak; x)\xi -\im \big((1+\tq_{\geq 1}(\zak;x)\big)\omega(\xi)+\ii \wt\tc_{\geq1}^{(-\frac12)}(\zak; x, \xi)}U_2\\
    &+\Opbw{ 
    \bB^{(0)}_{\geq N}(\zak; x, \xi) } U_2 + \bR_{\geq 1}(\zak) U_2\,,
\end{aligned}
\end{equation}
where $\omega(\xi)$ is defined in \eqref{omegaxi} and:\\ 
$\bullet$ $\check{\sfV}_{2}(\zak; x) $ is the real-valued function in $\wt \cF^\R_2$ defined in \eqref{V2.exp0}--\eqref{V2.exp1} 
and   $\check\sfV_{\geq 3}(\zak;x) $  a real-valued function in $\Sigma \cF^\R_{3}[r, N]$;\\
$\bullet$  
$\tq_{\geq 1}(\zak;x)$ is a  real-valued function in $\Sigma \cF^\R_{1}[r, N]$ and  $\wt\tc_{\geq1}^{(-\frac12)}(\zak; x, \xi)$ is a pluri-homogeneous real-valued symbol in $ \Sigma_1^{N-1}\wt \Gamma^{-\frac12}$;
\\
$\bullet$ $\bB_{\geq N}^{(0)}(\zak; x, \xi)$ is a real-to-real matrix of symbols in $ \Gamma_{\geq N}^0[r]$ and   $\bR_{\geq 1}(\zak) $ is a real-to-real matrix of smoothing remainders in $\Sigma \cR^{-\vr + {N}}_{1}[r, N]$.
\end{lemma}
\begin{proof}
We define the map $ {\bf \Psi}_1(\zak)$ as the time-1 flow 
$ {\bf \Psi}_1^\tau(\zak)\vert_{\tau = 1}$ of 
$$
\pa_\tau {\bf \Psi}_1^{\tau}(\zak)=
  \, \vOpbw{\ii \frac{\beta(\zak;x)}{1+\tau\beta_x(\zak;x)} \, \xi}
 {\bf \Psi}_1^{\tau}(\zak) , \qquad 
{ \bf \Psi}_1^0(\zak)={\rm Id} \, ,    
$$
and the real-valued function $\beta \in \widetilde\mF^\R_{1}$ is defined in \eqref{def:beta}.
By Lemma \ref{lem:flow.ad} $(i)$, 
 the map ${\bf \Psi}_1(\zak)$ is a $(0, 3)$- admissible transformation of arbitrary gain $ \tilde \vr > 3$.
We compute the equation solved by $U_2:= {\bf \Psi}_1(\zak)U_1$. As 
 $U_1$ solves equation \eqref{eq.diago.ord0}, we get
\begin{equation*}
    \partial_t U_2 = {\bf \Psi}_1(\zak)\Big(\bY(\zak)+ \Opbw{\bB_{\geq N}^{(0)}} + \bR_{\geq 1}(\zak) \Big)  {\bf \Psi}_1(\zak)^{-1} U_2 + (\partial_t {\bf \Psi}_1(\zak)) {\bf \Psi}_1(\zak)^{-1} U_2\,,
\end{equation*}
where
\begin{equation*}
    \bY(\zak):=   \vOpbw{-\ii \frac{\gamma}{2}  \sign(\xi)  - \ii \sfV_\gamma(\zak;x)\xi  - \ii \lambda(\zak;x) \omega(\xi) -\ii \tilde \td^{(-\frac12)}_{ \geq 1}(\zak; x, \xi) }\,.
\end{equation*}
We now study how each term is transformed. We start with $(\partial_t {\bf \Psi}_1(\zak)) {\bf \Psi}_1(\zak)^{-1}$. 
By Proposition \ref{prop:Egorov}-3 (with $p =1$) one has
$$
\big(\pa_t {\bf \Psi}_1(\zak)\big) {\bf \Psi}_1 (\zak)^{-1} =   \vOpbw{ \ii(\beta_t  - \beta_x \beta_t+ g_{\geq 3})\, \xi } + \bR_{\geq 2}(\zak)
$$
where $ g_{\geq 3}$ is a real-valued function in $ \Sigma \cF^{\R}_{3}[r, N]$  and 
$ \bR_{\geq 2}(\zak) $ is a real-to-real matrix of smoothing operators in  $  \Sigma\mathcal{R}^{-\vr}_{2}[r, N] $.

\smallskip

Next  consider   ${\bf \Psi}_1(\zak) \vOpbw{ -\ii \sfV_{\gamma} \,\x} {\bf \Psi}_1(\zak)^{-1}$. 
By applying Proposition  \ref{prop:Egorov}-1  (again with $p=1$)  with the real-valued symbol $ a^{(1)}_{\geq 1} = - \sfV_\gamma \xi $ one has (use formula \eqref{trans.transp})
\begin{equation}\label{conj.ord1.1}
{\bf \Psi}_1(\zak) \vOpbw{ -\ii \sfV_{\gamma} \,\x} {\bf \Psi}_1(\zak)^{-1} = 
  \vOpbw{ -\ii \big(  \sfV_{\gamma}  + \beta (\sfV_{\gamma})_x - \sfV_{\gamma} \beta_x +  h_{\geq 3}  \big) \xi }
  + \bR_{\geq 2}(\zak) 
\end{equation}
where $ h_{\geq 3} \in \Sigma \cF^{\R}_{ 3}[r, N]$ 
and $ \bR_{\geq 2}(\zak) $ is a real-to-real matrix of smoothing operators in  $ \Sigma {\mathcal R}^{-\vr+1}_{2} [r, N]$.

We now compute the remaining terms. 
By \cite[Lemma 3.21]{BD}, the   diffeomorphism of $ \T $ given by 
	$\Psi \pare{ \zak ;x } :=  x + \beta\pare{\zak; x}$
	is invertible and its  inverse has the form 
	$ \Psi^{-1}\pare{ \zak;y }  := y + \breve{\beta}\pare{\zak;y} $
	with $ \breve{\beta}\in  \Sigma \cF_{ 1}^\R\bra{r, N} $ (cf. also \eqref{eq:diffeo}).\\
So consider ${\bf \Psi}_1(\zak) \vOpbw{- \im \sign(\xi)} {\bf \Psi}_1(\zak)^{-1}$; by applying Proposition \ref{prop:Egorov}-1 with $ a_{\geq 0}^{(0)} = - \, \sign(\xi)$ one has 
\begin{equation}\label{conj.ord1.2}
\begin{aligned}
{\bf \Psi}_1(\zak) \vOpbw{- \im \, \sign(\xi)} {\bf \Psi}_1(\zak)^{-1} &= \vOpbw{-\ii \, \sign{(\xi)}\, \sign{\big({1+\breve{\beta}_y}_{|_{y=\Psi(\zak;x)}}\big)} + \ii g_{\geq 1}^{(-2)}} + \bR_{\geq 1}(\zak)
\\ 
&=\vOpbw{- \im \sign(\xi)} + \vOpbw{\ii g_{\geq 1}^{(-2)}} + \bR_{\geq 1}(\zak)
\end{aligned} 
\end{equation}
with $ g_{\geq 1}^{(-2)} \in \Sigma\Gamma^{-2}_{1}[r, N]$ a real-valued symbol 
and $ \bR_{\geq 1}(\zak)$ a  real-to-real matrix of smoothing operators in $\Sigma {\mathcal R}^{-\vr}_{1} [r, N]$, where in the last line we used that $\sign{({1+\breve{\beta}_y}_{|_{y=\Psi(\zak;x)}})} =1$ for $\beta$ small enough.\\
Next we consider ${\bf \Psi}_1(\zak) \vOpbw{-\ii \lambda(\zak)\omega(\xi) } {\bf \Psi}_1(\zak)^{-1}$. 
Writing $\lambda(\zak) = 1 + \lambda_{\geq 1}(\zak)$ with $\lambda_{\geq 1} \in \Sigma \cF^\R_1[r, N]$ (see \eqref{diag2})
and applying Proposition \ref{prop:Egorov}-1 with $ a^{(\frac12)}_{\geq 0} = - (1+ \lambda_{\geq 1}(\zak))\omega(\xi)$ one has 
\begin{equation}\label{conj.ord1.3}
\begin{aligned}
{\bf \Psi}_1(\zak) &\vOpbw{-\ii \big( 1 + \lambda_{\geq 1}\big) \omega(\xi) }{\bf \Psi}_1(\zak)^{-1} \\
&  = \vOpbw{
-\ii \big(1 + \lambda_{\geq 1}(\zak;y) \big) \, \,  
\omega\big(  ( 1+\breve\beta_y(\zak;y))\xi\big)\vert_{y = \Psi(\zak;x)} 
}+ \vOpbw{\im  {g}_{\geq 1 }^{(-\frac32)} } + \bR_{\geq 1}(\zak) 
\end{aligned}
\end{equation}
where  $  g^{(-\frac32)}_{\geq 1}$ is a real-valued symbol in $ \Sigma\Gamma^{-\frac32}_{ 1}[r, N]$ 
and $ \bR_{\geq 1}(\zak) $ a real-to-real matrix of smoothing operators in  $ \Sigma {\mathcal R}^{-\vr+\frac12}_{1}[r, N] $.
Arguing as in   \cite[Lemma 3.23]{BD}, one proves that 
$\omega\big(\xi (1+ \breve\beta_y(\zak;y))\big)$
is a symbol in $\Sigma \Gamma^{\frac12}_{0}[r, N]$. 
Then, using that $\omega^{(-\frac12)}(\xi):= \omega(\xi) - |\xi|^{\frac12} $
  is a real-valued Fourier multiplier in $ \wt \Gamma^{-\frac12}_0 $ we get 
 \begin{align}\notag
  \omega\Big(\xi (1+ \breve\beta_y(\zak;y)) \Big) 
  & = \big| 1+ \breve\beta_y(\zak;y) \big|^{\frac{1}{2}} \omega(\xi) + a_{\geq 1}^{(-\frac12)}(\zak; y, \xi)
 \end{align}
 where  
 \begin{align*}
a_{\geq 1}^{(-\frac12)}(\zak; y, \xi)
& = \omega^{(-\frac12)}\big(\xi(1+\breve\beta_y(\zak;y))\big)  - \omega^{(-\frac12)}\big(\xi\big)  + \big( 1-\big| 1+ \breve\beta_y(\zak;y) \big|^{\frac{1}{2}} \big)\omega^{(-\frac12)}(\xi)\
 \end{align*}
 is a real-valued symbol  in 
$ \Sigma\Gamma^{-\frac12}_{ 1}[r, N]$.
Then,   we get 
\begin{equation*}
{\bf \Psi}_1(\zak) \vOpbw{-\ii \big( 1 + \lambda_{\geq 1}\big) \omega(\xi) }{\bf \Psi}_1(\zak)^{-1} = 
\vOpbw{
-\ii \big(1 + \tq_{\geq 1} \big) \,  \omega(\xi) 
+\im {g}_{\geq 1 }^{(-\frac12)} } + \bR_{\geq 1}(\zak) 
\end{equation*}
where 
$$\tq_{\geq 1}(\zak;x) :=
\Big[\big( 1 + \lambda_{\geq 1}(\zak; y)\big) \big| 1+ \breve\beta_y(\zak;y) \big|^{\frac{1}{2}} - 1\Big]_{\vert_{y= x + \beta(\zak;x)}}
$$
is a  real-valued function in $\Sigma \cF^\R_{1}[r, N]$,   $  g^{(-\frac12)}_{\geq 1}$ is a real-valued symbol in $ \Sigma\Gamma^{-\frac12}_{ 1}[r, N]$ 
and $ \bR_{\geq 1}(\zak) $ a real-to-real matrix of smoothing operators in $ \Sigma {\mathcal R}^{-\vr+\frac12}_{1}[r, N] $.
Next, reasoning as above
\begin{equation}\label{conj.ord1.4}
\begin{aligned}
  &   {\bf \Psi}_1(\zak) \vOpbw{\ii \tilde \td^{(-\frac12)}_{ \geq 1}  }{\bf \Psi}_1(\zak)^{-1} =
    \vOpbw{ \ii \und{\td}^{(-\frac12)}_{ \geq 1}}+ \bR_{\geq 2}(\zak)   \\
 &    {\bf \Psi}_1(\zak) \Opbw{\bB^{(0)}_{\geq N}(\zak) }{\bf \Psi}_1(\zak)^{-1}= \Opbw{\hat{\bB}^{(0)}_{\geq N}(\zak) } + \bR_{\geq N}(\zak) 
\end{aligned}
\end{equation}
where $ \und{\td}^{(-\frac12)}_{ \geq 1}$ is a real-valued symbol in $\Sigma \Gamma_1^{-\frac12}[r,N]$, 
$\hat \bB^{(0)}_{\geq N}(\zak)$ is a real-to-real matrix of paradifferential operators with symbols in $\Gamma_{\geq N}^{0}[r]$, $ \bR_{\geq 2}(\zak)$ and $ \bR_{\geq N}(\zak)$ are real-to-real matrices of smoothing operators in in $  \Sigma {\mathcal R}^{-\vr}_{2}[r, N] $  respectively  $   {\mathcal R}^{-\vr}_{\geq N}[r] $.
Finally, by Proposition \ref{prop:Egorov}-2 the conjugated of the operator $\bR_{\geq 1}(\zak)$ in \eqref{eq.diago.ord0} is a real-to-real matrix of smoothing operators in 
$ \Sigma\mathcal{R}^{-\vr+N}_{1}[r, N]$.

    Therefore, by \eqref{conj.ord1.1}-\eqref{conj.ord1.4} we obtain
    \begin{align}
    \partial_t U_2 &=  \textup{Op}_{\tt vec}^{{\scriptscriptstyle{\mathrm BW}}}\Big(\im \big(\underbrace{ \beta_t- \sfV_{\gamma}  - \beta (\sfV_{\gamma})_x + \sfV_{\gamma} \beta_x  - \beta_x \beta_t }_{=  \check{\sfV}_2 + {\sfV}_{\geq 3}' \mbox{  by Lemma } \ref{lem:homo.eq.trans}}+g_{\geq 3} - h_{\geq 3}\big)\,  \xi\Big) U_2 \label{no.v1}\\
    &+ \vOpbw{-\ii \frac{\gamma}{2} \sign(\xi)-\ii \big(1 + \tq_{\geq 1} \big) \,  \omega(\xi) + \ii \wt\tc_{\geq1}^{(-\frac12)} } U_2 + \Opbw{\hat{\und{\bB}}^{(0)}_{\geq N}} U_2 + \bR_{\geq 1}(\zak) U_2\,,\nonumber
    \end{align}
    where
 $\wt\tc_{\geq1}^{(-\frac12)}:=\cP_{\leq N-1}\left[
 \und{\td}^{(-\frac12)}_{ \geq 1}+
 \frac \gamma 2 g_{\geq 1}^{(-2)}+ g^{(-\frac12)}_{\geq 1} \right]$
 is a pluri-homogeneous real-valued symbol in $\Sigma_1^{N-1}\wt\Gamma^{-\frac12}$, 
 $$
 \Opbw{\hat{\und{\bB}}^{(0)}_{\geq N}}:=\Opbw{\hat{{\bB}}^{(0)}_{\geq N}}+\cP_{\geq N}\left[
\vOpbw{\im g_{\geq 1}^{(-2)}+\im g^{(-\frac12)}_{\geq 1}+ \im \und{\td}^{(-\frac12)}_{ \geq 1}
}
\right]
 $$
 is a real-to-real matrix of paradifferential operators with symbols in
 $\Gamma_{\geq N}^{0}[r]$ 
 -- that to simplify notation we re-denote $\bB_{\geq N}^{(0)}(\zak)$ --
and $ \bR_{\geq 1}(\zak)$ is a real-to-real matrix of smoothing operators in  $ \Sigma {\mathcal R}^{-\vr+N}_{ 1}[r, N] $.
Equation \eqref{brand.new.w} follows with   $\check{\sfV}_{\geq 3}:= {\sfV}_{\geq 3}' + g_{\geq 3} - h_{\geq 3} $ which is a real-valued function in $\Sigma \cF^\R_3[r,N]$.
\end{proof}
\paragraph{Removal of symbols of non-negative order and  homogeneity $1$.} In this paragraph we remove the  degree $1$ homogeneous  component from the diagonal symbols in \eqref{brand.new.w} of order $\frac12$ and $0$. More precisely, we prove the following lemma.
\begin{lemma}[{\bf Normal form of degree 1 homogeneous components of order $\frac12$ and $0$}]
\label{lem:rid12}
Let $N\in\N$, $N\geq2$. Let $\vr>{2 N}$ and  $\wt \vr \geq 3$. 
There are $s_0, r  >0 $ and a $(0,\frac 3 2)$- admissible transformation $ {\bf \Psi}_2(\zak)$ with gain $\wt \vr$ (recall Definition \ref{admtra}) such that if $(\eta, \psi)(t) \in B_{s_0, \R}(I;r)$ solves \eqref{eq:etapsi}, and $\zak=\zak(t)$ is defined in \eqref{zak}, then the variable $U_3:= {\bf \Psi}_2(\zak)U_2$, with $U_2$ solving \eqref{brand.new.w}, solves
\begin{align}
    \partial_t U_3  = & {-\im  \bOmega(D) U_3 + \vOpbw{\im \check{\sfV}_{2}(\zak;x) \xi + \im \check{\sfV}_{\geq 3}(\zak; x)\xi + \im  {\tq}_{\geq 2}(\zak; x)\omega(\xi) + \ii \tb_{\geq 2}(\zak;x)\,\sign{(\xi)}+\ii \tilde \tc_{\geq1}^{(-\frac12)}(\zak; x, \xi)} U_3} \notag\\
    &+\bB_{\geq N}(\zak) U_3 + \bR_{\geq 1}(\zak) U_3\,,\label{eq.riduzione.ord12}
\end{align}
{where $\bOmega(D)$ is defined in \eqref{diaglin}-\eqref{relation},  
 $\omega(\xi) \in \wt\Gamma_0^\frac12$ in \eqref{omegaxi},  $\check{\sfV}_{2}\in \wt\cF^\R_2$ in 
\eqref{V2.exp0}--\eqref{V2.exp1},  
 $\check{\sfV}_{\geq 3}\in \Sigma \cF^\R_{\geq 3}[r,N]$  in Lemma \ref{ridlinord1} and}
\begin{itemize}
    \item $\tb_{\geq 2}(\zak;x)$, ${\tq}_{\geq 2}(\zak; x)$ are  real-valued functions in $ \Sigma \cF_{2}^\R[r, N]$ and $\tilde \tc_{\geq 1}^{(-\frac12)}(\zak;x, \xi)$ is a pluri-homogeneous real-valued symbol in $\Sigma_1^{N-1}\wt\Gamma^{-\frac12}$;
\item $\bB_{\geq N}(\zak)$ is a real-to-real matrix of spectrally localized maps in $\cS_{\geq N}^0[r]$ and  $\bR_{\geq 1}(\zak) $  a real-to-real matrix of smoothing remainders in $\Sigma \cR^{-\vr + {2N}}_{1}[r, N]$.
\end{itemize}
\end{lemma}
\begin{proof}
Let  $\tq_{\geq 1} (\zak; x)$ be the real-valued function in $  \Sigma \cF_1^\R[r, N]$ of \Cref{ridlinord1} and, recalling \eqref{pienne}, we denote by
$\tilde \tq_1:=\cP_1 [\tq_{\geq 1}]$ its real-valued 1-homogeneous component  in $\wt\cF^\R_1$. Then we  define the map $ {\bf \Psi}_2(\zak)$ as the time-1 flow ${\bf \Psi}_2(\zak):= {\bf \Psi}_2^{\tau}(\zak)|_{\tau=1}$ of 
$$
\pa_\tau {\bf \Psi}^{\tau}_2(\zak)=
  \, \vOpbw{ \ii\,  \wt g_1^{(\frac12)}(\zak;x,\xi)}
 {\bf \Psi}_2^{\tau}(\zak) , \qquad 
{\bf \Psi}_2^0(\zak)={\rm Id} \, ,   
$$
where $\tilde g_1^{(\frac12)}(\zak;x,\xi)$ is the real-valued $1-$homogeneous symbol 
\be\label{g112}
\begin{aligned}
  \tilde g_1^{(\frac12)}:=  &\tilde g_1^{(\frac12)}(\zak;x,\xi)  := \tg(\zak; x)\omega(\xi)+  \tf(\zak;x)\sign{(\xi)} \ , \\
    &\tg(\zak; x) := \tilde \tq_1((-\ii \vOmega(D))^{-1}\zak; x)= \sum_{\substack{j \in \Z_* \\ \sigma \in \{\pm 1\}} }\frac{\big(\tilde{\tq}_1\big)^\sigma_j}{-\ii \sigma \Omega_j(\gamma)} \zetina_j^\sigma e^{\ii \sigma j x}\,,\\
    &\tf(\zak; x) := -\frac12 \tg_x((-\ii \vOmega(D))^{-1}\zak; x)=-\frac12 \sum_{\substack{j \in \Z_* \\ \sigma \in \{\pm 1\}}} \frac{\big(\tg_x\big)^\sigma_j}{-\ii \sigma \Omega_j(\gamma)} \zetina_j^\sigma e^{\ii \sigma j x}
\end{aligned}
\ee
and $\Omega_j(\gamma)$  in \eqref{omegonejin}. 
The functions $\tg, \tf $ are well defined real-valued functions in $\wt \cF^\R_1$ as well as 
 $\tilde \tq_1 $ thanks to the  uniform lower bound \eqref{Omega>0} for $\Omega_j(\gamma)$.
Then  $\wt g_1^{(\frac12)}$ is a real-valued symbol in $\wt \Gamma_1^{\frac 1 2}$, and  by Lemma \ref{lem:flow.ad}--(iii),  the map ${\bf \Psi}_2(\zak)$ is a $(0, \frac 3 2)$- admissible transformation of arbitrary gain
{$\wt \vr \geq \frac32$}. 

We compute the equation solved by $U_3:= {\bf \Psi}_2(\zak)U_2$. As 
 $U_2$ solves equation \eqref{brand.new.w}, we get 
\begin{align*}
    \pa_t U_3 = & {\bf \Psi}_2(\zak)\Big(\bY(\zak)+ \Opbw{\bB_{\geq N}^{(0)}(\zak)} + \bR_{\geq 1}(\zak) \Big)  {\bf \Psi}_2^{-1}(\zak) U_3 + (\partial_t {\bf \Psi}_2(\zak)) {\bf \Psi}_2(\zak)^{-1} U_3 \ , \mbox{ where } \\
   &  \bY(\zak):= 
    \vOpbw{- \im 
    \frac{\gamma}{2}\sign(\xi) + 
    \im (\check{\sfV}_{2} +  \check{\sfV}_{\geq 3})\xi -\im \big((1+\tq_{\geq 1}\big)\omega(\xi)+\ii \wt\tc_{\geq1}^{(-\frac12)}(x,\xi)}
\end{align*}
Then by Proposition 
\ref{prop:FIO}-$(i)$ with $f=\tilde g_1^{(\frac12)}$, $p=1$, and $a = \frac{\gamma}{2}\sign(\xi)$
one has
\begin{equation}\label{mammolo}
    \begin{aligned}
        {\bf \Psi}_2(\zak) \vOpbw{- \im 
    \frac{\gamma}{2}\sign(\xi)} {\bf \Psi}_2(\zak)^{-1} 
        &=   \vOpbw{- \im  \frac{\gamma}{2}\sign(\xi) +\ii  a^{(-\frac 1 2)}_{\geq 1}}  + \bR_{\geq 1, 1}(\zak)\,,
    \end{aligned}
\end{equation}
with $a^{(-\frac 1 2)}_{\geq 1}$  a real-valued symbol in $\Sigma \Gamma_1^{- \frac 1 2}[r, N]$ and $ \bR_{\geq 1, 1}(\zak)$ a  real-to-real matrix of  smoothing operators in $\Sigma \cR_{ 1}^{-\vr}[r, N]$. Similarly, by applying  
\Cref{prop:FIO}-$(i)$, with $a=\check{\sfV}_2 + \check{\sfV}_{\geq 3}$,  one has
\begin{equation}\label{dottoceragia}
    \begin{aligned}
        {\bf \Psi}_2(\zak) &\,  \vOpbw{\ii (\check{\sfV}_2 + \check{\sfV}_{\geq 3}) \xi} \, {\bf \Psi}_2(\zak)^{-1}  = \vOpbw{\ii (\check{\sfV}_2 + \check{\sfV}_{\geq 3}) \xi +  \ii \{ \tilde g_1^{(\frac12)}\,,\  (\check{\sfV}_2 + \check{\sfV}_{\geq 3}) \xi\}}\\
        &+ \vOpbw{\frac{\ii}{2} \{ \tilde g_1^{(\frac12)}\,,\ \{ \tilde g_1^{(\frac12)}\,,\  (\check{\sfV}_2 + \check{\sfV}_{\geq 3}) \xi\} \} + \ii a^{(-\frac 1 2)}_{\geq 3} }  + \bR_{\geq 3}(\zak)\\
        &= \vOpbw{
        \ii \left((\check{\sfV}_2 + \check{\sfV}_{\geq 3}) \xi + q_{\geq 3} \omega(\xi) +  b_{\geq 3} \sign{(\xi)} +b_{\geq 3}^{(-\frac 1 2)} \right) } +  \bR_{\geq 3}(\zak)\,,
    \end{aligned}
\end{equation}
where\footnote{Here we used that $\xi\partial_\xi \omega(\xi)= \frac{1}{2}\omega(\xi)$ up to a Fourier multiplier in $\wt \Gamma^{-\frac12}_0$ and that $\omega(\xi)\partial_\xi \omega(\xi)=\frac12 \sign{\xi}$.}
\begin{align*}
& q_{\geq3} := q_{\geq 3}(\zak; x) := -(\check{\sfV}_2 + \check{\sfV}_{\geq 3}) \tg_x + \frac 1 2 \tg (\check{\sfV}_2 + \check{\sfV}_{\geq 3})_x \in \Sigma \cF^\R_{ 3}[r, N]\,,\\
& b_{\geq 3} := b_{\geq 3}(\zak; x) := \frac 1 4 (\tg (q_{\geq 3})_x - \tg_x  q_{\geq 3})- \tf_x (\check{\sfV}_2 + \check{\sfV}_{\geq 3}) \in \Sigma \cF^\R_{ 3}[r, N]\,,\\
& \qquad \quad a_{\geq 3}^{(-\frac 1 2)}, \, b_{\geq 3}^{(-\frac 1 2)} \in \Sigma \Gamma_{3}^{-\frac 1 2}[r, N] \text{ real-valued,}
\end{align*}
and $\bR_{\geq 3}(\zak)$ is a real-to-real matrix of smoothing operators in {$\Sigma \cR_{ 3}^{-\vr+1}[r, N]$}. 
Next, writing $\tq_{\geq 1} = \wt \tq_1 + \cP_{\geq 2}[\tq_{\geq 1}]$ and applying Proposition 
\ref{prop:FIO}-$(i)$ with $a=(1 + \tq_{\geq 1} ) \omega(\xi)$ ,
we  get
\begin{equation}\label{pisolo}
    \begin{aligned}
         {\bf \Psi}_2(\zak) &\,  \vOpbw{- \ii (1 + \tq_{\geq 1} ) \omega(\xi)} \, {\bf \Psi}_2(\zak)^{-1} \\
         & = 
         \vOpbw{-\ii (1 + \tq_{\geq 1})  \omega(\xi)
        +  \ii   \left\lbrace \tg\omega(\xi)\,,\ -(1 + \tq_{\geq 1} ) \omega(\xi)\right\rbrace
         +  \im c_{\geq 1}^{(-\frac 1 2)}
         } + \bR_{\geq 1,2}(\zak) \\
         &= \vOpbw{-\ii (1 +  
         {\wt\tq}_{1}+
         \cP_{\geq 2}[\tq_{\geq 1}]) \omega(\xi)+ \ii \big( 
         \frac{1}{2}\tg_x +
b_{\geq 2} \big)\sign{(\xi)}
+\ii c_{\geq 1}^{(-\frac 12)}
} + \bR_{\geq 1, 2}(\zak)
    \end{aligned}
\end{equation}
where 
$b_{\geq 2}:= b_{\geq 2}(\zak; x) = \frac{1}{2} \left( \tq_{\geq 1} \tg_x(\zak; x) - \tg \, (\tq_{\geq 1})_x\right)$
is a real-valued function in $\Sigma  \cF^\R_1[r,N]$,
$c_{\geq 1}^{(-\frac 1 2)}$ is a real-valued symbol in $\Sigma \Gamma_{1}^{-\frac 1 2}[r, N]$, and $\bR_{\geq 1, 2}(\zak) $ a real-to-real matrix of smoothing operators in 
$\Sigma \cR_{ 1}^{-\vr+\frac12}[r, N]$.
Analogously, by Proposition 
\ref{prop:FIO}-$(i)$ with $\tilde{\tc}_{\geq 1 }^{(-\frac 12)} $, one has
\begin{equation}\label{eolo}
\begin{aligned}
    {\bf \Psi}_2(\zak) \, \vOpbw{\ii \tilde{\tc}_{\geq 1 }^{(-\frac 12)} } \,  
    {\bf \Psi}_2(\zak)^{-1} &=  \vOpbw{\ii d_{\geq 1 }^{(-\frac 12)} } + \bR_{\geq 2, 1}(\zak)\,,
\end{aligned}
\end{equation}
with $d_{\geq 1 }^{(-\frac 12)}$ a real-valued symbol in $\Sigma \Gamma_{1 }^{-\frac 1 2}[r, N]$ and $\bR_{\geq 2, 1}(\zak) $ a real-to-real matrix of smoothing operators in $\Sigma \cR_{ 2}^{-\vr}[r, N]$.
Next, recalling the definition of spectrally localized map (see Definition \ref{def:specloc}) and that if $\bB_{\geq N}^{(0)}(\zak)$ is a real-to-real matrix of symbols in $\Gamma^0_{\geq N}[r]$ then it is also a real-to-real matrix of spectrally localized maps in $\cS_{\geq N}^0[r]$, by \Cref{prop:FIO}--3 one has
\begin{equation}\label{brontolo}
\begin{gathered}
{\bf \Psi}_2(\zak) \, \Opbw{\bB_{\geq N}^{(0)}(\zak) }\,  {\bf \Psi}_2(\zak)^{-1} 
= \check{\bB}_{\geq N}(\zak) \,,
\qquad 
{\bf \Psi}_2(\zak) \,  \bR_{\geq 1}(\zak) \,  {\bf \Psi}_2(\zak)^{-1} = \bR_{\geq 1, 3}(\zak)\,,
\end{gathered}
\end{equation}
with $\check{\bB}_{\geq N}(\zak)$ a real-to-real matrix of spectrally localized maps in $\cS_{\geq N}^0[r]$ and $\bR_{\geq 1,3}(\zak) $ a real-to-real matrix of smoothing operators in $\Sigma \cR^{-\vr+{2N}}_{ 1}[r, N]$. Finally,  Proposition \ref{prop:FIO}--4 one has
\begin{equation}\label{brontolo.parecchio}
\begin{aligned}
    \left(\pa_t {\bf \Psi}_2(\zak)\right) {\bf \Psi}_2 (\zak)^{-1}  =& \vOpbw{\ii  \tg(-\ii {\bf \Omega}(D) \zak ; x) \omega(\xi)+\ii \tf(-\ii {\bf \Omega}(D) \zak ; x) \sign(\xi)}\\
   & +\vOpbw{ \ii p_{\geq 2}  \omega(\xi) +\ii r_{\geq 2}\sign(\xi)+ \ii e_{\geq 2}^{(-\frac 1 2)}} + \bR_{\geq 2, 2}(\zak)\,,
\end{aligned}
\end{equation}
with 
\begin{equation*}
\begin{aligned}
    p_{\geq 2}:= p_{\geq 2} (\zak; x) &: = \tg(\cX_{\cH}(\zak); x) -  \tg(-\ii  {\bf \Omega}(D) \zak; x)\,,
    \\
    r_{\geq 2}:=r_{\geq 2} (\zak; x)&: = \tf(\cX_{\cH}(\zak); x) -  \tf(-\ii  {\bf \Omega}(D) \zak; x)+\frac12\{\tg(\zak; x) , \, \tg(\cX_{\cH}(\zak); x)  \}\,,
    \end{aligned}
\end{equation*} 
 both real-valued functions in $ \Sigma \cF_{ 2}^\R[r, N]$,  $e_{\geq 2}^{(-\frac 1 2)}$ a real-valued symbol in $\Sigma \Gamma_{2  }^{-\frac 1 2}[r, N]$, and $ \bR_{\geq 2, 2}(\zak)$ is a real-to-real matrix of smoothing operators in $\Sigma\cR_{2}^{-\vr}[r, N]$.
\\
Combining \eqref{mammolo}, \eqref{dottoceragia}, \eqref{pisolo}, \eqref{eolo}, \eqref{brontolo}, \eqref{brontolo.parecchio} and  that by the very definition \eqref{g112}
\begin{gather*}
        \ii \, \tg (-\ii {\bf \Omega}(D) \zak; x) - \ii \tilde{\tq}_{1} (\zak; x) = 0\,, \qquad 
        \ii \, \tf(-\ii \vOmega(D) \zak;x) +  \frac{\ii}{2}\tg_x (\zak; x)=0 \,,
\end{gather*}
one has that $U_3$ satisfies \eqref{eq.riduzione.ord12} with
\begin{equation}\label{esempioproiezioni}
\begin{aligned}
 &   \tq_{\geq 2}:= -\cP_{\geq 2}[\tq_{\geq 1}] + q_{\geq 3} +p_{\geq 2}
 \in\Sigma\cF_{\geq 2}^\R[r, N]\,, \quad
    \tb_{\geq 2} := b_{\geq 3} + b_{\geq 2}+ r_{\geq 2} \in\Sigma\cF_{ 2}^\R[r, N] \,, \\
  &  \tilde{\tc}_{\geq 1}^{(-\frac 1 2)} := \cP_{\leq N-1} \left[a_{\geq 1}^{(-\frac 1 2)} + b_{\geq 3}^{(-\frac 1 2)} + c_{\geq 1}^{(-\frac 1 2)} + d_{\geq 1}^{(-\frac 12)} + e_{\geq 2}^{(-\frac 1 2)}\right] \in\Sigma_1^{N-1}\wt\Gamma^{-\frac12}
  \,, \\
 &  \bB_{\geq N}(\zak) := \check{\bB}_{\geq N}(\zak) + \vOpbw{\cP_{\geq N} \left[ 
 a_{\geq 1}^{(-\frac 1 2)} + b_{\geq 3}^{(-\frac 1 2)} + c_{\geq 1}^{(-\frac 1 2)} + d_{\geq 1}^{(-\frac 12)} + e_{\geq 2}^{(-\frac 1 2)}
 \right]}\,,
\end{aligned}
\end{equation}
where $\bB_{\geq N}(\zak)$ is a real-to-real matrix of spectrally localized maps  in $\cS_{\geq N}^0[r]$ and  $\bR_{\geq 1}(\zak) $, obtained by collecting all the smoothing remainders, is a real-to-real matrix of smoothing remainders in $\Sigma \cR^{-\vr +{2N}}_{1}[r, N]$.
\\
Finally, by \eqref{omegonejin},  we  recombine the symbol of the dispersion relation writing 
$
  \vOpbw{ \im \frac{\gamma}{2}\sign(\xi)+  \im \omega(\xi) }  = \Opbw{\im \Omega(\xi)} =    \ii {\bf\Omega}(D)$.
\end{proof}
\paragraph{Removal of symbols of strictly-negative order and  homogeneity $1$.}
In this paragraph we remove the  degree $1$ homogeneous  component from the diagonal symbols in \eqref{eq.riduzione.ord12} of negative order. More precisely, we prove the following lemma.
\begin{lemma}[{\bf Normal form of degree 1 homogeneous components of negative order }]\label{riduzione.negativi}
Let $N\in\N$, $N\geq2$. Let $\vr>2N$ and $\wt\vr>3$. 
There exist $s_0, r  >0 $ and a $(0,0)$- admissible transformation $ {\bf \Psi}_3(\zak)$ with gain $\wt \vr$ (recall Definition \ref{admtra}) such that if $(\eta, \psi)(t) \in B_{s_0, \R}(I;r)$ solves \eqref{eq:etapsi}, and $\zak=\zak(t)$ is defined in \eqref{zak}, then the variable $U_4:= {\bf \Psi}_3(\zak)U_3$, with $U_3$ solving \eqref{eq.riduzione.ord12}, solves

\begin{equation}\label{eq.riduzione.ord.negativi}
\begin{aligned}
    \partial_t U_4  =
-\im  \bOmega(D) U_4 
    + \vOpbw{\im \sfD_{\geq 2}(\zak; x,\x) }U_4 +\bB_{\geq N}(\zak) U_4 + \bR_{\geq 1}(\zak) U_4\,,
\end{aligned}
\end{equation}
where $\bOmega(D)$ is defined in \eqref{diaglin}-\eqref{omegaxi} and 
\begin{itemize} 
    \item  $\sfD_{\geq 2}(\zak; x, \xi)$ is the real-valued symbol in $\Sigma \Gamma_2^1[r,N]$ defined in  \eqref{def:D};
\item  $\bB_{\geq N}(\zak)$ is a real-to-real matrix of spectrally localized maps in $\cS^0_{\geq N}[r]$ and $\bR_{\geq 1}(\zak) $ is a real-to-real matrix of smoothing remainders in $\Sigma \cR^{-\vr + {2N}}_{1}[r, N]$.
\end{itemize}
\end{lemma}
\begin{proof}
Define the map $ {\bf \Psi}_3(\zak)$ as the time-1 flow ${ \bf \Psi}_1(\zak):= { \bf \Psi}_1^{\tau}(\zak)|_{\tau=1}$ of 
$$
\begin{cases}
    \pa_\tau { \bf \Psi}_3^{\tau}(\zak)= \vOpbw{\ii \tilde g_1^{(-\frac12)}(\zak;x,\xi)}
  \, 
 { \bf \Psi}_3^{\tau}(\zak) \\
 { \bf \Psi}_3^0(\zak)={\rm Id}
\end{cases}, \qquad 
\quad\text{where } \quad  \tilde g_1^{(-\frac12)}:= \tilde g_1^{(-\frac12)}(\zak;x,\xi)\in \wt \Gamma_1^{-\frac12}
$$
is a symbol that we shall determine explicitly in \Cref{lem:g1-12} below.
By Lemma \ref{lem:flow.ad}--(iii),  the map ${\bf \Psi}_2(\zak)$ is a $(0,  0)$- admissible transformation of arbitrary gain
$\wt \vr \geq 0$.

We compute the equation solved by $U_4:= {\bf \Psi}_3(\zak)U_3$. As 
 $U_3$ solves equation  \eqref{eq.riduzione.ord12}, we get
\begin{gather}
    \partial_t U_4 = {\bf \Psi}_3(\zak)\Big(\vOpbw{
    - \im \Omega(\xi) + \im f^{(1)}_{\geq 1}} + \bB_{\geq N}(\zak) + \bR_{\geq 1}(\zak) \Big)  {\bf \Psi}_3^{-1}(\zak) U_4 + (\partial_t {\bf \Psi}_3(\zak)) {\bf \Psi}_3(\zak)^{-1} U_4\, , \\
    \text{where } \ \ 
    f^{(1)}_{\geq 1}:=  \big( \check{\sfV}_{2}  +  \check{\sfV}_{\geq 3}\big)\xi + {\tq}_{\geq 2}\, \omega(\xi) + \tb_{\geq 2}\,\sign{(\xi)}+ \tilde \tc_{\geq1}^{(-\frac12)} \in \Sigma \Gamma^{1}_1[r,N]  \,\notag
\end{gather}
is a real-valued symbol.
We  study how each term is transformed, starting with the   dispersion relation. By \Cref{prop:FIO}--2, we have
\begin{align}
 {\bf \Psi}_3(\zak)\vOpbw{-\ii\Omega(\xi)}   {\bf \Psi}_3(\zak)^{-1}= \vOpbw{-\ii\Omega(\xi)+ \tilde g_1^{(-\frac12)} \#_\vr \Omega(\xi)- \Omega(\xi)\#_\vr \tilde g_1^{(-\frac12)} +\ii d_{\geq 2}^{(-\frac52)}}+ \bR_{\geq 1}(\zak), 
 \label{disp3}
\end{align}
where $ d_{\geq 2}^{(-\frac52)}$ is a real-valued symbol in $ \Sigma \Gamma_{2}^{-\frac52}[r, N]$ and $\bR_{\geq 1}(\zak)$ is a real-to-real matrix of smoothing operators in $ \Sigma \cR^{-\vr}_1[r, N]$.
By \Cref{prop:FIO}--1,  we have 
\begin{align}
    {\bf \Psi}_3(\zak)\vOpbw{\im f^{(1)}_{\geq 1}}   {\bf \Psi}_3(\zak)^{-1}= \Opbw{\im f^{(1)}_{\geq 1}+ \ii d_{\geq 2}^{(-\frac12)}}+ \bR_{\geq 2}(\zak), 
    \label{fgeq13}
\end{align}
where $d_{\geq 2}^{(-\frac12)}$ is a real-valued symbol in $\Sigma \Gamma_2^{-\frac12}[r, N]$  and $\bR_{\geq 2}(\zak)$ is a real-to-real matrix of  smoothing operators in $ \Sigma \cR^{-\vr+1}_2[r, N]$.
Then, by \Cref{prop:FIO}--3, we have 
\begin{align}
     {\bf \Psi}_3(\zak)\bB_{\geq N}(\zak)   {\bf \Psi}_3(\zak)^{-1}=\check{\bB}_{\geq N}(\zak)  \ , \quad 
      {\bf \Psi}_3(\zak)\bR_{\geq 1}(\zak)  {\bf \Psi}_3(\zak)^{-1}=\bR'_{\geq 1}(\zak) , 
     \label{Bione3}
\end{align}
where $\check{\bB}_{\geq N}(\zak)$ is a real-to-real matrix of  spectrally localized maps in $\cS_{\geq N}^0[r]$
and
$\bR'_{\geq 1}(\zak)$ a real-to-real matrix of smoothing operators in 
$\Sigma \cR^{-\vr+{2N}}_1[r, N]$.
Finally, consider $(\partial_t {\bf \Psi}_1(\zak)) {\bf \Psi}_1(\zak)^{-1}$. By  \Cref{prop:FIO}--4, we get
\begin{align}
    (\partial_t {\bf \Psi}_3(\zak)) {\bf \Psi}_3(\zak)^{-1}=\vOpbw{\ii  \tilde g_1^{(-\frac12)}(-\ii \vOmega(D)\zak;x,\xi) +\ii d_{\geq 2}^{(-\frac12)}(\zak;x,\xi)}+\bR_{\geq 2}(\zak)
    \label{det3}
\end{align}
where $d_{\geq 2}^{(-\frac12)}\in\Sigma\Gamma^{-\frac12}_{2}[r, N]$ is real-valued and $\bR_{\geq 2}(\zak)$ is a real-to-real matrix of  smoothing operators in $\Sigma \cR_2^{-\vr}[r, N]$. 

Gathering \eqref{disp3}, \eqref{fgeq13}, \eqref{Bione3} and \eqref{det3}, denoting by $\wt\tc_{ 1}^{(-\frac12)} := \cP_{1}[\wt\tc_{\geq 1}^{(-\frac12)}]$ the real-valued, $1-$homogeneous 
 component in  $ \wtGamma_{1}^{-\frac12}$ of the symbol $\wt\tc_{ \geq 1}^{(-\frac12)} \in \Sigma_1^{N-1} \wt\Gamma^{-\frac12}$, 
and  projecting 
appropriately the symbols with respect to the homogeneity as done in \eqref{esempioproiezioni}, we get 
\begin{align}\notag
        \partial_t U_4  =& - \im \bOmega(D) U_4 +
       \textup{Op}_{\tt vec}^{{\scriptscriptstyle{\mathrm BW}}}\Big( \underbrace{ \im (\check{\sfV}_{2}+  \check{\sfV}_{\geq 3}) \xi + \im \, \tq_{\geq 2}\omega(\xi) + \ii \tb_{\geq 2}\,\sign{(\xi)}+\ii \tilde \tc_{\geq 2}^{(-\frac12)}}_{=: \im \, \sfD_{\geq 2}(\zak; \cdot) \mbox{ in } \eqref{def:D}} \Big) U_4\\
        &+ \vOpbw{\ii  \tilde g_1^{(-\frac12)}(-\ii \vOmega(D)\zak;x,\xi) +\tilde g_1^{(-\frac12)} \#_\vr \Omega(\xi)- \Omega(\xi)\#_\vr \tilde g_1^{(-\frac12)} +\ii \tilde \tc_1^{(-\frac12)}}U_4\label{homologichetta}\\
        \notag
    &+{\bB}_{\geq N}(\zak) U_4 + \bR_{\geq 1}(\zak) U_4\,
\end{align}
where $\bB_{\geq N}(\zak)$ is a real-to-real matrix of spectrally localized maps  in $\cS_{\geq N}^0[r]$ and  $\bR_{\geq 1}(\zak) $  is a real-to-real matrix of smoothing remainders in $\Sigma \cR^{-\vr +{2N}}_{1}[r, N]$.
Then \Cref{riduzione.negativi} follows by the following result.
\begin{lemma}\label{lem:g1-12}
Let $\vr >0$.
    There is a real-valued symbol $\tilde g_1^{(-\frac12)}\in \wt \Gamma_1^{-\frac12}$ such that 
    \begin{align}\label{sym.iter}
        \ii  \tilde g_1^{(-\frac12)}(-\ii \vOmega(D)\zak; \cdot) +\left(\tilde g_1^{(-\frac12)} \#_\vr \Omega(\xi)- \Omega(\xi)\#_\vr \tilde g_1^{(-\frac12)}\right)(\zak; \cdot ) +\ii \tilde \tc_1(\zak; \cdot) \in \wt \Gamma^{-\vr}_1.
    \end{align}
\end{lemma}
\begin{proof}
    Define $p\in \N$ such that $p+1> 2 \vr$ and 
$\tilde g_1^{(-\frac12)}= \tilde h_1^{(-\frac12)}+\ldots+\tilde h_1^{(-\frac{p}{2})}$ 
    where 
\begin{align*}
&\tilde h_1^{(-\frac12)}(\zak;x,\xi):= -\tilde \tc_1^{(-\frac12)}\Big((-\ii \vOmega(D))^{-1}\zak;x,\xi\Big) \in \widetilde \Gamma_1^{-\frac12},\\
&\tilde h_1^{(-1)}(\zak;x,\xi):= \ii \left(\tilde h_1^{(-\frac12)} \#_\vr \Omega(\xi)- \Omega(\xi)\#_\vr \tilde h_1^{(-\frac12)}\right)\Big((-\ii \vOmega(D))^{-1}\zak;x,\xi\Big)\in \widetilde \Gamma_1^{-1}\\
&{\scriptscriptstyle \vdots}\\
&\tilde h_1^{(-\frac{p}{2})}(\zak;x,\xi):= \ii \left(\tilde h_1^{(-\frac{p-1}{2})} \#_\vr \Omega(\xi)- \Omega(\xi)\#_\vr h_1^{(-\frac{p-1}{2})}\right)\Big((-\ii \vOmega(D))^{-1}\zak;x,\xi\Big)\in \widetilde \Gamma_1^{-\frac{p}{2}}.
\end{align*}
With this choice, by linearity and using that $p+1>2\vr $,  we obtain 
that the symbol in \eqref{sym.iter} equals
\begin{align*}
    \tilde h_1^{(-\frac{p}{2})} \#_\vr \Omega(\xi)- \Omega(\xi)\#_\vr h_1^{(-\frac{p}{2})} \in \wt \Gamma_1^{-\vr }.
\end{align*}
Finally, since $\tilde{\tc}_1^{(-\frac{1}{2})}$ is real-valued, it follows that $\tilde{h}_1^{(-\frac{1}{2})}$ is also real. Then, using \eqref{prop:ov}, we conclude that each $\tilde{h}_1^{(-\frac{j}{2})}$ is real-valued for $ j = 1, \ldots, p$, which in turn implies that $\tilde{g}_1^{(-\frac{1}{2})}$ is real.
\end{proof}
We conclude the proof including the operator  \eqref{homologichetta} in the smoothing remainder, thus getting \eqref{eq.riduzione.ord.negativi}.
\end{proof}

\subsubsection{Poincaré-Birkhoff normal form of the quadratic smoothing remainder}
In this  sub-section we remove the  $1$-homogeneous  component from the real-to-real matrix of smoothing operators $\bR_{\geq 1}(\zak)$ in 
\eqref{eq.riduzione.ord.negativi}. More precisely, we prove the following result.
\begin{proposition}[{\bf Poincaré-Birkhoff normal form of the water waves at quadratic degree}]\label{poinc}
Let $N \in \N$.
 Let $\vr> {2N+\frac32}$ and $\bar{\vr} = \vr - {2N}$. 
There exist $s_0, r  >0 $ and a $(0,0)$- admissible transformation ${\bf \Psi}_{4}(\zak)$
 with gain $\bar{\vr}$ (recall Definition \ref{admtra}) 
    such that if $(\eta, \psi)\in B_{s_0,\R}(I;r)$ solves \eqref{eq:etapsi}, and $\zak=\zak(t)$ is defined in \eqref{zak}, then the variable $Y={\bf \Psi}_4(\zak)U_4$, with $U_4$ solving \eqref{eq.riduzione.ord.negativi}, solves equation \eqref{BNF1}.
\end{proposition}

\begin{proof}
Define the map  ${\bf \Psi}_4(\zak)$ as the time-$1$  flow ${\bf \Psi}_4(\zak):={\bf \Psi}_4^\tau(\zak)_{|_{\tau=1}}$ of
 \be\label{BNFstep1}
\partial_{\tau} {\bf \Psi}_4^{\tau}(\zak)  = \bQ_1(\zak) {\bf \Psi}_4^{\tau}(\zak) \, , 
\quad {\bf \Psi}_4^{0}(\zak) = {\rm Id} \, ,  
\ee
with $ \bQ_1(\zak)$  a $1-$homogeneous real-to-real matrix of smoothing operators in $ \widetilde{\mathcal{R}}^{-\vr +{2N}}_1$   to be determined. 
    By Lemma \ref{flow.s.ad},  the map ${\bf \Psi}_4(\zak)$ is a $(0, 0)-$admissible transformation with gain $\bar\vr = \vr - {2N}$.

We compute the equation solved by $Y={\bf \Psi}_4(\zak)U_4$. As 
 $U_4$ solves equation  \eqref{eq.riduzione.ord.negativi}, we get
\begin{equation}\label{coniugazione.BNF1}
\begin{aligned}
\pa_t Y 
=  &  {\bf \Psi}_4(\zak)( -\ii {\bf\Omega}(D) ){\bf \Psi}_4(\zak)^{-1} Y
+ {\bf \Psi}_4(\zak)\vOpbw{ \ii \sfD_{\geq 2}}{\bf \Psi}_4(\zak)^{-1}Y \\
&+ (\pa_t {\bf \Psi}_4(\zak)) {\bf \Psi}_4(\zak)^{-1}Y + {\bf \Psi}_4(\zak) \bB_{\geq N} (\zak) {\bf \Psi}_4(\zak)^{-1}Y  + 
      {\bf \Psi}_4(\zak)  \bR_{\geq 1}(\zak)  {\bf \Psi}_4(\zak)^{-1} Y 
    \end{aligned}
\end{equation}
By \Cref{prop:Egorov_smoothing}--1 
  one has
   \begin{align}\label{core4}
   {\bf \Psi}_4(\zak)( -\ii {\bf\Omega}(D) ){\bf \Psi}_4(\zak)^{-1}  & =
 -\ii\bOmega(D)+ \big[ \bQ_1(\zak),-\ii \bOmega(D) \big]  + \bR_{\geq 2,1}(\zak) \,,\\
 \label{core2}
   {\bf \Psi}_4(\zak)\vOpbw{ \ii\sfD_{\geq 2})}{\bf \Psi}_4(\zak)^{-1} & =\vOpbw{ \ii\sfD_{\geq 2}}  + \bR_{\geq 2,2}(\zak) \,,
\end{align}
where $\bR_{\geq 2,1}(\zak)$ and $\bR_{\geq 2,2}(\zak)$ are  real-to-real matrices of smoothing operators  in $ \Sigma\mathcal{R}^{-\vr+ {2N+\frac12}}_{2}[r,N]$ respectively 
  $ \Sigma\mathcal{R}^{-\vr+{2N+1}}_{2}[r,N]$.
 Now denote by 
$$\bR_1(\zak):= \cP_{\leq 1}\Big[\bR_{\geq 1}(\zak) \Big]$$
the $1-$homogeneous component in $\wt\cR_1^{-\vr + 2N}$ of the real-to-real matrix of smoothing operators $\bR_{\geq 1}(\zak)  \in \Sigma \cR_1^{-\vr + 2N}[r,N]$ of \eqref{eq.riduzione.ord.negativi}; 
 by  \Cref{prop:Egorov_smoothing}--2, we obtain
 \begin{align}\label{core3}
&  {\bf \Psi}_4(\zak)\,  \bR_{\geq 1}(\zak)\, {\bf \Psi}_4(\zak)^{-1}=
  \bR_{ 1}(\zak) +  \bR_{\geq 2}(\zak)\,, \\
  \label{core} 
&  {\bf \Psi}_4(\zak)\,  \bB_{\geq N}(\zak) \,  {\bf \Psi}_4(\zak)^{-1}
= \bB_{\geq N}(\zak) + \bR_{\geq N}(\zak)\,,
 \end{align}
  where  $\bR_{\geq 2}(\zak)$  and 
  $\bR_{\geq N}(\zak)$ are real-to-real matrices of smoothing operators in
  $\Sigma\mathcal{R}^{-\vr+{2N}}_{ 2}[r, N]$ respectively
  $\mathcal{R}^{-\vr+2N}_{\geq N}[r]$.
Finally by \Cref{prop:Egorov_smoothing}--3,
 one gets
 \begin{align}
  (\pa_t {\bf \Psi}_4(\zak)) \, {\bf \Psi}_4(\zak)^{-1} &=  
    \bQ_1(-\ii \bOmega(D) \zak) + \bR_{\geq 2}(\zak) \label{core5}\,,
\end{align}
where $\bR_{\geq 2}(\zak)$ is a real-to-real matrix of smoothing operators in $\Sigma\mathcal{R}^{-\vr+{2N+\frac32}}_{2}[r, N]$.
Combining 
\eqref{coniugazione.BNF1}--\eqref{core5}, 
the variable $Y$ solves 
\be\label{Y.eq.R1}
\begin{aligned}
    \partial_t Y &  =
-\im  \bOmega(D) Y
    + \vOpbw{\im \sfD_{\geq 2} }Y  +\bB_{\geq N}(\zak) Y\\
    & \quad + \Big( \bR_{1}(\zak)  + \bQ_1(-\ii \bOmega(D) \zak)  + \big[ \bQ_1(\zak),-\ii \bOmega(D) \big] \Big) Y + 
    \bR_{\geq 2}(\zak)Y\,,
\end{aligned}
\ee
where $ \bR_{\geq 2}(\zak)$ is a real-to-real matrix of smoothing operators in $\Sigma \cR_2^{-\vr + 2N + \frac32}[r,N]$, obtained by collecting all the smoothing remainders. 
We choose   $\bQ_1(\zak)$ to solve the homological equation 
\begin{equation}\label{eq.omologica.linearmenorho}
\bR_1(\zak) +  \bQ_1(-\ii \bOmega(D) \zak)+ \big[  \bQ_1(  \zak), -\ii \bOmega(D) \big]=0 \, . 
\end{equation}
Expanding both vector fields in Fourier components as 
 \begin{equation}\label{formaQ1}
(\bR_1(\zak) Y)^{\sigma_1}_k = \sum_{ \fP_3
} R_{j, j', k}^{\sigma, \sigma', \sigma_1}
\zetina_{j}^{\sigma} y_{j'}^{\sigma'}\,,
\qquad (\bQ_1(\zak) Y)^{\sigma_1}_k = \sum_{ \fP_3
} Q_{j, j', k}^{\sigma, \sigma', \sigma_1}
\zetina_{j}^{\sigma} y_{j'}^{\sigma'}  
 \end{equation}
 where with the sum over $\fP_3$ we mean that the indexes $(j, j', k,\sigma, \sigma', -\sigma_1)$ belong to the set $\fP_3$ in \eqref{mom1}, 
the homological equation \eqref{eq.omologica.linearmenorho} is tantamount to 
\[
 \im \Big(\sigma_1 \Omega_k(\gamma) - \sigma' \Omega_{j'}(\gamma) - \sigma \Omega_j(\gamma) \Big)Q_{j, j', k}^{\sigma, \sigma', \sigma_1} + R_{j, j', k}^{\sigma, \sigma', \sigma_1}  = 0\,.
\]
By  \Cref{3onde} this equation has a unique solution $\bQ_1(\zak)$ of the form \eqref{formaQ1} with coefficients defined as 
\be\label{bQ1}
 Q_{j, j', k}^{\sigma, \sigma', \sigma_1}
 := 
  \dfrac{ \, R_{j, j', k}^{\sigma, \sigma', \sigma_1}}{ \im \big(  \sigma \Omega_j(\gamma) + \sigma'\Omega_{j'}(\gamma)  -\sigma_1 \Omega_k(\gamma)\big)}  \ , \qquad  (j, j', k,\sigma, \sigma', -\sigma_1)\in \fP_3\,.
 \ee
In conclusion, by  \eqref{Y.eq.R1} 
and the assumption that  $ \bQ_1(\zak) $ solves \eqref{eq.omologica.linearmenorho}  we deduce \eqref{BNF1}. 
To conclude the proof we only need to prove that $\bQ_1(\zak)$ is a real-to-real matrix of smoothing operators in $\wt\cR^{-\vr+{2N}}_1.$ This is done in the following lemma.
\end{proof}

 \begin{lemma}
  $\bQ_1(\zak)$ in \eqref{formaQ1}--\eqref{bQ1} is a real-to-real matrix of smoothing operators in  
 $\wt \cR^{-\vr + {2N}}_1 $.
 \end{lemma}
 \begin{proof}
 As $\bR_1(\zak)$ is a real-to-real matrix of  smoothing operators in $\wt \cR^{-\varrho+{2N}}_1$, its coefficients fulfill the estimate (see \eqref{smoocara}): for some $\mu\geq 0$, $C>0$, 
\begin{equation}\label{boundR1}
\abs{ R_{j, j',  k}^{ \sigma, \sigma',  \sigma_1} } \leq C 
\frac{
{\rm max}_2\{ \la j\ra,  \la j'\ra \}^{\mu}}{ \max\{ \la j\ra, \la j' \ra \}^{\varrho- {2N}}}
   \, , \quad \forall \, (j, j', k, \sigma, \sigma',  -\sigma_1) \in  \fP_3 \ ,
\end{equation}
 and satisfy the symmetric and reality properties  \eqref{M.coeff.p} and \eqref{M.realtoreal}.
 
 Consider now the coefficients $  Q_{j, j', k}^{\sigma, \sigma', \sigma_1}$ in \eqref{bQ1}. 
 Clearly, they satisfy the symmetric and reality properties \eqref{M.coeff.p} and \eqref{M.realtoreal}, proving that $\bQ_1(\zak)$ is a real-to-real matrix of operators. 
 We now bound the coefficients $  Q_{j, j', k}^{\sigma, \sigma', \sigma_1}$. By 
\eqref{boundR1},   \Cref{3onde}  and the momentum relation $\sigma_1 k= \sigma j + \sigma' j$, 
 $$
 \abs{   Q_{j, j', k}^{\sigma, \sigma', \sigma_1}} \leq \frac{C}{\tc_3(\gamma)} 
\frac{
{\rm max}_2\{ \la j\ra,  \la j'\ra \}^{\mu}}{ \max\{ \la j\ra, \la j' \ra \}^{\varrho-{2N}}}
   \, , \quad \forall (j, j', k, \sigma, \sigma',  -\sigma_1) \in  \fP_3 
\,.
 $$
This shows that  $\bQ_1(\zak)$ is a matrix of smoothing operators in  $\wt\cR_1^{-\varrho +{2N}}$.
 \end{proof}
\begin{proof}[Proof of Theorem \ref{thm:quadratic.nf}]\label{proof.thm:quadratic.nf}
Let $\vr > {2N+\frac92}$ and $\tilde \vr \geq  \vr+{\frac92-2N} $.  We define the map
\be
\bT(\zak):=  {\bf \Psi}_4(\zak) \circ  {\bf \Psi}_3(\zak) \circ  {\bf \Psi}_2(\zak)\circ 
 {\bf \Psi}_1(\zak)\circ \mathbf{\Psi}_0(\zak) \circ {\bf \Psi}(\zak)^{-1}
\ee
where ${\bf \Psi}(\zak)$ is defined in
 \Cref{diag}, $\mathbf{\Psi}_0(\zak)$ in \Cref{diag.ord0}, 
 $ {\bf \Psi}_1(\zak)$ in 
 \Cref{ridlinord1}, 
 ${\bf \Psi}_2(\zak)$ in 
 \Cref{lem:rid12}, 
$ {\bf \Psi}_3(\zak)$ in 
 \Cref{riduzione.negativi}  
 and finally
 ${\bf \Psi}_4(\zak)$ in 
 \Cref{poinc}.
Each 
transformation in the right hand side is admissible.
Hence, by Lemma \ref{lem:comp},  the map $\bT(\zak)$ is $(0, \frac92)$ admissible with gain 
$$ \min( \vr-{2N}, \tilde \vr-\frac92 ) = \vr-{2N} \geq \frac92 \  . $$
The explicit form of the transport term $\check{\sfV}_2(\zak; x)$ in \eqref{V2.exp0} follows by \Cref{lem:homo.eq.trans}.
\end{proof}
\subsection{The cubic normal form of Water Waves}\label{subsec:cubicNF}
So far we have conjugated the original water waves vector field in \eqref{eq:etapsi} to the quadratic paradifferential normal form  in \eqref{BNF1}. Such conjugation is valid {\em for every value of the vorticity} $\gamma \in \R$. 
Instead, in this section we are going to fix a resonant vorticity $\gamma <0$,  $\gamma^2 \in \Q$,  choose a 
set $\Lambda = \{\tm, \tn\}$, $\tm < 0 < \tn$ which is  $\gamma$-good  (cf. \Cref{g-good}), put the quadratic  transport term in its resonant normal form and 
the cubic smoothing remainder  in {\em strong $\Set$-normal form} (cf.  \Cref{def:wr}).
Remember that, given a resonant vorticity $\gamma$, there are countable many $\gamma$-good sets, as proved in \Cref{lem:goodset}.
Precisely, we  shall perform the following transformations: 
\begin{enumerate}
    \item First, in \Cref{thm:cubic.nf.1},  we put the quadratic  transport  velocity field  $\check{\sfV}_2(\zak;x) $ in \eqref{V2.exp0} in {\em resonant} normal form, removing  all its  non-resonant monomials and  reducing it to the quadratic transport  velocity field $\langle \sfV \rangle (\zak;x) $ in \eqref{VresZ}, that we compute explicitly in  \Cref{lem:coeffV};
    \item Second, in \Cref{thm:wnf}, we put the cubic part of the real-to-real  smoothing vector field $\bR_{\geq 2}(\zak)[Y]$ in \eqref{BNF1}  in {\em weak-$\Lambda$ normal form}, according to Definition \ref{def:wr}; in particular, we shall remove from the vector field all the non-resonant monomials with at most 2 indexes outside $\Lambda$. This is possible thanks to   \Cref{lem:wres} that characterize four-waves resonant interactions.
\end{enumerate}
\subsubsection{Cubic resonant normal form for the transport operator}

To reduce the  transport term 
$\vOpbw{\im \check{\sfV}_2 (\zak; x)\xi}$ in \eqref{def:D}  to its resonant normal form, we shall exploit the following result.
\begin{lemma}\label{lem:coeffV}

 Let $\gamma < 0 $, $\gamma^2 \in \Q$.
Let $N\in \N$, $N \geq 3$. There exist  $s_0, r>0$ such that the following holds true.
Let $(\eta,\psi) \in B_{s_0, \R}(I;r)$ be a solution of \eqref{eq:etapsi}, 
and define $\zak=\zak(t)$ by \eqref{zak}.
Let $\check{\sfV}_2$ be the real-valued function in $\wt\cF^\R_2$ 
defined in \eqref{V2.exp0}--\eqref{V2.exp1}. 
 There exists a real-valued homogeneous function $\beta_2 \in \wt\cF^\R_2$ such that 
\begin{equation}\label{be.small2}
   (\beta_2)_t + \check{\sfV}_{2} 
     =\la \sfV \ra + \beta_{\geq 3} \,, 
    \end{equation}
    where the resonant part $\la \sfV \ra$
    is the real-valued function in $\wt \cF^\R_2$  
    \begin{equation}
     \langle \, \sfV\, \rangle(\zak;x) := \hspace{-2em}
		\sum_{\sigma_1 \Omega_{k_1}(\gamma) + \sigma_2 \Omega_{k_2}(\gamma) = 0}
      \hspace{-2em}
        V_{k_1, k_2}^{\sigma_2, \sigma_2}(\gamma)\,   {\zetina_{k_1}^{\sigma_1}} \, {\zetina_{k_2}^{\sigma_2}} \, e^{\im (\sigma_1 k_1 + \sigma_2 k_2)x}  \, ,
        \qquad 
        V_{k_1,k_2}^{\sigma_1, \sigma_2}(\gamma) \mbox{  in } \eqref{V2.exp1}
        \label{Vres.int}
 \end{equation}
 and
 $ \beta_{\geq 3}$ is a real-valued function in $ \Sigma \cF^{\R}_{3}[r, N] $.
For any $\gamma$-good set $\Lambda = \{\tm, \tn\}$, $\tm <0<\tn$, the function
$\langle \, \sfV\, \rangle$ has the form  \eqref{VresZ}--\eqref{Vres.coeff2}.
\end{lemma}
\begin{proof}
The real-valued function $\beta_2 \in  \widetilde\cF^\R_{2}$ is chosen to remove the non-resonant quadratic monomials of $\check{\sfV}_2(\zak;x) =  \sum V_{k_1, k_2}^{\sigma_1, \sigma_2}(\gamma) \,  \zetina_{k_1}^{\sigma_1}\zetina_{k_2}^{\sigma_2} e^{\im (\sigma_1 k_1 + \sigma_2 k_2) x}$, namely the ones
 supported on indexes   $\sigma_1 \Omega_{k_1}(\gamma) + \sigma_2 \Omega_{k_2}(\gamma) \neq  0$.
We now prove that such $\beta_2$ 
is defined by 
\begin{align}
\beta_2(\zak;x) :=
\hspace{-2em}
    \sum_{\sigma_1\Omega_{k_1}(\gamma)+\sigma_1\Omega_{k_2}(\gamma)\neq 0} 
    \hspace{-2em}
    \beta_{k_1, k_2}^{\sigma_1, \sigma_2}(\gamma) \  \zetina_{k_1}^{\sigma_1}\zetina_{k_2}^{\sigma_2} e^{\im (\sigma_1 k_1 + \sigma_2 k_2) x}  \,,
    \quad \beta_{k_1, k_2}^{\sigma_1, \sigma_2}(\gamma):= \frac{V_{k_1, k_2}^{\sigma_1, \sigma_2}(\gamma)}{\ii \big(\sigma_1 \Omega_{k_1}(\gamma)+ \sigma_2 \Omega_{k_2}(\gamma)\big)} \,. 
\label{def:beta2}
\end{align}
First note that, by item  $(ii)$  of \Cref{lem:2wave} and since $\check{\sfV}_2 \in \wt\cF^\R_2$, there is $\mu \geq 0$ such that for any 
$(k_1, k_2, \sigma_1, \sigma_2) \in \Z_*^2 \times \{\pm \}^2$ fulfilling  $  \sigma_1\Omega_{k_1}(\gamma)+\sigma_1\Omega_{k_2}(\gamma)\neq 0$, the coefficients 
$\beta_{k_1, k_2}^{\sigma_1, \sigma_2}(\gamma)$ satisfy
$$
\abs{\beta_{k_1, k_2}^{\sigma_1, \sigma_2}(\gamma)} \lesssim (1+|k_1| + |k_2|)^{\mu} \,, \quad
\overline{\beta_{k_1, k_2}^{\sigma_1, \sigma_2}(\gamma)} = \beta_{k_1, k_2}^{-\sigma_1, -\sigma_2}(\gamma) \,,   
$$
 thus proving that $\beta_2 \in \wt\cF^\R_2$.
Then, since 
 $\zak$ solves \eqref{claim:zak.M}, 
 \begin{align*}
 (\beta_2)_t 
 &=
 \!\!\!\!\!\!\!\!\!\!\!\!\!\!
 \sum_{{\sigma_1\Omega_{k_1}(\gamma)+\sigma_1\Omega_{k_2}(\gamma)\neq 0}} 
 \!\!\!\!\!\!\!\!\!\!\!\!\!\!
 -\ii (\sigma_1 \Omega_{k_1}(\gamma)+\sigma_2 \Omega_{k_2}(\gamma))\  \beta_{k_1, k_2}^{\sigma_1, \sigma_2}(\gamma) \, \zetina_{k_1}^{\sigma_1} \ \zetina_{k_2}^{\sigma_2}e^{\ii(\sigma_1  k_1 + \sigma_2  k_2)x}  + \beta_{\geq 3}(\zak) 
 \end{align*} 
with $\beta_{\geq 3}(\zak):=  \beta_2(\bM_{\geq 1}(\zak)\zak, \zak) + \beta_2(\zak, \bM_{\geq 1}(\zak)\zak)\in \Sigma \cF^\R_3[r,N]$ (here $\bM_{\geq 1}(\zak)$ is the map appearing in \eqref{claim:zak.M}).
Then  identity  \eqref{be.small2} follows immediately by 
 \eqref{def:beta2} and 
defining the resonant part $\la \sf V \ra$ as in \eqref{Vres.int}.

Let us now fix a $\gamma$-good set $\Lambda = \{\tm, \tn\}$, with $\tm <0<\tn$ and 
$\gamma =- \dfrac{\tm + \tn}{ \sqrt{2(\tn - \tm)}} $  according to \eqref{Omega*},
and show that $\la \sfV \ra(\zak;x)$ writes explicitly as in \eqref{VresZ}--\eqref{Vres.coeff2}.
First note that $\sigma_1 \Omega_{k_1}(\gamma) + \sigma_2 \Omega_{k_2}(\gamma) = 0$ forces 
$\sigma_1 = - \sigma_2$ as 
$\Omega_k(\gamma) >0$ for any $k \in \Z_*$.
Then, by item $(i)$ of \Cref{lem:2wave},  either $k_1 = k_2$ or 
$\sign (k_1) \neq \sign (k_2)$. 
 Hence, starting from \eqref{Vres.int}, we split 
 \begin{align*}
  \langle \, \sfV\, \rangle(\zak;x) = & 
     \sum_{j \in \Z_* } \underbrace{\Big(V_{j,j}^{+, -}(\gamma) +V_{j,j}^{-,+}(\gamma)  \Big)}_{=:\sfV^{(\rm int)}_j } |\zetina_{j}|^2 \\
     & +\sum_{\substack{m<0<n \\  \Omega_{m}(\gamma) =   \Omega_{n}(\gamma)} }\,  
        \underbrace{\Big( V_{n,m}^{+, -}(\gamma) + V_{m,n}^{-,+}(\gamma)}_{=:\sfV^{(\res)}_{m,n}} \Big)\,   \bar{\zetina_m} \, \zetina_n  \, e^{\im (n-m)x}  + 
       \underbrace{ \Big( V_{n,m}^{-, +}(\gamma) + V_{m,n}^{+,-}(\gamma) \Big)}_{= \bar{\sfV^{(\res)}_{m,n}}}\,   {\zetina_m} \, {\bar \zetina_n}  \, e^{-\im (n-m)x} 
        \, ,
        \end{align*}
where we used $\bar{V_{k_1, k_2}^{\sigma_1, \sigma_2}}(\gamma) = V_{k_1, k_2}^{-\sigma_1, -\sigma_2}(\gamma)$ 
 being $ \la \sf V \ra (\zak; x)$  real-valued.
Then the last term is the complex conjugate of the second one and one obtains formula \eqref{VresZ}. In order to compute the coefficients $\sfV^{(\res)}_{m,n}$, we note that, by \eqref{2res}, for any $m<0<n$, the conditions $\gamma<0$ and $\Omega_{m}(\gamma) =   \Omega_{n}(\gamma)$ imply
\begin{equation}
    - \dfrac{m + n}{ \sqrt{2(n - m)}}=\gamma=- \dfrac{\tm + \tn}{ \sqrt{2(\tn - \tm)}}\,.
    \label{gamma_rel}
\end{equation}

Formulas \eqref{Vres.coeff} and \eqref{Vres.coeff2}  follow from \eqref{V2.exp1}, \eqref{numerelli} substituting 
the second equality in \eqref{gamma_rel} for $\sfV^{(\rm int)}_j$ and the first equality in \eqref{gamma_rel} for $\sfV^{(\res)}_{m,n}$  and using $\tm <0<\tn$, $\tn + \tm >0$, namely computing
$$
\sfV^{(\rm int)}_j=\left. \Big(V_{j,j}^{+, -}(\gamma) +V_{j,j}^{-,+}(\gamma)  \Big)\right|_{\gamma =- \frac{\tm + \tn}{ \sqrt{2(\tn - \tm)}}}, \qquad \sfV^{(\res)}_{m,n}= \left.\Big( V_{n,m}^{+, -}(\gamma) + V_{m,n}^{-,+}(\gamma)\Big)\right|_{ \gamma=- \frac{m + n}{ \sqrt{2(n - m)}}}\,.
$$

The computation is  verified in the  Mathematica notebook \texttt{lemma5\_25.nb} at the link \url{https://git.sissa.it/amaspero/transfer-ww-vorticity}.
\end{proof}
We now perform a paracomposition in order to put   $\vOpbw{\im \check{\sfV}_{2}(\zak;x)\xi}$ in resonant normal form, proving  the following result.
\begin{proposition}[{\bf Normal form  of  degree 2 homogeneous components in the transport term}]\label{thm:cubic.nf.1}
 Let $\gamma < 0 $, $\gamma^2 \in \Q$ and $\Set = \{\tm, \tn\}$ be a $\gamma$-good set according to Definition \ref{g-good}. Let $N\in \N$, $N \geq 3$ and $\vr >{3N +\frac32}$, $\wt\vr >3 $.
There exist $s_0, r  >0 $ and a $(0,3)$- admissible transformation ${\bf \Psi}_5(\zak)$
 with gain $\wt\vr$ (recall Definition \ref{admtra}) 
    such that if $(\eta, \psi)\in B_{s_0,\R}(I;r)$ solves \eqref{eq:etapsi}, then the variable $Y_1={\bf \Psi}_5(\zak)Y$, with $Y$ solving \eqref{BNF1}, solves 
\begin{equation}\label{BNF2}
 \begin{aligned}
 \pa_{t}Y_1&= -\ii\bOmega(D)Y_1+ \vOpbw {\ii \sfD_{\geq 2}^{(\res)}(\zak; x,\x)}Y_1  +  \bB_{\geq N}(\zak) Y_1+  \bR_{\geq 2}(\zak)Y_1
 \end{aligned}
 \end{equation}
 where $\bOmega(D)$ is defined in \eqref{diaglin}-\eqref{omegaxi},  
 $\sfD_{\geq 2}^{(\res)}(\zak; x, \xi)$ is a real-valued symbol in $\Sigma\Gamma^{1}_{ 2}[r, N]$ of the form 
     \begin{equation}\label{def:D2}
     \sfD_{\geq 2}^{(\res)}(\zak; x, \xi):= \big(\langle \, \sfV\, \rangle (\zak;x) + \und{\sfV}_{\geq 3}(\zak; x)\big) \xi + \und{\tq}_{\geq 2}(\zak; x) \omega(\xi) +  \und{\tb}_{\geq 2}(\zak;x)\,\sign{(\xi)}+ \und{\widetilde\tc}_{\geq 2}^{(-\frac12)}(\zak; x, \xi)
     \end{equation}
     with $\omega(\xi)$ defined in \eqref{omegaxi} and\\
     $\bullet$   
 $\langle \, \sfV\, \rangle(\zak;x)$ is the real-valued function in $ \wt\cF^{\R}_{2}$  defined in \eqref{VresZ}--\eqref{Vres.coeff2}, whereas  $\und{\sfV}_{\geq 3}(\zak;x)$ is a real-valued function in $\Sigma \cF^\R_{3}[r, N]$;
\\
$\bullet$ $\und{\tq}_{\geq 2}(\zak;x)$ and  $\und{\tb}_{\geq 2}(\zak;x)$ are real-valued functions in  $ \Sigma  \cF^\R_2[r, N]$ and  $\und{\widetilde \tc}^{(-\frac12)}_{\geq 2}(\zak;x, \xi)$ is a pluri-homogeneous real-valued symbol in  $  \Sigma_2^{N-1} \wt \Gamma^{-\frac 1 2}$;
\\
$\bullet$ $\bB_{\geq N}(\zak)$ is a real-to-real matrix of spectrally localized maps in $\cS_{\geq N}^{0}[r]$ and $ \bR_{\geq 2}(\zak)$ is a real-to-real matrix of smoothing  operators  
    in  $\Sigma\mathcal{R}^{-\vr+{3N+\frac32}}_{2}[r, N] $.
\end{proposition}
\begin{proof}
We define the map $ {\bf \Psi}_5(\zak)$ as the time-1 flow 
$ {\bf \Psi}_5^\tau(\zak)\vert_{\tau = 1}$ of 
$$
\pa_\tau {\bf \Psi}_5^{\tau}(\zak)=
  \, \vOpbw{\ii \frac{\beta_2(\zak;x)}{1+\tau(\beta_2)_x(\zak;x)} \, \xi}
 {\bf \Psi}_5^{\tau}(\zak) , \qquad 
{ \bf \Psi}_5^0(\zak)={\rm Id} \, ,    
$$
where the real-valued function $\beta_2 \in \widetilde\mF^\R_{2}$ is defined in \Cref{lem:coeffV}.
By item $(i)$ of Lemma \ref{lem:flow.ad}, 
 the map ${\bf \Psi}_5(\zak)$ is a $(0, 3)$- admissible transformation of arbitrary gain $ \tilde \vr > 3$.
We compute the equation solved by $Y_1:= {\bf \Psi}_5(\zak)Y$. As 
 $Y$ solves equation  \eqref{BNF1}, we get
\begin{gather}
    \partial_t Y_1 = {\bf \Psi}_5(\zak)\Big(\bA(\zak)+ \bB_{\geq N}(\zak) + \bR_{\geq 2}(\zak) \Big)  {\bf \Psi}_5^{-1}(\zak) Y_1 + (\partial_t {\bf \Psi}_5(\zak)) {\bf \Psi}_5(\zak)^{-1} Y_1\,,\quad \text{where}\label{conjprediffeo}\\
    \notag
    \bA(\zak):=   \vOpbw{
    \im  \big(\check{\sfV}_{2} + \check{\sfV}_{\geq 3}\big) \xi 
    + \ii (-1+ \tq_{\geq 2})\omega(\xi) 
    +\ii (- \frac{\gamma}{2} +  \tb_{\geq 2}) \sgn(\xi)+ \im \tilde \tc_{\geq 2}^{(-\frac12)}}  \,. 
\end{gather}
We now study how the various terms in \eqref{conjprediffeo} are  transformed.
We start with $(\partial_t {\bf \Psi}_5(\zak)) {\bf \Psi}_5(\zak)^{-1}$. 
By Proposition \ref{prop:Egorov}-3 (with $p=2$) one has
$$
\big(\pa_t {\bf \Psi}_5(\zak)\big) {\bf \Psi}_5 (\zak)^{-1} =   \vOpbw{ \ii((\beta_2)_t  + g_{\geq 4})\, \xi } + \bR_{\geq 4}(\zak)
$$
where $ g_{\geq 4}$ is a real-valued function in $ \Sigma \cF^{\R}_{ 4}[r,N]$  and 
$ \bR_{\geq 4}(\zak) $ is a real-to-real  matrix of smoothing operators in  $  \mathcal{R}^{-\vr}_{\geq 4}[r,N] $.
Then consider the transport term of order 1. 
By applying Proposition  \ref{prop:Egorov}-1, (use formula \eqref{trans.transp})
\begin{equation}
{\bf \Psi}_5(\zak) \vOpbw{ \ii \big(\check{\sfV}_{2} + \check{\sfV}_{\geq 3}\big) \xi} {\bf \Psi}_5(\zak)^{-1} = 
  \vOpbw{ \ii \big(\check{\sfV}_{2} + h_{\geq 3}\big) \xi}
  + \bR_{\geq 4}(\zak) 
  \label{conj.ord1.1_2}
\end{equation}
where $ h_{\geq 3}$ is a real-valued function in $ \Sigma \cF^{\R}_{ 3}[r,N]$ 
and $ \bR_{\geq 4}(\zak) $ is a real-to-real matrix of smoothing operators in  $ \Sigma {\mathcal R}^{-\vr+1}_{4} [r,N]$.
Next  consider the paradifferential operator of order $\frac12$; arguing as in the conjugation of \eqref{conj.ord1.3}, 
\begin{equation}
{\bf \Psi}_5(\zak) \vOpbw{\ii \big( -1+ \tq_{\geq 2}\big) \omega(\xi) }{\bf \Psi}_5(\zak)^{-1} = 
\vOpbw{
\ii \big(-1 + \und{\tq}_{\geq 2} \big) \,  \omega(\xi) 
+\im {g}_{\geq 2 }^{(-\frac12)} } + \bR_{\geq 2}(\zak) 
\end{equation}
where $\und{\tq}_{\geq 2}$
is a  real-valued function in $\Sigma \cF^\R_{2}[r,N]$,   $  g^{(-\frac12)}_{\geq 2}$ is a real-valued symbol in $ \Sigma\Gamma^{-\frac12}_{ 2}[r,N]$ 
and $ \bR_{\geq 2}(\zak) $ is a real-to-real matrix of smoothing operators in $ \Sigma {\mathcal R}^{-\vr+\frac12}_{2}[r, N] $.

Next consider the paradifferential operator of order $0$; arguing as in the conjugation of \eqref{conj.ord1.2} we get 
\begin{equation}
{\bf \Psi}_5(\zak) \vOpbw{\ii \big(-\frac{\gamma}{2}+\tb_{\geq 2}\big) \sgn(\xi)} {\bf \Psi}_5(\zak)^{-1}
= \vOpbw{\ii \left( - \frac{\gamma}{2} +\und{\tb}_{\geq 2}\right)\sgn(\xi)+\ii g_{\geq 2}^{(-2)}} + \bR_{\geq 2}(\zak)
\end{equation}
where  $ \und{\tb}_{\geq 2}$ is a real-valued function in $\Sigma \cF^{\R}_2[r,N]$,   $ g_{\geq 2}^{(-2)}$ is  a real-valued symbol in $  \Sigma\Gamma^{-2}_{2}[r,N]$ 
and $ \bR_{\geq 2}(\zak)$ is a real-to-real matrix of  smoothing remainders in $ \Sigma {\mathcal R}^{-\vr}_{2} [r,N]$.
Next, by Proposition 
\ref{prop:Egorov}-1 we get 
\begin{equation}
     {\bf \Psi}_5(\zak) \vOpbw{\ii \tilde \tc^{(-\frac12)}_{ \geq 2}  }{\bf \Psi}_5(\zak)^{-1} =
    \vOpbw{ \ii {\td}^{(-\frac12)}_{ \geq 2} }+ \bR_{\geq 4}(\zak)    
    \label{conj.ord1.4_2}
\end{equation}
where ${\td}^{(-\frac12)}_{ \geq 2}$ is  a real-valued symbol in $\Sigma \Gamma_2 ^{-\frac12}[r, N]$, and $\bR_{\geq 4}(\zak)$ is a real-to-real matrix of smoothing operators in  $ \Sigma {\mathcal R}^{-\vr}_{4}[r,N] $.
Finally by Proposition \ref{prop:Egorov}-2 
\begin{align}
 &    {\bf \Psi}_5(\zak) \bB_{\geq N}(\zak) {\bf \Psi}_5(\zak)^{-1}= \check{\bB}_{\geq N}(\zak) \ , 
 \quad 
  {\bf \Psi}_5(\zak) \bR_{\geq 2}(\zak) {\bf \Psi}_5(\zak)^{-1}= \check{\bR}_{\geq 2}(\zak)
\end{align}
where 
 $\check \bB_{\geq N}(\zak)$ is a real-to-real matrix of spectrally localized maps in $\cS_{\geq N}^{0}[r]$ and $\check{\bR}_{\geq 2}(\zak)$
 is a real-to-real matrix of smoothing operators  in 
$ \Sigma\mathcal{R}^{-\vr+{3N+\frac32}}_{2}[r, N]$.

    Therefore, by \eqref{conj.ord1.1_2}-\eqref{conj.ord1.4_2} one has
    \begin{align}
    \partial_t  Y_1 &= - \im \bOmega(D) Y_1 + 
    \textup{Op}_{\tt vec}^{{\scriptscriptstyle{\mathrm BW}}}\Big(\im ( 
    \hspace{-3em}
    \underbrace{(\beta_2)_t + \check{\sfV}_{2}}_{= \la \sfV \ra + \beta_{\geq 3} \mbox{ by Lemma } \ref{lem:coeffV} } 
    \hspace{-3em}
    +g_{\geq 4} + h_{\geq 3})\,  \xi \Big) Y_1  \label{no.v2}\\
    &+ \vOpbw{\ii  \und{\tq}_{\geq 2}  \,  \omega(\xi) + \ii  \und{\tb}_{\geq 2} \sign(\xi) + \ii \und{\wt\tc}_{\geq2 }^{(-\frac12)} } Y_1  + \bB_{\geq N}(\zak) Y_1 + \bR_{\geq 2}(\zak) Y_1\,,\nonumber
    \end{align}
   where
 $\und{\wt\tc}^{(-\frac12)}:= \cP_{\leq N-1}\left[ g_{\geq 2}^{(-2)}+ g^{(-\frac12)}_{\geq 2}+\td_{\geq 2}^{(-\frac12)} \right]$
 is a pluri-homogeneous real-valued symbol in $\Sigma_2^{N-1}\wt\Gamma^{-\frac12}$, 
 $$
 \bB_{\geq N}(\zak):=\check{{\bB}}_{\geq N}(\zak)+\cP_{\geq N}\left[
\vOpbw{\im ( g_{\geq 2}^{(-2)}+ g^{(-\frac12)}_{\geq 2}+ \td_{\geq 2}^{(-\frac12)} )
}
\right]
 $$
 is a real-to-real matrix of  spectrally localized maps in $\cS^{0}_{\geq N}[r]$
and $ \bR_{\geq 2}(\zak)$ is a real-to-real matrix of smoothing operators in  $ \Sigma {\mathcal R}^{-\vr+3N + \frac32}_{ 2}[r,N] $ obtained collecting all the matrices of smoothing operators.

 Equation \eqref{BNF2} follows with 
 $\sfD_{\geq 2}^{({\rm res})}$ as in 
 \eqref{def:D2}, and renaming  the 
real-valued function 
 $\beta_{\geq 3} + g_{\geq 4} + h_{\geq 3} $  in $\Sigma \cF^\R_3[r,N]$ simply as $\und{\sfV}_{\, \geq 3}$.
\end{proof}

 \subsubsection{The weak-$\Lambda$ normal form of water waves}\label{sub:weakL}
In this paragraph we put the cubic part of the smoothing  real-to-real vector field  $\bR_{\geq 2}(\zak)Y_1$ in 
\eqref{BNF2}
in weak-$\Lambda$ normal form, provided $\Lambda$ is a $\gamma$-good set.

\begin{proposition}[{\bf Weak-$\Set$ normal form of water waves}]\label{thm:wnf}

Let $\gamma < 0$, $\gamma^2 \in \Q$ and $\Set = \{\tm, \tn\}$ be a $\gamma$-good set according to Definition \ref{g-good}. Let $N\in \N$, $N \geq 3$ and $\vr>3N+\frac72$.
There exist $s_0, r >0$ and  a $(0,0)$ admissible transformation ${\bPsi}_w(\zak) \in \cM_{\geq 0}^{0}[r]$  with gain $\ov \vr = \vr-{3 N-2}$
such that
 if $(\eta,\psi) \in B_{s_0, \R}(I;r)$ solves \eqref{eq:etapsi} then
the variable 
\be \label{Def:ZW}
 Z:= \bPsi_w(\zak)Y_1 \,,
  \ee   
 with $Y_1$ solving \eqref{BNF2}, solves 
 \be \label{birk}
	 	\pa_t Z= 	 -\ii\bOmega(D)Z+ \vOpbw {\ii \sfD_{\geq 2}^{(\res)}(\zak; x,\x)}Z  +  \bB_{\geq N}(\zak)   Z   +  {\bR}_2^{(\Set)}(\zak)Z  + {\bR}_{\geq 3}(\zak)Z
	 \ee
	 where $\bOmega(D)$ is defined in \eqref{diaglin}-\eqref{omegaxi},    $\sfD_{\geq 2}^{(\res)}(\zak; x, \xi)$ is the real-valued symbol  in \eqref{def:D2} and 
\\
$\bullet$ ${\bR}_2^{(\Set)}(\zak)$ is a real-to-real matrix of smoothing operators in $\wt \cR_2^{-\vr+{3N+\frac32}}$ such that the  cubic vector field $X^{(\Set)}(Z):={\bR}_2^{(\Set)}(Z)Z$
(cf.  \eqref{Xlambda})
 is in weak-$\Set$ normal form (cf. Definition \ref{def:wr});\\
 
\noindent $\bullet$ $\bB_{\geq N}(\zak)$ is a  real-to-real matrix of spectrally localized maps in $\cS_{\geq N}^{0}[r]$ and ${\bR}_{\geq 3}(\zak)$  is a real-to-real matrix of smoothing operators in $\Sigma\cR_{3}^{-\varrho+{3N+\frac72}}[r,N]$ .
\end{proposition}

\begin{proof}
Define the map  ${\bf \Psi}_w(\zak)$ as the time-$1$  flow ${\bf \Psi}_w(\zak):={\bf \Psi}_w^\tau(\zak)_{|_{\tau=1}}$ of
  \be\label{BNFstep2}
\partial_{\tau} {\bf \Psi}_w^{\tau}(\zak)  = \bQ_2(\zak) {\bf \Psi}_w^{\tau}(\zak) \, , 
\quad {\bf \Psi}_w^{0}(\zak) = {\rm Id} \, ,  
\ee
with $ \bQ_2(\zak)$  a $2-$homogeneous real-to-real matrix of smoothing operators in $ \widetilde{\mathcal{R}}^{-\vr +{3N+2}}_2$    to be determined. 
    By Lemma \ref{flow.s.ad} the map ${\bf \Psi}_w(\zak)$ is a $(0, 0)-$admissible transformation with gain $\bar\vr = \vr - {3N-2}$.
    
    Recalling that $Y_1$ solves \eqref{BNF2}, the variable $Z={\bf \Psi}_w(\zak)Y_1$ fulfills
\begin{equation}\label{coniugazione.BNF2}
\begin{aligned}
\pa_t Z
= &  {\bf \Psi}_w(\zak)( -\ii {\bf\Omega}(D) ){\bf \Psi}_w(\zak)^{-1} Z
+ {\bf \Psi}_w(\zak)\vOpbw{ \ii \sfD_{\geq 2}^{(\res)}}{\bf \Psi}_w(\zak)^{-1}Z \\
&+ {\bf \Psi}_w(\zak)\bB_{\geq N}(\zak) {\bf \Psi}_w(\zak)^{-1}Z  + 
      {\bf \Psi}_w(\zak)  \bR_{\geq 2}(\zak)  {\bf \Psi}_w(\zak)^{-1} Z + (\pa_t {\bf \Psi}_w(\zak)) {\bf \Psi}_w(\zak)^{-1}Z \,.
    \end{aligned}
\end{equation}
By  \Cref{prop:Egorov_smoothing}-1 
  one has
   \begin{align}
   {\bf \Psi}_w(\zak)( -\ii {\bf\Omega}(D) ){\bf \Psi}_w(\zak)^{-1}  &=
 -\ii\bOmega(D)+ \big[ \bQ_2(\zak),-\ii \bOmega(D) \big]  + \bR_{\geq 4,1}(\zak) \,,\\
 \label{con1}
   {\bf \Psi}_w(\zak)\vOpbw{ \ii\sfD_{\geq 2}^{(\res)}}{\bf \Psi}_w(\zak)^{-1} & =\vOpbw{ \ii\sfD_{\geq 2}^{(\res)}}  + \bR_{\geq 4,2}(\zak) \,,
\end{align}
where $\bR_{\geq 4,1}(\zak), \bR_{\geq 4,2}(\zak) $ are  real-to-real matrices of smoothing operators  in $ \Sigma\mathcal{R}^{-\vr+ {3N+\frac52}}_{4}[r,N]$ respectively $ \Sigma\mathcal{R}^{-\vr+{3N+3}}_{{4}}[r,N]$.
Now denote by 
$$\bR_2(\zak):= \cP_{\leq 2}\Big[\bR_{\geq 2}(\zak) \Big] ;$$  by  \Cref{prop:Egorov_smoothing}-2  we obtain 
  \begin{align}
  {\bf \Psi}_w(\zak)\,  \bR_{\geq 2}(\zak) \, {\bf \Psi}_w(\zak)^{-1} & =
  \bR_{ 2}(\zak) +  \bR_{\geq 3}(\zak)\,, \\
   {\bf \Psi}_w(\zak) \bB_{\geq N}(\zak) {\bf \Psi}_w(\zak)^{-1}
& =\bB_{\geq N}(\zak) + \bR_{\geq N}(\zak)\,,
 \end{align}
where $\bR_{\geq 3}(\zak), \bR_{\geq N}(\zak)$ are real-to-real matrices of smoothing operators in $\Sigma\mathcal{R}^{-\vr+{3N+2}}_{ 3}[r,N]$ respectively 
$\mathcal{R}^{-\vr+{3N+2}}_{\geq N}[r]$.
Finally,  we consider the contribution coming from the conjugation of $ \pa_t $. 
By \Cref{prop:Egorov_smoothing}-3
 one gets 
 \begin{align}\label{core60}
  (\pa_t {\bf \Psi}_w(\zak)) \, {\bf \Psi}_w(\zak)^{-1} &= 
    \bQ_2(-\ii \bOmega(D) \zak, \zak ) +  \bQ_2(\zak, -\ii \bOmega(D) \zak )+ \bR_{\geq 3}(\zak) 
\end{align}
where $\bR_{\geq 3}(\zak)$ is a real-to-real matrix of smoothing operators in $\Sigma\mathcal{R}^{-\vr+{3N+\frac72}}_{3}[r,N]$.

Combining 
\eqref{coniugazione.BNF2}--\eqref{core60}, 
the variable $Z$ solves 
\be\label{Z.eq.R1}
\begin{aligned}
    \partial_t Z &  =
-\im  \bOmega(D) Z
    + \vOpbw{\im \, \sfD_{\geq 2}^{(\res)} }Z  +\bB_{\geq N}(\zak) Z+ 
    \bR_{\geq 3}(\zak)Z\\
    & \quad + \Big( \bR_{2}(\zak)  + \bQ_2(-\ii \bOmega(D) \zak, \zak ) +  \bQ_2(\zak, -\ii \bOmega(D) \zak ) +
     \big[ \bQ_2(\zak),-\ii \bOmega(D) \big] \Big) Z \,,
\end{aligned}
\ee
where $ \bR_{\geq 3}(\zak)$ is a real-to-real matrix of smoothing operators in $\Sigma \cR_3^{-\vr + 2N + \frac72}[r,N]$, obtained by collecting all the smoothing remainders.

We choose   $\bQ_2(\zak)$ to solve the homological equation 
\be\label{R2w}
  \bR_{2}(\zak)  + \bQ_2(-\ii \bOmega(D) \zak, \zak ) +  \bQ_2(\zak, -\ii \bOmega(D) \zak ) +
     \big[ \bQ_2(\zak),-\ii \bOmega(D) \big]  = \bR^{(\Lambda)}_2(\zak) \ , 
 \ee
where
$\bR^{(\Lambda)}_2(\zak)$ is defined by
 \be\label{RLambda.coeff}
 (\bR^{(\Lambda)}_2(\zak)Z)^\sigma_k:=  \sum_{\fC}  
 R_{j_1, j_2, j, k}^{\sigma_1, \sigma_2, \sigma', \sigma}
\zetina_{j_1}^{\sigma_1} \, \zetina_{j_2}^{\sigma_2} z_{j}^{\sigma'}
\ee
with $ R_{j_1, j_2, j, k}^{\sigma_1, \sigma_2, \sigma', \sigma}$ the Fourier coefficients of $\bR_2(\zak)$, i.e. 
$
( \bR_2(\zak) Z)^\sigma_k = \sum_{ \fP_4
} R_{j_1, j_2, j, k}^{\sigma_1, \sigma_2, \sigma', \sigma}
\zetina_{j_1}^{\sigma_1} \, \zetina_{j_2}^{\sigma_2} z_{j}^{\sigma'}$, 
and the sum is over the indexes $(j_1, j_2, j, k,\sigma_1, \sigma_2, \sigma', -\sigma)$ belonging to the set
\be\label{def:C}
\fC:= \fR^{(0)}_\Lambda \cup \fR^{(1)}_\Lambda \cup\fR^{(2)}_\Lambda \cup  \fP^{(3)}_\Lambda \cup \fP^{(4)}_\Lambda \subset \fP_4 \ , 
\qquad 
\fP^{(n)}_\Lambda, \  \fR^{(n)}_\Lambda  , \ \fP_4 \ \ \mbox{ in }  \eqref{momj}, \eqref{resj}, \eqref{mom1}  \,. 
\ee
Note that $\bR^{(\Lambda)}_2(\zak)$ is a 
 real-to-real matrix of smoothing operators in  $\wt \cR_2^{-\vr+{3N+\frac32}}$ as well as $\bR_2(\zak)$.
 
We claim that equation \eqref{R2w} is solved by
 \be\label{tQ2.def}
( \bQ_2(\zak)Z )^\sigma_k  := \sum_{ \fP_4} Q_{\vec{\jmath} , j,   k}^{\vec \sigma, \sigma',  \sigma} \, \zetina_{j_1}^{\sigma_1} \, \zetina_{j_2}^{\sigma_2} z_{j}^{\sigma'} \,,
\ee
where we denoted  $\vec{\jmath}=(j_1, j_2)$, $\vec \sigma=(\sigma_1, \sigma_2)$, with coefficients $Q_{\vec{\jmath}, j,   k}^{\vec \sigma, \sigma', \sigma}$ defined as
\be\label{tQ2.def2}
 Q_{\vec{\jmath}, j,   k}^{\vec \sigma, \sigma', \sigma}
 := 
 \begin{cases}
  \dfrac{R_{\vec{\jmath}, j,  k}^{\vec \sigma, \sigma',  \sigma}}{\im (\sigma_1 \Omega_{j_1}(\gamma) +\sigma_2 \Omega_{j_2}(\gamma) + \sigma' \Omega_{j}(\gamma) - \sigma \Omega_{k}(\gamma)   )}  \ , & (\vec{\jmath},j,  k, \vec\sigma , \sigma',  -\sigma) \in \fP_4\setminus \fC \\
  0  \ , & (\vec{\jmath}, j,  k, \vec \sigma, \sigma', -\sigma) \in  \fC
 \end{cases} \ , 
 \ee
as we prove in the next lemma.
 \begin{lemma}
  $\bQ_2(\zak)$ in \eqref{tQ2.def}--\eqref{tQ2.def2} is a real-to-real matrix of smoothing operators in  
 $\wt \cR^{-\vr + {3N+2}}_2 $ 
 fulfilling \eqref{R2w}.
 \end{lemma}
 \begin{proof}
 As $\bR_2(\zak)$ is a real-to-real matrix of smoothing operators in $\wt \cR^{-\varrho+{3N+\frac32}}_2$, its coefficients fulfill the estimate: for some $\mu\geq 0$, $C>0$, 
\begin{equation}\label{boundR}
\abs{ R_{\vec{\jmath}, j,  k}^{\vec \sigma, \sigma',  \sigma} } \leq C 
\frac{
{\rm max}_2\{ \la j_1\ra, \la j_2 \ra,  \la j\ra \}^{\mu}}{ \max\{ \la j_1\ra, \la j_2 \ra, \la j\ra \}^{\varrho-{3N-\frac32}}}
   \, , \quad \forall (\vec{\jmath}, j, k, \vec \sigma, \sigma',  -\sigma) \in  \fP_4 \ ,
\end{equation}
 and satisfy the symmetric and reality properties  \eqref{M.coeff.p} and \eqref{M.realtoreal}.
 Consider now the coefficients $ Q_{\vec{\jmath}, j,   k}^{\vec \sigma, \sigma', \sigma}$ in \eqref{tQ2.def2}. 
 Clearly, they satisfy the symmetric and reality properties \eqref{M.coeff.p} and \eqref{M.realtoreal}. 
 We now bound them. By 
\eqref{boundR},  items $(ii)$--$(iii)$ of Lemma \ref{lem:wres} (in particular \eqref{lower.R1},\eqref{lower.R2}) and the momentum relation $\sigma k= \sigma_1 j_1 + \sigma_2 j_2 + \sigma' j$, 
 $$
 \abs{  Q_{\vec{\jmath}, j,   k}^{\vec \sigma, \sigma', \sigma}} \leq C \frac{
{\rm max}_2\{ \la j_1\ra,\la j_2 \ra, \la j_3\ra \}^{\mu}}{ \max\{ \la j_1\ra,\la j_2 \ra, \la j_3\ra \}^{\varrho - {3N-2}}}\qquad  
\forall  (\vec{\jmath} , j,  k, \vec \sigma, \sigma', -\sigma) \in 
\fP^{(1)}_{\Lambda} \cup (\fP^{(2)}_\Lambda \setminus \fR^{(2)}_\Lambda ) \ ,
 $$
 (recall that  $\fR^{(1)}_\Lambda = \emptyset$), 
proving that   $\bQ_2(\zak)$ is a real-to-real matrix of smoothing operators in  $\wt\cR_2^{-\varrho +{3N+2}}$. 
Noting that  $\Pi_{\fP^{(0)}_\Lambda}(\bR_2(\zak)Z) =\Pi_{\fR^{(0)}_\Lambda}(\bR_2(\zak)Z) $  by item $(i)$ of Lemma \ref{lem:wres}, one has that
 $\bQ_2(\zak)$ fulfills \eqref{R2w}.
 \end{proof}
To conclude the proof of \Cref{thm:wnf}, we show that 
the vector field  $X^{(\Lambda)}(\zak)=  {\bR}_2^{(\Lambda)}(\zak)\zak $ is in weak-$\Lambda$ normal form, i.e. it fulfills 
  \eqref{wnf}.
  First note that  the coefficients of the vector field $X^{(\Lambda)}$,   obtained symmetrizing the coefficients in \eqref{RLambda.coeff} with respect to the first three multi-indexes as described in   \eqref{simmetrizzata}, are 
  $$
  X_{j_1, j_2, j_3, k}^{\sigma_1, \sigma_2, \sigma_3, \sigma}= \frac13 \Big(  R_{j_1, j_2, j_3, k}^{\sigma_1, \sigma_2, \sigma_3, \sigma} + R_{j_3, j_2, j_1, k}^{\sigma_3, \sigma_2, \sigma_1, \sigma}+ R_{j_1, j_3, j_2, k}^{\sigma_1, \sigma_3, \sigma_2, \sigma} \Big)  \ , \qquad (j_1, j_2, j_3, k,\sigma_1, \sigma_2, \sigma_3, -\sigma) \in \fC \ 
  $$
 where we used that  the set $ \fC$ is symmetric with respect to the first three multi-indexes.
 Then, by \eqref{def:C}, $\fP^{(n)}_\Set \cap \fC = \fR^{(n)}_\Set$, $ n = 0,1,2$, thus 
$\Pi_{\fP^{(n)}_\Set} X^{(\Lambda)} = \Pi_{\fR^{(n)}_\Set} X^{(\Lambda)}$ proving \eqref{wnf}.
\Cref{thm:wnf} is proved.
\end{proof}
 \subsection{ Proof of Theorem \ref{thm:nf}}\label{subsec:proofthm}
In the previous sections we have performed a paradifferential normal form that allowed us to conjugate the original system \eqref{eq:etapsi} to system 
\eqref{birk}. We shall now prove \Cref{thm:nf}.

\medskip 
\noindent{\sc Proof that $\bUpsilon(\zak)$ is an admissible transformation of order $(\frac12, 8)$ with gain {$\und{\vr}=\vr-3N-2$}:}
By hypothesis $\vr>3N+ 28 $. Fix $\wt \vr >\vr -3N+\frac{11}{2}$.
Define the operator
 \be\label{de.Upsilon}
 \bUpsilon(\zak):= \bPsi_w(\zak) \circ \bPsi_5(\zak) \circ \bT(\zak) \circ \pmb{\cG}_{\C}(\zak) 
 \ee
obtained combining the $(0,0)$-admissible transformation $\bPsi_w(\zak)$ with gain $\vr - 3N-2$  of \Cref{thm:wnf}, 
the $(0,3)$-admissible transformation $\bPsi_5(\zak)$ of gain $\wt\vr>3$  of \Cref{thm:cubic.nf.1}, 
the
$(0,\frac92)$-admissible transformation $\bT(\zak)$ with gain ${\vr - 2N}$  of \Cref{thm:quadratic.nf}, and 
 the $(\frac12, \frac12)$ admissible transformation 
 $ \pmb{\cG}_{\C}(\zak)$ of  gain $\wt \vr>1$ of \Cref{complex:good}.

In view of \Cref{lem:comp},
composing from left to right and using the fact that  $\wt \vr > \max(\vr-3N-2, \frac{15}{2})$, 
$\bUpsilon(\zak)$ is an admissible transformation of order $(\frac12, 8)$ and gain
\begin{align}
 \und{\vr}:=  \min\left(\vr-3N-2, \wt \vr -\frac{15}{2}\right) = \vr - 3N - 2  \geq 8 + \frac{1}{2} \ , 
\end{align}
where in the first passage we have used our largeness assumption on $\wt \vr$ and in the second we have used the one on $\vr$. Note that we need $\und{\vr} \geq 8 + \frac{1}{2}$ so that the map $\bUpsilon(\zak)$ satisfies assumption \eqref{restriction:vr_nu_m} of \Cref{loc.inv}, namely it is $(\frac 1 2, 8)$ admissible with gain $\und{\vr} \geq 8 + \frac 1 2$.

\medskip 
\noindent{\sc Proof that $Z(t)$ fulfills equation \eqref{Z.eq}:}
Define the nonlinear map
 \be
Z :=  \scF(\zak):=\bUpsilon(\zak)\zak ,  \quad  \bUpsilon(\zak) \in \cM^{\frac12}_{\geq 0}[r] \mbox{ in } \eqref{de.Upsilon}\,.
\label{de.Upsilon2}
 \ee
\Cref{loc.inv} ensures that there is $s_0'\geq 0$ such that for any $s \geq s_0'$, there is $r' = r'(s)>0$ and a local inverse map so that the nonlinear inverse map 
 $\scF^{-1}: B_{s_0'}(r') \cap \dot H^{s+\nu}_\R(\T;\C^2) \to \dot H^{s}_\R(\T;\C^2)$, having the  structure 
 \be\label{inv.F}
 \zak = \scF^{-1}(Z) = \bXi(Z) Z , \quad  \mbox{ with } 
 \bXi(Z) \in \cM_{\geq 0}^{\frac12}[r''] \mbox{ and }
 \bXi(Z) -  \uno \in \Sigma\cM^{ {24 + \frac 1 2}}_{ {1}}[r'', 3] \ , 
 \ee
 for some $r'' >0$.
 
 We then substitute $\zak$ in the internal variables of the operators in \eqref{birk}. 
 Consider first the symbol $\sfD_{\geq 2}^{(\res)}(\zak; x,\x)$ whose expression is given in \eqref{def:D2}.
 We have,  using  \Cref{lem:sostitution}--1, 
 \begin{multline}
    \vOpbw{\ii \sfD_{\geq 2}^{(\res)}(\scF^{-1}(Z); x, \xi)} = \bR'_{\geq 3}(Z)\\
     + \textup{Op}_{\tt vec}^{{\scriptscriptstyle{\mathrm BW}}}\Big(\im  
    \underbrace{\big(\langle \sfV \rangle (Z; x) \xi + \sfV_{\geq 3}(Z; x)\xi  + \td_{\geq 2}(Z; x)\omega(\xi) +  \tf_{\geq 2}(Z;x) \sign(\xi) +  \tg_{\geq2}^{(-\frac 1 2)}(Z; x, \xi)\big)}_{=:\sm^{(\res)}_{\geq 2}(Z;x) }\Big)\,,
 \end{multline}
where $\bR_{\geq3}'(Z)$ is a real-to-real matrix of smoothing operators in $\cR^{-\vr}_{\geq 3}[r'']$ and the real-valued symbols 
\begin{gather*}
{\sfV}_{\geq 3}(Z; x) := \langle \, \sfV\, \rangle(\scF^{-1}(Z);x)    -    \langle \, \sfV \, \rangle(Z;x)  + \und{\sfV}_{\geq 3}(\scF^{-1}(Z); x) \in  \cF^\R_{\geq 3}[r'']\,,\\
     \td_{\geq 2}(Z; x) := \und{\tq}_{\geq 2}(\scF^{-1}(Z);x) \in \Sigma\cF^\R_2[r'',3] \,, \qquad \tf_{\geq 2}(Z; x) := \und{\tb}_{\geq 2}(\scF^{-1}(Z); x)\in \Sigma\cF^\R_2[r'',3]\,,  \\
     \tg^{(-\frac 1 2)}_{\geq 2}(Z; x, \xi) := \und{\wt{\tc}}^{(-\frac 1 2)}_{\geq 2}(\scF^{-1}(Z); x, \xi) \in \Sigma\Gamma^{-\frac12}_2[r'',3] 
 \end{gather*}
Then using that, by item $(ii)$ of \Cref{loc.inv}, we have $ \| \scF^{-1}(Z)\|_{s}\lesssim \| Z\|_{s+\frac12}$ for any $s \geq s_0'$, we also get that
 $$
 \bB_{\geq N}(\scF^{-1}(Z)) \text{ is a matrix of real-to-real spectrally localized maps in } \cS^0_{\geq N}[r'']\,.
 $$
Then, using  \Cref{lem:sostitution}-2 and the fact that ${R}_2^{(\Lambda)}(Z) \in \wt \cR_2^{-\vr+3N+\frac32}$,\  ${\bR}_{\geq 3}(\zak) \in \Sigma\cR_{3}^{-\varrho+3N+\frac72}[r, N]$, (see \Cref{thm:wnf}) we obtain
 \begin{align*}
 {R}_2^{(\Lambda)}(\scF^{-1}(Z)) - {R}_2^{(\Lambda)}(Z) \in \cR_{\geq 3}^{-\vr+3N+26}[r'']\,, \qquad 
\bR_{\geq 3}(\scF^{-1}(Z)) \in \cR_{\geq 3}^{-\varrho  +3N +  28}[r'']\,.
 \end{align*}
 In conclusion we obtain that $Z(t)$ solves system \eqref{Z.eq}.

\medskip 
\noindent{\sc Proof that ${X}^{(\Lambda)}(Z) $ is in strong-$\Lambda$ normal form:}
We  exploit an identification argument to prove that   ${X}^{(\Lambda)}(Z) $ in \eqref{Xlambda} is actually in strong-$\Lambda$ normal form and compute explicitly 
$\Pi_{\fP^{(n)}_\Lambda}{X}^{(\Lambda)}$, $ n = 0,1,2$. 
We start with the following result
 about the projections on $\fR^{(n)}_{\Set}$, $n=0,1,2$, of cubic paradifferential vector fields.
\begin{lemma}\label{lem:proj.para}
Let $\gamma <0$, $\gamma^2\in \Q$ and  $\Set = \{\tm, \tn\}$ be a $\gamma$-good  set of the form \eqref{Lambda}.  Let $a(Z; x, \xi)$ be a real-valued, 2-homogeneous symbol in $\wt \Gamma_2^m$, $m \in \R$. Then 
 \be\label{ident.Op.ab}
\Pi_{\fR^{(n)}_{\Set} } \left[ \vOpbw{ \im a(Z; x, \xi)} Z \right] = 0 \ , \quad n = 0,1 \ ,
\ee
whereas 
$$\cD(Z):=\Pi_{\fR^{(2)}_{\Set} } \left[ \vOpbw{ \im a(Z; x, \xi)} Z \right]$$
is an integrable vector field fulfilling
\be\label{ident.Op.ab2}
  \Re (\cD(Z)^+_k \, \bar{z_k}) = 0 \quad  \forall k \in \Z_* \,.
 \ee
\end{lemma}
\begin{proof}
Recalling \eqref{vecop}, $\left( \vOpbw{ \im a(Z; x, \xi)} Z \right)^+ =  \Opbw{ \im a(Z; x, \xi)} z$. 
Using definition  \eqref{BW} specialized to quadratic  symbols, we get 
\begin{equation}\label{para2}
 \Opbw{ \ii  a(Z; x, \xi)} z = \sum_{ \sigma_1 j_1+ \sigma_2 j_2 +j=k} \im \chi_2 \left( j_1, j_2,\frac{j+k}{2}\right) a_{j_1, j_2}^{\sigma_1, \sigma_2}\left(\frac{j+k}{2}\right) z_{j_1}^{\sigma_1} z_{j_2}^{\sigma_2} z_j
{e^{\im k x}} \,.
\end{equation}
Recall that $\chi_2(\xi', \xi) \equiv 0$ whenever $|\xi'| > \delta_0 \la \xi \ra$ and $\delta_0 \ll 1$.   \\
\underline{Case $n = 0$:} 
In this case, recalling that $\delta_0{<\frac{1}{10}}$, we claim that  
\begin{equation}
    \chi_2 \left( j_1, j_2,\frac{j+k}{2}\right) = 0, \qquad \forall j_1, j_2, j, k\in \Set\,.
    \label{supporto4dentro}
\end{equation} 
 Indeed, as $\Lambda$ is $\gamma$-good, it has the form $\Lambda=\{\tm,\tn\}$ for some $\tn,\tm$ satisfying \eqref{Lambda} and \eqref{Omega*}--\eqref{G4}. The bounds in \eqref{Lambda} and the bound \eqref{G3'} in \Cref{lem:G3'} imply that 
\be\label{ciccio}
\tn>0, \qquad|\tm|< \tn\leq \upsilon|\tm|\leq 3 |\tm|\,.
\ee
Therefore, since $\min\{|j_1|,\, |j_2|\} \geq |\tm|$, $|j|+|k|\leq 2\tn$ and recalling that the number  $\delta_0 >0$ appearing in the cutoff in \eqref{supp.chi} satisfies $\delta_0\leq\tfrac{1}{10}$, on the support of $\chi_2 \left( j_1, j_2,\frac{j+k}{2}\right)$ one has
$$
|\tm|\leq \min\{|j_1|,\, |j_2|\} \leq \frac{\delta_0}{2} ( |j|+|k|)\leq  \delta_0 \tn \stackrel{\eqref{ciccio}}{\leq} 3\delta_0 |\tm|\leq \frac{3}{10}|\tm| \,, 
$$
which, together with $ |\tm|\geq 1$, yields a contradiction and proves the claim in \eqref{supporto4dentro}.

\noindent 
\underline{Case $n = 1$:} By Lemma \ref{lem:wres} $(ii)$, $\fR_\Set^{(1)} = \emptyset$ and there is nothing to prove.\\
\underline{Case $n = 2$:}  By Lemma \ref{lem:wres} $(iii)$, the indexes  $j_1,  j_2, j, k$ in $\fR^{(2)}_{\Set}$ are pairwise equal. 
In case $j_1 = j $ and $j_2 = k $ 
the cut-off vanishes since
$$
\chi_2\left(j_1, j_2, \frac{j+k}{2}\right)  = \chi_2 \left( j_1, k,\frac{j_1  + k}{2}\right)  \equiv  0 \quad \forall j_1, \ 
k \in \Z_*\,. 
 $$
The only case left is 
 $j_1 = j_2$ and $j = k$. By Lemma \ref{lem:wres}-$(iii)$ then we have $ \sigma_1 = -\sigma_2$;  in this case, using also \eqref{sym_sy}, we have
$$
\cD(Z)^+ = \sum_{(j_1,k)\in ( \Set\times  \Set^c)  \cup (\Set^c\times \Set)} 
\chi_2 \left( j_1, j_1,k\right) \im 2 a_{j_1, j_1}^{+, -}\left(k\right) |z_{j_1}|^2 \,  z_k
{e^{\im k x}}  \,, 
$$
 thus $\cD(Z)^+_k  \bar{z_k} = \sum_{j_1 } \chi_2 \left( j_1, j_1,k\right) \im \,2 a_{j_1, j_1}^{+, -}\left(k\right) |z_{j_1}|^2 \,  |z_k|^2 $ is purely imaginary since $a_{j_1, j_1}^{+, -}\in \R$ by \eqref{sym_sy} and \eqref{reality_cond}. Thus  \eqref{ident.Op.ab2} follows.
\end{proof} 
We are finally able to  prove that ${X}^{(\Lambda)}(Z) $ fulfills  the claimed properties in  \eqref{Y.str}.  This will be achieved thanks to the 
 the abstract identification argument of
 Proposition  \ref{prop.ident} and 
 the explicit form of the cubic  
  Hamiltonian vector field $\cX_{\cH_4^+}$  computed in 
 \Cref{prop:cH^{(4)}_+}.
 \begin{proposition}\label{prop.megliodemax2}
 The vector field ${X}^{(\Lambda)}(Z) $ in \eqref{Xlambda} is actually in strong-$\Lambda$ normal form (Definition \ref{def:wr}) and fulfills   \eqref{Y.str}. 
 \end{proposition}
\begin{proof}
 \Cref{thm:nf} guarantees that the variable $Z = \bUpsilon(\zak) \zak$ satisfies equation \eqref{Z.eq}, which has the form \eqref{eq:zetone} with the cubic vector field $$
 X_3^+(Z) := \vOpbw{\im  \sm_{2}^{(\res)}(Z; x, \xi)}Z
+  {X}^{(\Set)}(Z)\,, \qquad \sm_{2}^{(\res)}(Z; x, \xi) := \cP_2 \left[\sm_{\geq 2}^{(\res)}(Z; x, \xi) \right]\,,
$$
with $\sm_{\geq 2}^{(\res)}(Z; x, \xi)$ in \eqref{D.final}.
We are then in position to apply \Cref{prop.ident}, from which we deduce the identity (cf. \eqref{ide:res})
 \be\label{ident.X3}
\Pi_{\fR^{(n)}_{\Set} }
\left[ 
\vOpbw{ \im  \sm_{2}^{(\res)}(Z; x, \xi)}Z
+  {X}^{(\Set)} \right]
=
\Pi_{\fR^{(n)}_{\Set} } \cX_{\cH_4^+}   \ , \quad n = 0,1,2 \ , 
\ee
where $\cX_{\cH_4^+}$ is the Hamiltonian vector field of $\cH^+_4$ in \eqref{new.birkh}.\\ Consider first the cases $n=0, 1$.
We apply  Lemma \ref{lem:proj.para}  to the cubic vector field 
$\vOpbw{ \im  \sm_{2}^{(\res)}(Z; x, \xi)}Z$, getting that 
\begin{align}\label{ident.Op0}
& \Pi_{\fR^{(n)}_{\Set} } \left[ \vOpbw{ \im  \sm_{2}^{(\res)}(Z; x, \xi)}Z\right] = 0 \ ,  \quad n = 0,1 \ ,
\end{align}
from which, recalling that $X^{(\Lambda)}$ is in weak-$\Set$ normal form by \Cref{thm:wnf}, we deduce 
\begin{align*}
\Pi_{\fP^{(n)}_{\Set} }{X}^{(\Set)}  
& \stackrel{\eqref{wnf}}{=} 
\Pi_{\fR^{(n)}_{\Set} }{X}^{(\Set)} 
\stackrel{\eqref{ident.X3},  \eqref{ident.Op0}}{=}  \Pi_{\fR^{(n)}_{\Set}} {\cX}_{\cH_4^+}  \ ,  \quad n = 0,1 \,.
\end{align*}
By \Cref{prop:cH^{(4)}_+} the vector field $\cX_{\cH_4^+}$ is in strong-$\Set$ normal form, and its projections are computed in \eqref{H4+.str0}--\eqref{H4+.str}, proving that $X^{(\Set)}$ fulfills    the first and second of \eqref{Y.str}.
\\
It remains to consider the case $n = 2$.  Denoting 
\be\label{ident.Op2}
\cD(Z):=  \Pi_{\fR^{(2)}_{\Set} } \left[ \Opbw{ \im   \sm_{2}^{(\res)}(Z; x, \xi)} Z \right] \,,  
\ee
 we have 
\begin{align}
\Pi_{\fP^{(2)}_{\Set} }{X}^{(\Set)}  
& \stackrel{\eqref{wnf}}{=} 
\Pi_{\fR^{(2)}_{\Set} }{X}^{(\Set)} 
\stackrel{\eqref{ident.X3},  \eqref{ident.Op2}}{=}  \Pi_{\fR^{(2)}_{\Set}} {\cX}_{\cH_4^+}  - \cD(Z) \,. 
\label{2506:1412}
\end{align}
By   Lemma \ref{lem:proj.para},  
$
\Re (\cD(Z)^+_j \, \bar{z_j}) = 0$ for any  $j \in \Z_*$, but since also 
$$2 \Re \left( \left[\Pi_{\fR^{(2)}_{\Set}} {\cX}_{\cH_4^+}\right]_j^+ \bar{z}_j \right) 
\stackrel{\eqref{PP},\eqref{hvf}}{=}
 \{ \Pi_{\fR^{(2)}_{\Set}}\cH_4^+, |z_j|^2\}\stackrel{\eqref{HR12}}{=} 0 , 
$$
we deduce  the last of  \eqref{Y.str}.
\end{proof}
\section{The Effective Equation}\label{sec:eff.eq}
In this section we consider some special solutions, that we call {\em long-time controlled},   which are characterized by having a large apriori control on their $H^s$ norm and special smallness requirements on the  $L^2$-norm of the initial data, see Definition \ref{A}.
The goal of this section is to prove that  any long-time controlled solution fulfills the effective equation \eqref{eff.sys}. In the next section we shall prove that there exist long-time controlled solutions on timescales longer than those guaranteed by the local Cauchy theory.
\smallskip 

In order to derive the effective equation, we start regarding system \eqref{Z.eq} as a scalar equation in $z(t)$, thanks to the fact that the vector field at the right-hand side of  \eqref{Z.eq} is real-to-real. 
From now on we fix: \\
{\bf The set}
\be\label{fixLambda}
\Set = \{\tm, \tn\}\subset \Z_*
\ , \quad \tm < 0 < \tn \ , \quad \tm + \tn > 0 \ , \quad \gamma\mbox{-good}
\qquad  \mbox{(cf. \Cref{g-good})}
\ee
  and  {\bf the parameters}
\be\label{paraN}
N = 3\,, \quad \vr = 3N + 28 = 37\,, \quad \und{\vr} = \vr - 3 N - 2 = 26
\ee
in \Cref{thm:nf}, and define the parameters 
\be \label{s0r}
\begin{aligned}
& \fs_0\geq \max\{ s_0, \s_0', \s_0''\} +3 \,,  
\quad \fr:=\min\{r, \tr', \tr \}\,, \quad s> 4 \fs_0 \ , \\
& \theta \in (0, \theta_*), \quad \theta_* := \min\left\{\frac{s-4 \fs_0}{2(2s - \fs_0)},  \frac16 \right\} \ , \quad \delta_0 < \min \left( \frac{1}{10}, \frac{1}{2 \tn + 3} \right) \ , 
\end{aligned}
\ee
where $s_0, r$ are given in Theorem \ref{thm:nf} 
-- and $s_0$ is also the largest minimal regularity for which all the lemmas in Sections \ref{subsec:preliminary}-- hold--
whereas  $\s_0',\s_0''$, $\tr':= \tr'(\fs_0), \tr:= \tr(\fs_0)$ are respectively the low regularity thresholds and the low regularity radii in \Cref{UXequiv}, \Cref{UZequiv}. 
Finally  $\delta_0$ is the parameter appearing in the cutoff functions in \eqref{supp.chi} and in Definition \ref{quantizationtotale}.

\smallskip
We introduce the projectors $\Pi^\top$ on the set $\Set$ and  $\Pi^\perp$ on the complementary set $\Set^c$, namely 
\be\label{Pi}
z^\top:=\Pi^\top z := \sum_{j \in \Set} z_j \, e^{\im j x} = z_\tm e^{\ii \tm x} + z_\tn e^{\ii \tn x}  \ , \qquad z^\perp:=\Pi^\perp z := \sum_{j \not \in \Set} z_j \, e^{\im j x} \  
\ee
and use them to 
decompose any solution $z(t)$ as 
\be
\begin{aligned}\label{ztp}
z(t) = z^\top(t) + z^\perp(t) \ , \quad
&  z^\top(t) := \Pi^\top z(t) \ , \quad
  z^\perp(t):= \Pi^\perp z(t) \,.
\end{aligned}
\ee

In the next definition we describe the first  condition on the initial datum that we require.

\begin{definition}[{\bf Weakly well-prepared initial data}]\label{def:WWP}  
 Let $\e>0$. We say that an initial datum $z(0) \in L^2(\T;\C)$
is $\e$-{\em weakly well-prepared} if
\begin{equation}
\label{app.ass1}
\norm{ z^\top(0)}_{L^2} <  \e    \,, \qquad
\norm{ z^\perp(0)}_{L^2} < \e^{2}   \,, 
\qquad 
z^\top(0), \ z^\perp(0) \mbox{ in } \eqref{Pi} \,. 
\end{equation}
\end{definition}
Next we  introduce the notion of long-time controlled solutions. 
\begin{definition}[{\bf Long-time controlled solutions}]\label{A}
Let $\fs_0,\fr, s, \theta$ as in \eqref{s0r}.
Let also $T_\star>0$ and $ \e \in (0,\fr)$. 
We say that a solution $z(t) $ of system  \eqref{Z.eq} is {\em long-time controlled}  with parameters $(s, \theta, T_\star, \e)$ if 
\begin{itemize}
\item[$(A1)$] its initial datum $z(0)$ is $\e$-weakly well-prepared (\Cref{def:WWP});
\item[$(A2)$] $z(t)$
 belongs to the space $C^0([0, T_\star]; \dot H^s(\T;\C))$ 
and fulfills  
the large  a-priori bound
\begin{equation}
\label{app.ass2}
\sup_{0 \leq t \leq T_\star } 
\norm{z(t)}_s \leq  \e^{-\theta}    \,. 
\end{equation}
\end{itemize}
\end{definition}
\begin{remark}
  By the  Cauchy theory,  any initial data  $z(0)\in B_s(\e)$ which is also $\e$-weakly prepared gives rise to a solution which is long-time controlled with parameters $(s, \theta, T_\star, \e)$ and $T_\star \simeq \e^{-2}$. 
  The nontrivial fact is the existence of  long-time controlled solutions over longer time-scales, which are those needed to observe growth of Sobolev norms. 
  We shall prove this in \Cref{lem:cresce}. 
\end{remark}
As mentioned above, the  key property of long-time controlled solutions is that they  fulfill an {\em effective equation} with a very precise structure, as we now prove. 
\begin{proposition}[{\bf Effective equation}] \label{prop:eff}
Let
 $T_0 >0$, and let $s, \theta, \fr$ as in \eqref{s0r}. 
There exists $\e_\star = \e_\star(s, \theta, T_0, \fr) >0$ such that for any $\e \in (0, \e_\star)$ the following holds true. 
Let $z(t)$ be a solution  of  \eqref{Z.eq} which is long-time controlled  with parameters $(s, \theta, T_\star, \e)$ and with
\be\label{Tstar.T0}
T_\star\leq \frac{T_0}{\e^2}\log\left(\frac{1}{\e} \right) \,.
\ee
Then   $z(t, x) = z_{\tm}(t)e^{\ii \tm x} + z_\tn(t) e^{\ii \tn x} +  z^\perp(t, x) $ fulfills the system
\be\label{eff.sys}
\begin{cases}
 \pa_t z_\tm = - \im \big( \Omega_{\tm}(\gamma) +  2 \fa |z_\tm(0)|^2 + \fb|z_\tn(0)|^2\big) z_\tm  + \sfd_\tm(t) \\
 \pa_t z_{\tn}  = - \im (\Omega_\tn(\gamma) + 2\fc |z_{\tn}(0)|^2 + \fb |z_{\tm}(0)|^2 \big)  z_{\tn} + \sfd_{\tn}(t) \\
     \pa_t z^\perp  = - \im {\Omega(D)} z^\perp  + \im \Opbw{\fv(x - \tJ t) \xi + \sfW(t;x) \xi + \sfb(t;x, \xi)}z^\perp + \sfR(t) + {\sfY}(t)
\end{cases}
\ee
where $\fa, \fb , \fc$ are the numbers  in  \eqref{abc} and \\
$\bullet$ $\tJ$ is the real number
 \begin{equation}\label{J1}
\tJ:=  \frac{2\fc - \fb}{\tn - \tm}|z_\tn(0)|^2 + \frac{\fb-2\fa}{\tn-\tm}|z_\tm(0)|^2 \ ;  
\end{equation}
$\bullet$ the real-valued function  $ \fv(x)$ is given by
\begin{equation}\label{V0}
 \fv(x) := 
 \sfV^{(\rm int)}_\tm |z_\tm(0)|^2 + \sfV^{(\rm int)}_\tn |z_\tn(0)|^2 +
 2 \sfV^{(\rm res)}_{\tm, \tn} \, \Re \left(z_\tn(0) \, \bar{z_{\tm}(0)} \, e^{\im (\tn-\tm) x} \right) 
 \end{equation}
 with $\sfV^{(\rm int)}_\tm$, $\sfV^{(\rm int)}_\tn $, $\sfV^{(\rm res)}_{\tm, \tn}$ in \eqref{Vres.coeff};\\
$\bullet$ the real-valued, time dependent function $ \sfW(t;x)$ fulfills the estimate
\begin{equation}\label{Vpert}
\norm{ \sfW(t; \cdot)}_{W^{2,\infty}}  \leq C \e^{{3}- \theta} \ , \quad
\forall \, 0 \leq t \leq T_\star \,;
\end{equation}
$\bullet$ the real-valued symbol $\sfb(t;x , \xi) \in \Gamma^{\frac 1 2}_{W^{2, \infty}}$, recalling \eqref{seminorm}, fulfills the following estimate: for every $n \in \N_0$, there is $C_n >0$ such that 
\begin{equation}\label{bpert}
 \, \abs{ \sfb(t; \cdot) }_{\frac 12, W^{2, \infty}, n}  \leq C_n \e^2 \,, \quad
 \forall\ 0 \leq t \leq T_\star \,;
\end{equation}
$\bullet$ the functions $\sfd_{\tm}(t)$, $\sfd_\tn(t)$ fulfill the estimates
\be\label{d.bound}
\abs{\sfd_{j}(t)} \leq \e^{4-\theta}   \,, \quad
 \forall\ 0 \leq t \leq T_\star \ \qquad j \in \{\tm, \tn\} \,;
\ee
$\bullet$  $\sfR(t) \equiv \sfR(t, x)$ is a  smoothing vector field fulfilling  the estimate 
\begin{equation}\label{R(t,U)}
\norm{\sfR(t; \cdot)}_{s+1 } \leq C  \e^{2-\theta} \,, \quad
 \forall \, 0 \leq t \leq T_\star \,;
\end{equation}
$\bullet$ $\sfY(t) \equiv \sfY(t, x)$ is a bounded vector field fulfilling  the estimate 
\begin{equation}\label{Y(t,U)}
\norm{\sfY(t; \cdot)}_{s } \leq C  \e^{3-\theta} \ , \quad
 \forall \, 0 \leq t \leq T_\star \,.
\end{equation}
\end{proposition}
We pause to note that the  key features of \eqref{eff.sys} are the following:  
\begin{itemize}
    \item The dynamics of $z_\tm(t)$ and $z_\tn(t)$ is integrable up to the terms $\sfd_\tm(t)$, $\sfd_\tn(t)$ of much smaller size $\e^{4-\theta}$;
    \item the dynamics of $z^\perp(t)$ is led by the linear transport operator $ \im \Opbw{\fv(x - \tJ t) \xi }$ of size  $\e^2$, up to the linear operator $- \im \Omega(D)$ of larger size but lower order, a  para-differential operator with a lower order symbol $\sfb(t)$ of size $\e^2$,
    a order 1 symbol $\sfW(t,x)$ of smaller size $\e^{3-\theta}$, a vector field  $\sfR(t)$ of larger size $\e^{2-\theta}$ but more regular, and a bounded vector field $\sfY(t)$ of smaller size $\e^{3-\theta}$.
\end{itemize}
Before proving Proposition \ref{prop:eff}, we prove an immediate but important consequence. Note that the transport operator  $\fv(x-\tJ t) \xi$ in \eqref{eff.sys} is translating in space at a constant speed $\tJ$.
Therefore it is convenient to  put ourselves in a moving frame with  constant speed $\tJ$, in which the transport term is time independent and of the form \eqref{J+v}; the crucial sign property  \eqref{VmVn.def} --  fundamental to prove the positive  lower bound in \Cref{lem:pos.comm}--  
is a consequence of (G3) in  Definition \ref{g-good} of $\gamma$-good set.
\begin{corollary}\label{cor:dyn.eff}
    Let $s, \theta$ as in \eqref{s0r}. 
Fix also 
 $T_0 >0$.   
There exists $\e_\star = \e_\star(s, \theta, T_0,\fr) >0$ such that for any $\e \in (0, \e_\star)$ the following holds true. 
Let $z(t)$ be a solution  of  \eqref{Z.eq} which is long-time controlled  with parameters $(s, \theta, T_\star, \e)$  (see Definition \ref{A}) and with $T_\star$ fulfilling 
 \eqref{Tstar.T0}. 
 The variable 
\be
{
\zeg(t,x):= z^\perp \big( t,  x +\tJ t \big) }\ , \quad \tJ \mbox{ in } \eqref{J1} \ ,
\label{zeta}
\ee
with norm  $\norm{\zeg(t, \cdot)}_s = \norm{z^\perp(t, \cdot)}_s\ $ for any $t, s \in \R$, fulfills the equation
\begin{align}\label{eq.zeta}
 \pa_t \zeg  = &  - \im \Omega(D)\zeg   + \im \Opbw{{(\tJ+\fv(x ))} \xi }\zeg + \im \Opbw{ \wt \sfW(t; x) \xi + \wt \sfb(t; x,  \xi)}\zeg  + \wt{\sfR}(t)  + \wt{\sfY}(t)
 \end{align}
where\\
$\bullet$ the transport term 
\begin{equation}\label{J+v}
    \tJ+\fv(x )=  \tV_\tm |z_\tm(0)|^2 + \tV_\tn |z_\tn(0)|^2 +
 2 \sfV^{(\rm res)}_{\tm, \tn} \, \Re \left(z_\tn(0) \, \bar{z_{\tm}(0)} \, e^{\im (\tn-\tm) x} \right) 
\end{equation}
with $\sfV^{(\rm res)}_{\tm, \tn}$ given by \eqref{Vres.coeff} and 
    \be\label{VmVn.def}
\!\!\!\!\!\!\!\! \!\!\!\!\!\!\!\! \!\!\!\!\!\!\!\!  
\mbox{{\sc (Sign property)}} \qquad \quad 
\begin{aligned}
& \tV_\tm:= \frac{2 \tm^6-25 \tm^5 \tn+50 \tm^4 \tn^2+8 \tm^3 \tn^3-4 \tm^2 \tn^4+\tm \tn^5}{4 (3 \tn-\tm) (\tn-\tm)^3} >0  \ , \\
&
 \tV_\tn:=   \frac{\tn \left(19 \tm^4 - 33 \tm^3 \tn - 53 \tm^2 \tn^2 + 29 \tm \tn^3 - 2 \tn^4\right)}{4 (\tn-\tm)^2 (\tn-3\tm)} < 0\ ; 
 \end{aligned}
\ee
$\bullet$  the real-valued function $\wt \sfW(t; \cdot) \in W^{2,\infty}$, the real-valued symbol  
$\wt \sfb(t; \cdot) \in \Gamma^{\frac12}_{W^{2,\infty}}$ and the vector fields $\wt{\sfR}(t;\cdot) \in H^{{s+1}}(\T; \C)$, $ \wt{\sfY}(t; \cdot) \in H^{s}(\T; \C)$ 
fulfill, for any $0 \leq t \leq T_\star$, the estimates 
\be\label{est.706}
\begin{aligned}
& \norm{ \wt\sfW (t; \cdot)}_{W^{2,\infty}}  \leq  C \e^{{3}- \theta} \,, 
\quad
\vert \wt \sfb(t; \cdot) \vert_{\frac 1 2, W^{2, \infty}, n}  \leq C_n \e^2  \,, \ \ \forall n \in \N \,,  \\
&
\norm{\wt{\sfR}(t; \cdot)}_{s+1 } \leq C_s  \e^{2-\theta} \,, 
\quad 
\norm{\wt{\sfY}(t; \cdot)}_{s } \leq C_s  \e^{3-\theta}  
 \,.
\end{aligned}
\ee
\end{corollary}
\begin{proof}
After the translation \eqref{zeta},  the variable $\zeg(t, \cdot)$ fulfills \eqref{eq.zeta} with 
$\tJ + \fv(x)$ of the form \eqref{J+v} with (see \eqref{eff.sys}, \eqref{J1}, \eqref{V0})
\be
\begin{aligned}\label{VmVn}
& \tV_\tm:= \left(\sfV^{(\rm int)}_{\tm} + \frac{\fb-2\fa}{\tn-\tm}\right)  \,, \quad 
 & \tV_\tn:= \left(\sfV^{(\rm int)}_{\tn} + \frac{2\fc - \fb}{\tn-\tm}\right) \,, 
\end{aligned}
\ee 
and 
the real-valued function $\wt \sfW(t; x)$, the real-valued symbol  
$\wt \sfb(t; x,  \xi)$ and the vector fields $ \wt{\sfR}(t;x)$, $ \wt{\sfY}(t; x)$ given by  
\be
\wt \sfW(t; x) := \sfW(t, x+ \tJ t) \ , \quad 
\wt \sfb(t; x,  \xi):= \sfb(t; x+ \tJ t,  \xi) \,, \quad
 \wt{\sfR}(t; x) :=   \sfR(t; x+ \tJ t)  \,, \quad 
 \wt{\sfY}(t; x) :=   \sfY(t; x+ \tJ t) \,. 
\ee
The explicit formulas for $\tV_\tm$, $\tV_\tn$ follow from \eqref{VmVn} substituting  the values of $\fa, \fb, \fc$ in \eqref{abc} and of $\sfV^{(\rm int)}_{\tm}, \sfV^{(\rm int)}_{\tm}$ in \eqref{Vres.coeff}, see also
the  Mathematica notebook \texttt{corollary6\_5.nb} at the link \url{https://git.sissa.it/amaspero/transfer-ww-vorticity}.

The functions  $\tV_\tm$ and $\tV_\tn$ have the sign of their numerators $\tN(\tV_\tm),\  \tN(\tV_\tn)$ since $\tm < 0 < \tn$; to study them, note that
\[
\begin{aligned}
   &  \tN(\tV_\tm) = \tn^6 p\left(\frac{\tm}{\tn}\right) \ , \quad p(\lambda):=   2 \lambda^6-25 \lambda^5 +50 \lambda^4 +8 \lambda^3 -4 \lambda^2 +\lambda  \,, \\
     &  \tN(\tV_\tn) = \tn^5 q\left(\frac{\tm}{\tn}\right) \ , \quad q(\lambda):=  19 \lambda^4 - 33 \lambda^3  - 53 \lambda^2  + 29 \lambda  - 2 
     \,.
\end{aligned}
\]
Since the variable $\frac{\tm}{\tn} <0$, we are interested in their signs only for $\lambda <0$. Note that, since the set $\Lambda= \{\tn, \tm\}$ is $\gamma-$good, condition \eqref{G3} guarantees that $p(\tfrac{\tm}{\tn}) >0$. Moreover, a simple analysis (see again the notebook \texttt{corollary6\_4.nb}) gives that 
$$
q(\lambda) < 0 \Leftrightarrow   - {\und \upsilon } < \lambda < 0 \ , \ \ \  \und \upsilon \approx 1.3146
$$
which is fulfilled by the variable $\frac{\tm}{\tn} $ provided 
$- {\und \upsilon}^{-1} \tm  <\tn$.
Such condition is verified due to the fact that $\Set$ is $\gamma$-good; indeed, being $ {\und \upsilon}^{-1}\in (0,1)$ and $\tm < 0$, one has 
$\tn + {\und \upsilon}^{-1}\tm > \tn + \tm>0 $ (see \eqref{Lambda}).
\\
The estimates \eqref{est.706}   follow  by \eqref{Vpert}, \eqref{bpert}, \eqref{R(t,U)} and \eqref{Y(t,U)}.
\end{proof}
The remaining part of the section is devoted to the proof of Proposition \ref{prop:eff}.
\subsection{Dynamics of   $z^\top$ and $z^\perp$ variables}
In this subsection we compute the equations fulfilled by $z^\top$ and $z^\perp$ by projecting the vector field \eqref{Z.eq}. In particular, we 
 exploit the fact that the vector field in \eqref{Z.eq} is in strong $\Set$-normal form to prove that the projections of the cubic vector field 
 $X^{(\Set)}(z)$ are  the sum of an integrable vector  field -- that we denote by $X^{(\intt)}(z)$ respectively 
 $Y^{(\intt)}(z)$ --
 and a  vector field which is at least quadratic in $z^\perp$ -- that we denote by $X_{\pperp}(z)$ respectively $Y_{\pperp}(z)$. 

 Note that the analysis in this subsection remains valid independently of whether 
$z(t)$ is a long-time controlled solution. However, in the next subsection, we will demonstrate that, along a long-time controlled solution $z(t)$,  the vector fields  $X_{\pperp}(z(t))$ and $Y_{\pperp}(z(t))$  become perturbatively small.
We prove the following result.
\begin{lemma}\label{lem:sys}
Let $\fs_0, \fr$ as in \eqref{s0r} and $\s \geq \fs_0$. If $Z(t)=\vect{z(t)}{\bar z(t)} \in B_{\fs_0, \R}(I;\fr) \cap C^0(I; \dot H^\s_\R(\T; \C^2))$ 
 the variables
 $\big(z^\top(t), z^\perp(t)\big)$  defined in \eqref{ztp} fulfill  the system
\begin{align}
\label{eq.ztop}
\pa_t z^\top =  &- \im \Omega(D)  z^\top  + 
{X}^{(\intt)}(z) +  
 X_{\pperp} (z) + X_{\geq 4}(z) \\
\label{eq.zperp}
 \pa_t z^\perp  =  & -\im {\Omega(D)}   z^\perp  
+ \Opbw{  \im \,  \mathsf{m}(z; x, \xi) }z^\perp +Y^{(\intt)}_{\perp}(z) + Y_{\pperp}(z) 
+  Y_{\geq 4}(z)
\end{align}
where\\
$\bullet$ $X^{(\intt)}(z)$ is the integrable vector field
\begin{equation}
X^{(\intt)} (z) := X^{(\intt)} (z^\top) = -\im 
\left(2\fa |z_\tm|^2 + \fb |z_\tn|^2 \right) z_\tm e^{\im \tm x} 
-\im 
\left( 2\fc |z_\tn|^2  + \fb |z_\tm|^2 \right) z_\tn e^{\im \tn x} \ , 
\quad 
\fa, \fb , \fc \mbox{ in } \eqref{abc}     \ ;
\label{Y.int}
\end{equation}
$\bullet$ $\sm(z;x, \xi):= \sm^{(\res)}_{\geq 2}(Z; x, \xi)$  is the real-valued symbol in  $\Sigma\Gamma^1_{ 2}[\fr, 3]$ defined in \eqref{D.final}\,, i.e. 
\begin{equation}\label{tm}
    \sm(z; x, \xi):=  \big( \langle \,\sfV\, \rangle (Z;x) + {\sfV}_{\geq 3}(Z; x)\big) \xi + {\td}_{\geq 2}(Z; x) \omega(\xi) +  {\tf}_{\geq 2}(Z;x)\,\sign{(\xi)}+  \tg_{\geq 2}^{(-\frac12)}(Z; x, \xi)
\end{equation}
with $ \langle \, \sfV\, \rangle(Z;x)$ defined in \eqref{VresZ};\\
$\bullet$ $Y^{(\intt)}_{\perp}(z)$ 
is a cubic integrable  vector field fulfilling the energy identity
\begin{equation}\label{Y.perp}
\begin{gathered}
\Re \langle Y^{(\intt)}_{\perp}(z),\, z^\perp \rangle = 0\,
\end{gathered}
\end{equation}
and, for any $\s \geq \fs_0$, the smoothing estimate
\be
\label{Y.perp3}
\norm{Y^{(\intt)}_\perp(z)}_{\s +1} \lesssim \norm{z^\top}_{\fs_0}^2 \norm{z^\perp}_{\s} \ ; 
\ee
$\bullet$ $X_{\pperp}(z)$ and $Y_{\pperp}(z)$ are cubic vector fields, at least quadratic in the variable $z^\perp$, fulfilling for any $\s \geq \fs_0$ the smoothing estimates
\be\label{Y.perp2}
\norm{X_{\pperp}(z)}_{\s} \lesssim \left( \norm{z^\top}_{\fs_0} + \norm{z^\perp}_{\fs_0}\right) \norm{z^\perp}_{\fs_0}^2\,,
\qquad \norm{Y_{\pperp}(z)}_{\s + 1} \lesssim \left( \norm{z^\top}_{\fs_0} + \norm{z^\perp}_{\fs_0}\right) \norm{z^\perp}_{\fs_0} \norm{z^\perp}_{\s}\ ;
\ee
$\bullet$
$X_{\geq 4}(z)$ and $Y_{\geq 4}(z)$ are non-homogeneous vector fields fulfilling the estimate:  for any $\s\geq \fs_0$ there are  $C>0$, $\tr:=\tr(\s) \in (0,\fr)$ and for any  $z\in B_{\fs_0}(\tr)\cap \dot {H}^{\s}(\T;\C)$, 
\be\label{Y5}
\norm{X_{\geq 4}(z)}_{\s} + \norm{Y_{\geq 4}(z)}_{\s} \leq C  \norm{z}^3_{\fs_0} \norm{z}_{\s}  \,. 
\ee
\end{lemma}
To prove Lemma \ref{lem:sys}, we need to compute the projections of the cubic part of the vector field; therefore we shall use the following auxiliary lemma.
\begin{lemma}\label{lem:proiettori} Let $\Pi^\top$, $\Pi^\perp$ as in \eqref{Pi}.
The following holds true:
\begin{itemize}
    \item[$(i)$] Let  $\sm$ be the symbol defined in \eqref{tm}. If in the cutoffs appearing in \eqref{supp.chi} and in Definition \ref{quantizationtotale} one takes $\delta_0 < (2 \tn + 3)^{-1}$, then
    \begin{align}
\label{proj.op}
& \Pi^\top \Opbw{\im\, \sm } \Pi^\top  = \Pi^\top \Opbw{\im \, \sm_{\geq 3} } \Pi^\top \,, \\
\label{proj.op2}
&  \Pi^\top\Opbw{\im \,  \sm } \Pi^\perp = \Pi^\perp \Opbw{\im \, \sm } \Pi^\top = 0 \,, \\
\label{proj.op3}
& \Pi^\perp\Opbw{\im \,  \sm } \Pi^\perp  = \Opbw{\im \, \sm } \Pi^\perp  \,,
\end{align}
where $\sm_{\geq 3}$ in $\Gamma_{ \geq 3}^{1}[\fr]$ is given by
\be\label{smgeq3}
\sm_{\geq 3}(z; x, \xi):= 
 {\sfV}_{\geq 3}(z; x)\, \xi + \cP_{\geq 3} \left[{\td}_{\geq 2}(z; x) \omega(\xi) +  {\tf}_{\geq 2}(z;x)\,\sign{(\xi)}+  \tg_{\geq 2}^{(-\frac12)}(z; x, \xi)\right]\,.
\ee
\item[$(ii)$] Let $X^{(\Set)}(Z)$ be the vector field  in \eqref{Xlambda} in strong-$\Set$ normal form (Definition \ref{def:wr}),  then 
\begin{equation}\label{pi.bot.e.via}
\begin{aligned}
     \Pi^\top X^{(\Set)}(Z)^+ &=  X^{(\intt)}(z) + X_{\pperp}(z) \,, \\
     \Pi^\perp X^{(\Set)}(Z)^+ &= Y^{(\intt)}_{\perp}(z) + Y_{\pperp}(z)\,,
     \end{aligned}
\end{equation}
with $X^{(\intt)}(z)$ defined in \eqref{Y.int}, $Y^{(\intt)}(z)$ fulfilling \eqref{Y.perp}, \eqref{Y.perp3}, $X_{\pperp}(z)$ and $Y_{\pperp}(z)$ fulfilling \eqref{Y.perp2} for any $\s \geq \fs_0$, and where we recall the notation in \eqref{rtr}.
\end{itemize}
\end{lemma}
\begin{proof}
    $(i)$ We write $\sm = \sm_2 + \sm_{\geq 3}$, with $\sm_2 := \cP_{ 2}\left[ \sm\right] \in \wt{\Gamma}^{1}_{2}$ and $\sm_{\geq 3}:= \cP_{\geq 3} \left[\sm\right] \in \Gamma_{\geq 3}^1[\fr]$. 
    By  definition   \eqref{BW}, 
 \be\label{proj.op00}
\begin{aligned}
\Pi^\top  \Opbw{\im \, \sm_2 (Z; x, \xi)} \Pi^\top  z = & 
\!\!\!\!\!\! \sum_{\substack{ \sigma_1 j_1 + \sigma_2 j_2 +j=k  \\  j, k \in \Set} } \!\!\!  \chi_2 \left( j_1, j_2,\frac{j+k}{2}\right)  \im \, (\hat \sm_2)_{j_1, j_2}^{\sigma_1,\sigma_2}\left(\frac{j+k}{2}\right) z_{j_1}^{\sigma_1} \, {z_{j_2}^{\sigma_2}} \, z_j
\, {e^{\im k x}}  \,.
\end{aligned}
\ee
We now show that the cut-off  vanishes for any possible choice of $(j_1, j_2) \in \Z_*^2$ and  $j,k \in \Set$. 
Indeed, recalling that   $\chi_2(\xi', \xi) \equiv 0$ when $|\xi'| \equiv \max\{|\xi'_1|, |\xi'_2|\} \geq \delta_0 \langle \xi \rangle$, and using 
$\max\{|j_1|, |j_2|\} \geq 1$ (as $j_1, j_2 \neq 0$), $ j = k-\sigma_1 j_1 -\sigma_2 j_2$, $k \in \Set=\{\tm, \tn\}$ and provided $\delta_0 \leq 1/2(1+\tn)$,  one has 
\be\label{proj.op01}
\delta_0 {\big\langle \frac{j + k}{2} \big\rangle} = \delta_0 {\big\langle \frac{2 k-\sigma_1 j_1 -\sigma_2 j_2}{2} \big\rangle} 
\leq \delta_0 \frac{2\tn + 2\max{(|j_1|, |j_2|)}}{2}
 \leq \frac12 \max(|j_1|, |j_2|)   \,. 
\ee
Thus 
$\chi_2 \left( j_1, j_2,\frac{j+k}{2}\right) \equiv 0$  for any $j_1, j_2, j, k$ with $j,k \in \Set$ and 
 $\Pi^\top  \Opbw{\im \, \sm_2 } \Pi^\top  =0$, proving  \eqref{proj.op}.

\noindent \underline{Proof of \eqref{proj.op2}}. 
Again we write   explicitly the action of $\Pi^\top\Opbw{\im \, \sm_2 + \im\, \sm_{\geq 3} } \Pi^\perp $, using the quantization \eqref{BW} for the $2$-homogeneous symbol $\sm_2(z; \cdot)$ and \eqref{BWnon} for the non-homogeneous symbol  $\sm_{\geq 3}(z; \cdot)$, getting 
  \be\label{proj.op20}
\begin{aligned}
\Pi^\top  \Opbw{\im \, (\sm_2 + \sm_{\geq 3})} \Pi^\perp  z = & 
\!\!\!\!\!\! \sum_{\substack{ \sigma_1 j_1  + \sigma_2 j_2 +j=k  \\  \ j \in \Set^c,   k \in \Set} } \!\!\!  \chi_2 \left(j_1, j_2,\frac{j+k}{2}\right)  \im \, (\hat \sm_2)_{j_1, j_2}^{\sigma_1, \sigma_2}\left(\frac{j+k}{2}\right) z_{j_1}^{\sigma_1} \,  z_{j_2}^{\sigma_2} \,  z_j \, 
\, {e^{\im k x}} 
 \\
&  + \!\!\!\!\!\! \sum_{j \in \Set^c, \,   k \in \Set } \!\!\!  \chi \left(k-j,\frac{j+k}{2} \right)
\im \, \hat{\sm}_{\geq 3}\left(z; k-j, \frac{k+j}{2}\right) \, z_j \, {e^{\im k x}} \,.
\end{aligned}
\ee
Arguing as in \eqref{proj.op01},  the first line of \eqref{proj.op20} vanishes. 
To deal with the second line, recall that also  $\chi(\xi', \xi) \equiv 0$ when $|\xi'| \geq \delta_0 \langle \xi \rangle$. Hence for  $k \in \Set$, $j \in \Set^c$ (so $|j-k|\geq 1$) and $\delta_0 \leq 1/(2\tn + 3)$, recalling $|k| \leq \tn$, we get
\be\label{proj.op21}
\delta_0 {\big\langle \frac{j + k }{2} \big\rangle} 
\leq \delta_0 \frac{2+\tn + |j|}{2}  \leq \delta_0 \frac{2+2\tn + |j - k| }{2} \leq \frac{ |j - k| }{2}  \,.
\ee
This shows that  
$\chi \left( k-j, \frac{j+k}{2}\right) \equiv 0$ for any choice of $j \in \Set^c$, $k \in \Set$ and thus   the second line of \eqref{proj.op20} vanishes as well, proving the first of \eqref{proj.op2}.
The second identity is analogous exchanging the roles of $j$ and $k$.

\noindent \underline{Proof of \eqref{proj.op3}.} It follows writing $\Pi^\perp = \uno - \Pi^\top$ and using the first of \eqref{proj.op2}. 
\\[1mm]
    $(ii)$ For any cubic vector field $X(Z)$ one has the following identities, by the definition of the sets $\fP^{(n)}_\Set$ (see \eqref{momj}): 
    \begin{align}\label{proj.id1}
& \Pi_{\fP^{(0)}_\Set}X(z)^+  = \Pi^\top X(z^\top, z^\top, z^\top)^+  \,, \\
\label{proj.id2}
& \Pi_{\fP^{(1)}_\Set}X(z)^+  = {3}\Pi^\top X(z^\top, z^\top, z^\perp)^+ + \Pi^\perp X( z^\top, z^\top, z^\top)^+    \,, \\
\label{proj.id3}
&  \Pi_{\fP^{(2)}_\Set}X(z)^+  ={3} \Pi^\top X(z^\top, z^\perp, z^\perp)^+ + {3}\Pi^\perp X( z^\top, z^\top, z^\perp)^+ \,.
    \end{align}
    We now show the decomposition \eqref{pi.bot.e.via} for $ \Pi^\top X^{(\Set)}(z)^+$.  Writing $z = z^\top + z^\perp$ and exploiting the symmetry of $X^{(\Set)}$ in the internal variables,  we   decompose 
    \begin{align*}
         \Pi^\top X^{(\Set)}(z)^+ =  \underbrace{\Pi^\top X^{(\Set)}(z^\top, z^\top, z^\top)^+}_{ = \Pi_{\fP^{(0)}_\Set}X^{(\Set)}(z)^+ \mbox{ by } \eqref{proj.id1}} +
   \underbrace{ {3}\Pi^\top X^{(\Set)}(z^\top, z^\top, z^\perp)^+}_{ = \Pi^\top \Pi_{\fP^{(1)}_\Set}X^{(\Set)}(z)^+ \mbox{ by } \eqref{proj.id2} } + {3}\Pi^\top X^{(\Set)}(z^\top, z^\perp, z^\perp)^+
   + \Pi^\top X^{(\Set)}(z^\perp, z^\perp, z^\perp)^+
    \end{align*}
   Since $X^{(\Set)}(z)$ is in strong $\Set$-normal form,  by \eqref{Y.str}  we have that 
   \[
  X^{(\intt)}(z) :=\Pi_{\fP^{(0)}_\Set} X^{(\Set)}(z)^+ 
   \]
   has the claimed form \eqref{Y.int} and
that $  \Pi_{\fP^{(1)}_\Set} X^{(\Set)}(z)^+ =0$. 
Therefore  we have the first identity of \eqref{pi.bot.e.via} with 
   $$
   X_{\pperp}(z) := 3\Pi^\top X^{(\Set)}(z^\top, z^\perp, z^\perp)^+ + \Pi^\top X^{(\Set)}(z^\perp, z^\perp, z^\perp)^+ \,. 
   $$
   We prove that $X_{\pperp}$ fulfills the first estimate in \eqref{Y.perp2}.
   This is obtained recalling that by \eqref{Xlambda}, $X^{(\Set)}(Z) = \bR_2^{(\Set)}(Z)Z$ with $\bR_2^{(\Set)}(Z)$ a matrix of smoothing operators in $\wt{\cR}_{2}^{-\und{\vr}}$, using \eqref{smoothing} with $m\rightsquigarrow - \und{\vr} < -1$, and 
   recalling that $\Pi^\top$, being a projector on a finite dimensional subspace (see \eqref{Pi}) fulfills 
    \be \label{ptop.smooth}
    \|\Pi^\top z\|_{\s}\leq C_{\s, \fs_0} \|\Pi^\top z \|_{\fs_0}\,,
    \ee
    for some constant $C_{\s, \fs_0}>0$.
\\
   Consider next  $ \Pi^\perp X^{(\Set)}(z)^+$. As before we decompose it as 
    \begin{align*}
         \Pi^\perp X^{(\Set)}(z)^+ =  \underbrace{\Pi^\perp X^{(\Set)}(z^\top, z^\top, z^\top)^+}_{ = \Pi^\perp \Pi_{\fP^{(1)}_\Set}X^{(\Set)}(z)^+ \mbox{ by } \eqref{proj.id2}} +
 {3}\Pi^\perp X^{(\Set)}(z^\top, z^\top, z^\perp)^+
 + {3} \Pi^\perp X^{(\Set)}(z^\top, z^\perp, z^\perp)^+
   + \Pi^\perp X^{(\Set)}(z^\perp, z^\perp, z^\perp)^+\,.
    \end{align*}
    Again, the term  $\Pi_{\fP^{(1)}_\Lambda} X^{(\Set)}(z)^+ = 0$ by the second of \eqref{Y.str} and we define
    $$
    \begin{aligned}
    Y^{(\intt)}_{\perp}(z) &:= {3}\Pi^\perp X^{(\Set)}(z^\top, z^\top, z^\perp)^+ \stackrel{\eqref{proj.id3}}{=} \Pi^\perp \Pi_{\fP^{(2)}_\Set}X^{(\Set)}(z)^+  \,, \\
    Y_{\pperp}(z) &:= {3} \Pi^\perp X^{(\Set)}(z^\top, z^\perp, z^\perp)^+
   + \Pi^\perp X^{(\Set)}(z^\perp, z^\perp, z^\perp)^+\,.  
    \end{aligned}
    $$ 
    This gives the second identity of \eqref{pi.bot.e.via}.
    Then $Y^{(\intt)}_\perp$ fulfills the algebraic cancellation \eqref{Y.perp} by the last of \eqref{Y.str} and the smoothing estimate \eqref{Y.perp3}  by  \eqref{smoothing} with $m\rightsquigarrow - 1$ and \eqref{ptop.smooth}.
    Again, by \eqref{smoothing} with $m\rightsquigarrow - 1$ and \eqref{ptop.smooth}, the vector field $Y_{\pperp}$ satisfies \eqref{Y.perp2}.
\end{proof}

\begin{proof}[Proof of Lemma \ref{lem:sys}]
We compute the equations for $z^\top(t)$ and $z^\perp(t)$ by projecting the 
vector field in \eqref{Z.eq} via the projectors $\Pi^\top$ and $\Pi^\perp$ defined in \eqref{Pi}. 
\\[1mm]
\underline{Equation for $z^\top$.}  
By \eqref{Z.eq}, and recalling the notation in \eqref{tm},
\begin{align*}
   & \pa_t z^\top =   \Pi^\top \left[ (-\im \vOmega(D))^+ z + \Opbw{\im \sm}z + X^{(\Set)}(Z)^+ + [\bB_{\geq N}(Z)Z]^+ + [{\bR}_{\geq 3}(Z)Z]^+  \right] \\
   & \  \stackrel{ \eqref{proj.op}, \eqref{proj.op2}, \eqref{pi.bot.e.via}}{=} 
 -\im \Omega(D)z^\top +   X^{(\intt)}(z) + X_{\pperp}(z) +
   \underbrace{\Pi^\top \Opbw{\im \sm_{\geq 3}} \Pi^\top z + \Pi^\top \big[ [\bB_{\geq N}(Z)Z]^+ + [{\bR}_{\geq 3}(Z)Z]^+\big] }_{=: X_{\geq 4}(z)}
\end{align*}
 where $\sm_{\geq 3}$ in \eqref{smgeq3}.
By Lemma \ref{lem:proiettori},  $X^{(\intt)}$ and $X_{\pperp}$ fulfill the claimed properties. Finally, 
$X_{\geq 4}(z)$ fulfills \eqref{Y5} by Theorem \ref{thm:nf}, and combining the estimates \eqref{stimapar2}, \eqref{ptop.smooth}, \eqref{bound:specloc} (recall we fixed $N=3$), \eqref{piove}.
\\[1mm]
\underline{Equation for $z^\perp$.} Proceeding similarly for $z^\perp$ we obtain 
\begin{align*}
& \pa_t z^\perp  =  \Pi^\perp \left[(-\im \vOmega(D))^+ z + \Opbw{\im \sm}z + X^{(\Set)}(z)^+ + [\bB_{\geq N}(Z)Z]^+  + [{\bR}_{\geq 3}(Z)Z]^+\right] \\
   & \stackrel{ \eqref{proj.op2}, \eqref{proj.op3}, \eqref{pi.bot.e.via}}{=} 
  -\im \Omega(D)z^\perp +  \Opbw{\im \sm} z^\perp 
    +
   Y^{(\intt)}_{\perp}(z) + Y_{\pperp}(z)+
   \underbrace{ \Pi^\perp  \big[ [\bB_{\geq N}(Z)Z]^+  + [{\bR}_{\geq 3}(Z)Z]^+ \big]}_{=: Y_{\geq 4}(z)}\,.
\end{align*}

Again, by Lemma \ref{lem:proiettori}, $Y^{(\intt)}_{\perp}$ and $Y_{\pperp}$ fulfill the claimed properties, and 
$Y_{\geq 4}(z)$ fulfills \eqref{Y5}  by Theorem \ref{thm:nf}, the first  estimate in \eqref{piove} and estimate \eqref{bound:specloc}.
\end{proof}

\subsection{Time evolution of long-time controlled solutions}\label{sub:timeevo.lcs}

In this subsection we analyze system \eqref{eq.ztop}--\eqref{eq.zperp} when $z(t)$ is a long-time controlled solution. 
The first crucial property of long times controlled solutions, that we prove in the next lemma,  is that their  low norms $\norm{\ \cdot \ }_{\fs_0}$, $\norm{\ \cdot\ }_{L^2}$ remain small  over the whole time interval  $[0, T_\star]$.
To prove this,  it is crucial to exploit that the term in \eqref{eq.zperp}  linear in $z^\perp$, namely 
$ - \im \Omega(D) z^\perp + \Opbw{  \im \,  \mathsf{m}(z; x, \xi) }z^\perp +Y^{(\intt)}_{\perp}(z) $,
 vanishes in a $L^2$-energy estimate: indeed, the paradifferential part vanishes because its symbol is purely imaginary, 
 whereas  $Y^{(\intt)}_{\perp}(z)$ because of the energy identity \eqref{Y.perp}.
\begin{lemma}[\bf Bootstrap lemma]\label{lem:long.boot}
Let  $\fs_0, s, \theta$ as in \eqref{s0r}. Fix also 
 $T_0 >0$.   
There exist $\e_\star = \e_\star(\theta, T_0,\fr) >0$ and $C = C(\tm, \tn)>0$ such that for any $\e \in (0, \e_\star)$ the following holds true. 

Let $z(t)$ be a solution  of  \eqref{eq.ztop}--\eqref{eq.zperp} which is long-time controlled   with parameters $(s, \theta, T_\star, \e)$ (according to Definition \ref{A}) and with $T_\star$ fulfilling 
 \eqref{Tstar.T0}. 
Then 
  $z(t)$ fulfills  the  $L^2$-bound
\begin{equation}
\label{boot}
\norm{z^\top(t)}_{L^2} \leq 2 \e  \ , \qquad 
\norm{z^\perp(t)}_{L^2} \leq \e^{2- \frac32\theta}  \ ,  \quad \forall \,  0\leq t \leq  T_\star 
\end{equation}
and the low-norm bound 
\begin{equation}
\label{z.s0}
\norm{z(t)}_{\fs_0} \leq C \e \ , \qquad 
\norm{z^\perp(t)}_{\fs_0} \leq \e^{\frac32}  \ \ ,  \quad \forall \,  0 \leq t \leq  T_\star \  . 
\end{equation}
\end{lemma}

\begin{proof}
The proof is by a bootstrap argument. 
We assume the bounds 
\begin{equation}
\label{boot2}
 \norm{z^\top(t)}_{L^2} \leq 10 \e  , \qquad 
\norm{z^\perp(t)}_{L^2} \leq \e^{2- 2 \theta}  \ ,  \quad \forall \, 0 \leq t \leq T_\star
\end{equation}
and show that, provided $\e \in (0, \e_\star)$ with $\e_\star$  sufficiently small,   the better bounds in  \eqref{boot} hold.

First of all, we   bound $\norm{z^\perp(t)}_{\fs_0}$.
This is  done interpolating the bound on 
$\norm{z^\perp(t)}_{L^2}$ (that we have by the 
bootstrap assumption \eqref{boot2})
 and the   large a-priori bound on $\norm{z^\perp(t)}_{s}$  that we have by the long-time controlled assumption, see  \eqref{app.ass2}.
We obtain 
\begin{align}
\norm{z^\perp(t)}_{\fs_0} 
& \leq  \norm{z^\perp(t)}_{L^2}^{1-\frac{\fs_0}{s}} \, 
\norm{z^\perp(t)}_{s}^{\frac{\fs_0}{s}} \stackrel{\eqref{boot2}, \eqref{app.ass2} }{\leq}
 \e^{2 - \theta(2-\frac{\fs_0}{s}) - 2 \frac{\fs_0}{s}}  \leq \e^\frac32 \
 \label{z.perp.s0}
\end{align}
where in the last passage we have used that  $s, \theta$ are as in \eqref{s0r}. 
Recalling that, since $z^\top$ is supported on finitely many modes, $\|z^\top(t)\|_{\fs_0} \leq C_{\fs_0} \|z^\top(t)\|_{L^2}$ for some positive constant $C_{\fs_0}$ (which depends on $\fs_0$, $\tm$, and $\tn$ only), and combining estimates \eqref{boot2}, \eqref{z.perp.s0}, we get
\begin{equation}
\label{boot3}
\norm{z(t)}_{\fs_0} \leq (10 C_{\fs_0} + 1) \e \,, \qquad 
\forall\, 0\leq t  \leq T_\star \,. 
\end{equation}
Next we consider  $\norm{z^\top(t)}_{L^2}$ and prove the first improved estimate in \eqref{boot}.
Recall that the function $z^\top(t)$ fulfills equation \eqref{eq.ztop}; since $X^{(\intt)}(z)$ is integrable and using in particular \eqref{Y.int}, we get that
$\Re \big\langle X^{(\intt)}(z), z^\top \big \rangle = 0$. Therefore, for all times  $ 0 \leq t \leq T_\star$
\begin{align*}
\frac{\di}{\di t} \norm{z^\top(t)}_{L^2}^2  & = 2\underbrace{ \Re \la -\ii 
\Omega(D) z^\top+ X^{(\intt)}(z) , z^\top \ra }_{ = 0}
+ 2 \Re \la X_{\pperp}(z) + X_{\geq 4}(z), z^\top \ra \\
& \stackrel{\eqref{Y.perp2}, \eqref{Y5}}{ \lesssim}
\left( \norm{z^\top(t)}_{\fs_0} \norm{z^\perp(t)}_{\fs_0}^2  + \norm{z^\perp(t)}_{\fs_0}^3 +
\norm{z(t)}_{\fs_0}^4  \right) \, \norm{z^\top(t)}_{L^2}  \stackrel{\eqref{z.perp.s0}, \eqref{boot3},  \eqref{boot2}}{\lesssim} \e^{5}  \,. 
\end{align*}
Then, since $z(t)$ is long-time controlled, the norm  $\norm{z^\top(0)}_{L^2}$ of  its initial datum is bounded by   \eqref{app.ass1}; hence  for all times  $ 0 \leq t  \leq T_\star \leq \frac{T_0}{\e^2}\log\left(\frac{1}{\e} \right)$, 
\begin{equation}\label{boot.res}
 \norm{z^\top(t)}_{L^2}^2 \leq  \norm{z^\top(0)}_{L^2}^2 + |t| C \e^5 \leq \e^2 +  C T_0 \e^3 \log(\e^{-1}) \leq 4\e^2 \,,
\end{equation}
provided $ 0< \e \leq \e_\star$ and $\e_\star$ is sufficiently small. This proves the first estimate in \eqref{boot}.
\\
We now bound $\norm{z^\perp(t)}_{L^2}$.  
Since the paradifferential operator in  equation \eqref{eq.zperp} is   skew-adjoint and $Y^{(\intt)}_\perp(z)$ fulfills the energy identity \eqref{Y.perp}  we get, for all times $ 0 \leq t  \leq T_\star \leq  \frac{T_0}{\e^2}\log\left(\frac{1}{\e} \right)$, 
\begin{align*}
\frac{\di}{\di t} \norm{z^\perp(t)}_{L^2}^2 
& = 
2\underbrace{ \Re \la \Big(-\ii \Omega(D) +\Opbw{\im \sm(z;\cdot)}\Big) z^\perp + Y^{(\intt)}_{\perp}(z) , z^\perp \ra }_{ = 0}+
 2 \Re \la  Y_{\pperp}(z) + Y_{\geq 4}(z), z^\perp \ra \\
& \stackrel{\eqref{Y.perp2}, \eqref{Y5}}{ \lesssim}
%
 \left( \norm{z^\top(t)}_{\fs_0} \,\norm{z^\bot(t)}^2_{\fs_0}
+ \norm{z^\perp(t)}_{\fs_0}^3
+ \norm{z(t)}_{\fs_0}^4
\right)
\norm{z^\perp(t)}_{L^2} \\
& \hspace{-10pt}
 \stackrel{\eqref{boot3}, \eqref{z.perp.s0},\eqref{boot.res}, \eqref{boot2}}{\lesssim}
 (\e^4 + \e^{\frac{9}{2}})\e^{2-2\theta} \lesssim \e^{6-2\theta}\,.
\end{align*}
Again, being $z(t)$ long-time controlled, its initial datum $z^\perp(0)$ fulfills  \eqref{app.ass1}; hence, for all times 
$ 0 \leq t \leq T_\star  \leq \frac{T_0}{\e^2}\log\left(\frac{1}{\e} \right)$, we bound
\begin{equation}\label{boot.res2}
 \norm{z^\perp(t)}_{L^2}^2 \leq  \norm{z^\perp(0)}_{L^2}^2 + |t|  C 
 \e^{6-2\theta}
 \leq \e^{4} +  C T_0 \e^{4-2\theta} \log(\e^{-1})
\leq \e^{2(2- \frac32 \theta)} 
 \,,
\end{equation}
up to shrinking $\e_\star$. 
Estimates  \eqref{boot.res} and \eqref{boot.res2} prove  \eqref{boot}, thus closing the bootstrap argument. 
Then the second of \eqref{z.s0} follows by \eqref{z.perp.s0}, and the first of \eqref{z.s0} follows by \eqref{boot3}.
\end{proof}
We are now ready to prove Proposition \ref{prop:eff}, which is the main result of the section.
\begin{proof}[Proof of Proposition \ref{prop:eff}]
We shall use that, since $z(t)$ is long-time controlled  with parameters $(s, \theta, T_\star, \e)$ and with $T_\star$ fulfilling 
 \eqref{Tstar.T0}, by Lemma \ref{lem:long.boot} it satisfies the bounds \eqref{boot}, \eqref{z.s0}. 

 \medskip
\noindent{\sc Equations for $z_{\tm}(t)$ and $z_{\tn}(t)$.}
Writing equation  \eqref{eq.ztop} in components,  using the explicit expression of $X^{(\intt)}$ in \eqref{Y.int}, we get the coupled system
\be\label{z1.z-1}
\begin{cases}
 \pa_t z_\tm = - \im \Omega_{\tm}(\gamma) z_\tm - \im (2 \fa |z_\tm|^2 + \fb|z_\tn|^2) z_\tm  + \la X_{\pperp}(z) + X_{\geq 4}(z),  e^{\im \tm x} \ra \\
 \pa_t z_{\tn}  = - \im \Omega_{\tn}(\gamma) z_{\tn} - \im (2\fc |z_{\tn}|^2 + \fb |z_{\tm}|^2) z_{\tn} + \la  X_{\pperp}(z) + X_{\geq 4}(z), e^{\im \tn x} \ra
\end{cases}\,,
\ee
with $\fa, \fb , \fc$  in  \eqref{abc}.
Consider the equation for $z_\tm$. We write it as
\be\label{z1(t)}
\begin{gathered}
\pa_t z_\tm  = - \im ( \Omega_{\tm}(\gamma) +  2 \fa |z_\tm(0)|^2 + \fb|z_\tn(0)|^2)  z_\tm   +  \sfd_\tm(t) , \\
\sfd_\tm(t) := -\im \left(2 \fa |z_\tm(t)|^2 + \fb|z_\tn(t)|^2 - 2 \fa |z_\tm(0)|^2 - \fb|z_\tn(0)|^2\right) z_\tm(t) 
+ \la X_{\pperp}(z) + X_{\geq 4}(z), e^{\im \tm x} \ra \,,
\end{gathered}
\ee
giving the first equation in \eqref{eff.sys}. We prove now that $\sfd_\tm(t)$ fulfills the bound claimed in \eqref{d.bound}.  
First,  using the first of \eqref{z1.z-1}, we get for all times $ 0 \leq t \leq T_\star $ 
 \begin{align*}
\frac{\di}{\di t} |z_\tm(t)|^2 & = 2 \Re  \left( \la X_{\pperp}(z) + X_{\geq 4}(z), e^{\im \tm x} \ra \,  \bar z_\tm \right) \\
&
 \stackrel{\eqref{Y.perp2}, \eqref{Y5}}{ \lesssim}
 \big( \norm{z^\top(t)}_{\fs_0}\norm{z^\perp(t)}_{\fs_0}^2 + \norm{z^\perp(t)}_{\fs_0}^3 +
\norm{z(t)}_{\fs_0}^4  \big) \, \norm{z^\top(t)}_{L^2}  \stackrel{\eqref{z.s0},  \eqref{boot}}{\lesssim}  \e^5 \,,
 \end{align*}
 which implies, for all $t \in [0, T_\star]$ with $T_\star$ as in \eqref{Tstar.T0},
\be\label{controlz1}
\abs{ |z_\tm(t)|^2 - |z_\tm(0)|^2   } \lesssim |t|  \e^5 \lesssim  T_0 \, \e^3 \log(\e^{-1})  \,.
 \ee
 Analogously, one gets 
 \be\label{controlz2}
 \abs{  |z_\tn(t)|^2 - |z_\tn(0)|^2  } \lesssim  T_0 \, \e^3 \log(\e^{-1})  \,.
 \ee
 Hence, using estimates \eqref{controlz1}, \eqref{controlz2}, \eqref{Y.perp2}, \eqref{Y5}, and \eqref{boot}, we get that $\sfd_\tm(t)$ in \eqref{z1(t)} is bounded for $0 \leq t \leq T_\star$ by 
 \be\label{d1.bound}
 \begin{aligned}
\abs{\sfd_\tm(t)} &\lesssim 
   \abs{ \big( |z_\tm(t)|^2 - |z_\tm(0)|^2 \big) z_\tm(t) } +   \abs{  \big( |z_\tn(t)|^2 - |z_\tn(0)|^2 \big) z_\tm(t) } +
 \abs{\la X_{\pperp}(z) + X_{\geq 4}(z), e^{\im \tm x} \ra } \\
 &
 {\lesssim}  T_0 \, \e^4 \log(\e^{-1}) + \e^4 \,, 
 \end{aligned}
 \ee
 proving  \eqref{d.bound}  provided $\e_\star$ is sufficiently small.
  An analogous argument proves that $z_{\tn}(t)$ fulfills the second of \eqref{eff.sys} with $\sfd_\tn(t)$ satisfying \eqref{d.bound}\,.
\\
Hence, by Duhamel formula, $z_{\tm}(t)$ and $z_{\tn}(t)$ decompose as 
\begin{equation}
\label{z1z-1.new}
\begin{aligned}
& z_{\tm}(t) = \tz_{\tm}(t) + r_{\tm}(t) \,, \quad \mbox{ where} \quad 
\tz_{\tm}(t):= e^{- \im t\big(\Omega_{\tm}(\gamma) +  2 \fa |z_\tm(0)|^2 + \fb|z_\tn(0)|^2\big)} z_{\tm}(0)\,,\\
& z_{\tn}(t) = \tz_{\tn}(t) + r_{\tn}(t) \,, \quad \mbox{ where} \quad 
\tz_{\tn}(t):= e^{- \im t\big(\Omega_{\tn}(\gamma) +  2 \fc |z_\tn(0)|^2 + \fb|z_\tm(0)|^2\big)} z_{\tn}(0) \,,
\end{aligned}
\end{equation}
with the integral remainders 
\begin{equation}\label{r1}
\begin{aligned}
r_{\tm} (t) &:= \int_0^t e^{- \im (t - \tau) (\Omega_{\tm}(\gamma) +  2 \fa |z_\tm(0)|^2 + \fb|z_\tn(0)|^2)
}  \, \sfd_{\tm} (\tau)  \, \di \tau  \  \,,
\\
r_{\tn} (t) &:= \int_0^t e^{- \im (t - \tau) (\Omega_{\tn}(\gamma) +  2 \fc |z_\tn(0)|^2 + \fb|z_\tm(0)|^2)
}  \, \sfd_{\tn} (\tau)  \, \di \tau \,,
\end{aligned}
\end{equation}
 fulfilling, using  \eqref{d1.bound} (rather than \eqref{d.bound}), \eqref{Tstar.T0} and eventually shrinking again $\e_\star$,  the bounds 
\begin{equation}
\label{r.bound}
|r_{\tm} (t)|, \ |r_\tn(t)| \  \leq \e^{{2}-\theta} , \qquad \forall 0 \leq t \leq T_\star \,.
\end{equation}

\medskip
\noindent{\sc Equation for $z^\perp(t)$.} We start from equation \eqref{eq.zperp} and we substitute the evolution  of $z_\tm(t)$, $z_\tn(t)$ in \eqref{z1z-1.new}. 
Consider first the symbol $\sm (z;x, \xi)$ in \eqref{tm}. We shall extract from its component  $\la \,  \sfV \, \ra(Z; x)$, defined in \eqref{VresZ}, the main contribution which is the one supported on $\Set$.
More precisely, we have
\begin{equation}
\label{V.un.pezzo}
 \la \,  \sfV \, \ra(Z(t); x)  = \sfV_\Set(t;x) + \sfW_1(t; x)\,,
  \end{equation}
  with
 \begin{gather*}
      \sfV_\Set(t;x) := \sfV_\tm^{(\rm int)} |z_\tm(t)|^2 + \sfV_\tn^{(\rm int)} |z_\tn(t)|^2
 + 2 \sfV_{\tm, \tn}^{(\rm res)}\,  \Re \left( z_{\tn}(t) \, \bar{z_{\tm}(t)} \, e^{\im (\tn-\tm) x} \right)\,, \\
 \sfW_1(t;x):=  \sum_{j \not \in \Set} \sfV_j^{(\rm int)} |z_j(t)|^2 
 +    \sum_{\substack{m < 0 <n , \ m, n \in \Set^c \\ \Omega_m(\gamma) = \Omega_n(\gamma)} }
2\sfV_{m, n}^{(\rm res)}
\Re \big(  z_{n}(t) \, \bar{z_{m}(t)} \, e^{\im (n-m) x} \big)\,.
 \end{gather*}
 Note that we used the fact that, by item (i) of Lemma \ref{lem:2wave}, $n \notin \Set$ and $\Omega_m(\gamma) = \Omega_n(\gamma)$ implies $m \notin \Set$.
We now approximate the evolution of the modes in $\Set$ using the nonlinear oscillatory dynamics in \eqref{z1z-1.new}, and obtain
\begin{equation}\label{V.altro.pezzo}
    \begin{aligned}
   \sfV_\Set (t;x) & =  \sfV_\tm^{(\rm int)} |z_\tm(0)|^2 + \sfV_\tn^{(\rm int)} |z_\tn(0)|^2
 + 2 \sfV_{\tm, \tn}^{(\rm res)}\,  \Re \left( \tz_{\tn}(t) \, \bar{\tz_{\tm}(t)} \, e^{\im (\tn-\tm) x} \right) + \sfW_2(t;x)\\
 &= \sfV_\tm^{(\rm int)} |z_\tm(0)|^2 + \sfV_\tn^{(\rm int)} |z_\tn(0)|^2
 + 2 \sfV_{\tm, \tn}^{(\rm res)}\,  \Re \left( z_{\tn}(0) \, \bar{z_{\tm}(0)} \, e^{\im (\tn-\tm) (x - \tJ t)} \right)  + \sfW_2(t;x)\\
 & = \fv(x- \tJ t) + \sfW_2(t; x)\,,
    \end{aligned}
\end{equation}
where $\tJ$ is the real number defined in \eqref{J1}, $\fv$ is the function in \eqref{V0}, and we have defined the remainder
  \begin{align*}
  \sfW_2(t; x):= &       
 \sfV_\tm^{(\rm int)} (2\Re(\tz_\tm(t) \bar{r_\tm(t)}) + |r_\tm(t)|^2 )
+  \sfV_\tn^{(\rm int)} (2\Re(\tz_\tn(t) \bar{r_\tn(t)}) + |r_\tn(t)|^2 )\\
& 
 +2 \sfV_{\tm, \tn}^{(\rm res)}\,
 \Re 
 \left( \big(\tz_{\tn}(t) \, \bar{r_{\tm}(t)} + r_\tn(t) \bar{\tz_{\tm}(t)} +
 r_\tn(t) \bar{r_{\tm}(t)} \big)
 \, e^{\im (\tn - \tm) x}  \right)\,.
  \end{align*}
Using  \eqref{z1z-1.new}, \eqref{boot}, \eqref{z.s0},   \eqref{r.bound} and \eqref{Vres.coeff}--\eqref{Vres.coeff2}, $\sfW_1$ and $\sfW_2$ fulfill the bounds
\be\label{est.tV12}
\begin{aligned}
    \norm{\sfW_1(t; \cdot)}_{W^{2,\infty}} &\lesssim \norm{z^\perp(t)}_{\fs_0}^2 \lesssim \e^{3}\ ,  \qquad 
    \norm{\sfW_2(t; \cdot)}_{W^{2,\infty}} &\lesssim \e^{3-\theta} 
  \,, \quad
\forall 0 \leq t \leq T_\star \,. 
\end{aligned}
\ee
Then we decompose  $\sm(z;\cdot)$ in \eqref{tm} as 
$$
\sm(z(t); x, \xi) = \fv(x-\tJ t)\xi + \sfW(t;x)\xi + \sfb(t;x,\xi)\,,
$$
with (recall \eqref{V.un.pezzo} and \eqref{V.altro.pezzo})
\begin{gather*}
\sfW(t; x) :=  \sfW_1(t; x) +  \sfW_2(t; x) + {\sfV}_{\geq 3}(Z(t); x)\,, \\
\sfb(t; x, \xi) :=  {\td}_{\geq 2}(Z(t); x) \omega(\xi) +  {\tf}_{\geq 2}(Z(t);x)\,\sign{(\xi)}+  \tg_{\geq 2}^{(-\frac12)}(Z(t); x, \xi)\,.
\end{gather*}
The real-valued function   $\sfW(t;x)$ 
fulfills the claimed bound \eqref{Vpert} thanks to the 
 estimates  \eqref{est.tV12} for $\sfW_1$ and $\sfW_2$, and, recalling  that ${\sfV}_{\geq 3}\in \cF^\R_{\geq 3}[\fr]$, the bound 
$$
\norm{{\sfV}_{\geq 3}(Z(t); \cdot)}_{W^{2, \infty}} \stackrel{\eqref{nonhomosymbo}}{\lesssim}  \norm{z(t)}_{\fs_0}^3 \stackrel{\eqref{z.s0}}{\lesssim} \e^3 \ , \quad
\forall 0 \leq t \leq T_\star \ .
$$
The bound \eqref{bpert} for $\sfb(t;x, \xi)$ follows from \eqref{nonhomosymbo.homo}, \eqref{nonhomosymbo} and \eqref{z.s0}. It remains to analyze the term 
$Y^{(\intt)}_{\perp}(z) + Y_{\pperp}(z) 
+  Y_{\geq 4}(z)$
in \eqref{eq.zperp}.
We put
$$
\sfR(t, x):=Y^{(\intt)}_{\perp}(z(t)) + Y_{\pperp}(z(t)) \ , \quad \sfY(t;x):=  Y_{\geq 4}(z(t))
$$
which fulfill the estimates \eqref{R(t,U)}, \eqref{Y(t,U)} by \eqref{Y.perp3}, \eqref{Y.perp2}, \eqref{Y5} and using \eqref{z.s0} and \eqref{app.ass2}.
\end{proof}
 \section{Growth of Sobolev Norms for  Long-Time Controlled Solutions}\label{sec:mourre}
The goal of this section is to give sufficient conditions that guarantee that 
long-time controlled solutions defined up to a time $T_0 \e^{-2} \log(\e^{-1})$ undergo growth of Sobolev norms. 
As we shall prove, such conditions regard {\em only the initial data} $z(0)$, and we shall say that $z(0)$ is {\em strongly well-prepared} if it fulfills them, see \Cref{B}. 
Among such conditions, the two crucial ones are (B2) -- the no-sign condition of the transport term -- and (B3) -- the sign condition of the Mourre operator that we are going to  introduce.
\\
Given $\gamma <0$, $\gamma^2\in \Q$ and  $\Set = \{\tm, \tn\}$ the $\gamma$-good  set in  \eqref{fixLambda}, take  $s, \theta$ as in \eqref{s0r}, and for any $\e>0$ define the parameter
\be\label{def:R}
\tR := \e^{-2(3+\theta)} \,,
\ee
and  the symmetric {\em Mourre  operator} 
\begin{equation}
\label{fA}
\begin{aligned}
\fA := \fA_{s,\tR}&:= \Opbw{\fp (x, \xi)} \,, \qquad \fp(x, \xi):=\tp(x) \, |\xi| ^{2s} \, \varphi^2_{\tR}(\xi)  
 \in \Gamma^{2s}_{W^{2, \infty}}\,, \\ 
\tp(x)&:= - \sfV_{\tm,\tn}^{(\rm res)} \Im  \, \left( z_{\tn}(0) \,  \bar{z_{\tm}(0)} \, e^{\im (\tn-\tm)  x }\right)  \,, 
\quad \sfV_{\tm,\tn}^{(\rm res)}    < 0 \mbox{ in } \eqref{Vres.coeff}
\end{aligned}
\end{equation}
where
$\varphi_\tR(\xi)$ is  the
smooth step function 
 \be \label{etaR}
 \varphi_\tR(\xi):= \varphi\left(\frac{\xi}{\tR}\right),\quad \varphi(y):= 
 \begin{cases}
0& \mbox{ if } y\leq 1\\ 
 \dfrac{e^{-\frac{1}{y-1}}}{e^{-\frac{1}{y-1}} + e^{-\frac{1}{2-y}}}& \mbox{ if } y \in (1,2) \\
 1& \mbox{ if } y \geq 2
 \end{cases}
  \,.
 \ee

\medskip
Let us describe the further conditions  on the initial data of a long-time controlled solution that guarantee  Sobolev norm explosion. 
\begin{definition}[{\bf Strongly Well-prepared data}]\label{B} 
Let $\Set$ in \eqref{fixLambda}, 
 $s, \theta$ as in \eqref{s0r} and fix    $\nu_0, \e  >0$. 
We say that an initial datum $z(0) \in \dot H^s(\T; \C)$ is {\em strongly well-prepared} with parameters $(s, \theta, \nu_0, \e)$ if the following conditions hold true:
\begin{itemize} 
\item[{\rm (B1)}]  {\sc (Weakly well-prepared):}  $z(0)$ is $\e$-weakly well-prepared (cf.  \eqref{app.ass1}); in particular one has 
\be\label{zmzn0}
\abs{z_\tm(0)} < \e \,, \quad \abs{z_\tn(0)} < \e  \,;
\ee
\item[{\rm (B2)}]  {\sc (No-sign condition of the transport term):} The Fourier modes $z_\tm(0)$ and $z_\tn(0)$ supported on  $\Set$ fulfill the lower bound

\begin{equation}\label{def:kappa}
\upkappa :=  2 \abs{\sfV_{\tm, \tn}^{(\rm res)} \, z_\tm(0) \,  z_{\tn}(0) } - \abs{  \tV_\tm |z_{\tm}(0)|^2 
  + \tV_\tn |z_{\tn}(0)|^2 } > \nu_0 \e^2 \ 
\end{equation}
with $\tV_\tm, \tV_\tn$ in \eqref{VmVn.def} and $\sfV_{\tm, \tn}^{(\rm res)}$ in \eqref{Vres.coeff};
\item[{\rm (B3)}] {\sc (Sign condition of the Mourre operator):} The Mourre operator $\fA_{s,\tR }$ in \eqref{fA}, with $\tR$ in \eqref{def:R}, is positive on the  projected initial data  $z^\perp(0)$ supported on $\Set^c$ with the quantitative estimate
\be\label{tRbound}
 \la \fA_{s,\tR } \,  z^\perp(0),   z^\perp(0) \ra > \e^{3-4\theta} \,;
\end{equation}
\item[{\rm (B4)}] {\sc (High-Sobolev norm smallness):} $z(0)$ has a small $\dot H^s$-Sobolev norm:
	\be\label{z0.s}
	\norm{z(0)}_s < \e^\theta  \,. 
	\ee
\end{itemize}
\end{definition}

\begin{remark}
Recall that the transport operator $\Opbw{(\tJ + \fv(x))\xi}$ in \eqref{J+v} depends on the values of the  Fourier modes $z_\tm(0)$, $z_\tn(0)$ of the initial data.
The no-sign condition in (B2)  guarantees that the transport
$\tJ + \fv(x)$ changes sign on $\T$, and actually vanishes at some points. This  is the property that makes possible the existence of a Mourre operator $\fA$ whose commutator with 
$\Opbw{(\tJ + \fv(x))\xi}$ is positive (up to perturbative terms), as we shall prove in \Cref{lem:pos.comm}. 
\end{remark}
The next proposition is the main result of the section. It proves that
a solution $z(t)$ which is long-time controlled for times $ {T_0}{\e^{-2}} \log\left({\e^{-1}} \right) $ with $T_0$ sufficiently large and whose initial datum is strongly well-prepared undergoes arbitrary large growth of Sobolev norms. 
 \begin{proposition}[{\bf Growth for long-time controlled solutions with  strongly-well-prepared data}]\label{prop:instab} 
 Let $s, \theta, \fr$ as in \eqref{s0r}. Fix also   $T_0, \nu_0 >0$.
 There exists $\e_1 = \e_1(s, \theta, T_0, \nu_0, \fr)>0$ such that, for any  $\e \in (0, \e_1)$, the following holds true. 
 Let $z(t)\in \dot H^s(\T; \C)$ be  a solution of system  \eqref{eq.ztop}--\eqref{eq.zperp}  such that 
 \begin{itemize}
 \item[(i)] its initial datum $z(0)\in \dot H^s(\T; \C)$ is strongly  well-prepared with parameters $(s, \theta, \nu_0, \e)$ (cf. Definition \ref{B});
\item[(ii)] it is long-time controlled with parameters $(s, \theta, T_\star, \e)$ (see Definition \ref{A}), 
with $T_\star $ fulfilling \eqref{Tstar.T0}.
 \end{itemize}
 Then there exists $\tipac_s \in (0,1)$ such that 
\be\label{sopra.exp2}
\norm{z(t)}_s^2 \geq \tipac_s^2  \e^{1-4\theta} \, e^{(\tn - \tm){\nu_0}\e^2 t} \,, \quad \forall \, t \in [0, T_\star]\,. 
\ee
In particular, if 
\be\label{T0}
T_\star = \frac{{1-2\theta}}{(\tn - \tm)\nu_0 \e^2} \log\left(\frac{1}{\e} \right)   \,,
\ee
 then $z(t)$ undergoes growth of Sobolev norms: 
 \begin{equation}\label{zeta.grows}
 \norm{z(T_\star)}_s \geq  
{\tipac_s} \e^{-\theta}\,.
\end{equation}
 \end{proposition}
\begin{remark}
In fact, conditions (B1)--(B3) in Definition \ref{B} are sufficient to establish \eqref{zeta.grows}, while condition (B4) guarantees that the initial datum is small in $H^s$.
\end{remark}
The rest of the section is devoted to the proof of \Cref{prop:instab}.
We start by stating the  following two auxiliary lemmata proved in \cite[Appendix A]{MM}.
 \begin{lemma}[{\bf High-frequency symbols}]\label{lem:cutoff}
 Let $\tN \in \N_0$, $m \in \R$ and $\tR\geq 1$. If  $a \in \Gamma^m_{W^{\tN, \infty}}$, then 
 \be \label{2110:1812}
 a_\tR(x, \xi):= a(x, \xi)\,  \varphi_\tR(\xi) , \quad \varphi_\tR \mbox{ in }  \eqref{etaR}
 \ee is a symbol in $\Gamma^{m+\nu}_{W^{\tN, \infty}}$ for any $\nu \geq 0$ with quantitative bound
 \begin{equation}\label{a.etaR.sem}
| a_\tR |_{m + \nu, W^{\tN, \infty}, n} \leq C_{n}\,  \tR^{-\nu} \, | a |_{m , W^{\tN, \infty}, n} \quad \text{for any } n\in \N_0\,. 
\end{equation}
In addition, if $\tN \geq 2$ and  $b \in \Gamma^{m'}_{W^{2, \infty}}$, $m' \in \R$, one has for all $u \in H^s(\T; \C)$ the commutator estimate
\begin{equation}\label{est.aR.comp}
\norm{ [ \Opbw{ a_\tR}, \Opbw{b} ] u }_{s - m - m' - \nu +1 } \leq C_{{s}} \,
\tR^{- \nu} \,  | a |_{m , W^{2, \infty}, 7}  \, | b |_{m' , W^{2, \infty}, 7} \, \norm{u}_s  \,.
\end{equation}
 \end{lemma}
A consequence of this lemma is that 
the symbol 
 $\fp(x, \xi)$ in \eqref{fA}, having the high-frequency form \eqref{2110:1812},  fulfills the following estimates (recall  \eqref{seminorm}): for any $n \in \N_0$, there is $C_{{s, n}} = C_{{s, n}}(\tm,\tn) >0$ such that 
\begin{equation}\label{a.sem}
| \fp |_{2s , W^{2,\infty},  n } \leq  C_{s,n} \,   |z_\tm(0)|\,|z_{\tn}(0)|  \,, \qquad 
| \fp |_{ 2s +1 , W^{2,\infty},  n } \leq C_{s,n} \frac{  |z_\tm(0)|\,|z_{\tn}(0)| }{\tR}  \,,
\end{equation} 
 as it follows from its definition and from  Lemma \ref{lem:cutoff} with $a \leadsto \tp(x)|\xi|^{2s}\varphi_{\tR}(\xi)$, $m \leadsto 2s$, $N \leadsto 2$ and $\nu \leadsto 1$.
  \begin{lemma}[\bf Strong Garding's inequality]
\label{garding}
Let $\tR \geq 1$, $a(x) \in W^{3,\infty}$ and $a(x) \geq 0$. Let $\psi(\xi) \in \wt\Gamma^{m}_0$, $m > 0$,  a real-valued Fourier multiplier  with $\textup{supp } \psi \subseteq  [\tR, +\infty)$. Then there is $C >0$ such that for any $u \in H^m(\T;\C)$
\begin{equation}
\label{gard1}
\la \Opbw{a(x) \psi^2(\xi)} u , u \ra \geq - C \frac{\norm{a}_{W^{3,\infty}}}{\tR^2} \norm{u}_{m}^2  \,.
\end{equation}
\end{lemma}

The next crucial  lemma proves that, if the initial data is strongly well-prepared, the commutator between $\fA$ and $\Opbw{(\tJ + \fv(x))\xi}$ is positive up to a perturbative term. More precisely, one has the following result.
\smallskip

{\bf Notation:} Given two operators $\fA, \mathsf{B}$,  we will write $\fA \geq \mathsf{B}$ with the meaning $\la \fA u, u \ra \geq \la \mathsf{B} u, u \ra $ for any $u \in  H^s(\T;\C)$.
\begin{lemma}[\bf Positive commutator estimate]\label{lem:pos.comm}
Let $s, \theta$ as in \eqref{s0r} and  fix also   $\nu_0>0$. There is $ \e_0= \e_0(s,\theta,\nu_0)>0$ such that for any  $ \e \in(0,\e_0)$ and  $\tR$ as in \eqref{def:R}, the following holds true. 
Let $z(0)$ be strongly well-prepared with parameters $(s, \theta, \nu_0, \e)$ (see \Cref{B}), and consider 
the Mourre operator $\fA\equiv \fA_{s, \tR}$  in \eqref{fA} and the function $\tJ + \fv(x) $ in \eqref{J+v}. 
Then one has the  positive commutator estimate
    \begin{equation}\label{mourre}
\im \big[\fA,   \Opbw{ \big( \tJ + \fv(x)\big) \xi} \big] \geq (\tn-\tm) \, \upkappa  \, \fA  + \fR  \,, \quad 
\mbox{ where } \upkappa >\nu_0 \e^2  \mbox{ in } \eqref{def:kappa}
\end{equation}
 and  the operator $\fR\colon H^s(\T; \C) \to H^{-s}(\T; \C)$ fulfills the following estimate: there exists $C_s >0$ such that
\begin{equation}\label{mourre.R}
\norm{\fR  u}_{-s } \leq C_{s} \frac{\e^4}{\tR} \norm{u}_s  \,, 
\quad \forall \, u \in H^{s}(\T; \C)\,.
\end{equation}
\end{lemma}
\begin{proof}
 The symbol   $\big(\tJ +\fv(x)\big) \xi $ belongs to  $\Gamma^1_{W^{2, \infty}}$ with seminorm (recall \eqref{J+v})
\begin{equation}\label{fv.sem}
| \big( \tJ + \fv\big)  \xi |_{1, W^{2, \infty}, 7} \lesssim  \left(|z_\tm(0)|^2 +  |z_{\tn}(0)|^2 \right)  \stackrel{\eqref{zmzn0}}{\lesssim} \e^2 \,.
\end{equation}
We compute the commutator in \eqref{mourre}. 
Despite the fact  that $\fp(x, \xi)\in \Gamma^{2s}_{W^{2,\infty}}$, it is convenient to regard it as a high-frequency symbol in $\Gamma^{2s+1}_{W^{2, \infty}}$; in this way we lose the fact that the remainders operator in symbolic calculus are smoothing, but we gain  negative powers of $\tR$ in the  remainders estimates, thanks to the seminorms bounds in  \eqref{a.sem}.
Precisely, we apply   the symbolic calculus of \Cref{teoremadicomposizione} with 
 $a = \fp(x,\xi) \in \Gamma^{2s+1}_{W^{2,\infty}}$, $b = (\tJ + \fv(x))\xi\in \Gamma^{1}_{W^{2,\infty}}$ and  $\varrho  =2 $ 
getting (use also \Cref{rem_symbols})
\begin{equation}\label{mourre1}
\im  [\fA,   \Opbw{\big( \tJ +  \fv(x) \big) \xi} ] = \Opbw{ \{ \fp(x, \xi) \,   ,  \big( \tJ+ \fv(x)\big) \xi  \}} + \breve R \,,
\end{equation}
where
the remainder $\breve R:= \cQ(\fp,(\tJ + \fv)\xi )- \cQ((\tJ + \fv)\xi, \fp ) $ is bounded $  H^{s}(\T;\C) \to H^{ -s}(\T;\C)$ and, by  \eqref{comp020},  fulfill the quantitative estimate
\begin{equation}
\| \breve R u\|_{-s} \lesssim  | \fp |_{2s +1, W^{2,\infty}, 7}  \, 
| \big( \tJ + \fv \big) \xi |_{1, W^{2,\infty}, 7} \, 
\norm{u}_{s}    \stackrel{\eqref{a.sem}, \eqref{fv.sem}}{\lesssim_s} 
 \frac{  |z_\tm(0)|\,|z_{\tn}(0)| }{\tR} \e^2
\stackrel{\eqref{zmzn0}}{\lesssim_s} \frac{\e^4}{\tR} \norm{u}_s \,. 
\label{breveR}
\end{equation}
We now consider the  Poisson bracket  $\{ \fp(x, \xi) \,   ,  \big( \tJ+ \fv(x)\big) \xi  \}$.
We claim that 
\be\label{comm.av0}
\{ \fp(x,\xi) , \, \big(\tJ + \fv(x) \big)\xi \} = 
(\tn-\tm)\upkappa \, \fp(x, \xi)  + a(x,\xi)   \,, 
\quad \upkappa \mbox{  in } \eqref{def:kappa}
\ee
and 
$a(x, \xi)$ a  smooth, non-negative symbol having  the  structure  
\be\label{sym.a}
 a(x,\xi) =  a_1(x)  \psi_1(\xi)^2 +  a_2(x) \psi_2(\xi)^2 \,, 
\ee
with  $ a_j(x)$, $j=1,2$,   smooth, real-valued,  non-negative functions fulfilling
\be\label{est.aj}
\norm{a_j}_{W^{3,\infty}} \lesssim \left( |z_\tm(0)|^4 +  \, |z_{\tn}(0)|^4 \right) \stackrel{\eqref{zmzn0}}{\lesssim} \e^4 \,, 
\ee
and 
 $\psi_j(\xi)$, $j = 1,2$,  smooth, real-valued symbols in $\wt \Gamma^s_0$ with  support in $[\tR, +\infty)$.
\\
Assuming the claim \eqref{comm.av0}, we have that
\begin{equation}\label{mourre20}
\Opbw{ \{ \fp(x, \xi)\,   ,  \big(\tJ +  \fv(x) \big) \xi  \}}    = 
(\tn-\tm)\upkappa \fA  +  \Opbw{a(x,\xi)}    
\end{equation}
and 
we bound the operator $\Opbw{a}$ (with symbol $a$ in \eqref{sym.a}) from below using the  
 strong Garding inequality in \Cref{garding}, getting that there is $C>0$ such that for any $u \in H^s(\T;\C)$
\begin{equation}\label{mourre30}
\la\Opbw{a} u, u \ra   \geq - C   
\frac{\norm{a_1}_{W^{3,\infty}} + \norm{a_2}_{W^{3,\infty}}}{\tR^2} \norm{u}_s^2  \stackrel{\eqref{est.aj} }{\geq }
-C\frac{\e^4}{\tR^2} \la \la D\ra^{2s} u, u \ra \,. 
\end{equation}
By \eqref{mourre1}, \eqref{mourre20}, \eqref{mourre30}, we conclude that 
\begin{equation}\label{mourre2}
\im  [\fA,   \Opbw{ \big( \tJ + \fv(x) \big) \xi} ]   \geq   
(\tn - \tm) \upkappa \sfA + \fR \ , \quad 
\fR := \breve R  -C \frac{\e^4}{\tR^2}\la D\ra^{2s} \,,
\end{equation}
where 
the operator $\fR \colon H^{s} \to H^{-s }$ fulfills the  estimate \eqref{mourre.R}, also using \eqref{breveR}. 

\smallskip
\underline{Proof of \eqref{comm.av0}:} 
Using \eqref{poisson},  \eqref{fA} and denoting  $(\varphi')_\tR(\xi):= \varphi'(\xi/\tR)$, we compute
\begin{equation}
\label{pbcruz}
 \{ \fp(x, \xi)  ,    \big(\tJ + \fv(x)\big)\xi  \} 
 =
\big(   \tp \,  \fv_x - (\tJ + \fv) \, \tp_x   \big) |\xi|^{2s} \, \varphi_\tR^2 + (2s-1) \tp \fv_x   |\xi|^{2s} \varphi_\tR^2 
+  2 \tp  \fv_x \, |\xi|^{2s}\varphi_\tR \, \frac{\xi}{\tR}   \,  (\varphi')_\tR \,.
\end{equation}
Now, using the explicit definition of $\tp(x) $ in \eqref{fA}, of $\tJ + \fv(x)$ in \eqref{J+v} and that 
$$\tp_x(x)=- {\sfV_{\tm,\tn}^{(\rm res)}}(\tn-\tm) \Re  \, \left( z_{\tn}(0) \,  \bar{z_{\tm}(0)} \, e^{\im (\tn-\tm)  x }\right) \,, \quad 
\fv_x(x)=- 2(\tn-\tm) \sfV_{\tm,\tn}^{(\rm res)} \Im \left(z_\tn(0) \, \bar{z_{\tm}(0)} \, e^{\im (\tn-\tm) x} \right)  \,, $$
we get the lower bound
\begin{align}
\notag
\tp \fv_x - (\tJ + \fv )\tp_x  
& = 2 (\tn-\tm)(\sfV_{\tm,\tn}^{(\rm res)} )^2 |z_\tn(0)|^2 \, |z_\tm(0)|^2 - \tp_x \left( \tV_\tm |z_\tm(0)|^2 + \tV_\tn |z_\tn(0)|^2 \right) \\
\notag
&  \geq 
 (\tn-\tm)
 \abs{\sfV_{\tm,\tn}^{(\rm res)} \, 
  z_{\tn}(0) \,  z_{\tm}(0)}\,  
\Big( 2 \abs{\sfV_{\tm,\tn}^{(\rm res)} \, z_{\tn}(0) \, 
  z_{\tm}(0)} -  \abs{\tV_\tm |z_\tm(0)|^2 + \tV_\tn |z_\tn(0)|^2 } \Big) \\ \label{schiena}
 &    \geq (\tn - \tm) \tp \, \upkappa
\end{align}
where to pass from the first to the second line  we  used that
$
|\tp_x| \leq (\tn - \tm)\, |\sfV_{\tm,\tn}^{(\rm res)}| \, |z_\tm(0)| \, |z_{\tn}(0)|$, 
and to pass from the second to the third one we used $ \abs{\sfV_{\tm,\tn}^{(\rm res)} \, z_\tn(0) \, z_\tm(0)}  \geq \tp(x) $ and the definition of $\upkappa $ in \eqref{def:kappa}.
\\
Hence, adding and subtracting
$(\tn-\tm)\upkappa \fp $ 
in  \eqref{pbcruz}, we get the claimed formula \eqref{comm.av0} with 
\begin{equation*}
 a(x, \xi) := 
\underbrace{\left(\tp \fv_x - (\tJ + \fv) \tp_x   - (\tn-\tm)\upkappa \tp
+(2s-1) \tp \fv_x 
  \right)}_{=:a_1(x)} \underbrace{|\xi|^{2s} \varphi_\tR^2 }_{=:\psi_1(\xi)^2} +   \underbrace{2 \tp  \fv_x }_{=:a_2(x)}\, \underbrace{|\xi|^{2s}\varphi_\tR \, \frac{\xi}{\tR}   \,  (\varphi')_\tR}_{=:\psi_2(\xi)^2}   \,.
\end{equation*}
By \eqref{fA}, \eqref{J+v}, the function
\[
a_2:= 2\tp \fv_x = 2 (\tn-\tm) (\sfV_{\tm, \tn}^{(\rm res)})^2\,  \Im \left( z_{\tn}(0) \,  \bar{z_{\tm}(0)} \, e^{\im (\tn-\tm)  x } \right)^2 \geq 0 \,, 
\]
and, using also  \eqref{schiena}, also the   function
$$
a_1:=\tp \fv_x - (\tJ + \fv) \tp_x   - (\tn-\tm)\upkappa \tp 
+(2s-1) \tp \fv_x
\geq 0 \,.
$$
Therefore, both  functions  $a_j(x) \geq 0$, $j=1,2$; moreover, they   both are smooth and  fulfill estimate \eqref{est.aj}, which is easily checked  using  the definitions of $\tp(x)$ and $ \fv(x)$ respectively in \eqref{fA} and \eqref{V0}, of $\tJ$ in \eqref{J1} and $\upkappa$ in \eqref{def:kappa}.
To conclude, we show that 
   $\psi_1(\xi)= |\xi|^s \varphi_\tR(\xi)$ and $\psi_2(\xi) = |\xi|^s \sqrt{\varphi_\tR(\xi) \, \frac{\xi}{\tR} (\varphi')_\tR(\xi)}$ are smooth symbols in $\wt \Gamma^{s}_0$ supported in $[\tR, \infty)$.
We prove the claim only for $\psi_2$ since the one  for $\psi_1$ is trivial.
First of all, note that $\psi_2$ is well defined since, by  \eqref{etaR}, one has $\xi (\varphi')_\tR(\xi) \geq 0$.
Define
 $$
f(y):= \sqrt{\varphi(y) \, y \varphi'(y)} \,, \quad \textup{supp}(f) \subset [1,2] \,.
$$
 Then  $\psi_2(\xi)=|\xi|^s f(\xi/\tR)$ and is supported in $[\tR, 2\tR]$. So it remains to prove that 
$f(y)$ is a smooth function. It is easy to see that $\sqrt{y \varphi(y)}$ is smooth on its support while the function 
$$
\sqrt{\varphi'(y)} =
\begin{cases}
0 \,, & y\leq 1\\
\dfrac{\sqrt{2 y^2-6 y+5}}{e^{-\frac{1}{2-y}}+e^{-\frac{1}{y-1}} } \cdot \dfrac{e^{-\frac{1}{2(y-1)}}}{y-1} \cdot  
\dfrac{e^{-\frac{1}{2(2-y)}}}{2-y}   \,, & y \in (1,2) \\
0 \,,& y\geq 2
\end{cases}
$$
 is smooth by direct inspection.
\end{proof}
\begin{remark}
   Actually, only conditions (B1)--(B2) in Definition \ref{B} are needed to prove  \eqref{mourre}.
\end{remark}
\begin{proof}[Proof of \Cref{prop:instab}]
   Define the quadratic form 
\begin{equation}
\cA(t):= \langle \fA \,  \zeg(t),\, \zeg(t) \rangle, \quad \fA \equiv \fA_{s, \tR}  \mbox{ in } \eqref{fA} \,,\quad \zeg(t)  \mbox{ in } \eqref{zeta} \,.
\label{cA}
\end{equation}
We shall 
provide a positive lower bound on the 
  time derivative  $\frac{\di}{\di t} \cA(t)$, yielding a quantitative growth in time of the quadratic form $\cA(t)$.
Using  equation  \eqref{eq.zeta} for the evolution of $\zeg(t)$, we compute 
\begin{align}
\label{I}
\frac{\di}{\di t} \cA(t)  = & \la \im \big[ \fA, \Opbw{\big( \tJ + \fv(x) \big) \xi}\big]  \zeg,\, \zeg\ra
 \\
 \label{II}
 &  + \la \im  \big[ \fA,\Opbw{ \wt \sfW(t;x) \xi}\big] \zeg,\, \zeg \ra  
 \\
 \label{III}
 & + \la \im \big[ \fA,\Opbw{ - \Omega(\xi) + \wt \sfb(t; x, \xi)}\big] \zeg,\, \zeg \ra\\
 \label{IV}
 & + 2\Re\,  \la \fA \wt{\sfR}(t) ,\, \zeg \ra \\
 \label{V}
 & + 2\Re\,  \la \fA \wt{\sfY}(t) ,\, \zeg \ra  \,.
\end{align}
We now estimate each term.
Term \eqref{I} is estimated from below  using Lemma \ref{lem:pos.comm}, getting
 \begin{equation}\label{20250130:1908}
   \la \im [ \fA, \Opbw{(\tJ+\fv(x)) \xi}]  \zeg, \zeg \ra 
  {\geq }  (\tn - \tm) \, \upkappa  \la \fA \zeg, \zeg \ra  - C_s \frac{\e^{4}}{\tR} \norm{\zeg}_s^2\,.
 \end{equation}
Next  we estimate \eqref{II} from above. 
We consider $\fp(x,\xi)$ as a symbol in $\Gamma^{2s}_{W^{2,\infty}}$. By   estimates  \eqref{est.aR.comp} (with $\nu =0$, $m' = 1$, $m = 2s$),   we get 
\begin{equation}
\abs{\eqref{II}} \lesssim | \fp |_{2s, W^{2,\infty}, 7} \, |\wt \sfW(t,\cdot) \xi |_{1, W^{2,\infty}, 7} \, \norm{\zeg}_s^2
\stackrel{\eqref{a.sem}, \eqref{zmzn0},   \eqref{est.706}}{\lesssim } \ 
 \e^{5-\theta}\norm{\zeg}_s^2 \,. 
\label{II.est}
\end{equation}
Next we estimate \eqref{III}.
This time we exploit that  $\fp(x,\xi)$ is a  high-frequency symbol having the form  \eqref{2110:1812}.
We then use  estimate  \eqref{est.aR.comp} (with $\nu = \frac 1 2$, $m' = \frac 1 2$, $m = 2s$) to bound
\begin{equation}
\abs{\eqref{III}} \lesssim \frac{1}{\tR^{\frac12}} \, | \fp |_{2s, W^{2,\infty}, 7} \  | \Omega(\xi) + \wt\sfb(t, \cdot) |_{\frac12, W^{2,\infty}, 7} \, \norm{\zeg}_s^2
\stackrel{\eqref{a.sem}, \eqref{zmzn0},  \eqref{est.706} }{\lesssim } \ 
 \frac{\e^{2}}{\tR^{\frac12}}\norm{\zeg}_s^2 \,. 
\label{III.est}
\end{equation}
Next we estimate \eqref{IV}. 
We regard  $\fp(x, \xi)$ as  symbol in  $  \Gamma^{2s + 1}_{W^{2, \infty}}$ supported on high frequencies and exploit that $\wt \sfR(t)$ 
belongs to $H^{s+1}(\T;\C)$ (cf. \eqref{est.706}), and bound
\begin{equation}
  \abs{  \eqref{IV}} \lesssim 
  \norm{\fA \wt \sfR(t, \cdot)}_{-s} \norm{\zeg}_{s} 
\stackrel{  \eqref{cont00}}{\lesssim} 
 | \fp |_{2s+1, L^{\infty}, 7} \norm{ \wt \sfR(t, \cdot)}_{s+1}\norm{\zeg}_{s} 
\stackrel{  \eqref{a.sem}, \eqref{zmzn0},\eqref{est.706}}{\lesssim} \
  \frac{\e^{4-\theta}}{\tR} \norm{\zeg}_{s}  \,.
\label{IV.est}
\end{equation}
Finally, we estimate \eqref{V}. We regard $\fp(x,\xi)$ as a symbol in $\Gamma^{2s}_{W^{2,\infty}}$. By \eqref{cont00}
we get 
\begin{equation}\label{V.est}
\abs{\eqref{V}}   \lesssim 
\norm{\fA \wt \sfY(t, \cdot)}_{-s} \norm{\zeg}_{s} 
{\lesssim}
 | \fp |_{2s, L^{\infty}, 7} \norm{ \wt \sfY(t, \cdot)}_{s}\norm{\zeg}_{s} 
\stackrel{ \eqref{a.sem}, \eqref{zmzn0}, \eqref{est.706}}{\lesssim} \  \e^{5-\theta} \norm{\zeg}_{s} \,. 
\end{equation}
Altogether,  
 \eqref{20250130:1908}
 --\eqref{V.est} and the upper bound \be\label{20260128:2052}
 \norm{\zeg(t)}_s = \norm{z^\perp(t)}_s \leq \norm{z(t)}_s \leq  \e^{-\theta}
 \ee
 (recall that $z(t)$ is long-time controlled), imply that there are  a constant $C>0$ and $ \e_1=\e_1(s,\theta,\nu_0, \fr)>0$ such that, provided   $ \e \in (0,\e_1)$, 
 the  functional $\cA$ in  \eqref{cA}
  fulfills  for all times $0\leq t \leq T_\star$  the differential inequality
\begin{align}\label{dAt0}
\frac{\di}{\di t} \cA(t)\geq & \, (\tn-\tm) \, \upkappa  \,  \cA(t)
 -  C \left(   \frac{\e^{4-2\theta}}{\tR}  + \e^{5-3 \theta}  + \frac{\e^{2-2\theta}}{\tR^\frac12}    \right)\\
& \stackrel{\eqref{def:kappa}, \eqref{def:R}}{\geq }
(\tn - \tm) \nu_0 \e^2  \,  \cA(t) - 3 C \e^{5-3\theta}\,.
 \label{dAt}
\end{align}
Then $\cA(t)$ is a super-solution of the corresponding ODE, yielding the lower bound
\begin{equation*}
\cA(t) \geq e^{(\tn-\tm) \nu_0 \e^2 t} \, \left( \cA(0) - \frac{3C \e^{3-3\theta} }{(\tn-\tm) \nu_0}\right) +  \frac{3C \e^{3-3\theta} }{(\tn-\tm) \nu_0}  \ , \qquad  \forall \  0\leq t \leq T_\star  \,.
\end{equation*}
Thus, in order for $\cA(t)$ to grow in time, we need to select the initial data $z(0)$ so that  $\cA(0) >  \frac{3C \e^{3-3\theta} }{(\tn-\tm) \nu_0}$, a condition which is   guaranteed by (B3) of \Cref{B}. Indeed,   provided $\e$ is sufficiently small, 
\be
\cA(0) \stackrel{\eqref{cA}}{=} \la \fA \zeg(0), \zeg(0) \ra 
\stackrel{\eqref{zeta}}{=} \la \fA  z^\perp(0),  z^\perp(0) \ra  
\stackrel{\eqref{tRbound}}{>}
\e^{3-4 \theta} >
  \frac{6 C \e^{3-3\theta}}{(\tn-\tm) \nu_0}  \,. 
  \label{2110:1907}
 \ee
Then, using also the penultimate of the above inequalities,  $\cA(0) - \frac{3C \e^{3-3\theta} }{(\tn-\tm) \nu_0} >\e^{3-4 \theta}- \frac{3C \e^{3-3\theta} }{(\tn-\tm) \nu_0} > \frac12\e^{3-4\theta}$, and we get from \eqref{2110:1907}, the definition \eqref{cA} and the continuity Theorem \ref{thm:contS}, that, for any $t \in [0, T_\star]$
\be
\frac{1}{2} \e^{3-4\theta} \, e^{(\tn - \tm){\nu_0}\e^2 t} \leq \cA(t)   \leq  \norm{\fA_{s,\tR}\zeg(t)}_{-s} \norm{\zeg(t)}_s 
\stackrel{\eqref{fA}, \eqref{a.sem}}{\leq} C_s \e^2 \norm{\zeg(t)}_s^2 \stackrel{\eqref{20260128:2052}}{\leq} C_s \e^2 \norm{z(t)}_s^2 \,,
\ee
for some $C_s >1$. This proves \eqref{sopra.exp2}
with  $\tipac_s := 1/\sqrt{2 C_s} \in (0,1)$.
Then take    $ T_\star$  as in  \eqref{T0} to  get \eqref{zeta.grows}.
\end{proof}

 \section{Growth of Sobolev norms for $(\eta, \psi, \sfV, \sfB)(t)$}\label{sec:growth.z}
In the previous section we showed that long-time controlled solutions  with strongly well-prepared initial data and 
existing up to  the enhanced time scale \eqref{T0}
undergo Sobolev norm explosion, as proved in \Cref{prop:instab}.
The goal of this section is to leverage that analysis in order to prove the following result:
\begin{theorem}\label{thm:main}
Let the vorticity $\gamma $ fulfill \eqref{resonance.intro}.
There exists $\ts_0 > \frac32$  such that the following holds true.
    Let  $\ts> 5\ts_0$,  $0<\delta \ll 1$, $K \gg  1$   be given parameters. There exist a  time  $T>0$ and  a solution $(\eta, \psi, \sfV, \sfB) \in C^0([0,T]; X^\ts)$ of \eqref{eq:etapsi} 
    such that 
    \begin{equation}\label{explosion}
      \| (\eta,\psi,\sfV,\sfB)(0)\|_{X^\ts} \leq \delta  \,, 
      \qquad 
       \| (\eta,\psi,\sfV,\sfB)(T)\|_{X^\ts} \geq K
    \end{equation}
    and \begin{equation}\label{small_in_small_norm}
        \sup_{0 \leq t \leq T} \| (\eta,\psi,\sfV,\sfB)(t)\|_{X^{\ts_0}} \leq 2 \delta \,. 
    \end{equation}
\end{theorem}
The first step to  prove such result is to show that  strongly well-prepared initial data exist,  cf.  \Cref{lem:z0.wp}. 
Then, in \Cref{lem:cresce},
we show that {\em any} solution $z(t)$ of system 
\eqref{Z.eq} with a strongly well-prepared initial datum  undergoes Sobolev norm explosion.
We use this to prove, in \Cref{sub.proof_main},  that the couple $(\eta(t), \upomega(t))$ of the free surface profile and good unknown of Alinhac undergo growth of Sobolev norms and, by a dichotomy argument, we deduce \Cref{thm:main}.

 \subsection{Growth of Sobolev norms for the complex variable $z(t)$}

We start  proving that  strongly well-prepared initial data do exist, and exhibit an explicit example having support on  only four Fourier modes: two of them  in the $\gamma$-good set $\Set$, the other two at very high frequencies.
 \begin{lemma}[{\bf Existence of Strongly well-prepared initial data}]\label{lem:z0.wp}
Let $\Set$ in \eqref{fixLambda} and 
  $s, \theta$ as in \eqref{s0r},  $ \tV_\tm>0, \tV_\tn<0$ in  \eqref{VmVn.def}, $\sfV_{\tm, \tn}^{(\res)}\neq 0$ in \eqref{Vres.coeff} and define 
\be\label{sceltanufix}
\nu_0:= \abs{\sfV_{\tm, \tn}^{(\rm res)} } \frac{ \sqrt{-\tV_\tn} \sqrt{\tV_\tm}}{4(\tV_\tm - \tV_\tn)}>0  \,.
\ee
There exists  $ \e_0   >0$ such that, for any $ \e \in (0, \e_0)$, 
the initial datum 
\begin{align}
\label{z0}
 z(0,x) & := \e  \frac{\sqrt{-\tV_\tn}}{2\sqrt{\tV_\tm - \tV_\tn}} e^{\im \tm x }  + \e  \frac{\sqrt{\tV_\tm}}{2\sqrt{\tV_\tm - \tV_\tn}}  e^{\im \tn x} + \uprho\, e^{\im 3\tN x} - \im  \uprho \,  e^{\im (3\tN + \tn-\tm) x}  \\
 \label{param.z(0)}
 \mbox{where } 
 &  \uprho := \uprho(\e):= \e^{\theta + \und{\theta} +2s(3+\theta)} \ , \quad 
\und{\theta} = \frac{1-6\theta}{4}\,, \quad \tN:= \tN(\e):= \lceil \tR  \rceil = \lceil \e^{-2(3 + \theta)}\rceil\,, 
 \end{align}
is strongly well-prepared with parameters $(s, \theta, \nu_0, \e)$ (according  to \Cref{B}).
In particular, the 
 set of  strongly well-prepared initial data 
 is non-empty and open in $\dot H^s(\T; \C)$.
 \end{lemma}
 \begin{proof}
 The set of strongly well-prepared initial data is,
  open   since conditions (B1)--(B4) are open with respect to the $H^s_\R$ topology. 
  Then it is enough to show that 
  $z(0,x)$  is strongly well-prepared with parameters $(s, \theta, \nu_0, \e)$. 

\smallskip 

\noindent{\sc {Verification of (B1)}:} 
We show that  $z(0,x)$ in \eqref{z0} is 
 $\e$-weakly well-prepared, namely \eqref{app.ass1} holds true. Recalling the definitions of $z^\top(0)$ and $z^\perp(0)$ in \eqref{Pi}, we have
 $$
 \norm{z^\top(0)}_{L^2}^2 \stackrel{\eqref{z0}}{=} \frac14 \e^2  < \e^2 \ , \quad 
 \norm{z^\perp(0)}_{L^2}^2 \stackrel{\eqref{z0}}{=} 2\uprho^2  \stackrel{\eqref{param.z(0)}}{<} \e^4 
 $$
since $s > \frac23$ and provided $\e$ is small enough.\\
\noindent{\sc {Verification of  (B2)}:} 
Note that $z_\tm(0)$ and $z_\tn(0)$ are chosen so that
$$
\tV_\tm |z_{\tm}(0)|^2 
  + \tV_\tn |z_{\tn}(0)|^2 = 0 \,;
  $$
then 
\[
\upkappa \stackrel{\eqref{def:kappa} }{=} \e^2  2\abs{\sfV_{\tm, \tn}^{(\rm res)} } |z_\tm(0) \, z_\tn(0)|  \stackrel{\eqref{sceltanufix}}{=} 2 \nu_0 \e^2\,.
\]
{\sc Verification of (B3):} Recall that $\fA =\Opbw{\fp(x,\xi)}$,   with the symbol $\fp(x,\xi)$ defined in  \eqref{fA} which, in  view of the choice of  $z(0,x)$ in \eqref{z0}, has Fourier-coefficients in the $x$-variable given by
$$ 
\hat \fp(\ell, \xi) =
\begin{cases}
    \frac{\ii}{2} \e^2 \rho_\tn \rho_\tm \sfV^{(\res)}_{\tm, \tn} |\xi|^{2s} \eta_\tR^2(\xi) & \mbox{ if } \ell = \tn - \tm \\[0.5em]
    -  \frac{\ii}{2}\e^2 \rho_\tn \rho_\tm  \sfV^{(\res)}_{\tm, \tn} |\xi|^{2s} \eta_\tR^2(\xi)  & \mbox{ if } \ell = \tm - \tn\\[0.5em]
    0 & \mbox{ otherwise}
\end{cases}\,,
\qquad 
 \rho_\tm:= \frac{\sqrt{-\tV_\tn}}{2\sqrt{\tV_\tm - \tV_\tn}} \,, \ \ 
 \rho_\tn:= \frac{\sqrt{\tV_\tm}}{2\sqrt{\tV_\tm - \tV_\tn}}
$$
Hence,  using the definition of paradifferential operator in 
\eqref{BW},    we get 
 \begin{align*}
 \la \fA \Pi^\perp z(0), \Pi^\perp z(0) \ra & = 
 \sum_{k-j = \tn - \tm } \chi \big(\tn-\tm, j+\frac{\tn-\tm}{2} \big) \, 
 \hat\fp(\tn - \tm, j+\frac{\tn-\tm}{2}) \, z_j^\perp(0) \bar{z_k^\perp(0)}\\
 & +  \sum_{k-j = -(\tn - \tm) } \chi \big(\tm-\tn, j-\frac{\tn-\tm}{2} \big) \, 
 \hat\fp(\tm - \tn, j-\frac{\tn-\tm}{2}) \, z_j^\perp(0) \bar{z_k^\perp(0)}\,.
 \end{align*}
Since $z^\perp(0) = \uprho e^{\im 3\tN x} + \im \uprho e^{\im  (3\tN + \tn - \tm)x}$, the only term surviving in the first sum is the one with $k = 3\tN + \tn - \tm$, $j = 3\tN$ (hence $z_j^\perp(0) = \uprho$, $z_k^\perp(0) = -\im\, \uprho$), whereas in the second sum only the term with $k = 3\tN$ and $j = 3\tN + \tn - \tm$ survives (correspondingly  $z_j^\perp(0) = -\im \, \uprho$, $z_k^\perp(0) = \uprho$.
Plugging in these values  and since, by  \eqref{Vres.coeff} and $\tm <0<\tn$, $\sfV^{(\res)}_{\tm, \tn} ={\dfrac{\tm \tn (\tm + \tn)^3}{2(\tn-\tm)\sqrt{\tm \tn(\tm-3\tn)(\tn-3\tm)}}}<0$, we arrive at 
    \begin{equation*}
 \la \fA \Pi^\perp z(0), \Pi^\perp z(0) \ra  
 %
  =  {\e^2 \abs{\sfV^{(\res)}_{\tm, \tn}} \rho_{\tm} \rho_{\tn}} \, \abs{3\tN+\frac{\tn-\tm}{2}}^{2s} \, \eta_\tR^2\big(3\tN+\frac{\tn-\tm}{2}\big) \, \chi \big( \tn-\tm, 3\tN+\frac{\tn-\tm}{2} \big) \,  \uprho^2 \,.
 \end{equation*}
By taking $\tN$ sufficiently large -- i.e. in view of \eqref{param.z(0)}, $\e$ small enough --
the terms 
$\chi \big( \tn-\tm, 3\tN+\frac{\tn-\tm}{2} \big)  =  
\eta_\tR^2\big(3\tN+\frac{\tn-\tm}{2}\big) = 1$, 
 and we get
 $$
  \la \fA \Pi^\perp z(0), \Pi^\perp z(0) \ra \geq 
 \e^2 \abs{\sfV^{(\res)}_{\tm, \tn}} \rho_{\tn} \rho_{\tm} \abs{2\tN}^{2s} \,  \uprho^2  \stackrel{\eqref{param.z(0)}}{>} \e^{3-4\theta} 
 $$
eventually shrinking $\e$. 
Note that in the last inequality we used that   $\e^{\und{\theta}} \gg \e^{\frac{1-6\theta}{2}}$ (thanks to the choice of $\und{\theta}$), which also gives the second condition in the choice of $\theta$ in \eqref{s0r}.
\\
\noindent {\sc {Verification of (B4)}:} To impose  \eqref{z0.s} we require
$$
\norm{z(0)}_s^2 = \e^2 \la \tm \ra^{2s} \rho_\tm^2 + \e^2 \la \tn\ra^{2s} \uprho_\tn^2 + 
  \uprho^2 \la 3\tN \ra^{2s}+  \uprho^2 \la 3\tN+\tn - \tm \ra^{2s}\leq \e^{2\theta}\,.
  $$ 
  Each of the first two terms  is smaller than  $\frac14\e^{2\theta}$  by shrinking $\e$ enough. 
  The same is true for the third and fourth term, by the choice of $\uprho$ and $\tN$ in \eqref{param.z(0)}.
 \end{proof}
 \begin{remark}
 Let us summarize the compatibility constrains that we need on $\e, \tR, \uprho$ in order to impose all the required conditions. 
 First, in  \eqref{dAt0} we want the dominant term to be $\e^{5-3\theta}$, hence we impose 
 $
 \tR \geq \e^{-2(3-\theta)}$. In order to verify (B1) we need $\uprho \ll \e^2$. To verify (B3) we need $\tR^{2s} \uprho^2 \gg \e^{1-4\theta}$ and finally to verify (B4) $\tR^{2s} \uprho^2 \ll \e^{2\theta}$.
\end{remark}
The next proposition finally proves that any solution of system \eqref{Z.eq} with a well-prepared initial datum as in 
 Lemma \ref{lem:z0.wp} undergoes Sobolev norm explosion. 
 \begin{proposition}[{\bf Growth of Sobolev norms for the variable $z(t)$}]\label{lem:cresce}
  Let $s, \fs_0, \theta$ and $\fr$ as in \eqref{s0r}. 
 There exist $\e_2, C >0$ such that for any  $\e \in (0, \e_2)$ the following holds true. 
 Consider the solution $z(t)$ of system  \eqref{Z.eq}  with an initial datum $z(0)\in \dot H^s(\T; \C)$ 
  strongly well-prepared with   parameters $(s, \theta, \nu_0, \e)$,   $\nu_0$ in \eqref{sceltanufix}.
 Denote by 
 \be\label{T1}
0<\tT:= \tT(s,\e, z(0)) := \sup\left\lbrace t \geq 0 \colon \ \ \sup_{\tau \in [0,t]}\norm{z(\tau)}_s \leq \frac{\tipac_s}{2}\e^{-\theta} \quad \text{and} \quad  
\sup_{\tau \in [0,t]}
\norm{z(\tau)}_{\fs_0} \leq \fr \right\rbrace \,,
 \ee
with $\tipac_s \in (0,1)$ the constant in Proposition \ref{prop:instab}. Then $\tT$ is finite and bounded by 
$$
\tT \leq T_\star := \frac{1-2\theta}{(\tn - \tm)\nu_0} \e^{-2} \log(\e^{-1})  \,. 
$$
 Moreover, one has
 \begin{align}\label{cresce}
\sup_{0 \leq t \leq \tT}  \norm{z(t)}_{\fs_0} \leq C \e \ , 
 \quad  \norm{z(0)}_s \leq \e^\theta \,,  \\
 \label{cresce2}
\frac{\tipac_s}{4} \e^{-\theta} \leq \norm{z(t)}_{s} \leq \frac{\tipac_s}{2}\e^{-\theta} \,, \quad \forall \, t \in [\tT-1, \tT] \,. 
\end{align}
 \end{proposition}
 \begin{proof}
 The point of the proof is  that, since   the initial data is strongly well-prepared and  the solution $z(t)$ is long-time controlled with parameters $(s, \theta, \tT, \e)$, then by \eqref{sopra.exp2}
 its norm $\norm{z(t)}_{s}$ must be above the graph of an increasing exponential function which passes through $\tipac_s \e^{-\theta}$ at $t = T_\star$, hence it must itself reach the value $\frac12\tipac_s \e^{-\theta}$ at the  earlier time $t=\tT$. The subtlety is that we must guarantee that along the interval $[0,\tT]$ the solution exists and remains bounded by $\fr$ in the low norm $\norm{ \cdot }_{\fs_0}$ ; this will be achieved by a bootstrap argument. We give now the details.

Define $\wt \e_2:= \min(\e_\star, \e_0, \e_1, \fr)$ with
$\e_\star$ of Lemma \ref{lem:long.boot},  
  $\e_0$ of Lemma \ref{lem:z0.wp}, and $\e_1$ of Proposition \ref{prop:instab}.
  We shall take $0<\e < \e_2 <\wt \e_2 \ll \fr $, so that we can use the lemmata and proposition just mentioned.
  
First we observe that, due to local well-posedness of the autonomous Cauchy problem for $Z(t)$ in $B_{\fs_0}(\fr) \cap H^s$  (see \Cref{loc.ex}), on the interval $[0, \tT]$ the solution keeps existing and is continuous in  $B_{\fs_0}(\fr) \cap \dot H^s_\R(\T, \C^2)$.  Thus as $t \to \tT^-$
\be\label{eq:tocco.la.frontiera}
\text{either} \quad \|z(t)\|_{s} \to \frac{\tipac_s}{2}\e^{-\theta}\,, \quad \text{or} \quad \|z(t)\|_{\fs_0} \to \fr\,.
\ee
Assume now by contradiction that $\tT > T_\star$. Then the solution $z(t)$ is long-time controlled with parameters $(s, \theta, T_{\star}, \e)$ (according to \Cref{A}). 
Indeed,   condition (A1)  holds true since by assumption the data is strongly well-prepared, and thus in particular also $\e$-weakly well-prepared, whereas condition (A2) holds true by the very definition of $\tT$, with the continuity in time being granted again by \Cref{loc.ex}.
But then the assumptions of Proposition \ref{prop:instab} -- in particular \eqref{T0} -- are met, therefore 
 $$
  \sup_{0 \leq t \leq  T_\star } \norm{z(t)}_s \geq \tipac_s \e^{-\theta} \,,
 $$
and since $T_\star < \tT$, this contradicts the definition of $\tT$.
This proves that $\tT \leq {T_\star}$.
 
 We now prove \eqref{cresce}. We know that $z(t)$ is long-time controlled with parameters $(s, \theta, \tT, \e)$ and that  $\tT$ satisfies \eqref{Tstar.T0}, thus the bootstrap Lemma \ref{lem:long.boot} applies. Hence, by \eqref{z.s0}, and taking $\e$ small enough,
 \be\label{z-s0-venerdi}
 \sup_{0\leq t \leq \tT} \norm{z(t)}_{\fs_0} \leq C\e <  \fr \,,
 \ee
which is the first estimate of \eqref{cresce}.  
 The second inequality of \eqref{cresce} is satisfied by the
 strongly well-prepared assumption (B.4) on the initial datum.
Then, by \eqref{eq:tocco.la.frontiera}, \eqref{z-s0-venerdi} and  the very definition of $\tT$, we have
\be\label{eq:zisbig}
\norm{z(\tT)}_s = \frac{\tipac_s}{2}\e^{-\theta}\,.
\ee
It remains to prove \eqref{cresce2}. We differentiate the $\dot H^s_\R(\T;\C^2)$ norm of $Z(t) = \vect{z(t)}{\bar z(t)}$ using equation \eqref{Z.eq}, and recall  that we put $N = 3$ (cf. \eqref{paraN}).
We get the energy inequality
        \begin{align}\notag
       \frac{\di}{\di t}\| Z(t)\|_{s}^2 \leq  & \Big|\la\Big[| D|^{2s}, - \im \bOmega(D) Z + \vOpbw {\ii \sm_{\geq 2}^{\rm (res)}(Z; x,\x)}\Big]Z, \, Z \ra\Big| \\
       \notag
            &+ 2 \Big|\la | D|^{2 s} \left( X^{(\Set)}(Z) + \bB_{\geq N}(Z) Z+  \bR_{\geq 3}(Z)[Z]\right) ,\, Z  \ra\Big|\\
            &\lesssim_s \norm{Z(t)}_{\fs_0}^2
            \| Z(t)\|_s^2 \stackrel{\eqref{cresce}}{\lesssim_s}  \e^{2-2\theta}\,,
            \label{3001:1517}
        \end{align}
        where  in the first line we used that   $\sm_{\geq 2}^{\rm (res)}(Z; x,\x)$ in \eqref{D.final}  is  a real-valued symbol in $\Sigma \Gamma^1_2[r,3]$ and 
 to estimate the second line       that $\bB_{\geq N}(\zak)$, being a spectrally localized map in  $\cS_{\geq N}^0[r]$,  fulfills the bounds \eqref{bound:specloc}, and that $X^{(\Set)}(Z) = \bR_2^{(\Set)}(Z)Z$ with $\bR_2^{(\Set)}(Z)$, $\bR_{\geq 3}(\zak)$  real-to-real matrices of smoothing operators in 
 $\wt \cR_2^{-\vr +{3N+\frac32}}$ respectively $\cR_{\geq 3}^{-\varrho+{3N+28}}[r]$ (cf. \Cref{thm:nf}), thus 
 fulfilling \eqref{piove} with $m =0$.
Integrating the energy inequality in the interval $[t, \tT]$, for any $t \in [\tT - 1,\, \tT]$ we get
$$
\norm{Z(t)}^2_s \geq \norm{Z(\tT)}^2_s - |\tT - t| \sup_{\tau \in [t,\,\tT]} \frac{\di}{\di t} \|Z(\tau)\|_s^2 \stackrel{\eqref{eq:zisbig}}{\geq} \frac14\tipac_s^2 \e^{-2\theta} - C_s \e^{2-2\theta}\,,
$$
which gives \eqref{cresce2} up to shrinking again $\e$.
 \end{proof}
\begin{remark}\label{rem:enhanced}
  By the energy estimate \eqref{3001:1517} and the upper bound in \eqref{z-s0-venerdi},  there is $C_s >0$ such that
  $$
\frac{\tipac_s}{2}\e^{-\theta} \stackrel{\eqref{eq:zisbig}}{=} \norm{z(\tT)}_s \leq e^{C_s \e^2 \tT} \norm{z(0)}_s \stackrel{\eqref{z0.s}}{\leq} e^{C_s \e^2 \tT} \e^{\theta} 
\quad \Longrightarrow \quad  \tT \geq \frac{2 \theta}{C_s \e^2} \log \left( \left( \frac{\tipac_s}{2}\right)^{\frac{1}{2\theta}}\frac 1 \e\right)\,,
  $$
  showing that the solution  must indeed exist on timescales longer than $\e^{-2}$ given by the Cauchy theory. 
\end{remark}
\subsection{Growth of Sobolev norms  for $(\eta(t),\upomega(t))$
and $(\eta, \psi, \sfV, \sfB)(t)$ 
}\label{sub.proof_main}
By \Cref{lem:cresce} there exists a solution  $Z(t)$ 
of \eqref{Z.eq} with growing Sobolev norms. 
In view of \Cref{UZequiv}, 
 the complex variable $U(t)$ in \eqref{complexo} 
 has equivalent norms, hence its high $H^s$-Sobolev norms grow in time. 
In addition,  recalling the equivalence of norms of  $U = \vect{u}{\bar u}$ and of the couple $(\eta, \upomega)$ -- with $\eta$ the original free surface profile and $\upomega$ the good unknown of Alinhac -- see \eqref{eq:etaomegaUequiv}, 
also the high Sobolev norm $\norm{\eta(t)}_{s-\frac14} + \norm{\upomega(t)}_{s+\frac14}$ grows in time. 
Actually, we get the following   result regarding the  variables $(\eta, \psi, \sfV, \sfB)(t) \in X^{{s-\frac 3 4}}$ (recall \eqref{norm:Xs}).
\begin{proposition}\label{lem:partial_growth}
 Let $s, \theta,  \fs_0$ as in  \eqref{s0r}.  
There exist $ \tC, \tC_s  >1$, $\e_0 >0$  such that for any $\e \in (0, \e_0)$, there is $\tT >0$ and a solution $(\eta, \psi, \sfV, \sfB) \in C^0([0, \tT];  X^{s-\frac34})$ of \eqref{eq:etapsi}  fulfilling 
\begin{align}
&\norm{(\eta, \psi, \sfV, \sfB)(0)}_{X^{s-\frac34}} \leq \tC_s \e^\theta \,, 
\qquad \norm{(\eta, \psi, \sfV, \sfB)(\tT)}_{X^{s-\frac34}}  \geq \tC_s \e^{-\theta } \,.
\label{Xs_initial}  
\end{align}
In addition, recalling that $\upomega(t)$ is the good unknown in \eqref{GU}, one has the following estimates:
\begin{align}
\label{eta_omega_gr}
\tC_s^{-1} \e^{-\theta} \leq  & \norm{\eta(t)}_{s-\frac14} + \norm{\upomega(t)}_{s+\frac14} \leq \tC_s \e^{-\theta} \,, \quad \forall \, t \in [\tT - 1, \tT ] \,, \\
\label{upp.sol.low}
&  \sup_{t \in [0, \tT]} \left( \norm{(\eta, \psi, \sfV, \sfB)(t)}_{X^{\fs_0 -{\frac 3 4}}} + \norm{\upomega(t)}_{\fs_0 +{\frac{1}{4}}} \right) \leq \tC  \e \,, \\
\label{upp.sol.high}
&   \sup_{t \in [0, \tT]}  \left( \norm{(\eta, \psi, \sfV, \sfB)(t)}_{X^{s-\frac34}} + \norm{\upomega(t)}_{s+\frac14} \right) \leq \tC_{s}  \e^{-\theta} \,. 
  \end{align}
\end{proposition}
\begin{proof}
Apply \Cref{lem:cresce}  with initial datum $z(0,x)$ in \eqref{z0}, \eqref{param.z(0)} (which is ($s, \theta,  \nu_0, \e)$ strongly well-prepared), to construct a solution $Z(t) = \vect{z(t)}{\bar z(t)} \in \dot H^s_\R(\T, \C^2)$ fulfilling \eqref{cresce}, \eqref{cresce2} 
with  $\tT:=  \tT(s, \e, z(0))$ as in \eqref{T1}, so in particular 
$Z(t) \in B_{\fs_0, \R}(C\e)\cap \dot H^s$.
Now we wish to exploit the equivalences of norms 
given by \Cref{UXequiv} and \Cref{UZequiv}, 
which 
hold provided 
 $U(t)$ belongs to the small ball $B_{\s_0'}(\tr'(s))$ and $Z(t)$ belongs to the  small ball 
 $B_{\s_0''}(\tr(s))$; this clearly can be achieved by shrinking $\e$ and recalling that $\fs_0$ has been taken in \eqref{s0r} larger than both $\s_0', \s_0''$.
Hence,
for any $\s \in [\fs_0,  s]$,
\begin{equation*}
\norm{Z(t)}_{\s} 
 \simeq 
\norm{U(t)}_{\s} \stackrel{\eqref{eq:etaomegaUequiv}}{\simeq_\s }\norm{\eta(t)}_{\s-\frac14} + \norm{\upomega(t)}_{\s+\frac14} \simeq_\s \norm{(\eta, \psi, \sfV, \sfB)(t)}_{X^{\s - \frac 34}}  \,, \quad \forall\, t \in [0,\tT] \,,
\end{equation*}
and \eqref{Xs_initial}--\eqref{upp.sol.high} hold true due to the analogous properties of $z(t)$ in \eqref{T1}--\eqref{cresce2}.
\end{proof}

\begin{proof}[{\bf Proof of \Cref{thm:main}}]\label{proof:thm1}
Let the vorticity $\gamma <0$ such that $\gamma^2 \in \Q$ and $0<\delta \ll 1$, $K \gg  1$   be given parameters.
Let $ \fs_0, \fr$ be given by \eqref{s0r}, 
and put $\ts_0:= \fs_0 - \frac34$.
Since $\ts > 5 \ts_0$, the number $s := \ts +\frac34$
fulfill
$s > 4 \fs_0$ (use that $\fs_0 > 3$).
Then take $\theta = \theta(s)$ to fulfill \eqref{s0r}.\\
\Cref{lem:partial_growth} gives $\e_0>0$ such that, for any $\e \in (0, \e_0)$, there exists a solution 
 $(\eta, \psi, \sfV, \sfB) \in C^0([0, \tT];  X^\ts)$ of
 \eqref{eq:etapsi} fulfilling 
 \eqref{Xs_initial}--\eqref{upp.sol.high}.
 In particular, shrinking $\e$, we get
 \be\label{3001:1642}
\begin{aligned}
& \norm{(\eta, \psi, \sfV, \sfB)(0)}_{X^\ts} \leq \tC_\ts \e^\theta \leq \delta \ , \quad 
\sup_{t \in [0, \tT]} \norm{(\eta, \psi, \sfV, \sfB)(t)}_{X^{\ts_0}}   \leq \tC  \e \leq 2\delta\,,  \\
& \| (\eta,\psi,\sfV,\sfB)(T)\|_{X^\ts}   \geq \tC_\ts^{-1}\e^{-\theta} \geq K \,,
\end{aligned}
\ee
 proving  the theorem.
\end{proof}

 \section{Separated growth of Sobolev norms for the physical variables}\label{sec:growth.physical}
In the previous section, we constructed a solution to \eqref{eq:etapsi} for which the pair $(\eta(t), \upomega(t))$ exhibits arbitrarily large growth of Sobolev norms (see \Cref{lem:partial_growth}); in particular, at least one among $\eta(t)$ and $\upomega(t)$ undergoes Sobolev norm growth.
In this section, we show that all quantities of interest -- the free surface profile $\eta(t)$, the good unknown $\upomega(t)$, and the horizontal and vertical  components $\sfV(t)$ and $\sfB(t)$ of the irrotational part of the velocity field at the free interface -- undergo growth of Sobolev norms, thereby proving \Cref{thm:mainimprove}.

The key step in passing from the growth of the pair $(\eta(t), \upomega(t))$ to the growth of the individual variables is to exploit the wave-like structure of the paralinearized system \eqref{ParaWW}. Owing to the linear dynamics, the variables $\eta(t)$ and $\upomega(t)$ exchange energy over short time intervals, and we use this property to prove that actually both components separately undergo growth of Sobolev norm.

However, since the water waves system is quasilinear and the nonlinear terms constitute a singular perturbation of the linear part, the linear dynamics only provides a heuristic description of the true evolution. In particular, the quasilinear interactions prevent a clean decoupling of the linear wave mechanism from the nonlinear effects, and this difficulty becomes especially pronounced at high frequencies. This represents a genuine structural obstruction in the analysis.

To overcome this difficulty, we introduce an new upside-down virial argument, which allows us to recover the exchange of energy indirectly through lower bounds on the second order time derivatives of the Sobolev norms $\eta(t)$ and $\upomega(t)$, see \eqref{est.der2.eta}, \eqref{stima.der2.om}. Precisely we have:

\begin{proposition}[{\bf Upside-down virial argument}]\label{lem:UD}
Let $s$ and $\tC_s>1$ be as in \Cref{lem:partial_growth}, and set 
$
\ts := s - \frac{3}{4}
$.
We denote by $\tC_{\ts}:=\tC_{s - \frac{3}{4}}>0 $ the corresponding constant.
There exist  $\mu_\ts >0, $ and for any $\mu \in (0, \mu_\ts)$ there exists $\e_0 >0$ such that for any $\e \in (0,\e_0)$, the following holds true. Let $\tT$,  $(\eta(t), \psi(t), \sfV(t), \sfB(t)) \in C^0([0, \tT];  X^\ts)$ and  $\upomega(t)\in C^0([0, \tT];  H^{\ts+1}(\T;\R))$ be the time and the solutions  of \Cref{lem:partial_growth}, 
then one has: 
\begin{itemize}
    \item[(i)] {\bf Upside-down virial for $\eta(t)$:}
   Assume that $\norm{\eta(t)}_{\ts+\frac12}$ does not grow and quantitatively
\be\label{dicotomia1}
\sup_{t \in [\tT-1, \tT]}  \norm{\eta(t)}_{\ts+\frac12} < \mu \, \tC_\ts^{-1} \e^{-\theta}  \,,
\ee
then for every $t \in [\tT-1, \tT]$ one has 
    \begin{align}
    \label{est.der.eta}
        \left|\frac{\di}{\di t} \norm{\eta(t)}_\ts^2 \right| &\leq 3\mu \e^{-2\theta}
         \,, \\
     \frac{\di^2}{\di t^2} \| \eta(t)\|_\ts^2 
     &\geq \|\upomega(t)\|_{\ts + 1}^2 \geq (1-\mu)^2 \tC_{\ts}^{-2} \e^{-2\theta}\,. 
    \label{est.der2.eta}
    \end{align}
    \item[(ii)] {\bf Upside-down virial for $\upomega(t)$:} Assume that 
    $ \norm{\upomega(t)}_{\ts+1}$ does not grow and quantitatively
    \be\label{dicotomia3}
\sup_{t \in [\tT-1, \tT]}  \norm{\upomega(t)}_{\ts+1} \leq \mu \, \tC_\ts^{-1} \e^{-\theta}  \,,
\ee
then for every $t \in [\tT-1, \tT]$ one has
  \begin{align}
    \label{stima.der.om}
      \left|\frac{\di}{\di t} \norm{\upomega(t)}_{\ts+\frac12}^2 \right| &\leq 3 \mu \e^{-2\theta}\,,  \\
     \frac{\di^2}{\di t^2} \|\upomega(t)\|_{\ts+\frac12}^2 
     & \geq \norm{\eta(t)}^2_{\ts + \frac  12} \geq (1-\mu)^2 \tC_\ts^{-2} \e^{-2\theta}\,.  \label{stima.der2.om} 
    \end{align}
\end{itemize}
\end{proposition}
We shall prove  \Cref{lem:UD}$-(i)$ in  \Cref{UP.eta} and  \Cref{lem:UD}$-(ii)$ in 
\Cref{UP.omega}.  
First we collect some properties of the  paralinearized equation \eqref{ParaWW}, that we write compactly  as 
\be\label{eq.8}
\begin{cases}
& \pa_t \eta = |D| \upomega + \cB^{(-1)} \upomega + \cV \eta +f_1(t) \\
& \pa_t \upomega = - \eta + \gamma \Hilb \upomega + \cB^{(-2)} \upomega - \cV^\top \upomega + \cT \eta  + f_2(t)
\end{cases}\,, 
\ee
having denoted  by 
\be
\begin{aligned}
& \cV:= \Opbw{-\ii \xi \#_\vr \sfV_\gamma }, \qquad \sfV_\gamma \mbox{ in }  \eqref{def:Vgamma} \,, 
\qquad \vect{f_1}{f_2} := \bR(\eta, \psi) \vect{\eta}{ \upomega}\,, \\
&\cB^{(-1)}:= \Opbw{ b_{\geq 1}^{(-1)}(\eta;x,\xi)} \ , 
\qquad \cB^{(-2)}:= \Opbw{\gamma \frac{1}{\im \xi}\#_\vr  b_{\geq 1}^{(-1)}} \,, 
\qquad 
 \cT:= \Opbw{-\sfa} \,.
 \label{barcola}
\end{aligned}
\ee
The transpose of $\cV$,    recalling \eqref{A1b}, \eqref{realetoreale} and \eqref{prop:ov},  is given by
\be\label{Vtra}
\cV^\top= \Opbw{\sfV_\gamma\#_\vr\im  \x}   \,.
\ee
We now establish estimates for the operator norms of these operators and their time derivatives when evaluated along the solution $(\eta, \psi, \sfV, \sfB) \in C^0([0, \tT];  X^{\ts})$ of \eqref{eq:etapsi} constructed in \Cref{lem:partial_growth}, which satisfies the following bounds (recall $\ts = s-\frac34$ and $\ts_0 = \fs_0 - \frac34$ as stated in the beginning of the proof of \Cref{thm:main}):
\begin{align}\label{growth.final}
& \tC_\ts^{-1} \e^{-\theta} \leq   \norm{\eta(t)}_{\ts+\frac12} + \norm{\upomega(t)}_{\ts+1} \leq \tC_s \e^{-\theta} \,, \quad \forall \, t \in [\tT - 1, \tT ]\\
\label{upp.sol.low9}
&  \sup_{t \in [0, \tT]} \left( \norm{(\eta(t), \psi(t), \sfV(t), \sfB(t))}_{X^{\ts_0}} + \norm{\upomega(t)}_{\ts_0+1} \right) \leq \tC  \e \,, \\
\label{upp.sol.high9}
&   \sup_{t \in [0, \tT]}  \left( \norm{(\eta(t), \psi(t), \sfV(t), \sfB(t))}_{X^{\ts}} + \norm{\upomega(t)}_{\ts+1} \right) \leq \tC_{s}  \e^{-\theta} \,.
\end{align}
Exploiting the continuity \Cref{thm:contS}, the paracomposition estimate \eqref{sharp.est},  the bounds on the seminorms of the symbols
$\sfV_\gamma, \pa_t \sfV_\gamma, \sfa, \pa_t \sfa$ in \eqref{a.v.b.estimate}, \eqref{stima.der.BVa}, the seminorms of  the symbols $b^{(-1)}_{\geq 1}, \pa_t b^{(-1)}_{\geq 1}$ in \eqref{b-1dert}, and the low norm bound in \eqref{upp.sol.low9}, one gets that
for any $\s \in \R$
\begin{align}
    \label{E1.est.9}
&\norm{\cE_1(t)}_{\cL(H^\sigma, H^{\sigma-1})} + \norm{\pa_t \cE_1(t)}_{\cL(H^\sigma, H^{\sigma-1})} \lesssim_\s \e \ , \quad \forall \cE_1 \in \{\cV, \cV^\top\} \,, \\
    \label{T.est.9}
  &  \norm{\cT(t)}_{\cL(H^\sigma, H^{\sigma})} + \norm{\pa_t \cT(t)}_{\cL(H^\sigma, H^{\sigma})} \lesssim_\s \e \,, \\
  \label{E2.est.9}
   &   \norm{\cE_2(t)}_{\cL(H^\sigma, H^{\sigma+1})} + \norm{\pa_t \cE_2(t)}_{\cL(H^\sigma, H^{\sigma+1})} \lesssim_\s \e \ ,  \quad \forall \cE_2 \in \{\cB^{(-1)}, \cB^{(-2)}\} \,.
\end{align}
Analogously, using also the control on the norms of the smoothing remainders in \eqref{Rt_smooth}, one gets that the forcing terms $f_1(t), f_2(t) $ satisfy 
\begin{align}
    \label{f.est.9}
    &\norm{f(t)}_{\ts +1} + \norm{\pa_t f(t)}_{\ts+1} \lesssim \e^{1-\theta} \,, \quad \forall f \in \{f_1, f_2\} \,.
\end{align}
As an immediate application, using also 
\eqref{upp.sol.high9}, we get for any $t \in [0, \tT]$ the bounds
\be\label{de.t.om}
\norm{\pa_t \eta(t)}_{\ts - \frac12} + \norm{\pa_t \upomega(t)}_{\ts} \lesssim 
\norm{(\eta(t), \psi(t), \sfV(t), \sfB(t))}_{X^{\ts}} + \norm{\upomega(t)}_{\ts+1}
\lesssim
\e^{-\theta} \,.
\ee
Here we  establish the crucial cancellations underlying  the upside-down virial arguments for $\eta(t)$ and $\upomega(t)$.
 \begin{lemma}\label{lem:A2s}
    With the same assumptions as in \Cref{lem:UD}, the following holds. 
 The operators 
      \begin{equation}\label{lambda:A}
      A^{(2\ts)}:= |D|^{2\ts} \cV + \cV^\top |D|^{2\ts}  \,, 
      \qquad
       \wt A^{(2\ts+1)} := - |D|^{2\ts+1}\cV^\top - \cV |D|^{2\ts+1} 
        \end{equation}
        belong  to $\cL\left( H^\s;H^{\s-2\ts}\right)$ respectively $\cL\left( H^\s;H^{\s-2\ts-1}\right)$ for any 
$    \sigma \in \R$ and fulfill the estimates: 
 there is a constant $C_\s>0$ such that 
     \begin{align}
 \| A^{(2\ts)}\|_{\cL(H^\s, H^{\s-2\ts})} + \| \pa_t A^{(2\ts)}\|_{\cL(H^\s, H^{\s-2\ts})}  \leq C_\s \e \,; 
         \label{stima_Lam_A}\\
          \| \wt A^{(2\ts + 1)}\|_{\cL(H^\s, H^{\s-2\ts - 1})} + \| \pa_t \wt A^{(2\ts + 1)} \|_{\cL(H^\s, H^{\s-2\ts - 1})} \leq C_\s \e \,.
    \label{stima_Lam_B}
     \end{align}
 \end{lemma}
\begin{proof}
Consider first $A^{(2\ts)}$. Write 
\begin{align}
    A^{(2\ts)}=\left[ |D|^{2\ts}, \,  \cV \right]+ \left( \cV +\cV^\top\right)| D|^{2\ts}       
      \stackrel{\eqref{Vtra}}{=} \left[ |D|^{2\ts}, \,  \cV \right] - \Opbw{(\sfV_\gamma)_x} | D|^{2\ts}\,.
      \label{A_cancellata}
\end{align}
 We estimate the two terms appearing in \eqref{A_cancellata} separately.
Concerning the first one, recalling the definition of $\cV$ in \eqref{barcola} and using repeatedly  \eqref{comp01A} (with $\vr =1$) and \eqref{comp020}, we have that for any $\s \geq 0$
\begin{align}
\| [|D|^{2\ts},\ \cV]u\|_{\sigma - 2\ts } \lesssim_{\s, s}\| \sfV_\gamma\|_2 \| u\|_\s \lesssim_{s, \s}\e \|u\|_{\sigma} \quad \forall \, u \in H^\sigma\,,
\label{stima_A2s1}
\end{align}
where in the last passage we used 
\eqref{a.v.b.estimate} and the bound in 
\eqref{upp.sol.low9} (recall $\fs_0>3$). 
Similarly, concerning the second term, we apply \eqref{cont00} with $m = 0$, $s = \s - 2 \ts$, and,   recalling \eqref{upp.sol.low9}, we get
    \begin{align}
        \| \Opbw{(\sfV_\gamma)_x} |D|^{2\ts} u\|_{\s - 2\ts} \lesssim \|\sfV_\gamma\|_{2} \|u\|_{\sigma} \lesssim \e \| u\|_{\sigma}\,.
        \label{stima_A2s2}
    \end{align}
We now estimate $\pa_t A^{(2\ts)} = [|D|^{2\ts}, \pa_t \cV] - \Opbw{\pa_t (\sfV_\gamma)_x} |D|^{2\ts}.$ Arguing as to obtain \eqref{stima_A2s1}, \eqref{stima_A2s2}, one deduces that
$$
\norm{\pa_t A^{(2\ts)} u}_{\s - 2 \ts} \lesssim_{\s, \ts} \norm{\pa_t \sfV_\gamma}_2 \norm{u}_{\s}\lesssim \e \| u\|_{\sigma} \ ,   \quad \forall u \in H^\s\,,
$$
where we used \eqref{stima.der.BVa}, $\fs_0 > s_0 + 3$ and \eqref{upp.sol.low9}.

The proof for $\wt A^{(2\ts+1)}$ is analogous by noting that 
 $
     \wt A^{(2\ts+1)} 
     = 
     \left[ |D|^{2\ts+1},\  \cV \right]- | D|^{2\ts+1} \left( \cV +\cV^\top\right)
      = \left[ |D|^{2\ts+1}, \,  \cV \right] + | D|^{2\ts+1}\Opbw{(\sfV_\gamma)_x}$. 
\end{proof}

\subsection{Upside-down virial argument for  $\eta(t)$}\label{UP.eta}
\noindent
This section is devoted to the proof of Lemma \ref{lem:UD}$-(i)$.
We shall compute the first and second order derivatives of $\norm{\eta(t)}_\ts^2$ using the paralinearized equations \eqref{eq.8}.

\smallskip
\noindent{\sc Proof of \eqref{est.der.eta}:}
Using the first of \eqref{eq.8} we get
\begin{align}
     \frac{\di}{\di t} \norm{\eta(t)}_\ts^2 &=  \langle |D|^{2\ts} \pa_t \eta, \eta\rangle +  \langle |D|^{2\ts}  \eta, \pa_t\eta\rangle \,
     \notag\\
 =& 2\langle  |D|^{2\ts+1} \upomega, \eta\rangle+ 
 2\langle  |D|^{2\ts} \cB^{(-1)}   \upomega, \eta\rangle+
 \langle\underbrace{(|D|^{2\ts} \cV + \cV^\top |D|^{2\ts}) }_{=A^{(2\ts)} \mbox{ by } \eqref{lambda:A}}
 \eta,\eta\rangle+2 \langle |D|^{2\ts} f_1, \eta\rangle 
 \label{derivata_prima}
 \end{align}
We bound  each term at the right hand side of \eqref{derivata_prima}, for any $t \in [\tT -1, \tT]$. 
 The first term is bounded using the assumption \eqref{dicotomia1} and the upper bound in \eqref{upp.sol.high9}:
 \begin{align}\label{0402:1022}
     2 \left| \la |D|^{2\ts+1} \upomega, \eta \ra\right| \leq 2 \| \upomega(t)\|_{\ts+1}\| \eta(t)\|_\ts\leq 2  \mu \e^{-2\theta}\,.
 \end{align}
 We now bound the second term. 
 By \eqref{E2.est.9} and the upper bounds \eqref{dicotomia1}, \eqref{upp.sol.high9}, 
 \begin{equation}\label{cb0.est}
    2 \left| \langle |D|^{2\ts} \cB^{(-1)} \upomega, \eta \rangle\right|  \lesssim \|\cB^{(-1)} \upomega\|_{\ts + \frac 1 2} \|\eta \|_{\ts - \frac 1 2}\lesssim_\ts \e \|\upomega(t)\|_{\ts - \frac 1 2} \|\eta(t)\|_{\ts - \frac 1 2} \lesssim_\ts \e^{1-2\theta}\,. 
 \end{equation}
 To bound the third term we use  \eqref{stima_Lam_A} (with $\s\leadsto \ts$) and \eqref{dicotomia1} to get  
 \begin{align}
      \left| \langle A^{(2\ts)}\eta, \eta\rangle\right|\leq  \| A^{(2\ts)} \eta\|_{-\ts}\| \eta\|_\ts\lesssim_\ts  \e \| \eta(t)\|_{\ts}^2\lesssim_\ts  \e^{1-2\theta}\,.
      \label{est:A2s}
 \end{align}
 Finally, to bound  the last term we use \eqref{f.est.9} and \eqref{dicotomia1}
 to get 
 \begin{align}
      2&\left| \langle |D|^{2\ts} f_1, \eta\rangle\right|\leq  2 \| f_1(t)\|_{\ts}\| \eta(t)\|_\ts 
      \lesssim_\ts   \e^{1-2\theta}\,. 
      \label{est:R1}
 \end{align}
Gathering \eqref{0402:1022}--\eqref{est:R1} we proved that there is a constant $C_\ts>0$ such that 
$$
\abs{\frac{\di}{\di t} \norm{\eta(t)}_\ts^2 }\leq 2 \mu \e^{-2\theta}  + C_\ts \e^{1-2\theta}, \qquad \forall \, t \in [\tT-1, \tT] \,.
$$ 
Then \eqref{est.der.eta} follows taking $\e_0$ so small that  $ C_\ts \e\leq \mu$.
 
 \medskip
 \noindent{\sc Proof of \eqref{est.der2.eta}:}
We compute the time derivative of the identity \eqref{derivata_prima}, using also \eqref{eq.8}, 
and  further substitute
$\pa_t\upomega = - \cV^\top \upomega + (\pa_t \upomega + \cV^\top \upomega)$ and 
$\pa_t \eta = |D| \upomega + \cV \eta + (\pa_t \eta - |D|\upomega - \cV \eta)$ in the time derivative of the quadratic term $2\la |D|^{2\ts+1} \upomega, \eta \ra$, obtaining 
\begin{align}\notag
 \frac{\di^2}{\di t^2} \| \eta\|_\ts^2  = &  
2 \la |D|^{2\ts+2} \upomega, \upomega \ra  \\
 \label{L1}
 & \underbrace{+ 2 \la |D|^{2\ts} \cB^{(-1)} \pa_t \upomega, \eta \ra +
 2 \la |D|^{2\ts} \cB^{(-1)}  \upomega,  \pa_t \eta \ra 
 +
 2 \la |D|^{2\ts} (\pa_t \cB^{(-1)}) \upomega, \eta \ra }_{=: L_1(t)}\\
\label{L2}
& + \underbrace{\la A^{(2\ts)} \pa_t \eta, \eta\ra+ \la A^{(2\ts)}  \eta, \pa_t \eta\ra+ \la (\pa_t A^{(2\ts)})  \eta, \eta\ra}_{:=L_2(t)}\\ 
\label{L3}
&+ \underbrace{2  \la |D|^{2\ts} (\pa_t f_1) , \eta\ra+ 2  \la |D|^{2\ts}  f_1 , \pa_t\eta\ra}_{:= L_3(t)} \\
    \label{L4}
     &+\underbrace{2 \la \upomega, [|D|^{2\ts+1}, \cV] \eta \ra}_{=: L_4(t)}+ \underbrace{2\la |D|^{2\ts+1} (\pa_t \upomega + \cV^\top \upomega), \eta \ra}_{=: L_5(t)}\\
     &+ 
     \underbrace{2\la |D|^{2\ts+1} \upomega, (\pa_t\eta - |D| \upomega - \cV \eta) \ra}_{=: L_6(t)} \label{L6}
\end{align}
We claim that for any$\, t \in [\tT- 1, \tT]$ we have
\begin{align}
\label{estL1}
    |L_1(t)|  + |L_2(t)| + |L_3(t)| +|L_4(t)| + |L_6(t)| \lesssim_\ts  \e^{1-2\theta}  \ , 
    \qquad |L_5(t)| \leq (2+|\gamma|) \mu \e^{-2\theta}\,, 
\end{align}
implying that, for some  $C_\ts>0$, 
$$
\frac{\di^2}{\di t^2} \|\eta\|_\ts^2 \geq 2 \|\upomega\|_{\ts+1}^2 - (2+|\gamma|) \mu \e^{-2\theta} - C_\ts \e^{1-2\theta}\,.
$$
Now remark that, by assumption \eqref{dicotomia1}, the $H^{\ts+\frac12}$-norm of  $\eta(t)$ is quantitatively bounded,  thus  the $H^{\ts+1}$-norm of $\upomega(t)$ must have grown  in view of \eqref{growth.final}: quantitatively,
\be\label{dicotomia2-1}
(1-\mu) \tC_\ts^{-1} \e^{-\theta}\leq \norm{\upomega(t)}_{\ts+1} \leq \tC_\ts \e^{-\theta}   \,, \quad \forall\, t \in [\tT - 1, \tT ]
\,. 
\ee
Thus, taking $0<\mu <\min( \dfrac14, \dfrac{1}{4\tC_\ts^2(2+|\gamma|)}) $ and $\e$ so small that $C_\ts  \tC_\ts^2 \e\leq \frac14$  yields $\|\upomega\|_{\ts+1}^2 > (2+|\gamma|) \mu \e^{-2\theta} +C_\ts \e^{1-2\theta}$, proving 
\eqref{est.der2.eta}.
It remains to prove the  claim \eqref{estL1}. \\

\noindent {\em Estimate of $L_1(t)$ in \eqref{L1}:}
We start by estimating the first and second  term in $L_1(t)$.
Using twice  \eqref{cb0.est} with $\upomega \leadsto \pa_t \upomega$ respectively
$\eta \leadsto \pa_t \eta$ and 
then  the bounds \eqref{de.t.om}, \eqref{upp.sol.high9} we have, for every $t \in [\tT-1, \tT]$,
\begin{align}
\notag
    &\left| \langle |D|^{2\ts} \cB^{(-1)} \pa_t \upomega, \eta \rangle \right|
    +
     \left|\la |D|^{2\ts} \cB^{(-1)}  \upomega,  \pa_t \eta \ra \right|\\
    &\qquad  \lesssim_\ts  \e \|\eta(t)\|_{\ts-\frac 1 2}  \|\pa_t \upomega(t)\|_{\ts - \frac 1 2}
    +
    \e  \|\upomega(t)\|_{\ts -\frac 1 2} \|\pa_t \eta(t) \|_{\ts - \frac12}
    \lesssim_\ts \,  \e^{1-2\theta}\,. \label{stima:qui}
\end{align}
The third term in $L_1(t)$ is bounded by \eqref{E2.est.9} and  \eqref{upp.sol.high9} as 
\begin{align}
    \left| \la |D|^{2\ts} (\pa_t \cB^{(-1)}) \upomega, \eta\ra \right|  \leq \|  \eta\|_{\ts-\frac 1 2}  \|(\pa_t \cB^{(-1)}) \upomega\|_{\ts+\frac 1 2}  \lesssim_\ts \e \|  \eta(t)\|_{\ts-\frac 1 2} \|\upomega(t)\|_{\ts-\frac 12}\lesssim_\ts  \e^{1-2\theta}\,. \label{stima:qua}
\end{align}
Combining \eqref{stima:qui},   \eqref{stima:qua}  one gets the bound in \eqref{estL1} for $L_1(t)$.

\smallskip
\noindent {\em Estimate of $L_2(t)$ in \eqref{L2}:} 
 Using estimate \eqref{stima_Lam_A} -- which quantifies that 
 $A^{2\ts}$ is an operator of order $2\ts$ and not $2\ts+1$--
 and arguing as to estimate $L_1(t)$, we have
\begin{equation*}
   {\small \left| \la A^{(2\ts)} \pa_t \eta, \eta \ra \right| +
     \left| \la A^{(2\ts)} \eta, \pa_t \eta \ra \right|
    + \left| \la (\pa_t A^{(2\ts)}) \eta, \eta \ra \right|}
{\lesssim_\ts}  \e \|\eta(t)\|_{\ts+\frac 1 2} \|\pa_t \eta(t)\|_{\ts - \frac 1 2} 
    + \e \|\eta(t)\|_{\ts + \frac 1 2} \| \eta(t)\|_{\ts -\frac 1 2} \!\!\!\!\!\!\!\!\!
\stackrel{\eqref{de.t.om}, \eqref{upp.sol.high9}}{\lesssim_\ts} 
\!\!\!\!\!\!
\e^{1-2\theta}\,.
\end{equation*}

\noindent {\em Estimate of $L_3(t)$ in \eqref{L3}:} By \eqref{f.est.9}, \eqref{de.t.om}, and \eqref{upp.sol.high9}
\begin{equation*}
   \left| \la |D|^{2\ts} (\pa_t f_1), \eta\ra\right|
    +
      \left| \la |D|^{2\ts} f_1, \pa_t \eta\ra \right|
    \leq \|\eta(t)\|_{\ts+ \frac 1 2} \|\pa_t f_1(t) \|_{\ts - \frac 1 2} 
  +  \|  f_1(t)\|_{\ts + \frac 1 2} \|\pa_t \eta(t)\|_{\ts-\frac 1 2}  \lesssim_\ts  \e^{1-2 \theta}\,.
\end{equation*}

\noindent {\em Estimate of $L_4(t)$ in \eqref{L4}:} 
Arguing as to obtain \eqref{stima_A2s1} with $2s \leadsto 2\ts +1$
and using the  bounds \eqref{upp.sol.high9}, \eqref{dicotomia1}, one has
\begin{align*}
    \left|\la \upomega, [|D|^{2\ts+1}, \cV] \eta \ra \right| \leq \| \upomega\|_{\ts + 1} \| [|D|^{2\ts + 1},\, \cV] \eta \|_{-\ts - 1} \lesssim_\ts \e \| \upomega(t)\|_{\ts +1 } \|\eta(t)\|_{\ts} \lesssim_\ts \e^{1-2\theta}\,.
\end{align*}
\smallskip
\noindent {\em Estimate of $L_5(t)$ in \eqref{L4}:}  Recalling equation \eqref{eq.8} and using estimates \eqref{T.est.9}, \eqref{E2.est.9}, \eqref{f.est.9}, \eqref{dicotomia1} and \eqref{upp.sol.high9}  we have
$$
\begin{aligned}
    \Big| \Big\langle |D|^{2\ts + 1} & (\pa_t \upomega + \cV^\top \upomega), \eta \Big \rangle \Big|  \leq \|\eta\|_{\ts + \frac 1 2} \|\pa_t \upomega + \cV^\top \upomega\|_{\ts + \frac 1 2}= \|\eta\|_{ \ts + \frac 12} \| - \eta + \gamma \Hilb\upomega+ \cB^{(-2)} \upomega + \cT \eta + f_2\|_{\ts + \frac 1 2}\\
    &\leq \|\eta(t)\|_{\ts + \frac12}^2 + |\gamma| \|\eta(t)\|_{\ts + \frac12}\|\upomega(t)\|_{\ts + \frac 1 2} + \e  \|\eta(t)\|_{\ts + \frac12}\left(\norm{\upomega(t)}_{\ts - \frac 1 2}  +  \|\eta(t)\|_{\ts + \frac12} 
    +\e^{-\theta}
    \right)\\
     &\leq (1+|\gamma|)  \mu \e^{-2\theta}+ C_\ts\mu \e^{1-2\theta}\,,
\end{aligned}
$$
which implies the second bound in \eqref{estL1} provided $ \e C_\ts \leq 1$.

\smallskip
\noindent {\em Estimate of $L_6(t)$ in \eqref{L6}:}  Arguing as above and using \eqref{eq.8} and estimates  \eqref{E2.est.9}, \eqref{f.est.9} and \eqref{upp.sol.high9}, we obtain
\begin{align*}
\nonumber
    \left| \la |D|^{2\ts + 1} \upomega, \pa_t \eta - |D| \upomega - \cV \eta \ra \right| &
     \leq \| |D|^{2\ts + 1} \upomega\|_{-\ts} \|\pa_t \eta - |D| \upomega - \cV \eta \|_{\ts}
    \leq \|  \upomega\|_{\ts+1} \|\cB^{(-1)} \upomega + f_1\|_{\ts}\\
    &\lesssim_\ts 
    \e 
\|\upomega(t)\|_{\ts+1} \left( \|\upomega(t)\|_{\ts-1} + \e^{-\theta}\right) \lesssim \e^{1-2\theta}\,.
\end{align*}
The claim \eqref{estL1} is proved.
\qed
\subsection{Upside-down virial argument for  $\upomega(t)$}\label{UP.omega}
\noindent
This section is devoted to the proof of Lemma \ref{lem:UD}$-(ii)$.
We shall compute the first and second order derivatives of $\norm{\upomega(t)}_\ts^2$ exploiting the paralinearized equations \eqref{eq.8}.

\smallskip
\noindent{\sc Proof of \eqref{stima.der.om}}: using the second of  \eqref{eq.8} we have
\be\label{derivata_prima_om}
\begin{aligned} 
     \frac{\di}{\di t} \norm{\upomega(t)}_{\ts+\frac12}^2 &= 
-2\langle  |D|^{2\ts+1}\eta , \upomega\rangle +  
 2\langle  |D|^{2\ts + 1}  \cB^{(-2)}   \upomega, \upomega \rangle+ 2 \langle |D|^{2\ts+1} \cT \eta, \upomega \rangle \\
 & \ \ +
 \langle \underbrace{\left(- |D|^{2\ts+1}\cV^\top - \cV |D|^{2\ts+1}\right)  }_{=: \wt A^{(2\ts+1)} \mbox{ by } \eqref{lambda:A}}\upomega,\upomega\rangle +2 \langle |D|^{2\ts + 1} f_2, \upomega\rangle 
  \,.
 \end{aligned}
 \ee
 We bound each term at the right hand side of \eqref{derivata_prima_om},  for any $t \in [\tT -1, \tT]$. 
The first term is estimated by  \eqref{dicotomia3} 
and \eqref{upp.sol.high9}, getting
\begin{equation}
\label{stima:homer}
    2\left| \la |D|^{2\ts+1} \eta, \upomega \ra\right| \leq 2 \|\upomega(t)\|_{\ts+1} \| \eta(t) \|_{\ts}  \leq 2\mu \e^{-2\theta}\,.
\end{equation}
The second term of \eqref{derivata_prima_om} is estimated using 
\eqref{E2.est.9} and the
upper bound \eqref{upp.sol.high9}, 
which give
\begin{align}
\label{stima:marge}
    2 \left| \la |D|^{2\ts + 1} \cB^{(-2)} \upomega, \upomega \ra\right| \leq 2\|\upomega\|_{\ts +\frac 1 2} \|\cB^{(-2)} \upomega \|_{\ts + \frac 1 2} \lesssim_\ts \e \|\upomega(t)\|^2_{\ts + \frac 1 2} \lesssim_\ts  \e^{1-2\theta}\,.
\end{align}
The third term is estimated by  \eqref{T.est.9} and \eqref{upp.sol.high9}. One gets
\begin{align}
    \label{stima:bart}
    2 \left| \la |D|^{2\ts + 1} \cT \eta, \upomega \ra \right|  \leq 2 \|\upomega\|_{\ts + 1} \|\cT \eta\|_{\ts} \lesssim_{\ts} \e \|\upomega(t)\|_{\ts + 1} \|\eta(t)\|_{\ts} \lesssim_{\ts} \e^{1-2\theta}\,.
\end{align}
The fourth term is estimated using \eqref{stima_Lam_B} and \eqref{upp.sol.high9}, getting
\begin{align}
    \label{stima:lisa}
    \left| \la \wt{A}^{(2\ts + 1)} \upomega, \upomega \ra \right| \leq \|\upomega\|_{\ts + \frac 1 2} \| \wt{A}^{(2\ts+1)} \upomega\|_{-\ts-\frac 1 2} \lesssim_{\ts}  \e \|\upomega(t)\|_{\ts + \frac 1 2}^2 \lesssim \e^{1-2\theta}\,,
\end{align}
and  finally the last term is estimated by \eqref{f.est.9} and \eqref{upp.sol.high9}, getting
\begin{align}
    \label{stima:maggie}
    2 \left| \la |D|^{2\ts +1} f_2, \upomega \ra \right|
     \leq 2 \|\upomega(t)\|_{\ts + 1} \|f_2(t) \|_{\ts} \lesssim_\ts \e^{1-2\theta}\,.
\end{align}
Combining estimates \eqref{stima:homer}--\eqref{stima:maggie}, we proved that there is a constant $C_\ts>0$ such that 
$$
\abs{\frac{\di}{\di t} \norm{\upomega(t)}_{\ts+\frac12}^2} \leq 2 \mu \e^{-2\theta}  + C_\ts \e^{1-2\theta} \ .
$$ 
Then \eqref{stima.der.om} follows by taking $\e_0$ small enough so that $C_\ts \e \leq \mu$.

\smallskip
\noindent{\sc Proof of \eqref{stima.der2.om}}:
We compute the time derivative of the identity \eqref{derivata_prima_om}, using \eqref{eq.8} 
and  further substituting
$\pa_t\upomega = -\eta - \cV^\top \upomega + (\pa_t \upomega +\eta + \cV^\top \upomega)$ and 
$\pa_t \eta =   \cV \eta + (\pa_t \eta  - \cV \eta)$ in the time derivative of the quadratic term $-2\la |D|^{2\ts+1} \pa_t \eta, \upomega \ra - 2 \la |D|^{2\ts+1} \eta, \pa_t \upomega \ra$, obtaining 
\begin{align}\notag
    \frac{\di^2}{\di t^2} \|\upomega\|_{\ts+\frac 1 2}^2
&= 2 \la |D|^{2\ts + 1} \eta, \eta \ra \\
\label{M1}
& + \underbrace{2 \la |D|^{2\ts + 1} \cB^{(-2)} \pa_t \upomega, \upomega\ra + 2 \la |D|^{2\ts+1} \cB^{(-2)} \upomega, \pa_t \upomega \ra + 2 \la |D|^{2\ts+1} (\pa_t \cB^{(-2)} ) \upomega, \upomega \ra}_{:= M_1(t)} \\
\label{M2}
& + \underbrace{2 \la  |D|^{2\ts+1} \cT (\pa_t \eta - \cV \eta), \upomega\ra + 2 \la |D|^{2\ts+1} (\pa_t \cT) \eta, \upomega \ra + 2 \la |D|^{2\ts+1} \cT \eta, \pa_t \upomega  + \cV^\top \upomega \ra }_{:= M_2(t)}\\
\label{M3}
& + \underbrace{2 \la [|D|^{2\ts+1} \cT,\, \cV] \eta, \upomega \ra }_{:= M_3(t)}\\
\label{M4}
& + \underbrace{\la \wt{A}^{(2\ts+1)} \pa_t \upomega,  \upomega \ra  + \la \wt{A}^{(2\ts+1)} \upomega, \pa_t \upomega \ra + \la (\pa_t \wt{A}^{(2\ts+1)}) \upomega, \upomega \ra }_{:= M_4(t)}\\
\label{M56}
& + \underbrace{2 \la |D|^{2\ts+1} \pa_t f_2, \upomega \ra + 2 \la |D|^{2\ts+1} f_2, \pa_t \upomega \ra }_{:= M_5(t)} - \underbrace{2 \la |D|^{2\ts+1} (\pa_t \eta - \cV \eta), \upomega \ra}_{:= M_6(t)}\\
\label{M78}
& -\underbrace{2 \la |D|^{2\ts+1} \eta, \pa_t \upomega + \eta + \cV^\top \upomega \ra}_{:= M_7(t)} + \underbrace{ 2 \la [\cV,\ |D|^{2\ts+1}] \eta, \upomega \ra}_{:= -L_4(t) \mbox{ \footnotesize\  by\ }  \eqref{L4}}\,.
\end{align}
We claim that for any $t \in [\tT-1, \tT]$ we have 
\be
\label{estMs}
\begin{aligned}
&|M_1(t)|+|M_2(t)|+|M_3(t)|+|M_4(t)|+|M_5(t)| \lesssim_\ts  \e^{1-2\theta}  \,, 
    \\
&|M_6(t)|  + |M_7(t)| \leq( 2|\gamma| + 2) \mu \e^{-2\theta}  \ .
\end{aligned}
\ee
Using also the estimate $|L_4(t)| \lesssim \e^{1-2\theta}$ (cf. \eqref{estL1}) we get, for some  $C_\ts>0$, 
$$
\frac{\di^2}{\di t^2} \|\upomega(t)\|_{\ts+\frac12}^2 \geq 2 \|\eta(t)\|_{\ts+\frac12}^2 - (2+2|\gamma|) \mu \e^{-2\theta} - C_\ts \e^{1-2\theta}\,.
$$
As in the previous case,  by assumption \eqref{dicotomia3}, the $H^{\ts+1}$-norm of  $\upomega(t)$ is quantitatively bounded,  thus  the $H^{\ts+\frac12}$-norm of $\eta(t)$ must have grown  in view of \eqref{growth.final}: quantitatively,
\be\label{dicotomia4}
(1-\mu) \tC_\ts^{-1} \e^{-\theta}\leq \norm{\eta(t)}_{\ts+\frac12} \leq \tC_\ts \e^{-\theta}   \ , \quad \forall\, t \in [\tT - 1, \tT ]\,. 
\ee
Thus, taking 
$0<\mu <\min( \dfrac14, \dfrac{1}{4\tC_\ts^2(2+2|\gamma|)}) $ and $\e$ so small that $C_\ts  \tC_\ts^2 \e\leq \frac14$  one has 
$\|\eta\|_{\ts+\frac12}^2 > (2+2|\gamma|) \mu \e^{-2\theta} +C_\ts \e^{1-2\theta}$, proving 
\eqref{stima.der2.om}. It remains to
prove the  claim \eqref{estMs}. \\

\noindent {\em Estimate of $M_1(t)$ in \eqref{M1}:} 
by \eqref{E2.est.9} and \eqref{de.t.om} and the upper bounds 
in  \eqref{upp.sol.high9} we have
\begin{align}
\nonumber
    \big| \langle |D|^{2\ts + 1} \cB^{(-2)} \pa_t \upomega, \upomega \rangle &+ \langle |D|^{2\ts + 1} \cB^{(-2)} \upomega, \pa_t \upomega \rangle \big| \leq \|\upomega\|_{\ts + 1} \|\cB^{(-2)} \pa_t \upomega\|_{\ts} + \|\pa_t \upomega\|_{\ts} \| \cB^{(-2)} \upomega \|_{\ts + 1}\\
    &\lesssim  \e \left(\|\upomega(t)\|_{\ts + 1}   \norm{\pa_t \upomega(t)}_{\ts - 1} + \norm{\pa_t \upomega(t)}_{\ts} \norm{\upomega(t)}_{\ts}\right) \lesssim_{\ts} \e^{1-2\theta}\,.
    \label{stima:fred}
\end{align}
Moreover, using \eqref{E2.est.9} and \eqref{upp.sol.high9}
\begin{align}
    \left| \la |D|^{2\ts + 1} (\pa_t \cB^{(-2)} ) \upomega, \upomega \ra\right| &  \leq \|\upomega\|_{\ts+ \frac 1 2} \| (\pa_t \cB^{(-2)}) \upomega\|_{\ts +\frac 1 2} 
    \lesssim_\ts \e \|\upomega(t)\|_{\ts + \frac 1 2} \| \upomega(t)\|_{\ts -\frac 1 2} \lesssim_\ts  \e^{1-2\theta}\,.
    \label{stima:flintstone}
\end{align}
Combining estimates \eqref{stima:fred} and \eqref{stima:flintstone}, we get the bound \eqref{estMs} for $M_1(t)$.

\smallskip
\noindent {\em Estimate for $M_2(t)$ in \eqref{M2}:}
Arguing as above, by using the estimates
\eqref{T.est.9}, \eqref{E2.est.9}, \eqref{f.est.9} the upper bounds  \eqref{upp.sol.high9}, and the equations \eqref{eq.8},
we get
\begin{align}
\nonumber
    \left|\la |D|^{2\ts + 1} \cT (\pa_t \eta - \cV \eta), \upomega \ra \right| & \leq \|\upomega\|_{\ts + 1} \|\cT (\pa_t \eta - \cV \eta)\|_{\ts} 
\lesssim_\ts \e \|\upomega\|_{\ts+1} \|\pa_t \eta- \cV \eta\|_\ts\\
\notag
&\lesssim_\ts \e \|\upomega\|_{\ts+1} \| |D| \upomega + \cB^{(-1)} \upomega + f_1 \|_{\ts} \\
    &\lesssim_\ts \e \|\upomega(t)\|_{\ts+1} \left( \| \upomega(t)\|_{\ts+1} + \e \|\upomega(t)\|_{\ts-1} + \e^{1-\theta}\right) \lesssim_\ts \e^{1-2\theta}\,,
    \label{stima:leonardo}
\end{align}
as well as
\begin{align}
\nonumber
    \left| \la |D|^{2\ts + 1} \cT \eta, \pa_t \upomega + \cV^\top \upomega \ra \right| &\leq \|\pa_t \upomega + \cV^\top \upomega \|_{ \ts + \frac 1 2} \| \cT \eta\|_{\ts + \frac 1 2} \\
    \notag
    &\lesssim_\ts \e \|\eta\|_{\ts + \frac 1 2} \| - \eta + \gamma \Hilb \upomega+ \cB^{(-2)} \upomega + \cT \eta + f_2 \|_{\ts + \frac 1 2} \\
    &\lesssim_\ts \e \|\eta(t)\|_{\ts + \frac 1 2} \left( \| \eta(t)\|_{\ts + \frac 1 2} + \|\upomega(t)\|_{\ts + \frac 1 2} + \e^{1-\theta} \right) \lesssim_\ts \e^{1-2\theta}\,,
    \label{stima:donatello}
\end{align}
and finally 
\begin{align}
\left| \la |D|^{2\ts +1} (\pa_t \cT) \eta, \upomega \ra \right| \leq
\|\upomega\|_{ \ts + 1} \| (\pa_t \cT ) \eta\|_{\ts} \lesssim_\ts \e \|\upomega(t)\|_{\ts+1} \|\eta(t)\|_{\ts} \lesssim_\ts \e^{1-2\theta}\,.
\label{stima:raffaello}
\end{align}
The combination of \eqref{stima:leonardo}, \eqref{stima:donatello} and \eqref{stima:raffaello}
gives the upper bound for $M_2(t)$ in \eqref{estMs}.

\smallskip
\noindent{\em Estimate of $M_3(t)$ in \eqref{M3}:} by using \eqref{comp01A}, together with estimates \eqref{poi.est}, \eqref{E1.est.9} and \eqref{T.est.9} and the upper bounds \eqref{upp.sol.high9}, one has
\begin{equation*}
    \left| \la [|D|^{2\ts + 1} \cT,\ \cV] \eta, \upomega \ra \right| 
    \leq \|\upomega\|_{\ts + 1} \|[|D|^{2\ts + 1} \cT, \cV] \eta\|_{-\ts-1} \lesssim_\ts \e^2 \|\upomega(t)\|_{\ts+ 1} \|\eta(t)\|_{\ts} \lesssim_\ts \e^{2-2\theta}\,.    
\end{equation*}

\smallskip
\noindent {\em Estimate of $M_4(t)$ in \eqref{M4}:} using estimates \eqref{stima_Lam_B}, \eqref{upp.sol.high9} and \eqref{de.t.om}, one has
\begin{align}
\nonumber
    \Big| \Big \langle \wt{A}^{(2\ts +1 )} &\pa_t \upomega,  \upomega \Big \rangle \Big| + \left| \la  \wt{A}^{(2\ts +1 )} \upomega, \pa_t \upomega \ra \right|+ \left| \la (\pa_t \wt{A}^{(2\ts+1)} \upomega, \upomega \ra \right| \\
    \notag
    &\leq \| \upomega\|_{\ts+1} \|\wt{A}^{(2\ts+1)} \pa_t \upomega\|_{-\ts-1} + \|\pa_t \upomega\|_{\ts} \| \wt{A}^{(2\ts+1)} \upomega\|_{-\ts}+\|\upomega\|_{\ts+ 1} \| \pa_t \wt{A}^{(2\ts +1 )} \upomega\|_{-\ts-1} \\
    \notag
    &\lesssim_\ts \e (\|\upomega(t)\|_{\ts + 1} \|\pa_t \upomega(t)\|_{\ts} +\|\upomega(t)\|_{ \ts + 1} \|\upomega(t)\|_{\ts})\lesssim_\ts \e^{1-2\theta}\,.
\end{align}

\smallskip
\noindent {\em Estimate of $M_5(t)$ in \eqref{M56}:} using \eqref{f.est.9} and the upper bounds in  \eqref{upp.sol.high9} and \eqref{de.t.om}, we have 
\begin{equation} 
    \left| \la|D|^{2\ts +1} \pa_t f_2 , \upomega \ra \right|
    +
     \left| \la |D|^{2\ts +1} f_2, \pa_t \upomega \ra \right| 
    \leq \|\upomega(t)\|_{\ts + 1} \|\pa_t f_2(t)\|_{\ts}
    +
    \|\pa_t \upomega(t)\|_{\ts} \|f_2(t) \|_{\ts + 1}
 \lesssim_\ts \e^{1-2\theta}\,.
 \notag
\end{equation}

\smallskip
\noindent {\em Estimate of $M_6(t)$ in \eqref{M56}:} Using the equations for $\pa_t \eta$ in \eqref{eq.8}, and the upper bounds \eqref{dicotomia3}, \eqref{E2.est.9}, \eqref{f.est.9} and \eqref{upp.sol.high9}, we get
\begin{align}
\nonumber
    \left| 2 \la |D|^{2\ts + 1} (\pa_t \eta - \cV \eta), \upomega \ra \right| &\leq 2 \norm{\upomega}_{\ts + 1} \norm{\pa_t \eta - \cV \eta}_{\ts}  = 2 \|\upomega\|_{\ts+1} \| |D| \upomega + \cB^{(-1)} \upomega + f_1 \|_{\ts} \nonumber\\
    &\leq 2\|\upomega\|_{\ts+1}^2 + C_{\ts} \|\upomega\|_{\ts+1} \left( \e \|\upomega\|_{\ts-1} + \e^{1-\theta} \right) \leq 2 \mu^{2} \tC_\ts^{-2} \e^{-2\theta} + \tC_\ts C_\ts \e^{1-2\theta}\,,
\notag
\end{align}
which taking $\e$ so small that $2 \tC_{\ts} C_{\ts} \e \leq \mu \leq \frac 1 2$ and recalling $\tC_\ts > 1$, gives $M_6(t) \leq  \mu \e^{-2\theta}$\,.

\smallskip
\noindent {\em Estimate of $M_7(t)$ in \eqref{M78}:} Using the equations for $\pa_t \upomega$ in \eqref{eq.8} and estimates \eqref{E2.est.9}, \eqref{T.est.9} and \eqref{f.est.9} together with the upper bounds \eqref{dicotomia3} and \eqref{upp.sol.high9}, we have
\begin{align}
\nonumber
    \left| 2 \la |D|^{2\ts + 1} \eta, \pa_t \upomega + \eta + \cV^\top \upomega \ra \right|
    &  \leq 2 \|\pa_t \upomega + \eta + \cV^\top \upomega \|_{\ts + \frac 1 2} \| \eta\|_{\ts + \frac 1 2} = 2 \|\gamma \Hilb \upomega + \cB^{(-2)} \upomega + \cT \eta + f_2\|_{\ts + \frac 1 2} \|\eta\|_{\ts + \frac 1 2} \nonumber
    \\
    \notag
    &\leq 
2 |\gamma| \|\eta(t)\|_{\ts + \frac 1 2} \|\upomega(t)\|_{\ts + \frac 1 2}+
    C_\ts 
   \e  \|\eta(t)\|_{\ts + \frac 1 2}  \left(  \|\upomega(t)\|_{\ts-\frac12}+  \|\eta(t)\|_{\ts + \frac 1 2} + \e^{-\theta} \right) \\
    & \leq 2 \mu |\gamma|  \e^{-2\theta} + \tC_{\ts} C_\ts \e^{1-2\theta}\,.
\notag
\end{align}
Then by taking $\e$ small enough such that $\tC_{\ts} C_{\ts}\e <\mu$, one obtains $M_7(t) \leq (2 |\gamma| + 1) \mu \e^{-2\theta}$.

All the claimed estimates in  \eqref{estMs} are  proved.
\qed

\subsection{Proof of \Cref{thm:mainimprove}}
Let $\gamma<0$, $\gamma^2 \in \Q$, $\delta \in (0,1)$ and $K \geq 1$ be given. Fix $\ts, \ts_0, \theta$ as explained in the proof of \Cref{thm:main} at page \pageref{proof:thm1}.
\Cref{lem:partial_growth} gives $\e_0>0$ such that, for any $\e \in (0, \e_0)$, there exists a solution 
 $(\eta, \psi, \sfV, \sfB) \in C^0([0, \tT];  X^\ts)$ of
 \eqref{eq:etapsi} fulfilling \eqref{growth.final}--\eqref{upp.sol.high9}.
 
We first show that there exists $T_1 \in [\tT-1, \tT]$ such that 
 \be\label{ass:eta.gr}
\norm{\eta(T_1)}_{\ts+\frac12} \geq  \mu \, \tC_\ts^{-1} \e^{-\theta} \ .
\ee
We  argue  by contradiction. Assume that for any  $T_1
\in [\tT-1, \tT]$,  
\eqref{ass:eta.gr} is violated, namely  \eqref{dicotomia1} holds. 
By the  upside-down virial  \Cref{lem:UD} $(i)$ we deduce that the second derivative of $t \mapsto \|\eta(t)\|_{\ts }^2$ must have grown: 
 \be\label{1711:1629}
 \frac{\di^2}{\di t^2} \| \eta(t)\|_\ts^2 \geq \frac1 2 \tC_\ts^{-2}\e^{-2\theta} \,, 
 \quad \forall\, t \in [\tT-1, \tT] \ . 
 \ee
We deduce  a lower bound for the norm $\|\eta(t)\|_{\ts }^2$.
Indeed, by the Taylor expansion of $\|\eta(t)\|_{\ts }^2$ at $t = \tT -1$ and using \eqref{est.der.eta}, \eqref{1711:1629} we obtain
    \begin{equation}
    \begin{aligned}\label{eta.really.grow}
     \|\eta(\tT)\|_{\ts }^2  & \geq \|\eta(\tT-1)\|_{\ts }^2 +\frac{\di}{\di t} \|\eta(\tT-1)\|_{\ts }^2  + \frac 1 2  \inf_{\tau \in [\tT - 1, \tT] } \frac{\di^2}{\di t^2}\|\eta(\tau)\|_{\ts }^2 \\
     & \geq \left(\frac 1 4 \tC_\ts^{-2}  - 4 \mu  \right) \e^{-2\theta} \geq \frac 1 8 \tC_\ts^{-2} \e^{-2\theta}\,,
     \end{aligned}
    \end{equation}
    taking $0<\mu <      \dfrac{1}{32 \tC_\ts^2 } $, 
   which recalling also $\tC_\ts > 1$ contradicts \eqref{dicotomia1}, absurd.  
   Hence \eqref{ass:eta.gr} must be true.

\smallskip

Similarly we show that there is $T_2 \in [\tT-1, \tT]$ such that 
   \be\label{ass:omega.gr}
\norm{\upomega(T_2)}_{\ts+1} \geq \mu \, \tC_\ts^{-1} \e^{-\theta} \ .
\ee
Again we argue by contradiction.  Assume that for any  $t
\in [\tT-1, \tT]$,  
\eqref{ass:omega.gr} is violated, namely  \eqref{dicotomia3} holds. 
We deduce, by \Cref{lem:UD} $(ii)$,   that
\be\label{1711:1633}
 \frac{\di^2}{\di t^2} \| \upomega(t)\|_{\ts+\frac 1 2}^2
 \geq \frac12 \tC_\ts^{-2}\e^{-2\theta}   \ , \quad \forall\, t \in [\tT-1, \tT] \ . 
 \ee
     We Taylor expand 
      $\norm{\upomega(t)}_{\ts+\frac12}^2$ at $t = \tT -1$ and argue as in \eqref{eta.really.grow}, using this time the bounds in 
    \eqref{stima.der.om},  \eqref{dicotomia3} 
    and \eqref{1711:1633}
    to deduce the lower bound
    $$
    \norm{\upomega(\tT)}_{\ts+\frac12}^2 \geq \frac18 \tC_\ts^{-2} \e^{-2\theta}
    $$
    that contradicts \eqref{dicotomia3} for $\mu$ small enough, absurd. Then \eqref{ass:omega.gr} must hold true.
    
Finally, we  deduce that 
  \be\label{ass:VB.gr}
\norm{\sfV(T_2)}_{\ts}   \geq \frac12 \mu \, \tC_\ts^{-1} \e^{-\theta} \ , \quad 
\norm{\sfB(T_2)}_{\ts}   \geq \frac12 \mu \, \tC_\ts^{-1} \e^{-\theta} \ . 
\ee
Indeed,  by  \eqref{stima_V}, \eqref{ass:omega.gr}  and \eqref{upp.sol.high}, \eqref{3001:1642}, 
    $$
    \norm{\sfV(T_2)}_\ts 
    \geq \norm{\upomega(T_2)}_{\ts+1} - \norm{\sfV(T_2) - \upomega_x(T_2)}_\ts \geq \mu \tC_{\ts}^{-1} \e^{-\theta} - C_{\ts} \e^{1-\theta} \geq \frac12 \mu \, \tC_\ts^{-1} \e^{-\theta}
    $$
    provided $\e$ is sufficiently small.
    The lower bound for $\norm{\sfB(T_2)}_{\ts}$ follows analogously exploiting this time the estimate \eqref{stima_B} and \eqref{upp.sol.low9}--\eqref{upp.sol.high9}.
Again eventually shrinking $\e$, \eqref{ass:eta.gr} and \eqref{ass:VB.gr} give the claim in \eqref{explosion.improved}. 
The low Sobolev norm bound in \eqref{low.norm.intro0} is proved in \Cref{thm:main}.

\qed

\appendix 

\section{On the paralinearization of the Dirichlet-Neumann operator}\label{subsec:preliminary}
We  state a paralinearization formula for the Dirichlet--Neumann operator
$G(\eta)$ defined in~\eqref{DN}, including the pluri-homogeneous expansion of
the symbol and the smoothing remainders, as established in \cite[Chapter
7]{BD}. In addition to the paralinearization in \cite{BD}, we show that each homogeneous component of the
implicit negative-order symbol is real-valued, a feature required for the bootstrap
argument in \Cref{lem:long.boot}. This last property follows \emph{a
posteriori} by adapting the general strategy of \cite[Lemma 3.20]{BMM2},
and leveraging also the self-adjointness of the Dirichlet--Neumann operator; see also \cite[Section 5.3.2]{FMT}.
More precisely, we prove the following. 
\begin{lemma}\label{lemma:para.bella}
Let $N \in \N$ and $\vr>0$. There exists $ r >0$ such that the following hold.
\begin{enumerate}[label=(\roman*)]
    \item  {\bf  Paralinearization:} One has the paralinearization formula
\begin{align}
    G(\eta)\psi= |D|\upomega+ \Opbw{- \im \sfV(\eta,\psi;x)\, \xi - \tfrac{1}{2} \left(\sfV(\eta,\psi;x)\right)_x} \eta+ \Opbw{b_{\geq 1}^{(-1)}(\eta;x,\xi)}\upomega+ R(\eta)\psi\,,
    \label{para:DN}
\end{align}
where 
$\upomega$ is the good unknown of Alinhac  in \eqref{GU},
$\sfV(\eta,\psi; \cdot)$ is defined in \eqref{def:V},
$R(\eta)$ is a real smoothing remainder in $\Sigma\cR^{-\vr}_1[r, N]$ and $b_{\geq 1}^{(-1)}$ is a symbol in $\Sigma \Gamma^{-1}_1[r, N]$ expanding (cf. \eqref{symbols}) as 
\begin{align}
    b_{\geq 1}^{(-1)}= \sum_{p=1}^{N-1} \underbrace{b_{p}^{(-1)}}_{\in \wt \Gamma^{-1}_p}+ \underbrace{b_{\geq N}^{(-1)}}_{\in \Gamma_{\geq N}^{-1}[r]}, \qquad \begin{cases}
    b_{p}^{(-1)}=\overline{b_{p}^{(-1)}},\\
         b_{p}^{(-1)}=\big(b_{p}^{(-1)}\big)^\vee
    \end{cases}   \quad\text{for any } p=1,\dots, N-1.
    \label{real:b-1}
\end{align}
Moreover, one has that 
\be\label{G.mappa}
G(\eta) = |D| + M_{\geq 1}(\eta) \ , \quad M_{\geq 1}(\eta) \in \Sigma \cM_{1}^1[r,N] \ . 
\ee
\vspace{-0.5cm}
\item {\bf  Time derivative control:}
There exist  $ s_0 > 0$ such that for any   $\sigma \geq s_0$,  if 
 $(\eta, \psi, \sfV, \sfB) \in B_{X^{s_0}}(I;r) \cap X^{\s}$ is a solution of  \eqref{eq:etapsi},
 then for $t \in I$ and any $v \in C^1(I; H^{\s+\frac12})\cap B_{s_0}(r)$
\begin{equation}\label{pa_t:smooth}
\begin{gathered}
 |b^{(-1)}_{\geq 1}|_{-1, {W^{\s+ \frac12 -s_0,\infty}},M} + | \pa_t b^{(-1)}_{\geq 1}|_{-1, W^{\s + \frac 1 2-s_0,\infty},M}\lesssim_{\s, M} \| (\eta,\psi,\sfV,\sfB)\|_{ X^{\s}} \,, \quad \forall M \in \N \ , \\
    \| \pa_t (R(\eta) v) \|_{\sigma + \vr - \frac12}  \lesssim_{\s} 
    \| (\eta,\psi,\sfV,\sfB)\|_{X^{s_0}}
    \left( \| v\|_{\s+\frac12}  +   \|\pa_t v \|_{\s-\frac 1 2} \right)\\
    \qquad \qquad  + 
    \| (\eta,\psi,\sfV,\sfB)\|_{X^{\s}}
    \left( \| v\|_{s_0}  +   \|\pa_t v \|_{s_0-1} \right)
\end{gathered}
\end{equation}
In particular, when $ v=\psi$, one has 
\begin{align}\label{pa_t:smooth2}
        \| \pa_t (R(\eta) \psi) \|_{\sigma + \vr - \frac12} \lesssim_{\s} \| (\eta,\psi,\sfV,\sfB)\|_{X^{s_0}}\,  \| (\eta,\psi,\sfV,\sfB)\|_{X^{\s}}.
\end{align}
\end{enumerate}
\end{lemma}
To prove Lemma \ref{lemma:para.bella}, we shall use the following preliminary result.
\begin{lemma}\label{lem:tutto.X}
    There exist $r, s_0>0$ such that, if $\s \geq  s_0 $ and $(\eta, \psi, \sfV, \sfB)  \in C^0(I; X^{\s })$ is a solution   of \eqref{eq:etapsi} with 
    $\sup_{t \in I} \norm{(\eta(t), \psi(t), \sfV(t), \sfB(t))}_{X^{s_0}} \leq r$, 
 then 
\begin{equation}\label{pa_t:etapsi}
        \|\pa_t \eta\|_{\s-\frac12} + \|\pa_t \psi\|_{\s-\frac12} \lesssim_{\s} \|(\eta, \psi, \sfV, \sfB)\|_{X^{\s }}\,.
    \end{equation}
\end{lemma}
\begin{proof}
    We start with estimating $\pa_t \eta$. Using \eqref{eq:etapsi} and \eqref{espr_DN}, one has $\pa_t \eta = \sfB - \sfV \eta_x + \gamma \eta \eta_x$. Therefore, by standard tame estimates and recalling \eqref{norm:Xs}, one has
    \begin{align*}
    \| \pa_t \eta\|_{\s-\frac12}  & \lesssim_{\s} \| \sfB\|_{\s-\frac12} + \| \sfV\|_{\s-\frac12} \|\eta\|_{W^{1,\infty}} + \|\sfV\|_{L^\infty} \|\eta\|_{\s+\frac12} + \gamma \|\eta\|_{W^{1,\infty}} \|\eta\|_{\s+\frac12} \\
    & \lesssim_{\s} (1 + \|(\eta, \psi, \sfV, \sfB)\|_{X^{2}}
    ) \|(\eta, \psi, \sfV, \sfB)\|_{X^{\s }} \lesssim_{\s} \|(\eta, \psi, \sfV, \sfB)\|_{X^{\s}}\,.
    \end{align*}
    Then we estimate $\pa_t \psi$. We first note that, by \eqref{eq:etapsi} and recalling \eqref{form-of-B} and \eqref{espr_DN}, one has 
    \begin{align}
        \pa_t \psi= -\eta-\frac12 \psi_x^2 + \frac12\sfB^2 (1+\eta_x^2)+ \gamma \eta \psi_x+ \gamma \pa_x^{-1} (\sfB-\sfV \eta_x)\,.
    \end{align}
    Proceeding as done for the estimate on $\partial_t \eta$,  we deduce 
    \begin{align*}
        \|\pa_t \psi\|_{\s-\frac12}\lesssim_\s &  \|\eta\|_{\s-\frac12} +\|\psi\|_{W^{1,\infty}}\| \psi\|_{\s+\frac12} + \|\sfB\|_{L^\infty}\|\sfB\|_{\s-\frac12} +\|\sfB\|_{L^\infty}^2\|\eta\|_{W^{1,\infty}}\|\eta\|_{\s+\frac12}+ 
       \|\sfB\|_{\s-\frac12}\|\sfB\|_{L^\infty}\|\eta\|_{W^{1,\infty}}^2\\
        &+ \gamma \left(\| \eta\|_{L^\infty} \|\psi\|_{\s+\frac12}+  \| \psi\|_{W^{1,\infty}} \|\eta\|_{\s-\frac12}+  \| \sfB\|_{\s-\frac32}+ \| \sfV\|_{\s-\frac32} \| \eta\|_{W^{1,\infty}}+ \| \sfV\|_{L^\infty}\| \eta\|_{\s-\frac12}\right)\\
        \lesssim_\s &
      (1 + \|(\eta, \psi, \sfV, \sfB)\|_{X^{2}} 
      )\| (\eta,\psi,\sfV,\sfB)\|_{X^{\s}} \lesssim_\s \| (\eta,\psi,\sfV,\sfB)\|_{X^{\s}}\,.
    \end{align*}
  Estimate   \eqref{pa_t:etapsi} is proved.
\end{proof}
\begin{proof}[Proof of Lemma \ref{lemma:para.bella}]
$(i)$ {\sc Paralinearization:} Formula \eqref{para:DN} is proved in \cite[Proposition 7.5]{BD}; for its explicit expression, see formula (3.2) in \cite{BFP}.
The fact that the smoothing remainder has the form $R(\eta)\psi$ follows by observing that, since $G(\eta)\psi,\ \sfV,$ and $\upomega$ are linear in $\psi$, then by difference also the smoothing remainders $R_1(\eta)\upomega + R_2(\eta,\upomega)\eta$,
appearing in  \cite[Formula (3.2)]{BFP}, must be linear in $\psi$.
Note also that by \cite[Proposition 7.5]{BD} we have that $R(\eta)$ is a real operator.

 We now prove that $b_{\geq 1}^{(-1)}$ satisfies \eqref{real:b-1}. The expansion in \eqref{real:b-1} follows by definition of the class $\Sigma \Gamma_1^{-1}[r, N] \ni b^{(-1)}_{\geq 1}$. We now check the symmetry properties. Since $G(\eta)$, $|D|$ and $R(\eta)$ are real operators, it follows by difference that $\Opbw{b_{\geq 1}^{(-1)}}$ is so. Thus, from \eqref{realetoreale} it follows that $\overline{b_{\geq 1}^{(-1)}} = \big(b_{\geq 1}^{(-1)}\big)^\vee$, implying, by homogeneity, that each homogeneous component satisfies $\overline{b_p^{(-1)}} = (b_p^{(-1)})^\vee$, $p=1, \ldots, N-1$. Therefore, it remains to show that $\overline{b_{p}^{(-1)}} = b_{p}^{(-1)}$ for $p=1,\ldots,N-1$.
    Recalling that $\upomega = \psi - \Opbw{\sfB}\eta$ (see  \eqref{GU}) and using the composition formula \eqref{comp01A}, one has 
 \begin{align}
     G(\eta)\psi&= \Opbw{|\xi|+ b_{\geq 1}^{(-1)}}\upomega+ \Opbw{-\ii \sfV\xi- \tfrac{1}{2} \sfV_x}\eta+ R(\eta)\psi\notag\\
     &= \Opbw{|\xi|+ b_{\geq 1}^{(-1)}}\psi+ \mathfrak{g}(\eta)\psi+ R(\eta)\psi,
     \label{DN_para}
 \end{align}
where, in the second equality, we absorb into $R(\eta)$
the remainder arising from the composition of $\Opbw{b_{\geq 1}^{(-1)}}$ and $\Opbw{\sfB}$, 
whereas $\mathfrak{g}(\eta)$ denotes the linear operator defined by
 \begin{align*}
     \mathfrak{g}(\eta)\psi & := \Opbw{-\ii \sfV(\eta,\psi;x)\xi- \tfrac{1}{2} \left(\sfV(\eta,\psi;x)\right)_x-|\xi|\#_\vr \sfB(\eta,\psi;x)- b_{\geq 1}^{(-1)}(\eta;x, \xi)\#_\vr \sfB(\eta,\psi;x)}\eta \\
     & =: \Opbw{a_{\geq 1}^{(1)}(\eta,\psi;x,\xi)}\eta\,.
 \end{align*}
   Taylor expanding  the symbols 
   $b_{\geq 1}^{(-1)}$ and the smoothing remainder $R(\eta)\psi$ we obtain
 \begin{align}\label{decomposed}
 &\Opbw{ b_{\geq 1}^{(-1)}}= \sum_{p=1}^{N-1}\Opbw{b_{p}^{(-1)}}+\Opbw{ b_{\geq N}^{(-1)}}\,, \quad  R(\eta)= \sum_{p=1}^{N-1} R_p(\eta)+R_{\geq N}(\eta).
 \end{align}
 Then, expanding in homogeneity the identity $G(\eta) = G(\eta)^\top$ and using \eqref{A1b} and the already proved identity $\overline{b_p^{(-1)}} = (b_p^{(-1)})^\vee$, $p=1, \ldots, N-1$, we deduce 
\begin{align}
     \!\!\!\!\Opbw{\overline{b_p^{(-1)}} -b_p^{(-1)}}= \mathfrak{g}_p(\eta)- R^\top_p(\eta) +R_p(\eta)-\mathfrak{g}_p(\eta)^\top  \ , \quad p = 1,\ldots, N-1 \ ,
     \label{equ:brutta}
 \end{align}
 where, using that  $a_{\geq 1}^{(1)}$ is a symbol in $\Sigma \Gamma_{1}^1[r, N]$,
we defined
\begin{equation}
    \mathfrak{g}_p(\eta)\psi := \Opbw{a_p^{(1)}(\eta, \psi; \cdot)}\eta \,.
    \end{equation}
We claim that  the  r.h.s. of \eqref{equ:brutta} is a smoothing operator, namely that 
\be\label{tutto-smooth}
\mathfrak{g}_p(\eta) {-} \mathfrak{g}_p(\eta)^\top {+}R_{p}(\eta){-} R_{p}(\eta)^\top \in \wt \cR^{-\vr}_p\,.
\ee 
Then one eventually substitutes, if needed, $b_p^{(-1)}\leadsto \Re (b_p^{(-1)})$ for any $p=1,\ldots, N-1 $ including the imaginary part in the smoothing remainder. 

 \medskip 
 To prove  \eqref{tutto-smooth} 
we need to show that the Fourier-Taylor coefficients of the maps in \eqref{tutto-smooth}
fulfill the bounds \eqref{smoocara} with $m\leadsto - \vr$. 
We compute them using \eqref{equ:brutta}. First remark that
each homogeneous component  $ b_{p}^{(-1)}$   has Fourier expansion as in \eqref{espr.hom.sym}, namely 
 \begin{align*}
     &b_{p}^{(-1)}(\eta;x,\xi)=  \sum_{\substack{\vec{\jmath}=(j_1,\ldots,j_{p})\in \Z^p}} b_{\vec{\jmath}}(\xi) \eta_{{j_1}}\cdots \eta_{{j_p}} \, e^{\ii {(j_1+\ldots+j_p)} x}
 \end{align*}
where the coefficients $b_{\vec{\jmath}}(\xi) = b_{j_1,\ldots,j_p}(\xi)$ are symmetric in $(j_1,\ldots,j_p)$ (cf. \eqref{sym_sy}).
Then, by the quantization formula \eqref{BW},  we have
\begin{align}
   & \Opbw{\overline{b_p^{(-1)}} -b_p^{(-1)}}\psi = \sum_{j_1+\ldots+j_p+k=\ell}B_{j_1,\ldots,j_p,k,\ell}\eta_{j_1}\cdots \eta_{j_p}\psi_k e^{\ii \ell x}, \\
    & \quad   B_{j_1,\ldots,j_p,k,\ell}:=  \chi_p 
   \left(\vec{\jmath},\frac{k+\ell }{2}\right)\left[
   -b_{\vec{\jmath}}\left(\frac{k+\ell}{2}\right) +\overline{b_{-\vec{\jmath}}\left(\frac{k+\ell}{2}\right)}\right], \quad \vec{\jmath}:=(j_1,\ldots,j_p)\,. 
\end{align}
Thanks to the paradifferential spectral localization given by the cut-off $\chi_p$, we have the bounds 
\begin{align}
    B_{j_1,\ldots,j_p,k,\ell}\not=0\quad \implies \quad \begin{cases}
        j_1+\ldots+j_p+k=\ell\\
        \max\{ |j_1|, \ldots, |j_p|\} \ll |k|\sim |\ell|\\ 
        | B_{j_1,\ldots,j_p,k,\ell}|\lesssim \max\{ |j_1|, \ldots, |j_p|\}^\mu |k|^{-1}\,. 
    \end{cases}
    \label{spec_loc:para}
\end{align}
To improve them into smoothing bounds of the form \eqref{smoocara} (with $m\leadsto - \vr$), we shall prove finer estimates of  the Fourier-Taylor coefficients of the operators in  the r.h.s. of \eqref{tutto-smooth}.

We start by computing the Fourier-Taylor coefficients of 
\be\label{quant_g}
\mathfrak{g}_p(\eta)\psi := \Opbw{a_p^{(1)}(\eta, \psi; \cdot)}\eta = \sum_{j_1+\ldots+ j_{p}+k=\ell} G_{j_1,\ldots, j_{p},k, \ell } \eta_{j_1}\cdots\eta_{j_{p}} \psi_k e^{\ii \ell x} \ . 
\ee
Writing 
\begin{align*}
     &a_{p}^{(1)}(\eta,\psi;x,\xi)=  \sum_{\substack{\vec{\jmath}=(j_1,\ldots,j_{p-1},k)\in \Z^p}} a_{\vec{\jmath}}(\xi) \eta_{{j_1}}\cdots \eta_{{j_{p-1}}} \psi_k\, e^{\ii {(j_1+\ldots+j_p)} x},
 \end{align*}
where the  coefficients $a_{\vec{\jmath}}(\xi) = a_{j_1,\ldots,j_{p-1},k}(\xi)$ are symmetric in $(j_1,\ldots,j_{p-1})$,  and using the quantization formula \eqref{BW}, 
\begin{align}
    \Opbw{a_p^{(1)}}\eta=\sum_{\substack{\vec{\jmath}=(j_1,\ldots,j_{p-1},k)\in \Z^{p}\\ (j_p,\ell)\in \Z^2\\ j_1+\ldots+j_{p}+k=\ell}}\chi_p 
\left(\vec{\jmath},\frac{j+k}{2}\right)
 a_{\vec{\jmath}}\left(\frac{j_p+\ell}{2}\right) 
 \eta_{{j_1}}\cdots \eta_{{j_{p-1}}} \psi_k \eta_{j_p}
{e^{\im \ell x}} \ . 
\label{quant_a}
\end{align}
Comparing \eqref{quant_g} and  \eqref{quant_a}, after a symmetrization, we obtain
 \begin{align*}
 G_{j_1,\ldots, j_{p},k, \ell }
       = \chi_p 
 \left((\vec{\jmath}^{(2)},k),\frac{j^{(1)}+\ell }{2}\right)
 a_{\vec{\jmath}}\left(\frac{j^{(1)}+\ell}{2}\right), \quad |j^{(1)}|=\max\{ |j_1|, \dots, |j_p|\}, \quad \vec{\jmath}^{(2)}:= (j_1,\dots , \wh{j^{(1)}}, \dots, j_p)\,.
 \end{align*}
 In view of the  spectral localization given by the cut-off $\chi_p$, one has
 \begin{align}
 G_{j_1,\ldots, j_{p},k, \ell }\not=0 \quad \implies \quad \begin{cases}
 j_1+\ldots+ j_{p}+k=\ell\\
     {\rm max}\{ |\vec{\jmath}^{(2)}|, |k|\} \ll |j^{(1)}| \sim |\ell|, \\ |G_{j_1,\ldots, j_{p},k, \ell }|\lesssim \max\{ |\vec{\jmath}^{(2)}| , |k|\}^\mu |j^{(1)}|\,.
     \label{para_loc:g}
 \end{cases}
 \end{align}
 On the other hand, one has 
 \begin{align}
     \mathfrak{g}_p(\eta)^\top \psi= \sum_{j_1+\ldots+ j_{p}+k=\ell} \big(G^\top\big)_{j_1,\ldots, j_{p},k, \ell } \eta_{j_1}\cdots\eta_{j_p} \psi_k e^{\ii \ell x}, \qquad \big(G^\top\big)_{j_1,\ldots,j_p,k, \ell } := {G_{j_1,\ldots, j_{p},-\ell, -k }} \,.
 \end{align}
 Then, in view of \eqref{para_loc:g}, we have
 \begin{align}
 \big(G^\top\big)_{j_1,\ldots, j_{p},k, \ell }\not=0 \quad \implies \quad \begin{cases}
 j_1+\ldots+ j_{p}+k=\ell\\
     \max\{ |\vec{\jmath}^{(2)}|, |\ell|\} \ll |j^{(1)}| \sim |k|, \\ |G^\top_{j_1,\ldots, j_{p},k, \ell }|\lesssim \max\{ |\vec{\jmath}^{(2)}|, |\ell|\}^\mu|j^{(1)}| \,.
     \label{para_loc:gT}
 \end{cases}
 \end{align}
 Thus, one has
 \begin{gather}
          \max\{ |j_1|, \ldots, |j_{p}|,|k|\} \sim {\rm max}_2\{ |j_1|, \ldots, |j_{p}|, |k|\}, \quad \text{and}\notag\\ 
          |\big(G^\top\big)_{j_1,\ldots, j_{p},k, \ell }|\lesssim {\rm max}_2\{ |j_1|, \ldots, |j_{p}|,|k|\}^{\mu'} \max\{ |j_1|, \ldots, |j_{p}|,|k|\}^{-\vr}\,,\label{G:tras}
 \end{gather}
 proving that $ \mathfrak{g}_p(\eta)^\top$ is a smoothing remainder in $\wt \cR^{-\vr}_p$. Thus, the operator $G_{j_1,\ldots, j_p, k,\ell}=\big(G^\top\big)_{j_1,\ldots, j_p, -\ell,k} $ actually satisfy
 \begin{align}
     |G_{j_1,\ldots, j_{p},k, \ell }|\lesssim {\rm max}_2\{ |j_1|, \ldots, |j_{p}|,|\ell|\}^{\mu'} \max\{ |j_1|, \ldots, |j_{p}|,|\ell|\}^{-\vr}
     \label{nuova:G}\,.
 \end{align}
Consider now the remainder
$$
 R_p(\eta)\psi= \sum_{j_1+\ldots+ j_{p}+k=\ell} R_{j_1,\ldots, j_{p},k, \ell } \eta_{j_1}\cdots\eta_{j_{p}} \psi_k e^{\ii \ell x}\,. 
$$
Since $R_p(\eta)$ is a smoothing remainder, its coefficients satisfy
 \begin{align}
     R_{j_1,\ldots, j_{p},k, \ell }\not=0 \quad \implies \quad \begin{cases}
 j_1+\ldots+ j_{p}+k=\ell\\
     |R_{j_1,\ldots, j_{p}, k, \ell }|\lesssim {\rm max}_2\{ |j_1|, \ldots, |j_{p}|, |k|\}^\mu {\rm max}\{ |j_1|, \ldots, |j_{p}|, |k|\}^{-\vr} 
     \end{cases}
     \label{smoothing_coeff}
 \end{align}
 and its transpose $ R^\top(\eta)$ has coefficients $ \big( R^\top\big)_{j_1,\ldots, j_{p},k, \ell }$ satisfying 
      \begin{align}
    \big( R^\top\big)_{j_1,\ldots, j_{p},k, \ell }=R_{j_1,\ldots, j_{p},-\ell,-k }\not=0 \quad \implies \quad \begin{cases}
 j_1+\ldots+ j_{p}+k=\ell\\
     |R_{j_1,\ldots, j_{p}, k, \ell }|\lesssim \frac{{\rm max}_2\{ |j_1|, \ldots, |j_{p}|, |\ell|\}^\mu }{{\rm max}\{ |j_1|, \ldots, |j_{p}|, |\ell|\}^{\vr}} \,.
     \end{cases}
     \label{smoothing_coeff_tras}
 \end{align}
We now improve the estimate in \eqref{spec_loc:para} using \eqref{equ:brutta} and the support property $|k|\sim |\ell|$ that follows from \eqref{spec_loc:para}.
From \eqref{G:tras}, \eqref{nuova:G}, \eqref{smoothing_coeff}, and \eqref{smoothing_coeff_tras}, we get
\begin{align*}
    | B_{j_1,\ldots,j_p,k,\ell}|\leq & | G_{j_1,\ldots,j_p,k,\ell}|+|\big(R^\top\big)_{j_1,\ldots,j_p,k,\ell}|+ |R_{j_1,\ldots,j_p,k,\ell}|+  | \big(G^\top\big)_{j_1,\ldots,j_p,k,\ell}| \\
    \lesssim &\frac{{\rm max}_2\{ |j_1|, \ldots, |j_{p}|,|\ell|\}^{\mu''}}{ \max\{ |j_1|, \ldots, |j_{p}|,|\ell|\}^{\vr}}+ \frac{{\rm max}_2\{ |j_1|, \ldots, |j_{p}|,|k|\}^{\mu''} }{\max\{ |j_1|, \ldots, |j_{p}|,|k|\}^{\vr}}
    \stackrel{|k|\sim|\ell|}{\lesssim}\frac{{\rm max}_2\{ |j_1|, \ldots, |j_{p}|,|k|\}^{\mu''}}{ \max\{ |j_1|, \ldots, |j_{p}|,|k|\}^{\vr}}\,.
\end{align*}
This concludes the proof of \eqref{real:b-1}.\\
We now prove the expansion \eqref{G.mappa}.
Note that we have
\be\label{g.sviluppo}
\mathfrak{g}(\eta)= \sum_{p=1}^{N-1}\mathfrak{g}_p(\eta)+
 \mathfrak{g}_{\geq N}(\eta) \ , \quad 
 \mathfrak{g}_p(\eta) \in \wt \cM^{1}_p \,,  \quad 
\mathfrak{g}_{\geq N}(\eta)  \in \cM^1_{\geq N}[r] \,.
\ee
 Indeed, $\mathfrak{g}_p(\eta) \in  \wt \cM_p^1$ by \eqref{para_loc:g}, and the fact that $\mathfrak{g}_{\geq N} (\eta) \in \cM^{1}_{\geq N}[r]$ follows from \Cref{thm:contS}-Item(i) and observing that, since $a_{\geq N}^{(1)} \in \Gamma^{1}_{\geq N}[r]$ and $a_{\geq N}^{(1)}$ is linear in $\psi$, it satisfies the estimate
 $$
 |a_{\geq N}^{(1)}|_{1, L^\infty, 4} \lesssim \|\psi\|_{s_0} \|\eta\|^{N-1}_{s_0}\,.
 $$
Using \eqref{g.sviluppo} and identity  \eqref{DN_para}, estimate  \eqref{G.mappa} then follows
 with $M_{\geq 1}(\eta):=\Opbw{ b_{\geq 1}^{(-1)}}+ \mathfrak{g}(\eta)+ R(\eta)$ which is a map in  $\Sigma \cM_1^1[r,N]$, using that 
$\Opbw{ b_{\geq 1}^{(-1)}}$ is a map in $\Sigma \cM_1^1[r,N]$
by   the second bullet of \Cref{rmk:molto_bene}.

$(ii)$ {\sc Time derivative control:}
By Proposition 3.1 of \cite{BFP} and with the notations therein, one has that $b_{\geq 1}^{(-1)}$ and $R(\eta)$ are respectively a smoothing symbol in $\Sigma \Gamma^{-1}_{1, 0, 1}[r, N] \subset \Gamma^{-1}_{1, 0, 1}[r]$ and an operator in $\Sigma \cR_{1, 0, 1}^{-\vr}[r, N] \subset \cR^{-\vr}_{1,0, 1}[r]$; therefore, using Eq. (2.20) of \cite{BFP} with $k = K= 1$, $K'= 0$ and $N =1$ one has
\begin{equation}\label{pallavolo.2}
 \hspace{-0.28cm} \begin{aligned}
    \| \pa_t \left( R(\eta) v \right)\|_{\s + \vr - \frac 1 2} &\lesssim_\s \left(\|v\|_{\s + \frac 1 2} + \|\pa_t v\|_{\s - \frac 1 2} \right) \left( \|\eta\|_{s_0} + \|\pa_t \eta\|_{s_0}\right) \\
    & \quad  \quad + \left(\|v\|_{s_0} + \|\pa_t v \|_{s_0} \right)  \left( \|\eta \|_{\s + \frac 1 2} + \|\pa_t \eta\|_{\s - \frac 1 2}\right)\,,
\end{aligned}
\end{equation}
and, using Eq. (2.8) of \cite{BFP} with 
$k = K = 0$, $K'=0$,  $N=1$, and 
$k = K = 1$, $K'=0$, $N=1$, also
\begin{equation}
 | { b^{(-1)}_{{\geq 1}}}|_{-1, W^{\s +  1\!/\!2 - s_0}, M} + 
  | {\pa_t b^{(-1)}_{{\geq 1}}}|_{-1, W^{\s +  1\!/\!2 - s_0}, M} \lesssim_{\s, M} 
\left(\norm{\eta}_{\s+\frac12}
+ \norm{\pa_t\eta}_{\s-\frac12}\right)\,. 
\end{equation}
Then \eqref{pa_t:smooth} and \eqref{pa_t:smooth2} follow from \eqref{pa_t:etapsi} and recalling \eqref{norm:Xs}\,.

\end{proof}

\begin{proof}[Proof of \Cref{lem:B.V.a}]\label{proof.lemma5.3}
$(i)$   First consider the functions $\sfV, \sfB$ in \eqref{def:V}--\eqref{form-of-B} and   $\sfV_\gamma$ in \eqref{def:Vgamma}.
In view of the analyticity of the Dirichlet-Neumann operator $\eta \mapsto G(\eta)$ in   \eqref{eq:DNan},  there is $\s_0>\frac32$ so that for any $\s \geq \s_0$ there is $r = r(\s)>0$ so that the map
\be\label{an:DNX}
(\eta, \psi) \mapsto \Big(\sfB, \sfV, \sfV_\gamma \Big)  \, , \quad 
B_{\s}(r)\times \dot H^{\s}(\T; \R) \to \Big(H^{\s-1}(\T; \R) \Big)^3  \, , 
\quad \text{are { analytic}} \ , 
\ee
linear in $\psi$, and fulfilling 
\be
\label{B.eta.psi}
\norm{\sfB}_{\s- 1}  + \norm{\sfV}_{\s-1} \lesssim_\sigma  \norm{\psi}_{\s } \ ,  \quad 
 \norm{\sfV_\gamma}_{\s-1}\lesssim_\s  \norm{\psi}_{\s }+ \norm{\eta}_{\s}  \ ,   \quad  \forall (\eta, \psi) \in B_{\s}(r')\times \dot H^{\s}(\T; \R)
\ee
as stated in \eqref{a.v.b.estimate}.

Consider now $\sfa$ in \eqref{a.taylor}. Clearly $(\eta, \psi)\mapsto \pa_x \sfB(\eta, \psi)$ is real analytic 
$ B_\s(r) \times \dot H^{\s}(\T;\R) \to H^{\s-2}(\T;\R)$ for $\s \geq \s_0$ (enlarging $\s_0$). 

In view of the analyticity of the original water waves vector field $\cX_{\scH_\gamma} = (X^{(\eta)} , X^{(\psi)})$ (cf. \eqref{wwana}), also 
the map $(\eta, \psi) \mapsto \pa_t \sfB \equiv \di \sfB[\cX_{\scH_\gamma}]$, explicitly 
given by 
\begin{equation}
    \begin{aligned}\label{pat.B}
   (\pa_t \sfB) (\eta, \psi) &= \frac{- G(\eta)(\sfB X^{(\eta)}) - \pa_x(\sfV X^{(\eta)}) + G(\eta)X^{(\psi)} + [X^{(\eta)}]_x \psi_x + \eta_x [X^{(\psi)}]_x-2 \sfB \, \eta_x [X^{(\eta)}]_x}{1+\eta_x^2}
\end{aligned}
\end{equation}
is analytic as a map $  B_\s(r') \times  B_\s(r') \to H^{\s-2}(\T;\R) $. 
Note that in \eqref{pat.B} we also  used the ``shape-derivative formula'' 
$
\di G(\eta)[\wh \eta]\psi = - G(\eta)(\sfB \wh \eta) - \pa_x(\sfV \wh \eta)
$ (see \cite{LannesLivre}).

 From formula  \eqref{pat.B}, using the bounds in \eqref{wanna},\eqref{B.eta.psi} we also derive the bound
\begin{align}
\| \pa_t \sfB \|_{\s-2}
    \lesssim_\s \|\eta\|_{\s} +  \|\psi\|_{\s}\,, \quad 
    \forall (\eta, \psi) \in B_{\s}(r')\times  B_{\s}(r') \ . 
    \label{a.v.b.estimate.1}
\end{align}
The analyticity of $(\eta, \psi)\mapsto \sfa$ as stated in \eqref{BVa.ana} and the second estimate \eqref{a.v.b.estimate} follow. 


\smallskip

$(ii)$ The real-valued functions $\sfV, \sfB$ belong to $\Sigma \cF^\R_1[r, N]$ by \cite[Proposition 7.4 \& 7.5]{BD}. Consequently, also $\sfV_\gamma \in \Sigma \cF^\R_1[r, N]$. 
Consider now  $\pa_t \sfB$ and its explicit formula in   \eqref{pat.B}.
This is  a real-valued function in $\Sigma \cF^\R_1[r, N]$, as it follows since 
 by 
 \eqref{G.mappa} and \eqref{det.map},  $$
 G(\eta) - |D| \in \Sigma \cM_{1}^1[r,N]\ , \quad  
  X^{(\eta)} -  G(0) \psi\in \Sigma\cF^\R_2[r,N] \ , 
  \quad 
  X^{(\psi)}+\eta - \gamma G(0)\pa_x^{-1}\psi \in  \Sigma  \cF^\R_{2}[r,N]$$ 
  and  $\sfB, \sfV \in \Sigma \cF^\R_1[r, N]$. 
Then it follows that also $\sfa \in \Sigma \cF^\R_1[r, N]$.

\smallskip 

$(iii)$ The Taylor expansion of $\sfB$ in \eqref{esp.b.gamma.1.2} follows from its definition \eqref{form-of-B}, using the Taylor expansion of $G(\eta)$ given in \eqref{DNexp}--\eqref{DNexp2}.
Then  substitute \eqref{esp.b.gamma.1.2} in \eqref{def:V} to obtain that 
$\sfV(\eta, \psi) = \psi_x - \eta_x |D| \psi$ up to a function in $\Sigma \cF_{3}^\R[r,N]$, and then substitute such expansion in \eqref{def:Vgamma}, proving \eqref{esp.v.gamma.1.2}.

\smallskip 

$(iv)$ If $(\eta, \psi, \sfV, \sfB) \in B_{X^\s}(I;r)$ (cf. \eqref{norm:Xs}) then $(\eta, \psi) \in B_{\s+\frac12}(r)$. The estimate for $\pa_t \sfB$ follows from \eqref{a.v.b.estimate.1}.
Next consider $\sfV$. By \eqref{def:V}, we have $\pa_t \sfV = \pa_t \psi_x - \sfB \pa_t \eta_x - \eta_x \pa_t \sfB$, therefore using  Lemma \ref{lem:tutto.X},   \eqref{a.v.b.estimate.1}  one has
\begin{align*}
    \| \pa_t \sfV\|_{\s -\frac32 } \lesssim_{\s} \| \pa_t \psi\|_{\s - \frac12} + \|\sfB\|_{\s-\frac32}\| \pa_t \eta\|_{\s - \frac12} +\| \eta\|_{\s-\frac12} \|\pa_t \sfB\|_{\s-\frac32} \lesssim_{\s} \|(\eta, \psi, \sfV, \sfB)\|_{X^\s}
\end{align*}
and recalling the definition of $\sfV_\gamma$ (cf. \eqref{def:Vgamma}) analogous estimate follows for $\sfV_\gamma$. It remains to estimate $\pa_t \sfa\equiv \di \sfa[\cX_{\scH_\gamma}]$. Since by \eqref{BVa.ana} the map $\sfa: B_{\s-\frac 1 2}(r')\times B_{\s-\frac 1 2}(r') \to H^{\s-\frac 5 2}(\T;\R)$ is analytic with estimate \eqref{a.v.b.estimate}, one has by Cauchy estimate
\begin{align*}
    \|\di \sfa(\eta, \psi)[\cX_{\scH_\gamma}(\eta, \psi)]\|_{\s-\frac 5 2} & \lesssim_\s \|\di \sfa\|_
    {\cL(H^{\s- 1\!/\!2}\times H^{\s-1\!/\!2}; H^{\s-5\!/\! 2})}
    \|\cX_{\scH_\gamma}(\eta, \psi)\|_{\s-\frac 1 2}\\
    &\lesssim_\s (r')^{-1} \sup_{\|(\eta, \psi)\|_{\s-\frac 1 2}\leq r'} \|\sfa(\eta, \psi)\|_{\s-\frac52} \|\cX_{\scH_\gamma}(\eta, \psi)\|_{\s-\frac 1 2}\\
    &
    \lesssim_\s \|\eta\|_{\s + \frac 12} + \|\psi\|_{\s + \frac 1 2} \lesssim_\s \|(\eta, \psi, \sfV, \sfB)\|_{X^\s}\,,
\end{align*}
where in the second to last passage we have used also estimate \eqref{pa_t:etapsi}.
\end{proof}

\section{Flows and Conjugations}\label{sec:flows}

Following \cite{BD,FI2,BFP,BMM2,MMS}, in this section we collect several results concerning the conjugation of paradifferential operators and smoothing remainders under flows of the form
\begin{equation}\label{flusso1para0}
\begin{cases}
\partial_\tau {\bf \Phi}^\tau(\zak)= \bG(\zak;\tau),{\bf \Phi}^\tau(\zak)\,,\\
{\bf \Phi}^0(\zak)=\uno,
\end{cases}
\qquad {\bf \Phi}(\zak):={\bf \Phi}^\tau(\zak){|_{\tau=1}} \,,
\end{equation}
where the generator $\bG(\zak;\tau)$ is either a matrix of paradifferential operators or a matrix of smoothing remainders.
These conjugation rules are used in \Cref{sec:paraWW} to iteratively conjugate  vector fields of the form
\[
X(\zak;U)=-\ii \vOmega(D)\zak+ \Opbw{\bA(\zak;x,\xi)}U+ \bR(\zak)U+ \bB_{\geq N}(\zak)U\,,
\]
where $\vOmega(D)$ denotes the linear dispersion relation, $\bA(\zak;x,\xi)$ is a matrix-valued symbol,
$\bR(\zak)$ is a matrix of smoothing remainders, and $\bB_{\geq N}(\zak)$ is a matrix of bounded
remainders of homogeneity $\mathcal{O}(\zak^N)$, for some $N\in \N$.

Accordingly, each result in this section provides an expansion of each of the following conjugation:
\begin{enumerate}
    \item {\bf Conjugation of a paradifferential operator:} \(
    {\bf \Phi}^\tau(\zak)\,\Opbw{\bA(\zak;x,\xi)}\,{\bf \Phi}^\tau(\zak)^{-1}\,;
    \) occasionally the specific case of the conjugation of the dispersion relation \(    {\bf \Phi}^\tau(\zak)\,[-\ii \vOmega(D)]\,{\bf \Phi}^\tau(\zak)^{-1}\,
     \) is also treated separately;
    \item {\bf Conjugation of the remainder terms:} \(
    {\bf \Phi}^\tau(\zak)\,\bR(\zak)\,{\bf \Phi}^\tau(\zak)^{-1}\) and \(
    {\bf \Phi}^\tau(\zak)\,\bB_{\geq N}(\zak)\,{\bf \Phi}^\tau(\zak)^{-1};
    \)
    \item {\bf Time conjugation:} \(
    \left(\partial_t {\bf \Phi}^\tau(\zak)\right)\,{\bf \Phi}^\tau(\zak)^{-1}.
    \)
\end{enumerate}
The main difference with respect to the aforementioned papers \cite{BD,FI2,BFP,BMM2,MMS}
is that   we always need to  replace the time derivatives  with autonomous symbols and operators, substituting the equations of motions using \eqref{claim:zak.M}. 

\medskip

\noindent{\bf Conjugation by a flow generated by  a real-valued symbol of order one.} \label{sec:flow-1}
For $p=1,2$, consider the flow  ${\bf \Phi}^\tau(\zak)$, $\tau\in[-1,1]$, in \eqref{flusso1para0} with generator
\be\label{flusso1para1}
 \bG(\zak; \tau):= \vOpbw{\frac{\beta(\zak;x)}{1+\tau \beta_x(\zak;x)}\ii \xi}\,, \qquad
 \beta\in \wt \mF_p^\R \,. 
\ee
In \Cref{lem:flow.ad} it is shown that $\bPhi^\tau(\zak)$ is a $(0,3)$-admissible transformation. Moreover we recall that, by estimate \eqref{est.flow} (see also Lemma 3.22 in \cite{BD}), for any  
	$\zak \in B_{s_0,\R}(r)$ with $s_0>0$ sufficiently large and $r>0$ sufficiently small, 
	there is a constant $C_s >0$ such that, for any $ W\in H^{s}(\T; \C^2)$ and $\zak \in B_{s_0,\R}(r)$, one has
\begin{equation}\label{bound_flow_tras}
    \norm{ \bPhi^\tau (\zak) W }_{s}+\norm{ \bPhi^\tau (\zak)^{-1} W }_{s}
		 \leq   C_s  \| W \|_{s}\,.
    \end{equation}
\smallskip 
Following \cite{BD}, we define  the  path of diffeomorphisms of $ \T $ via 
\begin{equation}
	\label{eq:diffeo}
	\Psi \pare{ \zak, \tau ;x } :=  x + \tau \beta\pare{\zak; x}
	 \quad \text{with inverse}
	 \quad 
	 \Psi^{-1}\pare{ \zak, \tau;y }  := y + \breve{\beta}\pare{\zak, \tau;y} 
	\,, \quad 
	  \   \breve{\beta}\in  \Sigma \cF_{p}^\R\bra{r,N}\,, 
\end{equation}
for any $N\in \N$,  and set  $\Psi(\zak;x):=\Psi (\zak, 1;x)$.

\begin{proposition}[{\bf Conjugations for a transport flow}]
	\label{prop:Egorov}
	Let $m\in \R$, $\varrho>0$, $p=1,2$, $q\in \N_0$, $N\geq q+1$ and let ${\bf \Phi}(\zak)$ be the flow generated by \eqref{flusso1para0} with generator in \eqref{flusso1para1}. Then
	\begin{enumerate}
		
		\item \label{item:Egorovi} {\bf Conjugation of a para-differential operator:} Let $a_{\geq q}^{{(m)}} \in  \Sigma\Gamma^{m}_q\bra{r, N}$ be  a real-valued symbol and
		\begin{align}
		a^{\pare{m}}_\Psi\pare{\zak;x, \xi} :=  \left( a_{\geq q}^{{(m)}}\pare{\zak; y, \xi \ \partial_y\Psi^{-1}\pare{\zak;y}}\right)\Big|_{y=\Psi\pare{\zak;x}}
			\in  \Sigma\Gamma^{m}_q\bra{r, N} \,. 
            \label{simbolo:diffeo}
            \end{align}
		Then $a^{\pare{m}}_\Psi- a^{\pare{m}}_{\geq q}$ is  a real-valued symbol in $ \Sigma \Gamma_{p+q}^m[r,N]$ and 
		\begin{equation}
			\label{eq:conj_generic_symbol}
				\bPhi(\zak)  \vOpbw{\ii \, 
                a_{\geq q}^{{(m)}}\pare{\zak;x,\xi}}  \bPhi(\zak)^{- 1} =  \  \vOpbw{\ii a^{\pare{m}}_\Psi\pare{\zak;x, \xi} 
					+ \ii a^{\pare{m-2}}_{\geq (p+q)}  \pare{\zak;x, \xi}
				} + \bR_{\geq (p+q)}(\zak)\,,
		\end{equation} 
		where $a^{\pare{m-2}}_{\geq (p+q)}\pare{\zak;x,\xi}$
        is  a real-valued symbol  in $ \Sigma \Gamma^{m-2}_{p+q}\bra{r, N}$,
         whereas $ \bR_{\geq (p+q)}(\zak)$ is a real-to-real  matrix of smoothing operators in $ \Sigma \cR^{-\vr+m }_{p+q}\bra{r, N} $.\\
		In addition if $m=1$ and  $ a^{(m)}_{\geq q}\pare{\zak;x,\xi} = V\pare{\zak;x} \xi $ for some $ V \in  \Sigma \cF^\R_q\bra{r, N} $, 
        then in \eqref{eq:conj_generic_symbol} $ a^{\pare{ m-2 }}_{\geq (p+q)}\equiv 0 $   and 
        \begin{align}\label{trans.transp}
             a^{\pare{ m }}_{\Psi} \pare{\zak;x, \xi} = \Big(  V+ \beta V_x- V\beta_x+ {g_{\geq (q+2p)}}(\zak;x) \Big) \xi  \, \qquad  {g_{\geq (q+2p)}}\in \Sigma \cF^\R_{q+2p}[r, N]\,.
        \end{align}  

		\item{\bf  Conjugation of remainders:} If $ \bR_{\geq q}(\zak)$ is a real-to-real matrix of  smoothing operators in $ \Sigma  \cR^{-\vr}_q\bra{r, N} $  and $ \bB_{\geq N}(\zak)$ is a real-to-real matrix of spectrally localized maps in $\cS_{\geq N}^0[r]$, then 
		\begin{gather}
			\bPhi \pare{\zak}  \bR_{\geq q} \pare{\zak}  \bPhi\pare{\zak}^{-1} = \bR_{\geq q}(\zak)+ \bR_{\geq (p+q) }\pare{\zak},
            \label{conjtrasp_smoo}\\
            \bPhi \pare{\zak}  \bB_{\geq N}(\zak)  \bPhi\pare{\zak}^{-1}= \check\bB_{\geq N}(\zak) \label{conj:geq7_tras}
		\end{gather}
		where $\bR_{\geq (p+q) }\pare{\zak}$ is a real-to-real matrix of smoothing operators in $\cR^{-\vr+{N}}_{\geq (p+q)}\bra{r} $ and $\check \bB_{\geq N}(\zak)$  is a real-to-real matrix of spectrally localized maps in $\cS^{0}_{\geq N}\bra{r} $ . 
		
		\item {\bf Conjugation of $\partial_t$:} If $ \zak $ is the variable in  \eqref{zak} (thus solving \eqref{iniziale_zak}), then
		\begin{equation}\label{detconju}		
			(\partial_t \bPhi\pare{\zak})\,  \bPhi\pare{\zak}^{-1} = \vOpbw{\ii g(\zak;x)\xi } + {\bR_{\geq 2p}} \pare{\zak}, \qquad  g(\zak;x):= (\beta_t -\beta_t \beta_x ) + {g_{\geq 3p}(\zak;x)}
		\end{equation}
		where ${ g_{\geq 3p}\pare{\zak;x}}$  is  a real-valued function in $\Sigma \cF_{{3p}}^\R\bra{r,N}$ and $ {\bR_{\geq 2p}\pare{\zak}}$ is a real-to-real matrix of smoothing operators in $  \Sigma \cR^{-\vr}_{ {2p}}\bra{r,N}$.
	\end{enumerate} 
 \end{proposition}
	\begin{proof}
 During the proof we shall denote $b:= b(\zak;\tau, x):= \frac{\beta(\zak;x)}{1+\tau \beta_x(\zak;x)}$.
We first note that, as $ \bG(\zak;\tau)= \begin{pmatrix}
     1&0\\0& 1
 \end{pmatrix} \Opbw{\ii b \, \xi}$, one has $\bPhi^\tau(\zak)= \begin{pmatrix}
     1&0\\0& 1
 \end{pmatrix} \Phi^\tau(\zak)$ where $ \Phi^\tau(\zak)$ solves the scalar equation 
 \begin{align}
      \pa_\tau \Phi^\tau(\zak)= \Opbw{\ii b \, \xi}\Phi^\tau(\zak)\,.
      \label{scalar_flow}
 \end{align}
\noindent 	1. The thesis follows by Lemmas A.4 and A.5 in \cite{BFP} whose proofs are mainly contained in \cite{BD}.
The sole difference is that we claim that the symbol $ a^{(m-2)}_{\geq (p+q)}$ appearing in \eqref{eq:conj_generic_symbol} is real-valued, which can be checked following the algebraic steps of the proof of \cite{BD} and is a consequence of the reality of $a_{\geq q}^{(m)}$ and $\beta$.

\noindent 2. Equation \eqref{conjtrasp_smoo} follows as in \cite[Remark at pag. 89]{BD}  (see also \cite[Proposition A.2]{MMS} for details). We prove now \eqref{conj:geq7_tras}. 
Using the bound \eqref{bound_flow_tras} and estimate \eqref{bound:specloc} for $\bB_{\geq N}(\zak)$,  for any $\zak \in B_{s_0,\R}(r)$ and $V \in H^s(\T;\C^2)$ we get 
\begin{equation*}
    \| \bPhi(\zak)\bB_{\geq N}(\zak)\bPhi(\zak)^{-1}V\|_s\lesssim_s\|    \bB_{\geq N}(\zak)\bPhi(\zak)^{-1}V\|_s\lesssim_s  \| \zak\|_{s_0}^N \| \bPhi(\zak)^{-1}V\|_s\lesssim_s \| \zak\|_{s_0}^N \| V\|_s\,,
\end{equation*}
proving that $\breve \bB_{\geq N}(\zak)$ in \eqref{conj:geq7} belongs to $ \cS_{\geq N}^0[r]$.

\noindent 3.  We note that $\left(\pa_t \bPhi(\zak)\right)\bPhi(\zak)^{-1}= \begin{pmatrix}
    1&0\\ 0&1
\end{pmatrix} \pa_t \Phi(\zak) \Phi(\zak)^{-1}$, with $\Phi(\zak)$ in \eqref{scalar_flow}. 
The proof is the same of \cite[Proposition 3.28]{BD} and  \cite[Lemma A.4]{BFP}. The only difference is that we use an autonomous point of view, namely we substitute $\pa_t b $ using the equation \eqref{iniziale_zak} for $ \zak$ at any degree of homogeneity.
First, using \eqref{claim:zak.M} and that $\beta \in \wt\cF^\R_p$, $p=1,2$, we have 
\begin{align*}
  \beta_t (\zak;x)= \di \beta (\zak; x)[ - \im {\bf \Omega}(D)\zak + \bM_{\geq 1}(\zak) \zak] \in \Sigma\cF_{p}^\R[r, N] \,. 
\end{align*}
Then we compute the time derivative of the function $b=b(\tau,\zak;x)$ as
\begin{align}
    \pa_t b= \frac{\beta_t }{(1+{\tau} \beta_x)}- \tau \frac{{\beta}\beta_{tx}}{(1+\tau \beta_x)^2}= \beta_t - \tau \beta_t \beta_x - \tau \beta \beta_{t\, x}+ g_{\geq 3p}({\tau},\zak;x), \qquad  g_{\geq 3p}\in \Sigma \cF^{\R}_{3p}[r, N]\,, 
    \label{esp:dtb}
\end{align}
thus proving that $ \pa_t b $ belongs to $ \Sigma \cF^{\R}_{p}[r, N]$. 
Then  \eqref{detconju} follows arguing as in \cite{BD,BFP}.
	\end{proof}

\noindent{\bf Conjugation by a flow generated by  a real-valued symbol of order $<1$.} \label{sec:semiflow}
Given  $p=1,2$, $m<1$,  consider the flow  ${\bf \Phi}^\tau(\zak)$, $\tau\in[-1,1]$ in \eqref{flusso1para0} with generator 
\be\label{flusso2para}
\bG(\zak):= \vOpbw{\im f(\zak;x,\xi) } \,, \qquad 
  f(\zak; x, \xi) \mbox{ real-valued symbol in }  \wt \Gamma^m_p \,. 
\ee
In \Cref{lem:flow.ad} it is shown that $\Phi^\tau(\zak)$ is a $(0,\tfrac{3}{2})$-admissible transformation. Moreover, it is standard (see e.g.  Lemma 3.22 in \cite{BD}) that, for any  
	$\zak \in B_{s_0,\R}(r)$ with $s_0>0$ sufficiently large and $r>0$ sufficiently small, 
	the operator 
	$ \Phi^\tau (\zak)$  actually belongs to  $\cL\pare{H^s(\T; \C^2)} $ for any $s \in \R$. More precisely, there is a constant $C_s>0$
     such that	for any $ W\in H^{s}(\T; \C^2)$ and $\zak \in B_{s_0,\R}(r)$,
        \begin{align}
            \norm{ {\bf \Phi}^\tau (\zak) W }_{s}+\norm{ {\bf \Phi}^\tau (\zak)^{-1} W }_{s}
		 \leq   C_s  \| W \|_{s}\,.
         \label{analoga}
        \end{align}
We  set 
${\bf \Phi}(\zak):= {\bf \Phi}^1(\zak)$.
Following \cite[Lemma A.6]{BFP}, we have the following result.
\begin{proposition}[{\bf Conjugations for a semi-FIO}]
	\label{prop:FIO}
	Let $m<1$, $\varrho>0$, $p=1,2$, $q\in \N_0$, $N\geq q+1$ and let ${\bf \Phi}(\zak)$ be the flow generated by \eqref{flusso1para0} with generator in \eqref{flusso2para} at $\tau =1$. 
	\begin{enumerate}
		
		\item \label{item:FIO} {\bf  Conjugation of a para-differential operator:} Let $a\equiv a_{\geq q}^{{(m')}} \in  \Sigma\Gamma^{m'}_q\bra{r,N}$ be  a real-valued symbol, $m' \in\R$.
		Then 
		\begin{equation}
			\label{eq:conj_generic_symbol_semi}
			\begin{aligned}
				&{\bf \Phi}(\zak)  \vOpbw{\ii \, a}  {\bf \Phi}(\zak)^{- 1} = \\
                & \  \vOpbw{\ii \left[ a + \{ f, a\} + \frac12 \{f, \{f, a \}\}+
                a^{(m + m'-3)}_{q+p}
                + a^{(2m+m'-4)}_{q+2p}
                +a^{(3m+m'-3)}_{q+3p}
                \right]
				} + \bR_{\geq p+q}(\zak)\,,
			\end{aligned} 
		\end{equation} 
		where $a^{(m + m'-3)}_{q+p} \in \Sigma\Gamma^{m + m'-3}_{q+p}[r,N]$,  $a^{(2m+m'-4)}_{q+2p} \in \Sigma\Gamma^{2m+m'-4}_{q+2p}[r,N]$ and
       $ a^{(3m+m'-3)}_{q+3p}\in \Sigma\Gamma^{3m+m'-3}_{q+3p}[r,N]$ are real-valued, 
       and
       $ \bR_{\geq p+q}(\zak)$ is a matrix of smoothing operators in $\Sigma \cR^{-\vr +m' }_{ q+p}[r,N]$.
       \item \label{item:disp_rel}{\bf Conjugation of the dispersion relation:} If $\Omega(\xi)\in \wt \Gamma^{\frac12}_0$ is the Fourier multiplier in \eqref{omegonejin}
       \begin{align}
 {\bf \Phi}(\zak)\vOpbw{-\ii\Omega(\xi)}   {\bf \Phi}(\zak)^{-1}= \vOpbw{-\ii\Omega(\xi)+ f \#_\vr \Omega(\xi)- \Omega(\xi)\#_\vr f + \ii d_{\geq 2p}^{(2m - \frac32)}}+ \bR_{\geq p}(\zak)
 \label{conjfourier}
\end{align}
   where    $d_{\geq 2p}^{(2m - \frac32)}\in \Sigma \Gamma_{2p}^{2m-\frac32}[r, N]$  and $\bR_{\geq p}(\zak)$ is a matrix of smoothing remainders in $ \Sigma \cR^{-\vr+m+\frac12}_{p}[r, N]$.

		\item \label{item:smoo_rem}{\bf Conjugation of  remainders:} If $ \bR_{\geq q}(\zak)$ is a real-to-real matrix of  smoothing operators in $ \Sigma  \cR^{-\vr}_q\bra{r, N} $ and 
$\bB_{\geq N}(\zak)$ is a real-to-real  matrix of spectrally localized maps in $\cS_{\geq N}^0[r]$,
       then		
       \begin{align}\label{conj:rem}
			&{\bf \Phi}(\zak)  \bR_{\geq q}(\zak)  {\bf \Phi}(\zak)^{-1} \;= \bR_{\geq q}(\zak)+ \bR_{\geq (p+q) }\pare{\zak}  \,, \\
             &{\bf \Phi}(\zak)  \bB_{\geq N}\pare{\zak}  {\bf \Phi}(\zak)^{- 1} =  \hat{\bB}_{\geq N}(\zak)\,,
         \label{conj:geq7}
		\end{align}
		where $\bR_{\geq p+q }\pare{\zak}$ is a real-to-real matrix of smoothing operators in $\Sigma\cR^{-\vr+N^+}_{p+q}\bra{r,N} $, $N^+:= \begin{cases}
		    N & m>0\\
            0& m\leq 0.
		\end{cases}$,
        whereas
          $\hat{\bB}_{\geq N}(\zak)$ is a matrix of real-to-real  spectrally localized maps in $\cS_{\geq N}^0[r]$. 
		\item \label{item:det} {\bf Conjugation of $\partial_t$:} If $ \zak $ is the variable in  \eqref{zak} $($thus solving \eqref{iniziale_zak}$)$, then
		\begin{equation}\label{detconju.semi}		
			\left( \partial_t {\bf \Phi}(\zak)\right)\,  {\bf \Phi}(\zak)^{-1} = \vOpbw{\ii \left[ f_{\cX_\cH} +
           \frac12 \{f, f_{\cX_\cH}\} + 
            a^{(2m-3)}_{2p} + 
             a^{(3m-2)}_{3p}\right] 
            } + \bR_{\geq 2p} \pare{\zak}\,, 
		\end{equation}
		where the real-valued symbol 
        $f_{\cX_\cH}(\zak;\cdot) := \di_\zak  f(\zak)[\cX_\cH(\zak)]  \in \Sigma\Gamma^m_{\geq p}[r, N]$ fulfills 
        \be\label{fXH}
          f_{\cX_\cH}  -  \di_\zak  f(\zak)[- \im {\bf \Omega}(D)\zak] \in \Sigma\Gamma^m_{p+1}[r, N]\,, 
        \ee
       $ a^{(2m-3)}_{2p}\in \Sigma\Gamma^{(2m-3)}_{2p}[r,N] $, $a^{(3m-2)}_{3p} \in \Sigma\Gamma^{(3m-2)}_{3p}[r,N]$  are real-valued 
        and $ \bR_{\geq 2p}\pare{\zak}$ is a real-to-real matrix of smoothing operators in $ \Sigma \cR^{-\vr}_{2p}\bra{r,N}$.
	\end{enumerate} 
 \end{proposition}
\begin{proof}
\Cref{item:FIO,item:disp_rel}  follow as in \cite[Lemma A.6]{BFP}.
We prove 
    \Cref{item:smoo_rem}. For the spectrally localized operator $\bB_{\geq N}(\zak)$ it  follows by \eqref{analoga} and \eqref{bound:specloc} arguing as in the proof of \eqref{conj:geq7_tras} of Proposition \ref{prop:Egorov}.
    We prove now \Cref{item:smoo_rem} for the smoothing remainder $\bR_{\geq q}(\zak)$. The proof follows using the Lie expansion (see \cite[(A.3)]{BFP})
    \begin{align*}
    {\bf \Phi}^{1} (\zak) \bR_{\geq q}(\zak) ({\bf \Phi}^{1}(\zak))^{-1}   
&=  { \bR_{\geq q}(\zak) + 
\sum_{q=1}^{L}\frac{1}{q!}{\rm Ad}_{\vOpbw{\im f}}^{q}\Big[\bR_{\geq q}(\zak) \Big]} \\ 
& +\frac{1}{L!}\int_{0}^{1}
(1- \theta )^{L} 
{\bf \Phi}^{\theta}(\zak){\rm Ad}_{\vOpbw{\im  f}}^{L+1}\Big[ \bR_{\geq q}(\zak) \Big]({\bf \Phi}^{\theta}(\zak))^{-1} \di \theta \,,  
\end{align*}
{choosing $L = N-1$ so that $ (L+1)\,p+q\geq N$} and using also \Cref{comp-mappe-int} of  \Cref{lem:general_composition}. 

$4$. We follow \cite[Lemma A.7]{BFP}, using the  Lie expansion  (see \cite[(A.4)]{BFP}) to write
\begin{align*}
    (\pa_t {\bf \Phi}^{1} (\zak) )({\bf \Phi}^{1}(\zak))^{-1}   
&=  \vOpbw{ \im \pa_t f} +
\Big[\vOpbw{\im f}, \vOpbw{\im \pa_t f} \Big]
+
\sum_{q=3}^{L}\frac{1}{q!}{\rm Ad}_{\vOpbw{\im f}}^{q-1}\Big[\vOpbw{\im \pa_t f} \Big] \\ 
& +\frac{1}{L!}\int_{0}^{1}
(1- \theta )^{L} 
{\bf \Phi}^{\theta}(\zak){\rm Ad}_{\vOpbw{\im  f}}^{L}\Big[ \vOpbw{ \im \pa_t f} \Big]({\bf \Phi}^{\theta}(\zak))^{-1} \di \theta \,. 
\end{align*}
{The expansion in \eqref{detconju} follows transforming   $\pa_t f \equiv f_{\cX_\cH}(\zak):=  \di f(\zak)[\cX_{\cH}(\zak)]$ into an autonomous symbol by using the equations of motions.
To prove that
$f_{\cX_\cH}(\zak)$ is a symbol, use that  
  $f\in \wt\Gamma^m_p$, $p=1,2$, is a homogeneous symbol, than   Lemma \ref{lem:general_composition} and \eqref{claim:zak.M} give
 $$
 (\pa_t f)(\zak) =
 f_{\cX_\cH}(\zak):= 
 \di f(\zak)[- \im {\bf \Omega}(D)\zak + \bM_{\geq 1}(\zak) \zak] = \di f(\zak)[- \im {\bf \Omega}(D)\zak] + f_{\geq p+1}(\zak) \,, 
 $$
 where $ f_{\geq p+1}(\zak) :=\di f(\zak)[- \im {\bf \Omega}(D)\zak+ \bM_{\geq 1}(\zak)\zak] -   \di f(\zak)[- \im {\bf \Omega}(D)\zak] $ is  a real-valued symbol in 
$ \Sigma \Gamma^m_{p+1}[r, N] $.
}
\end{proof}
\noindent{\bf Conjugation by flows generated by linear smoothing operators.} 
Finally,  we study the conjugation rules for a flow 
$\bPhi(\zak):= \bPhi^{\tau}(\zak)\vert_{\tau =1}$ generated by 
\begin{equation}\label{flusso smoothing}
	\partial_\tau \bPhi^\tau(\zak) = \bQ(\zak) \circ \bPhi^\tau(\zak)\,,  \ \ 
	\bPhi^0(\zak) = \uno \,,
\end{equation}
with $\bQ(\zak)$ a matrix of smoothing operators in $\wt\cR^{-\varrho}_p$, $p=1,2$.
\smallskip

\begin{proposition}[{\bf Conjugation by flows generated by smoothing operators}]\label{prop:Egorov_smoothing}
	Let $m\in \R$, $p=1,2$, $ N\in \N$, $\vr\geq m$, 
    $\vr', r >0$. Let   $\bQ(\zak)$ be a matrix of smoothing operators in $\widetilde{\mathcal{R}}^{-\varrho}_p$ and $\bPhi(\zak)$ be the flow generated by $ \bQ(\zak) $ as in \eqref{flusso smoothing}. 
	Then the following holds:
	
	\begin{enumerate}
		\item {\bf Space conjugation:} \label{item:Egorov_smoothingi} If $ a\equiv a_{\geq 2}^{{(m)}} \in  \Sigma\Gamma^m_2[r,N] $   then 
\begin{align}\label{eq:conj_symbol_smoothing1}
			 \bPhi(\zak)\circ \vOpbw{a\pare{\zak;\cdot}} \circ  \bPhi(\zak)^{-1}
			-  \vOpbw{a \pare{\zak;\cdot}}    &\in  \Sigma\cR^{-\vr + \max\{m,0\}}_{ p+2 }[r,N]\,,\\
    \label{eq:conj_symbol_smoothing2}
			\bPhi(\zak)\circ  (-\im \vOmega(D) )  \circ \bPhi(\zak)^{-1}- \pare{- \im \vOmega(D)+ [\bQ(\zak), -\im \vOmega(D)] } & \in  \Sigma\cR^{-\vr + \frac12}_{ 2 p }[r,N] \,. 
		\end{align}
These matrices of operators are real-to-real provided $\bQ(\zak)$ is. 
		\item {\bf Conjugation of remainders:} \label{item:Egorov_smoothingii}If $ \bR(\zak)$ is a real-to-real matrix of smoothing operators in $\Sigma\cR^{-\vr'}_q[r,N] $, $q=1,2$,
         and $\bB_{\geq N}(\zak)$ is a matrix of spectrally localized maps in $\cS_{\geq N}^0[r]$, then 
		\begin{align} \label{bottadevita}
		 \bPhi(\zak)\circ \bR\pare{\zak}  \circ  \bPhi(\zak)^{-1} - \bR\pare{\zak} \in  & \Sigma\cR^{-\min\{\vr,\vr'\}}_{ q+p }[r,N] \,, \\
\label{eq:conj_symbol_smoothingzero}
			\bPhi(\zak)\circ \bB_{\geq N}\pare{\zak}\circ  \bPhi^{-1}(\zak)
			-  \bB_{\geq N}\pare{\zak}   &\in  \cR^{-\vr}_{\geq N}[r]\,. 
		\end{align}
		These matrices of operators are  real-to-real provided $\bQ(\zak)$ is. 
		\item  \label{item:Egorov_smoothing_iii} {\bf Conjugation of $\partial_t$:} If $ \zak $ is the variable in  \eqref{zak} (thus solving \eqref{zetone:formale}), then
		\begin{equation}\label{claimed}
			(\partial_t   \bPhi\pare{\zak})\circ \bPhi\pare{\zak}^{- 1} =
             \begin{cases}
               \bQ\pare{ -\im \vOmega(D) \zak} + \bR_{\geq 2}(\zak)\,, \quad &\mbox{ if } p=1\\
               \bQ(-\ii \bOmega(D) \zak, \zak ) +  \bQ(\zak, -\ii \bOmega(D) \zak ) + \bR_{\geq 3}(\zak)\,, \quad &\mbox{ if } p=2
             \end{cases}
		\end{equation}
        where $\bR_{\geq p}(\zak)$, $p=2, 3$, are matrices of smoothing operators in $\Sigma  \cR^{-\vr+{\frac32}}_{  p}\bra{r,N}$. Moreover
        $(\partial_t   \bPhi\pare{\zak})\circ \bPhi\pare{\zak}^{- 1}$
		 is real-to-real provided $\bQ(\zak)$ is. 
	\end{enumerate}
\end{proposition}

\begin{proof}
First, we note that by \Cref{flow.s.ad} the map $\bPhi(\zak)$ is a $(0,0)$-admissible transformation with gain $\vr$ and therefore belongs to $\cM^0_{\geq 0}[r]$. We now prove its claimed transformation rules.

1. We prove \eqref{eq:conj_symbol_smoothing1}. Denoting $\bA(\zak):=\vOpbw{a\pare{\zak;\cdot}}$ we have  
\be\label{UAU-1}
\bPhi(\zak)\circ \bA(\zak) \circ \bPhi(\zak)^{-1} = \bA(\zak) + \bR(\zak)\,,   
\quad \bR(\zak)= \bA(\zak)(\bPhi(\zak)^{-1}-\Id) + (\bPhi(\zak)- \Id ) \bA(\zak) \bPhi(\zak)^{-1}\,.
\ee
Therefore, the fact that $\bR(\zak)\in \Sigma\cR_{p+2}^{-\vr+\max\{m, 0\}}[r, N]$ follows from a repeated application of  \Cref{comp-mappe,comp-smo-pseudo} of  \Cref{lem:general_composition} using that $\bPhi^\tau(\zak) \in \cM^0_{\geq 0}[r]$ and the fact that, by \eqref{est:smotR},
\be\label{Ups-1}
\bPhi(\zak) - \Id = \int_0^1 \bQ(\zak) \bPhi^\tau(\zak) \, \di \tau \in \Sigma \cR_{p}^{-\vr}[r,N] \,. 
\ee
We prove now \eqref{eq:conj_symbol_smoothing2}. By a Taylor expansion we have that
\begin{equation}\label{agosto}
\begin{aligned}
\bPhi(\zak)\circ  (-\im \vOmega(D) )  \circ \bPhi(\zak)^{-1}&= - \im \vOmega(D)+ [\bQ(\zak), -\im \vOmega(D)] 
\\
&+ \frac{1}{2} \int_0^1 (1-\tau) \bPhi^\tau(\zak)  [\bQ(\zak),[\bQ(\zak), -\im \vOmega(D)] ]\bPhi^{-\tau}(\zak)\, \di \tau\,.
\end{aligned}
\end{equation}
By \Cref{comp-mappe} of \Cref{lem:general_composition} one has that 
$[\bQ(\zak),[\bQ(\zak), -\im \vOmega(D)] ]$ is a matrix of smoothing operators in $\Sigma\cR_{2p}^{-\vr + \frac{1}{2}}[r, N]$. 
Next, using \eqref{est:smotR},  we obtain that the integral term in \eqref{agosto}  is a matrix of smoothing operators in $\Sigma\cR_{2p}^{-\vr + \frac{1}{2}}[r, N]$.

2.  One follows the same lines as in the proof of \eqref{eq:conj_symbol_smoothing1}, using only \Cref{comp-mappe} of Lemma \eqref{lem:general_composition}.
Formula 
\eqref{eq:conj_symbol_smoothingzero} 
is proved arguing as  in \eqref{UAU-1} and using  estimates \eqref{bound:specloc}
 and \eqref{Ups-1}.

3. 
We introduce the conjugated quantity $\Theta^{\tau}(\zak) := (\bPhi^{\tau}(\zak))^{-1} \partial_t \bPhi^{\tau}(\zak)$ and, derivating that quantity w.r.t. the variable $\tau$ and using \eqref{claim:zak.M}, we obtain 
\begin{align*}
\partial_{\tau} \Theta^{\tau}(\zak) 
&= (\bPhi^{\tau}(\zak))^{-1} \di \bQ(\zak)[ - \im {\bf \Omega}(D)\zak + \bM_{\geq 1}(\zak) \zak]\, \bPhi^{\tau}(\zak) 
\end{align*}
with 
$\bM_{\geq 1}(\zak)\in \Sigma \cM^{\frac32}_1[r,N]$, 
whose solution, as $\Theta^0(\zak)=0$, is given by
\begin{equation*}
\Theta^\tau(\zak) \equiv  
(\bPhi^{\tau}(\zak))^{-1} \partial_t \bPhi^{\tau}(\zak) = \int_0^1 (\bPhi^{\tau}(\zak))^{-1}\, \di \bQ(\zak)[ - \im {\bf \Omega}(D)\zak + \bM_{\geq 1}(\zak) \zak] \, \bPhi^{\tau}(\zak)\, \di \tau\,.
\end{equation*}
Conjugating with 
$\bPhi^\tau(\zak)$ and performing the change of variables $\theta = 1-\tau$, we deduce
\begin{align*}
(\partial_t \bPhi(\zak)) \, \bPhi^{-1}(\zak) & = \int_0^1 \bPhi^{\tau}(\zak)\, \di \bQ(\zak)[ - \im {\bf \Omega}(D)\zak + \bM_{\geq 1}(\zak) \zak]\, \bPhi^{-\tau}(\zak)\, \di \tau\\
& = \begin{cases}
\int_0^1 
\bPhi^{\tau}(\zak)\, \bQ\big(- \im {\bf \Omega}(D)\zak + \bM_{\geq 1}(\zak) \zak\big)\, \bPhi^{-\tau}(\zak)\, \di \tau\,, \quad &p=1\\
\int_0^1 
\bPhi^{\tau}(\zak)\, 2\bQ\big(- \im {\bf \Omega}(D)\zak + \bM_{\geq 1}(\zak) \zak, \zak\big)\, \bPhi^{-\tau}(\zak)\, \di \tau\,, \quad &p=2
\end{cases}\,.
\end{align*}
Next, thanks to \Cref{lem:general_composition}-\Cref{smooth-comp-2}, we have that
the right hand side is a matrix of  smoothing operators in {$\Sigma\cR^{-\vr+\frac32}_p[r,N]$}. 
Moreover, by \eqref{bottadevita} (for simplicity we write the formula only for $p=1$, then for $p=2$ one argues analogously)
\begin{equation*}
 \int_0^1 \bPhi^{\tau}(\zak)\, \bQ\big(- \im {\bf \Omega}(D)\zak + \bM_{\geq 1}(\zak) \zak\big)\, \bPhi^{-\tau}(\zak)\, \di \tau 
 = 
 \bQ\big(- \im {\bf \Omega}(D)\zak\big) + \bR_{2}(\zak)\,, 
\end{equation*}
with 
$\bR_{2}(\zak)$ a  matrix of smoothing operators in  $\Sigma\cR_{2}^{-\vr+\frac32}[r,N]$, 
obtaining the claimed expansion \eqref{claimed}.
\end{proof}

\section{Local existence and continuity of $Z(t)$}\label{app:cauchy.z}
This section is devoted to the proof of \Cref{loc.ex}.
The first step is to derive \emph{a priori} bounds for the low and high Sobolev norms of the original variables $(\eta, \psi, \sfV, \sfB)$, as well as for the good unknown $\upomega$. To this end, we exploit in a crucial way the relations between $\eta$, $\sfV$, $\sfB$, and $\upomega = \psi - \Opbw{\sfB}\eta$ to gain an additional half-derivative of regularity for $\upomega$. Although $\upomega$ would naively belong to $H^{s+\frac12}(\T;\R)$ when $\eta, \psi \in H^{s+\frac12}(\T;\R)$, a cancellation mechanism ensures that $\upomega$ actually lies in the more regular space $H^{s+1}(\T;\R)$ provided   $\sfV$, $\sfB$ are taken in $H^{s}(\T;\R)$. 
The next lemma provides a quantitative statement of this fact.
\begin{lemma}\label{lem:magic_formula}
   Let $s_0$ as in \eqref{s0r}. For any $ \s\geq s_0$ there is $ r=r(\s)>0$ such that for any  $\eta,\psi\in B_{s_0}(r)\cap H^{\s+\frac12}(\T; \R)$ and $ \upomega, \sfV, \sfB$ defined in  \eqref{GU}, \eqref{def:V}, \eqref{form-of-B},  one has 
    \begin{align}        &\| \upomega\|_{\s+1} \lesssim \| \sfV\|_{\s}+ \| \sfB\|_{s_0}\|\eta\|_\s+ \| \eta\|_{s_0}\| \sfB\|_{\s}\,,\label{stima_omega}\\
        &\| \sfV - \upomega_x\|_{\s}\lesssim  \| \sfB\|_{s_0}\|\eta\|_\s+ \|\eta\|_{s_0}\| \sfB\|_{\s}\,,\label{stima_V}\\
        & \| \sfB- |D| \upomega\|_\s \lesssim
        (\| \sfV\|_{s_0} +
        \|\psi\|_{s_0} + 
        \| \sfB\|_{s_0} \| \eta\|_{s_0})\|\eta\|_\s+
        \left( 
        \|\sfV\|_\s+\|\psi\|_\s \right) \|\eta\|_{s_0} \,. 
        \label{stima_B}
    \end{align}
\end{lemma}
\begin{proof}
{\sc Proof of \eqref{stima_B}:}
    By equations \eqref{espr_DN}, \eqref{para:DN} and \eqref{bonyeq0}, together with \eqref{compsymbols} and \eqref{poisson}, one has 
    \begin{align}\notag
        \sfB& = G(\eta)\psi+ \sfV \eta_x \\
        &= |D|\upomega + \Opbw{
        -\ii \sfV \xi - \frac{1}{2}\sfV_x
        + \sfV\#_\vr (\ii \xi) }\eta + \Opbw{b_{\geq 1}^{(-1)}(\eta)}\upomega+ \Opbw{\eta_x }\sfV
        + R(\eta)\psi + \cR(\sfV, \eta_x)\notag\\
        &=|D|\upomega - \Opbw{ \sfV_x }\eta+\Opbw{\eta_x }\sfV
        +\Opbw{b_{\geq 1}^{(-1)}(\eta)}\upomega+ R(\eta)\psi + \cR(\sfV, \eta_x)
        \label{para:B}
    \end{align}
    where $R(\eta)\psi$ is the smoothing remainder of the Dirichlet-Neumann paralinearization formula \eqref{para:DN}, and  $\cR(\cdot, \cdot)$ is the smoothing remainder in the paralinearization formula
    \eqref{stima_BMM}.
    Then, using \eqref{para:B}, \eqref{cont00}, \eqref{stimapar2},  \eqref{piove},  \eqref{stima_BMM} (with $ \s_1=\s$ and $\s_2=\s_0$) and, by \eqref{GU}, $\norm{\upomega}_{\s} \lesssim \norm{\psi}_\s + \norm{\sfB}_{\s_0} \norm{\eta}_\s$, we obtain estimate \eqref{stima_B}.

{\sc Proof of \eqref{stima_V}:}
     On the other hand, 
 combining the definition \eqref{def:V} of $\sfV$, the definition  \eqref{GU} of $\upomega$ and using paraproduct formula \eqref{bonyeq0} and the symbolic calculus  \eqref{comp01A}, we have
    \begin{align}\label{BISMA1}
        \sfV = 
        \pa_x(\upomega  + \Opbw{\sfB} \eta) - \eta_x \sfB
       = \omega_x + \Opbw{\sfB_x}\eta-\Opbw{\eta_x}\sfB+ \cR(\sfB,\eta)  
    \end{align}
    where $\cR(\sfB,\eta)$ is the sum of the smoothing remainders appearing in \eqref{bonyeq0} and in \eqref{comp01A}. 
     Then, using \eqref{stimapar2},  \eqref{stima_BMM} (with $ \s_1=\s$ and $\s_2=s_0$) and \eqref{comp020}, we get 
    estimate \eqref{stima_V}.
    
{\sc Proof of \eqref{stima_omega}:}    Finally,  \eqref{stima_omega} is obtained by \eqref{BISMA1} putting in evidence $\upomega_x$ and taking the $H^\s$-norm.
\end{proof}
As a first application we  prove the equivalence between the norm of $U \in \dot H^\s_\R(\T; \C^2)$ and the norm of $(\eta, \psi, \sfV, \sfB) \in  X^{\s-\frac34}$.
\begin{lemma}\label{UXequiv}
There exists $\s'_0>0$ such that for any   $\s \geq \s'_0$, there exist $\tr'= \tr'(\s), C'>0$   such that if $U$  defined in \eqref{Ugrande} satisfies $U \in B_{\s'_0}(\tr') \cap \dot H^\s_\R(\T; \C^2)$ then  $(\eta, \psi, \sfV, \sfB) \in B_{X^{\s'_0-\frac34}}( C' \tr') \cap X^{\s-\frac34}$  and vice-versa
if $(\eta, \psi, \sfV, \sfB) \in B_{X^{\s'_0-\frac34}}( \tr') \cap X^{\s-\frac34}$ then $U \in B_{\s'_0}(C'\tr') \cap \dot H^\s_\R(\T; \C^2)$.
Moreover they have equivalent norms, i.e. there exists $C_\s>1$ such that 
 \be\label{eq:equiv0}
C_\s^{-1} \norm{(\eta, \psi, \sfV, \sfB)}_{X^{\s-\frac 3 4}} \leq \norm{U}_\s \leq C_\s  \norm{(\eta, \psi, \sfV, \sfB)}_{X^{\s - \frac 3 4}} \,,
 \ee
 where we recall the definition of $\| \cdot \|_{X^\s}$ in \eqref{norm:Xs}.
 \end{lemma}
\begin{proof}
 Assume first that  $ U \in B_{\s'_0}(\tr') \cap \dot H^\s_\R(\T; \C^2)$.    By \eqref{eq:etaomegaUequiv} with $\s + \frac 3 4\leadsto \s$, one has
    \be\label{U.e.om.1}
\norm{U}_{\s} \lesssim_\s \norm{\eta}_{\s - \frac 1 4} + \norm{\upomega}_{\s + \frac 1 4} \lesssim_\s \norm{U}_{\s}\,,
    \ee
    so in particular
    \be\label{0502:1722}
\norm{U}_{\s'_0} \simeq  \norm{\eta}_{\s'_0 - \frac 1 4} + \norm{\upomega}_{\s'_0 + \frac 1 4} \lesssim  \tr' \,. 
\ee
Thus it is sufficient to show that
    \be \label{equiv:pomeriggio}
    \norm{\eta}_{\s - \frac 1 4} + \norm{\upomega}_{\s + \frac 1 4} \lesssim_\s \norm{(\eta, \psi, \sfV, \sfB)}_{X^{\s - \frac 3 4}} \lesssim_\s  \norm{\eta}_{\s - \frac 1 4} + \norm{\upomega}_{\s + \frac 1 4}  \quad \forall \s \geq \s'_0\,.
    \ee
    By \eqref{0502:1722} and \eqref{B.eta.psi} with $\s = s_0 +1$, and choosing $\s'_0 \geq s_0 + \frac 5 4$,
    we have
    \be\label{v.b.control.1}
    \|\sfB\|_{s_0} + \|\sfV\|_{s_0} \lesssim \|\psi\|_{s_0+1} \lesssim \|\psi\|_{\s'_0 - \frac 1 4}\,.
    \ee
By the definition of $\upomega$ in  \eqref{GU} and  estimate  \eqref{cont00}, one has the  upper bound: for any $\s' \in \R$
\be\label{psi.norma.bassa}
\norm{\psi}_{\s'}\lesssim \| \upomega\|_{\s'}+ \norm{\sfB}_{L^\infty} \norm{\eta}_{\s'}\,.
\ee
Then, combining \eqref{psi.norma.bassa} 
with $\sigma' = \s_0'-\frac14$ and \eqref{v.b.control.1}
and using again the smallness of $\|\eta\|_{\s'_0-\frac14}$,
we deduce the low norm control
\be\label{small.n.small.r}
\| \sfB\|_{s_0} + \|\sfV\|_{s_0} \lesssim \|\upomega\|_{\s'_0 - \frac 1 4} \lesssim \tr'\,.
\ee
Then, by \eqref{stima_omega} and \eqref{0502:1722}, we deduce
\be
\|\eta\|_{\s -\frac14} + \|\upomega\|_{\s + \frac 1 4} \lesssim_\s \|\eta\|_{\s - \frac 1 4} + \|\sfV\|_{\s -\frac 3 4} + \|\sfB\|_{s_0} \|\eta\|_{\s - \frac 3 4} + \|\eta\|_{s_0} \|\sfB\|_{\s - \frac 3 4} \stackrel{\eqref{small.n.small.r}}{\lesssim_\s} \|(\eta, \psi, \sfV, \sfB)\|_{X^{\s - \frac 3 4}}\,,
\ee
which gives the first inequality in \eqref{equiv:pomeriggio}.

We show that also the second inequality holds true.
First, by \eqref{psi.norma.bassa} with $\sigma' = \sigma - \frac14$ and using the control on $\norm{\sfB}_{L^\infty}$ given by
\eqref{small.n.small.r}, we obtain
\be \label{psi.is.done}
\|\psi\|_{\s - \frac 1 4} \lesssim_\s \|\eta\|_{\s - \frac 1 4} + \|\upomega\|_{\s + \frac 1 4}\,, \quad \|\psi\|_{\s'_0 - \frac1 4} \lesssim \|\eta\|_{\s'_0 - \frac 1 4} + \|\upomega\|_{\s'_0 + \frac 1 4} \lesssim \tr'\,.
\ee
Furthermore, by \eqref{stima_V}--\eqref{stima_B} (using also the smallnesses in \eqref{small.n.small.r} and \eqref{0502:1722}) and \eqref{psi.is.done} we also have
\be\label{stima.v.b.327}
\|\sfV\|_{\s - \frac 3 4} + \|\sfB\|_{\s - \frac 3 4} \lesssim_\s \|\upomega\|_{\s + \frac 1 4} + \|\eta\|_{\s - \frac 3 4}  + \|\psi\|_{\s - \frac 3 4} \lesssim_\s \|\eta\|_{\s - \frac 1 4} + \|\upomega\|_{\s + \frac 1 4}\,.
\ee
Gathering the first of \eqref{psi.is.done} and \eqref{stima.v.b.327}, we get the second inequality in \eqref{equiv:pomeriggio}.
Note that
\eqref{equiv:pomeriggio} at $\sigma  = \sigma_0'$, together with \eqref{0502:1722}, gives also $(\eta, \psi, \sfV, \sfB) \in B_{\sigma_0'}(C'\tr')$ for some $C'>0$.

Vice-versa assume that  $(\eta, \psi, \sfV, \sfB) \in B_{X^{\s'_0-\frac34}}( \tr') \cap X^{\s-\frac34}$. By \eqref{stima_omega}, 
$\norm{\upomega}_{\s + \frac14} \lesssim   \norm{(\eta, \psi, \sfV, \sfB)}_{X^{\s - \frac34}}$, giving the left inequality of \eqref{equiv:pomeriggio}, the low-norm control 
$\norm{\upomega}_{\sigma_0' + \frac14} \lesssim \tr'$ 
and, arguing as above, also the right inequality. Then \eqref{eq:etaomegaUequiv} gives $U \in  B_{\s'_0}(C'\tr') \cap \dot H^\s_\R(\T; \C^2)$.
\end{proof}
In the following result we show the equivalence of norms between $Z$ solving \eqref{Z.eq} and the complex variable $U$ solving \eqref{complexo}. 
\begin{lemma}\label{UZequiv}
Let $N \in \N$, $\vr$, $\und{\vr}$ as in \Cref{thm:nf}. There exists $\s_0''>0$ such that for any   $\s \geq \s_0''$, there exist $\tr= \tr(\s)>0, C >0$   such that if  $Z \in B_{\s_0''}(\tr) \cap \dot H^\s_\R(\T; \C^2)$ then  $U$  defined in \eqref{Ugrande} satisfies
 $U \in B_{\s_0''}(C\tr) \cap \dot H^\s_\R(\T; \C^2) $, and vice-versa
 if $U \in B_{\s_0''}(\tr) \cap \dot H^\s_\R(\T; \C^2)$ then $Z \in B_{\s_0''}(C\tr) \cap \dot H^\s_\R(\T; \C^2)$.
 Moreover,   they have equivalent norms, i.e. there exists $C_\s>1$ such that 
 \be\label{eq:equiv}
C_\s^{-1} \norm{Z}_\s \leq \norm{U}_\s \leq C_\s  \norm{Z}_\s \,.
 \ee
\end{lemma}
\begin{proof}
Assume first that  $Z \in B_{\s_0''}(\tr) \cap \dot H^\s_\R(\T; \C^2)$. By \eqref{inv.F}, 
there exists $s'_0>0$ such that, provided $\s \geq \s_0'' \geq s_0'$,  the variable $\zak$ defined in \eqref{zak} fulfills 
\be\label{2410:1758}
\norm{\zak }_{s_0' -\frac12} \lesssim 
\norm{Z }_{s_0' } \lesssim \norm{Z }_{\s_0'' } \lesssim \tr \,, \qquad  \norm{\zak }_{\s - \frac 1 2} \lesssim_\s \norm{Z }_{\s} \,. 
\ee
We first prove the second inequality in \eqref{eq:equiv}. From \eqref{Def:ZU}, \eqref{de.Upsilon}  and \eqref{eq:complex_good} we deduce
\be
\label{Z.e.Us}
Z  = \Upsilon(\zak) \zak = \bPsi_w(\zak) \circ \bPsi_5(\zak) \circ \bT(\zak) U
\ee
where $\bPsi_w(\zak)$ (defined in  \Cref{thm:wnf}) is $(0, 0)$-admissible with gain $\ov{\vr} = \vr - 3N  -2$, $\bPsi_5(\zak)$ (defined in \Cref{thm:cubic.nf.1}) is $(0, 3)$-admissible with arbitrary gain $\wt \vr \geq \ov{\vr}$, and 
$\bT(\zak)$ (defined in \Cref{thm:quadratic.nf}) is $(0, \frac 9 2)$-admissible with gain $ \vr -2N$.
Linearly inverting such maps and substituting $\zak = \scF^{-1}(Z)$ (cf. \eqref{inv.F}) in the internal variable, we obtain
\begin{align}\label{0502:1753}
U = \bT^{-1}(\scF^{-1}(Z))\circ \bPsi_5^{-1}(\scF^{-1}(Z)) \circ  \bPsi_w^{-1}(\scF^{-1}(Z))   Z  \,. 
\end{align}
Then, by the  first of \eqref{2410:1758},    using repeatedly \eqref{lin.est.F} (with $\nu = 0$ and 
$\vr \leadsto  \ov{\vr}$) and eventually enlarging $\s_0''$, 
 for any $\s \geq \s_0 ''$ one gets
\begin{align}\nonumber
 \norm{U }_{\s} &\lesssim_{\s} \norm{Z }_{\s}+  \norm{\scF^{-1}(Z )}_{\s-\ov{\vr} }\norm{Z }_{\s_0''}
\stackrel{\eqref{2410:1758}}{\lesssim_{\s}}   \norm{Z }_{\s} + \norm{Z }_{\s - \ov{\vr} + \frac 1 2}\norm{Z }_{\s_0''}\nonumber
 \\
 &\stackrel{\ov{\vr} > \frac 1 2}{\lesssim_{\s}}  \norm{Z }_{\s} (1 + \norm{Z }_{\s_0''}) \lesssim_{\s}  \norm{Z }_{\s}\,,
 \label{UZ_bound}
\end{align}
where in the last passage we  used  $\norm{Z }_{\s_0''} < \tr$.
This gives the second inequality in \eqref{eq:equiv}. Note that
\eqref{UZ_bound} at $\sigma  = \sigma_0''$, together with the assumption  $Z \in B_{\sigma_0''}(\tr)$, gives also $U \in B_{\sigma_0''}(C \tr)$ for some $C>0$.
\\
We now prove the first inequality in \eqref{eq:equiv}. 
By  \eqref{Z.e.Us} and arguing as above 
$$
 \norm{Z }_{\s} \lesssim_{\s} \norm{U }_{\s}+  \norm{\zak }_{\s-\ov{\vr}}\norm{U }_{\s_0''}
 \stackrel{\eqref{2410:1758},\eqref{UZ_bound}}{\lesssim_\s} \norm{U }_{\s}+ \norm{Z }_{\s} \norm{Z }_{\s_0''}
$$
which implies the first inequality of \eqref{eq:equiv} using that  $\norm{Z }_{\s_0''}\leq \tr$.

Next assume  that  $U \in B_{\s_0''}(\tr) \cap \dot H^\s_\R(\T; \C^2)$. By \eqref{eq:complex_good}, $U = \mathscr{G}(\zak) := \pmb{\cG}_{\C}(\zak) \zak$ and in view of \Cref{loc.inv} this map is locally invertible, with inverse   $\zak = \mathscr{G}^{-1}(U)$, that, by \eqref{stima.inv.adm}, fulfills an estimate analogous to \eqref{2410:1758} with $Z \leadsto U$. Then the proof proceeds analogously to the previous case. 
\end{proof}

\begin{proof}[Proof of \Cref{loc.ex}]\label{proof.Lemma.loc.ex}
Let $\s_0 \geq \max(s_0, \s_0', \s_0'')$, with $s_0$ the regularity threshold in \Cref{thm:nf}, 
and $\s_0'$, $\s_0''$ those of \Cref{UXequiv}--\Cref{UZequiv}.
Let $Z_0 \in B_{\s_0}(\tr) \cap \dot H^\sigma_\R(\T; \C^2)$. 
Define $\zak_0 := \scF^{-1} (Z_0)$  with $\scF^{-1}$ as in \eqref{inv.F}.
By the properties of $\scF^{-1}$, $\zak_0\in B_{\s_0-\frac12}(C\tr)$. 
Then define
$\vect{\eta_0}{\psi_0} := \Lmap \zak_0 \in B_{\s_0 - \frac34}(C\tr)\times B_{\s_0 - \frac14}(C\tr)$ 
with $\Lmap$ in \eqref{zak}, $\sfB_0:= \dfrac{G(\eta_0) \psi_0 + (\eta_0)_x (\psi_0)_x}{1 + (\eta_0)_x^2}$ as in \eqref{form-of-B}, $\sfV_0:= (\psi_0)_x - (\eta_0)_x \sfB_0$ as in \eqref{def:V}
and with both of them in $B_{\s_0 - \frac 7 4}(C \tr)$.
Then define 
 $\upomega_0:= \psi_0 - \Opbw{\sfB_0} \eta_0 \in B_{\s_0 -\frac 3 4}(C \tr)$ as in \eqref{GU} and $U_0:= \Lmap^{-1} \vect{\eta_0}{\upomega_0} \in B_{\s_0 - 1}(C \tr)$ as in \eqref{Ugrande}.
 Actually, by \Cref{UXequiv} and  \Cref{UZequiv}, one has the improved regularity, for any $\s \geq \s_0$
 \be\label{0502:1822}
 \|Z_0\|_{\s} \simeq_\s \|U_0\|_{\s} \simeq_\s \|(\eta_0, \psi_0, \sfV_0, \sfB_0)\|_{X^{\s - \frac 3 4}}\,, \quad 
  \|Z_0\|_{\s_0} \simeq \|U_0\|_{\s_0} \simeq \|(\eta_0, \psi_0, \sfV_0, \sfB_0)\|_{X^{\s_0 - \frac 3 4}} \lesssim \tr\,. 
 \ee
 Consider now the initial data $(\eta_0, \psi_0, \sfV_0, \sfB_0)$ for equation \eqref{eq:etapsi}.
By the local Cauchy theory\footnote{See Alazard–Burq–Zuily \cite{ABZ1} for a local existence result for irrotational water waves using the same framework as in the present work, or Ifrim–Tătaru \cite{IFT} for a proof in the case of nonzero constant vorticity, employing holomorphic coordinates. See also \cite{CL} for a well-posedness result with general vorticity.}, given a small  initial datum as in \eqref{0502:1822} there are a time $T_{{\rm loc}}>0$ and a unique solution of \eqref{eq:etapsi} of the form
\begin{align}
    (\eta(t),\psi(t), \sfV(t), \sfB(t))\in C\left([-T_{{\rm loc}}, T_{{\rm loc}}];\ B_{\s_0-\frac 3 4}(C \tr)\cap {X^{\s-\frac34}}\right)\,.
\end{align}
Then, unfolding the changes of variables and using the equivalence of norms in \Cref{UXequiv,UZequiv},  we obtain the existence of a solution 
$$
Z(t) \in L^{\infty}\left( [-T_{{\rm loc}}, T_{{\rm loc}}];\ B_{\s_0}(C \tr) \cap {H^\s_\R(\T;\C^2)}\right) \,. 
$$
It remains to show that the map $t \mapsto Z(t)$ is continuous with respect to the norm $ \| \cdot \|_\s$. We first prove that the function $t \mapsto \norm{Z(t)}_\s$ is continuous.
Performing an energy estimate analogous to the one in \eqref{3001:1517}, we get 

     \begin{align*}
       \frac{\di}{\di t}\| Z(t)\|_{\s}^2 \lesssim_\s \norm{Z(t)}_{\s_0}^2
            \| Z(t)\|_\s^2 \lesssim_\s \tr^2  \| Z(t)\|_s^2\,
            \quad
            \Rightarrow
            \quad
            \| Z(t)\|_\s^2\leq e^{C_\s  \tr^2 t } \norm{Z(0)}_\s^2 \  \ \forall |t| \leq T_{{\rm loc}}\,.
        \end{align*}
        Consequently the time derivative $\frac{\di}{\di t}\| Z(t)\|_{\s}^2$ is bounded on the interval $[-T_{{\rm loc}},\ T_{{\rm loc}}]$ and the function $t \mapsto \norm{Z(t)}_\s$ is continuous. Moreover, since $Z(t)$ solves \eqref{Z.eq}, one has the bound 
    \begin{equation*}
        \| \pa_t Z(t)\|_{\s-1}\lesssim_\s \| Z
        (t)\|_\s\leq  e^{(C_\s  \tr^2 t)/2 } \norm{Z(0)}_\s \qquad \forall\, |t| \leq T_{{\rm loc}}\,,
    \end{equation*}
   proving that $Z(t) \in C\left( [-T_{{\rm loc}}, T_{{\rm loc}}];\ B_{\s_0}(C \tr) \cap {H^{\s-1}_\R(\T;\C^2)}\right)$, and therefore $t\mapsto Z(t)$ is weakly continuous in $H^{\s}_\R(\T;\C^2)$. Together with the continuity of $t\mapsto \norm{Z(t)}_\sigma$, a classical functional analysis argument gives 
$
Z(t) \in C\big([-T_{\rm loc},T_{\rm loc}];\, B_{\s_0}(C\,\tr) \cap H^\s_\R(\T;\C^2)\big)$, proving \Cref{loc.ex}.
\end{proof}

\section{Cubic Lifespan of Water-Waves with Constant Vorticity}\label{sec:cubic}
In this appendix, we apply the quadratic normal form result of \Cref{thm:quadratic.nf} to establish a cubic lifespan result for small and smooth initial data.
As a consequence, small solutions of size $\varepsilon$ remain small for long times $T \sim \varepsilon^{-2}$, as originally proved by Ifrim–Tataru \cite{IFT}. Since we do not keep track of the explicit Sobolev regularity threshold for the solutions, we recover the cubic lifespan of Ifrim–Tataru only for sufficiently regular initial data.
Precisely, recalling that for any $s>0$ the space $(X^s,\ \|\cdot\|_{X^s})$ is the Banach space defined by the norm in \eqref{norm:Xs}, we have the following result. 

\begin{proposition}\label{prop:e-2}
Fix an arbitrary  vorticity  $\gamma \in \R$.
    There is $\s_0>0$ such that for any $s\geq \s_0$ there is  $ \epsilon_0 = \epsilon_0(s)>0$
    such that for any initial datum
    \begin{align}
        (\eta_0,\psi_0, \sfV_0, \sfB_0) \in X^s \ \ \quad \text{with} \quad \varepsilon:= \|(\eta_0,\psi_0, \sfV_0, \sfB_0)\|_{X^{\s_0}}\leq \varepsilon_0\,,
        \label{initial:cond}
    \end{align}
where $ \sfV_0:= \sfV(\eta_0,\psi_0)$ and $ \sfB_0:= \sfB(\eta_0,\psi_0)$ are defined as in \eqref{def:V}-\eqref{form-of-B},  there is a unique classical solution of \eqref{eq:etapsi}
    \begin{align}
        (\eta(t),\psi(t), \sfV(t), \sfB(t))  \in C([-T_\varepsilon,T_\varepsilon];X^s), \quad \text{with} \quad T_\varepsilon \gtrsim \varepsilon^{-2}.
    \end{align}
    Moreover, such solution satisfies the bound 
    \begin{align}
       \sup_{t\in [0,T_\varepsilon]} \|(\eta(t),\psi(t), \sfV(t), \sfB(t))\|_{X^s}\lesssim_s 
    \, 
        \|(\eta_0,\psi_0, \sfV_0, \sfB_0)\|_{X^s}.
     \end{align}
\end{proposition}
\begin{proof}
The proof proceeds via a bootstrap argument. In addition to the local Cauchy theory, the main ingredient is a cubic energy estimate, which we derive as a consequence of \Cref{thm:quadratic.nf}.

\paragraph{Step 0: The local Cauchy theory.} By the local Cauchy theory (cf. the footnote in the proof of \Cref{loc.ex} at pag. \pageref{proof.Lemma.loc.ex}), given a small  initial datum as in \eqref{initial:cond} there are a time $T_{{\rm loc}}>0$ and a unique solution of \eqref{eq:etapsi} of the form
 $(\eta(t),\psi(t), \sfV(t), \sfB(t))\in C\left([-T_{{\rm loc}}, T_{{\rm loc}}];{X^s}\right)$.
\paragraph{Step 1: The bootstrap argument.} We fix 
\begin{align}
    \varepsilon:= \|(\eta_0,\psi_0, \sfV_0, \sfB_0)\|_{X^{\s_0}}\,, \qquad 
    \nu:= \|(\eta_0,\psi_0, \sfV_0, \sfB_0)\|_{X^{s}}\,.
\end{align}
We  prove the following
\begin{lemma}[{\bf Bootstrap}]\label{lem:boot}
    There is $\s_0>0$ such that for any $s\geq \s_0$ there are $c, \varepsilon_0>0$ and $\tK>1$ such that  for any $\varepsilon \in (0, \varepsilon_0)$, the following holds true: 
    if $T\leq c \varepsilon^{-2} $ and 
    \begin{align}
        \sup_{t\in [0,T]} \|(\eta(t),\psi(t), \sfV(t), \sfB(t))\|_{X^{\s_0}}\leq \tK \varepsilon \,, \qquad \sup_{t\in [0,T]} \|(\eta(t),\psi(t), \sfV(t), \sfB(t))\|_{X^{s}}\leq \tK \nu\,,
        \label{boot_ass}
    \end{align}
    then one has the improved bounds
    \begin{align}\label{improved}
        \sup_{t\in [0,T]} \|(\eta(t),\psi(t), \sfV(t), \sfB(t))\|_{X^{\s_0}}\leq \tfrac{\tK}{2} \varepsilon \,, \qquad \sup_{t\in [0,T]} \|(\eta(t),\psi(t), \sfV(t), \sfB(t))\|_{X^{s}}\leq \tfrac{\tK}{2} \nu \,. 
    \end{align}
\end{lemma}

\paragraph{Step 2: From the space $X^s$ to the good unknown of Alinhac.} 
We start by controlling  the $\dot H^{s+1}$-norm of the  Alinhac good unknown $\upomega= \psi- \Opbw{\sfB}\eta$. 
By \eqref{stima_omega} and the bootstrap assumption \eqref{boot_ass}, we have    $(\eta(t),\upomega(t)) \in H_0^{s+\frac12}(\T;\R)\times \dot H^{s+1}(\T;\R)$ and 
\begin{gather}
\label{stima_etaomega}
    \| \eta_0\|_{s+\frac12}+ \| \upomega_0\|_{s+1}\lesssim_s (1+\varepsilon)\nu \lesssim_s \nu, \qquad \sup_{t \in [0,T]}\| \eta(t)\|_{s+\frac12}+ \| \upomega(t)\|_{s+1}\lesssim_s \tK \nu  \,, \\
\label{stima_eta_omega_bassa}
    \| \eta_0\|_{\s_0+\frac12}+ \| \upomega_0\|_{\s_0+1}\lesssim \varepsilon \,,  
\end{gather}
choosing $\s_0 > s_0+1$
 and $\tK \varepsilon < 1$.

We now pass to the complex variable $ U$ in \eqref{Ugrande}.
By the equivalence of norms in \eqref{eq:etaomegaUequiv}, we have
\begin{align}
\label{u.e.om.boot}
 &    \| U(t)\|_{s+\frac34} \sim_s \| \eta(t)\|_{s+\frac{1}{2}} + \| \upomega(t)\|_{s+1} \stackrel{\eqref{stima_etaomega}}{\lesssim_s} \tK \nu, 
    \qquad \text{and }\  \| U(0)\|_{s+\frac34}\stackrel{\eqref{stima_etaomega}}{\lesssim_s} 
    \nu \,,\\
  &   \| U(0)\|_{\s_0+\frac34} \sim_s \| \eta_0\|_{\s_0+\frac{1}{2}} + \| \upomega_0\|_{\s_0+1} \stackrel{\eqref{stima_eta_omega_bassa}}{\lesssim }
    \varepsilon \,. 
      \label{stimaUdelta}
\end{align}
Consider now  the complex variable $\zak$ in \eqref{zak}: by the bootstrap assumption \eqref{boot_ass} one has  
\begin{align}
\sup_{t \in [0, T]}     \| \zak(t)\|_{\s_0}\lesssim \sup_{t \in [0, T]} \left(  \| \eta(t)\|_{\s_0-\frac14}+ \| \psi(t)\|_{\s_0+\frac14} \right) 
\lesssim  \tK \varepsilon\,.
    \label{zakeps}
\end{align}
In addition, by \Cref{complex:good}, one has also the  bounds
 \be\label{zak_and_U}
 \norm{\zak(t)}_{s} \lesssim_s \norm{U(t)}_{s  + \frac12}\lesssim_s \norm{\zak(t)}_{s  +1}  \,, \qquad \forall \,t \in[0, T]\,,
 \ee
 and arguing as in the proof of \Cref{UZequiv}, one obtains the following lemma.
 \begin{lemma}
Let $ Y=\bT(\zak)U$ be the the variable of  \Cref{thm:quadratic.nf}. For any $\s \in [\s_0, s + \frac34)$ and $\s_0$ sufficiently large, 
\be\label{BISMA_Y_U}
\norm{U(t)}_\s \sim_\s \norm{Y(t)}_\s \,, \quad \forall \, t \in [0, T] \,.  
\ee
\end{lemma}
\paragraph{Step 3: The energy estimate.} Using that $Y(t)$ solves the cubic equation \eqref{BNF1}, we prove the following 
\begin{lemma}
    Under the bootstrap assumption \eqref{boot_ass}, for any $\s\in [s_0, s+\frac34]$, we have 
    \begin{align}
        \Big|\frac{\di}{\di t}\| Y(t)\|_{\s}^2\Big|\lesssim_s \tK^2 \varepsilon^2 \| Y(t)\|_\s^2 \,,  \qquad 
        \forall \,t \in [0, T] \,. 
        \label{en:est}
    \end{align}
    \begin{proof}
        We differentiate the $H^\s$ norm using equation \eqref{BNF1} for $Y(t)$. We get 
        \begin{align*}
            \Big|\frac{\di}{\di t}\| Y(t)\|_{\s}^2\Big| =  \Big|2 \Re\la| D|^{2\s} \pa_t Y , Y  \ra\Big|\leq  & \Big|\la\Big[| D|^{2\s}, \vOpbw {\ii \sfD_{\geq 2}(\zak; x,\x)}\Big]Y, Y \ra\Big| \\
            &+ 2 \Big|\la | D|^{2\s} \left(\bB_{\geq N}(\zak) Y+  \bR_{\geq 2}(\zak)[Y]\right) , Y  \ra\Big|\\
            &\lesssim_s \tK^2 \varepsilon^2\| Y(t)\|_\s^2\,
        \end{align*}
        where, in addition to the bounds \eqref{zakeps}, in the first line we used that   $\sfD_{\geq 2}(\zak; x,\x)$ in \eqref{def:D}  is  a real-valued symbol in $\Sigma \Gamma^1_2[r,N]$ and         
 to estimate the second line       that $\bB_{\geq N}(\zak)$, being a spectrally localized map in  $\cS_{\geq N}^0[r]$,  fulfills the bounds \eqref{bound:specloc}, and that  $\bR_{\geq 2}(\zak)$ is a matrix of smoothing operators fulfilling \eqref{piove} with $m =0$ and $p=2$.
        This  proves \eqref{en:est}.
    \end{proof}
\end{lemma}
 Using repeatedly \eqref{BISMA_Y_U} and \eqref{en:est}, we are then able to
 perform the energy estimate
  \begin{align}
    \|U(t)\|_{ \s}^2\lesssim_\s \|Y(t)\|_{\s}^2\lesssim_\s  \|Y(0)\|_{\s}^2+ \int_0^T \Big|\frac{\di}{\di t}\| Y(t)\|_{\s}^2\Big|\, \di t 
    \lesssim_s \|U(0)\|_{\s}^2 + T \tK^2 \varepsilon^2 \|U(t)\|_{\s}^2\,. 
    \label{bobob}
\end{align}
 In conclusion, combining \eqref{bobob}, \eqref{u.e.om.boot}, and \eqref{stimaUdelta}, there are $C, C_s>0$ such that for any $t \in [0, T]$
\be\label{improved1}
\begin{aligned}
  & \| \eta(t)\|_{\s_0+\frac{1}{2}} + \| \upomega(t)\|_{\s_0+1}  \leq C \| U(t)\|_{\s_0+\frac34} \leq  C\left(\varepsilon^2 + T\tK^4 \varepsilon^4\right)^{\frac12}\,,  \\
  & \| \eta(t)\|_{s+\frac{1}{2}} + \| \upomega(t)\|_{s+1} \leq C_s  \| U(t)\|_{s+\frac34}  \leq  C_s\left(\nu^2 + T\tK^4 \varepsilon^2 \nu^2\right)^\frac12\,.
\end{aligned}
\ee
We are ready to prove the improved bounds \eqref{improved}. 
 By estimates \eqref{stima_V}, \eqref{stima_B}, \eqref{improved1} and by the  first bootstrap assumption \eqref{boot_ass} we have that for any $t \in [0, T]$
\begin{gather*}
    \| (\eta(t),\psi(t),\sfV(t),\sfB(t))\|_{X^s}\leq C_s\left( \| \eta(t)\|_{s+\frac12}+ \|\upomega(t)\|_{s+1} + \tK^2 \varepsilon \nu\right) \leq \und{C}_{\,s}\nu \left( 1+ \tK^2 \varepsilon T^{\frac12} +    \tK^2 \varepsilon\right)\leq \frac{\tK}{2} \nu\,,\\
    \| (\eta(t),\psi(t),\sfV(t),\sfB(t))\|_{X^{\s_0}}\leq C\left( \| \eta(t)\|_{\s_0+\frac12}+ \|\upomega(t)\|_{\s_0+1} + \tK^2 \varepsilon^2\right) \leq \und{C}\varepsilon \left( 1+ \tK^2 \varepsilon T^\frac12 +    \tK^2 \varepsilon\right)\leq \frac{\tK}{2} \varepsilon\,,
\end{gather*}
provided 
$\tK\geq 4 \max(\und{C}_{\,s},\ \und{C})$, $ \varepsilon\leq \varepsilon_0 \leq  (2 \tK^2)^{-1}$      and  $T\leq (4 \tK^4 \varepsilon^2)^{-1}$.
\\
Finally, the proof follows by a standard continuity argument, based on the bootstrap \Cref{lem:boot}.
\end{proof}

\footnotesize{
\noindent{\bf Acknowledgments:}
We thank Dario Bambusi, Ricardo Grande, Michela Procesi and Massimiliano Berti for useful comments.
 A. Maspero and S. Terracina are supported by the European Union  ERC CONSOLIDATOR GRANT 2023 GUnDHam, Project Number: 101124921.
 F. Murgante
is supported by the ERC STARTING GRANT 2021 
`Hamiltonian Dynamics, 
Normal Forms and Water Waves'' (HamDyWWa), Project Number: 101039762. 
 Views and opinions expressed are however those of the authors only and do not necessarily reflect those of the European Union or the European Research Council. Neither the European Union nor the granting authority can be held responsible for them.
 S. Terracina is also supported by GNAMPA-INdAM project "Growth of Sobolev norms and energy cascades in dispersive Hamiltonian PDE" CUP E53C25002010001
}


\end{document}